%% file: livre.tex
\documentclass[10pt,a4paper,reqno]{smfbook}
\usepackage{amsmath,amssymb,smfthm,stmaryrd,epic}
\usepackage[frenchb]{babel}
\usepackage[latin1]{inputenc}
\usepackage[all]{xy}

\usepackage{a4wide}

\usepackage{xypic}
\usepackage{rotating}

\NumberTheoremsIn{subsection}\SwapTheoremNumbers

\begin{document}

\input{macros2.tex}

\newcommand{\marque}{\addtocounter{smfthm}{1}
{\smallskip \noindent \textit{\thesmfthm}~---~}}

\renewcommand\atop[2]{\ensuremath{\genfrac..{0pt}{1}{#1}{#2}}}

\def\thechapter{\Roman{chapter}}

\newcommand{\clearemptydoublepage}{\newpage{\thispagestyle{empty} \cleardoublepage}}
\setcounter{secnumdepth}{3} \setcounter{tocdepth}{3}

\title{Faisceaux pervers des cycles évanescents des variétés de Drinfeld et groupes de
cohomologies du modèle de
Deligne-Carayol}

\alttitle{Perverse sheaf of vanishing cycles of Drinfeld varieties and cohomology group
of Deligne-Carayol
model}

\author{Boyer Pascal}
\email{boyer@math.jussieu.fr}
\address{Institut de mathématiques de Jussieu \\ UMR 7586, université Paris 6 \\
175 rue du Chevaleret Paris 13}

\frontmatter

\begin{abstract} Dans la première moitié du livre, on traduit, dans la situation
géométrique des variétés de Drinfeld,
les principaux résultats du livre d'Harris et Taylor. On explicite notamment la
restriction aux strates
ouvertes des faisceaux des cycles évanescents en fonction de certains systèmes locaux
dits d'Harris-Taylor
dont on calcule la somme alternée des groupes de cohomologie à supports compacts. Dans
la deuxième moitié du
livre, on décrit les gradués de la filtration de monodromie du faisceau pervers des
cycles évanescents ainsi
que la suite spectrale correspondante. D'après le théorème de comparaison de
Berkovich-Fargues, on obtient
alors une description de la filtration de monodromie-locale du modèle de
Deligne-Carayol.

\end{abstract}

\begin{altabstract} In the first half of the book, we translate in the geometric
situation of Drinfeld varieties, the principal results of the Harris and Taylor's book. We give in particular
the restriction to the open strata of the vanishing cycles sheaves in terms of some local systems named
Harris-Taylor's local systems which we calculate the alternated sum of the cohomology group with compact
supports. In the last half of the book, we describe the monodromy filtration of the vanishing cycles perverse
sheaf and the spectral sequence associated to it. Thanks to the Berkovich-Fargues' theorem, we obtain the
description of the local monodromy filtration of the Deligne-Carayol model.

\end{altabstract}



\altkeywords{Drinfeld varieties, formal modules, Langlands correspondences,
Jacquet-Langlands
correspondences, monodromy filtration, weight-monodromy conjecture, Shimura varieties,
perverse sheaves,
vanishing cycles}

\dedicatory{A mes parents, à Olivia et à Nathan}


\maketitle

\tableofcontents


\pagestyle{headings} \pagenumbering{arabic}


\mainmatter \clearemptydoublepage

\chapter*{Introduction}


\input{introduction-livre.tex}

\bigskip

\noindent \textit{Remerciements} Je tiens à exprimer ma profonde gratitude envers Gérard
Laumon pour son
soutien constant tout au long de ces années ainsi que pour sa perspicacité à dépister
les fausses bonnes
idées et à mettre en avant les autres. Merci aussi à Laurent Fargues de m'avoir fourni
une bonne notion de
monodromie locale qui m'a permis de simplifier grandement les preuves. J'adresse des
remerciements à
Jean-François Dat pour m'avoir fourni l'habillage catégoriel des faisceaux de Hecke
ainsi que pour ses
nombreux conseils ainsi qu'à Thomas Hausberger pour m'avoir indiqué une erreur dans mes
précédents travaux au
niveau des actions sur le modèle local.

\clearemptydoublepage

\chapter[Variétés d'Igusa]{Variétés d'Igusa, systèmes locaux d'Harris-Taylor et cycles
évanescents des variétés de Drinfeld}

\section*{Introduction}

\input{introduction1.tex}

\input{complements.tex}

\input{igusa1.tex}

\input{cycles.tex}



\input{calcul.tex}

\clearemptydoublepage

\chapter[Cas Iwahori]{Groupes de cohomologie du modèle local: cas Iwahori}

\section*{Introduction}

\input{introduction-chap2.tex}

\input{iwahori.tex}

\clearemptydoublepage

\chapter{Cohomologie des faisceaux de Harris-Taylor}

\section*{Introduction}

\input{introduction-chap3.tex}

\input{f-traces.tex}

\input{somme-alternee.tex}

\clearemptydoublepage

\chapter{Filtration de monodromie des cycles évanescents}

\section*{Introduction}


\input{introduction2.tex}


\input{enonces.tex}

\input{preuve.tex}

\clearemptydoublepage

\chapter[Compléments et applications]{Compléments sur la cohomologie globale et
applications}

\section*{Introduction}


\input{introduction-chap5.tex}

\input{applications.tex}



\clearemptydoublepage

\appendix


\clearemptydoublepage


\input{appendiceB.tex}


\clearemptydoublepage



\chapter{Récapitulatif des calculs}


\input{recap.tex}


\clearemptydoublepage

\chapter{Figures}




\input{figure7.tex}

\clearpage
\newpage

\input{figure8.tex}

\clearpage
\newpage

\textit{Légende pour les figures (\ref{figure5}) et (\ref{figure6}):}
\begin{itemize}

\item on indique les coordonnées $(p,q)$ telles que $E_2^{p,q} [\Pi^{\oo,o}]$ est non
nulle avec $\Pi$ une
représentation irréductible automorphe de $D_\Am^\times$ telle que $\Pi_\oo \simeq
\st_4(\pi_\oo)$ (resp.
$\speh_4(\pi_\oo)$) dans la figure 5 (resp. figure 6) pour $\pi_\oo$ une représentation
cuspidale de
$GL_2(F_\oo)$, et $\Pi_o \simeq \st_4(\pi_o)$ (resp. $\speh_4 (\pi_o)$) avec $\pi_o$
cuspidale de
$GL_2(F_o)$. Les représentations obtenues sont des représentations elliptiques de type
$\pi_o$. Pour chacune
de ces coordonnées, on indique qu'elle est la strate qui la donne dans la suite
spectrale de stratification
correspondante;

\item le nombre en dessous des $\overleftrightarrow{3}$ indique le poids;

\item toutes les flèches indiquées sont les seules non triviales et découlent des suites
exactes (\ref{suites-exactes});

\item on note avec un cercle plus grand les coordonnées qui contribuent à
l'aboutissement de la suite spectrale.
\end{itemize}

\begin{figure}
\input{figure5.tex} \caption{\label{figure5} $s=4$, $g=2$; $E_2^{p,q}[\Pi^{\oo,o}]$ de la suite spectrale des
cycles évanescents: cas $\Pi_o=\st_4(\pi_o)$. Le dessin du bas représente le terme $E_3$ et donc
l'aboutissement de la suite spectrale.}
\end{figure}

\newpage

\begin{figure}
\input{figure6.tex} \caption{\label{figure6}  $s=4$, $g=2$: $E_2^{p,q}[\Pi^{\oo,o}]$ de la suite spectrale des
cycles évanescents: cas $\Pi_o=\speh_s(\pi_o)$. Le dessin du bas représente le terme $E_3$ et donc
l'aboutissement de la suite spectrale.}
\end{figure}

\newpage

\clearpage

\textit{Légende pour les figures (\ref{figure1}), (\ref{figure2})}: le nombre qui accompagne les
représentations indique le décalage à donner à la représentation de Galois associée; ainsi quand, dans les
$E_2^{p,q}$ des suites spectrales (\ref{ss-dualite}) pour $\PC(g,l,\pi_o,k)$, on écrit
$\genfrac{}{}{0pt}{}{\overleftrightarrow{s-1}}{r}$, il faut lire
$$j^{\geq lg}_! HT(g,l,\pi_o,[\overleftrightarrow{s-1}]_{\pi_o}) \otimes
L_g(\pi_o)(-\frac{i(g-1)+l-1-k+r}{2}).$$

\bigskip

\begin{figure}[!h]
\input{figure1.tex} \caption{\label{figure1} $g=7$, $s=5$; $E_2^{p,q}$ de la suite spectrale
(\ref{ss-dualite}) pour $l=4$}
\end{figure}

\begin{figure}[!h]
\input{figure2.tex} \caption{\label{figure2}  $g=7$, $s=5$; $E_2^{p,q}$ de la suite spectrale
(\ref{ss-dualite}) pour $l=3$}
\end{figure}






\clearpage
\newpage

\clearemptydoublepage

\backmatter

\bibliographystyle{plain}
\bibliography{bib-ok}

\end{document}

%% file: macros2.tex
\newtheorem{prop-defi}[smfthm]{Proposition-Définition}
\newtheorem{theo-defi}[smfthm]{Théorème-définition}
\newtheorem{lem-defi}[smfthm]{Lemme-définition}
\newtheorem{notas}[smfthm]{Notations}
\newtheorem{nota}[smfthm]{Notation}
\newtheorem{defis}[smfthm]{Définitions}
\newtheorem{remas}[smfthm]{Remarques}

\newtheorem{theob}{Théorème}[section]
\def\thetheob{\Roman{chapter}.\arabic{section}.\arabic{theob}}
\newtheorem{propb}[theob]{Proposition}
\newtheorem{lemb}[theob]{Lemme}
\newtheorem{corob}[theob]{Corollaire}
\newtheorem{defib}[theob]{Définition}
\newtheorem{defisb}[theob]{Définitions}
\newtheorem{remab}[theob]{Remarque}

\renewcommand{\theequation}{\Roman{chapter}.\arabic{section}.\arabic{subsection}.\arabic{equation}}

\def\Am{{\mathbb A}}
\def\Fm{{\mathbb F}}
\def\Mm{{\mathbb M}}
\def\Nm{{\mathbb N}}
\def\Pm{{\mathbb P}}
\def\Qm{{\mathbb Q}}
\def\Zm{{\mathbb Z}}
\def\Dm{{\mathbb D}}
\def\Cm{{\mathbb C}}
\def\Rm{{\mathbb R}}

\def\AC{{\mathcal A}}
\def\CC{{\mathcal C}}
\def\DC{{\mathcal D}}
\def\EC{{\mathcal E}}
\def\FC{{\mathcal F}}
\def\GC{{\mathcal G}}
\def\HC{{\mathcal H}}
\def\IC{{\mathcal I}}
\def\JC{{\mathcal J}}
\def\KC{{\mathcal K}}
\def\LC{{\mathcal L}}
\def\MC{{\mathcal M}}
\def\NC{{\mathcal N}}
\def\OC{{\mathcal O}}
\def\PC{{\mathcal P}}
\def\UC{{\mathcal U}}
\def\VC{{\mathcal V}}
\def\XC{{\mathcal X}}
\def\YC{{\mathcal Y}}

\def\AF{{\mathfrak A}}
\def\GF{{\mathfrak G}}
\def\EF{{\mathfrak E}}
\def\CF{{\mathfrak C}}
\def\DF{{\mathfrak D}}
\def\JF{{\mathfrak J}}
\def\LF{{\mathfrak L}}
\def\MF{{\mathfrak M}}
\def\NF{{\mathfrak N}}
\def\XF{{\mathfrak X}}
\def\UF{{\mathfrak U}}

\def \longmapright#1{\smash{\mathop{\longrightarrow}\limits^{#1}}}
\def \mapright#1{\smash{\mathop{\rightarrow}\limits^{#1}}}
\def \lexp#1#2{\kern \scriptspace \vphantom{#2}^{#1}\kern-\scriptspace#2}
\def \linf#1#2{\kern \scriptspace \vphantom{#2}_{#1}\kern-\scriptspace#2}
\def \linexp#1#2#3 {\kern \scriptspace{#3}_{#1}^{#2} \kern-\scriptspace #3}

\def \a{\alpha}
\def \b{\beta}
\def \d{\delta}
\def \e{\epsilon}
\def \g{\gamma}
\def \k{\kappa}
\def \l{\lambda}
\def \m{\mu}
\def \n{\nu}
\def \o{\omega}
\def \r{\rho}
\def \s{\sigma}
\def \t{\tau}
\def \th{\theta}
\def \u {\upsilon}
\def \x{\chi}
\def \z {\dzeta}
\def \vphi {\varphi}

\let \leq=\leqslant
\let \geq=\geqslant
\def \lefto{\longleftarrow}
\def \fin{\hfill $\square$}
\let \DS=\displaystyle
\let \SS=\scriptstyle
\let \longto=\longrightarrow
\let \oo=\infty

\def \FH{\mathop{\mathrm{FH}}\nolimits}
\def \FPH{\mathop{\mathrm{FPH}}\nolimits}
\def \coh{\mathop{\mathrm{Coh}}\nolimits}
\def \res{\mathop{\mathrm{res}}\nolimits}
\def \op{\mathop{\mathrm{op}}\nolimits}
\def \rec {\mathop{\mathrm{rec}}\nolimits}
\def \art{\mathop{\mathrm{Art}}\nolimits}
\def \hyp {\mathop{\mathrm{Hyp}}\nolimits}
\def \cusp {\mathop{\mathrm{Cusp}}\nolimits}
\def \Iw {\mathop{\mathrm{Iw}}\nolimits}
\def \JL {\mathop{\mathrm{JL}}\nolimits}
\def \speh {\mathop{\mathrm{Speh}}\nolimits}
\def \isom {\mathop{\mathrm{Isom}}\nolimits}
\def \Vect {\mathop{\mathrm{Vect}}\nolimits}
\def \groth {\mathop{\mathrm{Groth}}\nolimits}
\def \lef {\mathop{\mathrm{Lef}}\nolimits}
\def \fix {\mathop{\mathrm{Fix}}\nolimits}
\def \hom {\mathop{\mathrm{Hom}}\nolimits}
\def \deg {\mathop{\mathrm{deg}}\nolimits}
\def \val {\mathop{\mathrm{val}}\nolimits}
\def \det {\mathop{\mathrm{det}}\nolimits}
\def \rep {\mathop{\mathrm{Rep}}\nolimits}
\def \spec {\mathop{\mathrm{Spec}}\nolimits}
\def \fr {\mathop{\mathrm{Fr}}\nolimits}
\def \frob {\mathop{\mathrm{Frob}}\nolimits}
\def \ker {\mathop{\mathrm{Ker}}\nolimits}
\def \im {\mathop{\mathrm{Im}}\nolimits}
\def \Red {\mathop{\mathrm{Red}}\nolimits}
\def \red {\mathop{\mathrm{red}}\nolimits}
\def \aut {\mathop{\mathrm{Aut}}\nolimits}
\def \diag {\mathop{\mathrm{diag}}\nolimits}
\def \spf {\mathop{\mathrm{Spf}}\nolimits}
\def \Def {\mathop{\mathrm{Def}}\nolimits}
\def \twist {\mathop{\mathrm{Twist}}\nolimits}
\def \supp {\mathop{\mathrm{Supp}}\nolimits}
\def \Id {{\mathop{\mathrm{Id}}\nolimits}}
\def \bar {\overline}
\def \ind {\mathop{\mathrm{Ind}}\nolimits}
\def \mod {\mathop{\mathrm{mod}}\nolimits}
\def \ker {\mathop{\mathrm{Ker}}\nolimits}
\def \coker {\mathop{\mathrm{Coker}}\nolimits}
\def \mult {\mathop{\mathrm{mult}}\nolimits}
\def \vide{\emptyset}
\def \bad {{\mathop{\mathrm{Bad}}\nolimits}}
\def \gal {{\mathop{\mathrm{Gal}}\nolimits}}
\def \Nr {{\mathop{\mathrm{Nr}}\nolimits}}
\def \rn {{\mathop{\mathrm{rn}}\nolimits}}
\def \vol {{\mathop{\mathrm{vol}}\nolimits}}
\def \ad {{\mathop{\mathrm{ad}}\nolimits}}
\def \tr {{\mathop{\mathrm{Tr~}}\nolimits}}
\def \Sp {{\mathop{\mathrm{Sp}}\nolimits}}
\def \lie {{\mathop{\mathrm{Lie}}\nolimits}}
\def \st {{\mathop{\mathrm{St}}\nolimits}}
\def \sp{{\mathop{\mathrm{Sp}}\nolimits}}
\def \card{{\mathop{\mathrm{card}}\nolimits}}
\def \sym{{\mathop{\mathrm{Sym}}\nolimits}}
\def \perv{\mathop{\mathrm{Perv}}\nolimits}
\def \const{\mathop{\mathrm{Const}}\nolimits}

\def \ele{élément }
\def \eles{éléments }
\def \cad{c'est à dire }
\def \rem{{\noindent\textit{Remarque:~}}}
\def \exem{{\noindent \textit{Exemple:~}}}
\def \ssi{~si et seulement si~}
\def \cl {{\mathop{\mathrm{cl}}\nolimits}}
\def \Tw {{\mathop{\mathrm{Tw}}\nolimits}}
\def \ob {{\mathop{\mathrm{Ob}}\nolimits}}
\def \ext {{\mathop{\mathrm{Ext}}\nolimits}}
\def \End {{\mathop{\mathrm{End}}\nolimits}}
\def \inv {{\mathop{\mathrm{inv}}\nolimits}}
\def \fix {{\mathop{\mathrm{Fix}}\nolimits}}

\def\semi{\mathrel{>\!\!\!\triangleleft}}


%% file: introduction-livre.tex
\noindent \textbf{0.1.} --- Soient $d$ un entier strictement positif et $K$ un corps local complet d'égale
caractéristique $p$, d'anneau des entiers $\OC_K$. On considère le groupe $D_{K,d}^\times$ (resp. $W_K$) des
éléments inversibles de ``l'''algèbre à division centrale sur $K$ d'invariant $1/d$ (resp. le groupe de Weil
de $K$). Pour un nombre premier $l \neq p$, Langlands (resp. Jacquet-Langlands) a (resp. ont) conjecturé
l'existence d'une bijection $L_d$ (resp. d'une injection $\JL$) entre les $\bar \Qm_l$-représentations
irréductibles admissibles de $GL_d(K)$ et les représentations $l$-adiques indécomposables de $W_K$ (resp.
entre les représentations admissibles irréductibles de $D_{K,d}^\times$ et les représentations
essentiellement de carré intégrable de $GL_d(K)$) qui sont compatibles à la formation des fonctions $L$ de
paires; on renvoie à \cite{he} pour des énoncés précis.

A l'aide de la cohomologie étale, Deligne a alors construit une série de représentations $\UC_{K}^{d,i}$ du
produit de ces trois groupes. Pour $d=2$ et pour $\rho$ une représentation irréductible de $D_{K,2}^\times$
telle que $\pi:=\JL(\rho)$ soit une représentation cuspidale de $GL_2(K)$, Carayol, dans \cite{ca}, montre
que la composante $\rho$-isotypique $\UC_{K}^{2,1}(\rho)$ de $\UC_{K}^{2,1}$ réalise les correspondances de
Langlands et de Jacquet-Langlands, i.e.
$$\UC_{K}^{2,1}(\JL^{-1}(\pi^\vee)) \simeq \pi \otimes L_d(\pi)(-\frac{1}{2}).$$
Le cas $d$ quelconque est traité dans \cite{boy}. En outre pour $d=2$, Carayol décrit également ce qui se
passe pour les autres représentations.

Le but premier de ce travail est de faire de même pour $d$ quelconque, i.e. calculer complètement les
$\UC_{K}^{d,i}(\rho)$ pour $\rho$ une représentation irréductible quelconque de $D_{K,d}^\times$. Dans le cas
où $\rho$ est la représentation triviale, le résultat se formule comme suit.

\medskip

\noindent \textbf{Théorème} \textit{Pour $0 \leq i \leq d-1$, on a
$$\UC_{K}^{d,i}(1)=\pi_{i} \otimes 1(-i)$$
où $\pi_{i}$ est l'unique quotient irréductible de l'induite parabolique
$$\ind_{P_{d-i,d}(K)}^{GL_{sg}(K)} 1(\frac{i(g-1)}{2}) \otimes \st_{i}(-\frac{(s-i)(g-1)}{2})$$
où $P_{d-i,d}$ est le parabolique standard associée aux $d-i$ premiers vecteurs et $\st_i$ est la
représentation de Steinberg de $GL_i(K)$.}

\medskip

L'énoncé du cas général, théorème (\ref{theo-ripsi-local}), s'énonce de manière similaire et fait intervenir
les correspondances de Langlands et Jacquet-Langlands. En particulier dans le cas Iwahori, via l'isomorphisme
de Faltings, cf. \cite{falt}, on retrouve le résultat principal de \cite{s-s}.

\medskip

\noindent \textbf{0.2.} --- La preuve est de nature globale et repose sur le théorème de comparaison de
Berkovich \footnote{en fait sur une version raffinée fournie par Fargues (cf. le théorème principal de
l'appendice de \cite{boy-duke})} des cycles évanescents locaux et globaux. Ainsi le deuxième résultat de ce
texte est la description explicite du complexe des cycles évanescents des variétés de Drinfeld-Stuhler ainsi
que de sa filtration de monodromie-poids et de la suite spectrale associée.

Par ailleurs les techniques s'appliquent dans le cadre de la caractéristique mixte pour les variétés de
Shimura associés à certains groupes unitaires étudiées dans \cite{h-t} ce qui fournit en particulier les
conjectures de monodromie-poids versions faisceautique et cohomologique. Ce cadre est l'objet de
\cite{boy-duke}, et nous donnons dans l'appendice A les grandes lignes des modifications formelles à
apporter.

\medskip

\noindent \textbf{0.3.} --- Soit $X$ une courbe projective lisse, irréductible et géométriquement connexe
définie sur "le" corps fini à $q=p^r$ \eles, $\Fm_q$ et soit $F$ son corps des fonctions. On fixe deux places
distinctes $\oo$ et $o$ de $X$ que l'on peut supposer par simplification, rationnelles sur $\Fm_q$ de sorte
que le complété $F_o$ du localisé de $F$ en $o$ est le corps local précédemment noté $K$. On note $A$
l'anneau des fonctions sur $X$, régulières en dehors de $\oo$. Étant donné un entier $d \geq 1$, on fixe une
algèbre à division centrale $D$ sur $F$ de dimension $d^2$, non ramifiée en $\oo$ et $o$, ainsi qu'un ordre
maximal $\DC$.

Dans \cite{lrs}, les auteurs construisent pour un idéal non trivial $I$ de $A$, un schéma $M_{I}$ défini sur
$F$, classifiant les $\DC$-faisceaux elliptiques sur $X$, munis d'une structure de niveau $I$. Pour $o \not
\in V(I)$, $M_I$ a un modèle entier $M_{I,o}$ lisse sur le complété $\OC_o$ de $A$ en la place $o$. Un tel
modèle non lisse dans le cas où $o \in V(I)$ est construit dans \cite{boy}. Les schémas $M_{I,o}$ sont
naturellement munis d'une action, par correspondances, de $(D_\Am^{\oo})^\times$.

\medskip

\noindent \textbf{0.4.} --- On s'intéresse alors à la fibre spéciale $M_{I,s_o}$ de $M_{I,o}$. Dans
\cite{boy}, je stratifie $M_{I,s_o}$ par des sous-schémas localement fermés $M_{I,s_o}^{=h}$ pour $1 \leq h
\leq d$, de pure dimension $d-h$ tels que l'on ait un équivalent du théorème de Serre-Tate pour les
$\DC$-faisceaux elliptiques à savoir: le complété de l'hensélisé strict de l'anneau local de $M_{I,o}$ en un
point géométrique de $M_{I,s_o}^{=h}$ est isomorphe à $\Def_n^h[[x_1,\cdots,x_{d-h}]]$ où $n$ est la
multiplicité de $o$ dans $I$ et $\Def_n^h$ représente le foncteur des déformations de niveau $n$ d'un
$\OC_o$-module formel de hauteur $h$ sur $\bar \Fm_p$.

Par ailleurs pour $1 \leq h < d$, il existe un sous-schéma fermé $M_{I,s_o,1}^{=h}$ de $M_{I,s_o}^{=h}$
stable sous les correspondances associées aux \eles du sous-groupe parabolique $P_{h,d}^{op}(F_o)$ de
$GL_d(F_o)$ (cf. la définition (\ref{defi-parab})) et tel que
$$M_{I,s_o}^{=h}= M_{I,s_o,1}^{=h} \times_{P_{h,d}^{op}(\OC_o/\MC_o^n)} GL_d(\OC_o/\MC_o^n)$$
où $n$ est la multiplicité de $o$ dans $I$: on dit que les strates non supersingulières sont géométriquement
induites.

\medskip

\noindent \textbf{0.5.} --- Dans le premier chapitre, en suivant \cite{h-t}, on étudie la restriction aux
strates $M_{I,s_o}^{=h}$ des faisceaux des cycles évanescents $R^i \Psi_{\eta_o}(\bar \Qm_l)$. Pour cela on
introduit selon loc. cit. les variétés d'Igusa de première et seconde espèce.

Les premières sont des revêtements galoisiens de $M_{I^o,s_o}^{=h}$ de groupe de Galois
$GL_{d-h}(\OC_o/\MC_o^n)$, définies par la donnée d'une structure de niveau $\MC_o^n$ sur la partie étale des
$\DC$-faisceaux elliptiques de $M_{I^o,s_o}^{=h}$: on dispose alors d'un morphisme radiciel de
$\IC_{I^o,n}^{=h}$ vers $M_{I,s_o,1}^{=h}$.

Les variétés d'Igusa de seconde espèce sont des revêtements galoisiens $\JC_{I^o,n}^{=h}(s) \longto
\IC_{I^o,n}^{=h}$ de groupe de Galois $(\DC_{o,h}/(\Pi_{o,h}^s))^\times$ où $\DC_{o,h}$ est l'ordre maximal
de l'algèbre à division $D_{o,h}$ centrale sur $F_o$ d'invariant $1/h$, avec $\Pi_{o,h}$ une uniformisante:
elles sont définies via une rigidification modulo $\Pi_{o,h}^s$ de la partie connexe des $\DC$-faisceaux
elliptiques de $\IC_{I^o,n}^{=h}$.

On dispose alors pour toute représentation $\tau$ irréductible et admissible de $D_{o,h}^\times$, d'un
système local $\FC_{\tau}$ dit de Harris-Taylor, sur $M_{I,s_o}^{=h}$. Suivant \cite{h-t}, la restriction à
$M_{I,s_o,1}^{=h}$ des faisceaux des cycles évanescents s'exprime alors en termes des systèmes locaux
d'Harris-Taylor et des faisceaux des cycles évanescents du modèle local.

\medskip

\noindent \textbf{0.6.} --- Au deuxième chapitre on prouve le cas Iwahori du théorème local. On commence par
des rappels sur les représentations elliptiques de $GL_d(F_o)$. On raisonne ensuite par récurrence en
supposant connu le cas Iwahori du théorème local en hauteur strictement inférieure à $d$. L'hypothèse de
récurrence et le fait que les strates non supersingulières soient induites, nous permettent alors d'étudier
combinatoirement la suite spectrale des cycles évanescents et de prouver le résultat.

\medskip

\noindent \textbf{0.7.} --- Au troisième chapitre, en suivant \cite{lrs}, il s'agit de calculer la somme
alternée des groupes de cohomologie des systèmes locaux d'Harris-Taylor. La démarche est classique: il s'agit
tout d'abord d'utiliser la formule des traces de Lefschetz et donc de compter les points fixes sous l'action
d'une correspondance de Hecke tordue par une puissance assez grande du Frobenius et ensuite de transférer les
intégrales orbitales obtenues afin de reconnaître le coté géométrique de la formule des traces de Selberg.

On en déduit alors un calcul de la somme alternée des groupes de cohomologie du modèle local de
Deligne-Carayol. En particulier, dans le cas Iwahori, des arguments de pureté nous redonnent les résultats
obtenus à la fin du chapitre précédent. Pour ce qui est du cas général, le théorème local correspond à dire
qu'il n'y a pas d'annulation dans la représentation virtuelle $\sum_{i=0}^{d-1} (-1)^i
[\UC_{F_o}^{d,i}(\rho_o)]$ où $\UC_{F_o}^{d,d-i}(\rho_o)$, pour $1 \leq i \leq d$, est donné par le $i$-ème
terme de plus haut poids.

\medskip

\noindent \textbf{0.8.} --- Dans le quatrième chapitre on étudie la filtration de monodromie du faisceau
pervers $R\Psi_{\eta_o}(\bar \Qm_l)[d-1]$ dont on notera $gr_k$ les gradués ainsi que la fibre en un point
supersingulier des termes $E_1^{i,j}$ de la suite spectrale
\begin{equation*}
E_1^{i,j}=h^{i+j} gr_{-i} \Rightarrow R^{i+j+d-1} \Psi_{\eta_o}(\bar{{\mathbb Q}_l})
\end{equation*}
et en particulier la fibre en un point supersingulier des termes $E_1^{i,j}$.

On décrit tout d'abord les $gr_k$ dans la catégorie des faisceaux pervers sur la tour des $(M_{I,s_o})_I$
munis d'une action par correspondances de $(D_\Am^\oo)^\times \times W_o$, en fonction des extension
intermédiaires $j^{\geq lg}_{!*} \FC_{\JL^{-1}(\st_l(\pi_o))}[d-lg]$ des systèmes locaux d'Harris-Taylor, où
$j^{\geq lg}$ désigne l'injection de la strate $M_{I,s_o}^{=lg}$ et $\pi_o$ est une représentation
irréductible cuspidale de $GL_g(F_o)$ avec $1 \leq l \leq d/g$.

On calcule ensuite les faisceaux de cohomologie de ces derniers ce qui donne le théorème local d'après le
théorème de comparaison de Berkovich. En fait Fargues, cf. l'appendice de \cite{boy-duke}, améliore les
résultats de Berkovich ce qui nous permet finalement de décrire la filtration de monodromie-locale du modèle
de Deligne-Carayol.

\medskip

\noindent \textbf{0.9.} --- Dans le dernier chapitre on donne divers compléments. On prouve en particulier la
conjecture (14.21) de \cite{lrs}, des cas particuliers de correspondances de Jacquet-Langlands par ailleurs
connus des experts. On montre aussi que les composantes locales de représentations automorphes de
$D_\Am^\times$ sont ce que l'on attend à savoir ce que l'on devrait obtenir à partir des résultats de
Moeglin-Waldspurger sur les composantes locales des carrés intégrables de $GL_d(\Am)$ via une correspondance
de Jacquet-Langlands globale. On montre enfin comment en caractéristique mixte, on peut prouver la conjecture
de monodromie-poids version cohomologique.


\medskip

\noindent \textbf{0.10.} --- Dans l'appendice A, on donne le dictionnaire entre les notations de \cite{h-t}
et les nôtres et on renvoie à loc. cit. pour les résultats géométriques et cohomologiques que l'on utilise.
Via ce dictionnaire, les résultats et les preuves du chapitre IV sont alors valables telles quelles; en
particulier on obtient une preuve de la conjecture de monodromie-poids dans le cadre des variétés de Shimura
de loc. cit.

\medskip

\noindent \textbf{0.11.} --- L'appendice B résume certains des résultats obtenus au cours de la preuve et
l'appendice C présente diverses figures des nombreuses suites spectrales que l'on utilise, essentiellement
pour $s=4$.

\bigskip

Les chapitres I et III sont essentiellement une traduction en égale caractéristique des résultats de
\cite{h-t}. Le chapitre II n'est pas utile mais il nous a semblé agréable de présenter une preuve du cas
Iwahori qui n'utilise pas toute la machinerie. Pour les lecteurs intéressés par les résultats du chapitre IV,
nous avons fait en sorte que celui-ci soit lisible indépendamment des deux premiers, sauf pour ce qui est des
rappels de \cite{ze} donnés \S \ref{rappel-ze-0}.


%% file: introduction1.tex
\noindent \textbf{0.1.} --- Rappelons la situation globale. Soit $X$ une courbe projective lisse,
irréductible et géométri\-quement connexe définie sur le corps fini à $q=p^r$ \eles, $\Fm_q$ et soit $F$ son
corps des fonctions. On fixe deux places $\oo,o$ distinctes de $X$ que l'on peut supposer par simplification,
rationnelles sur $\Fm_q$. On note $A$ l'anneau des fonctions sur $X$, régulières en dehors de $\oo$. Étant
donné un entier $d \geq 1$, on fixe une algèbre à division centrale $D$ sur $F$ de dimension $d^2$, non
ramifiée en $\oo$ et $o$, ainsi qu'un ordre maximal $\DC$.

Dans \cite{lrs}, les auteurs construisent pour un idéal non trivial $I$ de $A$, un schéma $M_{I}$ défini sur
$F$, classifiant les $\DC$-faisceaux elliptiques sur $X$, munis d'une structure de niveau $I$. Pour $o \not
\in V(I)$, $M_I$ a un modèle entier $M_{I,o}$ lisse sur le complété $\OC_o$ de $A$ en la place $o$. Dans
\cite{boy}, on construit un tel modèle non lisse dans le cas où $o \in V(I)$. Les schémas $M_{I,o}$ sont
naturellement munis d'une action, par correspondances, de $(D_\Am^{\oo})^\times$.

\medskip

\noindent \textbf{0.2.} --- Dans \cite{boy}, la fibre spéciale $M_{I,s_o}$ de $M_{I,o}$ est stratifiée par
des sous-schémas localement fermés $M_{I,s_o}^{=h}$ pour $1 \leq h \leq d$, de pure dimension $d-h$ tels que
l'on ait un équivalent du théorème de Serre-Tate pour les $\DC$-faisceaux elliptiques à savoir: le complété
de l'hensélisé strict de l'anneau local de $M_{I,o}$ en un point géométrique de $M_{I,s_o}^{=h}$ est un
anneau de séries formelles $\Def_n^h$ où ce dernier représente les déformations de niveau $n:=\mult_o(I)$
d'un $\OC_o$-module formel de hauteur $h$.

Par ailleurs ces strates pour $h \neq d$, sont géométriquement induites au sens où il existe un sous-schéma
fermé $M_{I,s_o,1}^{=h}$ de $M_{I,s_o}^{=h}$ stable sous les correspondances associées aux \eles du
sous-groupe parabolique $P_{h,d}^{op}(F_o)$ de $GL_d(F_o)$, opposé au parabolique standard associé aux $h$
premiers vecteurs de la base canonique et tel que
$$M_{I,s_o}^{=h}= M_{I,s_o,1}^{=h} \times_{P_{h,d}^{op}(\OC_o/\MC_o^n)} GL_d(\OC_o/\MC_o^n)$$

\medskip

\noindent \textbf{0.3.} --- Il s'avère que l'adhérence $M_{I,s_o,1}^{\geq h}$ des $M_{I,s_o,1}^{=h}$ est
lisse de sorte que la non lissité des strates $M_{I,s_o}^{=h}$ provient des intersections  entre ces
différentes composantes et donc de l'intrication des parties connexes et étales des $I$-structures de niveau.
Afin de démêler la situation on introduit les variétés d'Igusa de première et seconde espèce. Pour tout place
$o$ divisant $I$, on note $I=I^o \MC_o^n$ où $n$ est la multiplicité de $o$ dans $I$ de sorte que $o$ ne
divise pas $I^o$.

Les variétés d'Igusa de première espèce $\IC_{I^o,n}^{=h}$ sont des revêtements galoisiens de
$M_{I^o,s_o}^{=h}$ de groupe de Galois $GL_{d-h}(\OC_o/\MC_o^n)$, définis par la donnée d'une structure de
niveau $\MC_o^n$ sur la partie étale des $\DC$-faisceaux elliptiques de $M_{I^o,s_o}^{=h}$: on dispose alors
d'un morphisme radiciel de $\IC_{I^o,n}^{=h}$ vers $M_{I,s_o,1}^{=h}$.

Les variétés d'Igusa de seconde espèce sont des revêtements galoisiens $\JC_{I^o,n}^{=h}(s) \longto
\IC_{I^o,n}^{=h}$ de groupe de Galois $(\DC_{o,h}/(\Pi_{o,h}^s))^\times$ où $\DC_{o,h}$ est l'ordre maximal
de l'algèbre à division $D_{o,h}$ centrale sur $F_o$ d'invariant $1/h$, avec $\Pi_{o,h}$ une uniformisante:
elles sont définies via une rigidification modulo $\Pi_{o,h}^s$ de la partie connexe des $\DC$-faisceaux
elliptiques de $\IC_{I^o,n}^{=h}$.

On dispose alors pour toute représentation irréductible et admissible de $\DC_{o,h}^\times$, d'un système
local dit de Harris-Taylor, sur $M_{I,s_o,1}^{=h}$ qu'il est aisé de propager à toute la strate via l'action
de $GL_d(\OC_o)$. Afin de mieux suivre les actions, on partira plutôt d'une représentation irréductible
$\tau_o$ de $D_{o,h}^\times$ dont on prendra la restriction à $\DC_{o,h}^\times$ de sorte que le système
local obtenu $\FC_{\tau_o}$ ne sera pas irréductible.

\medskip

\noindent \textbf{0.4.} --- Suivant \cite{h-t}, on est en mesure de décrire assez simplement le complété
formel de $M_{I,o}$ le long de la strate $M_{I,s_o,1}^{=h}$. D'après les résultats de Berkovich, on en déduit
alors une expression de la restriction à $M_{I,s_o,1}^{=h}$ des faisceaux des cycles évanescents en termes
des systèmes locaux d'Harris-Taylor et des faisceaux des cycles évanescents $\Psi_{F_o,n}^{h,i}$ associés par
Berkovich aux anneaux $\Def_n^h$ tout en suivant l'action de $(D_\Am^{\oo,o})^\times \times W_o$.

\bigskip

\noindent \textbf{0.5.} --- Décrivons succinctement le contenu des différents paragraphes. Le premier
paragraphe est consacré aux rappels de \cite{lrs} et \cite{boy} sur les $\OC_o$-modules formels, \S
\ref{rapel-def}, sur les $\DC$-faisceaux elliptiques et leurs structures de niveau, \S \ref{rapel-fe}, les
variétés de modules $M_{I,o}$ munis d'une action par correspondances de $(D_\Am^\oo)^\times$, \S
\ref{retcor}. On rappelle ensuite, \S \ref{infini}, selon \cite{lrs}, la construction des systèmes locaux
$\LC_{\r_\oo}$ sur les variétés de $M_{I,o}$ associés à une représentation irréductible $\r_\oo$ des
inversibles de l'algèbre à division $\bar D_\oo$ centrale sur $F_\oo$ d'invariant $-1/d$. On rappelle enfin,
\S \ref{strates}, la stratification de la fibre spéciale et comme application du théorème de Serre-Tate, \S
\ref{completes}, on montre la lissité des $M_{I,s_o,1}^{\geq h}$.

\medskip

\noindent \textbf{0.6.} --- Le deuxième paragraphe s'intéresse plus en détail aux structures de niveau en y
démêlant parties étales et connexes. A tout $S$-point de $M_{I,s_o}$ est associé un $\phi$-faisceau: si $S
\longto M_{I,s_o}$ se factorise par $M_{I,s_o}^{=h}$, le $\phi$-faisceau associé se décompose alors en partie
étale et connexe. L'un des résultats clef de ce premier paragraphe est la proposition (\ref{scinde}) qui
montre qu'après extension des scalaires par une puissance assez grande du Frobenius, la décomposition en
composantes étale et connexe du $\phi$-faisceau d'un $S$-point de $M_{I,s_o}^{=h}$ est scindée modulo
$\pi_o^n$. On a aussi des résultats au niveau des déformations de points de $M_{I,s_o}^{=h}$, corollaire
(\ref{coro-defor}) et proposition (\ref{defor}). On montre enfin la relation de congruence, proposition
(\ref{cong}), qui permet de vérifier, \S \ref{cores-1}, la compatibilité des diverses définitions des
correspondances de Hecke.

\medskip

\noindent \textbf{0.7.} --- On introduit, \S \ref{igusa1}, suivant \cite{h-t}, \textit{les variétés d'Igusa
de première espèce} $\IC_{I^o,n}^{=h}$. On peut comprendre la démarche des définitions de celles-ci, de leurs
versions formelles $\hat \IC_{I^o,=h,m}(t)$ et de leurs liens aux variétés de Drinfeld, comme la différence
entre la donnée d'une structure de niveau sur un $\phi$-faisceau et la donnée de structures de niveaux sur
les parties étale et connexe du même $\phi$-faisceau. Les résultats dans ce sens sont les propositions
(\ref{igusa1-k}) et (\ref{iso-1}). En particulier on a un morphisme radiciel de $\IC_{I^o,n}^{=h}$ vers
$M_{I,s_o,a}^{=h}$ où $n$ est la multiplicité de $o$ dans $I$. Au niveau des complétés formels, le résultat
clef est la proposition (\ref{iso-1}) qui donne l'existence d'un isomorphisme canonique de $\hat
\IC_{I^o,=h,n}(n)$ vers $\widehat{M_{I,o,=h,a}}$. On définit ensuite, \S \ref{cores-1}, les correspondances
de Hecke associées aux \eles de $GL_h(F_o) \times GL_{d-h}(F_o)$ sur la tour des $\hat \IC_{I^o,=h,m}(t)$, de
manière compatible au théorème de Serre-Tate et à l'isomorphisme de la proposition (\ref{iso-1}).

\medskip

\noindent \textbf{0.8.} --- \textit{Les variétés d'Igusa de seconde espèce} $\JC_{I^o,m}^{=h}(s)$ de niveau
$s$, sont introduites \S \ref{igusa2}: le revêtement $\JC_{I^o,m}^{=h}(s) \longto \IC_{I^o,m}^{=h}$ est étale
de groupe de Galois $\DC_{o,h,s}^\times:=(\DC_{o,h}/(\Pi_{o,h}^{s+1}))^\times$ où $\DC_{o,h}$ est l'ordre
maximal de l'algèbre à division centrale sur $F_o$ d'invariant $1/h$ et où $\Pi_{o,h}$ est un uniformisante
de $\DC_{o,h}$. On définit \S \ref{cores-j}, les correspondances de Hecke sur la tour des
$\JC_{I^o,m}^{=h}(s)$ associées aux \eles de $(D_\Am^{\oo,o})^\times \times GL_{d-h}(F_o) \times
\widetilde{\NC_o}$ où $\widetilde{\NC_o}$ est le noyau du morphisme $GL_h(F_o) \times D_{o,h}^\times \times
W_{F_o} \longto \Zm$ qui au triplet $(g_o^c,\d_o,c_o)$ associe la valuation de $\det(g_o^c) \rn(\d_o)
\cl(c_o)$ où $\rn:D_{o,h}^\times \longto F_o^\times$ est la norme réduite et $\cl:W_{F_o} \longto F_o^\times$
est l'application de la théorie du corps de classe. Les résultats principaux sur les variétés d'Igusa sont
ceux de \S \ref{enonce-local}. Brièvement on se sert des $\JC_{I^o,m}^{=h}(s)$ pour ``détordre'' le schéma
formel $\hat \IC_{I^o,=h,m}(t)$ (cf. la proposition (\ref{iso-local}) et la remarque qui suit). Ce paragraphe
est très proche du texte de \cite{h-t}; en particulier en suivant loc. cit., on montre à la proposition
(\ref{niv-fini}), une version à niveau fini de ces résultats.

\medskip

\noindent \textbf{0.9.} --- Le cinquième paragraphe traite des cycles évanescents. \textit{A partir de ce
point et jusqu'à la fin, sauf mention expresse du contraire, on considère l'action naturelle, i.e. celle de
\cite{boy}, de $GL_d(F_o)$ sur $\Def_n^d$ tordue par l'application $g_o \mapsto \lexp t g_o^{-1}$; on
pointera cette modification par un tilde sur les espaces concernés.} On s'intéresse en premier lieu à
\textit{la restriction du faisceau des cycles évanescents} $R^i \Psi_{\eta_o}(\bar \Qm_l)$ de $M_{I,o}
\longto \spec (\OC_o)$ à la strate $M_{I,s_o,a}^{=h}$, en fonction du faisceau des cycles évanescents
$\widetilde{\Psi_{F_o,n}^{h,i}}$ du modèle local, i.e. de $\Def_n^h$ (cf. \S \ref{rapel-cycle}). L'énoncé
principal, théorème (\ref{theo-principal}), décrit cette restriction munie de son action de
$(\DC_\Am^{\oo,o})^\times \times P_{h,d}(F_o) \times W_{F_o}$. Formulons en rapidement l'énoncé: soit
$\twist_{\JC_{I^o,n}^{=h}(\oo)}(\widetilde{\Psi_{F_o,n}^{h,i}})$ le faisceau sur $\IC_{I^o,n}^{=h}$ tel que
sa restriction à $\JC_{I^o,n}^{=h}(s)$ muni de l'action naturelle de $\DC_{o,h}^\times$ est le faisceau
constant $\widetilde{\Psi_{F_o,n}^{h,i}}$ muni de l'action ``diagonale'' de $\DC_{o,h}^\times$. Ce faisceau
est naturellement muni d'une action de $(\DC_\Am^{\oo,o})^\times \times P_{h,d}^{op}(F_o) \times W_{F_o}$ et
on a un isomorphisme équivariant naturel
$$R^i \Psi_{\eta_o}(\bar \Qm_l)_{|M_{I,s_o,a}^{=h}} \simeq 
\twist_{\JC_{I^o,n}^{=h}(\oo)}(\widetilde{\Psi_{F_o,n}^{h,i}})$$
On relie ensuite \S \ref{action-cycle} ces faisceaux aux systèmes locaux\footnote{non irréductibles car la
restriction de $\tau_o$ à $\DC_{o,h}^\times$ est une somme directe de $e_{\tau_o}$ représentations
irréductibles.} $\FC_{\t_o}$ sur $M_{I,s_o,a}^{=h}$ associés aux revêtements d'Igusa de seconde espèce et à
une représentation admissible, irréductible $\t_o$ de $D_{o,h}^\times$. Cette description est donnée à la
proposition (\ref{prop-principal}).

\medskip

\noindent \textbf{0.10.} --- Enfin dans le dernier paragraphe, on donne, au niveau des complétés formels, des
équations explicites définissant les variétés d'Igusa de seconde espèce.


%% file: complements.tex
\section{Modules formels et variétés de Drinfeld: rappels} \label{defis}

\subsection{Rappels sur les déformations des $\OC_o$-modules formels}

\label{rapel-def}

\begin{defi}
Étant donnés une $\OC_o$-algèbre $R$ et $i:\OC_o \to R$ son morphisme structural, \textit{un $\OC_o$-module
formel sur $R$} est un couple $G=(F,(f_\l)_{\l \in \OC_o})$ où $F \in R[[X,Y]]$ et $f_\l \in R[[X]]$ avec
$f_\l(X)=i(\l)X$ modulo $(X^2)$, vérifiant les propriétés usuelles pour en faire un $\OC_o$-module, $F$
représentant l'addition et $f_\l$ la multiplication par $\l$.
\end{defi}

\noindent \textit{Remarques}:
- Si $R$ est un corps sur $\Fm_q$, il existe un
entier $d$, appelé \textit{la hauteur} de $G$, tel que $F(X,Y)=X+Y$ et $f_{\varpi_o}(X)$ est une série formelle en
$X^{q^d}$.\\
- Sur $\bar \Fm_q$, il existe à isomorphisme près un unique $\OC_o$-module formel de hauteur $d$
donnée que l'on notera $\Sigma_d$.
\\
- L'anneau des endomorphismes de ce $\OC_o$-module formel de
hauteur $d$ sur $\bar \k$, est l'ordre maximal $\DC_{o,d}$ de l'algèbre à division centrale sur $F_o$,
$D_{o,d}$, d'invariant $1/d$.

\begin{defis}
On note $\hat \OC_o^{nr}=\bar \Fm_q[[\varpi_o]]$ et on considère la catégorie $C$ dont les objets sont les
$\hat \OC_o^{nr}$-algèbres artiniennes. \textit{Une déformation} sur $R \in \ob(C)$ du $\OC_o$-module formel
de hauteur $d$ sur $\bar \Fm_q$, $\Sigma_d$, est un $\OC_o$-module formel $G=(F,(f_\l)_{\l \in \OC})$ sur $R$
dont la réduction modulo l'idéal maximal $\MC$ de $R$ est $\Sigma_d$. \textit{Une structure de niveau $n$ sur
$G$}, est la donnée d'un homomorphisme de $\OC_o$-module,
$$\iota_n:(\varpi_o^{-n}/\OC_o)^d \longto \MC,$$
tel que $f_{\varpi_o}(X)$ est divisible par $\prod_{\a \in (\varpi_o^{-1}\OC_o/\OC_o)^d} (X-\iota_n(\a))$.
\textit{Une déformation de niveau $n$} définie sur $R$, est par définition une déformation sur $R$ munie
d'une structure  de niveau $n$.
\end{defis}

\begin{prop} (cf. \cite{dr1} proposition 4.3)
Le foncteur des déformations de niveau $n$ du $\OC_o$-module formel de hauteur $d$ sur $\bar \Fm_q$,
$\Sigma_d$, est représenté par l'anneau $\Def_n^d$ vérifiant les propriétés suivantes:

\begin{itemize}

\item[-] pour $n=0$, le $\OC_o$-module formel universel sur $\Def_0^d$ est
de la forme $(X+Y,(f_\l)_{\l \in \OC_o})$ avec $f_{\varpi_o}(X)=\varpi_o X+a_1 X^q+ \cdots + a_{d-1}
X^{q^{d-1}}+X^{q^d}$ et $\Def_0^d \simeq \bar \Fm_q[[a_0,a_1,\cdots,a_{d-1}]]$ avec $a_0=\varpi_o$;

\item[-] soit pour $1 \leq i \leq d$, $e_i$ est le $i$-ème vecteur de la base canonique. Pour $n >0$, le morphisme
$$\phi: \bar \Fm_q[[v_1^n,\cdots,v_{d}^n]] \longto \Def_n^d$$
défini par $\phi(v_i^n)=\iota_n(e_i)$ où $\iota_n$ est la structure de niveau universelle, est un
isomorphisme.

\end{itemize}
\end{prop}

\rem Dans la suite quand on écrit
$$\Def_n^d \simeq \bar \Fm_q[[v_1^n,\cdots,v_d^n]] \quad (\hbox{resp. }
\Def_0^d \simeq \bar \Fm_q[[a_0,a_1,\cdots,a_{d-1}]], ~a_0=\varpi_o)$$ on sous-entend qu'il s'agit de
l'isomorphisme ci-dessus.

\medskip

\rem Le morphisme caractéristique $\bar \Fm_q[[\varpi_o]] \longto \Def_n^d$ peut se calculer de la façon
suivante. Le morphisme $\Def_0^d \longto \Def_1^d$ est donné par l'égalité de polynômes
$$\left | \begin{array}{ccccc}
v_1^1 & v_2^1 & \cdots & v_d^1 & X \\
(v_1^1)^q & (v_2^1)^q & \cdots & (v_d^1)^q & X^q \\
\vdots & \vdots & \vdots & \vdots & \vdots \\
(v_1^1)^{q^{d}} & (v_2^1)^{q^{d}} & \cdots & (v_d^1)^{q^{d}} & X^{q^{d}}
\end{array} \right |  = \l (\varpi_o X+a_1 X^q + \cdots + a_{d-1} X^{q^{d-1}} + X^{q^d})$$
où $\l$ est le coefficient dominant du polynôme de gauche. Ainsi l'image de $\varpi_o$ dans $\Def_1^d$ est
$$\left | \begin{array}{cccc}
v_1^1 & v_2^1 & \cdots & v_d^1 \\
(v_1^1)^q & (v_2^1)^q & \cdots & (v_d^1)^q \\
\vdots & \vdots & \vdots & \vdots \\
(v_1^1)^{q^{d-1}} & (v_2^1)^{q^{d-1}} & \cdots & (v_d^1)^{q^{d-1}}
\end{array} \right | ^{q-1}$$
En notant $a_i^{n}$ pour $0 \leq i <d$, les images de $a_i \in \Def_0^d$ dans $\Def_{n}^d$, les morphismes
$f_n:\Def_{n}^d \longto \Def_{n+1}^d$ se calculent de proche en proche par l'égalité suivante:
$$f_n(v_i^{n})=f_n(a_0^n) v_i^{n+1} + f_n(a_1^n) (v_i^{n+1})^q + \cdots + f_n(a_{d-1}^n) (v_i^{n+1})^{q^{d-1}} + 
(v_i^{n+1})^{q^d}$$

\begin{prop} (Deligne, Carayol: cf. \cite{ca} ou \cite{boy} 2.3)
L'action naturelle de $GL_d(\OC_o) \times \DC_{o,d}^\times \times
I_{o}$ sur $\Def_n^d$, se prolonge au noyau de l'application
$$\begin{array}{cll}
GL_d(F_o) \times D_{o,d}^\times \times W_{o} & \longto & \Zm \\
(g,\d,w) & \mapsto & \val(\det (g^{-1}) \rn(\d) cl (w))
\end{array}$$
où $\rn:D_{o,d}^\times \to F_o^\times$ est la norme réduite.
\end{prop}

On rappelle que l'action de $g^{-1} \in GL_d(F_o) \cap \Mm_d(F_o)$ se déduit de l'action naturelle à gauche
de $g$ sur la structure de niveau $\iota: (F_o/\OC_o)^d \longto \MC$; par ailleurs celle de
$\DC_{o,d}^\times$ est définie via l'identification avec $\aut(\Sigma_d)$.

\subsection{$\DC$-faisceau elliptiques et structures de niveau}

\label{rapel-fe}

\begin{defi} (cf. \cite{lrs})
Un $\DC$-faisceau elliptique $(\EC_i,j_i,t_i)$ sur un schéma $S$ est d'après \cite{lrs} un diagramme commutatif
$$
\xymatrix{
        \cdots \ar@{^{(}->}[r] & \EC_i \ar@{^{(}->}[r]^{j_i} & \EC_{i+1} \ar@{^{(}->}[r] & \cdots
        \\ \cdots \ar@{^{(}->}[r] & \lexp \t \EC_i \ar@{^{(}->}[ur]^{t_i} \ar@{^{(}->}[r]^{\lexp \t j_i} &
        \lexp {\t} {\EC}_{i+1} \ar@{^{(}->}[r] & \cdots
        }
$$
où:
\begin{itemize}

\item[-] $\EC_i$ est un $\DC_{X \times S}$-module à droite localement libre de rang $1$, et donc un $\OC_{X \times 
S}$-module
localement libre de rang $d^2$;

\item[-] $\lexp \t \EC_i$ est égal à $(\Id_X \times \frob_S)^* \EC_i$;

\item[-] $j_i$ et $t_i$ sont des injections $\DC_{X \times S}$-linéaires;

\item[-] $\EC_{i+d}\simeq \EC_i(\oo):= \EC_i \otimes_{\OC_X} \OC_X(\oo)$ et
le compos\'e $\EC_i \to \EC_{i+1} \to \cdots \to \EC_{i+d}$ est induit par l'injection canonique
$\OC_X \hookrightarrow \OC_X(\oo)$;

\item[-] $(pr_S)_*(\EC_i/\EC_{i-1})$ est un $\OC_S$-module localement libre de rang $d$ o\`u $pr_S:X \times S \to S$ 
est la
projection canonique. De manière équivalente, $\EC_i/\EC_{i-1}$ est isomorphe \`a
l'image directe $(\widetilde{i_{\oo}})_* (\Gamma_{\oo,i})$ d'un $\OC_S$-module
$\Gamma_{\oo,i}$ localement libre de rang $d$, par la section $\oo$:
$(\widetilde {i_{\oo}}):S \longto X \times S, ~ s \longmapsto (\oo,s);$

\item[-] l'image directe de $\coker t_i$ est un $\OC_S$-module localement libre de rang $d$. Le support de
$\coker t_i$ est disjoint de $\bad \cup \{ \oo \} \times S$, où $\bad$ désigne l'ensemble des places $x$ de
$X$ telles que $\DC_x$ n'est pas isomorphe à $\Mm_d(\OC_x)$. De manière \'equivalente, $\coker t_i$ est
isomorphe à l'image directe $(\widetilde{i_{0,i}})_* (\Gamma_{0,i})$ d'un $\OC_S$-module $\Gamma_{0,i}$
localement libre de rang $d$, par la section
$$ (\widetilde {i_{0,i}}):S \longmapright{(i_{0,i},id_S)} X \times S$$
induite par un morphisme $i_{0,i}:S \to X$ tel que $i_{0,i}(S) \subset |X'|$.
\end{itemize}
\end{defi}

\rem Les inclusions $\EC_i \hookrightarrow \EC_{i+1}$ étant des isomorphismes sur $(X \backslash \{ \oo
\})\times S$ et le support de $\coker t_i$ étant disjoint de $\oo \times S$, on en déduit que la donnée des
morphismes $(t_i)_i$ est équivalente à la donnée d'un seul $t_i$. Les morphismes $i_{0,i}$ sont indépendants
de $i$; on le note $i_0$, le morphisme caractéristique du $\DC$-faisceau elliptique.

\begin{defi} (cf. \cite{lrs})
Étant donnés $(\EC_i,j_i,t_i)$ un $\DC$-faisceau elliptique d\'efini sur $S$ et $I$ un idéal de $A$ tel que
$V(I) \cap i_0(S)=\emptyset$, \textit{une $I$-structure de niveau} sur $(\EC_i,j_i,t_i)$ est un isomorphisme
de $\DC_{I \times S}$-modules à droite, $\tilde \iota_I: \DC_I \boxtimes \OC_S \longmapright{\sim} \EC_{I
\times S}$ tel que le diagramme suivant est commutatif
$$
\xymatrix{
        \lexp \t \EC_{|I \times S} \ar[rr]^{t_{|I \times S}} && \EC_{|I \times S}\\
        & \DC_I \boxtimes \OC_S  \ar[ul]_{\lexp \t {\tilde \iota_I}} \ar[ur]_{\tilde \iota_I}
        }
$$
\end{defi}

Pour définir la notion de structure de niveau en une place $v \in i_0(S)$ d'après \cite{boy}, on est amené à
considérer les objets suivants. Pour alléger les notations on prend $v=o$ ce qui sera d'ailleurs le cas dans
la suite.

\begin{defis} \label{diverses-defis}
\begin{itemize}

\item[-] Le $(\OC_o \otimes \OC_S)$-module $\EC_i \otimes \OC_o$ est ind\'ependant de $i$; on le note
$\EC_o$. Un isomorphisme $\DC_o \simeq \Mm_d(\OC_o)$ étant fixé, soit $\FC_o$ le $(\OC_o \otimes
\OC_S)$-module libre de rang $d$ défini par $E_{1,1}.\EC_o$, où $E_{1,1}$ est l'idempotent de $\Mm_d(\OC_o)$
associé au premier vecteur de la base canonique. Par équivalence de Morita, on a un isomorphisme
$\DC_o$-équivariant $\EC_o \simeq \FC_o^d$. Les morphismes $t_i$ induisent alors un morphisme $t'_o:\lexp \t
\FC_o \longto \FC_o$. \textit{Le $\phi$-faisceau associé} est par définition le couple $(\FC_o,\phi_o)$ avec
$\phi_o=t_o': \lexp \t \FC_o \longto \FC_o$.

\item[-] \textit{Le $\OC_o$-module de Dieudonné sur $B$ associé à $(\EC_i,j_i,t_i)$} est le couple $(V_o,\vphi_o)$ 
o\`u
$V_o$ est le $(\OC_o \hat \otimes_{\k(o)} \OC_S)$-module localement libre de rang $d$,
$$\FC_o \otimes_{(\OC_o \otimes_{\k(o)} \OC_S)} (\OC_o \hat \otimes_{\k(o)} \OC_S)$$
et $\vphi_o:(\Id_{\OC_o} \hat \otimes_{\k(o)} \frob_{\k(o)})^* V_o \longto V_o$ est l'application induite par 
$\phi_o$
Pour tout $n$, on notera $V_{o,n}:=V_o \varpi_o^n \backslash V_o$.

\item[-] Pour tout $n$, on note $\FC_{o,n}=\FC_o \otimes_{\OC_o} (\OC_o \backslash \MC_o^{-n})$ et on considère
$Gr(\FC_{o,n})$, le $S$-schéma en $\OC_o$-modules finis, d'ordre $q^{nd}$ tués par $\varpi_o^n$, qui
représente le foncteur de la catégorie des $S$-schémas dans celle des ensembles:
$$S' \mapsto \{ u \in \hom_{\OC_S}(\FC_{o,n},S')~/~u(\phi_{o,n}(x))=u(x)^q ~ \forall x \in \FC_{o,n} \}$$
On a alors la suite exacte
$$0 \longto Gr(\FC_{o,n}) \longmapright{i_n} Gr(\FC_{o,n+1}) \longmapright{\varpi_o^n} Gr(\FC_{o,n+1}).$$

\item[-] Si $P$ est un $R$-point d'un schéma $Y$, on note $[P]$ le sous-$R$-schéma de $Y$ qu'il définit. Pour $(P_i)$ 
une
famille finie
de tels points, on note $\sum [P_i]$ le sous-schéma de $Y$ défini par le faisceau d'idéaux produit des faisceaux
d'idéaux définissant les $[P_i]$.

\item[-] \textit{Une $\MC_o^n$-structure de niveau sur $(\EC_i,j_i,t_i)/S$} est un homomorphisme de $\OC_o$-modules
$\iota_{o,n}':(\MC_o^{-n}/\OC_o)^d  \longto Gr(\FC_{o,n})(S)$ tel que le sous-schéma
$${\DS \sum_{z \in (\OC_o/\MF_o^n)^d}}[\iota_{o,n}'(z)]$$
de $Gr(\FC_{o,n})$ coïncide avec $Gr(\FC_{o,n})$.
\end{itemize}
\end{defis}

\rem Pour tout \ele $z$ de $(\OC_o/\MC_o^n)^d$, $\iota_{o,n}'(z)$ est un \ele de $\FC_{o,n}^*$ tel que
$\phi_{o,n}^*(\iota_{o,n}'(z))=\iota_{o,n}'(z)^{q}.$
Le morphisme $\iota_{o,n}'$ fournit alors un homomorphisme de $\OC_o$-modules $(\OC_o/\MC_o^n)^d \longto \FC_{o,n}^*$ 
qui
après équivalence de Morita, donne un homomorphisme de $\DC_o$-modules
$\iota_{o,n}:\DC_{o,n} \longto \EC_{o,n}^*$ tel que le diagramme ci-dessous commute
$$
\xymatrix{
        \DC_{o,n} \ar[rr]^{\iota_{o,n}} \ar[dr]_{\lexp \t \iota_{o,n}} & & \EC_{o,n}^* \ar[dl]_{t_{o,n}^*} \\
        & \lexp \t \EC_{o,n}^*
        }
$$
Dans la suite, $\iota_{o,n}$ désignera la structure de niveau après équivalence de Morita.

\begin{prop} (cf. \cite{boy} propositions (7.1.3) et (7.1.4)) \label{niveau-so}
Soit
$$\iota'_{o,n}:(\MC_o^{-n}/\OC_o)^d \longto Gr(\FC_{o,n})(S)$$
une $\MC_o^n$-structure de niveau définie sur un schéma $S$ tel que $S^o$ l'ouvert complémentaire de
$i_0^{-1}(\spec \k(o))$, est non vide. Alors $\iota_{o,n}' \otimes_{\OC_S} \OC_S(S^o)$ induit un isomorphisme
de $\OC_o$-modules
$$(\OC_o/\MC_o^n)^d \boxtimes \OC_S(S^o) \longto \left ( \FC_{o,n} \otimes_{\OC_S} \OC_S(S^o) \right ) ^*.$$
Réciproquement soit $S$ est un schéma intègre tel que $S^o$ est non vide. Si l'homomorphisme de $\OC_o$-modules
$\iota_{o,n}':(\MC_o^{-n}/\OC_o/)^d  \longto Gr(\FC_{o,n})(S)$ induit un isomorphisme sur $S^o$ $(\MC_o^{-n}/\OC_o)^d
\boxtimes \OC_S(S^o) \longmapright{\sim} (\FC_{o,n} \otimes_{\OC_S} \OC_S(S^o))^*,$ alors $\iota_{o,n}'$ est une
$\MC_o^n$-structure de niveau.
\end{prop}

\subsection{Schémas de modules et correspondances de Hecke}

\label{retcor}

\begin{prop} (cf. \cite{lrs})
Le classifiant des classes d'équivalence des $\DC$-faisceaux elliptiques munis d'une $I$-structure de niveau,
définit un schéma régulier (un champ si $I=A$), $M_I \to X'$ de dimension relative $d-1$.
\end{prop}

\begin{defi} On notera
$M_{I,o}:= M_I \times_{X'} \spec (\OC_o)$ et $M_{I,s_o}$ (resp. $M_{I,\eta_o}$) la fibre spéciale (resp.
générique) de $M_{I,o}$.
\end{defi}

\marque Soient $\Am$ l'anneau des adèles de $F$ et $D^\times_\Am:=D^\times \otimes_F \Am$. La limite
projective $M_o:={\DS \lim_{\genfrac{}{}{0pt}{}{\longleftarrow}{I}}} ~M_{I,o}$ est munie d'une action par
correspondance de Hecke de $(D_\Am^\oo)^\times$. Au dessus de la fibre générique $M_{\eta_o}$ de $M_o$ et à
niveau fini, pour $K_\Am^\oo$ un sous-groupe compact ouvert de $(D_\Am^\oo)^\times$, l'action d'un \ele
$g^{\oo}$ de $(D_\Am^{\oo})^\times$ se décrit par la correspondance géométrique
$$
\xymatrix{
        & (M_{\eta_o})^{(K^\oo_\Am \cap (g^{\oo})^{-1} K^\oo_\Am g^{\oo})} \ar[dl]_{c_1} \ar[dr]^{c_2} \\
        (M_{\eta_o})^{K^\oo_\Am} \ar[dr] \ar@{-->}[rr] & & (M_{\eta_o})^{K^\oo_\Am} \ar[dl] \\ & \spec (F_o)
        }
$$
où le morphisme $c_1$ (resp. $c_2$) est induit par l'inclusion\footnote{cf. \cite{lrs} \S 7}
$$(K^\oo_\Am \cap (g^\oo)^{-1} K^\oo_\Am g^\oo) \subset K^\oo_\Am \hbox{ (resp. } (K^\oo_\Am \cap (g^\oo)^{-1} 
K^\oo_\Am g^\oo) \longmapright{Ad(g^\oo)} K^\oo_\Am)$$
En particulier on note $K^\oo_{\Am,I}$ le noyau de $(\DC^\oo_\Am)^\times \longto \prod_{x\neq \oo}
\DC_x^\times \otimes_{\OC_x} \OC_x/\MC_x^{n_x}$, o\`u $n_x$ est la multiplicité de $x$ dans $I$. On a
$(M_{\eta_o})^{K^\oo_{\Am,I}}=M_{I,\eta_o}$.

\marque Dans \cite{boy} \S 7.3, on montre que ces correspondances s'étendent de manière unique sur $M_o$ tout
entier. Si la composante $g^\oo_o$ de $g^\oo$ en $o$ est triviale, le résultat est évident. Précisons la
situation pour $g_o \in GL_d(F_o) \cap \Mm_d(\OC_o)$. On pose $S:=M_o$ et $S^o:=M_{\eta_o}$ et on considère
le $\DC$-faisceau elliptique universel sur $S$ muni de sa $I$-structure de niveau universelle pour tout idéal
$I$ de $A$. La structure de niveau fournit un morphisme
$$\iota'_o:(F_o/\OC_o)^d \times S \longto \FC_o^* \otimes_{\OC_o} (F_o/\OC_o)$$
selon les notations précédentes. Sur l'ouvert $S^o$, $\iota_o^{'o}:=\iota'_o \otimes_{\OC_o} F_o$ est un
isomorphisme. A $g_o$ on associe, d'après loc. cit., la multiplication à gauche de $\lexp t g_o \times$ sur
$(F_o/\OC_o)^d$ et donc le morphisme
$$[g_o^*]:\FC_o^* \otimes_{\OC_o} (F_o/\OC_o) \longto \FC_o^* \otimes_{\OC_o} (F_o/\OC_o)$$
tel que sa restriction $[g_o^*]^o$ à $S_o$ est par définition induite, via l'isomorphisme $\iota_o^{'o}$, par la
multiplication à gauche de $\lexp t g_o$ sur $(F_o/\OC_o)^d$. On note $[g_o]$ l'endomorphisme de $\FC_o$ déduit de
$[g_o^*]$. L'image par $g_o$ du $\DC$-faisceau elliptique universel $((\EC_i,j_i,t_i),\iota_\oo)$ est le 
$\DC$-faisceau
$((\EC_i',j_i',t_i'),\iota_\oo')$ défini via les diagrammes cartésiens

$$\diagram
\EC_i' \rrto \dto & & \FC_o \dto^{[g_o]} \\
\EC_i \rrto_i & & \FC_o
\enddiagram$$
où $i$ est l'injection canonique, ou de manière duale par $(\EC')^*:=\EC^* \otimes_{\EC_o^*,[g_o^*]} \EC_o^*$. La
structure de niveau est alors définie de manière naturelle comme suit. Soient $n$ et $m$ tels que $\ker \lexp t g_o
\subset (\MC_o^{-n}/\OC_o)^d$ et $(\MC_o^{-m}/\OC_o)^d \subset \im g_o$. Si $\iota_{o,n}$ est la structure de niveau 
$n$
sur $\FC_{o,n}$, la composée
$$(\MC_o^{-m}/\OC_o)^d \longto (\MC_o^{-n}/\OC_o)^d/\ker g_o \longmapright{\iota_{o,n}} \FC_{o,n}^*
\longmapright{[g_o^*]} \FC_{o,m}^*$$ définit une flèche $\iota'_{o,m}:(\MC_o^{-m}/\OC_o)^d \longto (\FC_{o,n}')^*$ 
dont
on vérifie aisément en utilisant la proposition (\ref{niveau-so}), qu'elle définit une structure de niveau $m$.

\subsection{Structures de niveau à l'infini}
\label{infini}

On résume simplement ici le paragraphe 8 de \cite{lrs}.

\marque Étant donné un $\Fm_q$-schéma $S$ et $(\EC_i,j_i,t_i)$ un $\DC$-faisceau elliptique sur $S$, on
considère pour tout $i \in \Zm$,  les $\OC_\oo \hat \otimes \OC_S$-modules localement libre de rang $d^2$:
$$\check M_i:= \EC_i^\vee | (X \times S)^\vee_\oo$$
où $\EC_i^\vee$ est le faisceau dual de $\EC_i$ et $(X \times S)^\vee_\oo$ est la complétion de $X \times S$
le long de $\{ \oo \} \times S$. Les $j_i$ définissent alors un système inductif $\check M:=\cdots
\hookrightarrow \check M_i \hookrightarrow \check M_{i+1} \hookrightarrow \cdots$ et l'on identifie les
$\check M_i$ avec leur image dans $\check M$ qui est donc un $F_\oo \hat \otimes \OC_S$-module de rang $d^2$.
Les $t_i$ induisent alors un isomorphisme
$$\check M \longmapright{\sim} \lexp \t {\check M}=(F_\oo \hat \otimes \frob_S)^* \check M$$
qui envoie $\check M_i$ surjectivement sur $\lexp \t {\check M}_{i+1}$. On définit alors $\check \psi: \lexp
\t {\check M} \longto \check M$ comme l'inverse de l'isomorphisme ci-dessus; $\check \psi(\lexp \t {\check
M}_i)=\check M_{i-1} \subset \check M_i$.

\marque L'action à droite de $\DC$ sur les $\EC_i$ induit une action à gauche sur $\check M$ qui commute à
$\check \psi$ et qui stabilise $\check M_i \subset \check M$. Un isomorphisme $\DC_\oo \simeq \Mm_d(\OC_\oo)$
étant fixé, on a une équivalence de Morita
$$(\check M, \check \psi)=(\check N,\check \psi)^d \hbox{ et } \check M_i=(\check N_i)^d$$
Ainsi $\check N$ est un $F_\oo \hat \otimes \OC_S$-module de rang $d$ et $\check N_i \subset \check N$ en est un
$\OC_\oo \hat \otimes \OC_S$-sous-module localement libre de rang $d$ vérifiant
$$\varpi_\oo \check N_i = \check N_{i- (\deg \oo) d} \subset \cdots \subset \check N_{i-1} \subset \check N_i, \qquad
\check \psi (\lexp \t {\check N}_i) = \check N_{i-1}$$  tel que le quotient $\check N_i / \varpi_\oo \check
N_i \longto \check N_i/ \check \psi (\lexp \t {\check N}_i)$ est un $\OC_{\{ \oo \} \times S}$-module
supporté par le graphe d'un $\Fm_q$-morphisme de schéma $i_{\oo,i}:S \longto \{ \oo \} $ et est localement
libre de rang $1$ sur son support avec $i_{\oo,i+1}=i_{\oo,i} \circ \frob_S$.

\marque Soit $S$ un schéma muni d'un $\Fm_q$-morphisme de schéma $i_{\oo,0}:S \longto \{ \oo \}$. Tout autre
$\Fm_q$-morphisme de schéma de $S$ vers $\{ \oo \}$ est alors de la forme $i_{\oo,i}:=i_{\oo,0} \circ
\frob_S^i$ pour un unique $i \in \Zm / \deg (\oo) \Zm$. On associe au couple $(S,i_{\oo,0})$ le triplet
$(N_{d,1},\psi_{d,1},\NC_{d,1})$ où
$$N_{d,1}=(F_\oo \hat \otimes \OC_S)^d= \bigoplus_{i=0}^{\deg(\oo)-1} (F_\oo \hat \otimes_{\k(\oo),i_{\oo,i}^*}
\OC_S)^d$$
de base canonique $(e_{i,j})_{0 \leq i < \deg (\oo), 1 \leq j \leq d}$, $\psi_{d,1}: \lexp \t N_{d,1}
\longmapright{\sim} N_{d,1}$ est défini par $\psi_{d,1}(e_{i,j})=e_{i+1,j}$ pour $i \neq \deg(\oo) -1$ et
$\psi_{d,1}(e_{\deg(\oo)-1,j})=e_{0,j-1}$ pour $j \neq 1$ et $\varpi_\oo e_{0,d}$ pour $j=1$, et
$$\NC_{d,1}=(\OC_\oo \hat \otimes \OC_S)^d =\bigoplus_{i=0}^{\deg(\oo)-1} (\OC_\oo \hat \otimes_{\k(\oo),i_{\oo,i}^*} 
\OC_S)^d$$

\marque Soit $\k(\oo)_d$ l'extension de degré $d$ de $\k(\oo)$ et soient $F_{\oo,d}=F_\oo \hat
\otimes_{\k(\oo)} \k(\oo)_d$, $\OC_{\oo,d}$ son anneau des entiers et $\s_{\oo,d}=F_\oo \hat
\otimes_{\k(\oo)} \frob_q^{\deg(\oo)}$. Soit $F_{\oo,d}[\t_\oo]$ l'algèbre de polynômes sur $F_{\oo,d}$, non
commutative avec la règle de commutation $\t_\oo a =\s_{\oo,d}^{-1} (a) \t_\oo$ pour tout $a \in F_{\oo,d}$.
L'\ele $\t_\oo^d - \varpi_\oo$ est alors central dans $F_{\oo,d}[\t_\oo]$ et $\bar
D_\oo:=F_{\oo,d}[\t_\oo]/(\t_\oo^d - \varpi_\oo)$ est alors ``l'''algèbre à division centrale sur $F_\oo$
d'invariant $-1/d$, d'ordre maximal $\bar \DC_\oo:= \OC_{\oo,d}[\t_\oo]/(\t_\oo^d-\varpi_\oo)$.

\marque Si $\l:S \longto \spec (\k(\oo)_d)$ est un $\Fm_q$-morphisme de schémas, on peut construire une
injection de $F_\oo$-algèbres
$$\l^*:\bar D_\oo \hookrightarrow \End(N_{d,1},\psi_{d,1})$$
de la façon suivante où $i_{\oo,0}:S \longmapright{\l} \spec (\k(\oo)_d) \longmapright{can} \{ \oo \}$.
Pour $\a \in \k(\oo)_d$, l'image de $1 \hat \otimes \a \in F_{\oo,d} \subset \bar D_\oo$ est donnée par
$$\l^*(1 \hat \otimes \a)(e_{i,j})=(1 \hat \otimes \l^*(\frob_q^{i-j \deg(\oo)}(\a)))e_{i,j}$$
et l'image de $\t_\oo \in \bar D_\oo$ par
$$\l^*(\t_\oo)(e_{i,j})= \left \{ \begin{array}{l} \varpi_\oo e_{i,d} \hbox{ si } j=1 \\ e_{i,j-1} \hbox{ sinon} 
\end{array} \right. $$

Si $\l$ et $\l'$ sont deux $\Fm_q$-morphismes de schémas de $S$ vers $\spec (\k(\oo)_d)$ qui relèvent
$i_{\oo,0}$, on a $\l'=\l \circ \frob_S^{n \deg (\oo)}$ pour $n \in \Zm/d \Zm$ et alors $(\l')^* =\l^* \circ
Ad(\t_\oo^{-n})$. On vérifie aussi que $\l^*(\bar D_\oo)$ laisse $\NC_{d,1}$ stable.

\begin{defi} Une structure de niveau à l'infini sur $(\EC_i,j_i,t_i)$ est une paire $(\l,\a)$ où
$\l:S \longto \spec (\k(\oo)_d)$ est un $\Fm_q$-morphisme de schémas qui relève le pôle $i_{\oo,0}$ et où
$\a:\NC_{d,1} \longmapright{\sim} \check N_0$ est un isomorphisme de $\OC_\oo \hat \otimes \OC_S$-modules tel
que le diagramme suivant soit commutatif
$$\diagram
\lexp \t {\check \NC}_{d,1} \rrto^{\genfrac{}{}{0pt}{}{\lexp \t \a}{\sim}} \dto^{\psi_{d,1}} & & \lexp \t {\check N_0} 
\dto_{\check \psi} \\
\NC_{d,1} \rrto^{\genfrac{}{}{0pt}{}{\a}{\sim}} & & \check N_0 \\
\enddiagram$$
\end{defi}

\marque On a une notion évidente d'isomorphisme entre $\DC$-faisceaux elliptiques sur $S$ munis d'une
structure de niveau à l'infini; on note $\MC_I$ (resp. $\widetilde{\MC_I} (S)$) la catégorie des classes
d'isomorphismes des $\DC$-faisceaux elliptiques sur $S$ munis d'une structure de niveau $I$ (resp. et d'une
structure de niveau à l'infini). On obtient ainsi une catégorie fibrée $\widetilde{\MC_I}$ qui est un
pro-champ; en effet se donner un isomorphisme
$$\a:\NC_{d,1} \longmapright{\sim} \check N_0$$
de $\OC_\oo \hat \otimes \OC_S$-modules, revient à se donner un système projectif
$$(\a_n= \a \mod \varpi_\oo^{n+1})_{n \geq 0}$$
d'isomorphismes de $\OC_\oo / (\varpi_\oo^{n+1}) \hat \otimes \OC_S$-modules et $\a$ commute avec les $\psi$ \ssi les 
$\psi_n$
commutent avec les $\psi$ modulo $\varpi_\oo^{n+1}$.

\marque On a en outre un morphisme d'oubli
$$r_{\oo,I}:\widetilde{\MC_I} \longto \MC_I$$
qui envoie $((\EC_i,j_i,t_i),\iota_I,(\l,\a))$ sur $((\EC_i,j_i,t_i),\iota_I)$, ainsi qu'un morphisme de structure
$$\l_I:\widetilde{\MC_I} \longto \spec (\k(\oo)_d)$$
qui envoie $((\EC_i,j_i,t_i),\iota_I,(\l,\a))$ sur $\a$ et qui relève $i_{\oo,0} \circ r_{\oo,I}$.
Pour $I \subset J \subset X \backslash \{ \oo \}$, on a des $2$-diagrammes commutatifs
$$\diagram
\widetilde{\MC_J} \rto^{\tilde r_{J,I}} \dto^{r_{\oo,J}} & \widetilde{\MC_I} \dto^{r_{\oo,I}} \\
\MC_J \rto^{r_{J,I}} & \MC_I
\enddiagram$$

\marque Sur $\widetilde{\MC_I}$, on a des actions continues à droites du groupe profini $\bar \DC_\oo^\times$
et de $\Zm/d\Zm$: $\d \in \bar \DC_\oo^\times$ (resp. $n \in \Zm/d\Zm$) envoie
$((\EC_i,j_i,t_i),\iota_I,(\l,\a))$ sur
$$((\EC_i,j_i,t_i),\iota_I,(\l, \a \circ \l^*(\d))) \quad (resp. ~((\EC_i,j_i,t_i),\iota_I,(\l \circ \frob_S^{n 
\deg(\oo)},\a))).$$
Comme on a
$$(\l \circ \frob_S^{n \deg(\oo)})^*=\l^* \circ Ad(\t_\oo^{-n})$$
ces deux actions induisent une action à droite continue du groupe profini
$$\bar \DC_\oo^\times \semi \Zm/d \Zm$$
où $n \in \Zm/d \Zm$ agit sur $\bar \DC_\oo^\times$ par $Ad(\t_\oo^{-n})$. On identifie alors ce produit semi-direct 
avec
$$\bar D_\oo^\times / \varpi_\oo^\Zm$$
en envoyant $(\d,n)$ sur $\d \t_\oo^{-n}$. On obtient ainsi une action à droite de $\bar D_\oo^\times /\varpi_\oo^\Zm$ 
qui commutent
aux $r_{\oo,I}$, $\tilde r_{J,I}$ et $\l_I$ si on fait agir  $\bar D_\oo^\times /\varpi_\oo^\Zm$ sur $\spec 
(\k(\oo)_d)$ à travers
son quotient
$$- \oo \circ \rn: \bar D_\oo^\times \longto \Zm/d\Zm$$
($\gal(\k(\oo)_d/\k(\oo))=\Zm/d\Zm$).

\begin{theo} (cf. \cite{lrs} théorème (8.10) et proposition (8.8))

(i) Le morphisme de pro-champs
$$r_{\oo,I}:\widetilde{\MC_I} \longto \MC_I$$
est représentable et est un revêtement pro-galoisien de groupe de Galois $\bar D_\oo^\times /\varpi_\oo^\Zm$;
$\widetilde{\MC_I}$ est alors un schéma que l'on notera $\widetilde{M_I}$.

(ii) Les correspondances de Hecke sur les $M_I$ associées aux \eles de $(D_\Am^\oo)^\times \times \Zm$ se relèvent 
sur
$\widetilde{M_I}$ et commutent
à l'action de $\bar D_\oo^\times /\varpi_\oo^\Zm$. De plus $F^\times$ envoyé diagonalement dans $(\bar D_\oo^\times 
/\varpi_\oo^\Zm)
\times (D_\Am^\oo)^\times$ agit trivialement sur les $\widetilde{M_I}$.
\end{theo}

\rem La caractéristique étant disjointe de $\oo$, les arguments de loc. cit. s'appliquent sans modification sur les 
variétés $M_I$
même dans le cas de mauvaise réduction.

\subsection{Stratification de la fibre spéciale}
\label{strates}

\begin{defi} \label{defi-phd}
Pour tout $1 \leq h \leq d$, on notera $P_{h,d}$ le parabolique standard associé aux $h$
premières coordonnées et soit $P_{h,d}^{op}$ le parabolique opposé
\end{defi}

\begin{prop} \label{existence-et}
(cf \cite{lau} proposition 2.4.6 et \cite{boy} lemme 6.2.3)
\begin{itemize}
\item[(a)] Tout $\OC_o$-module de Dieudonné $(V_o,\vphi_o)$ sur une extension $k$ de $\k(o)$ se décompose en une somme 
directe
$(V_o^c,\vphi_o^c) \oplus (V_o^{et},\vphi_o^{et})$
où $\vphi_o^c$ est topologiquement nilpotent et $\vphi_o^{et}$ est bijective.

\item[(b)] Tout $\OC_o$-module de Dieudonné $(V_o,\vphi_o)$ sur une $\OC_o$ algèbre $R$ dans laquelle l'image de 
$\varpi_o$
 est nilpotente, se dévisse en une suite exacte
$$0 \longto (V_o^{et},\vphi_o^{et}) \longto (V_o,\vphi_o) \longto (V_o^c,\vphi_o^c) \longto 0$$
où $\vphi_o^c$ est topologiquement nilpotent et $\vphi_o^{et}$ est bijective.
\end{itemize}
\end{prop}

\begin{defi}
Pour tout $1 \leq h \leq d$, on définit dans \cite{boy} un sous-schéma $M_{I,s_o}^{=h}$ de $M_{I,s_o}$ de
pure dimension $d-h$, stabilisé par les correspondances de Hecke et caractérisé par l'une des propriétés
équivalentes suivantes: en tout point géométrique de $M_{I,s_o}^{=h}$,

\begin{itemize}
\item[-] le $\OC_o$-module de Dieudonné qui lui est associé, a sa composante connexe de hauteur $h$.

\item[-] le $\phi$-faisceau $(\FC_o,\phi_o)$ qui lui est associé est tel que si on note $M_{o,1}$ la matrice de
$\phi_{o,1}:=\phi_o \otimes_{\OC_o} (\OC_o \backslash \MC_o^{-1})$ dans une base quelconque de $\FC_{o,1}^*:=\FC_o^*
\otimes_{\OC_o} (\MC_o \backslash \OC_o)$, alors d'après loc. cit.,
$$M_{o,1}^{!h}:=M_{o,1} (\lexp \t M_{o,1}) \cdots (\lexp {\t^{h-1}} M_{o,1})$$
a un mineure d'ordre $d-h$ inversible et tous ses mineures d'ordre $r>d-h$ sont nuls.
\end{itemize}
\end{defi}

\rem Le schéma $M_{I,s_o}^{=h}$ n'est pas réduit sauf si la multiplicité de $o$ dans $I$ est nulle; on notera
$M_{I,s_o,red}^{=h}$ le réduit associé.

\begin{defis} Soit $I=I^o \MC_o^n$ un idéal de $A$, $n$ étant la multiplicité de $o$ dans $I$. L'ensemble
$$G/P(d,h,n):=GL_d(\OC_o/\MC_o^n)/P_{h,d}(\OC_o/\MC_o^n)$$
classifie les facteurs directs $a$ de rang $h$ de $(\MC_o^{-n}/\OC_o)^d$. Pour tout $a \in G/P(d,h,n)$, on
notera
$$P_{h,d;a}(\OC_o/\MC_o^n):=aP_{h,d}(\OC_o) a^{-1}$$
$$K_{o^n}:=\ker (GL_d(\OC_o) \longto GL_d(\OC_o/\MC_o^n))$$
$$K_{o,h,d,n}(a)=K_{o^n} \cap P_{h,d;a}^{op}(\OC_o/\MC_o^n)$$
\end{defis}

\begin{defi} (cf. \cite{boy} définition (10.4.1))
Pour tout $a \in G/P(d,h,n)$, il existe un sous-schéma fermé
$M_{I,s_o,a}^{=h}$ de $M_{I,s_o}^{=h}$ stable sous
l'action des correspondances de Hecke associées aux éléments de
$$K_\Am^\oo/K_{\Am,I}^{\oo,o} \times K_{o,h,d,n}(a)$$
\end{defi}

\rem Concrètement si $\iota_{o,n}$ est la
$\MC_o^n$-structure de niveau universelle sur $M_{I,o}$ alors le noyau de $\iota_{o,n}^{et}:=(\iota_{o,n}
\times_{M_{I,o}} M_{I,s_o,a}^{=h})_{|\FC_o^{et}}:(\MC_o^{-n}/\OC_o)^d \longto \FC_o^{et,*}$, où $\FC_o^{et}$
est la partie étale du $\phi$-faisceau universel sur $M_{I,s_o,a}^{=h}$ donnée par la proposition
(\ref{existence-et}), est le facteur direct $a$ de $(\MC_o^{-n}/\OC_o)^d$.

\marque On en déduit alors le fait fondamental suivant.

\noindent \textbf{Propriété géométrique fondamentale:} \textit{les strates non supersingulières sont
induites, i.e. \footnote{On rappelle que l'application canonique $GL_d(\OC_o) \longto GL_d(\OC_o/\MC_o^n)$
est surjective.}
$$M_{I,s_o}^{=h}=M_{I,s_o,1}^{=h} \times_{P_{h,d}^{op}(\OC_o/\MC_o^n)} GL_d(\OC_o/\MC_o^n)$$
où $M_{I,s_o,1}^{=h}$ est la composante associée à la classe de $I_d$ dans $G/P(d,h,n)$.}

\begin{prop} L'ensemble des points géométriques de $M_{I,s_o,a}^{=h}$ au dessus d'un point géométrique donné
de $M_{I^o,s_o}^{=h}$ est en bijection avec les isomorphismes

$$\iota_{o,n}^{et}:(\MC_o^{-n}/\OC_o)^d/a \longto \FC_{o,n}^{et,*}$$
tels que $\phi_{o,n}^{et,*} \circ \iota_{o,n}^{et}=\lexp \t \iota_{o,n}^{et}$ où $n$ est la multiplicité de
$o$ dans $I$. Cet ensemble est alors de cardinal $\#GL_{d-h}(\OC_o/\MC_o^n)$.
\end{prop}

\begin{defi} On notera $M_{I,s_o}^{\geq h} \hbox{ (resp. }M_{I,s_o,a}^{\geq h})$
l'adhérence schématique de $M_{I,s_o}^{=h}$ (resp. $M_{I,s_o,a}^{=h}$) dans $M_{I,s_o}$.
\end{defi}

\subsection{Théorème de Serre-Tate et conséquences}
\label{completes}

D'après la théorie du module de coordonnées, cf. \cite{gen}, $\Def_n^d$ représente aussi les déformations de
niveau $n$, du $\OC_o$-module de Dieudonné sur $\bar \Fm_q$ connexe de hauteur $d$. \footnote{\label{note}
Plus généralement, on montre qu'il y a une équivalence de catégories, des déformations de niveau $n$ du
$\OC_o$-module divisible de type $(h,j)$ (la composante connexe est de hauteur $h$ et la partie étale de
hauteur $j$ (cf. \cite{dr1})) sur $\bar \Fm_q$ vers les déformations de niveau $n$ du $\OC_o$-module de
Dieudonné sur $\bar \Fm_q$ dont la partie connexe est de rang $h$ et la partie étale de rang $j$. Le foncteur
de ces déformations est représentable (cf. \cite{dr1}) par l'anneau $\Def_n^{h,j} \simeq \Def_n^h
[[w_1^n,\cdots,w_j^n]].$}

\begin{defi} Pour tout $1 \leq h \leq d$, soit $\Def_n^{d;=h}$ la $\hat \OC_o^{nr}$-algèbre
qui classifie les déformations de niveau $n$ de $\Sigma_d$ telle que la
composante connexe du $\OC$-module de Dieudonné associé est de dimension $h$. On note $\Def_n^{d;\geq h}$ son
adhérence.
\end{defi}

\rem Concrètement, cf. \cite{boy} \S 9, en notant $M_1$ la matrice modulo $\varpi_o \otimes 1$, de $\vphi$
dans une base du $\OC_o \hat \otimes_{\k} \Def_n^d$-module $V$, $\Def_n^{d;\geq h}$ est le quotient de
$\Def_n^d$ par l'idéal engendré par tous les mineures d'ordre supérieur strictement à $d-h$ de
$$M^{!h}_1:=M_1 (\lexp \t M_1) \cdots (\lexp {\t^{h-1}} M_1).$$

\begin{prop-defi} (cf. \cite{boy} proposition (9.3.2)) On a une suite exacte
$$0 \longto (V_h^{et},\vphi_h^{et}) \longto (V,\vphi) \otimes_{\Def_n^d} \Def_n^{d;=h} \longto (V_h^c,\vphi_h^c) 
\longto 0$$
où $(V_h^{et},\vphi_h^{et})$ (resp. $(V_h^c,\vphi_h^c)$) est un $\OC$-module de Dieudonné étale (resp.
topologiquement nilpotent) de dimension $d-h$ (resp. $h$). La structure de niveau $n$
universelle sur $\Def_n^d$,
$$\iota_n:(\MC^{-n}/\OC_o)^d \longto \MC_n^d$$
où $\MC_n^d$ est l'idéal maximal de $\Def_n^d$, est alors telle que le noyau de
$$\iota_n^{et}:=(\iota_n \otimes_{\Def_n^d} \Def_n^{d;=h})_{|V_h^{et}}:(\MC^{-n}/\OC)^d \longto \MC_n^d 
\otimes_{\Def_n^d} \Def_n^{d;=h}$$
est un facteur direct $a$ de $(\MC_o^{-n}/\OC_o)^d$, de sorte que
$$\Def_n^{d;=h}= \prod_{a \in G/P(d,h,n)} \Def_{n,a}^{d;=h}$$
où $\Def_{n,a}^{d;=h}$ est muni d'une action de $K_{o^n} \cap P_{h,d;a}(\OC_o/\MC_o^n)$.
\end{prop-defi}

\begin{prop} \label{prop-st} (\textbf{Théorème de Serre-Tate} cf. \cite{boy} théorème 7.4.4)
Le morphisme naturel du foncteur des déformations d'un point géométrique de $M_{I,o}$
vers le foncteur des déformations de niveau $n=\mult_o(I)$ du $\OC_o$-module de Dieudonné qui lui est associé, est 
une
équivalence de catégories. En outre l'action d'un \ele $(g_o^c,\frob_o^r) \in GL_h(F_o) \times W_o$ tel que
$r=-\val(\det g_o^c)$ sur $M_{I,o}$, induit l'action de $(\lexp t(g_o^c)^{-1},\frob_o^r) \in \NC_o$ sur $\Def_n^d$.
\end{prop}

\begin{proof} Par rapport à loc. cit., remarquons en effet que l'action de $g_o^c$ sur le $\OC_o$-module de
Dieudonné $(V_o,\vphi_o)$ est donnée par la multiplication à gauche par $g_o^c$, de sorte qu'étant donnée une
structure de niveau $\iota_{o,n}:(\MC_o^{-n}/\OC_o)^d \longto Gr(V_{o,n},\vphi_{o,n})$, l'action de $g_o^c$
est donnée par la multiplication à gauche de $\lexp t g_o^c$ sur $(\MC_o^{-n}/\OC_o)^d$ à comparer avec
l'action définie au paragraphe précédent qui est donnée par la multiplication de $(g_o^c)^{-1}$ à gauche sur
$(\MC_o^{-n}/\OC_o)^d$.

\end{proof}

\begin{coro} \label{coro-def}

\begin{itemize}

\item[(a)] Le complété formel de l'anneau local de $M_{I,o}$ en un point géométrique quelconque de $M_{I,s_o}^{=h}$,
est non canoniquement isomorphe\footnote{Au paragraphe (\ref{enonce-local}), on reviendra plus précisément sur ce 
point.}
 à l'anneau qui représente les déformations d'un $\OC_o$-module divisible de hauteur $d$,
extension d'un $\OC_o$-module formel de hauteur $h$ par sa partie étale,
$\Def_n^{h,d-h}$ où $n$ est la multiplicité de $o$ dans $I$.

\item[(b)] Le complété formel de l'anneau local de
$M_{I,s_o}^{\geq h}$ en un point géométrique de $M_{I,s_o}^{=h'}$ pour $h' \geq h$ est isomorphe par le théorème
de Serre-Tate à $\Def_n^{h',d-h';\geq h}$ défini comme suit avec
$j'=d-h'$:
$$\Def_n^{h',j';\geq h} \simeq \Def_n^{h';\geq h} \hat \otimes_{\bar \Fm_q} \bar \Fm_q[[w_1^n,\cdots,w_{j'}^n]]$$
où $\Def_n^{h';\geq h}$ est le réduit du quotient de $\Def_n^{h'}$ par l'idéal engendré par $\vphi_{0 \to n}(u_i)$ 
pour $1 \leq i < h$ où
$$\vphi_{0 \to n}: \Def_0^{h'}=\bar \Fm_q[[\varpi_o,a_1,\cdots,a_{h'-1}]] \longto \Def_n^{h'}$$
est le morphisme d'oubli du niveau.

\item[(c)] Le complété formel de l'anneau local de $M_{I,s_o,a_h}^{\geq h}$ en un point géométrique de 
$M_{I,s_o,a_{h'}}^{=h'}$ pour
$h' \geq h$ et $a_h \subset a_{h'}$
 est isomorphe par le théorème de Serre-Tate à $\Def_{n,a_h \subset a_{h'}}^{h',d-h';\geq h}$ défini comme suit avec 
$j'=d-h'$:
$$\Def_{n,a_h \subset a_{h'}}^{h',j';\geq h} \simeq \Def_{n,a'_h}^{h';\geq h} \hat \otimes_{\bar \Fm_q} \bar 
\Fm_q[[w_1^n,\cdots,w_{j'}^n]]$$
où $a'_h$ est le facteur direct de rang $h'-h$ de $(\MC_o^{-n}/\OC_o)^{h'}$ associé au quotient $a_{h'}/a_h$ et où
$\Def_{n,a'_h}^{h';\geq h}$ est le réduit
du quotient de $\Def_n^{h'}$ par l'idéal engendré par les $\iota_n(e_i)$ pour $1 \leq i \leq h$ et $(e_1,\cdots,e_h)$
une base de $a'_h$, où  $\iota_n$ est la structure de niveau $n$ universelle sur $\Def_n^{h'}$.

\end{itemize}
\end{coro}

\begin{prop} \label{def-reg}
Pour tout $n \geq 0$ et tout $a \in G/P(d,h,n)$, $(\Def_{n,a}^{d;\geq h})_{red}$ est régulier de dimension $d-h$.
\end{prop}

\begin{proof} Soit $(w_1,\cdots,w_h)$ une base de $a$; pour $1 \leq i \leq h$, on écrit $w_i=\l_1^i
e_1+\cdots+\l_d^i e_d$, où $(e_1,\cdots,e_d)$ est la base canonique de $(\MC_o^{-n}/\OC_o)^d$ et $\l_k^i \in
\OC_o$. Pour tout $1 \leq i \leq h$, l'idéal $(\l_1^i,\cdots ,\l_d^i)$ est égal à $\OC_o$. Notons
$(F,(f_\l)_{\l \in \OC_o})$ le $\OC_o$-module formel de hauteur $d$ universel avec
$$F(X,Y)=X+Y+\cdots \hbox{ et pour tout } \l \in \OC_o, ~f_\l(X)=i(\l)X+\cdots$$
où $i$ est l'injection naturelle $\OC_o \to \Def_n^d$. Soit $\JF_a$ l'idéal de $\bar
\Fm_q[[v_1^n,\cdots,v_d^n]]$ engendré par les \eles
$$F \left ( f_{\l_1^i}(v_1^n),\cdots,F(f_{\l_{d-1}^j}(v_{d-1}^n,f_{\l_d^i}(v_d^n)) \cdots \right ),\qquad 1 \leq i 
\leq h$$
Pour tout $1 \leq i \leq h$, les équations
$$F \left ( f_{\l_1^i}(v_1^n),F(f_{\l_2^i}(v_2^n),\cdots,F(f_{\l_{d-1}^i}(v_{d-1}^n,f_{\l_h^i}(v_d^n)))) \cdots \right 
)=0.$$
s'écrivent sous la forme
$$i(\l_1^j) v_1^n + \cdots +i(\l_d^j) v_d^n + \hbox{ termes de degré } >1.$$
La matrice $(d \times h)$ des $(\bar \l_k^j)_{1 \leq \genfrac{}{}{0pt}{}{k}{j} \leq \genfrac{}{}{0pt}{}{d}{h}
}$ étant de rang $h$, $\bar \Fm_q[[v_1^n,\cdots,v_h^n]]/\JF_A$ est alors isomorphe à $\bar
\Fm_q[[u_1,\cdots,u_{d-h}]]$. En notant $\IC_a$ l'idéal de définition de $\Def_{n,a}^{d;h}$, on a clairement
$\IC_a \subset \JC_a$. En outre $\JC_a \subset \sqrt{\IC_a}$ car d'après la proposition (9.3.3) de
\cite{boy}, pour tout $z \in a$, $(\iota_n \otimes_{\Def_n^d} (\Def_{n,a}^{=h})_{red})(z)$ est nilpotent,
d'où le résultat.

\end{proof}

\begin{coro} \label{coro-reg}
 Les schémas fermés réduits $(M_{I,s_o,a}^{\geq h})_{red}$ sont réguliers pour tout $a \in G/P(d,h,n)$ où $n$ est la
multiplicité de $o$ dans $I$.
\end{coro}

\begin{rema} Soit $Y^h$ une composante irréductible, donc connexe, de $M_{I^o,s_o}^{\geq h}$.
On note\footnote{On ne sait pas à priori si
$Y^{h+i}$ est vide ou pas.} $Y^{h+i}=Y^h \times_{M_{I^o,s_o}^{\geq h}} M_{I^o,s_o}^{\geq h+i}$
et pour tout $n$ la multiplicité de $o$ dans $I$,
$Y^h_n=Y^h \times_{M_{I^o,s_o}^{\geq h}} M_{I,s_o}^{\geq h}=\bigcup_{a \in G/P(d,h,n)} Y_{n,a}^h$.
Pour tout $a,a' \in G/P(d,h,n)$, on a $$Y_{n,a}^h \cap Y_{n,a'}^h=Y_{n,a+a'}^{\dim (a+a')}$$
où l'on considère $a,a'$ comme des facteurs directs de rang $h$ de $(\MC_o^{-n}/\OC_o)^d$.
En particulier si $Y^d$ est non vide, alors les $Y_{n,a}^h$ sont connexes, lisses et donc irréductibles.
\end{rema}

\begin{proof} Au dessus de $M_{I^o,s_o}^{=h}$, on a $Y_n^{=h}=\coprod_{a \in G/P(d,h,n)} Y_{n,a}^{=h}$. En outre au 
dessus de tout point géométrique de
$Y^{=h}$, il y a exactement $|G/P(d,h,n)| \times |GL_{d-h}(\OC_o/\MC_o^n)|$ points géométriques de $Y_n^{=h}$, de 
sorte que
pour tout $a \in G/P(d,h,n)$, $Y_{n,a}^{=h}$ est la réunion disjointes de $|GL_{d-h}(\OC_o/\MC_o^n)|$ composantes 
irréductibles.
L'inclusion
$Y_{n,a}^h \cap Y_{n,a'}^h \subset Y_{n,a+a'}^{\dim (a+a')}$ est évidente. Montrons l'inclusion réciproque.
Soit donc $z \in Y_{n,a+a'}^{\dim (a+a')}$. Le complété de l'anneau local de $Y_n^h$ en $z$ est
$$\Def_{n,a+a'}^{\dim(a+a'),d-\dim (a+a');h}=\bigcup_{b \in G/P(\dim (a+a'),h,n)} \Def_{n,b \subset a+a'}^{\dim 
(a+a'),d-\dim(a+a');h}$$
de sorte que $z \in Y_{n,\tilde a}^h$, pour tout $\tilde a \subset a+a'$.

\end{proof}

\section{Compléments sur la géométrie des strates}

Nous verrons \S \ref{igusa1} et \S \ref{igusa2}, que les relations entre les variétés d'Igusa de première et
seconde espèce avec les variétés $M_{I,o}$ de Drinfeld-Stuhler, sont directement en rapport avec les liens
entre la structure de niveau $\iota_{o,n}$ et la donnée d'une structure de niveau sur la partie étale,
$\iota_{o,n}^{et}$, et connexe, $\iota_{o,n}^c$, séparément, ce qui sera abordée au paragraphe suivant. Pour
définir précisément les parties connexes et étales de $\iota_{o,n}$, on est amené à étudier plus précisément
le scindage en parties connexe et étale des $\OC_o$-modules de Dieudonné sur un $\k(o)$-schéma ou sur un
épaississement infinitésimal.

\subsection{Scindage en partie étale et connexe des $\OC_o$-modules de Dieudonné}

\begin{prop} \label{scinde}
Soit $S$ un $\k(o)$-schéma. Pour tout $S$-point de $M_{I,s_o}^{=h}$, on note $(V_o,\vphi_o)$, le $(\OC_o \hat
\otimes_{\k(o)} \OC_S)$-module de Dieudonné qui lui est associé. Il existe alors un sous-module canonique
$V_o^{et}$ de $V_o$ de rang $(d-h)$ stable sous l'action de $\vphi_o$ tel que

- $\vphi_o^{et}:=(\vphi_o)_{|V_o^{et}}: \lexp \t V_o^{et} \longto V_o^{et}$ est inversible;

- $\vphi_o^c: \lexp \t V_o^c \longto V_o^c:=V_o/V_o^{et}$, l'application induite par $\vphi_o$,
 est topologiquement nilpotente;

- pour tout $n \geq 0$, la suite exacte \footnote{On rappelle que l'indice $n$ signifie que l'on prend la
réduction modulo $\varpi_o^n$: $V_{o,n}=V_o/(\varpi_o^n)$, cf. la définition (\ref{diverses-defis}).}
$$0 \longto (V_{o,n}^{et},\vphi_{o,n}^{et}) \longto
(V_{o,n},\vphi_{o,n}) \longto (V_{o,n}^c,\vphi_{o,n}^c) \longto 0$$
se scinde après le changement de base $S \longmapright{\fr_o^{nh}} S$.
\end{prop}

\begin{proof} Pour tout entier $i \geq 0$, on note $\lexp {\t^i} \vphi_o: \lexp {\t^{i+1}} V_o \longto \lexp {\t^i}
V_o$
l'application $(\OC_o \hat \otimes_{\k(o)} \OC_S)$-linéaire déduite de $\vphi_o: \lexp \t V_o \longto V_o$. On pose
$$\vphi_o^{!h}:=\vphi_o \circ \lexp \tau\vphi_o \circ \cdots \circ \lexp {\tau^{h-1}} \vphi_o: \lexp {\t^h} V_o 
\longto V_o.$$
Soit $ \spec R \longto S$ un ouvert affine de $S$ tel que $V_o \times_S \spec R$ est libre. Si $M_{o,R}$ est
la matrice de $\vphi_o \times_S \spec R$ par rapport à une base $(b_i)_i$ de $V_o$ alors
$$M_{o,R}^{!h}=M_{o,R} ( \lexp \tau M_{o,R}) \cdots (\lexp {\tau^{h-1}} M_{o,R})$$
est celle de $\vphi_{o,R}^{!h}:=\vphi_o^{!h} \times_S \spec R$ par rapport à cette même base $(b_i)_i$.

D'après (\ref{existence-et}) (b), pour tout $h' \geq h$, $M_{o,R,1}^{!h'}$ a un mineure d'ordre $d-h$ inversible
de sorte que\footnote{l'indice $1$ désigne la réduction modulo $\varpi_o$} $\vphi_{o,R,1}^{!h'}(\lexp {\t^{h'}} 
V_{o,R,1})$
 est un facteur direct de $V_{o,R,1}$ de rang $d-h$ car tous les mineures d'ordre supérieure à
$d-h+1$ sont nuls. En remarquant que pour $h_1 \leq h_2$, on a
$$\vphi_{o,R}^{!h_2}(\lexp {\t^{h_2}} V_{o,R}) \subset \vphi_{o,R}^{!h_1}(\lexp {\t^{h_1}} V_{o,R}),$$
on pose pour tout $h' \geq h$
$$V_{o,R,1}^{et}: = \vphi_{o,R,1}^{!h'}(\lexp {\t^{h'}} V_{o,R,1})=\vphi_{o,R,1}^{!h}(\lexp {\t^h} V_{o,R,1}).$$

\begin{lemm} L'application $\vphi_{o,R,1}$ induit une application $\vphi_{o,R,1}^{et}:\lexp \t V_{o,R,1}^{et} \longto 
V_{o,R,1}^{et}$
qui est bijective.
\end{lemm}

\begin{proof} Le diagramme commutatif suivant
$$\diagram
\lexp {\t^{h+1}} V_{o,R} \rrto^{\lexp \tau \vphi_{o,R}^{!h}}
\drrto_{\vphi_{o,R}^{!(h+1)}} & & \lexp \t V_{o,R} \dto^{\vphi_{o,R}}\\
 & & V_{o,R}
\enddiagram $$
permet de définir $\vphi_{o,R,1}^{et}:\lexp \t V_{o,R,1}^{et} \longto V_{o,R,1}^{et}$.
En effet étant donné $v \in V_{o,R,1}^{et}$, il existe $w \in \lexp {\t^h} V_{o,R,1}$ tel que 
$v=\vphi_{o,R,1}^{!h}(w)$,
de sorte que
$$\vphi_{o,R,1}(v \otimes \l)=\vphi_{o,R,1}\circ \lexp \tau \vphi_{o,R,1}^{!h} (w \otimes \l)=\vphi_{o,R,1}^{!h+1}(w 
\otimes \l)\in V_{o,R,1}^{et},$$
d'où l'existence de $\vphi_{o,R,1}^{et}$.

En ce qui concerne la surjectivité de $\vphi_{o,R,1}^{et}$, soit $v=\vphi_{o,R,1}^{!h+1}(w) \in V_{o,R,1}^{et}$
avec $w \in \lexp {\t^{h+1}} V_{o,R,1}$. On a alors $v=\vphi_{o,R,1} \circ \lexp \tau \vphi_{o,R,1}^{!h} (w)$ avec
$\lexp \tau \vphi_{o,R,1}^{!h}(w) \in \lexp \t V_{o,R,1}^{et}$ d'où le résultat.

Pour l'injectivité soit $(e_{1,R}, \cdots,e_{j,R})$ une base de $V_{o,R,1}^{et}$ en tant que $R$-module:
$(e_{i,R} \otimes 1)_{i=1 \cdots j}$ est alors une base de $\lexp \t V_{o,R,1}^{et}$ et $(\vphi(e_{i,R}
\otimes 1))_{i=1 \cdots j}$ est d'après la surjectivité de $\vphi_{o,R,1}^{et}$ une base de $V_{o,R,1}^{et}$.
Ainsi si $v=\sum_{i=1}^j e_{i,R} \otimes \l_j$ appartient à $\ker \vphi_{o,R,1}^{et}$, on en déduit que
$\sum_{i=1}^j \l_i \vphi_{o,R,1}(e_{i,R} \otimes 1)=0$ soit $\l_i=0$ pour $i=1, \cdots, j$.

\end{proof}

On remarque par ailleurs que $(\vphi_{o,R,1}^c)^{!h}:\lexp {\t^h} V_{o,R,1}^c \longto V_{o,R,1}^c$ est nulle.
On définit alors globalement sur $S$, $V_{o,1}^{et}:=\vphi_{o,1}^{!h}(V_{o,1})$ ainsi que
$V_{o,1}^c=V_{o,1}/V_{o,1}^{et}$. Pour tout entier $n$, on définit de même
$$V_{o,n}^{et}:=\vphi_{o,n}^{!nh}(\lexp{\t^{nh}} V_{o,n})$$
qui est un facteur direct de $V_{o,n}$ de rang $n(d-h)$. De la même façon $\vphi_{o,n}^{et}: \lexp \t V_{o,n}^{et} 
\longto V_{o,n}^{et}$
est bijective et $(\vphi_{o,n}^c)^{!nh}$ est nulle.
Pour $n \geq n'$, on a une flèche naturelle
$$V_{o,n}^{et}=\vphi_{o,n}^{!nh}(\lexp{\t^{nh}} V_{o,n})
\longmapright{r_{n',n}}  \vphi_{o,n'}^{!nh}(\lexp{\t^{nh}}
V_{o,n'})=\vphi_{o,n'}^{!n'h}(\lexp{\t^{n'h}} V_{o,n'})=V_{o,n'}^{et}$$
telle que pour $n_1 \geq n_2 \geq n_3$, on a $r_{n_3,n_1}=r_{n_3,n_2} \circ
r_{n_2,n_1}$. On pose alors
$$V_o^{et}:= \lim_{\genfrac{}{}{0pt}{}{\lefto}{n}} V_{o,n}^{et}$$
qui est ainsi stable sous l'action de $\vphi_o$. Soit $V_o^c:=V_o/V_o^{et}$ de sorte que $\vphi_o^c: \lexp \t
V_o^c \longto V_o^c$ est topologiquement nilpotente.

Montrons que le noyau $N_{o,n}$ de $\vphi_{o,n}^{!nh} : \lexp {\t^{nh}} V_{o,n} \longto V_{o,n}$ est un relèvement de
$\lexp {\t^{nh}} V_{o,n}^c$. L'application $\vphi_{o,n}^{et}$ étant injective, on a  $N_{o,n} \cap \lexp {\t^{nh}} 
V_{o,n}^{et} = \{ 0 \}$.
Soit alors $x \in \lexp {\t^{nh}} V_{o,n}$ et $v=\vphi_{o,n}^{!nh}(x)$. D'après la surjectivité de $\vphi_{o,n}^{et}$, 
soit
$w \in \lexp {\t^{nh}} V_{o,n}^{et}$ tel que $v= \vphi_{o,n}^{et,!nh}(w)$. On a alors $x-w \in N_{o,n}$ de sorte que
$$\lexp {\t^{nh}} V_{o,n}= N_{o,n} \oplus \lexp {\t^{nh}} V_{o,n}^{et}.$$
Par ailleurs l'égalité
$\vphi_{o,n}^{!nh} \circ \vphi_{o,n}^{nh}=\vphi_{o,n}^{!nh+1}= \vphi_{o,n} \circ \lexp \tau \vphi_{o,n}^{!nh}$, montre 
que
$\lexp {\tau^{nh}} \vphi_{o,n}: \lexp {\t^{nh+1}} V_{o,n} \longto \lexp {\t^{nh}} V_{o,n}$
induit une application $\vphi_{o,n}^c:\lexp \t N_{o,n} \longto N_{o,n}$, d'où le résultat.

\end{proof}

En ce qui concerne les déformations, on a la proposition suivante.

\begin{coro} \label{coro-defor}
Soit $S$ un $\spec \OC_o$ schéma artinien de réduit $\bar S$ et $S \to M_{I,o}$ tel que $\bar S \to M_{I,o}$ se 
factorise par
$M_{I,s_o}^{=h}$. En notant $(V_{o},\vphi_{o})$
le $\OC_o$-module de Dieudonné associé à $S \to M_{I,o}$, il existe alors un sous-$\OC_o$-module de Dieudonné étale
$(V_{o}^{et},\vphi_{o}^{et})$ tel que $V_{o,n}^{et}$ soit l'image de $\vphi_{o,n}^{!nh}$. De plus, si on définit 
$(V_o^c,\vphi_o^c)$
par la suite exacte
$$0 \longto (V_{o}^{et},\vphi_{o}^{et}) \longto (V_{o},\vphi_{o}) \longto (V_{o}^c,\vphi_{o}^c) \longto 0$$
$\vphi_{o}^c$ est topologiquement nilpotente.
\end{coro}

\begin{proof}  On raisonne localement sur un ouvert affine $\spec R \to S$, pour $R$
une $\OC_o$-algèbre artinienne et $\MC$ un idéal nilpotent de $R$. Classiquement on peut se ramener à $\MC^2=(0)$.
Sur $\bar R$, on choisit une décomposition
$\bar V_o=\bar V_o^{et} \oplus \bar V_o^c$ ainsi que des bases $(\bar e_1, \cdots, \bar e_{d-h})$ de $\bar V_o^{et}$ 
et
$(\bar f_1,\cdots, \bar f_h)$ de $\bar V_o^c$ de sorte que la matrice de $\bar \vphi_o$ par rapport à cette base est
$$\left ( \begin{array}{cc} \bar M_{et} & \bar M_{ext} \\ 0 & \bar M_c \end{array} \right ) $$
On fixe des relèvements quelconques $(e_1,\cdots,e_{d-h})$ et $(f_1,\cdots,f_h)$; la matrice de $\vphi_o$ dans ces 
bases est alors de la forme
$$\left ( \begin{array}{cc} M_{et} & M_{ext} \\ \MC_1 & M_c \end{array} \right )$$
où $\MC_1$ est à coefficient dans $\MC$.
En effectuant le changement de base via la matrice $\left ( \begin{array}{cc} Id & 0 \\ P_0 & Id \end{array} \right 
)$, où $P_0$ est
à coefficients dans $\MC$ de sorte que $\lexp \tau P_0$ est la matrice nulle, on obtient
$$\left ( \begin{array}{cc} M_{et} & M_{ext} \\ P_0 M_{et}+ \MC_1 & P_0 M_{ext}+ M_c \end{array} \right )$$
de sorte qu'en choisissant $P_0=-\MC_1 M_{et}^{-1}$, on obtient une
décomposition $V_o=V_o^{et} \oplus V_o^c$ dans laquelle la matrice de $\vphi_o$ est de la forme
$$\left ( \begin{array}{cc} M_{et} & M_{ext} \\ 0 & M'_c \end{array} \right )$$
avec $M'_c$ topologiquement nilpotente. On remarque par ailleurs que
$V_{o,n}^{et}$ est à nouveau défini globalement sur $S$ comme l'image de $\vphi_{o,n}^{!nh}$.

\end{proof}

\begin{defi}
Dans la suite, dans la situation du corollaire (\ref{coro-defor}), on notera $\FC_o^{*,c}$ le sous-$\OC_o$-module de 
$\FC_o^*$ des
formes qui s'annulent sur $\FC_o^{et}$. Celui ci est stable sous $\phi_o^*$ et on notera $\phi_o^{*,c}$ la 
restriction
de $\phi_o^*$ à $\FC_o^{*,c}$. On remarque ainsi que
$\FC_o^{et,*}$ est isomorphe au quotient $\FC_o^*/\FC_o^{*,c}$ et on note $\phi_o^{et,*}$ l'application induite par 
$\phi_o^*$.
On notera parfois $\FC_o^{*,et}$ un supplémentaire de $\FC_o^{*,c}$, qui est isomorphe à $\FC_o^{et,*}$ mais qui n'est 
pas
stable par $\phi_o^*$.
\end{defi}

\subsection{Retour sur les structures de niveau} \label{struct}

Soit $S$ un $\spec \OC_o$-schéma artinien, $\IC$ un faisceau d'idéaux nilpotents et $\bar S$ le sous-schéma fermé 
associé.
Soit $S \to M_{I,o}$ un $S$-point tel que $\bar S \to M_{I,o}$ se factorise par $M_{I,s_o,a}^{=h}$. On note
$(\FC_{o,n},\phi_{o,n})$ le $\phi$-faisceau
associé sur $S$. D'après le paragraphe précédent, on écrit $\FC_{o,n}^*$ comme une somme directe
$\FC_{o,n}^{*,et} \oplus \FC_{o,n}^{*,c}$ où la matrice de $\phi_{o,n}^*$ dans une base adaptée à cette décomposition; 
est de la forme
$$M_{o,n}^*= \left ( \begin{array}{cc} M_{o,n}^{*,et} & 0 \\ M_{o,n}^{*,ext} &
M_{o,n}^{*,c} \end{array} \right )$$
On note  $\iota_{o,n}:(\MC_o^{-n}/\OC_o)^d \longto \FC_{o,n}^*$ la structure de
niveau $n$ associée.
Pour tout \ele $z \in (\MC_o^{-n}/\OC_o)^d$, on a
$$M_{o,n}^* \iota_{o,n}(z)=\iota_{o,n}(z)^q,$$
\cad en écrivant $\iota_{o,n}(z)=\iota_{o,n}^{et}(z)+\iota_{o,n}^c(z)$:
\begin{eqnarray}
M_{o,n}^{*,et} \iota_{o,n}^{et}(z)=\iota_{o,n}^{et}(z)^q \label{1} \\
M_{o,n}^{*,ext} \iota_{o,n}^{et}(z) + M_{o,n}^{*,c} \iota_{o,n}^c(z)=\iota_{o,n}^c(z)^q \label{2}
\end{eqnarray}

\marque Par définition, l'ensemble $\{ z \in (\MC_o^{-n}/\OC_o)^d ~/~ \iota_{o,n}^{et}(z)=0 \}$ est égal au
sous-module $a$. La condition de Drinfeld s'exprime alors comme suit:

\begin{itemize}

\item[-] $\iota_{o,n}^{et}:(\MC_o^{-n}/\OC_o)^d/a \times S \simeq (\MC_o^{-n}/\OC_o)^{d-h} \times S \longto 
\FC_{o,n}^{et,*}
\simeq \FC_{o,n}^{*,et}$ est un isomorphisme vérifiant (\ref{1});

\item[-] $\iota_{o,n}^c: a \longto \FC_{o,n}^{*,c}$ est une structure de niveau au sens de Drinfeld.

\end{itemize}

\rem Réciproquement si on se donne $\iota_{o,n}^{et}$ et $\iota_{o,n}^c$ comme ci-dessus, le choix d'un supplémentaire
de $a$ permet
de définir une application
$$\iota_{o,n}=\iota_{o,n}^{et} \oplus \iota_{o,n}^c.$$
Le problème est qu'elle ne vérifiera pas obligatoirement la relation (\ref{2}) sauf par exemple si
$M_{o,n}^{ext}$ est nulle. En outre il n'est pas certain que l'on puisse définir $\iota_{o,n}^c$ sur $b$
telle que (\ref{2}) soit vérifiée.

\begin{prop} \label{prop-nulle}
Si $R$ est un $\k(o)$-anneau réduit et $((\EC_i,j_i,t_i),\iota_I)$ un $R$-point de $M_{I,s_o,a}^{=h}$, alors
$\iota_{o,n}^c:a \longto \FC_{o,n}^{*,c}$ est nulle.
\end{prop}

\begin{proof} Soit $\MC$ un idéal maximal de $R$; $\iota_{o,n}^c \otimes_R R/\MC$ est alors nulle, de sorte que 
$\iota_{o,n}^c$ est à coefficients
dans $\bigcap_{\MC \hbox{ maximal}} \MC$ qui est l'idéal nul car $R$ est réduit.
\end{proof}

\begin{prop} \label{degre}
Le morphisme d'oubli du niveau $M_{I,s_o,red}^{=h} \longto M_{I^o,s_o}^{=h}$ est fini et plat de degré
$$\#GL_d(\OC_o/\MC_o^n) / \# GL_h(\OC_o/\MC_o^n)$$
où $n$ est la multiplicité de $o$ dans $I$.
\end{prop}

\begin{proof} La platitude découle de l'équivalent du théorème de Serre-Tate, proposition (\ref{prop-st}), et des
résultats de Drinfeld rappelés au paragraphe (\ref{rapel-def}). Calculons alors les degrés. L'ensemble des
points géométriques de $M_{I,\bar s_o}^{=h}:=M_{I,s_o}^{=h} \times_{\spec \k(o)} \spec \bar \k(o)$ au dessus
d'un point géométrique donné de $M_{I^o,\bar s_o}^{=h}$, est de cardinal
$$\#G/P(d,h,n) \#GL_{d-h}(\OC_o/\MC_o^n).$$
Soit $s_I$ un tel point géométrique de $M_{I,s_o}^{=h}$ au dessus d'un point $s$
de $M_{I^o,s_o}^{=h}$. Soient alors $(M_{I,s_o,red}^{=h})_{s_I}^{\widehat{~}}$ et
$(M_{I^o,s_o}^{=h})_s^{\widehat{~}}$ les complétés formels de respectivement $M_{I,s_o,red}^{=h}$ et
$M_{I^o,s_o}^{=h}$ aux points $s_I$ et $s$. D'après l'équivalent du théorème de Serre-Tate, et avec les
notations du paragraphe (\ref{rapel-def}) et plus particulièrement de la note (\ref{note}), on a
$$(M_{I^o,s_o}^{=h})_s^{\widehat{~}} \simeq \bar \Fm_q [[w_1^0,\cdots ,w_{d-h}^0]] \quad \hbox{ et } \quad
(M_{I,s_o,red}^{=h})_{s_I}^{\widehat{~}} \simeq \bar \Fm_q [[w_1^n,\cdots,w_{d-h}^n]].$$
D'après \cite{dr1}, le degré de
$\bar \Fm_q[[w_1^0,\cdots,w_{d-h}^0]] \longto \bar \Fm_q[[w_1^n,\cdots,w_{d-h}^n]]$ est égal à
$$q^{nh(d-h)}=\#(\OC_o/\MC_o^n)^{h(d-h)},$$
d'où le résultat.

\end{proof}

\begin{coro} \label{degre-a}
Le morphisme $M_{I,s_o,a,red}^{=h} \longto M_{I^o,s_o}^{=h}$ est fini et plat de degré
$$\#(\OC_o/\MC_o^n)^{h(d-h)} \#GL_{d-h}(\OC_o/\MC_o^n)$$
où $n$ est la multiplicité de $o$ dans $I$.
\end{coro}

\marque La proposition suivante justifie l'existence du morphisme de la proposition (\ref{iso-1}). En langage
clair, étant donné un point $\bar S \longto M_{I,s_o,a}^{=h}$ et une déformation de $\bar S \longto
M_{I^o,s_o}^{=h}$, obtenu via le morphisme d'oubli $M_{I,s_o,a}^{=h} \longto M_{I^o,s_o}^{=h}$, ainsi que des
déformations des parties étale et connexe de la structure de niveau en la place $o$, on construit une
déformation canonique de $\bar S \longmapright{\frob_o^{nh}} \bar S \longto M_{I,s_o,a}^{=h}$, où $n$ est la
multiplicité de $o$ dans $I$.

\begin{prop} \label{defor}
Soit $S$ un $\spec \OC_o$-schéma artinien, un sous-faisceau d'idéaux nilpotents $\IC_S$ de $\OC_S$ et $\bar
S$ le fermé de $S$ associé. Soit $\bar S \longto M_{I,s_o,a}^{=h}$ et $S \longto M_{I^o,o}$ une déformation
du morphisme composé
$$\bar S \longto M_{I,s_o,a}^{=h} \longto M_{I^o,s_o}^{=h}.$$
Soit $(V_o^{et},\vphi_o^{et})$ le sous-$\OC_o$-module de Dieudonné étale de $(V_o,\vphi_o)$ associé à $S
\longto M_{I^o,o}$, donné par le corollaire (\ref{coro-defor}). On note $V_o^{*,c}$ le sous-$\OC_o$-module
des \eles de $V_o^*$ qui s'annulent sur $V_o^{et}$. Soient $b$ un supplémentaire de $a$ dans
$(\MC_o^{-n}/\OC_o)^d$, où $n$ est la multiplicité de $o$ dans $I$, et
$$\iota_{o,n}^{et}:b \longto V_{o,n}^{et,*} = V_{o,n}^*/V_{o,n}^{*,c} \quad \hbox{ et } \quad \iota_{o,n}^c:a \longto 
V_{o,n}^{*,c}$$
des déformations respectives de
$$\bar \iota_{o,n}^{et}:b \longto \bar V_{o,n}^{et,*} = \bar V_{o,n}^*/\bar V_{o,n}^{*,c} \quad \hbox{ et }
\quad \bar \iota_{o,n}^c:a \longto \bar V_{o,n}^{*,c}$$
Associée à toutes ces données, il existe alors une déformation canonique de
$\bar S \longmapright{\frob_o^{nh}} \bar S \longto M_{I,s_o,a}^{=h}$
telle que son $\OC_o$-module de Dieudonné $(\tilde V_o,\tilde \vphi_o)$ se scinde modulo $\varpi_o^n$ en partie étale 
et connexe.
\end{prop}

\begin{proof} D'après la proposition (\ref{prop-st}), il suffit de donner une déformation de
$(\lexp {\t^{nh}} {\bar V_o},\lexp {\t^{nh}} {\bar \vphi_o})$
munie de sa structure de niveau $n$, $\lexp {\t^{nh}} {\bar \iota_{o,n}}$.
Soit $(V_o,\vphi_o)$ le $\OC_o$-module de Dieudonné associé à $S \longto M_{I^o,o}$.
D'après le corollaire (\ref{coro-defor}), on a $V_o=V_o^{et} \oplus V_o^c$, où $V_o^{et}$ est défini canoniquement et 
est stable par $\vphi_o$ dont
la matrice dans une base associée à cette décomposition est de la forme
$$ \left ( \begin{array}{cc} M_{o}^{et} & M_o^{ext} \\ 0 & M_{o}^c \end{array} \right )$$
avec $M_o^{et}$ inversible. On considère alors $\tilde V_o=V_o^{et} \oplus V_o^c \varpi_o^n$ défini canoniquement et
muni de l'application $\tilde \vphi_o$ déduite de $\vphi_o$, dont une matrice dans une base adaptée à l'écriture
ci-dessus est
$$\left ( \begin{array}{cc} M_{o}^{et} & \varpi_o^n M_{o}^{ext} \\ 0  & M_{o}^{c} \end{array} \right )$$
de sorte que $(\tilde V_{o,n},\tilde \vphi_{o,n})$ se scinde en partie étale et connexe. En outre comme $\bar
\vphi_o^{!nh}$ induit un isomorphisme $\lexp {\t^{nh}} {\bar V_o} \simeq \bar V_o^{et} \oplus \bar V_o^c
\varpi_o^n$, on en déduit que $(\tilde V_o,\tilde \vphi_o)$ est une déformation de $(\lexp {\t^{nh}} {\bar
V_o},\bar \vphi_o^{nh})$. La structure de niveau $\tilde \iota_{o,n}$ est alors définie comme la composée de
$\iota_{o,n}^c \oplus \iota_{o,n}^{et}$ et de l'inclusion $V_o^* \hookrightarrow \tilde V_o^*$, définition
licite car $(\tilde V_{o,n},\tilde \vphi_{o,n})$ se scinde en partie étale et connexe.

\end{proof}

\begin{rema}
On reprend les notations de la proposition précédente. A $f:S \longto M_{I,s_o,a}^{=h}$, on associe $S
\longto M_{I^o,o}$ ainsi que des déformations $\iota_{o,n}^{et}$ et $\iota_{o,n}^c$ de respectivement $\bar
\iota_{o,n}^{et}$ et $\bar \iota_{o,n}^c$. La proposition précédente fournit alors une déformation $S \to
M_{I,o}$ de $\bar S  \longmapright{\frob_o^{nh}} \bar S \longto M_{I,s_o,a}^{=h}$. Le morphisme ainsi défini
n'est alors rien d'autre que $f \circ \frob_o^{nh}$.
\end{rema}

\begin{defi} \label{defi-f}
On note $\widehat{M_{I,o,=h,a}}$ l'ouvert au dessus de $M_{I,s_o,a}^{=h}$ du complété formel
$\widehat{M_{I,o,h,a}}$ de $M_{I,o}$ le long de $M_{I,s_o,a}^{\geq h}$.
\end{defi}

\subsection{Relation de congruence}

\begin{prop} \label{cong}
L'action d'un \ele $g_o=(g_o^c,1) \in GL_h(F_o) \times GL_{d-h}(F_o)$ (resp. $(\varpi_o,1)$) sur ${\DS
\lim_{\genfrac..{0pt}{1}{\lefto}{I}} M_{I,s_o,1}^{=h}}$ (resp. $\widehat{M_{I,o,=h,1}}$) est donnée par
$\frob_o^{\val(\det g_o^c)}$ (resp. $\frob_o^{h}$).
\end{prop}

\begin{proof} On reprend les notations du paragraphe (\ref{retcor}). Remarquons en premier lieu qu'il suffit de
montrer le résultat pour $g_o^c \in GL_h(F_o) \cap M_h(\OC_o)$. On rappelle que l'action de $W_o$ sur
$M_{I,o}$ est telle que l'image du $\DC$-faisceau elliptique $(\EC_i,j_i,t_i)$ par un frobenius géométrique
est $(\lexp \t \EC_i,\lexp \t j_i,\lexp \t t_i)$. On raisonne sur le $\phi$-faisceau universel
$(\FC_o,\phi_o)$ sur ${\DS \lim_{\genfrac..{0pt}{1}{\lefto}{I}} S_I}$ avec $S_I=M_{I,o}$. A la décomposition
$(F_o/\OC_o)^d= (F_o/\OC_o)^h \oplus (F_o/\OC_o)^{d-h}$ on associe la décomposition
$$(\FC_o^* \otimes_{\OC_o} (F_o/\OC_o)) \otimes_{\OC_S} \OC_{S^o}=\FC_{o,h}^{o} \oplus \FC_{o,d-h}^{o}.$$
L'action de $g_o$ sur $S^o$ est alors induite par celle à droite de $(\lexp t g_o^c,\Id)$ sur $\FC_{o,h}^{o} \oplus
\FC_{o,d-h}^{o}$.

\begin{lemm} Il existe des faisceaux en $\OC_S$-modules $\FC_{o,h}$ et
$\FC_{o,d-h}$ contenus respectivement dans $\FC_{o,h}^{o}$ et $\FC_{o,d-h}^{o}$
tels que
$$\FC_o^* \otimes_{\OC_o} (F_o/\OC_o)=\FC_{o,h} \oplus \FC_{o,d-h}.$$
\end{lemm}

\begin{proof} On rappelle que $S$ étant régulier, $\FC_o^* \otimes_{\OC_o} (F_o/\OC_o)$ s'injecte canoniquement dans 
$(\FC_o^*
\otimes_{\OC_o}(F_o/\OC_o))\otimes_{\OC_S} \OC_{S^o}$. On note $\FC_{o,h}$ (resp. $\FC_{o,d-h}$) le faisceau en
$\OC_S$-modules défini localement pour tout ouvert affine $\spec R \to S$, comme le $R$-module engendré par 
l'ensemble
des $f_{h,R}$ (resp. $f_{d-h,R}$) tels qu'il existe $f_R \in \FC_o^* \otimes_{\OC_o} (F_o/\OC_o)$ tel que $f_R 
\otimes
1=f_{h,R}+f_{d-h,R}$ dans $(\FC_o^* \otimes_{\OC_o} (F_o/\OC_o)) \otimes_{\OC_S} \OC_{S^o}$. Pour tout entier $n$,
$g_o^n(f_R)$ est un \ele de $\FC_{o,R}^*$ qui se décompose sous la forme $f_{d-h,R}+\varpi_o^n f_{h,R}$. Or pour $n$
assez grand, $\varpi_o^n f_{h,R}$ appartient à $\FC_{o,R}^*$ de sorte que $\FC_{o,d-h}=\FC_o^* \cap \FC_{o,d-h}^{o}$ 
et
donc $\FC_{o,h}=\FC_o^* \cap \FC_{o,h}^{o}$ et finalement $\FC_o^*=\FC_{o,h} \oplus \FC_{o,d-h}$.

\end{proof}

\marque On note $((\EC'_i,j'_i,t'_i),\iota')$ le $\DC$-faisceau elliptique muni de sa structure de niveau
infinie, obtenu comme l'image par $g_o$ du $\DC$-faisceau elliptique universel sur $S$,
$(\EC_i,j_i,t_i,\iota)$. Il est alors défini par le diagramme commutatif
$$\diagram
(\FC_{o,h} \dto_{\lexp t g_o^c} & \oplus & \FC_{o,d-h})^d \dto^{\Id} \rto & \EC^* \dto \\
(\FC_{o,h} & \oplus & \FC_{o,d-h})^d \rto & (\EC')^*
\enddiagram$$
En termes de modules de Dieudonné, d'après la proposition (\ref{scinde}), sur $M_{I,s_o,1}^{=h}$, $\FC_{o,d-h}$
correspond au dual de la composante étale $V_o^{et}$ de $V_o$. On note $m= \val(\det g_o^c)$. On rappelle que
$\vphi_o^{!m}: \lexp {\t^m} V_o \longto V_o$ est injective et a pour image $V_o^{et} \oplus (g_o^c. V_o^c)$, de sorte
que sur $M_{I,s_o,1}^{=h}$, $t_i^{!m}: \lexp {\t^m} \EC_i \longto \EC_i$ se factorise en un isomorphisme $\lexp {\t^m} 
\EC_i
\longto \EC_i'$.

\marque En ce qui concerne les structures de niveaux, le seul problème se situe à la place $o$. Soient donc
$\iota_{o,n+r}:(\MC_o^{-n-r}/\OC_o)^d \longto \FC_{o,n+r}^*$ la structure de niveau $n+r$ en $o$ sur
$(\EC_i,j_i,t_i)$ où $r$ est assez grand. On rappelle qu'alors la structure de niveau $n$ en $o$ sur
$(\EC'_i,j'_i,t'_i)$ est définie comme la composée de l'inclusion
$$(\MC_o^{-n}/\OC_o)^d \hookrightarrow (\MC_o^{-n-r}/\OC_o)^h/(\MC_o^{-r}/\OC_o)^h \oplus (\MC_o^{-n}/\OC_o)^{d-h}$$
avec $\iota_{o,n+r} \circ (\times \lexp t g_o)$ et de l'identification de $[g_o^*] \FC_o^*$ avec
$(\FC_o')^*$. Sur la fibre spéciale, la partie connexe de la structure de niveau est nulle, il n'y a donc
rien de plus à vérifier. Pour les déformations, on remarque que la matrice de $\phi_o^{!h}$ est de la forme
$\left (
\begin{array}{cc} \varpi'_o & 0 \\ 0 & \Id \end{array} \right )$ pour une uniformisante $\varpi'_o$ de
$\OC_o$, d'où le résultat.

\end{proof}

\begin{coro}
Soient $(g_o^c,g_o^{et}) \in GL_h(F_o) \times GL_{d-h}(\OC_o)$ et $n \geq m$ assez grand tel que $(g_o^c,g_o^{et})$ 
définit un morphisme
$${\DS M_{I^o\MC_o^n,s_o,1}^{=h} \longmapright{(g_o^c,g_o^{et})} M_{I^o\MC_o^m,s_o,1}^{=h}}.$$
Le diagramme suivant est alors commutatif
$$\diagram
M_{I^o \MC_o^n,s_o,1}^{=h} \rrto^{(g_o^c,g_o^{et})} \dto_{c_1} & & M_{I^o\MC_o^m,s_o,1}^{=h} \dto^{c_1} \\
M_{I^o,s_o}^{=h}  \rrto_{\frob_o^{\val(\det g_o^c)}} & & M_{I^o,s_o}^{=h}
\enddiagram$$
où $c_1$ est le morphisme de restriction du niveau.
\end{coro}

\begin{coro} L'action d'un \ele de $GL_d(\OC_o)$ de la forme
$$\left ( \begin{array}{cc} I_h & 0 \\ * & I_{d-h} \end{array} \right )$$
est triviale sur $(M_{I,s_o,1}^{=h})_{red}$
\end{coro}


%% file: igusa1.tex
\section{Variétés d'Igusa de première espèce}

\label{igusa1}

\subsection{Définition}

Soit $\IC_{I^o,n}^{=h}(S)$ le $\k(o)$-schéma dont les $S$-points sont les $S$-points de $M_{I^o,s_o}^{=h}$
munis d'un isomorphisme
$$\iota_{o,n}^{et}:(\MC_o^{-n}/\OC_o)^{d-h} \times S \longmapright{\sim} \FC_{o,n}^{*,et} \simeq 
\FC_{o,n}^{et,*}=\FC_{o,n}^*/\FC_{o,n}^{*,c}$$
tel que $\phi_{o,n}^{*,et} \circ \iota_{o,n}=\lexp \t \iota_{o,n}$, selon les notations habituelles.

\marque On a alors un morphisme d'oubli de la partie connexe de la structure de niveau
$$i: M_{I,s_o,a}^{=h} \longto \IC_{I^o,n}^{=h}$$
où $n$ est la multiplicité de $o$ dans $I$ et où $i$ est donné par
$$((\EC_i,j_i,t_i),\iota_{I^o},\iota_{o,n}) \longmapsto ((\EC_i,j_i,t_i),\iota_{I^o},\iota_{o,n}^{et}).$$

\begin{prop} \label{igusa1-k}
Pour tout \ele $a$ de $G/P(d,h,n)$,
il existe un morphisme $g_{n,a}:\IC_{I^o,n}^{=h} \longto M_{I,s_o,a}^{=h}$ qui rend
le diagramme ci-dessous commutatif
$$\diagram
M_{I,s_o,a}^{=h} \drto^i \rrto^{\frob_o^{nh}} \ddrto & & M_{I,s_o,a}^{=h} \drto^i \ddrto \\
 & \IC_{I^o,n}^{=h} \urto^{g_{n,a}} \rrto^{\frob_o^{nh}} \dto & & \IC_{I^o,n}^{=h} \dto \\
& M_{I^o,s_o}^{=h} \rrto_{\frob_o^{nh}} & & M_{I^o,s_o}^{=h}
\enddiagram$$
où $n$ est la multiplicité de $o$ dans $I$.
Ainsi le morphisme radiciel $g_{n,a}$ se factorise en un isomorphisme $g_{n,a,red}:\IC_{I^o,n}^{=h} \longto 
M_{I,s_o,a,red}^{=h}$.
\end{prop}

\begin{proof}  Le morphisme de schéma $g_{m,a}$ est défini de la façon suivante. Étant
donné un $S$-point de $\IC_{I^o,m}^{=h}$, il lui correspond un $\DC$-faisceau
elliptique $(\EC_i,j_i,t_i)/S$ muni d'une $I^o$-structure de niveau et
de l'isomorphisme $\iota_{o,m}^{et}/S$. On lui associe alors le $S$-point de
$M_{o^m,a}$:
$$\left \{ \begin{array}{l}
        \lexp {\t^{mh}} (\EC_i,j_i,t_i) \\
        \lexp {\t^{mh}} \iota_{I^o} \\
        \iota_{o,m} \hbox{ défini comme suit.}
\end{array} \right .$$
D'après la proposition \ref{scinde}, $\lexp {\t^{mh}} (\FC_{o,m}^*,\phi_{o,m}^*)$ se scinde en partie étale
et con\-nexe. Soit $b$ un supplémentaire de $a$ dans $(\MC_o^{-m}/\OC_o)^d$. On définit alors la structure de
niveau $m$:
$$\iota_{o,m}:(\MC_o^{-m}/\OC_o)^d=a \oplus b \longto \lexp {\t^{mh}} (\FC_{o,m})^* = \lexp {\t^{mh}} (\FC_{o,m}^c)^* 
\oplus \lexp {\t^{mh}} (\FC_{o,m}^{et})^*$$
comme étant triviale sur $a$ et telle que sa restriction à $b$ soit donnée par $\lexp {\t^{mh}}
(\iota_{o,m}^{et})$. Vu la stabilité de $\lexp {\t^{mh}} (\FC_{o,m}^c)^*$ et de $\lexp {\t^{mh}}
(\FC_{o,m}^{et})^*$ par $\lexp {\t^{mh}} \phi_{o,m}^*$, on définit bien ainsi une $\MC_o^m$-structure de
niveau sur $\lexp {\t^{mh}} (\EC_i,j_i,t_i)$ de telle sorte que la classe d'équivalence de
$((\EC_i,j_i,t_i),\iota_I)$ est indépendante du choix de $b$ et de $(\FC_{o,m}^{et})^*$.

Vérifions alors la commutativité du triangle
$$\diagram
M_{I,s_o,a}^{=h} \drto^i \rrto^{\frob_o^{nh}}  & & M_{I,s_o,a}^{=h} \\
 & \IC_{I^o,n}^{=h} \urto^{g_{n,a}}
\enddiagram$$
de l'énoncé, les autres commutativités en découlant de manière immédiate. Soit $S'=\spec R \to S$ un ouvert
affine tel que $\FC_{o,R}$ est libre et soit $\iota_{o,m}^c / S'$ une structure de niveau $m$ sur
$(\FC_{o,R,m}^c)^*$. On peut choisir une base de $(\FC_{o,R,m}^c)^*$ de sorte que la matrice de
$\phi_{o,m}^{c,*}$ relativement à ce choix soit de la forme
$$\left ( \begin{array}{cccc} 0 & \cdots & \cdots & 0 \\ 1 & \ddots &  & \vdots \\ 0 & \ddots & \ddots & \vdots \\ 0 & 
\cdots & 1 & 0
\end{array} \right )$$
Soit $z \in (\MC_o^{-m}/\OC_o)^h$, on pose $\iota_{o,m}^c(z)=\left ( \begin{array}{c} x_1 \\ \vdots \\ x_h \end{array} 
\right )$ avec
$x_i=\sum_{k=0}^{m-1} x_i^k \varpi_o^k$. La relation
$\left ( \begin{array}{c} x_1^q \\ \vdots \\ x_h^q \end{array} \right ) =\phi_{o,m}^{c,*} \left ( \begin{array}{c} x_1 
\\ \vdots \\ x_h \end{array} \right )$
permet d'exprimer tous les $x_i^k$ en fonction des $(x_1^0)^{q^i}$ pour $0
\leq i \leq mh$. De plus, comme $(M_{o,n}^{ext,!nh})^*$ est nulle d'après (\ref{defor}),
la condition de Drinfeld s'exprime par $(x_1^0)^{q^{mh}}=0$ de sorte que $\lexp {\t^{mh}} (\iota_{o,m}^c)$ est 
triviale.

Le schéma $\IC_{I^o,n}^{=h}$ étant réduit, $g_{n,a}$ se factorise en un morphisme
$g_{n,a,red}:\IC_{I^o,n}^{=h} \longto M_{I,s_o,a,red}^{=h}$ qui est, d'après ce qui précède, une bijection au niveau 
des points
géométriques et tel que le diagramme suivant est commutatif
$$\diagram
\IC_{I^o,n}^{=h} \rrto^{g_{n,a,red}} \dto & & M_{I,s_o,a,red}^{=h} \dto \\
M_{I^o,s_o}^{=h} \rrto^{\frob_o^{nh}} & & M_{I^o,s_o}^{=h}
\enddiagram$$
D'après le corollaire (\ref{degre-a}), $M_{I,s_o,a,red}^{=h} \longto M_{I^o,s_o}^{=h}$ est fini et plat de
degré
$$\#(\OC_o/\MC_o^n)^{h(d-h)} \#GL_{d-h}(\OC_o/\MC_o^n)$$
où $n$ est la multiplicité de $o$ dans $I$. De même
$\IC_{I^o,n}^{=h} \longto M_{I^o,s_o}^{=h}$ est fini et plat de degré $\#GL_{d-h}(\OC_o/\MC_o^n)$ et
$M_{I^o,s_o}^{=h}$ étant régulier de dimension $d-h$, $M_{I^o,s_o}^{=h} \longmapright{\frob_o^{nh}} M_{I^o,s_o}^{=h}$
est de degré $\#(\OC_o/\MC_o^n)^{h(d-h)}$. Comme
$\IC_{I^o,n}^{=h}$ et $M_{I,s_o,a,red}^{=h}$ sont lisses, on en déduit que $g_{n,a,red}$ est un isomorphisme.

\end{proof}

\begin{coro}
Pour tout $m \geq m'$,  on a des morphismes de transitions
$$\IC_{I^o,m}^{=h} \longmapright{\frob_o^{(m-m')h}} \IC_{I^o,m'}^{=h}$$
et des diagrammes commutatifs
$$\diagram \IC_{I^o,m}^{=h} \rrto^{\frob_o^{(m-m')h}} \dto^{g_{m,a}} & & \IC_{I^o,m'}^{=h} \dto^{g_{m',a'}} \\
M_{I^o\MC_o^m,s_o,a}^{=h} \rrto^{c_1} & & M_{I^o\MC_o^{m'},s_o,a'}^{=h} \enddiagram$$ où $c_1$ est le
morphisme de restriction du niveau.
\end{coro}

\subsection{Expression du complété formel de $M_{I,o}$ le long d'une strate}
\label{defi-complete1}

Pour $I$ un idéal de $A$ et $m$ un entier, on considère l'extension étale
$$\hat \IC_{I^o,=h,m} \longto \widehat{M_{I^o,=h}}$$
de fibre spéciale $\IC_{I^o,m}^{=h} \longto M_{I^o,s_o}^{=h}$.
On considère aussi
$$\hat \IC_{I^o,=h,m}(t) \longto \hat \IC_{I^o,=h,m}$$
l'espace classifiant des structures de niveau $t$ sur la partie connexe du $\vphi$-faisceau universel sur $\hat 
\IC_{I^o,=h,m}$.

\begin{prop} \label{iso-1}
Il existe un isomorphisme canonique
$$\hat g_{n,a}:\hat \IC_{I^o,=h,n}(n) \longto \widehat{M_{I,o,=h,a}}$$
qui prolonge le morphisme $g_{n,a}:\IC_{I^o,n}^{=h} \longto M_{I,s_o,a}^{=h}$ défini au paragraphe précédent et tel 
que le diagramme ci-dessous
soit commutatif
$$\diagram
\hat \IC_{I^o,=h,n}(n) \rto^{\hat g_{n,a}} \dto & \widehat{M_{I,o,=h,a}} \dto \\
\widehat{M_{I^o,=h}} \rto^{\frob_o^{nh}} & \widehat{M_{I^o,=h}}
\enddiagram$$
\end{prop}

\begin{proof} La définition de $\hat g_{m,a}$ découle directement de la proposition (\ref{defor}), il ne reste plus
qu'à voir qu'il s'agit d'un isomorphisme. Soit donc $S$ un $\spec \OC_o$-schéma muni d'un sous-faisceau
d'idéaux nilpotents de $\OC_S$ dont le fermé associé est $\bar S$ et $(S \longmapright{f_i}
\IC_{I^o,m}^{=h},\linf i \iota_{o,n}^c)$ pour $i=1,2$, dont les images par $\hat g_{m,a}$ sont des
déformations isomorphes de $\bar S \longto \IC_{I^o,m}^{=h} \longto M_{I,s_o,a}^{=h}$. Notons $(\linf i V_o,
\linf i \vphi_o)$ les $\OC_o$-modules de Dieudonné sur $S$ associés à $f_i$ pour $i=1,2$. Par hypothèse, on a
un isomorphisme
$$h: \linf 1 {\tilde V_o}=\linf 1 V_o^{et} + \varpi_o^n \linf 1 V_o^c \longmapright{\sim} \linf 2 {\tilde V_o}=\linf 2 
V_o^{et} + \varpi_o^n \linf 2 V_o^c$$
tel que $\bar h$ est l'identité et commute aux actions de $\linf i \vphi_o$. Clairement $h$ induit un
isomorphisme $h^{et}: \linf 1 V_o^{et} \longmapright{\sim} \linf 2 V_o^{et}$. Soit donc $v_1 \in \linf 1
V_o^c$: $h(\varpi_o^n v_1)=\varpi_o^n v_2 + w$ avec $v_2 \in \linf 2 V_o^c$ et $w \in \linf 2 V_o^{et}$. On a
$$\begin{array}{rlr}
h(\linf 1 \vphi_o^{!nh}(\varpi_o^n v_1)) & =h(\varpi_o^n(\varpi_o^n u_1+u_2)) & u_1 \in \linf 1 V_o^c,~ u_2 \in \linf 
1 V_o^{et} \\
 & = \varpi_o^n w' & w' \in \linf 2 V_o \\
 & = \linf 2 \vphi_o^{!nh}(\varpi_o^n v_2 + w)
\end{array}$$
soit $\linf 2 (\vphi_o^{et})^{!nh}(w) \in \varpi_o^n \linf 2 V_o^{et}$, soit $w \in \varpi_o^n \linf 2
V_o^{et}$. On en déduit donc un morphisme $\linf 1 V_o^c \longto \linf 2 V_o$ et finalement $h$ provient d'un
isomorphisme $\linf 1 V_o \simeq \linf 2 V_o$. Si en outre $\iota_{o,n,1}$ et $\iota_{o,n,2}$ sont
isomorphes, il est clair du fait que les $\linf i {\tilde V_{o,n}}$, pour $i=1,2$, se scindent en partie
étale et connexe, que les triplés $(\linf i V_o, \linf i \iota_{o,n}^{et}, \linf i \iota_{o,n}^c)$ pour
$i=1,2$, sont isomorphes. Pour des raisons de dimension, $\hat g_{m,a}$ induit alors un isomorphisme sur les
espaces tangents, c'est donc un isomorphisme.

\end{proof}

\medskip

\rem Soit $i:\widehat{M_{I,o,=h,a}} \longto \hat \IC_{I^o,=h,n}(n)$ la flèche donnée par
$$i((\EC_i,j_i,t_i),\iota_{I^o},\iota_{o,n})=((\EC_i,j_i,t_i),\iota_{I^o},\iota_{o,n}^{et},\iota_{o,n}^c)$$
qui prolonge celle définie au paragraphe précédent. Le composé
$$\widehat{M_{I,o,=h,a}} \longto \hat \IC_{I^o,=h,n}(n) \longmapright{\hat g_{n,a}} \widehat{M_{I,o,=h,a}}$$
est $\frob_o^{nh}$.

\subsection{Correspondances de Hecke}
\label{cores-1}

\begin{prop} La tour $(\IC_{I^o,n}^{=h})_n$ est munie de correspondances de Hecke associées aux \eles de 
$GL_{d-h}(F_o) \times \Zm$,
compatibles aux morphisme $g_{m,a}$, \cad que
pour tout $n \geq 0$, et tout \ele $(g_o^{et},r) \in GL_{d-h}(F_o) \times \Zm$, il existe $m_0 \geq n$ ainsi que des 
morphismes
$\IC_{I^o,m}^{=h} \longmapright{(g_o^{et},r)} \IC_{I^o,n}^{=h}$ pour tout $m \geq m_0$ compatibles aux morphismes de 
restriction
du niveau, tels que le diagramme suivant soit commutatif
$$\diagram
\IC_{I^o,m}^{=h} \rrto^{(g_o^{et},r)} \dto^{g_{m,a_m}} & & \IC_{I^o,n}^{=h} \dto^{g_{n,a_n}} \\
M_{I^o \MC_o^m,s_o,a_m}^{=h} \rrto^{(g_o^{et},g_o^c)} & & M_{I^o\MC_o^n,s_o,a_n}^{=h}
\enddiagram$$
pour $a_n,~a_m$ des \eles de $G/P(d,h,n)$ et $G/P(d,h,m)$ respectivement tels que $a_n$ soit l'image de $a_m$ par la 
surjection
canonique, et où $g_o^c$ est un \ele quelconque de $GL_h(F_o)$ tel que $r=\val(\det g_o^c)$.
\end{prop}

\begin{proof} On définit les correspondances de Hecke sur $\IC_{I^o,m}^{=h}$
associées aux \eles de $GL_{d-h}(F_o) \times \Zm$, en procédant comme suit.

- Soit $g_o^{et} \in GL_{d-h}(F_o)$ tel que $g_o^{et} \in \Mm_{d-h}(\OC_o)$. On choisit $m$ et $n$ tels que
le noyau de $\lexp t g_o^{et}:(F_o/\OC_o)^{d-h} \longto (F_o / \OC_o)^{d-h}$ est contenu dans
$(\MC_o^{-m}/\OC_o)^{d-h}$ et tels que $(\MC_o^{-n}/\OC_o)^{d-h}$ est contenu dans
$$\im (\lexp t g_o^{et}:(\MC_o^{-m}/\OC_o)^{d-h}) \to (\MC_o^{-m}/\OC_o)^{d-h}).$$
Soient $S=\IC_{I^o,m}^{=h}$ et $((\EC_i,j_i,t_i),\iota_{I^o},\iota_{o,m}^{et})$ l'objet universel de
$\IC_{I^o,m}^{=h}$. A $g_o^{et}$ on associe via
la structure de niveau $\iota_{o,m}^{et}$, le morphisme $[g_o^{et}]: \EC_{o,m}^{et} \to \EC_{o,m}^{et}$ défini comme 
au
paragraphe (\ref{retcor}) du chapitre précédent.

- L'action de $(g_o^{et},-k)$ pour $g_o^{et} \in \Mm_{d-h}(\OC_o)$ et $k$ un entier naturel assez grand tel
que $\varpi_o^{[k/h]} (g_o^{et})^{-1} \in \Mm_{d-h}(\OC_o)$, de sorte que $\varpi_o^{[k/h]} \EC_{o,m}^{et}
\subset [g_o^{et}] \EC_{o,m}^{et}$, est alors définie comme suit. Soit $\EC_i'$ défini par le carré cartésien
$$\diagram \EC_i' \rto \dto & \EC_{o,n}^{et} + (\vphi_o^c)^{!k}(\lexp {\t^k} \EC_{o,n}^c)
\dto^{[g_o^{et}] + i_o^c}  \\ \EC_i \rto & \EC_{o,m} =\EC_{o,m}^{et} \oplus \EC_{o,m}^c
\enddiagram$$
où $i_o^c$ est l'inclusion canonique.

On rappelle que modulo $\EC_{o,m}^{et}$, on a $(\vphi_o^c)^{!h}(\lexp {\t^h} \EC_o^c) =\varpi_o \EC_o^c$,
et donc l'image par $\vphi_o$ de $\lexp \t ((\vphi_o^c)^{!k}(\lexp {\t^k} \EC_{o}^c)$
est incluse dans $[g_o^{et}](\EC_{o,n}^{et}) + (\vphi_o^c)^{!k}(\lexp {\t^k} \EC_{o,n}^c))$, ce qui permet de définir 
des
applications $t'_i:\lexp \t \EC_i' \longto \EC_{i+1}'$ déduites des $t_i$. Les
$\EC_i'$ peuvent ainsi être organisés en un $\DC$-faisceau elliptique  et on
peut de plus le munir de la $\MC_o^n$-structure de niveau $\iota_{o,n}^{'~et}$
sur la partie étale déduite de $\iota_{o,m}^{et}$. Le $\DC$-faisceau
elliptique image par $(g_o^{et},k)$ de
$((\EC_i,j_i,t_i),\iota_{I^o},\iota_{o,m}^{et})$ est alors
$((\lexp {\t^{m-n}} \EC_i',\lexp {\t^{m-n}} j_i',\lexp {\t^{m-n}} t_i'),\lexp {\t^{m-n}} \iota_{I^o}', \lexp 
{\t^{m-n}} \iota_{o,n}^{'~et})$.

- On fait agir un \ele $\varpi_o^k,~k \geq 0$ du centre de $GL_{d-h}(F_o)$ en tordant $\EC_o^{et}$ par le
faisceau inversible $\OC_X(k.o)$, \cad avec des notations similaires à celles introduites ci-dessus,
$\EC_o':=\varpi_o^{k} \EC_o^{et} \oplus \EC_o^c$. On peut de manière évidente définir des applications $t'_i$
et donc organiser les $\EC_i'$ en un $\DC$-faisceau elliptique muni d'une $\MC_o^n$-structure de niveau sur
la partie étale.

Soient $r \in \Zm$ et $m \geq n$ tels que $m-n-r \geq 0$. On définit le
morphisme $[r]: \IC_{I^o,m}^{=h} \longto \IC_{I^o,n}^{=h}$ comme étant égal à $\frob_o^{m-n-r}$.

En remarquant que pour $k \geq 0$, $\vphi_o^{!k}$ induit un isomorphisme entre $\lexp {\t^k} (\EC_o^{et}
\oplus \EC_o^c)$ et $\EC_o^{et} + (\vphi_o^c)^k (\EC_o^c)$, il est facile\footnote{essentiellement il s'agit
de vérifier que l'action de $(g_o^{et},k+s)$ est égale à celle de $(g_o^{et},k)(1,s)$ pour $(g_o^{et})^{-1}
\in \Mm_d(\OC_o)$, $s \geq 0$ et $k$ assez grand} de vérifier que l'on définit bien ainsi une action de
$GL_{d-h}(F_o) \times \Zm$ sur ${\DS \lim_{\genfrac{}{}{0pt}{}{\lefto}{m}}} (\IC_{I^o,m}^{=h})$.

\medskip

En vertu de la proposition (\ref{cong}) du chapitre précédent, il est immédiat de vérifier la commutativité du 
diagramme de l'énoncé.

\end{proof}

\begin{prop} La tour $(\hat \IC_{I^o,=h,n}(n))_n$ est munie de correspondances de Hecke associées aux \eles de
$GL_{d-h}(F_o) \times GL_h(F_o)$,
compatibles aux morphisme $\hat g_{m,a}$, \cad que
pour tout $n \geq 0$, et tout \ele $(g_o^{et},g_o^c) \in GL_{d-h}(F_o) \times GL_h(F_o)$, il existe $m_0 \geq n$ ainsi 
que des morphismes
$$\hat \IC_{I^o,=h,m}(m) \longmapright{(g_o^{et},g_o^c)} \hat \IC_{I^o,=h,n}(n)$$
pour tout $m \geq m_0$ compatibles aux morphismes de restriction du niveau, tels que les diagrammes suivants
soient commutatifs
$$\diagram \hat \IC_{I^o,=h,m}(m) \rrto^{(g_o^{et} \times g_o^c)} \dto & & \hat \IC_{I^o,=h,n}(n) \dto \\
\IC_{I^o,m}^{=h} \rrto_{(g_o^{et},\val (\det g_o^c))} & & \IC_{I^o,n}^{=h} \\
\enddiagram$$
$$\diagram \hat \IC_{I^o,=h,m}(m) \rrto^{(g_o^{et} \times g_o^c)} \dto^{\hat g_{m,a}} & & \hat \IC_{I^o,=h,n}(n) 
\dto^{\hat g_{n,a'}} \\
\widehat{M_{I,o,h,a}} \rrto_{\left ( \begin{array}{cc} g_o^{et} & 0 \\ 0 & g_o^c \end{array} \right ) } & & 
\widehat{M_{I,o,h,a'}}
\enddiagram$$
\end{prop}

\begin{proof} On définit une correspondance de Hecke associée à un \ele $(g_o^{et},g_o^c)$ de
$GL_{d-h}(F_o) \times GL_h(F_o)$ sur $\hat \IC_{I^o,=h,m}(m)$, de telle sorte
que l'action induite sur $\IC_{I^o,m}^{=h}$ est donnée par
$(g_o^{et},\val (\det g_o^c)) \in GL_{d-h}(F_o) \times \Zm$ telle qu'elle est définie ci-dessus.

- Soit $g_o^c$ (resp. $g_o^{et}$) un \ele de $GL_h(F_o)$ (resp. de $GL_{d-h}(F_o)$) tel que
$$g_o^c \in \Mm_h(\OC_o)$$
$$(resp.~g_o^{et} \in \Mm_{d-h}(\OC_o)) \hbox{ et } \varpi_o^{\val(\det g_o^c)/h} (g_o^{et})^{-1} \in 
\Mm_{d-h}(\OC_o)$$
de sorte que $\varpi_o^{[\val(\det g_o^c)/h]} \EC_{o,m}^{et} \subset [g_o^{et}] \EC_{o,m}^{et}$. On choisit $m$ et 
$n$
tels que les noyaux de $\lexp t g_o^c:(F_o/\OC_o)^h \longto F_o/\OC_o)^h$ et $\lexp t g_o^{et}:(F_o/\OC_o)^{d-h} 
\longto
(F_o/\OC_o)^{d-h}$ sont contenus dans $(\MC_o^{-m}/\OC_o)^h$ et tels que $(\MC_o^{-n}/\OC_o)^h$ (resp.
$\MC_o^{-n}/\OC_o)^{d-h}$) est contenu dans
$$\im (\lexp t g_o^c:(\MC_o^{-m}/\OC_o)^h \longto \MC_o^{-m}/\OC_o)^h)$$
(resp. dans $\im (\lexp t g_o^{et}:(\MC_o^{-m}/\OC_o)^{d-h} \longto (\MC_o^{-m}/\OC_o)^{d-h})$). Soient $S$ un $(\hat
\OC_o^{nr})$-schéma et ainsi qu'un $S$-point $((\EC_i,j_i,t_i),\iota_{I^o},\iota_{o,m}^c,\iota_{o,m}^{et})$ de $\hat
\IC_{I^o,=h,m}(m)$. Pour définir l'action du couple $(g_o^{et},g_o^c)$, il suffit de définir une structure de niveau 
$n$
sur $(\vphi_o^c)^{!(\val (\det g_o^c))}(\EC_o^c) \simeq \lexp {\t^{\val (\det g_o^c)}} \EC_o^c$. Soit $[g_o^c]:
\EC_{o,m}^c \to \EC_{o,m}^c$  le morphisme associé à $g_o^c$ défini via $\iota_{o,m}^c$, la structure de niveau $m$ 
(cf.
le paragraphe (\ref{retcor})). On a déjà vu que $[g_o^c](\EC_{o,m}^c) \simeq \lexp {\t^{\val (\det g_o^c)}} 
\EC_{o,m}^c$
de sorte que la structure de niveau $n$ déduite de $\iota_{o,m}^c$ convient.

- L'action d'un \ele du centre de $GL_{d-h}(F_o)$ est définie de manière évidente. De même l'action de
$\varpi_o^{k}$ vu comme \ele du centre de $GL_h(F_o)$ se définit de manière naturelle $\hat \IC_{I^o,=h,m}(m)
\longto \hat \IC_{I^o,=h,n}(n)$ pour $m$ et $n$ tels que $m-n+hk \geq 0$, et relève
$$\frob_o^{m-n+hk}: \IC_{I^o,m}^{=h} \longto \IC_{I^o,n}^{=h}.$$
On vérifie aisément que l'on définit bien ainsi une action de $GL_{d-h}(F_o) \times GL_h(F_o)$ sur ${\DS
\lim_{\genfrac{}{}{0pt}{}{\lefto}{m}}} ~\hat \IC_{I^o,=h,m}(m)$ de telle sorte que l'on a le diagramme
commutatif suivant
$$\diagram \hat \IC_{I^o,=h,m}(m) \rrto^{(g_o^{et} \times g_o^c)} \dto & & \hat \IC_{I^o,=h,n}(n) \dto \\
\IC_{I^o,m}^{=h} \rrto_{(g_o^{et},\val (\det g_o^c))} & & \IC_{I^o,n}^{=h} \\
\enddiagram$$
En vertu de la proposition (\ref{defor}) et de la définition de $\hat g_{m,a}$, le diagramme suivant est alors 
commutatif
$$\diagram \hat \IC_{I^o,=h,m}(m) \rrto^{(g_o^{et} \times g_o^c)} \dto^{\hat g_{m,a}} & & \hat \IC_{I^o,=h,n}(n) 
\dto^{\hat g_{n,a'}} \\
\widehat{M_{I,o,h,a}} \rrto_{\left ( \begin{array}{cc} g_o^{et} & 0 \\ 0 & g_o^c \end{array} \right ) } & & 
\widehat{M_{I,o,h,a'}}
\enddiagram$$

\end{proof}

\section{Variétés d'Igusa de seconde espèce}

\label{igusa2}

On rappelle que pour un idéal $I$ de $A$, on note $I=I^o \MC_o^n$ avec $o \not \in V(I^o)$. Soit $\Pi_{o,d}$
l'\ele suivant de $GL_d(\OC_o)$
$$\left ( \begin{array}{cccc}
0 & \cdots & 0 & \varpi_o \\ 1 & 0 & \vdots & 0 \\ 0 & \ddots & & \vdots \\ 0 &  0 & 1 & 0
\end{array} \right ) $$

\subsection{Définition}
\label{defi-2}

Étant donnés un $\k(o)$-schéma $S$, ainsi qu'un $S$-point de $\IC_{I^o,m}^{=h}$, on a un $\OC_o \otimes_{\k(o)} 
\OC_S$-module
localement libre de rang $d$, $\FC_o$ ainsi qu'un morphisme linéaire injectif $\phi_o: \lexp \t \FC_o \longto \FC_o$.

Considérons la catégorie fibrée $\JF_{I^o,m}^{=h}(s)$ sur la catégorie des
$\k(o)$-schémas dont les objets sont
$$S \mapsto \left \{ \begin{array}{l}
.~ S \longto \IC_{I^o,m}^{=h} \\
.~ \hbox{une section globale } \s \hbox{ de } \FC_o/ \phi_o^{!(s)}(\lexp {\t^{s}} \FC_o) \hbox{ telle que } \\
\phi_{o}^{!h}(\s \otimes 1)=\varpi_o \s \hbox{ et telle que l'application } \OC_S \longto \FC_o/\phi_o(\lexp \t 
\FC_o)\\
a \mapsto a \bar \s, \hbox{ est inversible}
\end{array} \right. $$
Dans la suite, on notera $\FC_{o,/s}$ le quotient $\FC_o/ \phi_o^{!s}(\lexp {\t^s} \FC_o)$.

\begin{prop} \label{propdefi2}
 Le morphisme naturel $\JF_{I^o,m}^{=h}(s) \longto \IC_{I^o,m}^{=h}$ est relativement représentable par un schéma
$\JC_{I^o,m}^{=h}(s) \longto \IC_{I^o,m}^{=h}$ étale galoisien de groupe de Galois
$\DC_{o,h,s}^\times:=(\DC_{o,h}/(\Pi_{o,h}^{s}))^\times$. En outre le diagramme suivant est cartésien
$$\diagram \JC_{I^o,m}^{=h}(s) \rto \dto & \JC_{I^o,0}^{=h}(s) \dto \\ \IC_{I^o,m}^{=h} \rto & M_{I^o,s_o}^{=h} 
\enddiagram$$
\end{prop}

\begin{proof} La représentabilité est immédiate et le diagramme est clairement cartésien.
Soit alors $\spec R \to S$ un ouvert affine tel que $\FC_{o,R}:=\FC_o \otimes_{\OC_S} R$
est libre.
D'après la proposition (\ref{scinde}), il existe une décomposition
$$\FC_{o,R}=\FC_{o,R}^{c} \oplus \FC_{o,R}^{et}$$
dans laquelle la matrice de $\phi_o$ est
de la forme $\left (  \begin{array}{cc} M_c & 0 \\ M_{ext} & M_{et} \end{array} \right ) $
avec $M_{et}$ inversible de taille $d-h$ de sorte que $\FC_{o,R}^{et} \subset \phi_o^{!(s+1)}(\lexp {\t^{s+1}} 
\FC_{o,R})$
pour tout $s >0$. En outre $\FC_{o,R}/\phi_o(\lexp \t \FC_{o,R})$ étant un $R$-module libre de rang $1$,
il existe donc $e_0 \in \FC_{o,R}^{c}$ tel que
$$(e_0, \phi_o(e_0\otimes 1), \cdots, \phi_o^{!(h-1)}(e_0 \otimes 1))$$
soit une base de $\FC_{o,R}^c$ avec
$$\phi_o^{!h}(e_0 \otimes 1)= \varpi_o\otimes 1 \sum_{i=0}^{h-1} \g_i \phi_o^{!i}(e_0 \otimes 1) \quad \hbox{ modulo } 
\FC_o^{et}$$
où $\g_i \in \OC_o \hat \otimes_{\k(o)} R$. Il s'agit ainsi de trouver $\s=\sum_{i=0}^{h-1} \a_i.
\phi_o^{!i}(e_0 \otimes 1)$ avec $\a_i ={\DS \sum_{\genfrac{}{}{0pt}{}{0 \leq j}{i+j.h \leq s}} \b_i^j
\varpi_o^j}$ telle que $\phi_o^{!h}(\s \otimes 1)=\varpi_o \s$ et $\b_0^0$ inversible. Notons
$$(W_{h,S,o},\phi_{h,S,o})$$
le $\phi$-faisceau $W_{h,S,o}:=(\OC_o \otimes_{\k(o)} \OC_S)^h$ tel que la matrice de
$$\phi_{h,S,o}:~\lexp \t W_{h,S,o} \longto W_{h,S,o}$$
dans la base canonique est donnée par
$\phi_{h,S,o}(e_i \otimes 1)=e_{i+1}$ pour $1 \leq i < h$ et $\phi_{h,S,o}(e_h \otimes 1)=\varpi_o.  e_1$.
Le groupe $\aut(W_{h,S,o},\vphi_{h,S,o})$ est
$$\DC_{o,h}^\times \simeq \{ P \in GL_h(\OC_o)~/~P \Pi_{o,h}=\Pi_{o,h} P \}$$
Trouver $\s$ revient alors à \textit{rigidifier la partie connexe de $\FC_o$}, i.e. à donner un isomorphisme
$$(\FC_{o}^c,\phi_{o}^c) \simeq (W_{h,S,o},\phi_{h,S,o})$$
modulo $\Pi_{o,h}^{s}$. On peut ainsi obtenir aisément des équations du revêtement à partir des \eles de la
matrice $M_c$. On en déduit en outre que le revêtement est étale de groupe de Galois:
$$\aut(W_{h,S,o,/s},\phi_{h,S,o,/s}) \simeq \{ P \in GL_h(\OC_o)/(1+\Pi_{o,h}^{s} \Mm_h(\OC_o)~/~P
\Pi_{o,h}= \Pi_{o,h}P \}$$

\end{proof}

\rem Si $S$ est le spectre d'un corps ou d'un anneau artinien, il correspond, par la théorie du module de coordonnées 
(cf. \cite{gen}), à
$(V_o^c,\vphi_o^c)$ un $\OC_o$-module formel de hauteur $h$, de sorte que les $S$-points de $\JC_{I^o,m}^{=h}(s)$ sont 
ceux de
$\IC_{I^o,m}^{=h}$ muni d'une rigidification à l'ordre $s$ du $\OC_o$-module formel de hauteur $h$ qui leur est 
associé. Ainsi
$\JC_{I^o,m}^{=h}(s)$ est la variété d'Igusa de seconde espèce telle qu'elle est définie dans \cite{h-t}.

\subsection{Correspondances de Hecke sur $\JC_{I^o,n}^{=h}(s)$}
\label{cores-j}

\begin{prop} La tour $(\JC_{I^o,n}^{=h}(s))_{s,n}$ est munie de correspondances de Hecke associées aux \eles
du noyau $\widetilde{\NC_o}$ de l'application
$$\begin{array}{cl} GL_h(F_o) \times D_{o,h}^\times \times W_o & \longto \Zm \\
(g_o^c,\d_o,\s_o) & \mapsto \val (\det(g_o^c) \rn(\d_o) \cl(\s_o))
\end{array}$$
où $\rn:D_{o,h}^\times \longto F_o^\times$ est la norme réduite, et $\cl:W_{F_o} \longto F_o^\times$ est
l'application de la théorie du corps de classe, de manière compatible aux morphismes $\JC_{I^o,n}^{=h}(s)
\longto \IC_{I^o,n}^{=h}$, i.e. pour tout $n,s$, il existe $m_0 \geq n$ et $t_0 \geq s$ tels que pour tout $m
\geq m_0$ et $t \geq t_0$, on ait des morphismes
$$\diagram \JC_{I^o,m}^{=h}(t) \rrto^{g_o^{et} \times g_o^c \times \d_o \times \s_o} & & 
\JC_{I^o,n}^{=h}(s)\enddiagram$$
compatibles aux morphismes de restriction du niveau et tels
que le diagramme suivant soit commutatif:
$$\xymatrix {
\JC_{I^o,m}^{=h}(t) \ar^{g_o^{et} \times g_o^c \times \d_o \times \s_o}[rrr] \ar[d] & & & \JC_{I^o,n}^{=h}(s) \ar[d] 
\\
\IC_{I^o,m}^{=h} \ar_{g_o^{et} \times (\val(\det(g_o^c)))}[rrr] & & & \IC_{I^o,n}^{=h} } $$
\end{prop}

\begin{proof}
- Les correspondances géométriques sur $\IC_{I^o,n}^{=h}$ associées aux \eles de $(D_\Am^{\oo,o})^\times
\times GL_h(\OC_o) \times GL_{d-h}(F_o) $ se remontent aisément sur $\JC_{I^o,n,o}^{=h}(m)$.

- Celles associées aux \eles de $\DC_{o,h}^\times$ sont données via l'isomorphisme $\aut(V_o^c,\vphi_o^c)
\simeq \DC_{o,h}^\times$, de sorte que
$$\JC_{I^o,n}^{=h}(s)= ({\DS \lim_{\genfrac{}{}{0pt}{}{\lefto}{s'}}}~ 
\JC_{I^o,n}^{=h}(s'))^{(1+\Pi_{o,h}^{s}\DC_{o,h})}.$$

- Soit $(g_o^c,\d_o,1)$ un \ele de $\widetilde{\NC_o}$ avec $g_o^c \in \Mm_h(\OC_o)$. On choisit $m \geq n$
tels que
$$\ker (\lexp t g_o^{et}) \subset (\MC_o^{-m}/\OC_o)^{d-h},$$
$$(\MC_o^{-n}/\OC_o)^{d-h} \subset  \im (\lexp t g_o^{et})$$
$$\varpi_o^{\val (\det g_c)/h} (g_o^{et})^{-1} \in \Mm_{d-h}(\OC_o),$$
de sorte que $(g_o^{et},g_o^c)$ définisse une correspondance de Hecke
$$\IC_{I^o,m}^{=h} \longto \IC_{I^o,n}^{=h}.$$
Il reste alors à définir, pour
$t$ assez grand, à partir d'un isomorphisme
$$\a_{t}:(\FC_{o,/t}^c,\phi_{o,/t}^c) \simeq (W_{h,S,o,/t},\phi_{h,S,o,/t}),$$
un isomorphisme
$$\a'_{s}:(\FC_{o,/s}^{c,'},\phi_{o,/s}^{c,'}) \simeq (W_{h,S,o,/s},\phi_{h,S,o,/s}),$$
où $(\FC_o^{c,'},\phi_o^{c,'})$ est le $\phi$-faisceau connexe associé à $(\FC_o^c,\phi_o^c)$ par $\lexp t
g_o^c$. Rappelons que l'application $(\phi_o^{c})^{!\val(\det(g_o^c))}$ (resp.
$\phi_{h,S,o}^{!\val(\rn(\d_o^{-1}))}$) réalise un isomorphisme de $\phi$-faisceau de $\lexp
{\t^{\val(\det(g_o^c))}} \FC_o$ (resp. $\lexp {\t^{\val(\rn(\d_o^{-1}))}} W_{h,S,o}$) vers $\FC_o^{c,'}$
(resp. $\d_o^{-1}(W_{h,S,o})$), de sorte que
$$\d_o \circ \phi_{h,S,o}^{!\val(\rn(\d_o))} \circ \a_{t} \circ (\phi_o^{!\val(\det(g_o^c))})^{-1}$$
restreint à $\FC_{o,/t-\val(\det(g_o^c))}^{c,'}$ définit une rigidification à l'ordre $s:=t-\val(\det(g_o^c))$ de la
partie connexe du $\phi$-faisceau $(\FC_o',\phi_o')$.

- pour $z \in F_o^\times \cap \OC_o$, on fait agir $(z,z^{-1},1)$ de manière évidente, \cad avec les
notations ci-dessus, $\FC_o^{c,'}=z \FC_o^c$ avec $\a'_s=z \circ \phi_{h,S,o}^{!h \val(z)} \circ \a_{t} \circ
(\phi_o^{!h \val(z)})^{-1}$.

- Le procédé est exactement identique pour les \eles de $\widetilde{\NC_o}$ de la forme $(1,\d_o,c_o)$, en
remarquant que $c_o \d_o$ est un automorphisme de $(W_{h,S,o,/s},\phi_{h,S,o,/s})$.

\end{proof}

\subsection{Complétés formels}
\label{enonce-local}

On se sert des $\JC_{I^o,m}^{=h}(s)$ pour ``détordre'' le schéma formel $\hat \IC_{I^o,=h,m}(t)$. Le lemme suivant est 
immédiat:

\begin{lemm} \label{lem-iso}
En notant $\JC_{I^o,m}^{=h}(\oo):={\DS \lim_{\genfrac{}{}{0pt}{}{\lefto}{s}} \JC_{I^o,m}^{=h}(s)}$, un point
fermé sur $\k$ de $\JC_{I^o,m}^{=h}(\oo)$ correspond à la donnée d'un point fermé de $\IC_{I^o,m}^{=h}$ et
d'un isomorphisme de la partie connexe $(V_o^{*,c},\vphi_o^{*,c})$ de son $\OC_o$-module de Dieudonné
associé, avec $(W_{h,\k,o},\vphi_{h,\k,o})$.
\end{lemm}

\begin{defi} Pour tout $s$, en accord avec les notations de \cite{h-t}, on notera $\Tw_{I^o,h,m,t}(s)$ le quotient
$$\JC_{I^o,m}^{=h}(s) \times_{\spec \k(o)} \spf \Def_t^h/\DC_{o,h}^\times$$
où $\DC_{o,h}^\times$ agit diagonalement. L'espace topologique sous-jacent aux $\Tw_{I^o,h,m,t}(s)$ est
$\IC_{I^o,m}^{=h}$. On notera aussi $\Tw_{I^o,h,m,t}(\oo)$ la limite ${\DS
\lim_{\genfrac{}{}{0pt}{}{\lefto}{s}} \Tw_{I^o,h,m,t}(s)}$.
\end{defi}

\begin{prop} \label{comute}
La tour $\Tw_{I^o,h,m,t}(\oo)$ est munie de correspondances de Hecke associées aux \eles de
$$(D_\Am^{\oo,o})^\times \times GL_{d-h}(F_o) \times \widetilde{\NC_o},$$
vérifiant la propriété suivante:
étant donné un point géométrique $x_{I^o,m}(\oo)$ de $\JC_{I^o,m}^{=h}(\oo)$
au dessus d'un point $x_{I^o,m}$ de $\IC_{I^o,m}^{=h}$, les morphismes canoniques
$$f_{x_{I^o,m}(\oo)}:(\hat \IC_{I^o,=h,m}(t))_{x_{I^o,m}} \longto \spf \Def_t^h$$
ainsi que
$$h_{x_{I^o,m}(\oo)}:\Tw_{I^o,h,m,t}(\oo)_{x_{I^o,m}} \longto \spf \Def_t^h$$
sont tels que pour tout \ele $(g^{\oo,o},g_o^{et},(g_o^c,\d_o,\s_o)) \in (D^{\oo,o})^\times\times GL_{d-h}(F_o) 
\times
\widetilde{\NC_o}$, les diagrammes suivant sont commutatifs:\footnote{pour un choix convenables de $J^o,m',t'$, comme
ci-avant}
$$\xymatrix{
(\hat \IC_{J^o,=h,m'}(t'))_{x_{J^o,m'}} \ar^{f_{x_{J^o,m'}(\oo)}}[rrr] \ar_{(g^{\oo,o},g_o^{et},g_o^c)}[d] & & &
\spf \Def_{t'}^h \ar^{(\lexp t (g_o^c)^{-1},\d_o,\s_o)}[d] \\
(\hat \IC_{I^o,=h,m}(t))_{x_{I^o,m}} \ar^{f_{(g^{\oo,o},g_o^{et},(g_o^c,\d_o,\s_o))x_{J^o,m'}(\oo)}}[rrr] & & &
\spf \Def_t^h
}$$
$$\xymatrix{
\Tw_{J^o,h,m',t'}(\oo)_{x_{J^o,m'}} \ar^{h_{x_{J^o,m'}(\oo)}}[rrr] \ar_{(g^{\oo,o},g_o^{et},(g_o^c,\d_o,\s_o))}[d]
& & & \spf \Def_{t'}^h \ar^{(\lexp t (g_o^c)^{-1},\d_o,\s_o)}[d] \\
\Tw_{I^o,h,m,t}(\oo)_{x_{I^o,m}} \ar^{h_{(g^{\oo,o},g_o^{et},(g_o^c,\d_o,\s_o))x_{J^o,m'}(\oo)}}[rrr] & & & \spf 
\Def_t^h
}$$
où $(g^{\oo,o},g_o^{et},g_o^c,\s_o)x_{J^o,m'}=x_{I^o,m}.$
\end{prop}

\begin{proof} Pour définir les correspondances sur la tour $\Tw_{I^o,h,m,t}(\oo)$, il suffit de faire agir
$\widetilde{\NC_o}$ diagonalement\footnote{$(g_o^c,\d_o,\s_o) \in \widetilde{\NC_o}$ agit sur $\Def_t^h$ via
$(\lexp t (g_o^c)^{-1},\d_o,\s_o) \in \NC_o$ laquelle est définie par Deligne et Carayol, cf. le paragraphe
(\ref{rapel-def})}. Concrètement, étant donné
$$(g^{\oo,o},g_o^{et},g_o^c,\d_o,\s_o) \in (D_\Am^{\oo,o})^\times \times GL_{d-h}(F_o) \times \widetilde{\NC_o},$$
tel que $g_o^{et} \in \Mm_{d-h}(\OC_o)$, $g_o^c \in \Mm_h(\OC_o)$ et $\d_o \in \DC_{o,h}^\times$, soient
$J^o$, $I^o$ des idéaux de $A$ de multiplicité nulle en $o$, et $m,m',t,t',s,s'$ des entiers tels que:

\begin{itemize}

\item[-] $K^{\oo,o}_{\Am,J^o} \subset K^{\oo,o}_{\Am,I^o} \cap (g^{\oo,o})^{-1} K^{\oo,o}_{\Am,I^o} g^{\oo,o}$;

\item[-]  $\ker (\lexp t g_o^{et}) \subset (\MC_o^{-m'}/\OC_o)^{d-h}$ et $(\MC_o^{-m}/\OC_o)^{d-h} \subset  \im (\lexp 
t
g_o^{et})$;

\item[-]  $\ker (\lexp t g_o^{c}) \subset (\MC_o^{-t'}/\OC_o)^{h}$ et $(\MC_o^{-t}/\OC_o)^{h} \subset  \im (\lexp t
g_o^{c})$;

\item[-] $\varpi_o^{\val(\det(g_o^c))/h} (g_o^{et})^{-1} \in \Mm_{d-h}(\OC_o)$;

\item[-] $s'=s-\val(\det(\d_o))$.

\end{itemize}
On a alors une correspondance de Hecke $(g^{\oo,o},g_o^{et},g_o^c,\d_o,\s_o)$:
$$\JC_{J^o,m'}^{=h}(s') \times_{\spec \k(o)} \spf \Def_{t'}^h \longto \JC_{I^o,m}^{=h}(s) \times_{\spec \k(o)} \spf 
\Def_{t}^h$$
Les morphismes canoniques $f_{x_{I^o,m}(\oo)}$ et $h_{x_{I^o,m}(\oo)}$ sont définis à partir de l'isomorphisme fixé 
par le lemme
(\ref{lem-iso}) et la commutativité des diagrammes est alors évidente.

\end{proof}

\rem En particulier, on notera que l'on peut munir les
$\Tw_{I^o,h,m,t}(s)$ de correspondances de Hecke associées aux \eles
$$(g^{\oo,o},g_o^{et},g_o^c,\s_o) \in (D_\Am^{\oo,o})^\times \times GL_{d-h}(F_o) \times GL_h(F_o) \times W_{F_o},$$
à partir de celles définies ci-dessus, en choisissant un \ele quelconque $\d_o \in D_{o,h}^\times$ tel que
$(g_o^c,\d_o,\s_o) \in \widetilde{\NC_o}$, le résultat ne dépendant pas de ce choix car on a quotienté par
$\DC_{o,h}^\times$.

\begin{prop} \label{iso-local}
Au dessus de tout ouvert affine $S$ de $\IC_{I^o,m}^{=h}$, il existe des isomorphismes
$$\Tw_{I^o,h,m,t}(\oo) \times_{\IC_{I^o,m}^{=h}} S \longmapright{\sim} \hat \IC_{I^o,=h,m}(t) 
\times_{\IC_{I^o,m}^{=h}} S$$
tels qu'en tout point géométrique $x$ de $\IC_{I^o,m}^{=h}$, on a le diagramme commutatif
$$\diagram
(\Tw_{I^o,h,m,t}(\oo) \times_{\IC_{I^o,m}^{=h}} S)_x \rto^{\sim} \dto_{h_{x(\oo)}}
& (\hat \IC_{I^o,=h,m}(t) \times_{\IC_{I^o,m}^{=h}} S)_x \dto^{f_{x(\oo)}} \\
\spf \Def_t^h \rto^{=} & \spf \Def_t^h
\enddiagram$$
où $x(\oo)$ est un point géométrique de $\JC_{I^o,m}^{=h}(\oo)$ au dessus de $x$.
\end{prop}

\rem Cet énoncé peut paraître imprécis car les isomorphismes en question ne sont pas canoniques et dépendent
comme on le verra du choix d'une extension de la partie connexe par la partie étale des $\OC_o$-modules de
Dieudonné. A priori il n'y a aucun moyen de les recoller, et donc bien sur il n'y a pas de sens à essayer de
les rendre compatibles aux actions des correspondances de Hecke. La propriété essentielle est en fait la
commutativité du diagramme de l'énoncé, qui nous permettra comme on le verra au paragraphe (\ref{iso-unique})
de montrer que tous ces isomorphismes coïncident au niveau des cycles évanescents\footnote{Moralement les
cycles évanescents ne dépendent que de la partie connexe du $\OC_o$-module de Dieudonné et pas de l'extension
de cette partie par la partie étale.} et se recollent de façon compatible aux actions de Hecke, en vertu de
la proposition (\ref{comute}). On peut toutefois suspecter que le schéma tout entier est affine.

\medskip

\begin{proof} Commençons par construire la flèche. Soit donc $R$ une $\OC_o$-algèbre artinienne, $\MC$ un idéal
maximal de carré nul et $\bar R=R/\MC$. Supposons donné un $\bar R$-point de $\JC_{I^o,m}^{=h}(\oo)$. Étant
donnée $(W_{h,R,o},\vphi_{h,R,o},\iota_{h,R,t})$ une déformation de niveau $t$ de $(W_{h,\bar
R,o},\vphi_{h,\bar R,o})$, il faut construire une déformation sur $R$, du $\bar R$-point correspondant de
$\IC_{I^o,m}^{=h}$. D'après l'équivalent du théorème de Serre-Tate (cf. la proposition (\ref{prop-st})), il
suffit de construire une déformation de niveau $(m,t)$\footnote{structure de niveau $m$ (resp. $t$) sur la
partie étale (resp. connexe)} du $\OC_o$-module de Dieudonné $(\bar V_o,\bar \vphi_o)$ associé sur $\bar R$,
cette construction devant être compatible à l'action diagonale de $\DC_{o,h}^\times$ sur le produit
$\JC_{I^o,m}^{=h}(\oo) \times_{\spec \k(o)} \spf \Def_t^h$.

\begin{lemm} \label{lem-extension}
Soit $R$ une $\OC_o$-algèbre artinienne, $\MC$ un idéal maximal de carré nul et $\bar R=R/\MC$. Soit $(\bar
V_o,\bar \vphi_o)$ un $\OC_o$-module de Dieudonné sur $\bar R$, ainsi qu'une suite exacte
$$0 \longto (\bar V_o^{et},\bar \vphi_o^{et}) \longto (\bar V_o,\bar \vphi_o) \longto (\bar V_o^c,\bar \vphi_o^c) 
\longto 0$$
avec $(\bar V_o^{et},\bar \vphi_o^{et})$ étale et $(\bar V_o^c,\bar \vphi_o^c)$ connexe. Supposons en outre donnée une 
déformation
$(V_o^c,\vphi_o^c)$ sur $R$ de $(\bar V_o^c,\bar \vphi_o^c)$. Il existe alors une déformation $(V_o,\vphi_o)$ définie 
sur $R$
de $(\bar V_o,\bar \vphi_o)$ ainsi qu'une suite exacte
$$0 \longto (V_o^{et},\vphi_o^{et}) \longto (V_o,\vphi_o) \longto (V_o^c,\vphi_o^c) \longto 0$$
dont la réduction modulo $\MC$ est la suite exacte précédente sur $\bar R$.
\end{lemm}

\begin{proof} Soit $(V_o^{et},\vphi_o^{et})$ la déformation sur $R$ de $(\bar V_o^{et},\bar \vphi_o^{et})$.
Le problème est alors de construire une extension de $(V_o^{et},\vphi_o^{et})$ par $(V_o^c,\vphi_o^c)$.
La question est classique et découle des résultats de \cite{ill} chapitre 4 \S 3. L'obstruction à l'existence d'une
extension $V_o$ de $\bar V_o$ par $\MC \otimes_{\bar R} \bar V_o$ réside\footnote{Cette obstruction
est nulle et de plus il n'y a qu'une seule extension possible car en utilisant l'argument de la fin de la preuve on 
peut montrer que
$\ext^1_{\OC_o \hat \otimes_{\k(o)} \bar R}(\bar  V_o, \MC \otimes_{\bar R} \bar V_o)$ est nul. Par contre les 
flèches
$V_o^{et} \longto V_o \longto V_o^c$ ne sont pas uniquement définies.} dans
$$\ext^2_{\OC_o \hat \otimes_{\k(o)} \bar R}(\bar  V_o, \MC \otimes_{\bar R} \bar V_o)$$
Si cette obstruction est nulle, une extension étant choisie,
l'obstruction à l'existence d'une flèche $V_o^{et} \longto V_o$ (resp. $V_o \longto V_o^c$) réside dans le groupe
$$\ext^1_{\OC_o \hat \otimes_{\k(o)} \bar R}(\bar  V_o^{et}, \MC \otimes_{\bar R} \bar V_o)
\quad (\hbox{resp. } \ext^1_{\OC_o \hat \otimes_{\k(o)} \bar R}(\bar  V_o, \MC \otimes_{\bar R} \bar V_o^c))$$
et l'ensemble de ces flèches est un torseur sous
$$\hom_{\OC_o \hat \otimes_{\k(o)} \bar R}(\bar  V_o^{et}, \MC \otimes_{\bar R} \bar V_o)
\quad (\hbox{resp. }\hom_{\OC_o \hat \otimes_{\k(o)} \bar R}(\bar  V_o, \MC \otimes_{\bar R} \bar V_o^c)).$$
Si ces obstructions sont nulles,
$V_o^{et} \longto V_o$ est injective et $V_o \longto V_o^c$ est surjective. En outre l'ensemble des flèches $V_o^{et} 
\longto V_o^c$
est un torseur sous
$$\hom_{\OC_o \hat \otimes_{\k(o)} \bar R}(\bar  V_o^{et}, \MC \otimes_{\bar R} \bar V_o^c),$$
de sorte qu'il existe des flèches $V_o^{et} \longto V_o$ et $V_o \longto V_o^c$ tel que la suite
$0 \longto V_o^{et} \longto V_o \longto V_o^c \longto 0$ soit exacte. Montrons que toutes ces obstructions
sont nulles. En fait on va montrer que les groupes dans lesquelles elles vivent sont nuls. L'argument est identique 
pour tous
ces groupes, traitons par exemple le cas de $\ext^1_{\OC_o \hat \otimes_{\k(o)} \bar R}(\bar  V_o^{et}, \MC 
\otimes_{\bar R} \bar V_o)$.
Les $\OC_o \hat \otimes_{\k(o)} \bar R$-modules $\bar V_o$ et $\bar V_o^{et}$ étant localement libres, les faisceaux
$\underline{Ext}^i_{\OC_o \hat \otimes_{\k(o)} \bar R}(\bar  V_o^{et}, \MC \otimes_{\bar R} \bar V_o)$ sont nuls pour 
$i >0$.
La suite spectrale locale-globale pour le $\ext^1$ montre que
$\ext^1_{\OC_o \hat \otimes_{\k(o)} \bar R}(\bar  V_o^{et}, \MC \otimes_{\bar R} \bar V_o)$
est égal à
$$H^1(\spec (\bar R),\underline{Ext}^0_{\OC_o \hat \otimes_{\k(o)} \bar R}(\bar  V_o^{et}, \MC \otimes_{\bar R} \bar 
V_o)$$
qui est nul car $\spec (\bar R)$ est affine.

Il faut alors construire $\vphi_o: \lexp \t V_o \longto V_o$ qui soit compatible à $\vphi_o^{et}$ et $\vphi_o^c$.
On commence par remarquer que $\frob_o:R \longto R$ se factorise par $\bar \frob_o:\bar R \longto R$, de sorte
$\lexp \t V_o=\bar V_o \otimes_{\bar R,\bar \frob_o} R$. Comme précédemment il n'y a pas d'obstruction à l'existence
d'une application $\vphi'_o:\lexp \t V_o \longto V_o$ et l'ensemble des telles applications est un torseur sous
$\hom_{\OC_o \hat \otimes_{\k(o)} \bar R}(\lexp \t {\bar V_o}, \MC \otimes_{\bar R} \bar V_o)$. Fixons une telle 
application
$\vphi'_o$. Les applications $V_o^{et} \longto V_o$, $\vphi'_o$ et $\vphi_o^{et}$ fournissent alors un \ele
de $\hom_{\OC_o \hat \otimes_{\k(o)} \bar R}(\lexp \t {\bar V_o^{et}},\MC \otimes_{\bar R} \bar V_o)$. En appliquant 
le foncteur
$\hom_{\OC_o \hat \otimes_{\k(o)} \bar R}(\bullet,\MC \otimes_{\bar R} \bar V_o)$ à la suite exacte courte
$$0 \to \lexp \t {\bar V_o^{et}} \longto \lexp \t {\bar V_o} \longto \lexp \t {\bar V_o^c} \to 0,$$
on obtient la suite exacte
\begin{multline*}
0 \to \hom_{\OC_o \hat \otimes_{\k(o)} \bar R}(\lexp \t {\bar V_o^c},\MC \otimes_{\bar R} \bar V_o) \longto
\hom_{\OC_o \hat \otimes_{\k(o)} \bar R}(\lexp \t {\bar V_o},\MC \otimes_{\bar R} \bar V_o) \to \\
\longto \hom_{\OC_o \hat \otimes_{\k(o)} \bar R}(\lexp \t {\bar V_o^{et}},\MC \otimes_{\bar R} \bar V_o) \longto 0
\end{multline*}
car, de manière identique à ce qui précède,
$\ext^1_{\OC_o \hat \otimes_{\k(o)} \bar R}(\lexp \t {\bar V_o^c},\MC \otimes_{\bar R} \bar V_o)$ est nul.
On modifie donc $\vphi'_o$ en un $\vphi''_o$ qui soit compatible à $\vphi_o^{et}$; l'ensemble des telles $\vphi''_o$ 
est alors
un torseur sous $\hom_{\OC_o \hat \otimes_{\k(o)} \bar R}(\lexp \t {\bar V_o^c},\MC \otimes_{\bar R} \bar V_o)$.
De même, $V_o \longto V_o^c$, $\vphi''_o$ et $\vphi_o^c$ définissent un \ele de
$\hom_{\OC_o \hat \otimes_{\k(o)} \bar R}(\lexp \t {\bar V_o^c},\MC \otimes_{\bar R} \bar V_o^c)$. En appliquant le 
foncteur
$\hom_{\OC_o \hat \otimes_{\k(o)} \bar R}(\lexp \t {\bar V_o^c},\bullet)$, à la suite exacte
$$0 \to \MC \otimes_{\bar R} \bar V_o^{et} \longto \MC \otimes_{\bar R} \bar V_o \longto \MC \otimes_{\bar R} \bar 
V_o^c \to 0,$$
on obtient la suite exacte
\begin{multline*}
0 \to \hom_{\OC_o \hat \otimes_{\k(o)} \bar R}(\lexp \t {\bar V_o^c},\MC \otimes_{\bar R} \bar V_o^{et}) \longto
\hom_{\OC_o \hat \otimes_{\k(o)} \bar R}(\lexp \t {\bar V_o^c},\MC \otimes_{\bar R} \bar V_o) \to \\
\to \hom_{\OC_o \hat \otimes_{\k(o)} \bar R}(\lexp \t {\bar V_o^c},\MC \otimes_{\bar R} \bar V_o^c) \longto 0
\end{multline*}
car $\ext^1_{\OC_o \hat \otimes_{\k(o)} \bar R} (\lexp \t {\bar V_o^c}, \MC \otimes_{\bar R} \bar V_o^{et})$ est nul.
Il est alors possible de modifier $\vphi''_o$ en un $\vphi_o$ compatible à $\vphi_o^{et}$ et $\vphi_o^c$.

\end{proof}

Soit donc comme dans le lemme ci-dessus, $(V_o^{et},\vphi_o^{et})$ le relèvement sur $R$ de $(\bar V_o^{et},\bar 
\vphi_o^{et})$.
Soit aussi $\iota_{o,m}^{et}$ le relèvement sur $R$ de $\bar \iota_{o,m}^{et}$:
c'est une structure de niveau $m$ sur $(V_o^{et},\vphi_o^{et})$.
De l'isomorphisme sur $\bar R$:
$$(\bar V_o^{c*},\bar \vphi_o^{c*}) \simeq (W_{h,\bar R,o},\vphi_{h,\bar R,o}),$$
et de la déformation $(W_{h,R,o},\vphi_{h,R,o})$, on en déduit une déformation $(V_o^c,\vphi_o^c)$ définie
sur $R$, de $(\bar V_o^c,\bar \vphi_o^c)$, munie de plus d'une structure de niveau $t$. Il suffit alors
d'appliquer le lemme précédent. Il est en outre immédiat que cette construction est invariante par l'action
diagonale de $\DC_{o,h}^\times$ sur le produit $\JC_{I^o,m}^{=h}(\oo) \times_{\spec \k(o)} \spf \Def_t^h$ et
que le diagramme de l'énoncé est bien commutatif.

\medskip

Montrons ensuite que la flèche en question est un isomorphisme.
Le pro\-blème se ramène immédiatement au cas $t=0$. Il suffit alors de montrer que l'on obtient un isomorphisme au 
niveau
des espaces tangents. Pour des raisons de dimension, il suffit de vérifier que l'on a une injection ce qui est 
clairement le cas.

\end{proof}

\marque Dans la suite nous proposons suivant \cite{h-t}, une autre façon de voir cet isomorphisme en un cran
fini. Notons
$$\hat \JC_{I^o,=h,m}(s) \longto \hat \IC_{I^o,=h,m}$$
l'extension étale de fibre spéciale $\JC_{I^o,m}^{=h}(s) \longto \IC_{I^o,m}^{=h}$. On note de même
$$\hat \JC_{I^o,=h,m}(s,t) \longto \hat \JC_{I^o,=h,m}(s)$$
le classifiant des structures de niveau $t$ sur la partie
connexe du $\vphi$-faisceau universel sur $\hat \JC_{I^o,=h,m}(s)$. On définit de manière identique au
paragraphe (\ref{cores-j}), des correspondances de Hecke sur $\hat \JC_{I^o,=h,m}(s,t)$ associées aux \eles
de $(D_\Am^{\oo,o})^\times \times GL_{d-h}(F_o) \times \widetilde{\NC_o}$.

On note $\MC_t^h$ l'idéal maximal de $\Def_t^h$ et soient pour $S$ un ouvert affine de $\IC_{I^o,m}^{=h}$:
$$S_s:= \JC_{I^o,m}^{=h}(s) \times_{\IC_{I^o,m}^{=h}} S, \qquad \hat S_{s,t}:= \hat \JC_{I^o,=h,m}(s,t) 
\times_{\IC_{I^o,m}^{=h}} S.$$
Étant donné un point géométrique $x(\oo)$ de $\JC_{I^o,m}^{=h}(\oo)$, on a des morphismes canoniques
$$(S_s \times_{\spec \k(o)} \spf (\Def_t^h/(\MC_t^h)^N))_{x(\oo)} \longto \spf (\Def_t^h/(\MC_t^h)^N)$$
et
$$(\hat S_{s,t})_{x(\oo)} \longto \spf \Def_t^h$$
ainsi qu'une compatibilité aux correspondances de Hecke comme dans la proposition (\ref{comute}).

\begin{prop} \label{niv-fini}
Pour tout $N \in \Nm$, il existe un $s_0$ assez grand tel que pour tout $s \geq s_0$, et $S$ un ouvert affine
de $\IC_{I^o,m}^{=h}$, on ait un morphisme
$$f_{N,s}: S_s \times_{\spec \k(o)} \spf \Def_t^h \longto \hat S_{s,t}$$
tel qu'en tout point géométrique $x(\oo)$ de $\JC_{I^o,m}^{=h}(\oo)$, on ait
le diagramme suivant:
$$\xymatrix{
(S_s \times_{\spec \k(o)} \spf (\Def_t^h/(\MC_t^h)^N))_{x(\oo)} \ar^{f_{N,s,t}}[rrr] \dto & & & (\hat 
S_{s,t})_{x(\oo)} \ar[d] \\
\spf (\Def_t^h/(\MC_t^h)^N) \ar@{^{(}->}[rrr]  & & & \spf \Def_t^h
}$$
\end{prop}

\begin{proof} On commence par prouver le lemme suivant qui est l'équivalent pour les $\OC_o$-modules de Dieudonné du 
lemme (1.1.3)
de \cite{kat}.

\begin{lemm}
Soit $R$ une $\OC_o$-algèbre artinienne de morphisme structural $i:\OC_o \to R$, et $r \in \Nm$ tel que
$i(\varpi_o^r)=0$ dans $R$. Soient $(V_o,\vphi_o)$ et $(W_o,\psi_o)$ deux $\OC_o$-modules de Dieudonné
connexes sur $R$. Soit $\MC$ l'idéal maximal de $R$ et $\bar R=R/\MC$: on suppose que $\MC^k=(0)$. On notera
$\bar V_o$ et $\bar W_o$ la réduction modulo $\MC$ de $V_o$ et $W_o$. On a alors

\begin{itemize}
\item[-] $\hom((V_o,\vphi_o),(W_o,\psi_o)) \longto \hom((\bar V_o,\bar \vphi_o),(\bar W_o,\bar \psi_o))$ est 
injective;

\item[-] pour toute application $\bar f:(\bar V_o,\bar \vphi_o) \longto (\bar W_o,\bar \psi_o)$, on peut trouver un
relèvement de $\varpi_o^{N} \bar f$, où $N \geq r (\ln k / \ln q)$.
\end{itemize}
\end{lemm}

\begin{proof} Commençons par montrer que $\varpi_o^N \otimes 1$ annule $\MC W_o$ pour $N \geq r (\ln k /\ln q)$. Dans 
$W_o$,
en considérant $\psi_o:W_o \longto W_o$ comme une application $\t$-linéaire, on a
$(\varpi_o^r \otimes 1) \Id =(1 \otimes i(\varpi_o^r)) \Id + \sum_{i=1}^{rh} a_i (\psi_o)^i$. Comme $i(\varpi_o^r)$ 
est nul dans $R$, on en
déduit que pour $v \in \MC W_o$, $(\varpi_o^r \otimes 1) v \in \MC^p W_o$, d'où le résultat.

Soit donc $g:(V_o,\vphi_o) \longto \MC (W_o,\psi_o)$: $(\varpi_o^N \otimes 1)v \in \ker g$. Or $(V_o,\vphi_o)$ étant 
connexe,
on en déduit qu'il existe $w \in V_o/ \vphi_o(V_o)$ tel que $(\varpi_o^N \otimes 1)v=(\vphi_o)^i(w)$. On a alors
$(\psi_o)^i(g(w))=0$ soit $g(w)=0$ et donc $g=0$.

Pour le dernier point, on définit l'image d'un \ele $v \in V_o$ comme $\varpi_o^N \otimes 1$ fois
le relèvement quelconque de $\bar f (\bar v)$. Cette définition ne dépend pas du choix de ce relèvement car $\MC W_o$ 
est tué
par $\varpi_o^N \otimes 1$.

\end{proof}

\begin{lemm} \label{lem-alpha}
(cf. \cite{h-t} lemme III.2.6)
Il existe une fonction $\a: \Nm \longto \Nm$ croissante de limite l'infini, telle que:

\begin{itemize}

\item[-] tout \ele de $\varpi_o^{s-t} \DC_{o,h}$ se relève en un endomorphisme de
$$(W_{h,\Def_t^h/(\MC_t^h)^{\a(s)},o},\vphi_{h,\Def_t^h/(\MC_t^h)^{\a(s)}})$$

\item[-] tout \ele de $(1+\varpi_o^s \DC_{o,h})$ agit trivialement sur $\Def_t^h/(\MC_t^h)^{\a(s)}$.

\end{itemize}
\end{lemm}

\begin{proof} On remarque tout d'abord que pour $t >0$, $\varpi_o \in (\MC_t^h)^{(q-1)q^{(t-1)h}}$. Le premier point 
découle alors du lemme
précédent en prenant $\a(s)$ tel que
$$\a(s) \ln (\a(s)) \leq (s-t) (\ln q) (q-1)q^{(t-1)h}.$$
Le deuxième point est alors immédiat: en effet soit $\d \in \DC_{o,h}^\times$ tel que $\d \equiv 1 \mod \varpi_o^s$. 
Alors
$(\d-1)/\varpi_o^m$ se relève en un endomorphisme de
$$(W_{h,\Def_t^h/(\MC_t^h)^{\a(s)},o},\vphi_{h,\Def_t^h/(\MC_t^h)^{\a(s)},o}).$$
Il en est de même de $(\d^{-1}-1)/\varpi_o^m$, de sorte que $\d$ se relève en un automorphisme de
$$(W_{h,\Def_t^h/(\MC_t^h)^{\a(s)},o},\vphi_{h,\Def_t^h/(\MC_t^h)^{\a(s)},o})$$
qui est l'identité sur les points de $\varpi_o^m$-torsion.

\end{proof}

Soit alors $S_{s,t}(N):=S_s \times_{\spec \k(o)} \spf (\Def_t^h/(\MC_t^h)^{N})$. D'après le lemme ci-dessus,
on a une action de $(1+\varpi_o^s \DC_{o,h})$ sur
$(W_{h,\Def_t^h/(\MC_t^h)^{\a(s)},o},\vphi_{h,\Def_t^h/(\MC_t^h)^{\a(s)},o})$ telle que pour tout $a \geq t$,
$(1+\varpi_o^{s-t+a} \DC_{o,h})$ agit trivialement sur
$$(W_{h,\Def_t^h/(\MC_t^h)^{\a(s)},o,a},\vphi_{h,\Def_t^h/(\MC_t^h)^{\a(s)},o,a}).$$
Ainsi $(1+\varpi_o^s \DC_{o,h})/(1+\varpi_o^{s+a-t} \DC_{o,h})$ agit diagonalement sur
$$(W_{h,\Def_t^h/(\MC_t^h)^{\a(s)},o,a},\vphi_{h,\Def_t^h/(\MC_t^h)^{\a(s)},o,a}) \times S_{s+a-t,t}(\a(s))$$
et le quotient est
un $\OC_o$-module de Dieudonné tronqué à l'ordre $a$
$$(V_{o,a}^{c,s},\vphi_{o,a}^{c,s})$$
connexe sur $S_{s,t}(\a(s))$. La limite projective
des $(V_{o,a}^{c,s},\vphi_{o,a}^{c,s})$ est alors un $\OC_o$-module de Dieudonné $(V_{o}^{c,s},\vphi_{o}^{c,s})$ sur 
$S_{s,t}(\a(s))$
tel que
$$(V_{o,t}^{c,s},\vphi_{o,t}^{c,s}) \simeq (W_{h,\Def_t^h,o,t},\vphi_{h,\Def_t^h,o,t})$$
de sorte que $(V_{o}^{c,s},\vphi_{o}^{c,s})$ est muni d'une structure de niveau $t$. Étant donné un $\bar
R$-point de $\JC_{I^o,m}^{=h}(s)$, une $\OC_o$-algèbre artinienne $R_s$, $\MC$ un idéal maximal de $R_s$ tel
que $\bar R=R_s/\MC$, ainsi qu'un $R_s$-point de $\spf (\Def_t^h/(\MC_t^h)^{\a(s)})$, on construit une
déformation sur $R_s$ de ce $\bar R$-point, en donnant une déformation de niveau $(m,t)$ de son
$\OC_o$-module de Dieudonné. On procède comme dans la preuve de la proposition précédente. Soit
$(V_{o}^{et,s},\vphi_{o}^{et,s},\iota_{o,m}^{et,s})$ le relèvement sur $R$ de $(\bar V_o^{et},\bar
\vphi_o^{et},\bar \iota_{o,m}^{et})$. D'après le lemme (\ref{lem-extension}), il existe une extension
$(V_o^s,\vphi_o^s)$ de $(V_o^{et,s},\vphi_o^{et,s})$ par $(V_o^{c,s},\vphi_o^{c,s})/R$ ainsi qu'une suite
exacte
$$0 \longto (V_o^{et,s},\vphi_o^{et,s}) \longto (V_o^s,\vphi_o^s) \longto (V_o^{c,s},\vphi_o^{c,s}) \longto 0$$
sur $R_s$ qui se restreigne sur $R_s \otimes_{\Def_t^h/(\MC_t^h)^{\a(s)}} \Def_t^h/(\MC_t^h)^{\a(s-1)}$ au pull-back 
de
$$0 \longto (V_o^{et,s-1},\vphi_o^{et,s-1}) \longto (V_o^{s-1},\vphi_o^{s-1}) \longto (V_o^{c,s-1},\vphi_o^{c,s-1}) 
\longto 0$$
On en déduit donc l'existence d'un\footnote{non unique} morphisme pour tout ouvert affine $S$ de
$\IC_{I^o,m}^{=h}$,
$$f_{s,t}: S_{s,t}(\a(s)) \longto \hat S_{s,t}$$
tel que le diagramme ci-dessous est commutatif
$$\diagram
S_{s,t}(\a(s)) \rto^{f_{s,t}} & \hat S_{s,t} \\
S_{s+1,t}(\a(s)) \uto \rto^{f_{s+1,t}} & \hat S_{s+1,t} \uto
\enddiagram$$
En particulier si $x(\oo)$ est un point fermé de $\JC_{I^o,m}^{=h}(\oo)$,
on a le diagramme commutatif
$$\xymatrix{
(S_{s,t}(\a(s)))_{x(s)} \ar^{f_{s,t}}[r] \ar[d] & (\hat S_{s,t})_{x(s)} \ar^{j_{x(\oo)}}[d] \\
\spf (\Def_t^h/(\MC_t^h)^{\a(s)}) \ar@{^{(}->}[r] & \spf \Def_t^h }$$ L'existence du morphisme $f_{N,s,t}$ de
l'énoncé découle alors d'un lemme de Berkovich (cf. lemme 1 de l'annexe de \cite{h-t}, ou bien le lemme
II.5.8 de loc. cit.)

\end{proof}


%% file: cycles.tex
\section{Cycles évanescents et systèmes locaux d'Harris-Taylor}

On rappelle que le but des variétés d'Igusa est en particulier de décrire la restriction à $M_{I,s_o}^{=h}$
du faisceau des cycles évanescents de $M_{I,o} \longto \spec (\OC_o)$. Cette question se pose dans le
contexte suivant : la cohomologie de la fibre générique de la tour des $(M_{I,o})_I$ est décrite dans
\cite{lrs} de manière ``assez'' précise. Cette cohomologie est aussi l'aboutissement de la suite spectrale
des cycles évanescents
$$E_2^{p,q}(I)=H^p(M_{I,\bar s_o},R^q \Psi_{\eta_o}(\bar \Qm_l) \Rightarrow H^{p+q}(M_{I,\eta_o},\bar \Qm_l)$$
Afin de séparer la contribution due aux différentes strates, entre en jeu la suite spectrale associée à la 
stratification
$$E_1^{p,q}(I,i)=H_c^{p+q}(M_{I,\bar s_o}^{=d-p},R^i \Psi_{\eta_o}(\bar \Qm_l)) \Rightarrow
H^{p+q}(M_{I,\bar s_o},R^i \Psi_{\eta_o}(\bar \Qm_l)$$ On cherche donc à avoir des renseignements sur ${\DS
\lim_{\genfrac..{0pt}{1}{\lefto}{I}}} H^i_c(M_{I,\bar s_o,1}^{=h},R^j \Psi_{\eta_o}(\bar \Qm_l))$ en tant que
représentation de $(D_\Am^{\oo,o})^\times \times P_{h,d}(F_o) \times W_{F_o}$.

\subsection{Rappels sur les cycles évanescents associés à $\Def_n^h$}

\label{rapel-cycle}

Pour tout $n \geq 0$, $\Def_n^d$ représente le foncteur des déformations de niveau $n$, par isogénies, du
$\OC_o$-module formel de hauteur $d$ sur $\bar \Fm_q$, cf. \cite{dr1}. Soit alors $\Psi_{F_o,n}^{d,i}$ le
$\bar \Qm_l$-espace vectoriel de dimension finie obtenu via la théorie de Berkovich comme le $i$-ème foncteur
des cycles évanescents associé au morphisme structural
$$\spf \Def_n^d \longto \spf \hat \OC_o^{nr}.$$
Cet espace vectoriel est muni entre autre d'une action de $GL_d(\OC_o)$ qui se factorise par le morphisme
surjectif naturel $GL_d(\OC_o) \longto GL_d(\OC_o/\MC_o^n)$ et on pose
$$\Psi_{F_o}^{d,i}= {\DS \lim_{\genfrac{}{}{0pt}{}{\longto}{m \geq n}}} \Psi_{F_o,m}^{d,i}$$
de sorte que pour $K_{o,n}:=\ker (\OC_o^\times \longto (\OC_o/\MC_o^n)^\times)$,
$\Psi_{F_o,n}^{d,i}=(\Psi_{F_o}^{d,i})^{K_{o,n}}$. On introduit le groupe $\NC_o$ (resp. $\NC_o'$) défini
comme le noyau de
$$(g_o,\d_o,c_o) \in GL_d(F_o) \times D_{o,d}^\times \times W_o \mapsto \val(\det(g_o^{-1})\rn(\d_o)\cl(c_o))
\in \Zm$$ (resp. composé avec la projection canonique $\Zm \longto \Zm/d\Zm$. Comme rappelé au paragraphe
(\ref{rapel-def}), pour $\xi_o$ un caractère d'ordre fini de $F_o^\times$, $\Psi_{F_o}^{d,i}(\xi_o')$ (resp.
$\Psi_{F_o}^{d,i}$), où $\xi_o'$ est la restriction de $\xi_o$ à $\OC_o^\times$, est muni d'une action de
$\NC_o'$ (resp. de $\NC_o$).

\marque Dans la définition de $\Def_n^d$, il est agréable de considérer plutôt les déformations par
quasi-isogénies de sorte que la construction précédente fourni des $\bar \Qm_l$-espaces vectoriels
$\UC_{F_o,n}^{d,i} \simeq (\UC_{F_o}^{d,i})^{K_{o,n}}$ où
$$\UC_{F_o,\xi_o}^{d,i}:=\ind_{\NC_o'}^{GL_d(F_o) \times D_{o,d}^\times \times W_{F_o}}\Psi_{F_o,\xi_o}^{d,i}$$
est une représentation de $GL_d(F_o) \times D_{o,d}^\times \times W_{F_o}$.

Pour toute représentation admissible irréductible $\t_o$ de $D_{o,d}^\times$, la réciprocité de Frobenius
donne un isomorphisme
$$\UC_{F_o}^{d,i}(\t_o) \simeq \hom_{\DC_{o,d}^\times}(\res^{D_{o,d}^\times}_{\DC_{o,d}^\times} 
\t_o,\Psi_{F_o,\xi_o}^{d,i})$$
où $\xi_o$ est le caractère central de $\t_o$ et l'action de $(g_o^c,\s_o)$ est donnée par celle de
$(g_o,\d_o,\s_o) \in \NC_o'$ pour $\d_o \in D_{o,h}^\times$ quelconque. Afin d'avoir un théorème de
Serre-Tate équivariant on introduit les modifications suivantes.

\begin{defi} \label{action-tilde}
Pour tout $h,i$, $\widetilde{\UC_{F_o}^{h,i}}$ est la représentation $\UC_{F_o}^{h,i}$ où l'action de
$GL_h(F_o)$ est tordue par $g_o \mapsto \lexp t g_o^{-1}$. On définit de même $\widetilde{\Psi_{F_o}^{h,i}}$
comme étant l'espace $\Psi_{F_o}^{h,i}$ muni d'une action de
$$\widetilde{\NC_o}:= \ker( (g_o,\d_o,c_o) \in GL_d(F_o) \times D_{o,d}^\times \times W_o \mapsto
\val(\det(g_o)\rn(\d_o)\cl(c_o)) \in \Zm$$ via l'action de $(\lexp t g_o^{-1},\d_o,c_o) \in \NC_o$ sur
$\Psi_{F_o}^{h,i}$.
\end{defi}

\subsection{Restriction aux strates ouvertes et variétés d'Igusa}

\label{iso-unique}

On rappelle le lemme suivant bien connu.

\begin{lemm} \label{lem-lisse}
Soit $\LC$ un $\bar \Qm_l$-faisceau lisse sur $M_{I,o}$. On a alors
$$R^i \Psi_{\eta_o}(\LC) \simeq (R^i \Psi_{\eta_o}(\bar \Qm_l)) \otimes \LC_{s_o}$$
où $\LC_{s_o}$ est la restriction de $\LC$ à la fibre spéciale $M_{I,s_o}$.
\end{lemm}

\marque Appliqué au faisceau $\LC_{\r_\oo}$ cela nous conduit à étudier $R^i \Psi_{\eta_o}(\bar
\Qm_l)_{M_{I,o}}$ le $i$-ème faisceau des cycles évanescents du morphisme structural
$$M_{I,o} \longto \spec \OC_o.$$
On notera aussi $R^i \Psi_{\eta_o}(\bar \Qm_l)_{M_{I,s_o,a}^{=h}}$ sa restriction à $M_{I,s_o,a}^{=h}$. On
considère le tiré en arrière $R^i \Psi_{\eta_o}(\bar \Qm_l)_{\IC_{I^o,n}^{=h}}$ de  $R^i \Psi_{\eta_o}(\bar
\Qm_l)_{M_{I,s_o,a}^{=h}}$ par le morphisme $\IC_{I^o,n}^{=h} \longto M_{I,s_o,a}^{=h}$, où $n$ est la
multiplicité de $o$ dans $I$. On rappelle qu'étant donné un point géométrique $x(\oo)$ de
$\JC_{I^o,m}^{=h}(\oo)$ au dessus d'un point $x$ de $M_{I,s_o}^{=h}$ (resp. $y$ de $\IC_{I^o,n}^{=h}$), on a
des isomorphisme canoniques
$$\widetilde{\Psi_{F_o,n}^{h,i}} \simeq (R^i \Psi_{\eta_o}(\bar \Qm_l))_{(\widehat{M_{I,o,=h,a}})_x} \simeq
(R^i \Psi_{\eta_o}(\bar \Qm_l)_{M_{I,s_o,a}^{=h}})_x $$
$$\widetilde{\Psi_{F_o,n}^{h,i}} \simeq (R^i \Psi_{\eta_o}(\bar \Qm_l))_{(\hat \IC_{I^o,=h,n}(n))_y} \simeq
(R^i \Psi_{\eta_o}(\bar \Qm_l)_{\IC_{I^o,n}^{=h}})_y $$
que l'on note dans les deux cas $j_{x(\oo)}$.
Le but de ce paragraphe et finalement des variétés d'Igusa est de décrire le faisceau $R^i \Psi_{\eta_o}(\bar 
\Qm_l)_{M_{I,s_o,a}^{=h}}$.
Selon loc. cit. nous introduisons les notations suivantes:

- soient $\bar Y$ un schéma lisse sur $\bar \Fm_q$, et $(\bar Y_n \to \bar Y)$ un système projectif de
revêtements étales de groupe de Galois $G={\DS \lim_{\genfrac{}{}{0pt}{}{\lefto}{n}}} G_n$;

- soit $\LF$ un $\bar \Qm_l$-faisceau lisse sur $\bar Y$, muni d'une action de
$G$ se factorisant par un quotient fini $G_N$ de $G$.

Le faisceau $\LF_{|\bar Y_n}$ est ainsi muni d'une action de $G_n \times G_n$ et l'action diagonale de $G_n$
est compatible à son action sur $\bar Y_n$ de sorte qu'il existe un faisceau
$$\twist_{(\bar Y_n)} (\LF)$$
sur $\bar Y$ tel que sa restriction à $\bar Y_n$ munie de l'action naturelle
de $G_n$ coïncide avec $\LF_{|\bar Y_n}$ muni de l'action diagonale de $G_n$.

\marque Rappelons l'énoncé suivant dû à Berkovich (cf. loc. cit. appendice). Supposons donné sur $\hat
\OC_o^{nr}$, un schéma formel $\XC$ muni d'une action d'un groupe $G$ agissant trivialement sur la fibre
spéciale $\XC_s$. On note $\XC(n)$ le schéma formel $(\XC_s,\OC_{\XC}/\IC_{\XC}^n)$ où $\IC_{\XC}$ est un
idéal de définition de $\XC$. On suppose que $G={\DS \lim_{\genfrac{}{}{0pt}{}{\lefto}{n}} G(n)}$ tel que
l'action de $G$ sur $\XC(n)$ se factorise par $G(n)$. Soit alors $\XC_n \longto \XC$ des revêtements
galoisiens de groupe de Galois $G(n)$ et on fait les hypothèses suivantes:

\begin{itemize}

\item[-] le quotient de $\XC_n(n)$ par l'action diagonale de $G(n)$ est un schéma formel $\YC^n$ (c'est le cas si par 
exemple
$\XC$ est quasi-projectif sur $\spf \hat \OC_o^{nr}$);

\item[-] il existe un schéma formel spécial $\YC$ tel que $\YC(n)=\YC^n$;

\item[-] la famille des complétions formelles de $\XC$ le long d'un point fermé de $\XC_s$, a un nombre fini de 
classes d'isomorphismes.

\end{itemize}

\marque On note $\Psi_m^i(\XC)$ (resp. $\Psi_m^i(\YC)$) le $i$-ème faisceau des cycles évanescents du
faisceau constant $\Zm/p^m\Zm$ associé à $\XC$ (resp. $\YC$). D'après les travaux de Berkovich (cf.
\cite{berk}), il existe un entier $N$, tel que tout automorphisme sur $\XC$ trivial sur $\XC(N)$, agit
trivialement sur $\Psi_m^i(\XC)$. On fixe donc pour tout $m$, un tel entier $N(m)$, de sorte que l'action de
$G$ sur $\Psi_m^i(\XC)$ se factorise par $G(N(m))$. En particulier les faisceaux $\twist_{\XC_n}
(\Psi_m^i(\XC))$ sont indépendants de $n \geq N(m)$; on le note $\twist_{\XC_\oo} (\Psi_m^i(\XC))$. Berkovich
démontre alors le résultat suivant.

\begin{theo} \label{theo-berk} (cf. l'appendice de Berkovich dans \cite{h-t})
Il existe un système compatible d'isomorphismes canoniques de faisceaux:
$$\Psi_m^i(\YC) \simeq \twist_{\XC_\oo} (\Psi_m^i(\XC))$$
\end{theo}

\marque On cherche à appliquer le théorème ci-dessus à $\YC=\Tw_{I^o,h,m,t}(\oo)$. On pose donc $\XC=
\IC_{I^o,m}^{=h} \times_{\spec \k(o)} \spf \Def_t^h$ qui est muni d'une action de $G=\DC_{o,h}^\times$
agissant trivialement sur $\XC_s=\IC_{I^o,m}^{=h}$. On note $\MC_t^h$ l'idéal maximal de $\Def_t^h$ et soit
$\XC(n):=\IC_{I^o,m}^{=h} \times_{\spec \k(o)} \spf (\Def_t^h/(\MC_t^h)^n)$. On considère une fonction $N:\Nm
\longto \Nm$, telle que $\a \circ N=\Id$, où $\a$ est la fonction définie au lemme (\ref{lem-alpha}), de
sorte que l'action de $G$ sur $\XC(n)$ se factorise à travers $G(n):= \DC_{o,h,N(n)}$. On considère ensuite
$\XC_s:=\JC_{I^o,m}^{=h}(N(s)) \times_{\spec \k(o)} \spf \Def_t^h$ de sorte  que $\XC_s \longto \XC$ est
galoisien de groupe de Galois $G(s)=\DC_{o,h,N(s)}$. On est ainsi dans les conditions d'application du
théorème ci-dessus d'où le résultat suivant.

\begin{theo} \label{theo-principal}
Il existe, sur $\IC_{I^o,m}^{=h} \otimes_{\k(o)} \bar \k(o)$, un isomorphisme canonique de $\bar
\Qm_l$-faisceaux
$$R^i \Psi_{\eta_o}(\bar \Qm_l)_{\Tw_{I^o,h,m,t}(\oo)} \longmapright{\sim} \twist_{\JC_{I^o,m}^{=h}(\oo)}
(\widetilde{\Psi_{F_o,t}^{h,i}})$$ vérifiant les propriétés suivantes

\begin{itemize}

\item[(1)] si $x$ est un point géométrique de $\IC_{I^o,m}^{=h}$ et si
$x(\oo)=(x(s))_s$ est un système projectif de points de $\JC_{I^o,m}^{=h}(s)$ au dessus d'un point $x$ de 
$\IC_{I^o,m}^{=h}$,
on a un isomorphisme
$$(\Tw_{I^o,h,m,t}(\oo))_x \longto (\IC_{I^o,m}^{=h} \times_{\spec \k(o)} \spf \Def_t^h)_x .$$
Par le théorème de changement de base lisse le faisceaux des cycles évanescents sur $\IC_{I^o,m}^{=h}
\times_{\spec \k(o)} \spf \Def_t^h$ est le faisceau constant $\widetilde{\Psi_{F_o,t}^{h,i}}$. On en déduit
donc un isomorphisme
$$f_{x(\oo)}: \widetilde{\Psi_{F_o,t}^{h,i}} \longto R^i \Psi_{\eta_o}(\bar \Qm_l)_{\Tw_{h,m,t}(\oo)_x} .$$
De même la donnée de $x(\oo)$, donne un isomorphisme
$$h_{x(\oo)}: \widetilde{\Psi_{F_o,t}^{h,i}} \longto (\twist_{J_{I^o,m}^{=h}(\oo)} \widetilde{\Psi_{F_o,t}^{h,i}})_x 
$$
et l'on a le diagramme commutatif
$$\diagram (R^i \Psi_{\eta_o}(\bar \Qm_l)_{\Tw_{I^o,h,m,t}(\oo)})_{x} \rrto^{\sim} \dto &
&  (\twist_{\JC_{I^o,m}^{=h}(\oo)} (\widetilde{\Psi_{F_o,t}^{h,i}}))_{x}  \\
R^i \Psi_{\eta_o}(\bar \Qm_l)_{\Tw_{I^o,h,m,t}(\oo)_x} & & \widetilde{\Psi_{F_o,t}^{h,i}} \uto^{h_{x(\oo)}}
\llto^{f_{x(\oo)}}
\enddiagram$$

\item[(2)] l'action par correspondances d'un \ele
$$(g^{\oo,o},g_o^{et},g_o^c) \in (D_\Am^{\oo,o})^\times \times GL_{d-h}(F_o) \times GL_h(F_o) $$
sur $R^i \Psi_{\eta_o}(\bar \Qm_l)_{\Tw_{I^o,h,m,t}(\oo)}$ donne par l'isomorphisme de l'énoncé, une action
par correspondances sur $\twist_{\JC_{I^o,m}^{=h}(\oo)} (\widetilde{\Psi_{F_o,t}^{h,i}})$ qui découle de
l'action naturelle par correspondances d'un \ele
$$(g^{\oo,o},g_o^{et},(g_o^c,\s_o,\d_o)) \in (D^{\oo,o})^\times/K_{\Am,I}^{\oo,o} \times GL_{d-h}(F_o) \times 
\widetilde{\NC_o}$$
\footnote{le résultat ne dépend pas du choix de $\d_o$} sur le produit $\JC_{I^o,m}^{=h}(s) \times \spf \Def_t^h$ où
l'on rappelle que $(g_o^c,\s_o,\d_o) \in \widetilde{\NC_o}$ agit diagonalement.

\end{itemize}
\end{theo}

\begin{proof} Pour montrer que les correspondances de Hecke se décrivent comme indiqué au point (2), il suffit de
travailler sur les germes aux points géométriques de $\IC_{I^o,m}^{=h}$ en remarquant que les diagrammes
suivant sont commutatifs, cf. la proposition (\ref{comute}):

$$\xymatrix{
\Psi_{F_o,t'}^{h,i} \ar_{(\lexp t (g_o^c)^{-1},\d_o,\s_o)}[d] \ar^{f_{x_{J^o,m'}(\oo)}}[rrrrrr] & & & & & &
R^i \Psi_{\eta_o}(\bar \Qm_l)_{\Tw_{J^o,h,m',t'}(\oo)_{x_{J^o,m'}}} \ar_{(g^{\oo,o},g_o^{et},(g_o^c,\d_o,\s_o))}[d] 
\\
\Psi_{F_o,t}^{h,i} \ar^{f_{(g^{\oo,o},g_o^{et},(g_o^c,\d_o,\s_o))x_{J^o,m'}(\oo)}}[rrrrrr] & & & & & &
R^i \Psi_{\eta_o}(\bar \Qm_l)_{\Tw_{I^o,h,m,t}(\oo)_{x_{I^o,m}}}
}$$

$$\xymatrix{
\Psi_{F_o,t'}^{h,i} \ar_{(\lexp t (g_o^c)^{-1},\d_o,\s_o)}[d] \ar^{h_{x{J^o,m'}(_\oo)}}[rrrrrr] & & & & & &
(\twist_{\JC_{J^o,m'}^{=h}(\oo)} \Psi_{F_o,t'}^{h,i})_{x_{j^o,m'}} \ar_{(g^{\oo,o},g_o^{et},(g_o^c,\d_o,\s_o))}[d] \\
\Psi_{F_o,t}^{h,i} \ar^{h_{(g^{\oo,o},g_o^{et},(g_o^c,\d_o,\s_o))x_{J^o,m'}(\oo)}}[rrrrrr] & & & & & &
(\twist_{\JC_{I^o,m}^{=h}(\oo)} \Psi_{F_o,t}^{h,i})_{x_{I^o,m}} }$$ où $\d_o$ est un \ele quelconque de 
$D_{o,h}^\times$
tel que $(g_o^c,\d_o,\s_o) \in \widetilde{\NC_o}$, et $J^o,I^o$, $m,m',t,t'$ sont comme au paragraphe
(\ref{enonce-local}), avec
$$(g^{\oo,o},g_o^{et},(g_o^c,\d_o,\s_o))(x_{J^o,m'}(\oo))=x_{I^o,m}(\oo).$$

\end{proof}

\begin{coro} \label{coro-iso}
Il existe un isomorphisme canonique
$$R^i \Psi_{\eta_o}(\bar \Qm_l)_{\hat \IC_{I^o,=h,m}(t)} \longto \twist_{\JC_{I^o,m}^{=h}(\oo)}
(\widetilde{\Psi_{F_o,t}^{h,i}})$$ tel qu'en tout point fermé $x$ de $\IC_{I^o,m}^{=h}$, le diagramme suivant
soit commutatif
$$\diagram
(R^i \Psi_{\eta_o}(\bar \Qm_l)_{\hat \IC_{I^o,=h,m}(t)})_x \rto & (\twist_{\JC_{I^o,m}^{=h}(\oo)}
(\widetilde{\Psi_{F_o,t}^{h,i}}))_x \\
\widetilde{\Psi_{F_o,t}^{h,i}} \rto^{\Id} \uto^{j_{x(\oo)}} & \widetilde{\Psi_{F_o,t}^{h,i}}
\uto^{h_{x(\oo)}}
\enddiagram$$
où les flèches verticales dépendent du choix d'un point fermé $x(\oo)$ de $\JC_{I^o,m}^{=h}(\oo)$ au dessus
de $x$. En outre ces isomorphismes sont compatibles aux correspondances de Hecke associées aux \eles de
$(D_\Am^{\oo,o})^\times \times GL_{d-h}(F_o) \times GL_h(F_o) \times W_{F_o}$ de la manière suivante: si
$(g^{\oo,o}, g_o^{et}, g_o^c)$ est un \ele de $(D_\Am^{\oo,o})^\times \times GL_{d-h}(F_o) \times GL_h(F_o)$,
qui induit donc, en particulier, une correspondance de Hecke $[g]$ sur $\IC_{I^o,m}^{=h}$, on a le diagramme
commutatif suivant
$$\diagram
([g],(\fr_o^{\deg \s_o})^*)^* R^i \Psi_{\eta_o}(\bar \Qm_l)_{\hat \IC_{I^o,=h,m'}(t')}
\rrto^{~(g^{\oo,o},g_o^{et},g_o^c,\s_o)} \dto  & &
 R^i \Psi_{\eta_o}(\bar \Qm_l)_{\hat \IC_{I^o,=h,m}(t)} \dto \\
([g],(\fr_o^{\deg \s_o})^*)^* \twist_{\JC_{I^o,m'}^{=h}(\oo)}(\widetilde{\Psi_{F_o,t'}^{h,i}})
\rrto^{\genfrac{}{}{0pt}{}{\hspace{1.7cm} \s_o^* \circ (g^{\oo,o},g_o^{et},g_o^c,\s_o)}{~}}
 & & \twist_{\JC_{I^o,m}^{=h}(\oo)}(\widetilde{\Psi_{F_o,t}^{h,i}})
\enddiagram$$
où l'action de $(g^{\oo,o},g_o^{et},g_o^c,\s_o)$ sur $\twist_{\JC_{I^o,m}^{=h}(\oo)}
(\widetilde{\Psi_{F_o,t}^{h,i}})$ est précisée au point (2) du théorème précédent.
\end{coro}

\begin{proof} La démonstration découle directement du résultat précédent et de la proposition (\ref{iso-local})
notamment la description des correspondances de Hecke se montre comme dans la preuve du théorème ci-dessus.
On peut toutefois en donner une preuve directe à partir de la proposition (\ref{niv-fini}). D'après le
théorème 4.1 de \cite{berk}, pour tout entier $r$, on peut choisir un entier $N$, tel que deux morphismes de
schéma formels sur $\hat \OC_o^{nr}$
$$\spf \Def_t^h[[X_1,\cdots,X_h]] \longto \spf \Def_t^h$$
qui coïncide sur $\spf (\Def_t^h/(\MC_t^h)^N)[[X_1,\cdots,X_h]]$, induisent la même application sur les cycles 
évanescents
du faisceau constant $\Zm/ l^r \Zm$. De la proposition (\ref{niv-fini}), on en déduit donc un morphisme
$$f_s: R^i \Psi_{\eta_o}(\Zm / l^r \Zm)_{S_s \times_{\spec \k(o)} \spf \Def_t^h}=\widetilde{\Psi_{F_o,t}^{h,i}} 
\longto R^i
\Psi_{\eta_o}(\Zm/l^r \Zm)_{\hat S_s}$$
tel que pour tout point fermé $x(\oo)$ de $\JC_{I^o,m}^{=h}(\oo)$, le morphisme
$$\widetilde{\Psi_{F_o,t}^{h,i}} \longto (R^i \Psi_{\eta_o}(\Zm/l^r \Zm)_{\hat S_s})_{x_s}$$
est le morphisme $j_{x(\oo)}$, c'est en particulier un
isomorphisme. Ce dernier se descend en un isomorphisme
$$\twist_{\JC_{I^o,m}^{=h}(s)} \widetilde{\Psi_{F_o,t}^{h,i}} \longto R^i \Psi_{\eta_o}(\Zm/l^r \Zm)_{\hat 
\IC_{I^o,=h,m}(t)}$$
tel qu'en tout point fermé $x(\oo)$ de $\JC_{I^o,m}^{=h}(\oo)$ au dessus d'un point $x$ de $\IC_{I^o,m}^{=h}$, on ait 
le diagramme
commutatif suivant
$$\diagram
(\twist_{\JC_{I^o,m}^{=h}(s)} \widetilde{\Psi_{F_o,t}^{h,i}})_x \rto & (R^i \Psi_{\eta_o}(\Zm/l^r \Zm)_{\hat 
\IC_{I^o,=h,m}(t)})_x \\
\widetilde{\Psi_{F_o,t}^{h,i}} \rto^{\Id} \uto^{h_{x(\oo)}} & \widetilde{\Psi_{F_o,t}^{h,i}}
\uto^{j_{x(\oo)}}
\enddiagram$$
où les flèches verticales définies plus haut, dépendent du choix de $x(\oo)$.
En travaillant au niveau des germes, on voit que ces isomorphismes sont compatibles lorsque $r$ varie, de sorte que 
l'on
peut les recoller en un isomorphisme vérifiant les propriétés de l'énoncé.

\end{proof}

\subsection{Systèmes locaux d'Harris-Taylor}

\label{action-cycle}

Soit $\t_o$ une représentation admissible irréductible de $D_{o,h}^\times$. L'espace que l'on souhaite
étudier est le $(GL_h(F_o) \times W_{F_o})$-module
$$\hom_{D_{o,h}^\times}(\t_o,\UC_{F_o}^{h,i})=\UC_{F_o}^{h,i}(\t_o) \simeq 
\Psi_{F_o}^{h,i}(\t_o):=\hom_{\DC_{o,h}^\times}(\t_o,\Psi_{F_o}^{h,i})$$
ou de manière équivalente leur version avec un tilde, $\widetilde{\UC_{F_o}^{h,i}(\t_o)}$ qui désigne
l'espace $\UC_{F_o}^{h,i}(\t_o)$ muni de l'action tordue de $GL_h(F_o) \times W_o$ par $g_o^c \mapsto \lexp t
(g_o^c)^{-1}$.

\begin{defi} Soit $\t_o$ une représentation irréductible de $D_{o,h}^\times$, sa restriction à
$\DC_{o,h}^\times$ est une somme de représentations irréductibles
$$\r_{o,1} \oplus \cdots \oplus \r_{o,e_{\t_o}}$$
et on notera $e_{\t_o}$ le nombre de celles ci. Étant donnée une représentation irréductible $\r_o$ de
$\DC_{o,h}^\times$, soient alors $\t_o$ et $\t_o'$ des sous-représentations irréductibles de l'induite de
$\DC_{o,h}^\times$ à $D_{o,h}^\times$ de $\r_o$: d'après la réciprocité de Frobenius, ce sont exactement
celles telles que leur restriction à $\DC_{o,h}^\times$ contienne $\r_o$. On en déduit alors que $\t_o$ et
$\t_o'$ sont inertiellement équivalentes, i.e. $\t_o' \simeq \t_o \otimes \xi_o$ avec $\xi_o:\d \mapsto
x^{v(\det \d)}$ pour $x \in \bar \Qm_l^\times$. On note  $\CF_h$ l'ensemble des classes d'équivalences
inertielles des représentations admissibles et irréductibles du groupe $D_{o,h}^\times$.
\end{defi}

\marque Pour tout $\t_o$, on a un morphisme naturel de $\NC_o$-modules:
$$\Psi_{F_o}^{h,i}(\t_o) \otimes \t_o \longto \Psi_{F_o}^{h,i}$$
qui envoie $f \otimes v$ sur $f(v)$. On note $\Psi_{F_o}^{h,i}[\t_o]$ l'image de ce morphisme et soit
$\Psi_{F_o,m}^{h,i}[\t_o]$ la préimage de $\Psi_{F_o}^{h,i}[\t_o]$ dans $\Psi_{F_o,m}^{h,i}$.
Le sous-module $\Psi_{F_o}^{h,i}[\t_o]$ ne dépend que de la classe d'équivalence inertielle de $\t_o$.
Le groupe $\DC_{o,h}^\times$ étant compact, on a
$$\Psi_{F_o}^{h,i} = \bigoplus_{\t_o \in \CF_h} \Psi_{F_o}^{h,i}[\t_o]$$

\marque Soit $e_{\t_o}$ le nombre de composantes irréductibles de $(\t_o)_{|\DC_{o,h}^\times}$; c'est aussi
le nombre de caractère $\xi: \Zm \longto \bar \Qm_l^\times$ tels que $\t_o \simeq \t_o \otimes (\xi \circ
\val \circ \det)$. Comme $\Pi_{o,h}^h$ est dans le centre de $D_{o,h}^\times$, $e_{\t_o}$ divise $h$. Soit
$\Delta_{\t_o}$ un ensemble d'\eles de $D_{o,h}^\times$ tel que les congruences des $\val (\det \d)$ pour $\d
\in \Delta_{\t_o}$ forment un système de représentants de $\Zm/e_{\t_o} \Zm$. On a alors le lemme suivant
dont la preuve est claire.

\begin{lemm} L'application
$$\begin{array}{rl}
\Psi_{F_o}^{h,i}(\t_o) \otimes \t_o & \longto \bigoplus_{\d \in \Delta_{\t_o}} \Psi_{F_o}^{h,i}[\t_o]^\d \\
f \otimes v & \mapsto (f(\d^{-1}v))_\d
\end{array}$$
où $\Psi_{F_o}^{h,i}[\t_o]^\d$ est l'espace $\Psi_{F_o}^{h,i}[\t_o]$ muni de la structure de $\NC_o$-module où 
$(g,d,w)$ agit via
$(g,\d^{-1}d \d ,w)$, est un isomorphisme de $\NC_o$-modules.
\end{lemm}

\begin{lemm}  \label{Fto}
Soit $\FC_{\t_o,1}$ le système local sur $M_{I,s_o,1}^{=h}$ associé à $\t_o$ et au revêtement d'Igusa de
seconde espèce. On a alors un isomorphisme naturel de $(D_\Am^{\oo,o})^\times \times GL_{d-h}(F_o) \times
GL_h(F_o) \times W_{F_o}$-modules
$$(\twist_{\JC_{I^o,m}^{=h}(\oo)}(\widetilde{\Psi_{F_o,t}^{h,i}}))^h \simeq \bigoplus_{\t_o \in \CF_h} (\FC_{\t_o,1} 
\otimes
\widetilde{\Psi_{F_o,t}^{h,i}(\t_o)})^{h/e_{\t_o}}$$
\end{lemm}

\begin{proof} D'après le lemme précédent, on a
$$\twist_{\JC_{I^o,m}^{=h}(\oo)} (\widetilde{\Psi_{F_o,t}^{h,i}})=\bigoplus_{\t_o \in \CF_h} 
\twist_{\JC_{I^o,m}^{=h}(\oo)}
(\widetilde{\Psi_{F_o,t}^{h,i}}[\t_o])$$
et
$$\FC_{\t_o,1} \otimes \widetilde{\Psi_{F_o,t}^{h,i}(\t_o)} \simeq \bigoplus_{\d \in \Delta_{\t_o}} 
\twist_{\JC_{I^o,m}^{=h}(\oo)}
(\widetilde{\Psi_{F_o,t}^{h,i}}[\t_o]^\d)$$ D'après le point (2) du théorème (\ref{theo-principal}), tous les
$\twist_{\JC_{I^o,m}^{=h}(\oo)}(\widetilde{\Psi_{F_o,t}^{h,i}}[\t_o]^\d)$ sont isomorphes en tant que
$(D_\Am^{\oo,o})^\times \times GL_{d-h}(F_o) \times GL_h(F_o) \times W_{F_o}$-modules, d'où le résultat.

\end{proof}

\rem Pour tout $t \leq t' \leq \oo$, on a
$$\twist_{\JC_{I^o,m}^{=h}(\oo)}(\widetilde{\Psi_{F_o,t}^{h,i}})=(\twist_{\JC_{I^o,m}^{=h}(\oo)}
(\widetilde{\Psi_{F_o,t'}^{h,i}}))^{K_{o,t}}$$
$$\twist_{\JC_{I^o,m}^{=h}(\oo)}(\widetilde{\Psi_{F_o,t}^{h,i}}[\t_o])=(\twist_{\JC_{I^o,m}^{=h}(\oo)}
(\widetilde{\Psi_{F_o,t'}^{h,i}}[\t_o]))^{K_{o,t}}$$ où l'on rappelle que $K_{o,t}=\ker (GL_h(\OC_o) \longto
GL_h(\OC_o/\MC_o^t))$.

D'après le corollaire (\ref{coro-iso}) et le théorème (\ref{theo-principal}), on en déduit le résultat suivant.

\begin{prop} \label{prop-principal}
\begin{itemize}

\item[(i)] $(R^i \Psi_{\eta_o}(\bar \Qm_l)_{\IC_{I^o,m}^{=h}})^h$ est isomorphe à
$$\bigoplus_{\t_o \in \CF_h} (\FC_{\t_o,1} \otimes \widetilde{\UC_{F_o,m}^{h,i}(\t_o)})^{h/e_{\t_o}}.$$

\item[(ii)] L'action de $(D_\Am^{\oo,o})^\times \times GL_h(F_o) \times GL_{d-h}(F_o) \times W_{F_o}$
sur ${\DS \lim_{\genfrac{}{}{0pt}{}{\lefto}{m}}} ~\IC_{I^o,m}^{=h}$ par les correspondances de Hecke, induit
des morphismes $(g^{\oo,o} \times g_o^c \times g_o^{et} \times \s_o)$:
$$(g^{\oo,o},\frob_o^{\val(\det g_o^c)+\deg \s_o},g_o^{et})^* R^i \Psi_{\eta_o}(\bar
\Qm_l)_{\IC_{I^o,m}^{=h}} \longto R^i \Psi_{\eta_o}(\bar \Qm_l)_{\IC_{J^o,m'}^{=h}}$$
(pour $m,m'$ et $I^o,J^o$ convenables) qui induisent par (i) le morphisme
$$\diagram
\bigoplus_{\t_o \in \CF_h}[(\frob_o^{\deg \s_o})^* \circ (g^{\oo,o},g_o^{et},(g_o^c,\d_o,\s_o))^* ]
(\FC_{\t_o,1} \otimes \widetilde{\UC_{F_o,m}^{h,i}(\t_o)}) \dto^{(\s_o^* \circ (g^{\oo,o},g_o^{et},g_o^c) \otimes 
(g_o^c,\d_o,\s_o))} \\
\bigoplus_{\t_o \in \CF_h} (\FC_{\t_o,1} \otimes \widetilde{\UC_{F_o,m'}^{h,i}(\t_o)})
\enddiagram$$
où $\d_o$ est un \ele quelconque de $D_{o,h}^\times$ tel que $(g_o^c,\d_o,\s_o) \in \widetilde{\NC_o}$
\footnote{On rappelle que $\widetilde{\NC_o}$ est le noyau de $(g_o,d_o,c_o) \in GL_h(F_o) \times
D_{o,h}^\times \times W_{F_o} \mapsto \val(\det g_o \rn (d_o) \cl(c_o)) \in \Zm$, on peut par exemple prendre
$\d_o=\Pi_{o,h}^{-\deg \s_o - \val(\det g_o^c)}$, où $\Pi_{o,h}$ est une uniformisante de $\DC_{o,h}$ telle
que $\Pi_{o,h}^h=\varpi_o$.}

\item[(iii)] on a
$$R^i \Psi_{\eta_o}(\bar \Qm_l)_{\IC_{I^o,m}^{=h}} \simeq (\twist_{\JC_{I^o,m}^{=h}(\oo)}
(\widetilde{\Psi_{F_o}^{h,i}}))^{1+\varpi_o^m \Mm_h(\OC_o)}.$$
\end{itemize}
\end{prop}

\begin{remas} \label{rema-action-Z}
 (i) Les systèmes locaux $\FC_{\t_o,1}$ ne sont pas irréductibles mais plutôt une somme directe de
$e_{\t_o}$ systèmes locaux irréductibles. La complexité de l'écriture du point (i), est la contre-partie de
la simplicité de la description de l'action de $GL_d(F_o) \times W_o$ qui tient au fait que l'on a fait
apparaître $\widetilde{\UC_{F_o}^{h,i}}(\t_o)$ plutôt que
$\hom_{\DC_{o,h}^\times}(\r_o,\widetilde{\Psi_{F_o}^{h,i}})$ qui est seulement une représentation de
$\widetilde{\NC_o} \cap (GL_h(F_o) \times D_{o,h}^{\times})$.

(ii) Soit $\Zm \longmapright{\sim} \widetilde{\NC_o}/\DC_{o,h}^\times$ l'isomorphisme défini par $n \mapsto
\Pi_{o,h}^n$. On considère ainsi
$${\DS \lim_{\genfrac{}{}{0pt}{}{\longto}{I^o,n}}} H^j_c(\IC_{I^o,n}^{=h},\FC_{\t_o,1})$$
comme une représentation de
$$(D_\Am^{\oo,o})^\times \times GL_{d-h}(F_o) \times \Zm \hbox{ ou de } (D_\Am^{\oo,o})^\times \times GL_{d-h}(F_o)
\times(D_{o,h}^\times/ \DC_{o,h}^\times)$$
\end{remas}

En vertu du lemme (\ref{lem-lisse}), on en déduit alors le corollaire suivant

\begin{coro} \label{corofond1}
On a un isomorphisme canonique
$$H^j_c(M_{I,s_o,1}^{=h}, R^i\Psi_{\eta_o}(\LC_{\r_\oo}))^h \simeq \bigoplus_{\t_o \in \CF_h}
(H^j_c(M_{I,s_o,1}^{=h},\FC_{\t_o,1} \otimes \LC_{\r_\oo}) \otimes
\widetilde{\UC_{F_o,m}^{h,i}(\t_o)})^{h/e_{\t_o}}$$ tel que l'action de $(g^{\oo,o},g_o^c,g_o^{et},\s_o) \in
(D_\Am^{\oo,o})^\times \times GL_h(F_o) \times GL_{d-h}(F_o) \times W_{F_o}$ sur la limite inductive indexée
par les idéaux $I=I^o\MC_o^m$ du membre de gauche, induit l'action de
$$(g^{\oo,o},g_o^{et},-\val(\det g_o^c)-\deg(\s_o)) \otimes (g_o^c,\s_o)$$
sur la limite inductive du membre de droite.
\end{coro}

\begin{defi} On note $\FC_{\tau_o}:= \FC_{\tau_o,1} \times_{P_{h,d}^{op}(F_o)} GL_d(F_o)$ le faisceau sur 
$M_{I,s_o}^{=h}$
induit par $\FC_{\tau_o,1}$.
\end{defi}

On déduit du corollaire précédent la proposition suivante dont la proposition (\ref{prop-segment}) nous
permettra au paragraphe suivant de donner le cas Iwahori du théorème local.

\begin{prop} \label{prop-poids}
Soient $1 \leq lg <d$, $\pi_o$ une représentation irréductible cuspidale de $GL_g(F_o)$ et $\Pi_l$ une
représentation de $GL_{lg}(F_o)$. Pour $1 \leq i \leq l$ et pour tout $j$,
$$\lim_{\genfrac..{0pt}{1}{\to}{I}} H^j_c(M_{I,s_o}^{=lg}, \FC_{\JL^{-1}(\st_l(\pi_o))} \otimes \LC_{\rho_\oo}) 
\otimes \Pi_l
\otimes L_g(\pi_o)$$ est muni d'une action naturelle de $(D_\Am^{\oo,o})^\times \times P_{lg,d}(F_o)^{op}
\times W_o$ telle que $(g^{\oo,o},g_o^c \times g_o^{et},c_o)$ agisse via l'action naturelle de
$(g^{\oo,o},-\val(\det(g_o^c))-\deg(c_o)) \in (D_\Am^{\oo,o})^\times \times \Zm$ et celle de $(g_o^c,c_o)$
sur $\Pi_l \otimes L_g(\pi_o)$. Cet espace, en tant que représentation de $GL_d(F_o) \times W_o$ est alors de
la forme
$$\bigoplus_\xi(\ind_{P_{lg,d}^{\op}(F_o)}^{GL_d(F_o)} \Pi_l \circ \xi(\val(\det))\otimes \pi_\xi)\otimes
L_g(\pi_o \circ \xi(\val(\det)))$$ où $\xi$ décrit les caractères $\Zm \longto \bar \Qm_l^\times$ et
$\pi_\xi$ est une représentation de $GL_{d-lg}(F_o)$.
\end{prop}

\begin{proof} Le faisceau $\FC_{\JL^{-1}(\st_l(\pi_o))}$ étant induit à partir de $\FC_{\JL^{-1}(\st_l(\pi_o)),1}$
on a:
\begin{multline} \label{iso-induite1}
\lim_{\genfrac..{0pt}{1}{\to}{I}} H^j_c(M_{I,s_o}^{=lg},\FC_{\JL^{-1}(\st_l(\pi_o))} \otimes \LC_{\rho_\oo}
\otimes \Pi_l)
\otimes L_g(\pi_o)) \simeq \\
\ind_{P_{lg,d}^{\op}(F_o)}^{GL_d(F_o)} \lim_{\genfrac..{0pt}{1}{\to}{I}}
H^j_c(M_{I,s_o,1}^{=lg},\FC_{\JL^{-1}(\st_l(\pi_o)),1} \otimes \LC_{\rho_\oo} \otimes \Pi_l) \otimes
L_g(\pi_o)
\end{multline}
tel que l'action de $((g_o^c,g_o^{et}),c_o) \in P_{lg,_d}(F_o)^{op} \times W_o$ sur le membre de droite de
(\ref{iso-induite1}), soit donnée par l'action de $(g_o^{et},-\val(\det g_o^c)-\deg(c_o)) \times (g_o^c,c_o)
\in GL_{d-lg}(F_o) \times \Zm \times GL_{lg}(F_o) \times W_o$ où ${\DS \lim_{\genfrac..{0pt}{1}{\to}{I}}}
H^j_c(M_{I,s_o,1}^{=lg},\FC_{\JL^{-1}(\st_l(\pi_o))} \otimes \LC_{\rho_\oo})$ (resp. $\Pi_l \otimes
L_g(\pi_o)$) est vue comme un $GL_{d-lg}(F_o) \times \Zm$-module (resp. $GL_{lg}(F_o) \times W_o$-module). Le
résultat découle alors de l'écriture
$$\lim_{\genfrac..{0pt}{1}{\to}{I}} H^j_c(M_{I,s_o,1}^{=lg},\FC_{\JL^{-1}(\st_l(\pi_o))} \otimes \LC_{\rho_\oo}) 
\simeq
\bigoplus_{\xi} \xi \otimes \pi_\xi$$
où $\xi$ décrit les caractères $\Zm \longto \bar \Qm_l^\times$ et où $\pi_\xi$ est une représentation de
$GL_{d-lg}(F_o)$, à priori réductible.

\end{proof}


%% file: calcul.tex
\section{Equations formelles des variétés d'Igusa de seconde espèce}

On propose de donner au niveau des complétés formels des équations explicites du revêtement d'Igusa de
seconde espèce, qui fournissent les systèmes locaux $\FC_{\t_o}$ associés à une représentation irréductible
$\t_o$ de $D_{o,d}^\times$.

\begin{lemb} \label{ext-def}
Soit $z$ un point géométrique de $M_{I^o,s_o}^{h'}$ avec $h' > h$. Le revêtement $\JC_{I^o,0}^{=h} \longto 
M_{I^o,s_o}^{=h}$ donne
au dessus du complété formel $\Def_0^{h'd-h';h}$ de $M_{I^o,s_o}^h$ en $z$, une extension
$$\Def_0^{h',d-h';=h} \hookrightarrow J_0^{h',d-h';=h}(s)$$
qui provient d'une extension $\Def_0^{h';=h} \longto J_0^{h';=h}(s)$; i.e.
$J_0^{h',d-h';=h}(s)$ est isomorphe à
$$J_0^{h';=h}(s) \hat \otimes_{\Def_0^{h';=h}} \Def_0^{h',d-h';=h}$$
où l'on rappelle que $\Def_0^{h',d-h';=h}=\Def_0^{h';=h}[[w_1^0,\cdots,w_{d-h'}^0]]$.
\end{lemb}

\begin{proof} Par définition des variétés d'Igusa de seconde espèce (cf. le paragraphe (\ref{defi-2})), l'extension
$\Def_0^{h',d-h';=h} \hookrightarrow J_0^{h',d-h';=h}(s)$ est donnée par la rigidification à l'ordre $s$ de la partie 
connexe
du $\OC_o$-module de Dieudonné universel sur $\Def_0^{h',d-h';=h}$. Le résultat découle alors de la forme de ce 
$\OC_o$-module
de Dieudonné sur $\Def_0^{h',d-h';=h}$ donnée à la proposition (9.9.1) de \cite{boy}.

\end{proof}

\marque Soit alors $z$ un point géométrique de $M_{I^o,s_o}^{=h'}$ avec $h' >h$. Le complété de l'anneau
local de $M_{I^o,s_o}^h$ en $z$ est $\Def_0^{h',d-h';h} \simeq \Def_0^{h';h}[[w_1,\cdots,w_{d-h'}]]$ avec
$\Def_0^{h';h} \simeq \hat \OC_o^{nr} [[a_h,\cdots,a_{h'-1}]].$ L'ouvert correspondant à $M_{I^o,0}^{=h}$ est
$$\Def_0^{h',d-h';=h} \simeq \Def_0^{h';=h}[[w_1,\cdots,w_{d-h'}]]$$
avec$\Def_0^{h';=h} \simeq \hat \OC_o^{nr}((a_h)) [[a_{h+1},\cdots,a_{h'-1}]].$
D'après le lemme précédent, l'extension
$\Def_0^{h';=h} \hookrightarrow J_0^{h';=h}(s)$ est donnée par l'équation
$$\t^h \circ \a=\a \circ (a_h \t^h+ \cdots + a_{h'-1} \t^{h'-1}+\t^{h'}) \qquad \hbox{ modulo } \t^{h+m}$$
avec $\a=\sum_{i=0}^{m-1} \a_i \t^i$ avec $\a_0$ non nul. On obtient alors le système suivant
$$\left \{ \begin{array}{l}
\a_0^{p^h-1}=a_h \\
\a_1^{p^h}=\a_0 a_{h+1} + \a_1 a_h^p \\
\cdots \\
\a_{h'-h-1}^{p^h}=\a_0 a_{h'-1}+\a_1 a_{h'-2}^p+ \cdots + \a_{h'-h-1} a_h^{p^{h'-h-1}}\\
\a_{h'-h}^{p^h}=\a_0+ \a_1 a_{h'-1}^p+ \cdots + \a_{h'-h} a_h^{p^{h'-h}} \\
\cdots \\
\a_{m}^{p^h}=\a_{m-h'+h}+ \a_{m-h'+h+1} a_{h'-1}^{p^{m-h'+h}}+ \cdots + \a_{m} a_h^{p^{m-1}}
\end{array} \right.
$$

On peut effectuer le changement de variable $t_i=\a_i t_0^{-p^i}$ pour $i>0$ avec $t_0=\a_0$. Le système précédent 
s'écrit alors
$$\left \{
\begin{array}{l}
t_0^{p^h-1} =a_h \\
t_1^{p^h}-t_1 =a_{h+1} t_0^{1-p^{h+1}} \\
t_2^{p^h}-t_2 = a_{h+2} t_0^{1-p^{h+2}} + (a_{h+1}t_0^{1-p^{h+1}})^p t_1 \\
\cdots  \\
t_{h'-h-1}^{p^h}-t_{h'-h-1} = a_{h'-1} t_0^{1-p^{h'-1}}+\sum_{k=1}^{h'-h-2} t_k (a_{h'-1-k} t_0^{1-p^{h'-1-k}})^{p^k} 
\\
t_{h'-h}^{p^h}-t_{h'-h} = t_0^{1-p^{h'}} + \sum_{k=1}^{h'-h-1} t_k (a_{h'-k}t_0^{1-p^{h'-k}})^{p^k} \\
t_{h'-h+1}^{p^h}-t_{h'-h+1} = t_1 t_0^{p-p^{h'+1}}+\sum_{k=1}^{h'-h-1} t_{h'-h+1-k} (a_{h+k} 
t_0^{1-p^{h+k}})^{p^{h'-h+1-k}} \\
\cdots  \\
t_s^{p^h}-t_s = t_{s-h'+h} t_0^{p^{s-h'+h}-p^{s+h}} + \sum_{k=1}^{h'-h-1} t_{s-k} (a_{h+k} t_0^{1-p^{h+k}})^{p^{s-k}} 
\\
\end{array}
\right .
$$
Pour $i \leq h'-h$ (resp. $i \geq h'-h$), on peut éliminer dans le membre de droite de l'équation dont le
membre de gauche est $t_i^{p^h}-t_i$, les termes
$t_0^{1-p^{h+k}}a_{h+k}$ pour $k < i$ (resp. $1 \leq k \leq h'-h-1$ et $t_0^{1-p^{h'}}$); on obtient alors le système 
suivant:

$$\left \{
\begin{array}{l}
t_0^{p^h-1} = a_h \\
t_1^{p^h}-t_1 = a_{h+1} t_0^{1-p^{h+1}} \\
t_2^{p^h}-t_2 = a_{h+2} t_0^{1-p^{h+2}} +r_2 \\
\cdots  \\
t_{h'-h-1}^{p^h} -t_{h'-h-1} = a_{h'-1} t_0^{1-p^{h'-1}}+r_{h'-h-1} \\
t_{h'-h}^{p^h} - t_{h'-h} = t_0^{1-p^{h'}}+r_{h'-h} \\
t_{h'-h+1}^{p^h} - t_{h'-h+1} = r_{h'-h+1} \\
\cdots  \\
t_{s}^{p^h}-t_s =r_s \\
\end{array}
\right . $$
avec $r_i \in \Fm_q[[t_1,\cdots,t_{i-1}]]$ pour $1 \leq i \leq s$ définis par les relations suivantes:
\begin{itemize}
\item $r_1=0$ et pour $2 \leq i \leq h'-h$, $r_i= \sum_{k=1}^{i-1} t_{i-k} (t_k^{p^h}-t_k-r_k)^{p^{i-k}}$;

\item pour $i > h'-h$, $r_i=\sum_{k=1}^{h'-h} t_{i-k}(t_k^{p^h}-t_k-r_k)^{p^{i-k}}$.
\end{itemize}


%% file: introduction-chap2.tex
\noindent \textbf{0.1.} --- On commence par des rappels sur les  induites paraboliques d'après \cite{ze}. Les
paraboliques que l'on considère sont les paraboliques opposés aux paraboliques standards. On s'intéresse tout
particulièrement aux représentations elliptiques, d'après une terminologie de Dat, i.e. aux sous-quotients
irréductibles de l'induite normalisée du parabolique $P_{g,2g,\cdots,sg}^{\op}(F_o)$ à $GL_{sg}(F_o)$,
$$\pi_o(\frac{s-1}{2}) \times \cdots \times \pi_o(\frac{1-s}{2})$$
où $\pi_o$ est une représentation irréductible cuspidale de $GL_g(F_o)$, celles-ci seront dites de type
$\pi_o$. Ainsi on notera $[\overleftarrow{s-1}]_{\pi_o}$ (resp. $[\overrightarrow{s-1}]_{\pi_o}$) l'unique
quotient (resp. sous-espace) irréductible de cette induite, noté plus habituellement $\st_s(\pi_o)$ (resp.
$\speh_s(\pi_o)$). On rappelle que d'après \cite{ze}, ces sous-quotients irréductibles sont paramétrés par
les orientations $\vec \Gamma^s$ du graphe linéaire complet à $s$ sommets
$$\Gamma^s: \circ^{\frac{s-1}{2}} - \circ \cdots - \circ^{\frac{1-s}{2}}$$
Étant données deux représentations $\pi_1$ et $\pi_2$, elliptiques de type $\pi_o$, de respectivement
$GL_{s_1g}(F_o)$ et $GL_{s_2g}(F_o)$, les induites normalisées du parabolique
$P_{s_1g,(s_1+s_2)g}^{\op}(F_o)$ à $GL_{(s_1+s_2)g}(F_o)$,
$$\pi_1(-\frac{s_2(g-1)}{2}) \times \pi_2(\frac{s_1(g-1)}{2}) \hbox{ et }\pi_1(-\frac{s_2(g+1)}{2}) \times
\pi_2(\frac{s_1(g+1)}{2})$$ sont de longueur $2$, chacun des constituant étant une représentation elliptique
de type $\pi_o$, dont l'orientation de $\Gamma^{s_1+s_2}$ est une des deux qui prolonge celle du graphe
partiellement orienté obtenu par la concaténation du graphe de $\pi_1$ avec celui de $\pi_2$ (resp. du graphe
de $\pi_2$ avec celui de $\pi_1$)
$$\vec \Gamma^s:\overbrace{\circ \lefto \circ \cdots \to \circ}^{\vec \Gamma^{s_1}} \longleftrightarrow
\overbrace{\circ \to \circ \cdots \circ}^{\vec \Gamma^{s_2}}$$

$$(\hbox{resp. } \vec \Gamma^s:\overbrace{\circ \lefto \circ \cdots \to \circ}^{\vec \Gamma^{s_2}} 
\longleftrightarrow
\overbrace{\circ \to \circ \cdots \circ}^{\vec \Gamma^{s_1}})$$ On notera alors ces induites sous la forme
$\pi_1 \overrightarrow{\times} \pi_2$ (resp. $\pi_1 \overleftarrow{\times} \pi_2$). Ces notations sont
compatibles au calcul des foncteurs de Jacquet pour les paraboliques opposés au paraboliques standards tandis
qu'elles sont inversement orientées pour les foncteurs de Jacquet pour les paraboliques standards, cf. \S
\ref{jacquet}.

\medskip

\noindent \textbf{0.2.} --- On donne ensuite une preuve, théorème (\ref{theo-iwahori}), de
(\ref{theo-ripsi-local}) pour les représentations qui ont des vecteurs invariants sous le sous-groupe
d'Iwahori standard $\Iw_o$, i.e. les représentations elliptiques de type $1_o$. Par un argument
d'accouplement cohomologique, proposition (\ref{torsion2}), on obtient le résultat pour les représentations
elliptiques de type $\xi_o$, où $\xi_o$ est un caractère de $(\bar{{\mathbb Q}_l})^\times$.

Pour prouver le cas Iwahori, on utilise que $\Def_{st}^d:=(\Def_n^d)^{\Iw_o/K_{o,n}}$ muni de son morphisme
structurel, est semi-stable (cf. la proposition (\ref{model-st})) et possède une interprétation modulaire
très simple en termes de drapeaux (cf. (\ref{drapeaux})). On dispose alors d'une description explicite du
complexe des cycles évanescents associé, rappelée au théorème (\ref{ill}), et on en déduit en particulier au
corollaire (\ref{dim-poids}) que pour toute représentation irréductible $\t_o$ de $D_{o,d}^\times$,
$(\Psi_{F_o}^{d,i}(\t_o))^{\Iw_o}$ est un $\bar \Qm_l$-espace vectoriel de dimension le coefficient binomial
$(\genfrac{}{}{0pt}{}{i}{d-1})$, qui est pur de poids $2i$. Le reste de la preuve est alors purement
combinatoire à partir de la proposition (\ref{prop-segment}).

%% file: iwahori.tex
\section{Rappels sur les représentations de $GL_d(F_o)$}
\label{rappel-ze-0}

Tous les énoncés ainsi que leur preuve peuvent être trouvés dans \cite{ze}.

\subsection{Induites paraboliques}
\label{defi-induite}

\begin{defi} \label{defi-parab}
Pour une suite $0< r_1 < r_2 < \cdots < r_k=d$, on note $P_{r_1,r_2,\cdots,r_k}$ le sous-groupe parabolique
de $GL_d$ standard associé au sous-groupe de Levi $GL_{r_1}(F_o) \times GL_{r_2-r_1}(F_o) \times \cdots
\times GL_{r_k-r_{k-1}}(F_o)$ et soit $P_{r_1,\cdots,r_k}^{\op}$ le parabolique opposé dont on note
$N_{r_1,\cdots,r_k}^{\op}$ le radical unipotent.
\end{defi}

\marque Soient $\pi_1$ et $\pi_2$ des représentations de respectivement $GL_{n_1}(F_o)$ et $GL_{n_2}(F_o)$;
on note selon la coutume, $\pi_1 \times \pi_2$ l'induite parabolique
$$\pi_1 \times \pi_2:= \ind_{P_{n_1,n_1+n_2}^{\op}(F_o)}^{GL_{n_1+n_2}(F_o)} \pi_1(-n_2/2)
\otimes \pi_2(n_1/2)$$

\rem Le symbole $\times$ est associatif, i.e. $\pi_1 \times (\pi_2 \times \pi_3)=(\pi_1 \times \pi_2) \times
\pi_3$ que l'on notera donc $\pi_1 \times \pi_2 \times \pi_3$.

\begin{defis} Soit $g$ un diviseur de $d=sg$ et $\pi_o$ une représentation cuspidale irréductible de $GL_g(F_o)$:

- les sous-quotients irréductibles de
$$V(\pi_o,s):=\pi_o(\frac{s-1}{2}) \times \pi_o(\frac{s-3}{2}) \times \cdots \times \pi_o(\frac{1-s}{2})$$
seront dits elliptiques de type $\pi_o$;

- $V(\pi_o,s)$ possède un unique quotient (resp. sous-espace) irréductible que l'on notera
$$[\overleftarrow{s-1}]_{\pi_o} \qquad (\hbox{resp. } [\overrightarrow{s-1}]_{\pi_o});$$
c'est une représentation de Steinberg (resp. de Speh) généralisée notée habituellement $\st_s(\pi_o)$ (resp.
$\speh_s(\pi_o)$;

- pour $\pi_1$ et $\pi_2$ des représentations respectivement de $GL_{l_1g}(F_o)$ et $GL_{l_2g}(F_o)$, on
notera
$$\pi_1 \overrightarrow{\times} \pi_2 \qquad (\hbox{resp. } \pi_1 \overleftarrow{\times} \pi_2)$$
l'induite parabolique $\pi_1(l_2/2) \times \pi_2(-l_1/2)$ (resp. $\pi_1(-l_2/2) \times \pi_2(l_1/2)$).
\end{defis}

La propriété habituelle de transitivité des induites paraboliques donne le lemme suivant.

\begin{lemm} Pour $(\pi_i)_{1 \leq i \leq 3}$ des représentations de $GL_{l_ig}(F_o)$, on a les égalités suivantes:
$$\begin{array}{rl}
(\pi_1 \overrightarrow{\times} \pi_2)^\vee \simeq & \pi_1^\vee \overleftarrow{\times} \pi_2^\vee \\
(\pi_1 \overrightarrow{\times} \pi_2) \overrightarrow{\times_g} \pi_3 = & \pi_1 \overrightarrow{\times}
(\pi_2
\overrightarrow{\times} \pi_3) \\
(\pi_1 \overleftarrow{\times} \pi_2) \overleftarrow{\times_g} \pi_3= & \pi_1 \overleftarrow{\times} (\pi_2
\overleftarrow{\times} \pi_3)
\end{array}$$
En outre si $\pi_1$ et $\pi_2$ sont elliptiques de type $\pi_o$, il en est de même de $\pi_1
\overrightarrow{\times} \pi_2$ et donc de $\pi_1 \overleftarrow{\times} \pi_2$.
\end{lemm}

\rem Pour $\pi_1$, $\pi_2$ et $\pi_3$ des représentations elliptiques de type $\pi_o$, on a aussi
$$\pi_1 \overrightarrow{\times} (\pi_2 \overleftarrow{\times} \pi_3)= \pi_2 \overleftarrow{\times} (\pi_1
\overrightarrow{\times} \pi_3)$$ En effet d'après \cite{ze}, on a $\pi_1 \times \pi_2 \simeq \pi_2 \times
\pi_1$ si les supports cuspidaux de $\pi_1$ et $\pi_2$ sont disjoints \footnote{En fait loc. cit. donne un
critère plus précis; dans le cas où $\pi_1$ et $\pi_2$ sont irréductibles, il faut et il suffit que leurs
supports cuspidaux ne soient pas liés.}.

\begin{prop-defi} Pour $g$ un diviseur de $d=sg$ et $\pi_o$ une représen\-tation irréductible cuspidale de 
$GL_g(F_o)$, on a
pour $1 \leq l \leq s$:

\begin{itemize}
\item[-] $[\overleftarrow{l-1}]_{\pi_o} \overrightarrow{\times} [\overrightarrow{s-l-1}]_{\pi_o}$ (resp.
$[\overrightarrow{l-1}]_{\pi_o} \overrightarrow{\times} [\overleftarrow{s-l-1}]_{\pi_o}$) est de longueur
$2$; on notera
$$[\overleftarrow{l-1},\overrightarrow{s-l}]_{\pi_o} \qquad (\hbox{resp. }
[\overrightarrow{l},\overleftarrow{s-l-1}]_{\pi_o})$$ son unique sous-espace irréductible, et
$$[\overleftarrow{l},\overrightarrow{s-l-1}]_{\pi_o} \qquad (\hbox{resp. }
[\overrightarrow{l-1},\overleftarrow{s-l}]_{\pi_o})$$ son unique quotient irréductible;

\item[-] par dualité $[\overleftarrow{l-1}]_{\pi_o} \overleftarrow{\times} [\overrightarrow{s-l-1}]_{\pi_o}$ (resp.
$[\overrightarrow{l-1}]_{\pi_o} \overleftarrow{\times} [\overleftarrow{s-l-1}]_{\pi_o}$) est de longueur $2$
avec
$$[\overrightarrow{s-l-1},\overleftarrow{l}]_{\pi_o} \qquad (\hbox{resp. }
[\overleftarrow{s-l},\overrightarrow{l-1}]_{\pi_o})$$ pour unique sous-espace irréductible et
$$[\overrightarrow{s-l},\overleftarrow{l-1}]_{\pi_o} \qquad (\hbox{resp. }
[\overleftarrow{s-l-1},\overrightarrow{l}]_{\pi_o})$$ pour unique quotient irréductible;

\item[-] on notera $\lfloor \pi \rfloor$ (resp. $\lceil \pi \rceil$) l'unique, s'il existe, sous-espace (resp. 
quotient)
irréductible de $\pi$. Pour $\pi_1$ et $\pi_2$ des représentations irréductibles elliptiques de type $\pi_o$,
$\pi_1 \overrightarrow{\times} \pi_2$ et $\pi_1 \overleftarrow{\times} \pi_2$ ont un unique sous-espace et un
unique quotient irréductible et on a:
$$\lfloor \pi_1 \overrightarrow{\times} \pi_2 \rfloor = \lceil \pi_2^\vee \overleftarrow{\times} \pi_1^\vee 
\rceil^\vee$$
$$\lceil \pi_1 \overrightarrow{\times} \pi_2 \rceil = \lfloor \pi_2^\vee \overleftarrow{\times} \pi_1^\vee 
\rfloor^\vee$$

Pour $\pi_3$ une troisième représentation irréductible elliptique de type $\pi_o$, on a
$$\lfloor \lfloor \pi_1 \overrightarrow{\times} \pi_2 \rfloor \overrightarrow{\times} \pi_3 \rfloor=\lfloor
\pi_1 \overrightarrow{\times} \lfloor \pi_2 \overrightarrow{\times} \pi_3 \rfloor \rfloor=\lfloor \pi_1
\overrightarrow{\times} \pi_2 \overrightarrow{\times} \pi_3 \rfloor$$

$$\lceil \lceil \pi_1 \overrightarrow{\times} \pi_2 \rceil \overrightarrow{\times}
\pi_3 \rceil=\lceil \pi_1 \overrightarrow{\times} \lceil \pi_2 \overrightarrow{\times} \pi_3 \rceil
\rceil=\lceil \pi_1 \overrightarrow{\times} \pi_2 \overrightarrow{\times} \pi_3 \rceil$$

\item[-] Soient $r \geq 1$ et $\Gamma^s=(a_i,\e_i)_{1 \leq i \leq r}$ tel que les $a_i$ sont des entiers strictement 
positifs
avec $s-1=a_1+\cdots + a_r$ et $\e_i= \pm 1$; on notera $\Gamma^s$ sous la forme $(\overleftarrow{a_1},
\cdots , \overrightarrow{a_r})$ où pour tout $i$ la flèche au dessus de $a_i$ est $\overleftarrow{a_i}$
(resp. $\overrightarrow{a_i}$) si $\e_i=-1$ (resp. $\e_i=1$). On associe à $\Gamma^s$ un sous-quotient
irréductible $[\Gamma^s]$ de $V(s,\pi_o)$ que l'on note aussi sous la forme $[\overleftarrow{a_1}, \cdots ,
\overrightarrow{a_r}]_{\pi_o}$. On convient par ailleurs des égalités suivantes:
$$[\cdots,\overleftarrow{a},\overleftarrow{b},\cdots]_{\pi_o}=[\cdots,\overleftarrow{a+b},\cdots]_{\pi_o}
\qquad
[\cdots,\overrightarrow{a},\overrightarrow{b},\cdots]_{\pi_o}=[\cdots,\overrightarrow{a+b},\cdots]_{\pi_o}$$
D'après \cite{ze}, on obtient alors une bijection, modulo les identifications ci-dessus, entre les
sous-quotients irréductibles de $V(s,\pi_o)$ et l'ensemble des $[\Gamma^s]$ telle que l'on ait les relations
suivantes:
$$\lfloor [\cdots,\overleftrightarrow{a}]_{\pi_o} \overrightarrow{\times} 
[\overleftrightarrow{b},\cdots]_{\pi_o}\rfloor
= [\cdots,\overleftrightarrow{a},\overrightarrow{1},\overleftrightarrow{b},\cdots]_{\pi_o}$$
$$\lceil [\cdots,\overleftrightarrow{a}]_{\pi_o} \overrightarrow{\times} 
[\overleftrightarrow{b},\cdots]_{\pi_o}\rceil
= [\cdots,\overleftrightarrow{a},\overleftarrow{1},\overleftrightarrow{b},\cdots]_{\pi_o}$$
$$\lfloor [\cdots,\overleftrightarrow{a}]_{\pi_o} \overleftarrow{\times} 
[\overleftrightarrow{b},\cdots]_{\pi_o}\rfloor
= [\cdots,\overleftrightarrow{b},\overleftarrow{1},\overleftrightarrow{a},\cdots]_{\pi_o}$$
$$\lceil [\cdots,\overleftrightarrow{a}]_{\pi_o} \overleftarrow{\times} [\overleftrightarrow{b},\cdots]_{\pi_o}\rceil
= [\cdots,\overleftrightarrow{b},\overrightarrow{1},\overleftrightarrow{a},\cdots]_{\pi_o}$$ où
$\overleftrightarrow{c}$ désigne arbitrairement $\overleftarrow{c}$ ou $\overrightarrow{c}$.

\end{itemize}
\end{prop-defi}

\begin{rema} \label{rema-numero} Dans \cite{ze}, $(\overleftarrow{a_1}, \cdots , \overrightarrow{a_r})$ est représenté 
sous la forme
$$\circ^{\frac{s-1}{2}} \lefto \circ \cdots \circ \lefto \cdots \circ \to \circ^{\frac{1-s}{2}}$$
où les $a_1$ premières flèches vont dans le sens de la flèche au dessus de $a_1$, les $a_2$ suivantes dans le
sens de la flèche au dessus de $a_2$ ... Ainsi graphiquement on a
$$[\Gamma^{s_1}_1] \overrightarrow{\times} [\Gamma^{s_2}_2]:\overbrace{\circ \lefto \circ \cdots \to
\circ}^{\Gamma^{s_1}_1} \longleftrightarrow \overbrace{\circ \to \circ \cdots \circ}^{\Gamma^{s_2}_2}$$

$$(\hbox{resp. } [\Gamma^{s_1}_1] \overleftarrow{\times} [\Gamma_2^{s_2}]:
\overbrace{\circ \lefto \circ \cdots \to \circ}^{\Gamma_2^{s_2}} \longleftrightarrow \overbrace{\circ \to
\circ \cdots \circ}^{\Gamma_1^{s_1}})$$
\end{rema}

\medskip

\noindent \textit{Exemples}: la représentation cuspidale irréductible $\pi_o$ de $GL_g(F_o)$ étant fixé, avec
nos notations, on a:
$$\begin{array}{rl}
{[}\overleftarrow{s-1}{]}_{\pi_o} &= \lfloor \overbrace{[\overrightarrow{0}]_{\pi_o} \overleftarrow{\times}
\cdots \overleftarrow{\times} [\overrightarrow{0}]_{\pi_o}}^{s} \rfloor = \lceil
\overbrace{[\overleftarrow{0}]_{\pi_o} \overrightarrow{\times} \cdots \overrightarrow{\times}
[\overleftarrow{0}]_{\pi_o}}^s
\rceil \\
{[}\overrightarrow{s-1}{]}_{\pi_o} &= \lfloor \overbrace{[\overleftarrow{0}]_{\pi_o} \overrightarrow{\times}
\cdots \overrightarrow{\times} [\overleftarrow{0}]_{\pi_o}}^s \rfloor =\lceil
\overbrace{[\overrightarrow{0}]_{\pi_o} \overleftarrow{\times} \cdots \overleftarrow{\times}
[\overrightarrow{0}]_{\pi_o}}^s
 \rceil \\
{[} \overleftarrow{l-1},\overrightarrow{s-l} {]}_{\pi_o} &= \lfloor [\overleftarrow{l-1}]_{\pi_o}
\overrightarrow{\times} [\overrightarrow{s-l-1}]_{\pi_o} \rfloor =\lceil [\overleftarrow{l-2}]_{\pi_o}
\overrightarrow{\times}
[\overrightarrow{s-l}]_{\pi_o} \rceil \\
&= \lceil [\overrightarrow{s-l-1}]_{\pi_o} \overleftarrow{\times} [\overleftarrow{l-1}]_{\pi_o} \rceil=
\lfloor [\overrightarrow{s-l}]_{\pi_o} \overleftarrow{\times} [\overleftarrow{l-2}]_{\pi_o} \rfloor
\end{array}$$

\noindent \textbf{Notation}: \textit{Dans la suite $[\overleftrightarrow{s-1}]_{\pi_o}$ désignera une
représentation elliptique quelconque de type $\pi_o$ de $GL_{sg}(F_o)$.}

\subsection{Foncteur de Jacquet}
\label{jacquet}

Soit $P=MN$ un parabolique de $GL_d$ de Lévi $M$ et de radical unipotent $N$.

\begin{defi}
Pour $\pi$ une représentation admissible de $GL_d(F_o)$, l'espace des vecteurs $N(F_o)$-coinvariants est
stable sous l'action de $M(F_o) \simeq P(F_o)/N(F_o)$. On notera $J_N(\pi)$ cette représentation tordue par
$\d_P^{-1/2}$.
\end{defi}

\begin{lemm} Soit $g$ un diviseur de $d=sg$ et $\pi_o$ une représentation irréductible cuspidale de
$GL_g(F_o)$. Pour $1 \leq h \leq d$, le foncteur de Jacquet vérifie les propriétés suivantes:
\begin{itemize}
\item si $g$ ne divise pas $h$, alors
$$J_{N_{h,d}^{op}}([\overleftarrow{s-1}]_{\pi_o})=J_{N_{h,d}^{op}}([\overrightarrow{s-1}]_{\pi_o})=
J_{N_{h,d}}([\overleftarrow{s-1}]_{\pi_o})=J_{N_{h,d}}([\overrightarrow{s-1}]_{\pi_o}) =(0)$$

\item si $h=lg$ alors
$$J_{N_{lg,sg}}([\overleftarrow{s-1}]_{\pi_o})=[\overleftarrow{l-1}]_{\pi_o((s-l)/2)} \otimes
[\overleftarrow{s-l-1}]_{\pi_o(-l/2)}$$
$$J_{N_{lg,sg}}([\overrightarrow{s-1}]_{\pi_o})=[\overrightarrow{l-1}]_{\pi_o((l-s)/2)} \otimes
[\overrightarrow{s-l-1}]_{\pi_o(l/2)}$$
$$J_{N_{lg,sg}^{op}}([\overleftarrow{s-1}]_{\pi_o})=[\overleftarrow{l-1}]_{\pi_o((l-s)/2)} \otimes
[\overleftarrow{s-l-1}]_{\pi_o(l/2)}$$
$$J_{N_{lg,sg}^{op}}([\overrightarrow{s-1}]_{\pi_o})=[\overrightarrow{l-1}]_{\pi_o((s-l)/2)} \otimes
[\overrightarrow{s-l-1}]_{\pi_o(-l/2)}$$
\end{itemize}
\end{lemm}

\begin{proof} soit $\Gamma=(a_i,\e_i)_{1 \leq i \leq r}$, les résultats découlent alors des propriétés suivantes
que l'on trouve dans \cite{ze}:

- $J_{N_{g,2g,\cdots,sg}^{op}}([\Gamma^s])$ est de la forme $\pi_o(\frac{1-s}{2}+\s(0)) \otimes \cdots
\pi_o(\frac{1-s}{2}+\s(s-1))$ où $\s$ est une permutation de l'ensemble $\{ 0,\cdots,s-1 \}$ soumise à la
règle suivante: soit $1 \leq i \leq r$ avec $\e_i=1$ (resp. $\e_i=-1$): pour tout $a_1+\cdots + a_{i-1} \leq
r < r' \leq a_1+ \cdots +a_i$ alors $\s^{-1}(r) < \s^{-1}(r')$ (resp. $\s^{-1}(r) > \s^{-1}(r')$).
\footnote{Autrement dit $\s$ est compatible aux orientations des flèches.}

- en ce qui concerne $J_{N_{g,2g,\cdots,sg}}([\Gamma^s])$ la règle est inversée, i.e. $\s^{-1}(r)>
\s^{-1}(r')$ (resp. $\s^{-1}(r) < \s^{-1}(r')$).

\end{proof}

\section{Preuve du cas Iwahori}

On rappelle le lemme suivant:

\begin{lemb} \label{lem-torsion-x}
Pour tout caractère $\psi:\Zm \longto \bar \Qm_l^\times$, on a un isomorphisme
$$\widetilde{\UC_{F_o}^{d,i}} (\t_o \otimes (\psi \circ \val \circ \det)) \simeq \widetilde{\UC_{F_o}^{d,i}}
(\t_o) \otimes (\psi \circ d_d)$$ où $d_d:(g,\s) \in GL_d(F_o) \times W_{F_o} \longto -\val(\det g)-\deg(\s)
\in \Zm$.
\end{lemb}

Le but est alors d'étendre le lemme précédent en la proposition suivante:

\begin{propb} \label{torsion2}
Pour tout caractère $\eta_o$ de $F_o^\times$, on a
$$\widetilde{\UC_{F_o}^{d,i}}(\t_o \otimes \eta_o) \simeq \widetilde{\UC_{F_o}^{d,i}}(\t_o)\otimes(\eta_o \otimes 
\eta_o
\circ \cl^{-1})$$ en tant que $(GL_d(F_o) \otimes W_o)$-module.
\end{propb}

\begin{proof} L'isomorphisme de l'énoncé découle d'un accouplement cohomologique que l'on va expliciter dans les
lignes qui suivent. On considère la catégorie des faisceaux en $\bar \Qm_l$-modules sur les espaces
considérés. Soit $\XF_n$ l'espace rigide analytique associé au schéma formel $\spf \Def_n^d$ sur $\hat
\OC_o^{nr}$ muni de la topologie étale de Berkovich (cf. \cite{ber1}); on note en particulier $\XF_s$ (resp.
$\XF_\eta$) sa fibre spéciale (resp. générique). Dans loc. cit., l'auteur construit le foncteur des cycles
évanescents $\Psi_\eta$ tel que pour tout faisceau $\FC$ sur $\XF_n$, $R^q \Psi_\eta(\FC)$ est le faisceau
associé au préfaisceau qui à une extension étale $\NF_s \to \XF_s$ associe le $\bar \Qm_l$-espace vectoriel
de dimension finie $H^q(\NF_\eta,\FC)$.

Ainsi pour tous faisceaux $\FC_1$ et $\FC_2$, on a une flèche:
$$T^0:R^0 \Psi_\eta(\FC_1) \otimes_{\bar \Qm_l} R^0 \Psi_\eta(\FC_2) \longto R^0 \Psi_\eta(\FC_1 \otimes_{\bar \Qm_l} 
\FC_2)$$
En outre étant donné un \ele de $\NC_o$, la correspondance $(c_1,c_2)$ associée fournit
$$c_1^* R^0 \Psi_\eta(\FC) \longto R^0 \Psi_\eta(c_1^* \FC) \longto R^0 \Psi_\eta(c_2^! \FC) \longto c_2^! R^0 
\Psi_\eta(\FC)$$
compatible avec $T^0$, \cad que le diagramme suivant est commutatif:
$$\diagram c_1^* R^0 \Psi_\eta(\FC_1) \otimes c_1^* R^0 \Psi_\eta(\FC_2) \rrto^{c_1^* T^0} \dto & &
c_1^* R^0 \Psi_\eta(\FC_1 \otimes \FC_2) \dto \\
c_2^! R^0 \Psi_\eta(\FC_1) \otimes c_2^! R^0 \Psi_\eta(\FC_2) \rrto^{c_2^! T^0} & & c_2^! R^0 \Psi_\eta(\FC_1
\otimes \FC_2)
\enddiagram$$

D'après, par exemple, \cite{bre} théorème II 6.2, on en déduit pour tout $i \geq 0$, des flèches
$$T^i: R^0 \Psi_\eta(\bar \Qm_l) \otimes_{\bar \Qm_l} R^i \Psi_\eta(\bar \Qm_l) \longto R^i \Psi_\eta(\bar \Qm_l)$$
qui, d'après la commutativité du diagramme ci-dessus et l'unicité des flèches $T^i$, sont compatibles à
l'action de $\NC_o$, i.e. $T^i (n s_1,n s_2)=n T^i (s_1,s_2)$ pour tout $n \in \NC_o$.

On fixe $s \in R^0 \Psi_\eta(\bar \Qm_l)$ et on considère l'application
$$T^i_s:R^i \Psi_\eta(\FC) \longto R^i \Psi_\eta(\FC)$$
définie par $T^i_s(t)=T^i(s,t)$. Pour tout $\FC$, on vérifie aisément que $T^0_s$ est un isomorphisme, de
sorte que d'après loc. cit. théorème 6.2 (c), qui n'est autre qu'une application du lemme des cinq, $T^i_s$
est un isomorphisme pour tout $\FC$. On note encore $T^i$ l'application
$$\widetilde{\UC_{F_o}^{d,0}} \otimes \widetilde{\UC_{F_o}^{d,i}} \longto \widetilde{\UC_{F_o}^{d,i}}$$
En notant $R^0 \Psi_\eta^1 (\bar \Qm_l)$ les cycles évanescents sur $\spf \Def_n^1$, et
$\widetilde{\UC_{F_o}^{1,0}}$ la représentation de $F_o^\times \times F_o^\times \times W_o$ associée, on
rappelle que l'on a
$$\widetilde{\UC_{F_o}^{1,0}}= \eta_o \otimes \eta_o \otimes \eta_o \circ \cl^{-1}$$
La restriction de $T^i$ à
$$\widetilde{\UC_{F_o}^{1,0}} \otimes \widetilde{\UC_{F_o}^{d,i}} \longto \widetilde{\UC_{F_o}^{d,i}}$$
est d'après ce qui précède un isomorphisme qui fournit l'isomorphisme de l'énoncé.

\end{proof}

\begin{propb} \label{model-st}
L'anneau $\Def_n^d$ est muni d'une action de $GL_d(\OC_o/\MC_o^n)$ et on définit
$$\Def_{st}^d:=(\Def_n^d)^{\Iw_o/K_{o,n}}$$
où $\Iw_o$ est le sous-groupe d'Iwahori standard, i.e. l'ensemble
des matrices de $GL_d(\OC_o)$ triangulaires supérieures modulo $\MC_o$. L'anneau ainsi défini ne dépend pas
de l'entier $n >0$ choisi, il est régulier et semi-stable sur $\hat \OC_o^{nr}$, i.e. il existe des
coordonnées $\a_1,\cdots,\a_d$ tels que $\Def_{st}^d \simeq \bar \Fm_p [[\a_1,\cdots,\a_d]]$ et l'image de
$\pi_o$ par le morphisme structural $\hat \OC_o^{nr} \simeq \bar \Fm_p [[\pi_o]] \longto \Def_{st}^d$ est
égale à $\prod_{i=1}^d \a_i$.
\end{propb}

\begin{proof} L'indépendance de la définition de $\Def_{st}^d$ relativement à $d$, découle simplement du fait que
$K_{o,n}$ est un sous-groupe distingué de $\Iw_o$ pour tout $n >0$. Ainsi on a $\Def_{st}^d =
(\Def_1^d)^{\Iw_o/K_{o,1}}$ où $\Iw_o/K_{o,1}$ est isomorphe au sous-groupe des matrices triangulaires
supérieures de $GL_d(\k(o))$. On rappelle en outre que $\Def_0^d \hookrightarrow \Def_1^d$ représente le
foncteur des structures de niveau $1$ sur le $\OC_o$-module formel universel sur $\Def_0^d=\bar
\Fm_p[[a_0,\cdots,a_{d-1}]]$ à savoir $\t^d+a_{d-1} \t^{d-1} + \cdots + a_0$ où $\t=x^q$; la structure de
niveau universelle sur $\Def_1^d$ est donnée par une application $\Fm_q$-linéaire $\iota_{o,1}:\Fm_q^d
\longto \Def_1^d$ avec $\iota_{o,1}(e_i)=\a_i^1 \in \Def_1^d$ où $(e_i)_{1 \leq i \leq d}$ est la base
canonique de $\Fm_q^d$ qui vérifient la condition de Drinfeld
$$\t^d+a_{d-1} \t^{d-1} + \cdots + a_0=\prod_{(\l_i)_{1 \leq i \leq d} \in \k(o)^d} (X-\sum_{i=1}^d \l_i \a_i^1)$$
Par ailleurs, on a $\Def_1^d \simeq \bar \Fm_p[[\a_1^1,\cdots,\a_d^1]]$, l'action de $M \in GL_d(\k(o))$
étant donnée par la multiplication à droite de $M^{-1}$ sur $\Fm_q^d$. L'orbite sous l'action du Borel
standard de $GL_d(\Fm_q)$ de la base canonique de $\Fm_q^d$ est le drapeau complet
$$(0) \subset \Vect(e_1) \subset \Vect(e_1,e_2) \subset \cdots \subset \Vect(e_1,\cdots,e_d)=\Fm_q^d$$
Ainsi l'inclusion $\Def_0^d \hookrightarrow \Def_{st}^d$ représente le foncteur des ``drapeaux complets'' du
noyau $\Sigma_d[\pi_o]$ de la multiplication par $\pi_o$ du $\OC_o$-module formel universel sur $\Def_0^d$,
i.e.
$$(0) \subset G_1 \subset G_2 \subset \cdots \subset G_d=\Sigma[\pi_o]$$
où les $G_i$ sont des sous-$\Fm_q$-modules de $\Sigma[\pi_o]$ de rang $i$. En posant
$$P_i(X)=\prod_{\sum_i \l_i \a_i^1 \in G_i} (X-\sum_i \l_i \a_i^1),$$
dans $\Def_1^d [X]$, on a $P_i(X)=(X^q+\a_i) \circ P_{i-1}(X)$ et donc $(\t+\a_d) \circ \cdots \circ
(\t+\a_1)=\t^d+a_{d-1} \t^{d-1}+\cdots + a_0$. En outre la donnée des $\a_i$ définit complètement les $G_i$
de sorte que $\Def_{st}^d \simeq \bar \Fm_p[[\a_1,\cdots,\a_d]]$ avec en particulier $a_0=\pi_o=\prod_{i=1}^d
\a_i$.

\end{proof}

\begin{corob} \label{drapeaux}
 L'inclusion $\Def_0^d \hookrightarrow \Def_{st}^d$ représente le foncteur des ``drapeaux complets'' du noyau
$\Sigma_d[\pi_o]$ de la multiplication par $\pi_o$ du $\OC_o$-module formel sur $\Def_0^d$, i.e. la donnée de
sous-$\Fm_q$-modules de rang $i$, $G_i$ telles que $(0) \subset G_1 \subset G_2 \subset \cdots \subset
G_d=\Sigma_d[\pi_o]$. Il existe alors des indéterminées $\a_i$ pour $1 \leq i \leq d$ tels que $\Def_{st}^d
\simeq \bar \Fm_p[[\a_1,\cdots,\a_d]]$ avec
$$(\t+\a_d) \circ \cdots \circ (\t+\a_1)= \t^{d}+a_{d-1} \t^{d-1}+ \cdots +a_0.$$
\end{corob}

\begin{theob} \label{theo-iwahori}
Pour tout $0 \leq i <d$, on a
$$\widetilde{\UC_{F_o}^{d,i}}(1_o) \simeq [\overleftarrow{i},\overrightarrow{d-i-1}]_{1_o} \otimes |\cl|^{-i}$$
et pour $g \neq 1$ un diviseur de $d=sg$ et $\pi_o$ une représentation cuspidale de $GL_g(F_o)$, on a
$$(\widetilde{\UC_{F_o}^{d,i}}(\JL^{-1}([\overleftarrow{s-1}]_{\pi_o}))^{\Iw_o}=(0)$$
où $\Iw_o$ désigne l'Iwahori standard de $GL_d(\OC_o)$.
\end{theob}

\begin{proof} On raisonne par récurrence sur $d$; on suppose donc le résultat acquis pour tout $h \leq d$ ce qui est
vérifié pour $h=1$. Commençons par rappeler un résultat sur les cycles évanescents sur un schéma semi-stable.

\begin{theob} \label{ill}
(cf. \cite{ill}) Soit $X \longto S=\spec \OC_o$ un schéma semi-stable, i.e. tel que localement pour la
topologie étale $\OC_X$ est de la forme $\OC_o[t_1,\cdots,t_d]/(\pi_o- \prod_{i=1}^d t_i)$; en notant comme
d'habitude $\Psi_{\eta_o}$ le foncteur des cycles proches, on a
$$R^i \Psi_{\eta_o} (\bar \Qm_l) \simeq \Lambda^i R^1 \Psi_{\eta_o}(\bar \Qm_l)$$
$$R^1 \Psi_{\eta_o} (\bar \Qm_l) \simeq \bigoplus_{i=1}^d (\bar \Qm_l)_{\bar X_{s,i}}/(\bar \Qm_l \textrm{ 
diagonal})(-1)$$
avec $X_{i,s}$ localement défini par $t_i=0$ dans la fibre spéciale $X_s$ de $X$.
\end{theob}

\begin{corob} \label{dim-poids}
Pour tout $0 \leq i < d$ et toute représentation irréductible admissible $\t_o$ de $D_{o,d}^\times$,
$(\Psi_{F_o}^{d,i})^{\Iw_o}$ est un $\bar \Qm_l$-espace vectoriel de dimension le coefficient binomial
$\genfrac{(}{)}{0pt}{}{i}{d-1}$ qui est pur de poids $2i$.
\end{corob}

\begin{proof} En effet avec les notations du théorème précédent, en tout point géométrique $z$ de $X_s$, les germes
de $R^i \Psi_{\eta_o}(\bar \Qm_l)_z$ sont purs de poids $2i$, poids indépendant du relèvement choisi pour le
Frobenius. On applique le théorème avec
$$X=\spec (\bar \Fm_p[\a_1,\cdots,\a_d]) \longto \spec \bar \Fm_p [\pi_o],$$
l'image de $\pi_o$ étant donnée par le produit des $\a_i$. Le théorème de comparaison de Berkovich donne
alors que le $i$-ème faisceau des cycles évanescents $\Psi_{F_o,st}^{d,i}$ de $\spf \Def_{st}^d$ est pur de
poids $2i$. En outre en notant $\pi_{n \to st}:\spf \Def_n^d \longto \spf \Def_{st}^d$, le faisceau constant
$\bar \Qm_{l,st,\eta_o}$ sur la fibre générique de l'espace analytique $\spf \Def_{st}^d$ est isomorphe à
$\pi_{n \to st,\eta_o,*} \bar \Qm_{l,n,\eta_o}$ avec des notations évidentes. On a
$$R^i \Psi_{st,\eta_o}(\pi_{n \longto st,\eta_o,*} \bar \Qm_l) \simeq \pi_{n \to st,s_o,*} R^i \Psi_{n,\eta_o}
(\bar \Qm_{l,n,\eta_o})$$ soit donc
$$R^i \Psi_{st,\eta_o} (\bar \Qm_{l,st,\eta_o}) \simeq \pi_{n \to st,s,*} R^i \Psi_{n,\eta_o}(\bar
\Qm_{l,n,\eta_o})^{\Iw_o}$$ Or comme les fibres spéciales sont réduites à un point, on obtient que
$\Psi_{F_o,st}^{d,i} \simeq (\Psi_{F_o}^{d,i})^{\Iw_o}$, d'où le résultat.

\end{proof}

\begin{lemb} (cf. \cite{ze}) \label{inv-iwa}
L'espace des invariants sous l'Iwahori standard de
$$[\overleftarrow{h-1},\overrightarrow{d-h}]_{1_o}$$
est de dimension le coefficient binomial $\binom{h-1}{d-1}$.
\end{lemb}

\begin{proof} D'après un argument classique (cf. par exemple \cite{boy} lemme 4.7), cette dimension est égale au
nombre de fois que $J_{N_d} ([\overleftarrow{h},\overrightarrow{d-h-1}]_{1_o})$ contient la représentation
triviale du tore maximal de $GL_d(F_o)$, où $N_d$ est le sous-groupe unipotent maximal du Borel standard et
$J$ est le foncteur de Jacquet. Or d'après \cite{ze}, $J_{N_d}
([\overleftarrow{h},\overrightarrow{d-h-1}]_{1_o})$ est la représentation triviale avec la multiplicité égale
au cardinal de $\Lambda^+$, le sous-ensemble de l'ensemble des permutations, cf. la remarque
(\ref{rema-numero}), de $\{ 1,\cdots ,d \}$ telles que $\l(i) < \l(i+1)$ pour tout $1 \leq i \leq h$ et
$\l(i)> \l(i+1)$ pour tout $h \leq i \leq d$. Le cardinal de $\Lambda^+$ est alors $\binom{h-1}{d-1}$; en
effet une fois choisi $h-1$ entiers distincts entre $2$ et $d$, on les classe par ordre décroissant, $\a_1
\geq \cdots \geq \a_{h-1}$; on classe de même par ordre croissant le complémentaire $\b_1 \leq \cdots \leq
\b_{d-h}$ et on pose $\l(i)=\a_i$ pour $1 \leq i <h$, $\l(h)=1$, $\l(h+i)=\b_i$ pour $1 \leq i \leq d-h$. On
vérifie aisément que l'on définit bien ainsi une bijection sur $\Lambda^+$.

\end{proof}

\marque On commence dans un premier temps par traquer les sous-quotients irréductibles de
$$n-\ind_{B^{op}(F_o)}^{GL_d(F_o)} (|-|^{(d-1)/2} \times \cdots \times |-|^{(1-d)/2})$$
i.e. les représentations elliptiques de type $1_o$. On raisonne par récurrence, le cas $d=1$ correspondant à
la théorie de Lubin-Tate. \footnote{Le cas $d=2$ est prouvé dans \cite{ca}.}

\begin{propb} \label{prop-segment}
Les représentations elliptiques de type $1_o$ de
$$\lim_{\genfrac..{0pt}{1}{\to}{I}} H^j_c(M_{I,\bar s_o}^{=h},R^i \Psi_{\eta_o}(\bar 
\Qm_l))=\ind_{P_{h;d}^{op}(F_o)}^{GL_d(F_o)}
\lim_{\genfrac..{0pt}{1}{\to}{I}} H^j_c(M_{I,\bar s_o,1}^{=h},R^i \Psi_{\eta_o}(\bar \Qm_l))$$ sont, avec les
notations du \S \ref{defi-induite}, de la forme
$[\overleftrightarrow{t},\overleftarrow{i},\overrightarrow{h-i-1},\overleftrightarrow{d-h-t}]_{1_o} \otimes
|\cl|^{-i-t}$ avec $0 \leq t \leq d-h$; plus précisément pour $t=0$ (resp. $t=d-h$, resp. $0 < t < d-h$) on
obtient
\begin{multline*}
[\overleftarrow{i},\overrightarrow{h-i-1}]_{1_o} \overrightarrow{\times}
[\overleftrightarrow{d-h-1}]_{1_o} \\
(resp.~ [\overleftarrow{i},\overrightarrow{h-i-1}]_{1_o} \overleftarrow{\times}
[\overleftrightarrow{d-h-1}]_{1_o}, \\
resp.~ ([\overleftarrow{i},\overrightarrow{h-i-1}]_{1_o} \overrightarrow{\times}
[\overleftrightarrow{d-h-t-1}]_{1_o}) \overleftarrow{\times} [\overleftrightarrow{t-1}]_{1_o})
\end{multline*}
où par exemple l'unique quotient est de la forme
\begin{multline*}
[\overleftarrow{i},\overrightarrow{h-i},\overleftrightarrow{d-h-1}]_{1_o} \\
(resp. ~
[\overleftrightarrow{d-h-1},\overrightarrow{1},\overleftarrow{i},\overrightarrow{h-i-1}]_{1_o}, \\
resp. ~ [\overleftrightarrow{t-1},\overrightarrow{1},\overleftarrow{i},\overrightarrow{h-i},
\overleftrightarrow{d-h-t-1}]_{1_o}).
\end{multline*}
Par cela on entend tout d'abord que dans le groupe de Grothendieck correspondant, les 2 ou 4 représentations
en question apparaissent simultanément; ensuite de manière plus précise, les induites de l'énoncé
apparaissent comme sous-espace et comme quotient (quitte à changer les flèches non précisées).
\end{propb}

\begin{proof} Elle est similaire à celle de la proposition (\ref{prop-poids}). Soit $\xi_o$ un caractère (non
ramifié) de $F_o^\times$. D'après l'hypothèse de récurrence on a
$$\widetilde{\UC_{F_o}^{h,i}}(\xi_o)\simeq [\overleftarrow{i},\overrightarrow{h-i-1}]_{\xi_o} \otimes \xi_o(-i)$$
en tant que représentation de $GL_h(F_o) \times W_o$. Les représentations elliptiques de type $1_o$
s'obtiennent exclusivement comme sous-quotient des induites normalisées
$$n-\ind_{B^{op}(F_o)}^{GL_d(F_o)} |-|^{i_1} \times \cdots \times |-|^{i_d}$$
avec $\{ i_1,\cdots ,i_d \}=\{ \frac{1-d}{2},\cdots , \frac{d-1}{2} \}$ ou
encore avec des induites non normalisées comme les sous-quotients de
$$\ind_{B^{op}(F_o)}^{GL_d(F_o)} |-|^{i_1+(1-d)/2} \times \cdots \times |-|^{i_d + (d-1)/2}.$$
Dans le groupe de Grothendieck des représentations admissibles de $GL_{d-h}(F_o) \times \Zm$, on découpe
${\DS \lim_{\genfrac..{0pt}{1}{\to}{I}}} H^j_c(M_{I,s_o}^{=h},\FC_{\xi_o})$ suivant les caractères $\chi_o$
de $\Zm$ et ensuite suivant les supports cuspidaux pour $GL_{d-h}(F_o)$ ce qui s'écrit $\sum_{\eta_o,\chi_o}
\eta_o \otimes \chi_o$ où $\eta_o$ (resp. $\chi_o$) décrit les $\bar \Qm_l$-représentations irréductibles de
$GL_{d-h}(F_o)$ (resp. $\Zm$). On en déduit alors que pour $0 \leq i \leq h-1$, ${\DS
\lim_{\genfrac..{0pt}{1}{\to}{I}}} H^j_c(M_{I,s_o}^{=h},R^i\Psi_{\eta_o}(\bar \Qm_l))$ est de la forme
$$\sum_{\xi_o,\eta_o,\chi_o} \ind_{P_{h,d}^{op}(F_o)}^{GL_d(F_o)} ([\overleftarrow{i},\overrightarrow{h-i-1}]_{\xi_o 
\otimes \chi_o}
\otimes \eta_o) \otimes (\xi_o \chi_o)(-i)$$
En ce qui concerne les représentations elliptiques de type $1_o$, on considère
les $\chi_o$ de la forme $\xi_o^{-1}(-t)$, ce qui donne avec les notations ci-dessus 
$i_1=\frac{d-1}{2}-t,\cdots,i_h=\frac{d-1}{2}-h+1-t$,
d'où le résultat.

Si on ne veut plus simplement raisonner dans le groupe de Grothendieck, on choisit $\eta_o \otimes \chi_o$ un
sous-espace (resp. un quotient) irréductible du groupe de cohomologie précédent. Le reste du raisonnement est
alors identique.

\end{proof}

\marque On globalise $\st_\oo \otimes \st_o$ en une représentation automorphe $\Pi$ de $D_\Am^\times$ telle
qu'il existe deux places distinctes de $\oo,o$ et des places de ramification $\bad$ de $D$, pour lesquelles
la composante locale de $\Pi$ est cuspidale; pour l'existence d'une telle globalisation cf. par exemple
\cite{badu}. D'après le théorème (14.12) de \cite{lrs}, la composante isotypique $H^i_{\eta_o}[\Pi^\oo]$ est
nulle pour $i \neq d-1$ et pour $i=d-1$ elle est isomorphe, en tant que représentation du groupe de
Weil-Deligne local en $o$, à $\sp(d) \otimes |\cl|^{(1-d)/2}$ où $\sp_d$ est la représentation spéciale i.e.
$\sp(d)=\bar \Qm_l((1-d)/2) \oplus \cdots \oplus \bar \Qm_l ((d-1)/2)$ en tant que représentation de $W_o$
avec l'action de la monodromie donnée par $N:\bar \Qm_l(k/2) \simeq \bar \Qm_l(k/2+1)$. \footnote{Dans le
contexte de \cite{h-t}, ce résultat est prouvé dans \cite{y-t}.}

\begin{propb} Pour tout $0 \leq i <d$, $[\overleftarrow{d-1}]_{1_o} \otimes |\cl|^{-i}$ est un sous-quotient de
$$[\lim_{\genfrac..{0pt}{1}{\to}{I}} H^i(M_{I,\bar s_o},R^{d-1-i} \Psi_{\eta_o}(\bar \Qm_l))][\Pi^\oo].$$
\end{propb}

\begin{proof}  On utilise simplement que pour un schéma à réduction semi-stable, la filtration aboutissement de la
suite spectrale des cycles évanescents coïncide avec la filtration des noyaux (cf. \cite{ill}); le résultat
découle alors directement de la description de la cohomologie générique donnée dans \cite{lrs}.

\end{proof}

\begin{lemb} Soit $h_0$ le plus grand (s'il existe) $0< h \leq d$ tel qu'il existe $t,j$ pour lesquels
$$[\overleftrightarrow{t},\overleftarrow{i},\overrightarrow{d-h-i-1},\overleftrightarrow{h-t}]_{1_o} \otimes
|\cl|^{-i-t}$$ est un sous-quotient de ${\DS \lim_{\genfrac..{0pt}{1}{\to}{I}} H^j_c (M_{I,\bar
s_o}^{=d-h},R^i \Psi_{\eta_o}(\bar \Qm_l)) [\Pi^\oo]}$. Alors pour tout $h'>h_0$,
$$[\overleftrightarrow{t},\overleftarrow{i},\overrightarrow{d-h_0-i-1},\overleftrightarrow{h_0-t}]_{1_o}
\otimes |\cl|^{-i-t}$$ n'est pas un sous-quotient de ${\DS \lim_{\genfrac..{0pt}{1}{\to}{I}}
H^{j+1}_c(M_{I,\bar s_o}^{=d-h'},R^i \Psi_{\eta_o}(\bar \Qm_l)) [\Pi^\oo]}$.
\end{lemb}

\begin{proof} Si $\pi_{1_o} \otimes |\cl|^{-i-t}$ est un sous-quotient de
$$\lim_{\genfrac..{0pt}{1}{\to}{I}} H^{j+1}_c(M_{I,\bar s_o}^{=d-h'},R^i \Psi_{\eta_o}(\bar \Qm_l)) [\Pi^\oo]$$
pour $\pi_{1_o}$ une représentation elliptique de type $1_o$, pour $h'>h_0$ alors pour des raisons de poids
$\pi_{1_o}$ est, d'après la proposition (\ref{prop-segment}) de la forme
$[\overleftrightarrow{t},\overleftarrow{i},\overrightarrow{d-h'-i-1},\overleftrightarrow{h'-t}]_{1_o}$. On
conclut alors par la maximalité de $h_0$.

\end{proof}

\begin{lemb} Pour tout $h'<h$ et tout $0 \leq t <h$ (resp. $t=h$),
$$[\overleftrightarrow{t},\overleftarrow{i},\overrightarrow{d-h-i-1},\overleftarrow{\a},
\overleftrightarrow{h-t-\a}]_{1_o} \otimes |\cl|^{-i-t},$$
$$(\hbox{resp. }
[\overleftrightarrow{h},\overleftarrow{i},\overrightarrow{d-h-i-1}]_{1_o} \otimes |\cl|^{-h-i}),$$ pour $0 <
\a \leq h-t$, n'est pas un sous-quotient de
$${\DS \lim_{\genfrac..{0pt}{1}{\to}{I}} H^{j-1}_c(M_{I,\bar s_o}^{=d-h'},R^i \Psi_{\eta_o}(\bar \Qm_l)) 
[\Pi^\oo]}.$$
\end{lemb}

\begin{proof} En effet d'après la proposition (\ref{prop-segment}), les représentations elliptiques de type $1_o$ de
poids $2(i+t)$ de ${\DS \lim_{\genfrac..{0pt}{1}{\to}{I}} H^{j-1}_c(M_{I,\bar s_o}^{=d-h'},R^i
\Psi_{\eta_o}(\bar \Qm_l)) [\Pi^\oo]}$, pour $0 \leq t \leq h'$, sont de la forme
$$[\overleftrightarrow{t},\overleftarrow{i},\overrightarrow{d-h'-i-1},\overleftrightarrow{h'-t}]_{1_o}$$
de sorte qu'en numérotant comme d'habitude les sommets du graphe $\vec \Gamma^d$ de $\frac{d-1}{2}$ à
$\frac{1-d}{2}$, on a toujours
$$\cdots \bullet^{\frac{d-1-2(t+d-h-1)}{2}} \longto \bullet^{\frac{d-1-2(t+d-h)}{2}} \cdots$$
alors que dans $[\overleftrightarrow{t},\overleftarrow{i},\overrightarrow{d-h-i-1},\overleftarrow{\a},
\overleftrightarrow{h-t-\a}]_{1_o}$ cette flèche est orientée dans l'autre sens car $\a>0$ (pour $t=h$ un
simple argument de poids suffit).

\end{proof}

La proposition suivante montre que dans (\ref{prop-segment}), on doit avoir $t=0$.

\begin{propb} \label{prop-t}
Pour $i>0$ fixé et pour tout $j$, $t>0$ et $h \neq 0$,
$$[\overleftrightarrow{t},\overleftarrow{i},\overrightarrow{d-h-i-1},\overleftrightarrow{h-t}]_{1_o}$$
n'est pas un sous-quotient de ${\DS \lim_{\genfrac..{0pt}{1}{\to}{I}}} H^j_c (M_{I,\bar s_o}^{=d-h},R^i
\Psi_{\eta_o}(\bar \Qm_l)) [\Pi^\oo]$.
\end{propb}

\begin{proof} On raisonne par l'absurde; soit donc $h_0$ le plus grand, l'idée est de montrer que cette 
représentation
se retrouverait, à travers la suite spectrale de stratification et celle des cycles évanescents dans la
cohomologie de la fibre générique ce qui n'est pas d'après \cite{lrs}.

La suite spectrale associée à la stratification de la fibre spéciale $M_{I,s_o}$:
$$E_1(I,i)^{p,q}=H^{p+q}_c(M_{I,\bar s_o}^{=d-p},R^i \Psi_{\eta_o}(\bar \Qm_l))
\Rightarrow H^{p+q}(M_{I,\bar s_o},R^i \Psi_{\eta_o}(\bar \Qm_l)),$$ le fait que $H^0(M_{I,s_o}^d,R^i
\Psi_{\eta_o}(\bar \Qm_l))$ soit de poids $2i$, la proposition (\ref{prop-segment}), et les deux lemmes
précédents, montrent que pour $0 < t < h_0$ (resp. $t=h_0$) et tout $r \geq 1$,
$$\pi_o:=
[\overleftrightarrow{t-1},\overrightarrow{1},\overleftarrow{i},\overrightarrow{d-h_0-i-1},\overleftarrow{1},
\overleftrightarrow{h_0-1}]_{1_o} \otimes |\cl|^{-i-t}$$
$$(\hbox{resp. } \pi_o:=[\overleftrightarrow{h_0-1},
\overrightarrow{1},\overleftarrow{i},\overrightarrow{d-h_0-i-1}]_{1_o} \otimes |\cl|^{-i-t})$$ est un
sous-quotient de ${\DS (\lim_{\genfrac..{0pt}{1}{\to}{I}} E_r^{h_0,j-h_0}(I,i))[\Pi^\oo]}$. \footnote{En fait
c'est un quotient.}

Montrons que pour tout $k \geq 2$, il n'existe pas de $j$ et $h$ tels que $\pi_o \otimes |\cl|^{-i-t}$ soit
un sous-quotient de ${\DS \lim_{\genfrac..{0pt}{1}{\to}{I}}} H^{j+k}_c (M_{I,\bar s_o}^{=d-h},R^{i+1-k}
\Psi_{\eta_o}(\bar \Qm_l)) [\Pi^\oo]$. D'après la proposition (\ref{prop-segment}), une représentation
elliptique de type $1_o$ de cette dernière représentation de poids $2(i+t)$, est de la forme
$$[\overleftrightarrow{t+k-1},\overleftarrow{i-k+1},\overrightarrow{d-h-i+k-2},
\overleftrightarrow{h-t-k+1}]_{1_o}$$ en particulier on doit avoir $d-h-i+k-2 \leq d-h_0-i-1$ soit $h \geq
h_0+k-1$ ce qui contredit la maximalité de $h_0$.

De même pour tout $k \geq 2$, il n'existe pas de $j,h$ tels que $\pi_o \otimes |\cl|^{-i-t}$ soit un
sous-quotient de ${\DS \lim_{\genfrac..{0pt}{1}{\to}{I}}} H^{j-k}_c (M_{I,\bar s_o}^{=d-h},R^{i+k-1}
\Psi_{\eta_o}(\bar \Qm_l)) [\Pi^\oo]$. En effet pour $h>0$, un tel sous-quotient elliptique de type $1_o$ de
poids $2(i+t)$ de cet espace, est de la forme
$$[\overleftrightarrow{t-k+1},\overleftarrow{i+k-1},\overrightarrow{d-h-i+k-2},
\overleftrightarrow{h-t+k-1}]_{1_o},$$ de sorte que $i+k-1 \leq i$, contradiction. Pour $h=0$, pour des
raisons de poids, il faut $t=j-1$; le résultat découle alors du lemme suivant.

\begin{lemb} La dimension des invariants sous $\Iw_o$ de
$$[\overleftrightarrow{j-2},\overrightarrow{1},\overleftarrow{i},\overrightarrow{d-h-i-1},
\overleftrightarrow{h-j+1}]_{1_o}$$ est strictement supérieure à $\binom{i+j-2}{d-1} \binom{i-1}{i+j-2}$.
\end{lemb}

\begin{proof} On procède comme rappelé au lemme (\ref{inv-iwa}). Pour donner une numérotation qui induise
l'orientation donnée, on peut commencer par donner le numéro $1$ au sommet $(i+j-1)$-ème sommet, puis on
choisit $i+j-2$ parmi $d-1$ qui serviront à numéroter les $i+j-2$ premiers sommets. Parmi ces $i+j-2$, on en
choisit à nouveau $i-1$ pour numéroter les sommets de $j-1$ à $i+j-2$. On peut aussi numéroter avec $1$ un
des premiers sommets dans une configuration $\leftarrow \bullet \to $, il y en a forcément un parmi les
$i+j-2$ premiers, d'où l'inégalité stricte.

\end{proof}

Ainsi $\pi_o \otimes |\cl|^{-i-t}$ apparaît dans l'aboutissement de la suite spectrale des cycles
évanescents, ce qui n'est pas d'après \cite{lrs}.

\end{proof}

\begin{corob} \label{corob-poids}
Tous les ${\DS \lim_{\genfrac..{0pt}{1}{\to}{I}}} H^j_c (M_{I,\bar s_o}^{=d-h},R^i \Psi_{\eta_o}(\bar \Qm_l))
[\Pi^\oo]$ sont purs de poids $2i$.
\end{corob}

\begin{lemb} Si $[\overleftarrow{i},\overrightarrow{d-h-i},\overleftarrow{h-1}]_{1_o} \otimes |\cl|^{-i}$
est un sous-quotient de
$$\lim_{\genfrac..{0pt}{1}{\to}{I}} H^{j}_c(M_{I,\bar s_o}^{=d-h'},R^i \Psi_{\eta_o}(\bar \Qm_l))[\Pi^\oo]$$
alors $h'=h$ ou $h-1$.
\end{lemb}

\begin{proof} A nouveau si $[\overleftarrow{i},\overrightarrow{d-h-i},\overleftarrow{h-1}]_{1_o} \otimes |\cl|^{-i}$
est un sous-quotient de l'espace en question pour $h'>h$, alors d'après (\ref{prop-segment})
$$[\overleftarrow{i},\overrightarrow{d-h'-i-1},\overleftarrow{1},\overrightarrow{h'-h},
\overleftarrow{h-1}]_{1_o} \otimes |\cl|^{-i}$$ aussi. Soit alors $h_0$ le plus grand $h'$ tel que
$\lim_{\genfrac..{0pt}{1}{\to}{I}} H^{j}_c(M_{I,\bar s_o}^{=d-h'},R^i \Psi_{\eta_o}(\bar \Qm_l))[\Pi^\oo]$
contienne
$$[\overleftarrow{i},\overrightarrow{d-h'-i-1},\overleftarrow{\a_1},\overrightarrow{\a_2},
\overleftrightarrow{h-\a_1-\a_2}]_{1_o} \otimes |\cl|^{-i}$$ avec $\a_1$ et $\a_2$ strictement positifs.
Celle-ci se retrouve alors par maximalité de $h_0$ et par un lemme analogue au précédent, dans
l'aboutissement de la suite spectrale de stratification à savoir: $\lim_{\genfrac..{0pt}{1}{\to}{I}}
H^j(M_{I,\bar s_o},R^i \Psi_{\eta_o}(\bar \Qm_l))[\Pi^\oo]$. Pour des raisons de poids, elle se retrouve
aussi dans l'aboutissement de la suite spectrale des cycles évanescents; en effet elle ne peut pas être
compensée par une contribution des points supersinguliers en vertu de (\ref{dim-poids}), ni d'après
(\ref{prop-segment}) par des $R^{i'} \Psi$ pour $i' > i$ et ni d'après (\ref{prop-t}) pour $i'<i$. On en
déduit donc $h' \leq h$.

Si $h'<h-1$, d'après (\ref{prop-segment})
$[\overleftarrow{i},\overrightarrow{d-h-i-1},\overleftarrow{h-1}]_{1_o}$ ne peut pas être un sous-quotient de
${\DS \lim_{\genfrac..{0pt}{1}{\to}{I}} H^{j}_c(M_{I,\bar s_o}^{=d-h'},R^i \Psi_{\eta_o}(\bar
\Qm_l))[\Pi^\oo]}$ car $d-h'-i-1>d-h-i$; d'où le résultat.

\end{proof}

\begin{propb} \label{prop-iwa}
Pour tout $0<h \leq d-i-1$,
$${\DS \lim_{\genfrac..{0pt}{1}{\to}{I}} H^h_c(M_{I,\bar s_o}^{=d-h},R^i \Psi_{\eta_o}(\bar \Qm_l))[\Pi^\oo]}$$
admet $[\overleftarrow{i},\overrightarrow{d-h-i},\overleftarrow{h-1}]_{1_o} \otimes |\cl|^{-i}$ comme
sous-quotient.
\end{propb}

\begin{proof} On raisonne par récurrence descendante sur $h$; pour $h=d-i-1$, on utilise le fait que ${\DS 
\lim_{\genfrac..{0pt}{1}{\to}{I}}} H^{d-i-1} (M_{I,\bar s_o},R^i \Psi_{\eta_o} (\bar \Qm_l)) [\Pi^\oo]$ admet $\st_o 
\otimes
|\cl|^{-i}$ comme sous-quotient. On rappelle que d'après (\ref{corob-poids})
$$\lim_{\genfrac..{0pt}{1}{\to}{I}} H^{d-i-1}_c(M_{I,\bar s_o}^{=i+1},R^i \Psi_{\eta_o}(\bar \Qm_l))[\Pi^\oo]$$
se calcule comme l'induite $\ind_{P_{i+1}(F_o)}^{GL_d(F_o)} [\overleftarrow{i}]_{1_o} \otimes \pi'$ pour une
certaine représentation admissible $\pi'$ de $GL_{d-i-1}(F_o)$ et que, d'après l'hypothèse de récurrence sur
les $\Psi_{F_o}^{h,i}(1_o)$, pour tout $h \neq i+1$ la partie de poids $2i$ de
$$\lim_{\genfrac..{0pt}{1}{\to}{I}} H^{d-i-1}_c(M_{I,\bar s_o}^{=i+1},R^i \Psi_{\eta_o}(\bar \Qm_l))[\Pi^\oo]$$
est soit nulle pour $h \leq i$, soit, pour $h > i+1$ de la forme
$[\overleftarrow{i},\overrightarrow{h-1-i},\overleftrightarrow{d-h}]_{1_o} \otimes |\cl|^{-i}$. La suite
spectrale associée à la stratification impose alors que $\pi'$ admet la représentation de Steinberg
$[\overleftarrow{d-i-2}]_{1_o}$ comme sous-quotient et donc ${\DS \lim_{\genfrac..{0pt}{1}{\to}{I}}
H^{d-i-1}_c(M_{I,\bar s_o}^{=i+1},R^i \Psi_{\eta_o}(\bar \Qm_l))[\Pi^\oo]}$ admet
$$[\overleftarrow{i},\overrightarrow{1},\overleftarrow{d-i-2}]_{1_o} \otimes |\cl|^{-i}$$
comme sous-quotient.

Supposons donc la propriété vérifiée pour $1<h <d-i$, i.e.
$$\lim_{\genfrac..{0pt}{1}{\to}{I}} H^h_c(M_{I,\bar s_o}^{=d-h},R^i \Psi_{\eta_o}(\bar \Qm_l))[\Pi^{\oo}]$$
admet $[\overleftarrow{i},\overrightarrow{d-h-i},\overleftarrow{h-1}]_{1_o} \otimes |\cl|^{-i}$ comme
sous-quotient. En outre on sait cette dernière n'apparaît pas dans l'aboutissement de la suite spectrale des
cycles évanescents de sorte qu'il existe $j,i',h'$ tels qu'elle soit un sous-quotient de
$$\lim_{\genfrac..{0pt}{1}{\to}{I}} H^j_c(M_{I,\bar s_o}^{=d-h'},R^{i'} \Psi_{\eta_o}(\bar \Qm_l))[\Pi^{\oo}]$$
Pour des raisons de poids et d'après le corollaire (\ref{corob-poids}), il faut $i'=i$ de sorte que
l'annulation doit se faire dans la suite spectrale associée à la stratification avec $j=h \pm 1$ et $h' \neq
h$, i.e. $[\overleftarrow{i},\overrightarrow{d-h-i},\overleftarrow{h-1}]_{1_o} \otimes |\cl|^{-i}$ est un
sous-quotient de
$${\DS \lim_{\genfrac..{0pt}{1}{\to}{I}} H^{h \pm 1}_c(M_{I,\bar s_o}^{=d-h'},R^i \Psi_{\eta_o}(\bar
\Qm_l))[\Pi^\oo]}.$$ D'après le lemme précédent, on doit avoir $h'=h-1$, soit en remarquant encore que
$${\DS \lim_{\genfrac..{0pt}{1}{\to}{I}} H^{h-1}_c(M_{I,\bar s_o}^{=d-h+1},R^i \Psi_{\eta_o}(\bar \Qm_l))[\Pi^\oo]}$$
est de la forme
$\ind_{P_{d-h+1}(F_o)}^{GL_d(F_o)} [\overleftarrow{i},\overrightarrow{d-h-i}]_{1_o} \otimes \pi',$
on obtient que $\pi'$ admet $[\overleftarrow{h-1}]_{1_o}$ comme sous-quotient et donc que
$${\DS \lim_{\genfrac..{0pt}{1}{\to}{I}} H^{h -1}_c(M_{I,\bar s_o}^{=d-h+1},R^i \Psi_{\eta_o}(\bar 
\Qm_l))[\Pi^\oo]}$$
admet $[\overleftarrow{i},\overrightarrow{d-h-i+1},\overleftarrow{h-2}]_{1_o} \otimes |\cl|^{-i}$ comme
sous-quotient; d'où le résultat.

\end{proof}

\begin{corob} Pour tout $0 \leq i <d$, $\Psi_{F_o}^{d,i}(1_o)$ admet
$[\overleftarrow{i},\overrightarrow{d-i-1}]_{1_o} \otimes |\cl|^{-i}$ comme sous-quotient.
\end{corob}

\begin{proof} D'après la proposition ci-dessus, ${\DS \lim_{\genfrac..{0pt}{1}{\to}{I}} H^1_c(M_{I,\bar 
s_o}^{=d-1},R^i
\Psi_{\eta_o}(\bar \Qm_l))[\Pi^\oo]}$ admet comme sous-quotient $[\overleftarrow{i},\overrightarrow{d-i-1}]_{1_o} 
\otimes
|\cl|^{-i}$. Cette dernière représentation ayant des vecteurs fixes sous l'Iwahori
$\Iw_o$, et comme
$${\DS \lim_{\genfrac..{0pt}{1}{\to}{I}} H^1(M_{I,\bar s_o},R^i \Psi_{\eta_o}(\bar \Qm_l))[\Pi^\oo]}$$
n'admet pas comme sous-quotient celle-ci, la suite spectrale associée à la stratification implique que ${\DS
\lim_{\genfrac..{0pt}{1}{\to}{I}} H^0 (M_{I,\bar s_o}^d,R^i \Psi_{\eta_o}(\bar \Qm_l))[\Pi^\oo]}$ admet
$[\overleftarrow{i},\overrightarrow{d-i-1}]_{1_o} \otimes |\cl|^{-i}$ comme sous-quotient. On utilise alors
la proposition (15.2) de \cite{boy}, cf. (\ref{h0-ss}), qui donne un isomorphisme $(D_\Am^\oo)^\times \times
W_o$-équivariant

$$\lim_{\genfrac..{0pt}{1}{\to}{I}} H^0(M_{I,s_o}^d,R^i \Psi_{\eta_o}(\bar \Qm_l))\simeq \hom_{\bar D_o^\times}
((\CC^\oo_{\bar D})^\vee,\widetilde{\UC_{F_o}^{d,i}})$$ où $\bar D$ est l'algèbre à division centrale sur $F$
dont les invariants sont ceux de $D$ excepté en les places $\oo,o$ où ils sont respectivement égaux à $-1/d$
et $1/d$, soit en particulier $\bar D_o \simeq D_{o,d}$, et où $\CC^\oo_{\bar D}$ est l'algèbre de
convolution des fonctions localement constantes sur $\bar D^\times \backslash (\bar D_\Am^\oo)^\times$. Le
résultat découle alors directement d'une correspondance de Jacquet-Langlands globale entre $\Pi$ et une
unique sous-représentation de $\CC^\oo_{\bar D}$ de multiplicité $1$; concrètement l'ensemble des
sous-représentations $\t^\oo$ de $\CC^\oo_{\bar D}$ telles que $\t^{\oo,o} \simeq \Pi^{\oo,o}$ est réduit à
un \ele de multiplicité $1$. Ainsi le lemme de Schur donne
$$\widetilde{\UC_{F_o}^{d,i}}(1_o)= \lim_{\genfrac..{0pt}{1}{\to}{I}} H^0(M^d_{I,s_o},R^i \Psi_{\eta_o} (\bar \Qm_l)) 
[\Pi^{\oo,o}]$$
en tant que représentation de $GL_d(F_o) \times W_o$, d'où le résultat.

\end{proof}

Ainsi d'après (\ref{dim-poids}) et (\ref{inv-iwa}), il ne reste plus de place dans
$\widetilde{\UC_{F_o}^{d,i}}(1_o)$ pour d'autres représentations ayant des vecteurs invariants sous l'Iwahori
$\Iw_o$, autres que les
$$[\overleftarrow{i},\overrightarrow{d-i-1}]_{1_o} \otimes |\cl|^{-i}.$$
De la même façon pour toute représentation irréductible $\t_o$ de $D_{o,d}^\times$ qui n'est pas dans la
classe d'équivalence inertielle de la représentation triviale, $\Psi_{F_o}^{d,i}(\t_o)$ n'admet comme
sous-quotient aucune représentation ayant des vecteurs invariants sous $\Iw_o$.

\end{proof}


%% file: introduction-chap3.tex
\noindent \textbf{0.1.} --- Le but est de calculer la somme alternée des groupes de cohomologie des systèmes
locaux d'Harris-Taylor. La démarche est classique: il s'agit tout d'abord d'utiliser la formule des traces de
Lefschetz et donc de compter les points fixes sous l'action d'une correspondance de Hecke tordue par une
puissance assez grande du Frobenius et ensuite de transférer les intégrales orbitales obtenues afin de
reconnaître le coté géométrique de la formule des traces de Selberg. On en déduit alors un calcul de la somme
alternée des groupes de cohomologie du modèle local de Deligne-Carayol. Dans le cas Iwahori, des arguments de
pureté nous redonne les résultats obtenus à la fin du chapitre précédent.

\medskip

\noindent \textbf{0.2.} --- Dans un premier temps, \S \ref{adel-h}, on donne une description adélique des
points géométriques des variétés d'Igusa de seconde espèce: on ne fait ici qu'adapter les résultats de
\cite{lrs}. On compte alors, \S \ref{fixe}, les points fixes sous l'action des correspondances de Hecke et
d'une puissance arbitraire du frobenius en $o$. Ce comptage se fait en termes d'intégrales orbitales; à
nouveau ces résultats sont une réadaptation de ceux de \cite{lrs}, les fonctions de transfert de la
proposition (\ref{transfert}) étant prises dans \cite{h-t}. Par une application de la formule des traces de
Lefschetz et de Selberg, on en déduit alors, au théorème (\ref{strate-alt}), le calcul de
$[H^*_{h,\r_\oo,\t_o}]$.
\medskip

\noindent \textbf{0.3.} --- Le dernier paragraphe est consacré au calcul, théorème (\ref{theo-calcul-psi}),
de la représentation virtuelle $[\widetilde{\Psi_{F_o}^{d,*}}(\t_o)],$ où $\tau_o$ est une représentation
irréductible quelconque de $D_{o,d}^\times$. Ce résultat est l'équivalent dans notre cadre du théorème
VII.1.5 de \cite{h-t}. En particulier comme d'après le corollaire (\ref{dim-poids}), pour tout $0 \leq i < d$
et toute représentation irréductible admissible $\t_o$ de $D_{o,d}^\times$, $(\Psi_{F_o}^{d,i})^{\Iw_o}$ est
pur de poids $2i$, le calcul virtuel (\ref{theo-calcul-psi}) redonne bien le théorème (\ref{theo-iwahori}).

%% file: f-traces.tex

\section{Cohomologie des systèmes locaux d'Harris-Taylor}

L'essentiel des résultats est déjà présent dans \cite{lrs}, il suffit simplement de les adapter aux variétés
d'Igusa.

\subsection{Description adélique de $M_{I,\bar s_o}(\bar \Fm_p)$}

\label{adel-h} \label{formule-trace}

\begin{defi} Un $\vphi$-espace $(V,\vphi)$ sur $\bar \k(o)$, est un $F \otimes_{\Fm_q} \bar k(o)$-espace vectoriel $V$ 
de dimension finie,
muni d'une application $F \otimes_{\Fm_q} \frob_q$-semi-lin\'eaire bijective,
$$\vphi:V \longto V.$$
\end{defi}

\marque Étant donné un $\DC$-faisceau elliptique $(\EC_i,j_i,t_i)$, soit $V_i=\EC_{i,\eta}$ la fibre de
$\EC_i$ au point générique $\eta \otimes \bar \Fm_q$ de $X \otimes_{\Fm_q} \bar \Fm_q$; pour tout $i$, les
$j_i$ induisent des isomorphismes $V_i \simeq V_{i+1}$. On note $V=V_0$ et $\vphi:V \longto V$ l'application
bijective $F \otimes \frob_q$-semi-linéaire induite par $t_0$: $(V,\vphi)$ est un $\vphi$-espace. L'action de
$D$ sur $V$ commute à $\vphi$ et donne un homomorphisme de $F$-algèbre: $\iota:D^{op} \longto \End(V,\vphi)$.
Le triplet $(V,\vphi,\iota)$ est appelé \textit{la fibre générique} du $\DC$-faisceau elliptique
$(\EC_i,h_i,t_i)$.

\begin{defi}
Deux $\DC$-faisceaux elliptiques de caractéristique $o$ sur $\bar \Fm_q$ sont dits isogènes si leurs
fibres génériques sont isomorphes.
\end{defi}

\marque Si $x$ est une place de $F$, on considère le $F_x$-module de Dieudonné
$$(V_x,\vphi_x):=(F_x \hat \otimes_F V, F_x \hat \otimes_F \vphi)$$
muni du morphisme de $F_x$-algèbre $\l_x:D_x^{op} \longto \End (V_x,\vphi_x).$ On pose
$$M_x=H^0(\spec(\OC_x \hat \otimes \bar \k(o)),\EC_0),$$
qui est un $\DC_x$-réseau de $V_x$ stable sous $\l_x(D_x^{op})$.

\begin{prop} \label{reseau}
(cf. \cite{lrs} proposition 9.4) La construction ci-dessus d\'efinit une bijection entre l'ensemble des classes 
d'isomorphismes
des $\DC$-faisceaux elliptiques sur $\bar \k(o)$ et l'ensemble des classes d'isomorphismes des paires
$$((V,\vphi,\l), (M_x)_{x \in |X|})$$
où $(V,\vphi)$ est un $\vphi$-espace de rang $d^2$ sur $F \otimes \bar \k(o)$, $\l:D^{op} \to \End(V,\vphi)$ est
un morphisme de $F$-algèbres et $(M_x)_{x \in |X|}$ est une collection de $\DC_x$-réseaux des $F_x$-modules
de Dieudonné $(V_x,\vphi_x)=(F_x \hat \otimes_F V, F_x \hat \otimes_F \vphi)$ qui vérifient les propriétés suivantes:

\begin{itemize}

\item (i) si $x=\oo$, on a\footnote{où l'on a supposé pour simplifier $\deg (\oo)=1$}
 $$\vphi_\oo(M_\oo) \supset M_\oo, $$
$$\vphi_\oo^{d}(M_\oo)=\varpi_\oo^{-1}M_\oo,$$
$$\dim_{\bar \k(o)} (\vphi_\oo(M_\oo)/M_\oo)=d,$$

\item (ii) si $x=o$, on a $$\varpi_o M_o \subset \vphi_o(M_o) \subset M_o,$$
le $\k(o) \otimes \bar \k(o)$-module $M_o/\vphi_o(M_o)$ est de longueur $d$ et il est supporté par la composante 
connexe de
$\spec (\k(o) \otimes \bar \k(o))$ qui correspond à l'inclusion $\k(o) \hookrightarrow \bar \k(o)$;

\item (iii) si $x \neq o,\oo$, on a $$\vphi_x(M_x)=M_x;$$

\item (iv) toute base du $F \otimes \bar \k(o)$-espace vectoriel $V$ appartient et engendre le
$\OC_x \hat \otimes \bar \k(o)$-sous-module $M_x$ de $V_x$ pour presque toutes les places $x \neq o,\oo$ de $F$.

\end{itemize}

\end{prop}

\begin{defi} Une $\vphi$-paire $(\tilde F,\tilde \Pi)$ est un couple formé d'une $F$-algèbre $\tilde F$, commutative 
de dimension finie
et d'un \ele $\tilde \Pi \in \tilde F^\times \otimes \Qm$ qui satisfait à la propriété suivante: pour toute 
$F$-sous-algèbre propre $F'$
de $\tilde F$, $\tilde \Pi$ n'appartient pas à $F^{'\times} \otimes \Qm \subset \tilde F^\times \otimes\Qm$.
\end{defi}

A tout $\vphi$-espace $(V,\vphi)$, Drinfeld associe une $\vphi$-paire (cf. \cite{lrs} A.4).

\begin{prop} \label{paire} (cf. \cite{lrs} proposition 9.9)
Soit $(\tilde F,\tilde \Pi)$ la $\vphi$-paire associée au $\vphi$-espace $(V,\vphi)$. On a alors
les propriétés suivantes:

\begin{itemize}

\item (i) $\tilde F$ est un corps et $[\tilde F:F]$ divise $d$;

\item (ii) $F_\oo \otimes_F \tilde F$ est un corps et si $\tilde \oo$ est l'unique place de $\tilde F$ divisant $\oo$, 
on
a l'égalité $\deg(\tilde \oo) \tilde \oo(\tilde \Pi)=-[\tilde F:F]/d$;

\item (iii) il existe une unique place $\tilde o \neq \tilde \oo$ de $\tilde F$ telle que $\tilde o(\tilde \Pi) \neq 
0$; de plus
$\tilde o$ divise $o$;

\item (iv) on a l'égalité $h=d[\tilde F_{\tilde o} :F_o]/[\tilde F:F]$, où $h$ est l'indice de la strate à laquelle 
$(\EC_i,j_i,t_i)$
appartient.

\end{itemize}

\end{prop}

\begin{coro} \label{algebre}
L'algèbre $\End (V,\vphi,\l)$ est une algèbre à division centrale sur $\tilde F$ de dimension
$(d/[\tilde F:F])^2$ dont les invariants sont donnés comme suit:
$$\inv_{\tilde x}(\End(V,\vphi,\l))= \left \{ \begin{array}{ll}
{-[\tilde F:F]/d }& si ~\tilde x=\tilde \oo \\
{[\tilde F:F]/d }& si ~\tilde x=\tilde o \\
{[\tilde F_{\tilde x}:F_x]} \inv_x(D) & sinon \end{array} \right. $$
pour tout place $x$ de $F$ et toute place $\tilde x$ de $\tilde F$ divisant $x$.

\end{coro}

\begin{defi} \label{defphipaire}
Un $(D,\oo,o)$-type est une $\vphi$-paire $(\tilde F,\tilde \Pi)$ telle que:

\begin{itemize}

\item (i) $\tilde F$ est un corps et $[\tilde F:F]$ divise $d$:

\item (ii) $F_\oo \otimes_F \tilde F$ est un corps et si $\tilde \oo$ est l'unique place de $\tilde F$ divisant $\oo$, 
on
a $$\deg(\tilde \oo) \tilde \oo(\tilde \Pi)=-[\tilde F:F]/d;$$

\item (iii) il existe une unique place $\tilde o \neq \tilde \oo$ de $\tilde F$ telle que $\tilde o(\tilde \Pi) \neq 
0$; de plus
$\tilde o$ divise $o$;

\item (iv) pour toute place $x$ de $F$ et toute place $\tilde x$ de $\tilde F$ divisant $x$, on a
$$(d[\tilde F_{\tilde x}:F_x]/[\tilde F:F])\inv_x(D) \in \Zm.$$
\end{itemize}
\end{defi}

\begin{theo} \label{isogene} (cf. \cite{lrs} théorème 9.13)
L'application composée
$$(\EC_i,j_i,t_i) \longmapsto (V,\vphi,\l) \longmapsto (\tilde F,\tilde \Pi),$$
qui à un $\DC$-faisceau elliptique défini sur $\bar \k(o)$ associe son $(D,\oo,o)$-type,
induit une bijection de l'ensemble des classes d'isogénie des $\DC$-faisceaux elliptiques
définis sur $\bar \k(o)$ sur l'ensemble des classes d'isomorphismes des $(D,\oo,o)$-types.
\end{theo}

\marque Toujours selon \cite{lrs} (9.12), la bijection inverse est la suivante. Étant donné $(\tilde F,
\tilde \Pi)$, soit $(W,\psi)$ ``le'' $\vphi$-espace sur $\bar \k(o)$ qui lui correspond et soit $\Delta$
``l'''algèbre à division centrale sur $\tilde F$ dont les invariants sont
$$\inv_{\tilde x} \Delta = \left \{ \begin{array}{ll} [\tilde F:F]/d & \hbox{ si } \tilde x= \tilde \oo \\
-[\tilde F:F]/d & \hbox{ si } \tilde x =\tilde o \\
{[} \tilde F_{\tilde x}:F_x {]} \inv_x (D) & \hbox{ sinon} \end{array} \right.
$$
pour tout place $x$ de $F$ et toute place $\tilde x$ de $\tilde F$ divisant $x$. En particulier $\Delta$ est de
dimension $(d/[\tilde F:F])^2$ sur $\tilde F$ et $D^{op} \otimes_F \Delta$ et $\Mm_d(\End(W,\psi))$ sont des algèbres
simples centrales sur $\tilde F$ de même dimension et possédant les mêmes invariants en toute place $\tilde x$ de
$\tilde F$. En vertu du théorème de Skolem-Noether, on choisit un isomorphisme
$$\a: D^{op} \otimes_F \Delta \longto \Mm_d(\End(W,\psi))$$
et on pose $(V,\vphi):=(W,\psi)^d$. On obtient alors un homomorphisme de $F$-algèbre
$$\iota:D^{op} \longmapright{\delta \mapsto \delta \otimes 1} D^{op} \otimes_F \Delta \longmapright{\a} 
\Mm_d(\End(W,\psi))=\End(V,\vphi)$$
tel que le commutant de $\iota(D^{op})$ dans $\End(V,\vphi)$ est l'image de $\Delta$ par l'homomorphisme de $\tilde 
F$-algèbre
$$\Delta \longmapright{\delta \mapsto 1 \otimes \delta} D^{op} \otimes_F \Delta \longmapright{\a} \End(V,\vphi)$$

\rem Si $(\tilde F,\tilde \Pi)$ est un $(D,\oo,o)$-type associé à la strate $h$, on a
les injections d'algèbres $\Delta^{\oo,o} \hookrightarrow D_\Am^{\oo,o}$,
$\Delta_o^{et}:=\Delta_o^{\tilde o} \hookrightarrow \Mm_{d-h}(F_o)$ et $\Delta_o^c:=\Delta_{\tilde o} \hookrightarrow 
D_{o,h}$.

\begin{prop} Soit $(\tilde F,\tilde \Pi)$ un $(D,\oo,o)$-type associé à la strate $h$.
Il existe pour tout id\'eal $I$ de $A$, une bijection de $\JC_{I^o,m}^{=h}(s)(\bar \k(o))_{(\tilde F,\tilde \Pi)}$ 
avec le quotient
$$\Delta^\times \backslash \left [ (D_\Am^{\oo,o})^\times / K_{\Am,I^o}^{\oo,o} \times GL_{d-h}(F_o)/K_{o,m} \times 
D_{o,h}^\times/(1+\Pi_{o,h}^{s+1} \DC_{o,h}) \right ]$$
De plus si $(\JC_{J^o,m'}^{=h}(s'),c_1,c_2)$ est une correspondance de Hecke naturelle\footnote{i.e. comme dans
\cite{lrs}} sur $\JC_{I^o,m}^{=h}(s)$ associée à un \ele $(g^{\oo,o},g_o^{et},(g_o^c,\d_o,\s_o))$ de
$(D_\Am^{\oo,o})^\times \times GL_{d-h}(F_o) \times \widetilde{\NC_o}$ \footnote{$J^o$ est donc tel que
$K_{\Am,J^o}^{\oo,o} \subset K_{\Am,I^o}^{\oo,o} \cap (g^{\oo,o})^{-1} K_{\Am,I^o}^{\oo,o} g^{\oo,o}$ et $m'$ tel que
$K_{o,m'} \subset K_{o,m} \cap (g_o^{et})^{-1} K_{o,m} g_o^{et}$} la correspondance induite par ces bijections est
$$\xymatrix{
\Delta^\times \backslash [\frac{(D_\Am^{\oo,o})^\times}{K_{\Am,I^o}^{\oo,o}} \times
\frac{GL_{d-h}(F_o)}{K_{o,m}} \times
\frac{D_{o,h}^\times}{(1+\Pi_{o,h}^{s+1} \DC_{o,h})} ] \\
\Delta^\times \backslash [\frac{(D_\Am^{\oo,o})^\times}{K_{\Am,J^o}^{\oo,o}} \times
\frac{GL_{d-h}(F_o)}{K_{o,m'}} \times
\frac{D_{o,h}^\times}{(1+\Pi_{o,h}^{s'+1} \DC_{o,h})}] \dto_{c_1} \uto^{c_2} \\
\Delta^\times \backslash [\frac{(D_\Am^{\oo,o})^\times}{K_{\Am,I^o}^{\oo,o}} \times
\frac{GL_{d-h}(F_o)}{K_{o,m}} \times \frac{D_{o,h}^\times}{(1+\Pi_{o,h}^{s+1} \DC_{o,h}) }] }
$$
où $c_1$ est induit par les inclusions
$$K_{\Am,J}^{\oo,o} \subset K_{\Am,I}^{\oo,o},\quad K_{o,m'} \subset K_{o,m},\quad (1+\Pi_{o,h}^{s'+1} \DC_{o,h}) 
\subset (1+\Pi_{o,h}^{s+1} \DC_{o,h})$$
et où $c_2$ est induit par la multiplication à droite de $(g^{\oo,o})^{-1}$ sur $(D_\Am^{\oo,o})^\times$, de
$(g_o^{et})^{-1}$ sur $GL_{d-h}(F_o)$, et la multiplication à gauche de $\d_o$ sur $D_{o,h}^\times$.
\end{prop}

\begin{proof} En vertu de \cite{lrs} \S 9 et \S 10, $\IC_{I^o,m}^{=h}(\bar \k(o))$ est en bijection avec
$$\Delta^\times \backslash [(D_\Am^{\oo,o})^\times/K_{\Am,I^o}^{\oo,o} \times GL_{d-h}(F_o)/K_{o,m} \times \Zm ]$$
où l'action de $(D_\Am^{\oo,o})^\times \times GL_{d-h}(F_o)$ se décrit comme dans l'énoncé. L'action d'un
\ele $g_o^c$ de $GL_h(F_o)$ n'intervient que sur la composante $\Zm$. Le groupe $PSL_d(F_o)$ \'etant simple,
l'endomorphisme de $\Zm$ associé \`a $g_o^c$ est la translation de valeur $k.\val(\det(g_o))$, pour un
certain entier $k$: en prenant pour $g_o^c$, l'\ele $\varpi_o$ du centre, en utilisant la proposition B.10 de
\cite{lrs}, on obtient $k=-1$. En outre d'après loc. cit., le frobenius géométrique en $o$ agit par
translation de valeur $-1$ sur la composante $\Zm$. Le résultat en découle alors de manière immédiate, en
identifiant $\Zm \times \DC_{o,h}^\times$ avec $D_{o,h}^\times$ où l'on envoie $(n,\d_o)$ sur $\Pi_{o,h}^n
\d_o$.

\end{proof}

\marque On rappelle, selon loc. cit., que le couple $(\tilde F,\tilde \Pi)$ est construit de la manière
suivante. Soit $D^\times_\natural$ l'ensemble des classes de conjugaisons d'\eles de $D^\times$. Soit alors
$\gamma \in D^\times_\natural$ et $F'=F[\gamma]$ tel qu'il existe une unique place $\oo'$ au dessus de $\oo$
ainsi qu'une unique place $o'$ au dessus de $o$ vérifiant $o'(\gamma) \neq 0$ avec $h=\frac{d
[F'_{o'}:F_o]}{[F':F]}$. Soit donc $\Pi' \in F'$ tel que $\oo'(\Pi') \neq 0$, $o'(\Pi') \neq 0$ et
$x'(\Pi')=0$ pour tout $x' \neq \oo',o'$. On définit alors $\tilde F= \cap_{\genfrac{}{}{0pt}{}{n \in \Zm}{n
\neq 0}} F[(\Pi')^n]$ et $\tilde \Pi$. Le couple $(\tilde F,\tilde \Pi)$ est ainsi un $(D,\oo,o)$-type
associé à la strate $h$ et tous ceux-ci sont obtenus par ce procédé. Un tel \ele $\gamma \in D^\times$ est
dit elliptique en $\oo$ et de type $h$ en $o$: son image dans $\Delta$ définit un \ele $\d \in
\Delta^\times_\natural$. On a ainsi $\tilde F \subset F'=F[\gamma]=\tilde F[\d] \subset \Delta$ ainsi qu'une
inclusion naturelle: $\Delta_\Am^{\oo,o} \hookrightarrow D_\Am^{\oo,o}$. Précisons la situation en la place
$o$. A conjugaison près, on peut supposer que
$$\gamma=(\gamma_o^{et},\gamma_o^c) \in GL_{d-h}(F_o) \times GL_h(F_o) \subset GL_d(F_o)$$
avec $F_o[\gamma_o^{et}]=(F')_o^{o'} \subset \Mm_{d-h}(F_o)$ et $F_o[\gamma_o^c]=F'_{o'} \subset \Mm_h(F_o)$. On 
obtient ainsi des
injections naturelles
$$(\Delta_o^{et})^\times:=(\Delta_o^{\tilde o})^\times \hookrightarrow GL_{d-h}(F_o) \hbox{ et }
(\Delta_o^c)^\times:=\Delta_{\tilde o}^\times \hookrightarrow D_{o,h}^\times.$$

\subsection{Description adélique de $\widetilde{M_{I,\bar s_o}}(\bar \Fm_p)$}

On considère désormais les structures de niveau à l'infini comme au paragraphe (\ref{infini}).

\marque En termes de la description (\ref{reseau}), l'application sur les variétés sans niveau
$$r_{\oo}: \widetilde{M_\emptyset} \longto M_\emptyset$$
se décrit comme suit. Pour tout
$$((V,\phi,\l),(M_x)_{x \in | X | }) \in M_\emptyset(\bar \k(o))$$
soit $i_{\oo,0}: \spec (\bar \k(o)) \longto \spec (\k(\oo))$
son pôle, i.e. le support du $\k(\oo) \otimes \bar \k(o)$-module $\vphi_\oo(M_\oo)/M_\oo$.
Via l'identification $\DC_\oo \simeq \Mm_d(\OC_\oo)$,
l'équivalence de Morita donne
$$(V_\oo,\vphi_\oo)=(V'_\oo,\vphi'_\oo)^d \hbox{ et } M_\oo=(M'_\oo)^d.$$
On note $\check M'_\oo$ le dual
du $\OC_\oo \hat \otimes \bar \k(o)$-module libre $M_\oo'$, de rang $d$, et soit $\check \vphi_\oo': \check M'_\oo 
\longto \check M'_\oo$
la restriction de l'application $\check \psi'_\oo: \check V'_\oo \longto \check V'_\oo$ duale de $\vphi'_\oo$.

\begin{prop} (cf. \cite{lrs} page 274) L'ensemble $\widetilde{M_\emptyset}(\bar \k(o))$
est en bijection avec l'ensemble des classes d'isomorphismes des triplets
$$((V,\phi,\l),(M_x)_{x \in | X | },(\nu,\a))$$
où $((V,\phi,\l),(M_x)_{x \in | X | })$ appartient à $M_\emptyset(\bar \k(o))$ et
$$\nu:\spec (\bar \k(o)) \longto \spec (\k(\oo)_d)$$
est un relèvement de $i_{\oo,0}$ et
$$\a:\NC_{d,1} \longmapright{\sim} \check M'_\oo$$
est un isomorphisme de $\OC_\oo \hat \otimes \bar \k(o)$-modules qui commute avec les $\psi$.
L'application $r_\oo$ envoie $((V,\phi,\l),(M_x)_{x \in | X | },(\nu,\a))$ sur $((V,\phi,\l),(M_x)_{x \in | X | })$
et $\bar D_\oo^\times \simeq \bar \DC_\oo^\times \semi \Zm/d \Zm$ agit comme décrit au paragraphe (\ref{infini}).
\end{prop}

\marque On rappelle qu'étant donné un $(D,\oo,o)$-type, $(\tilde F,\tilde \Pi)$, on a une inclusion
d'algèbres
$$\Delta \hookrightarrow \End(N_{d,1},\phi_{d,1})$$
On identifie $\End(N_{d,1},\phi_{d,1})$ avec $(\bar D_\oo)^{op}$ de sorte que l'on obtient une inclusion
$\Delta^\times \hookrightarrow ((\bar D_\oo)^{op})^\times \simeq \bar D_\oo^\times$.
On note $\widetilde{\JC_{I^o,m}^{=h}}(s)$ le produit fibré
$$\widetilde{\JC_{I^o,m}^{=h}}(s):= \JC_{I^o,m}^{=h}(s) \times_{M_{I^o}} \widetilde{M_{I^o}}.$$

\begin{coro}
Soit $(\tilde F,\tilde \Pi)$ un $(D,\oo,o)$-type associé à la strate $h$.
Il existe pour tout idéal $I$ de $A$, une bijection de
$\widetilde{\JC_{I^o,m}^{=h}}(s)(\bar \k(o))_{(\tilde F,\tilde \Pi)}$ avec le quotient
$$\Delta^\times \backslash \left [ (\bar D_\oo^\times /\varpi_\oo^\Zm) \times
(D_\Am^{\oo,o})^\times / K_{\Am,I^o}^{\oo,o} \times GL_{d-h}(F_o)/K_{o,m} \times D_{o,h}^\times/(1+\Pi_{o,h}^{s+1}
\DC_{o,h}) \right ]$$ De plus si $(\widetilde{\JC_{J^o,m'}^{=h}}(s'),c_1,c_2)$ est une correspondance de Hecke
naturelle\footnote{i.e. comme dans \cite{lrs}} sur $\widetilde{\JC_{I^o,m}^{=h}}(s)$ associée à un \ele
$$(\bar g_\oo,g^{\oo,o},g_o^{et},(g_o^c,\d_o,\s_o)) \in
(\bar D_\oo^\times/\varpi_\oo^\Zm) \times (D_\Am^{\oo,o})^\times \times GL_{d-h}(F_o) \times \widetilde{\NC_o}$$
\footnote{$J^o$ est donc tel que $K_{\Am,J^o}^{\oo,o} \subset K_{\Am,I^o}^{\oo,o} \cap (g^{\oo,o})^{-1}
K_{\Am,I^o}^{\oo,o} g^{\oo,o}$ et $m'$ tel que $K_{o,m'} \subset K_{o,m} \cap (g_o^{et})^{-1} K_{o,m} g_o^{et}$} la
correspondance induite par ces bijections est
$$\xymatrix{
\Delta^\times \backslash [(\bar D_\oo^\times/\varpi_\oo^\Zm) \times
\frac{(D_\Am^{\oo,o})^\times}{K_{\Am,I^o}^{\oo,o}} \times \frac{GL_{d-h}(F_o)}{K_{o,m}} \times
\frac{D_{o,h}^\times}{(1+\Pi_{o,h}^{s+1} \DC_{o,h})} ] \\
\Delta^\times \backslash [(\bar D_\oo^\times/\varpi_\oo^\Zm) \times
\frac{(D_\Am^{\oo,o})^\times}{K_{\Am,J^o}^{\oo,o}} \times \frac{GL_{d-h}(F_o)}{K_{o,m'}} \times
\frac{D_{o,h}^\times}{(1+\Pi_{o,h}^{s'+1} \DC_{o,h})}] \dto_{c_1} \uto^{c_2} \\
\Delta^\times \backslash [(\bar D_\oo^\times/\varpi_\oo^\Zm) \times
\frac{(D_\Am^{\oo,o})^\times}{K_{\Am,I^o}^{\oo,o}} \times \frac{GL_{d-h}(F_o)}{K_{o,m}} \times
\frac{D_{o,h}^\times}{(1+\Pi_{o,h}^{s+1} \DC_{o,h}) }] }
$$
où $c_1$ est induit par les inclusions
$$K_{\Am,J}^{\oo,o} \subset K_{\Am,I}^{\oo,o}, \quad K_{o,m'} \subset K_{o,m}, \quad
(1+\Pi_{o,h}^{s'+1} \DC_{o,h}) \subset (1+\Pi_{o,h}^{s+1} \DC_{o,h})$$ et où $c_2$ est induit par la
multiplication à gauche par $g_o$ sur $(\bar D_\oo^\times/\varpi_\oo^\Zm)$, à droite de $(g^{\oo,o})^{-1}$
sur $(D_\Am^{\oo,o})^\times$, de $(g_o^{et})^{-1}$ sur $GL_{d-h}(F_o)$, et la multiplication à gauche de
$\d_o$ sur $D_{o,h}^\times$.
\end{coro}

\rem L'application
$$\l(\k)_{(\tilde F,\tilde \Pi)}: \widetilde{\JC_{I^o,m}^{=h}}(s)(\k)_{(\tilde F,\tilde \Pi)} \longto \spec 
(\k(\oo)_d)$$
est induite par l'application
$$\bar D_\oo^\times / \varpi_\oo^\Zm \longto \Zm/d \deg(\oo) \Zm$$
qui envoie $(n,\d)$ sur $n + \deg (\oo) \oo(\rn(\d)) \mod d \deg(\oo) \Zm$.

\subsection{Comptage des points fixes}

\label{fixe}

On fixe un sous-groupe normal $\bar K_\oo \subset (\bar D_\oo^\times / \varpi_\oo^\Zm)$ et
on note
$$M(s)=\widetilde{\JC_{I^o,m}^{=h}}(s)/ \bar K_\oo$$
qui est donc muni d'une action du produit $\DC_{o,h}^\times / (1+\Pi_{o,h}^{s+1} \DC_{o,h}) \times (\bar
D_\oo^\times/ \varpi_\oo^\Zm)/ \bar K_\oo$ et d'une action par correspondances de Hecke associées aux \eles
$g^{\oo,\tilde o}=(g^{\oo,o},g_o^{et}) \in (D_\Am^{\oo,o})^\times \times GL_{d-h}(F_o)$:
$$\xymatrix{
& M(g^{\oo,\tilde o},s)_{(\tilde F,\tilde \Pi)}(\bar \Fm_q) \dlto_{c_1} \drto^{c_2} \\
M(s)_{(\tilde F,\tilde \Pi)}(\bar \Fm_q) & & M(s)_{(\tilde F,\tilde \Pi)}(\bar \Fm_q) \ar@{-->}[ll]
}$$
où $M(g^{\oo,\tilde o},s)$ est le quotient
$$(\lim_{\genfrac..{0pt}{1}{\lefto}{I}^o,m} \widetilde{\JC_{I^o,m}^{=h}}(s))/((K_{\Am,I^o}^{\oo,o} \cap 
(g^{\oo,o})^{-1} K_{\Am,I^o}^{\oo,o}
g^{\oo,o}) \times K_{o,m} \cap (g_o^{et})^{-1} K_{o,m} g_o^{et}$$
et où la correspondance est induite par
$$\Delta^\times \backslash \left [ \frac{\bar D_\oo^\times /\varpi_\oo^\Zm}{\bar K_\oo} \times
\frac{(D_\Am^{\oo,o})^\times}{K_{\Am,I^o}^{\oo,o}} \times \frac{GL_{d-h}(F_o)}{K_{o,m}} \times
\frac{D_{o,h}^\times}{(1+\Pi_{o,h}^{s+1} \DC_{o,h}^\times)} \right ]$$
$$ \uparrow ~ c_2$$
\begin{multline*}
\Delta^\times \backslash \{ \frac{\bar D_\oo^\times / \varpi_\oo^\Zm}{\bar K_\oo} \times
\frac{(D_\Am^{\oo,o})^\times}{K_{\Am,I^o}^{\oo,o} \cap (g^{\oo,o})^{-1} K_{\Am,I^o}^{\oo,o} g^{\oo,o}}
\times \\
\frac{GL_{d-h}(F_o)}{K_{o,m} \cap (g_o^{et})^{-1} K_{o,m} g_o^{et}} \times
\frac{D_{o,h}^\times}{(1+\Pi_{o,h}^{s+1} \DC_{o,h})} \}
\end{multline*}
$$\downarrow ~c_1$$
$$\Delta^\times \backslash \left [ \frac{\bar D_\oo^\times / \varpi_\oo^\Zm}{\bar K_\oo} \times
\frac{(D_\Am^{\oo,o})^\times}{K_{\Am,I^o}^{\oo,o}} \times \frac{GL_{d-h}(F_o)}{K_{o,m}} \times
\frac{D_{o,h}^\times}{(1+\Pi_{o,h}^{s+1} \DC_{o,h})} \right ] $$ où
$$c_1[h^{\oo,o}(K_{\Am,I^o}^{\oo,o} \cap (g^{\oo,o})^{-1} K_{\Am,I^o}^{\oo,o} g^{\oo,o}),h_o^{et}(K_{o,m}\cap 
(g_o^{et})^{-1} K_{o,m} g_o^{et}),d_o (1+\Pi_{o,h}^{s+1}\DC_{o,h})]$$
$$~ \hfill =[h^{\oo,o} K_{\Am,I^o}^{\oo,o},h_o^{et} K_{o,m},d_o (1+\Pi_{o,h}^{s+1} \DC_{o,h})]$$
et
$$c_2[h^{\oo,o}(K_{\Am,I^o}^{\oo,o} \cap (g^{\oo,o})^{-1} K_{\Am,I^o}^{\oo,o} g^{\oo,o}),h_o^{et}(K_{o,m}\cap 
(g_o^{et})^{-1} K_{o,m} g_o^{et}),d_o (1+\Pi_{o,h}^{s+1}\DC_{o,h})]$$
$$~ \hfill =[h^{\oo,o}(g^{\oo,o})^{-1} K_{\Am,I^o}^{\oo,o},h_o^{et}(g_o^{et})^{-1} K_{o,m},d_o (1+\Pi_{o,h}^{s+1} 
\DC_{o,h})]$$

\begin{rema} \label{rema-action}
 Dans cette description l'action de
 $$(g_o^c,\d_o,\s_o) \in \NC_o \subset GL_h(F_o) \times D_{o,h}^\times \times W_{F_o}$$
sur $M_{(\tilde F,\tilde \Pi)}$ est donnée par la multiplication à
gauche de $\d_o$ sur la composante $D_{o,h}^\times/(1+\Pi_{o,h}^{s+1}\DC_{o,h})$.
\end{rema}

\begin{defi} Pour $r>0$, $\bar g_\oo \in (\bar D_\oo^\times/ \varpi_\oo^\Zm)$ et $\d_o \in \DC_{o,h}^\times$, on note
$$\fix_r^{=h}(\bar g_\oo,\d_o,g_o^c,g^{\oo,\tilde o},s)$$
l'ensemble $\{ m \in M(g^{\oo,\tilde o},s)(\bar \Fm_q)~/~(g_o^c,\Pi_{o,h}^{-\val(\det(g_o^c))-r}
\d_o,\frob_o^r)(c_1(m)).\bar g_\oo \bar K_\oo=c_2(m) \}.$
\end{defi}

\begin{lemm} (cf. \cite{lrs} lemme 11.1) Les ensembles $\fix_r^{=h}(\d_o,g_o^c,g^{\oo,\tilde o},s)$ pour $r>0$ sont 
finis et
chacun de ces points fixes est de multiplicité $1$.
\end{lemm}

\marque On décompose ces ensembles, suivant leur $(D,\oo,o)$-type
$$\fix_r^{=h}(\bar g_\oo,\d_o,g_o^c,g^{\oo,\tilde o},s)=\coprod_{(\tilde F,\tilde \Pi)} \fix_r^{=h}(\bar 
g_\oo,\d_o,g_o^c,g^{\oo,\tilde o},s)_{(\tilde F,\tilde \Pi)}$$
où la réunion porte sur les $(D,\oo,o)$-type associés à la strate $h$.

\noindent Soit $(\tilde F,\tilde \Pi)$ un $(D,\oo,o)$-type associé à la strate $h$ et soit $\Delta$ l'algèbre à 
division centrale sur $F$
correspondante. On note $\Delta_\natural^\times$ un système de représentants des classes de conjugaisons dans 
$\Delta^\times$ et
pour $\d \in \Delta^\times$, soit $\Delta_\d^\times$ le centralisateur de $\d$ dans $\Delta^\times$.

\begin{prop} L'ensemble $\fix_r^{=h}(\bar g_\oo,\d_o,g_o^c,g^{\oo,\tilde o},s)_{(\tilde F,\tilde \Pi)}$
est la réunion disjointe des doubles classes
$$\coprod_{\d \in \Delta_\natural^\times} \Delta_\d^\times
[\bar h_\oo \bar K_\oo,h^{\oo,o}K_{\Am,I^o}^{\oo,o},h_o^{et} K_{o,m},d_o (1+\Pi_{o,h}^{s+1} \DC_{o,h})]$$
telles que
$$\left \{ \begin{array}{l}
(\bar h_\oo)^{-1} \d \bar h_\oo \in \bar g_\oo \bar K_\oo \\
(h^{\oo,o})^{-1} \d h^{\oo,o} \in K_{\Am,I^o}^{\oo,o} g^{\oo,o} K_{\Am,I^o}^{\oo,o} \\
(h_o^{et})^{-1} \d h_o^{et} \in K_{o,m} g_o^{et} K_{o,m} \\
\d \in \Pi_{o,h}^{r+\val(\det(g_o^c))} \d_o^{-1} (1+\Pi_{o,h}^{s+1} \DC_{o,h})
\end{array} \right.$$
\end{prop}

\begin{proof} D'après la description des correspondances de Hecke donnée ci-dessus,
\begin{multline*}
\Delta^\times [\bar h_\oo \bar K_\oo,h^{\oo,o}(K_{\Am,I^o}^{\oo,o} \cap
(g^{\oo,o})^{-1} K_{\Am,I^o}^{\oo,o} g^{\oo,o}), \\
h_o^{et}(K_{o,m}\cap (g_o^{et})^{-1} K_{o,m} g_o^{et}),d_o (1+\Pi_{o,h}^{s+1}\DC_{o,h})]
\end{multline*}
appartient à $\fix_r^{=h}(\d_o,g_o^c,g^{\oo,\tilde o},s)_{(\tilde F,\tilde \Pi)}$
\ssi il existe $\d \in \Delta^\times$, $\bar k_\oo \in \bar K_\oo$, $k^{\oo,o} \in K_{\Am,I^o}^{\oo,o}$,
$k_o^{et} \in K_{o,m}$ et $k_s \in (1+\Pi_{o,h}^{s+1} \DC_{o,h})$ tels que
$$\left \{ \begin{array}{l}
\bar h_\oo \bar g_\oo = \d \bar h_\oo \bar k_\oo \\
h^{\oo,o}=\d h^{\oo,o} (g^{\oo,o})^{-1} k^{\oo,o} \\
h_o^{et}=\d h_o^{et} (g_o^{et})^{-1} k_o^{et} \\
d_o=\d \Pi_{o,h}^{-r-\val(\det(g_o^c))} \d_o d_o k_s
\end{array} \right. $$
soit, \ssi, il existe $\d \in \Delta^\times$ tel que
$$\left \{ \begin{array}{l}
(\bar h_\oo)^{-1} \d \bar h_\oo \in \bar g_\oo \bar K_\oo \\
(h^{\oo,o})^{-1} \d h^{\oo,o} \in K_{\Am,I^o}^{\oo,o} g^{\oo,o} \\
(h_o^{et})^{-1} \d h_o^{et} \in K_{o,m} g_o^{et} \\
\d \in \Pi_{o,h}^{r+\val(\det(g_o^c))} \d_o^{-1} (1+\Pi_{o,h}^{s+1} \DC_{o,h})
\end{array} \right.$$
en particulier, on a $o(\rn(\d))=r + \val(\det g_o^c)$.
L'application
$$h^{\oo,o}(K_{\Am,I^o}^{\oo,o} \cap (g^{\oo,o})^{-1} K_{\Am,I^o}^{\oo,o} g^{\oo,o}) \mapsto h^{\oo,o} 
K_{\Am,I^o}^{\oo,o}$$
$$(\hbox{resp. } h_o^{et}(K_{o,m}\cap (g_o^{et})^{-1} K_{o,m} g_o^{et}) \mapsto h_o^{et} K_{o,m})$$
de l'ensemble des classes satisfaisant
$$(h^{\oo,o})^{-1} \d h^{\oo,o} \in K_{\Am,I^o}^{\oo,o} g^{\oo,o} \hbox{ (resp. }(h_o^{et})^{-1} \d h_o^{et} \in 
K_{o,m} g_o^{et})$$
dans l'ensemble des classes vérifiant
$$(h^{\oo,o})^{-1} \d h^{\oo,o} \in K_{\Am,I^o}^{\oo,o} g^{\oo,o} K_{\Am,I^o}^{\oo,o} \hbox{ (resp. } (h_o^{et})^{-1} 
\d h_o^{et} \in K_{o,m} g_o^{et} K_{o,m})$$
est clairement bijective. Le résultat découle alors du lemme suivant.

\begin{lemm} \label{unicite}
Pour tout $h^{\oo,o} \in (D_\Am^{\oo,o})^\times$, $h_o^{et} \in GL_{d-h}(F_o)$ et $d_o \in D_{o,h}^\times$, le seul
\ele $\d \in \Delta^\times$ tel que
$$\left \{ \begin{array}{l}
(h^{\oo,o})^{-1} \d h^{\oo,o} \in K_{\Am,I^o}^{\oo,o} \\
(h_o^{et})^{-1} \d h_o^{et} \in K_{o,m}\\
o(\rn(\d))=0
\end{array} \right. $$
est l'identité.
\end{lemm}

\begin{proof} cf. \cite{lau} (3.2.6).

\end{proof}

\subsection{Intégrales orbitales}

Pour tout $\d \in \Delta^\times$, soit
$$(D_\Am^{\oo,o})^\times_\d \quad (\hbox{resp. }GL_{d-h}(F_o)_\d, \quad \hbox{resp. }(D_{o,h})^\times_\d, \quad 
\hbox{resp. }
(\bar D_\oo^\times)_\d)$$ le centralisateur de $\d$ dans $(D_\Am^{\oo,o})^\times$ (resp. $GL_{d-h}(F_o)$,
resp. $D_{o,h}^\times$, resp. $\bar D_\oo^\times$). On fixe une mesure de Haar,
$$dh^{\oo,o} \quad (\hbox{resp. }dh_o^{et}, \quad resp.~d \bar h_\oo)$$
sur $(D_\Am^{\oo,o})^\times$ (resp. $GL_{d-h}(F_o)$, resp. $\bar D_\oo^\times/ \varpi_\oo^\Zm$) normalisée
par
$$\vol(K_{\Am,I}^{\oo,o},dh^{\oo,o})=1 \quad (\hbox{resp. }\vol(K_{o,m},dh_o^{et})=1 \hbox{resp. }\vol(\bar
K_\oo,d \bar h_\oo)=1).$$ Soit
$$dh^{\oo,o}_\d \quad (\hbox{resp. } dh_{o,\d}^{et}, \quad \hbox{resp. } d\bar h_{o,\d},\quad \hbox{resp. }d \bar 
h_{\oo,\d})$$
une mesure de Haar sur $(D_\Am^{\oo,o})^\times_\d$ (resp. $GL_{d-h}(F_o)_\d$, resp. $(D_{o,h})^\times_\d$, resp. 
$(\bar
D_\oo^\times)_\d / \varpi_\oo^\Zm$) et soit $d \d'$ la mesure de comptage sur $\Delta^\times_\d$. Soit alors
$$f^{\oo,o} \quad (\hbox{resp. } f_o^{et},\quad \hbox{resp. } \bar f_\oo)$$
la fonction caractéristique de
$$K_{\Am,I^o}^{\oo,o} g^{\oo,o}K_{\Am,I^o}^{\oo,o} \quad (\hbox{resp. }K_{o,m} g_o^{et} K_{o,m}, \quad
\hbox{resp. }\bar g_\oo \bar K_\oo)$$ dans $(D_\Am^{\oo,o})^\times$ (resp. $GL_{d-h}(F_o)$, resp. $\bar
D_\oo^\times / \varpi_\oo^\Zm$). On introduit alors les intégrales orbitales
$$O_\d(\bar f_\oo,d \bar h_{\oo,\d})=\int_{(\bar D_\oo^\times)_\d \backslash \bar D_\oo^\times} \bar f_\oo((\bar 
h_\oo)^{-1} \d \bar h_\oo)
\frac{d \bar h_\oo}{d \bar h_{\oo,\d}},$$
$$O_\d(f^{\oo,o},dh^{\oo,o}_\d)=\int_{(D_\Am^{\oo,o})^\times_\d \backslash (D_\Am^{\oo,o})^\times} 
f^{\oo,o}((h^{\oo,o})^{-1} \d h^{\oo,o}) \frac{dh^{\oo,o}}{dh^{\oo,o}_\d},$$
$$O_\d(f_o^{et},dh_{o,\d}^{et})=\int_{GL_{d-h}(F_o)_\d \backslash GL_{d-h}(F_o)} f_o^{et}((h_o^{et})^{-1} \d h_o^{et}) 
\frac{dh_o^{et}}{dh^{et}_{o,\d}},$$
Ces intégrales sont absolument convergentes et on introduit le volume $V$ défini par
$$\vol(\Delta_\d^\times \backslash [((\bar D_\oo^\times)_\d / \varpi_\oo^\Zm) \times
(D_\Am^{\oo,o})^\times_\d \times GL_{d-h}(F_o)_\d],\frac{dh^{\oo,o}_\d dh_{o,\d}^{et}}{d \d'})$$

Pour tout $\d \in \Delta^\times$, l'inclusion d'algèbres $\Delta^{\oo,o} \hookrightarrow D_\Am^{\oo,o}$
(resp. $\Delta_o^{et} \hookrightarrow \Mm_{d-h}(F_o)$, resp. $\Delta_o^c \hookrightarrow D_{o,h}$, resp.
$\Delta \hookrightarrow \bar D_\oo^{op}$), induit un isomorphisme de groupe $(\Delta^{\oo,o})^\times_\d
\simeq (D_\Am^{\oo,o})^\times_\d$ (resp. $(\Delta_o^{et})^\times_\d \simeq GL_{d-h}(F_o)_\d$, resp.
$(\Delta_o^c)^\times_\d \simeq (D_{o,h})^\times_\d$, resp. $(\Delta_\oo^\times)_\d \simeq (\bar
D_\oo^\times)_\d$). En particulier la mesure de Haar $d \bar h_{\oo,\d} \times dh^{\oo,o}_\d \times
dh_{o,\d}^{et}$ induit une mesure de Haar $d h_\d$ sur $(\Delta_\Am^\times)_\d/ \varpi_\oo^\Zm$.

\begin{lemm} Le volume $V$ est égal au volume
$$\vol(\Delta_\d^\times \backslash (\Delta_\Am^\times)_\d/ \varpi_\oo^\Zm, \frac{d h_\d}{d \d'})$$
\end{lemm}

\begin{proof} cf. \cite{lau} (3.3.4).
\end{proof}

\marque Soit $\lef_r^{=h}(\bar g_\oo,\d_o,g_o^c,g^{\oo,\tilde o},s)_{(\tilde F,\tilde \Pi)}$ le cardinal de
$\fix_r^{=h}(\bar g_\oo,\d_o,g_o^c,g^{\oo,\tilde o},s)_{(\tilde F,\tilde \Pi)}$. D'après ce qui précède
$\lef_r^{=h}(\bar g_\oo,\d_o,g_o^c,g^{\oo,\tilde o},s)_{(\tilde F,\tilde \Pi)}$ est égal à

\begin{multline*}
\sum_{\genfrac{}{}{0pt}{}{\d \in \Delta^\times_\natural}{\d \in \Pi_{o,h}^{r+\val(\det g_o^c)}
\d_o^{-1}(1+\Pi_{o,h}^{s+1}\DC_{o,h})}} \vol(\Delta^\times_\d \backslash (\Delta_\Am^\times)_\d/
\varpi_\oo^\Zm, \frac{d h_\d}{d \d'}) \\
 O_\d(\bar f_\oo, d \bar h_{\oo,\d}) O_\d(f^{\oo,o},dh^{\oo,o}_\d) O_\d(f_o^{et},dh_{o,\d}^{et})
\end{multline*}

\begin{lemm} \label{bij} (cf. \cite{lau} (3.4))
La construction de la fin du paragraphe (\ref{adel-h}) définit une bijection naturelle
$$\{ \gamma \in D^\times_\natural ~/~ \gamma \hbox{ est elliptique en } \oo \hbox{ et de type } h \hbox{ en } o \}
\longmapright{\sim} \coprod_{(\tilde F,\tilde \Pi)} \{ \d \in \Delta^\times_\natural \}$$
où $(\tilde F,\tilde \Pi)$ décrit l'ensemble des $(D,\oo,o)$-type associé à la strate $h$.
\end{lemm}

\marque Soit alors $\gamma \in D^\times_\natural$ elliptique en $\oo$ et de type $h$ en $o$ et soit $((\tilde
F,\tilde \Pi),\d)$ le triplet correspondant.

\begin{itemize}
\item On peut voir $\bar D_\oo^\times$ comme une forme intérieure de $D_\oo^\times=GL_d(F_\oo)$ et si $\bar \gamma \in 
\bar D_\oo^\times$
est le transfert de $\gamma \in D_\oo^\times$ défini à conjugaison près, $\bar \gamma$ et l'image de $\d \in 
\Delta^\times$ dans
$\bar D_\oo^\times$ sont conjugués dans $\bar D_\oo^\times$; on peut alors identifier le centralisateur $(\bar 
D_\oo^\times)_{\bar \gamma}$
de $\bar \gamma$ dans $\bar D_\oo^\times$ avec $(\Delta_\oo^\times)_\d$.

\item Comme $\Delta_\Am^{\oo,o}$ est le centralisateur de $\tilde F$ dans $D_\Am^{\oo,o}$ avec $\tilde F \subset
F'=F[\gamma]=\tilde F[\d] \subset \Delta$, le centralisateur $(D_\Am^{\oo,o})^\times_\gamma$ de $\gamma$ dans
$(D_\Am^{\oo,o})^\times$ coïncide avec $(\Delta_\Am^{\oo,o})^\times_\d$.

\item De même comme $\Delta_o^{et}$ (resp. $\Delta_o^c$) est le centralisateur de $\tilde F_o^{\tilde o}$ (resp.
$\tilde F_{\tilde o}$) dans $\Mm_{d-h}(F_o)$ (resp. $D_{o,h}$), le centralisateur $GL_{d-h}(F_o)_{\gamma_o^{et}}$ 
(resp.
$(D_{o,h})^\times_{\gamma_o^c}$) de $\gamma_o^{et}$ (resp. $\gamma_o^c$) dans $GL_{d-h}(F_o)$ (resp. 
$D_{o,h}^\times$)
coïncide avec $(\Delta_o^{et})^\times_\d$ (resp. $(\Delta_o^c)^\times_\d$).
\end{itemize}

\marque La mesure de Haar $d h_\d=d \bar h_{\oo,\d} d h^{\oo,o}_\d \times d h^{et}_{o,\d} \times d \bar
h_{o,\d}$ sur $(\Delta^\times)_\d$ définit alors une mesure de Haar $d h^{o}_\gamma$ (resp. $d
h_{o,\gamma}^{et}$, resp. $d \bar h_{o,\gamma}$) sur $(D_\Am^{o})^\times_\gamma$ (resp.
$GL_{d-h}(F_o)_\gamma$, resp. $(D_{o,h})^\times_\gamma$). Soit comme précédemment $\gamma \in
D^\times_\natural$ elliptique en $\oo$ et de type $h$ en $o$ et soit $((\tilde F,\tilde \Pi),\d)$ le triplet
associé. Le centralisateur $D^\times_\gamma$ de $\gamma$ dans $D^\times$ est une forme intérieure du
centralisateur $\Delta^\times_\d$ de $\d$ dans $\Delta^\times$.

\begin{lemm} On a
$$\vol(\Delta^\times_\d \backslash (\Delta_\Am^\times)_\d/ \varpi_\oo^\Zm ,
\frac{d \d_{\Am,\d}}{d \d'})= \vol(D_\gamma^\times \backslash (D_\Am^\times)_\gamma/ \varpi_\oo^\Zm ,\frac{d
h_{\Am,\gamma}}{d h_\gamma})$$
\end{lemm}

\begin{proof} cf. \cite{lau} (3.5).

\end{proof}

Soit $D^\times_{\natural_h}$ un système de représentants des classes de conjugaisons de $D^\times$ elliptiques en 
$\oo$ et
de type $h$ en $o$.

\medskip

\rem - Soit $\gamma \in D^\times_{\natural_h}$ et $\d \in \Delta^\times_\natural$ l'\ele qui lui est associé
d'après le lemme (\ref{bij}). L'image de $\d$ par l'inclusion $(D^{\oo,o})^\times \hookrightarrow
(D_\Am^{\oo,o})^\times$ est conjuguée à $\gamma \in D^\times \subset (D_\Am^{\oo,o})^\times$. De même en
notant $\d_o^{et}$ (resp. $\bar \d_o$) l'image de $\d$ dans $\Delta_o^{et} \hookrightarrow \Mm_{d-h}(F_o)$
(resp. $\Delta_o^c \hookrightarrow D_{o,h}$), si $g_o^c$ est un \ele semi-simple de $GL_h(F_o)$ de même
polynôme caractéristique que $\d_o^c$, alors $(g_o^{et},g_o^c) \in GL_d(F_o)$ est conjugué à $\gamma_o \in
D_o^\times$.

- Soit $g_o^{et} \in GL_{d-h}(F_o)$ (resp. $\bar g_o \in D_{o,h}^\times$) dans le support de $f_o^{et}$
(resp. de $\bar f_o$). Si $g_o^c$ est un \ele semi-simple de $GL_h(F_o)$ de même polynôme caractéristique que
$\bar g_o$, le centralisateur de $(g_o^{et},g_o^c)$ dans $GL_d(F_o)$ est le même que dans $GL_{d-h}(F_o)
\times GL_h(F_o)$; il suffit en effet de remarquer que l'on a pris $r>0$ pour la définition de $\bar f_o$.
Ainsi étant donnée une mesure de Haar $dh_{(g_o^{et},g_o)}$ sur $GL_d(F_o)_{(g_o^{et},g_o^c)}$, on lui
associe les mesures de Haar $dh_{g_o^{et}}$ et $dh_{g_o^c}$. On obtient alors la proposition suivante.

\begin{prop} On considère, comme précédemment, les mesures de Haar $d \bar h_\oo$, $d h^\oo$. Pour $\gamma \in 
D^\times_{\natural_h}$,
on fixe une mesure de Haar $d h_{\Am,\gamma}=d h_{\oo,\gamma} \times d h^\oo_\gamma$ sur le centralisateur
$(D_\Am^\times)_\gamma / \varpi_\oo^\Zm$ dans $D_\Am^\times / \varpi_\oo^\Zm$ et soit $dh_\gamma$ la mesure
de comptage sur $D_\gamma^\times$. Soit $\bar \gamma \in \bar D_\oo^\times$, défini à conjugaison près, le
transfert de $\gamma$; son centralisateur $(\bar D_\oo^\times)_{\bar \gamma}$ dans $\bar D_\oo^\times$ est
une forme intérieure de $(D_\oo^\times)_\gamma$; soit alors $d \bar h_{\oo,\bar \gamma}$ le transfert  de $d
h_{\oo,\gamma}$, la mesure de Haar sur $(\bar D_\oo^\times)_{\bar \gamma} / \varpi_\oo^\Zm$. On note
$$\e_\oo(\gamma)=(-1)^{\frac{d}{[ F_\oo[\gamma]:F_\oo ] } -1}$$
le signe de Kottwitz à l'infini de $\gamma$.

En posant $f^\oo=f^{\oo,o} f_o^{et}$, le nombre de Lefschetz $\lef_r^{=h}(\bar g_\oo,\d_o,g_o^c,g^{\oo,\tilde
o},s)$ est alors donné par

\begin{multline*}
\sum_{\gamma \in D^\times_{\natural_h} \cap C_h(\Pi_{o,h}^{\val(\det g_o^c)+r} \d_o^{-1})} \e_\oo(\gamma)
\vol(D^\times_\gamma \backslash (D_\Am)^\times_\gamma/\varpi_\oo^\Zm,\frac{d h_{\Am,\gamma}}{d h_\gamma}) \\
O_{\bar \gamma}(\bar f_\oo,d \bar h_{\oo,\bar \gamma}) O_{\gamma^o}(f^o,dh^o_\gamma)
O_{g_o^{et}}(f_o^{et},dh_{g_o^{et}})
\end{multline*}

où: - $C_h(\Pi_{o,h}^\a \d_o^{-1}))$ désigne l'ensemble des $\gamma$ conjugués à un élément de la forme
$(g_o^{et},g_o^d)$ avec $g_o^{et} \in GL_{d-h}(F_o)$ et $g_o^d \in GL_h(F_o)$ semi-simple de même polynôme
caractéristique que $\Pi_{o,h}^\a \d_o^{-1}$;

- $dh_{g_o^{et}}$ et $dh_{g_o^c}$ sont respectivement les mesures de Haar sur $GL_{d-h}(F_o)_{g_o^{et}}$ et
$GL_h(F_o)$ associées à la mesure de Haar $dh_{o,\gamma}$ par le procédé rappelé ci-dessus.
\end{prop}

\subsection{Formule des traces de Lefschetz}

On fixe une représentation irréductible
$$\r_\oo:\bar D_\oo^\times/ \varpi_\oo^\Zm \longto GL(L)$$
sur un $\bar \Qm_l$-espace vectoriel de dimension finie $L$, qui est définie sur une extension finie $E_\l$
de $\Qm_l$ dans $\bar \Qm_l$ et qui est continue pour la topologie pro-finie sur $\bar D_\oo^\times$ et la
topologie $l$-adique sur $GL(L)$. Ainsi $\r_\oo$ se factorise à travers un quotient fini $(\bar D_\oo^\times
/ \varpi_\oo^\Zm)/ \bar K_\oo$, pour $\bar K_\oo$ un sous-groupe ouvert normal. Le revêtement
$\widetilde{\JC_{I^o,m}^{=h}}(s) \longto \JC_{I^o,m}^{=h}(s)$ et $\r_\oo$ définissent un $\bar
\Qm_l$-faisceau localement constant $\LC_{\r_\oo}$ sur $\JC_{I^o,m}^{=h}(s)$.

\marque Soit de même $\t_o$ une représentation irréductible admissible de $D_{o,h}^\times$ et $\FC_{\t_o}$ le
système local sur $\IC_{I^o,m}^{=h}$ associé à la restriction de $\t_o$ à $\DC_{o,h}^\times$ et au revêtement
d'Igusa de seconde espèce $\JC_{I^o,m}^{=h}(\oo) \longto \IC_{I^o,m}^{=h}$. Soit $\HC_{I^o,m}^{\oo,\tilde o}$
la $\Qm$-algèbre des fonctions localement constantes à support compact
$$f^\oo:(D_\Am^{\oo,o})^\times \times GL_{d-h}(F_o) \longto \Qm$$
qui sont $K_{I^o}^{\oo,o} \times K_{o,m}$-bi-invariantes par translations à droite et à gauche; le produit est donné
par le produit de convolution avec la mesure de Haar $d h^{\oo,o} \times d h_o^{et}$.
Une base de $\HC_{I^o,m}^{\oo,\tilde o}$ est donnée par les fonctions caractéristiques
$$1_{K_{I^o}^{\oo,o} g^{\oo,o} K_{I^o}^{\oo,o}} \times 1_{K_{o,m}g_o^{et} K_{o,m}}$$
des doubles classes $K_{I^o}^{\oo,o} g^{\oo,o} K_{I^o}^{\oo,o} \subset (D_\Am^{\oo,o})^\times$
(resp. $K_{o,m} g_o^{et} K_{o,m} \subset GL_{d-h}(F_o)$), où
$g^{\oo,o}$ (resp. $g_o^{et}$) décrit un système de représentants de ces doubles classes. Pour tout
$(g^{\oo,o},g_o^{et}) \in (D_\Am^{\oo,o})^\times \times GL_{d-h}(F_o)$, on a une correspondance de Hecke
$$\xymatrix{
  & {\DS \lim_{\genfrac{}{}{0pt}{}{\lefto}{J^o,n}}} \widetilde{\JC_{J^o,n}^{=h}}(s) \dlto_{c_1(s)} \drto^{c_2(s)} \\
\widetilde{\JC_{I^o,m}^{=h}}(s) & & \widetilde{\JC_{I^o,m}^{=h}}(s) \ar@{-->}[ll] }$$ Cette correspondance
agit sur chaque
$$H^i_{h,I^o,m,\r_\oo,\t_o}:=H_c^i(\IC_{I^o,m}^{=h} \otimes_{\k(o)} \bar \k(o), \LC_{\r_\oo} \otimes \FC_{\t_o})$$
et cette action dépend seulement des doubles classes
$$K_{I^o}^{\oo,o} g^{\oo,o} K_{I^o}^{\oo,o} \times K_{o,m}g_o^{et} K_{o,m}.$$

\marque On obtient ainsi une action de $\HC_{I^o,m}^{\oo,\tilde o}$ sur $H^i_{h,I^o,m,\r_\oo,\t_o}$. D'après
la remarque (\ref{rema-action}), $C_c(D_{o,h}^\times/\DC_{o,h}^\times)$ agit sur $H^i_{h,I^o,m,\r_\oo,\t_o}$
et son action commute à celle de $\HC_{I^o,m}^{\oo,\tilde o}$. Dans le groupe de Grothendieck de
$\HC_{I^o,m}^{\oo,\tilde o} \times C_c(D_{o,h}^\times/\DC_{o,h}^\times)$, on note
$$[H_{h,I^o,m,\r_\oo,\t_o}^*]:=\sum_{i=0}^{2d-2h} (-1)^{d-h-i} [H^i_{h,I^o,m,\r_\oo,\t_o}]$$
Avec $\bar K_\oo$ choisi comme ci-dessus, on note aussi
$$[\widetilde{H_{h,I^o,m,s}}^*]:=\sum_{i=0}^{2d-2h} (-1)^{d-h-i} [H_c^i(\widetilde{\JC_{I^o,m}^{=h}(s)} 
\otimes_{\k(o)} \bar \k(o), \bar \Qm_l)]$$
dans le groupe de Grothendieck de $(\bar D_\oo^\times / \varpi_\oo^\Zm) \times \HC_{I^o,m}^{\oo,\tilde o}
\times \widetilde{\NC_o}$. On notera bien que $\widetilde{\NC_o}/\DC_{o,h}^\times$ est isomorphe à $GL_h(F_o)
\times W_{F_o}$. Soit $W$ (resp. $V$) l'espace vectoriel sous-jacent à $\t_o$ (resp. $\r_\oo$) et soit $s
\geq 0$ tel que $\t_o$ est trivial sur $(1+\Pi_{o,h}^{s+1} \DC_{o,h})$, alors par définition des systèmes
locaux,
$$H_c^i(\IC_{I^o,m}^{=h} \otimes_{\k(o)} \bar \k(o), \LC_{\r_\oo} \otimes \FC_{\t_o})$$
est isomorphe à
$$(H_c^i(\widetilde{\JC_{I^o,m}^{=h}}(s),\bar \Qm_l) \otimes (V \otimes W))^{((\bar D_\oo^\times / \varpi_\oo^\Zm)/
\bar K_\oo) \times \DC_{o,h}^\times/(1+\Pi_{o,h}^{s+1} \DC_{o,h})}$$
où $\DC_{o,h}^\times/(1+\Pi_{o,h}^{s+1} \DC_{o,h})$ (resp. $(\bar D_\oo^\times / \varpi_\oo^\Zm)/\bar K_\oo$)
agit sur $W$ (resp. $v$) via $\t_o$ (resp. $\r_\oo$). En particulier, on a
\begin{multline*}
\tr (1_n \times 1_{K_{\Am,I^o}^{\oo,o} g^{\oo,o} K_{\Am,I^o}^{\oo,o}} \times 1_{K_{o,m} g_o^{et} K_{o,m}}, 
[H^*_{h,I^o,m,\r_\oo,\t_o}])= \\
\frac{1}{|\DC_{o,h}^\times:(1+\Pi_{o,h}^{s+1}\DC_{o,h})| |(\bar D_o^\times / \varpi_\oo^\Zm)/\bar K_\oo | } \\
\sum_{\d_o \in \DC_{o,h}^\times/(1+\Pi_{o,h}^{s+1} \DC_{o,h})} \sum_{\bar g_\oo \bar K_\oo \in (\bar D_\oo^\times /
\varpi_\oo^\Zm)/\bar K_\oo}  \tr(\r_\oo(\bar g_\oo)) \tr(\t_o(\Pi_{o,h}^{n} \d_o)) \\
\tr(\bar g_\oo \times (g_o^c,\Pi_{o,h}^{n}\d_o,\frob_o^r) \times 1_{K_{\Am,I^o}^{\oo,o} g^{\oo,o} K_{\Am,I^o}^{\oo,o}} 
\times
1_{K_{o,m} g_o^{et} K_{o,m}},[\widetilde{H_{h,I^o,m,s}}^*])
\end{multline*}
où $g_o^c \in GL_h(F_o)$ et $r \in \Zm$ sont tels que $(g_o^c,\Pi_{o,h}^n,\frob_o^r) \in \widetilde{\NC_o}$.
En appliquant la formule des traces de Lefschetz, on obtient la proposition suivante

\begin{prop} \label{prop-1-trace}
Avec les notations précédentes,
$$\tr (1_n \times 1_{K_{\Am,I^o}^{\oo,o} g^{\oo,o} K_{\Am,I^o}^{\oo,o}} \times 1_{K_{o,m} g_o^{et} 
K_{o,m}},[H^*_{h,I^o,m,\r_\oo,\t_o}])$$
est donné par la formule
\begin{multline*}
\frac{1}{|\DC_{o,h}^\times:(1+\Pi_{o,h}^{s+1}\DC_{o,h})| |(\bar D_\oo^\times / \varpi_\oo^\Zm)/\bar K_\oo| }
\sum_{\d_o \in \DC_{o,h}^\times/(1+\Pi_{o,h}^{s+1} \DC_{o,h})} \\
\sum_{\bar g_\oo \bar K_\oo \in (\bar D_\oo^\times / \varpi_\oo^\Zm)/\bar K_\oo} \lef_r^{=h}(\bar
g_\oo,\d_o,g_o^c,g^{\oo,\tilde o},s) \tr (\t_o(\Pi_{o,h}^{n}\d_o)) \tr (\r_\oo(\bar g_\oo))
\end{multline*}
soit en remplaçant $\lef_r^{=h}(\bar g_\oo,\d_o,g_o,g^{\oo,\tilde o},s)$ par sa valeur:
\begin{multline*}
\sum_{\gamma \in D^\times_{\natural_h}\cap C_h(\Pi_{o,h}^{\val(\det g_o^c)+r} \d_o^{-1})} \e_\oo(\gamma)
\vol(D_\gamma^\times \backslash (D_\Am)^\times_\gamma/(\varpi_\oo^\Zm), \frac{d h_{\Am,\gamma}}{d h_\gamma})
\\ O_{\bar \gamma} (\bar f_\oo, d \bar h_{\oo,\gamma}) O_{\gamma}(f^{\oo,o},d h_{\gamma}^{\oo,o})
O_{g_o^{et}}(f_o^{et},d h_{g_o^{et}}) \frac{\tr \t_o(\Pi_{o,h}^{-\a}\d_o)}{\vol(\DC_{o,h}^\times,d \bar h_o)}
\end{multline*}
où:

\noindent - les mesures de Haar $dh^\oo$, $dh_{\Am,\gamma}$, $dh^\oo_\gamma$, $dh_\gamma$, et $d \bar
h_{\oo,\gamma}$ sont choisis comme précédemment,

\noindent - $d \bar h_\oo$ est arbitraire avec
$$\bar f_\oo =\frac{\xi_{\r_\oo}}{\vol(\bar D_\oo^\times/ \varpi_\oo^\Zm, d \bar h_\oo)}$$
où $\xi_{\r_\oo}$ est le caractère de $\r_\oo$.
\end{prop}

\rem Le produit $\bar f_\oo d \bar h_\oo$ est indépendant du choix de $d \bar h_\oo$.

\subsection{Transfert des intégrales orbitales}

\begin{defi} Une représentation $\pi_x$ de $GL_d(F_x)$ est dite de carré intégrable (resp. essentiellement de carré 
intégrable)
si $\pi_x$ (resp. s'il existe un caractère $\psi_x$ tel que
$\psi_x \otimes \pi_x$) a un coefficient matriciel de carré intégrable sur $GL_d(F_x)/ F_x^\times$.
\end{defi}

\textit{A la place $\oo$}: On suit de très près le paragraphe (13.8) de \cite{lrs}. D'après \cite{badu}, il existe une 
unique
représentation irréductible $\pi_\oo$ de $D_\oo^\times / \varpi_\oo^\Zm$ essentiellement de carré intégrable qui est 
déterminée
par la relation suivante sur la restriction des caractères aux \eles elliptiques réguliers; si $\gamma \in 
D_\oo^\times / \varpi_\oo^\Zm$
est elliptique régulier\footnote{i.e. $F_\oo[\gamma]$ est une extension séparable de degré $d$ de $F_\oo$}
correspondant à un \ele $\bar \gamma \in \bar D_\oo^\times / \omega_\oo^\Zm$ par transfert, alors
$$\xi_{\pi_\oo}(\gamma)=(-1)^{d-1} \xi_{\r_\oo}(\bar \gamma)$$
D'après \cite{ze}, $\pi_\oo$ est de la forme $\st_{2t+1}(\pi_\oo')$, où $\pi_\oo'$ est une représentation
cuspidale de $GL_{d'}(F_\oo)$ avec $d=d' (2t+1)$ et $t \in \frac{1}{2} \Zm$, où l'on rappelle que
$\st_{2t+1}(\pi'_\oo)$ est l'unique sous-représentation irréductible de la représentation induite normalisée
(cf. ci-après)
$$\pi'_\oo(t) \times \cdots \times \pi'_\oo(-t)$$
En outre cette dernière induite possède un unique quotient irréductible que l'on note $\speh_{2t+1}(\pi'_\oo)$.

La proposition suivante est bien connue (cf. par exemple \cite{h-t} lemme I.3.4 (3))

\begin{prop} Il existe une fonction $f_{\pi_\oo}$ sur $D_\oo^\times / \varpi_\oo^\Zm$, localement constante, à support 
compact,
vérifiant les propriétés suivantes:
\begin{itemize}
\item[(i)] Les intégrales orbitales non elliptiques de $f_{\pi_\oo}$ sont nulles; pour $\gamma \in D_\oo^\times$ 
elliptique, on a
$$O_\gamma(f_{\pi_\oo},dh_{\oo,\gamma})=\e_\oo(\bar \gamma) O_{\bar \gamma}(\bar f_\oo, d \bar h_{\oo,\gamma})$$
où $\bar f_\oo$ et les mesures de Haar sont définies comme au paragraphe précédent.

\item[(ii)] Pour une représentation irréductible $\tilde \pi_\oo$ de $D_\oo^\times / \varpi_\oo^\Zm$, on a
$$\tr \tilde \pi_\oo(f_{\pi_\oo})= \left \{
\begin{array}{l} 1 \hbox{ si } \tilde \pi_\oo \simeq \pi_\oo=\st_{2t+1}(\pi'_\oo) \\
(-1)^{2t} \hbox{ si } \tilde \pi_\oo \simeq \speh_{2t+1}(\pi'_\oo) \\
0 \hbox{ sinon}
\end{array} \right. $$
\end{itemize}
\end{prop}

\textit{A la place $o$}: On reprend textuellement les résultats de \cite{h-t} VI.5 On rappelle que
$P_{h,d}\subset GL_d$ est le sous-groupe parabolique constitué des matrices triangulaires supérieures par
blocs de Levi $GL_{h} \times GL_{d-h}$. Soit alors $N_{h,d}$ son radical unipotent. On note $N_{h,d}^{op}$ le
radical unipotent du parabolique opposé $P_{h,d}^{op}$ de $P_{h,d}$. Pour $\pi_o$ une représentation
irréductible admissible de $GL_d(F_o)$ et $P$ un parabolique de radical unipotent $N$, le module de Jacquet
$J_N(\pi_o)$ est la représentation admissible $\pi_{o,N} \otimes \d_P^{1/2}$ du groupe $(P/N)(F_o)$ dont
l'espace sous-jacent est l'espace des $N(F_o)$-coinvariants de l'espace de $\pi_o$ et
$$\d_P(h)=|\det(\ad(h)_{|\lie N})|_{F_o}.$$
Si $\pi_o$ est une représentation admissible de $(P/N)(F_o)$, on note
$$n-\ind_{P(F_o)}^{GL_d(F_o)} (\pi_o):=\ind_{P(F_o)}^{GL_d(F_o)}(\pi_o \otimes \d_P^{1/2}).$$

Pour $\pi_o$ une représentation irréductible admissible essentiellement de carré intégrable de $GL_d(F_o)$,
de caractère central $\psi_{\pi_o}$,
Deligne, Kazhdan et Vigneras (cf. \cite{badu}) montrent l'existence d'une fonction $\phi_{\pi_o} \in 
C_c^\oo(GL_d(F_o),\psi_{\pi_o}^{-1})$,
que l'on appelle un pseudo-coefficient de $\pi_o$, vérifiant les propriétés suivantes:
\begin{itemize}
\item $\tr \pi_o(\phi_{\pi_o})=\vol(D_{o,d}^\times/F_o^\times)$ \footnote{où l'on considère sur $D_{o,d}^\times$ le 
transfert de la
 mesure de Haar sur $GL_d(F_o)$};

\item si $P$ est un sous-groupe parabolique de $GL_d$ de Levi $GL_{d_1} \times \cdots \times GL_{d_s}$, $s>1$. 
Supposons donnés des
représentations irréductibles admissibles essentiellement de carré intégrable $\pi_{o,i}$ de $GL_{d_i}(F_o)$ tels que
$\psi_{\pi_{o,1}} \cdots \psi_{\pi_{o,s}}=\psi_{\pi_o}$, alors
$$\tr (n-\ind_{P(F_o)}^{GL_d(F_o)} (\pi_{o,1} \times \pi_{o,s})) (\phi_{\pi_o})=0;$$

\item si $\gamma \in GL_d(F_o)$ est un \ele semi-simple non elliptique alors
$$O_\gamma^{GL_d(F_o)}(\phi_{\pi_o})=0$$
(cf. \cite{h-t} \S I.3)

\item si $\gamma \in GL_d(F_o)$ est elliptique semi-simple et si $\d \in D_{o,d}^\times$ a le même polynôme
caractéristique que $\gamma$ alors (cf. \cite{h-t} lemme I.3.1)
$$O_\gamma^{GL_d(F_o)}(\phi_{\pi_o})=(-1)^{d(1-[F_o(\gamma):F_o]^{-1})} \vol(D_{o,d}^\times/Z_{D_{o,d}^\times}(\d)) 
\tr \JL^{-1}(\pi_o^\vee)(\d)$$

\end{itemize}

Soient $d h_o^{et}$ et $d h_o$ des mesures de Haar sur respectivement $GL_{d-h}(F_o)$ et $GL_d(F_o)$.

\begin{prop} \label{transfert}
(cf. \cite{h-t} lemmes VI.5.1 et VI.5.2) Soient $\t_o$ une représentation irré\-ductible de $D_{o,h}^\times$ et
$\phi_o^{et} \in C_c^\oo(GL_{d-h}(F_o))$; il existe alors une fonction
$$IPC_{\t_o}(\phi_o^{et};d h_o^{et}) \in C_c^\oo(GL_d(F_o))$$ telle que

(1) $O_{\gamma}(IPC_{\t_o}(\phi_o^{et};d h_o^{et}))$ est nulle si $\gamma$ n'appartient pas à
$C_h(\Pi_{o,h}^{-\a} \d_o)$ et sinon est égal à

$$(-1)^{h-1}\\ O_{g_o^{et}}^{GL_{d-h}(F_o)}(\phi_o^{et},d h_{o,g_o^{et}}^{et})
 \frac{\tr \t_o(\Pi_{o,h}^{-\a}\d_o)}{\vol(\DC_{o,h}^\times,d \bar h_o)}$$

(2) En outre si $\pi_o$ est une représentation irréductible admissible de $GL_d(F_o)$ et si
$$[J_{N_{h,d}}(\pi_o) \otimes \d_{P_{h,d}}^{1/2}]= \sum_{\a,\b} m_{\a,\b}[\a \otimes \b]$$
où $\a$ (resp. $\b$) décrit l'ensemble des représentations irréductibles admissibles de $GL_h(F_o)$ (resp. 
$GL_{d-h}(F_o)$), alors
\begin{multline*}
\tr \pi_o(IPC_{\t_o}(\phi_o^{et}, d h_o^{et})= \\
\sum_{\a,\b,\psi} \tr \psi(\bar \phi_o) \frac{m_{\a,\b}}{\vol(D_{o,h}^\times/F_o^\times,d \bar h_o)^{-1}} \tr
\a (\phi_{\JL(\t_o \otimes \psi^{-1})}) \tr \b(\phi_o^{et})
\end{multline*}
où $\a$ (resp. $\b$, resp. $\psi$) décrit les représentations irréductibles admissibles de $GL_h(F_o)$ (resp.
$GL_{d-h}(F_o)$, resp. $D_{o,h}^\times/\DC_{o,h}^\times$); la somme sur $\psi$ porte sur les $\psi$ tels que
$\a$ et $\t_o \otimes \psi^{-1}$ ont le même caractère central.
\end{prop}

\begin{proof} Celle-ci est strictement similaire à celle de loc. cit. avec les modifications suivantes. En premier
lieu on note que $\gamma$ est associé à $\Pi_{o,h}^{\a} \d_o^{-1}$ alors que l'on considère $\tr
\t_o(\Pi_{o,h}^{-\a} \d_o)$; en outre $g_o^c$ induit la multiplication par $-\val(\det g_o^c)$ sur $\Zm$, à
comparer avec $\val(\det g_o^c)$ dans loc. cit.. On se retrouve alors avec $\JL(\t_o \otimes \psi^{-1})$ au
lieu de $\JL(\t_o^\vee \otimes \psi)$.
\end{proof}

\marque On introduit alors
$$\Red_{\t_o}^h: \groth(GL_d(F_o)) \longto \groth(D_{o,h}^\times/\DC_{o,h}^\times \times GL_{d-h}(F_o))$$
défini comme la composition des deux homomorphismes suivant.

\begin{itemize}
\item En premier lieu, on a un homomorphisme
$$\begin{array}{l}
\groth(GL_d(F_o)) \longto \groth(GL_h(F_o) \times GL_{d-h}(F_o)) \\
~ [ \pi_o ] \mapsto [ J_{N_{h,d}}(\pi_o) \otimes \d_{P_{h,d}}^{1/2} ]
\end{array}$$

\item Ensuite on a un homomorphisme
$$\begin{array}{l}
\groth(GL_h(F_o) \times GL_{d-h}(F_o)) \longto \groth(D_{o,h}^\times/\DC_{o,h}^\times \times GL_{d-h}(F_o)) \\
~ [\a \otimes \b] \mapsto \sum_{\psi} \vol(D_{o,h}^\times/F_o^\times,d \bar h_o)^{-1} \tr \a (\phi_{\JL(\t_o
\otimes \psi^{-1})})[\psi \otimes \b],
\end{array}$$
où $\psi$ décrit les caractères de $D_{o,h}^\times/\DC_{o,h}^\times$ tels que $\a$ et $\t_o \otimes
\psi^{-1}$ ont le même caractère central et où l'on considère des mesures de Haar associées sur $GL_h(F_o)$
et $D_{o,h}^\times$.
\end{itemize}

\begin{coro} Pour $\pi_o$ une représentation irréductible admissible de $GL_d(F_o)$, on a
$$\tr \pi_o(IPC_{\t_o}(\phi_o^{et},d h_o^{et}))= \tr \Red_{\t_o}^h(\pi_o)(\phi_o^{et})$$
où $\phi_o^{et} \in C_c^\oo(GL_{d-h}(F_o))$.
\end{coro}

D'après la proposition (\ref{prop-1-trace}) et en remarquant que $h[F'_\oo:F_\oo]=d[F'_o:F_o]$ avec les
notations précédentes, on obtient le résultat suivant.

\begin{theo} \label{theo-trace}
Avec les notations précédentes, on pose
$$f=f_{\pi_\oo} f^{\oo,o} IPC_{\t_o}(f_o^{et}).$$
On a alors
\begin{multline*}
\tr (1_n \times 1_{K_{\Am,I}^{\oo,o} g^{\oo,o} K_{\Am,I}^{\oo,o}} \times 1_{K_{o,m} g_o^{et} 
K_{o,m}},H^*_{h,I^o,m,\r_\oo,\t_o}) = \\
\sum_{\gamma \in D^\times_{\natural}} \vol(D_\gamma^\times \backslash (D_\Am)^\times_\gamma/(\varpi_\oo^\Zm),
\frac{d h_{\Am,\gamma}}{d h_\gamma}) O_\gamma(f,d h_{\Am,\gamma})
\end{multline*}
\end{theo}

\subsection{Formule des traces de Selberg}

Soit $\AC(D^\times \backslash D_\Am^\times/(\varpi_\oo^\Zm))$
l'espace des fonctions localement constantes muni de la représentation régulière
à droite de $D_\Am^\times/(\varpi_\oo^\Zm)$. Comme $D$ est une algèbre à division, $D^\times \backslash 
D_\Am^\times/(\varpi_\oo^\Zm)$ est
compact, de sorte que
$$\AC(D^\times \backslash D_\Am^\times/(\varpi_\oo^\Zm))= \bigoplus_{\Pi} m(\Pi) \Pi$$
avec $m(\Pi)$ fini et où $\Pi$ décrit l'ensemble des représentations irréductibles admissibles de $D_\Am^\times/ 
(\varpi_\oo^\Zm)$.
Si la multiplicité $m(\Pi)$ n'est pas nulle, la représentation $\Pi$ est dite automorphe. L'opérateur induit par une 
fonction localement
constante à support compact sur $D_\Am^\times / (\varpi_\oo^\Zm)$, a une trace:
$$\begin{array}{rl}
\tr (f;\AC(D^\times \backslash D_\Am^\times/(\varpi_\oo^\Zm)))& = \sum_\Pi m(\Pi) \tr \Pi(f)\\
& =\sum_{\gamma \in D^\times_\natural} \vol(D^\times_\gamma \backslash
(D_\Am)^\times_\gamma/(\varpi_\oo^\Zm),\frac{d h_{\Am,\gamma}}{d h_\gamma}) O_\gamma(f,d h_{\Am,\gamma})
\end{array}$$

\begin{prop} Pour tout $r>0$, on a l'égalité
$$\tr (1_n \times 1_{K_{\Am,I}^{\oo,o} g^{\oo,o} K_{\Am,I}^{\oo,o}} \times 1_{K_{o,m} g_o^{et} 
K_{o,m}},H^*_{h,I^o,m,\r_\oo,\t_o}) =
\sum_\Pi m(\Pi) \tr \Pi(f)$$
\end{prop}

\marque On note $[H^i_{h,\r_\oo,\t_o}]$ l'\ele du groupe de Grothendieck
$$\groth((D_\Am^{\oo,o})^\times GL_{d-h}(F_o) \times D_{o,h}^\times/\DC_{o,h}^\times)$$
défini par
la limite inductive sur $I^o$ et $m$ des $[H^i_{h,I^o,m,\r_\oo,\t_o}]$.
En vertu de la proposition précédente, on a le corollaire suivant.

\begin{coro} \label{strate-alt}
Soit $\Pi$ une représentation automorphe de $D_\Am^\times$, alors dans le groupe de Grothendieck
$\groth(GL_{d-h}(F_o) \times (D_{o,h}^\times /\DC_{o,h}^\times))$, on a
$$[H^*_{h,\r_\oo,\t_o}(\Pi^{\oo,o})]= \left \{
  \begin{array}{ll}
    m(\Pi) \Red_{\t_o}^h (\Pi_o) & \hbox{si } \Pi_\oo \simeq \st_{2t+1}(\pi'_\oo) \hbox{ ou }
 \speh_{2t+1}(\pi'_\oo) \\
0 & \hbox{sinon}
  \end{array}
\right .
$$
\end{coro}


%% file: somme-alternee.tex
\section{Application à la cohomologie du modèle local}

\subsection{Calcul de $\sum_i (-1)^i \widetilde{\Psi_{F_o}^{d,i}}(\JL^{-1}(\st_s(\pi_o)))$}

\label{calculht}

Dans \cite{boy}, pour une représenta\-tion cuspidale $\pi_o$ de $GL_d(F_o)$, on montre que
$$\widetilde{\UC_{F_o}^{d,i}}(\pi_o)= \left \{ \begin{array}{l}
\JL^{-1}(\pi_o) \otimes L_d(\pi_o) (\frac{1-d}{2}) \hbox{ pour } i=d-1 \\
0 \hbox{ si } i \neq d-1
\end{array} \right. $$
où $L_d$ (resp. $\JL$) est la bijection de Langlands (resp. de Jacquet-Langlands) de l'ensemble des classes
d'équivalences des représentations irréductibles essentiellement de carré intégrables de $GL_d(F_o)$ dans
l'ensemble des classes d'équivalences des représentations de dimension $d$ de $W_{F_o}$ (resp. des
représentations irréductibles admissibles de $D_{o,d}^\times$). Le but de ce paragraphe est d'obtenir le
pendant du théorème VII.1.5 de \cite{h-t}, à savoir.

\begin{theob} \label{theo-calcul-psi} (cf. \cite{h-t} théorème VII.1.5)
Soit $d=sg$ avec $s$ et $g$ des entiers positifs, et soit $\pi_o$ une représentation irréductible cuspidale
de $GL_g(F_o)$, alors
$$[\widetilde{\UC_{F_o}^{d,*}}(\JL^{-1}(\st_s(\pi_o)))]= \sum_{j=1}^s (-1)^{s-j}
[\overleftarrow{j-1},\overrightarrow{s-j}]_{\pi_o} \otimes L_g(\pi_o) (-\frac{d-s+2(j-1)}{2}).$$
\end{theob}

\begin{proof} Dans le groupe de Grothendieck de $(D_\Am^\oo)^\times \times W_{F_o}$, on note
$$H^i_{\eta_o,\r_\oo}=\sum_{\Pi^\oo} \Pi^\oo \otimes W_{\r_\oo,i}(\Pi^\oo)$$
$$[H_{\eta_o,\r_\oo}^*]:= \sum_{i=0}^{2d-2} (-1)^{d-1-i} [\lim_{\genfrac..{0pt}{1}{\lefto}{I}} H^i(M_{I,\eta_o} 
,\LC_{\r_\oo})]=
\sum_{\Pi^\oo} [\Pi^\oo][W_{\r_\oo}^*(\Pi^\oo)]$$ où la somme porte sur les représentations irréductibles
automorphes de $(D_\Am^\oo)^\times$. On rappelle alors l'un des résultats principaux de \cite{lrs}.

\begin{theob} \label{theo-lrs}
Si $\Pi$ est une représentation automorphe de $D_\Am^\times$ telle que $\Pi_\oo \simeq \st_\oo$ et
telle qu'il existe deux places $x_1,x_2$ de $X'$ distinctes de $\oo$ et $o$ avec $\Pi_{x_i}$ cuspidales pour $i=1,2$, 
alors
$$[W_{1_\oo}^*(\Pi^\oo)_{|W_o}] \simeq [L_d(\Pi_o)] (\frac{1-d}{2})$$
où $L_d$ désigne la correspondance de Langlands. \footnote{On a même la nullité des $W_{1_\oo,i}(\Pi^\oo)$
pour $i \neq d-1$ et alors l'égalité en tant que représentations et pas seulement des semi-simplifiées.}
\end{theob}

\marque La suite spectrale des cycles évanescents donne alors l'égalité
$$[H_{\eta_o,\r_\oo}^*]=\sum_{0 \leq j \leq d-1} (-1)^j [\lim_{\genfrac..{0pt}{1}{\lefto}{I}} H^*(M_{I,\bar s_o},R^j 
\Psi_{\eta_o}(\LC_{\r_\oo}))]$$
Cette égalité combinée à la suite spectrale associée à la stratification donne l'égalité suivante
$$[H_{\eta_o,\r_\oo}^*]=\sum_{\genfrac{}{}{0pt}{}{0 \leq j \leq d-1}{1 \leq h \leq d}} (-1)^j 
[\lim_{\genfrac..{0pt}{1}{\lefto}{I}} H^*_c(M_{I,\bar s_o}^{=h}, R^j \Psi_{\eta_o}(\LC_{\r_\oo}))]$$
D'après \cite{boy}, pour $\FC$ un faisceau sur $M_{I,s_o}^{=h}$, on a
$$\lim_{\genfrac..{0pt}{1}{\lefto}{I}} H^i_c(M_{I,\bar s_o}^{=h},\FC)= \ind_{P_{h,d}^{op}(F_o)}^{GL_d(F_o)} 
\lim_{\genfrac..{0pt}{1}{\lefto}{I}} H^i_c(M_{I,\bar s_o,1}^{=h},\FC)$$
On obtient ainsi l'égalité
\begin{eqnarray} \label{egalite}
[H_{\eta_o,\r_\oo}^*]=\sum_{\genfrac{}{}{0pt}{}{0 \leq j \leq d-1}{1 \leq h \leq d}} (-1)^j
[\ind_{P_{h,d}^{op}(F_o)}^{GL_d(F_o)} \lim_{\genfrac..{0pt}{1}{\lefto}{I}} H_c^*(M_{I,\bar s_o,1}^{=h}, R^j
\Psi_{\eta_o}(\LC_{\r_\oo}))]
\end{eqnarray}
D'après le corollaire (\ref{corofond1}) et le lemme (\ref{lem-lisse}) du chapitre 3, on a
\begin{equation} \label{egalite2}
h[H^i_c(M_{I,\bar s_o,1}^{=h},R^j \Psi_{\eta_o}(\LC_{\r_\oo}))]=\bigoplus_{\t_o \in \CF_h} \frac{h}{e_{\t_o}}
{[} H^i_c(M_{I,\bar s_o,1}^{=h},\LC_{\r_\oo} \otimes \FC_{\t_o}) {]} *_{d_h}
[\widetilde{\UC_{F_o,n}^{h,j}}(\t_o)]
\end{equation}
où $n$ est la multiplicité de $o$ dans $I$, où
$$d_h: \left \{ \begin{array}{l}
GL_h(F_o) \times W_{F_o} \longto (D_\Am^{\oo,o})^\times \times GL_{d-h}(F_o) \times D_{o,h}^\times/\DC_{o,h}^\times 
\\
(g_o^c,\s_o) \longmapsto (1,1,\d_o)
\end{array} \right .$$
où $\d_o$ est tel que $\val(\rn(\d_o))=-\val(\det(g_o^c)) - \deg(\s_o)$ et où pour $\pi_i$ une représentation
de $G_i$ pour $i=1,2$ et $d:G_2 \longto Z(G_1)$, on note $\pi_1 \otimes_d \pi_2$ la représentation de $G_1
\times G_2$ définie par
$$(\pi_1 \otimes_d \pi_2) (g_1,g_2)=\pi_1(g_1 d(g_2)) \otimes \pi_2(g_2);$$
dans le groupe de Grothendieck de $G_1 \times G_2$, on note $[\pi_1]*_d [\pi_2]$ l'image de $[\pi_1 \otimes_d 
\pi_2]$.
Cette égalité est vue dans
$$\groth((D_\Am^{\oo,o})^\times \times GL_{d-h}(F_o) \times GL_h(F_o) \times W_{F_o})$$
via l'application évidente
$$\groth((D_\Am^{\oo,o})^\times \times GL_{d-h}(F_o) \times D_{o,h}^{\times}/\DC_{o,h}^\times) \longto 
\groth((D_\Am^{\oo,o})^\times \times GL_{d-h}(F_o))$$

Soient alors $g$ un diviseur de $d=sg$ et $\pi_o$ une représentation irréductible cuspidale de $GL_g(F_o)$.
On considère $\Pi$ une représentation irréductible de $D_\Am^\times$ vérifiant les conditions du théorème
(\ref{theo-lrs}) et telle que $\Pi_o \simeq \st_s(\pi_o)$ \footnote{pour l'existence d'une telle
globalisation cf. \cite{he2}}.

\begin{lemm} On a
$$\Red_{\JL^{-1}(\st_s(\pi_o))}([\overleftarrow{s-1}]_{\pi_o})=e_{\JL^{-1}(\st_s(\pi_o))}
[\overleftarrow{s-l-1}]_{\pi_o(\frac{l(g-1)}{2})} \otimes \Xi^{\frac{(s-l)(g-1)}{2}}$$ où $\Xi$ est le
caractère multiplicatif de $\Zm$ tel que $\Xi(1)=\frac{1}{p}$.
\end{lemm}

\begin{proof} Le résultat découle directement de \cite{ze} 2.2. En effet on a
$$J_{N_{lg,d}}([\overleftarrow{s-1}]_{\pi_o}=[\overleftarrow{l-1}]_{\pi_o(\frac{(s-l)}{2})} \otimes
[\overleftarrow{s-l-1}]_{\pi_o(\frac{-l}{2})}$$ de sorte qu'après multiplication par $\d_{P_{lg,d}}^{1/2}$ on
obtient
$$[\overleftarrow{l-1}]_{\pi_o(-\frac{(s-l)(g-1)}{2})} \otimes [\overleftarrow{s-l-1}]_{\pi_o(\frac{l(g-1)}{2})}$$
d'où le résultat.

\end{proof}

Ainsi en combinant (\ref{egalite}), (\ref{egalite2}) avec les lemmes précédents, on obtient
\begin{multline} \label{egalite3}
[\overleftarrow{s-1}]_{\pi_o} \otimes [L_g(\pi_o) \otimes \sp_s]= \\ \sum_{l=1}^s
[\widetilde{\UC_{F_o}^{lg,*}}(\JL^{-1}(\st_l(\pi_o)))](-\frac{(s-l)(g-1)}{2}) \overrightarrow{\times}
[\overleftarrow{s-l-1}]_{\pi_o}
\end{multline}
où $(-\frac{(s-l)(g-1)}{2})$ est la torsion sur la partie galoisienne \footnote{Celle sur $GL_{lg}(F_o)$ est
contenu dans le symbole $\overrightarrow{\times}$.} et $[\sp_s]= 1(\frac{1-s}{2})+ \cdots +
1(\frac{s-1}{2})$. On suppose alors, par récurrence sur $s$, que pour tout $1 \leq l < s$, on ait
$$[\widetilde{\UC_{F_o}^{lg,*}}(\JL^{-1}(\st_l(\pi_o)))]=\sum_{r=0}^{l-1} (-1)^r
[\overleftarrow{l-1-r},\overrightarrow{r}]_{\pi_o} \otimes L_g(\pi_o)(-\frac{l(g-1)+2(l-r-1)}{2})$$ On
réinjecte ces égalités dans (\ref{egalite3}) ce qui donne en rassemblant selon les poids:
\begin{multline*}
[\widetilde{\UC_{F_o}^{d,*}}(\JL^{-1}(\st_s(\pi_o)))]=\sum_{l=0}^{s-1} L_g(\pi_o)(-\frac{s(g-1)+2l}{2}) \\
[\overleftarrow{s-1}]_{\pi_o} - \sum_{k=0}^{s-l-1} (-1)^k [\overleftarrow{l-1},\overrightarrow{k}]_{\pi_o}
\overrightarrow{\times} [\overleftarrow{s-l-k-1}]_{\pi_o}
\end{multline*}
d'où le résultat en utilisant le lemme suivant.

\end{proof}

\begin{lemm} \label{lem-combi-ze}
Pour tout $0<l<s$, on a l'égalité suivante dans le groupe de Grothendieck des représentations admissibles de
$GL_d(F_o)$:
$$\sum_{r=0}^{s-l-1} (-1)^r [\overleftarrow{l-1},\overrightarrow{r}]_{\pi_o} \overrightarrow{\times}
[\overleftarrow{s-l-r-1}]_{\pi_o}=[\overleftarrow{s-1}]_{\pi_o} + (-1)^{s-l-1}
[\overleftarrow{l-1},\overrightarrow{s-l}]_{\pi_o}$$
\end{lemm}

\begin{proof} Soit $a_k=\sum_{r=0}^{k} (-1)^r [\overleftarrow{l-1},\overrightarrow{r}]_{\pi_o}
\overrightarrow{\times} [\overleftarrow{s-l-r-1}]_{\pi_o}$ et montrons par récurrence sur $k$ de $0$ à
$s-l-1$ que $a_k=[\overleftarrow{s-1}]_{\pi_o} + (-1)^{k}
[\overleftarrow{l-1},\overrightarrow{k+1},\overleftarrow{s-l-k-1}]_{\pi_o}$. Le résultat est clairement
vérifié pour $k=0$ supposons le vrai au rang $k-1$ et traitons le cas de $k$ soit
$$\begin{array}{ll} a_k = & a_{k-1}+ (-1)^k [\overleftarrow{l-1},\overrightarrow{k}]_{\pi_o}
\overrightarrow{\times}
[\overleftarrow{s-l-k-1}]_{\pi_o} \\
= & [\overleftarrow{s-1}]_{\pi_o} + (-1)^{k-1}(
[\overleftarrow{l-1},\overrightarrow{k},\overleftarrow{s-l-k}]_{\pi_o} - \\ &
[\overleftarrow{l-1},\overrightarrow{k},\overleftarrow{s-l-k}]_{\pi_o}-[\overleftarrow{l-1},\overrightarrow{k+1},
\overleftarrow{s-l-k-1}]_{\pi_o}) \\
= & [\overleftarrow{s-1}]_{\pi_o} + (-1)^{k}
[\overleftarrow{l-1},\overrightarrow{k+1},\overleftarrow{s-l-k-1}]_{\pi_o}
\end{array}$$

\end{proof}

\subsection{Retour sur le cas Iwahori}

D'après le corollaire (\ref{dim-poids}), pour tout $0 \leq i < d$ et toute représentation irréductible
admissible $\t_o$ de $D_{o,d}^\times$, $(\Psi_{F_o}^{d,i})^{\Iw_o}$ est pur de poids $2i$. Ainsi d'après
(\ref{theo-calcul-psi}), il est égal à $[\overleftarrow{i},\overrightarrow{d-i-1}]_{1_o} \otimes |\cl|^{-i}$,
ce qui redonne bien le théorème (\ref{theo-iwahori}).


%% file: introduction2.tex
\noindent \textbf{0.1.} --- Soit $K$ un corps local complet d'égale caractéristique $p$, d'anneau des entiers
$\OC_K$. Pour un entier $d$ strictement positif fixé, on introduit le groupe $D_{K,d}^\times$ (resp. $W_K$)
des éléments inversibles de ``l'''algèbre à division centrale sur $K$ d'invariant $1/d$ (resp. le groupe de
Weil de $K$). Pour un nombre premier $l \neq p$, Langlands (resp. Jacquet-Langlands) a (resp. ont) conjecturé
l'existence d'une bijection $L_d$ (resp. d'une injection $\JL$) entre les $\bar \Qm_l$-représentations
irréductibles admissibles de $GL_d(K)$ et les représentations $l$-adiques indécomposables de $W_K$ (resp.
entre les représentations admissibles irréductibles de $D_{K,d}^\times$ et les représentations
essentiellement de carré intégrable de $GL_d(K)$) qui sont compatibles à la formation des fonctions $L$ de
paires \footnote{Pour un énoncé précis, cf. \cite{he}}.

A l'aide de la cohomologie étale Deligne a construit une série de représentations $\UC_{K}^{d,i}$ du produit
de ces trois groupes. Pour $d=2$ et $\rho$ une représentation irréductible admissible de $D_{K,d}^\times$
telle que $\pi:=\JL(\rho)$ est une représentation cuspidale de $GL_d(K)$, Carayol, dans \cite{ca}, montre que
la composante $\rho$-isotypique $\UC_{K}^{2,1}(\rho)$ de $\UC_{K}^{2,1}$ réalise les correspondances de
Langlands et de Jacquet-Langlands, i.e.
$$\UC_K^{2,1}(\JL^{-1}(\pi)^\vee) \simeq \pi \otimes L_d(\pi)^\vee (-\frac{d-1}{2})$$
Dans \cite{boy}, en égale caractéristique $p$, je traite le cas $d$ quelconque.

En outre pour $d=2$, Carayol décrit également les $\UC_K^{2,1}(\rho)$ pour $\rho$ quelconque. Le but premier
de ce travail est de faire de même pour $d$ quelconque, i.e. calculer les $\UC_{K}^{d,i}(\rho)$ pour $\rho$
une représentation irréductible admissible de $D_{K,d}^\times$. Dans le cas où $\rho$ est la représentation
triviale, rappelons que d'après le théorème (\ref{theo-iwahori}), on a

\medskip

\noindent \textbf{Théorème 1} \textit{Pour $0 \leq i \leq d-1$, on a
$$\UC_{K}^{d,i}(1)=\pi_{i} \otimes 1(-i)$$
où $\pi_{i}$ est l'unique quotient irréductible de l'induite parabolique
$$\ind_{P_{d-i,d}(K)}^{GL_{d}(K)} 1 \otimes \st_{i}$$
où $P_{d-i,d}$ est le parabolique standard associée aux $d-i$ premiers vecteurs et $\st_i$ est la
représentation de Steinberg de $GL_i(K)$.}

\medskip

L'énoncé du cas général, théorème (\ref{theo-ripsi-local}), s'énonce de manière similaire en faisant
intervenir les correspondances de Langlands et Jacquet-Langlands. Une autre formulation du résultat revient à
dire qu'il n'y a pas d'annulation dans l'expression, cf. le théorème (\ref{theo-calcul-psi}), de la
représentation virtuelle $\sum_{i=0}^{d-1} (-1)^i [\UC_{K}^{d,i}(\rho)]$ où $\UC_{K}^{d,d-i}(\rho)$, pour $1
\leq i \leq d$, y sera alors donné par le $i$-ème terme de plus haut poids.

\medskip

\noindent \textbf{0.2.} --- La preuve du théorème 1 dans le cas général procède par globalisation via l'étude
des variétés de Drinfeld-Stuhler et le théorème de comparaison de Berkovich\footnote{en fait sur une version
raffinée fournie par Fargues, cf. le théorème de l'appendice de \cite{boy-duke}} des cycles évanescents
locaux et globaux.

Soit donc $X$ une courbe projective lisse, irréductible et géométriquement connexe définie sur le corps fini
à $q=p^r$ \eles $\Fm_q$, et soit $F$ son corps des fonctions. On fixe deux places distinctes $\oo$ et $o$ de
$X$ que l'on peut supposer par simplification, rationnelles sur $\Fm_q$, de sorte que le complété $F_o$ du
localisé en $o$ de $F$ est isomorphe au corps local précédemment noté $K$. On note $A$ l'anneau des fonctions
sur $X$, régulières en dehors de $\oo$. Étant donné un entier $d \geq 1$, on fixe une algèbre à division
centrale $D$ sur $F$ de dimension $d^2$, non ramifiée en $\oo$ et $o$, ainsi qu'un ordre maximal $\DC$. Dans
\cite{lrs}, les auteurs construisent pour un idéal non trivial $I$ de $A$, un schéma $M_{I}$ défini sur $F$,
classifiant les $\DC$-faisceaux elliptiques sur $X$, munis d'une structure de niveau $I$. Pour $o \not \in
V(I)$, $M_I$ a un modèle entier $M_{I,o}$ lisse sur le complété $\OC_o$ de $A$ en la place $o$. Un tel modèle
non lisse dans le cas où $o \in V(I)$ est construit dans \cite{boy}. Les schémas $M_{I,o}$ sont naturellement
munis d'une action, par correspondances, de $(D_\Am^{\oo})^\times$.

\medskip

\noindent \textbf{0.3.} --- On s'intéresse alors à la fibre spéciale $M_{I,s_o}$ de $M_{I,o}$. Dans
\cite{boy}, je stratifie $M_{I,s_o}$ par des sous-schémas localement fermés $M_{I,s_o}^{=h}$ pour $1 \leq h
\leq d$, de pure dimension $d-h$ tels que l'on ait un équivalent du théorème de Serre-Tate pour les
$\DC$-faisceaux elliptiques à savoir: le complété de l'hensélisé strict de l'anneau local de $M_{I,o}$ en un
point géométrique de $M_{I,s_o}^{=h}$ est isomorphe à $\Def_n^h[[x_1,\cdots,x_{d-h}]]$ où $n$ est la
multiplicité de $o$ dans $I$ et $\Def_n^h$ représente le foncteur des déformations de niveau $n$ d'un
$\OC_o$-module formel de hauteur $h$ sur $\bar \Fm_p$. Par ailleurs pour $1 \leq h < d$, il existe un
sous-schéma fermé $M_{I,s_o,1}^{=h}$ de $M_{I,s_o}^{=h}$ stable sous les correspondances associées aux \eles
du sous-groupe parabolique $P_{h,d}^{op}(F_o)$ de $GL_d(F_o)$ (cf. la définition (\ref{defi-parab})) et tel
que
$$M_{I,s_o}^{=h}= M_{I,s_o,1}^{=h} \times_{P_{h,d}^{op}(\OC_o/\MC_o^n)} GL_d(\OC_o/\MC_o^n)$$
où $n$ est la multiplicité de $o$ dans $I$: on dit que les strates non supersingulières sont géométriquement
induites.

\medskip

\noindent \textbf{0.4.} --- Suivant \cite{h-t}, on introduit sur chacune des strates $M_{I,s_o,1}^{=h}$, des
systèmes locaux $\FC_{\tau_o}$ associés aux représentations $\tau_o$ du groupe des inversibles
$D_{o,h}^\times$ de l'algèbre à division centrale sur $F_o$ d'invariant $1/h$. On décrit alors la restriction
des cycles évanescents $R^i\Psi_{\eta_o}(\bar \Qm_l)$ à la strate $M_{I,s_o}^{=h}$ en fonction des
$\FC_{\t_o}$ et des cycles évanescents locaux $\Psi_{F_o,n}^{h,i}$, cf. (\ref{iso1}). D'après le théorème de
comparaison de Berkovich, le théorème local se déduit de la connaissance de la fibre en un point
supersingulier des $R^i\Psi_{\eta_o}(\bar \Qm_l)$.

\medskip

\noindent \textbf{0.5.} --- Le complexe $R\Psi_{\eta_o}(\bar \Qm_l)[d-1]$ est vu comme un faisceau pervers
muni d'une filtration de monodromie dont on notera $gr_k$ les gradués. Le deuxième résultat concerne la
description des gradués $gr_k$ dans la catégorie des faisceaux pervers sur $M_{I,s_o}$ munis d'une action
compatible de $(D_\Am^\oo)^\times \times W_o$. Pour $1 \leq lg \leq d$, et pour $\pi_o$ (resp. $\Pi_l$) une
représentation irréductible cuspidale de $GL_g(F_o)$ (resp. quelconque de $GL_{lg}(F_o)$), on introduit le
faisceau $HT(g,l,\pi_o,\Pi_l)$ sur la strate $M_{I,s_o}^{=lg}$, "induit" à partir du système local
$\FC_{\JL^{-1}(\st_l(\pi_o))} \otimes \Pi_l$ sur la composante $M_{I,s_o,1}^{=lg}$.

Les composantes $\pi_o$-isotypiques $gr_{k,\pi_o}$ des $gr_k$, cf. la proposition (\ref{prop-so}), se
décrivent alors au moyen des faisceaux pervers $j^{\geq lg}_{!*} HT(g,l,\pi_o,\st_l(\pi_o))[d-lg]$, cf. le
théorème (\ref{theo-global1}), où $j^{\geq lg}$ désigne l'injection de la strate $M_{I,s_o}^{=lg}$; en ce qui
concerne la partie associée à $\pi_o$ triviale l'énoncé est le suivant:

\medskip

\noindent\textbf{Théorème 2} \textit{Les faisceaux pervers $gr_{k,1_o}$ sont nuls pour $k \geq d$ et sinon on
a}
$$gr_{k,1_o,}=\bigoplus_{\genfrac{}{}{0pt}{}{|k| < l \leq  d}{l \equiv k-1 \mod 2}} j^{\geq l}_{!*} HT(1,l,1_o,\st_l)
(-\frac{lg+k-1}{2})[d-l]$$

\medskip

\noindent L'énoncé pour $\pi_o$ quelconque est similaire et fait intervenir les correspondances de Langlands
et Jacquet-Langlands.

\medskip

\noindent \textbf{0.6.} --- En utilisant le théorème de monodromie-poids, la preuve du théorème 2 découle de
la connaissance du théorème 1 pour toutes les hauteurs $h < d$ ainsi que de la description des restrictions
aux strates des faisceaux des cycles évanescents en fonction des systèmes locaux $\FC_{\tau_o}$ comme rappelé
en (0.4.). Par ailleurs le théorème 1 en hauteur $d$, découle d'après le théorème de comparaison de
Berkovich-Fargues, du calcul des germes aux points supersinguliers des faisceaux de cohomologie des $gr_k$.

On raisonne alors par récurrence en supposant connus\footnote{En fait on suppose plutôt connu les gradués
locaux $gr_{h,k,loc}$ de la filtration de monodromie-locale du complexe des cycles évanescents
$\Psi_{F_o}^{h,\bullet}$.} les $\UC_{F_o}^{h,i}$ du modèle de Deligne-Carayol de hauteur $h$ pour tout $1
\leq h <d$. On en déduit alors le théorème 2. Par ailleurs on sait déterminer tous les faisceaux de
cohomologie $h^i gr_k$ des $gr_k$ en dehors des points supersinguliers. La technique repose sur l'étude de la
suite spectrale associée à la filtration de monodromie:
\begin{equation} \label{ss-intro}
E_1^{i,j}=h^{i+j} gr_{-i} \Rightarrow R^{i+j+d-1} \Psi_{\eta_o}(\bar{{\mathbb Q}_l})
\end{equation}
en essayant d'en deviner les termes initiaux, l'aboutissement étant connu. En outre en utilisant la
perversité des $gr_k$ ainsi que la compatibilité de $R\Psi_{\eta_o}$ à la dualité de Verdier, on obtient un
contrôle sur les germes aux points supersinguliers des $h^i gr_k$. On étudie ensuite la suite spectrale des
cycles évanescents dont à nouveau on essaie de deviner les termes initiaux alors que l'aboutissement est
connu, en utilisant en particulier le théorème de Lefschetz difficile. Le contrôle obtenu précédemment nous
permet alors de prouver le théorème 1 et on explicite la suite spectrale de monodromie (\ref{ss-intro}).

\medskip

\noindent \textbf{0.7.} --- En ce qui concerne les résultats globaux que l'on obtient, citons

\begin{itemize}
\item la description ``explicite'' des gradués pour la filtration de monodromie du faisceau pervers des cycles 
évanescents
ainsi que de la suite spectrale associée;

\item la détermination des extensions intermédiaires des systèmes locaux d'Harris-Taylor;

\item le calcul de tous les groupes de cohomologie des différents faisceaux ou complexes de faisceaux qui sont 
introduits.
\end{itemize}

\medskip

\noindent \textbf{0.8.} --- Les résultats obtenus dans le cadre des variétés de Drinfeld s'adaptent aussi en
caractéristique mixte dans le cadre des variétés de Shimura de type PEL étudiées dans \cite{h-t}; c'est ce
travail qui est effectué dans \cite{boy-duke}. Le fait est que tous les arguments reposent de façon formelle
sur les propriétés géométriques du paragraphe (\ref{rappels-globaux}) et les propriétés cohomologiques du
paragraphe (\ref{ht-prop}); nous avons ainsi donné dans l'appendice A, le dictionnaire entre nos notations et
celles de \cite{h-t} en caractéristique mixte et rappelé où se trouvaient dans loc. cit. les propriétés
cohomologiques que l'on utilise dans le paragraphe (\ref{ht-prop}). La différence essentielle réside dans le
fait qu'en caractéristique mixte, nous ne disposons pas à priori de monodromie-poids. Ainsi dans le texte
nous avons, chaque fois que nécessaire, indiqué des preuves qui n'utilisent pas cette propriété. Au final
dans la situation de \cite{h-t}, on obtient alors une preuve de la conjecture de monodromie-poids versions
faisceautique et cohomologique.

\medskip

\noindent \textbf{0.9.} --- Comme conséquence directe de ces calculs on obtient une correspondance de
Jacquet-Langlands globale entre les algèbres à divisions $D$ et $\bar D$, résultat qui d'après I. Badulescu,
peut s'obtenir aisément à partir de la formule des traces simples. En outre on montre que les composantes
locales des représentations automorphes de $D_\Am^\times$ sont celles prévues par l'existence, en général
encore conjecturale, d'une correspondance de Jacquet-Langlands globale avec $GL_d(\Am)$ à partir des
résultats de Moeglin-Waldspurger sur les représentations de carré intégrable modulo le centre de $GL_d(\Am)$
et la conjecture de Ramanujan-Peterson prouvée par Lafforgue. Ainsi dans les cas défavorables où il n'existe
pas encore une telle correspondance de Jacquet-Langlands, on obtient la description des composantes locales
des représentations automorphes de $D_\Am^\times$ vérifiant $\hyp(\oo)$ et donc en particulier la conjecture
de Ramanujan-Peterson.

\bigskip

\noindent \textbf{0.10.} --- Décrivons succinctement le contenu des divers paragraphes. On commence par des
rappels sur les données géométriques locales et globales tirées  des chapitres précédents dans le cas d'égale
caractéristique, et de \cite{h-t} dans celui de la caractéristique mixte, cf. l'appendice \ref{appendiceB}.

En ce qui concerne les données globales, outre le fait que les strates non supersingulières soient
géométriquement induites, la donnée fondamentale est celle des systèmes locaux d'Harris-Taylor
$\FC(g,l,\pi_o)$ sur la strate $lg$, attachés aux représentations $\JL^{-1}(\st_l(\pi_o))$ des inversibles
$D_{o,lg}^\times$ de l'algèbre à division centrale sur $F_o$ d'invariant $1/lg$. La propriété essentielle que
l'on utilisera sur ces systèmes locaux est l'isomorphisme (\ref{iso1}).

On redonne en outre, d'après \cite{boy}, la description de l'ensemble des points supersinguliers et on
définit pour tout diviseur $g$ de $d=sg$, le faisceau $\FC(g,s,\pi_o)$ à support sur les points
supersinguliers. On notera par ailleurs que les $\FC(g,l,\pi_o)$ ne sont pas à priori irréductibles car ils
proviennent de la restriction à $\DC_{o,lg}^\times$ de $\JL^{-1}(\st_l(\pi_o))$.

\medskip

\noindent \textbf{0.11.} --- Au deuxième paragraphe, on explicite le lien entre $\Psi_{F_o}^{h,i}$ et
$\UC_{F_o}^{h,i}$ et on introduit l'entier $e_{\pi_o}$ qui est le cardinal de la classe d'équivalence
inertielle, définition (\ref{defi-inert}), de $\pi_o$. En particulier $e_{\pi_o}$ est égal au nombre de
sous-représentations irréductibles de la restriction à $\DC_{o,lg}^\times$ de $\JL^{-1}(\st_l(\pi_o))$.

On rappelle alors le théorème (3.2.4) de \cite{boy} qui décrit les parties cuspidales des groupes de
cohomologie des modèles locaux de Deligne-Carayol. A ce propos, pour avoir un énoncé exact dans loc. cit., il
faut considérer l'action naturelle de $GL_d(F_o)$ sur $\UC_{F_o}^{d,i}$, tordue par $g_o \mapsto \lexp t
g_o^{-1}$; on introduira un tilde pour marquer cette modification. On rappelle alors le théorème (VII.1.5) de
\cite{h-t} qui donne le calcul de la somme alternée, dans le groupe de Grothendieck des $GL_d(F_o) \times
W_o$-modules, $\sum_i (-1)^i \widetilde{\UC_{F_o}^{d,i}}(\JL^{-1}(\st_s(\pi_o)))$.

En global l'isomorphisme (\ref{iso1}) se traduit cohomologiquement par la proposition (\ref{prop-hic}). On
introduit selon \cite{lrs}, pour toute représentation $\rho_\oo$ du groupe des inversibles $\bar
D_\oo^\times$ de ``l'''algèbre à division centrale sur $F_\oo$ d'invariant $-1/d$, le système local
$\LC_{\rho_\oo}$ sur les variétés $M_{I,o}$. On considère alors les groupes de cohomologie du produit
tensoriel de ce dernier avec un système local d'Harris-Taylor. Le résultat fondamental qui résulte des
arguments de comptage de points est la proposition (\ref{prop-somme-alternee}) qui calcule la somme alternée
des groupes de cohomologie à support compact des systèmes locaux d'Harris-Taylor, dans le groupe de
Grothendieck des $(D_\Am^\oo)^\times \times \Zm$-modules, où pour tout $h$, $\Zm$ est identifié au quotient
$D_{o,h}^\times/ \DC_{o,h}^\times$ via la valuation de la norme réduite. On introduit pour cela le caractère
$\Xi$ de $\Zm \longto \bar \Qm_l^\times$ défini par $\Xi(1)=\frac{1}{p}$.

Ainsi les représentations automorphes $\Pi$ qui interviennent dans cette description ainsi que dans celle de
la cohomologie de la fibre générique à valeurs dans $\LC_{\rho_\oo}$, vérifient des conditions $\hyp(\oo)$ à
la place $\oo$ que l'on donne à la définition (\ref{defi-hyp}): en résumé pour
$\rho_\oo=\JL^{-1}(\st_s(\pi_\oo))$, il faut que $\Pi_\oo$ soit isomorphe à $\st_s(\pi_\oo)$ ou à
$\speh_s(\pi_\oo)$.

\medskip

\noindent \textbf{0.12.} --- Au troisième paragraphe le théorème (\ref{theo-ripsi-local}) est la motivation
initiale de ce travail, i.e. décrire chacun des $\widetilde{\UC_{F_o}^{d,i}}(\JL^{-1}(\st_s(\pi_o)))$.
Finalement on obtient un résultat plus précis, théorème (\ref{theo-local-fil}), qui est la description des
gradués de la filtration de monodromie-locale définie par Fargues dans l'appendice de \cite{boy-duke} et de
la suite spectrale correspondante. Le fait primordial démontré dans par Laurent Fargues est le théorème de
comparaison à la Berkovich qui implique qu'en tout point géométrique de $M_{I,s_o}^{=h}$, la filtration et la
suite spectrale correspondante, induite par la monodromie globale sur le germe en ce point du complexe des
cycles évanescents, donne la suite spectrale de filtration de monodromie-locale du complexe des cycles
évanescents du modèle de Deligne-Carayol pour la hauteur $h$. On référera à ce résultat comme le théorème de
comparaison de Berkovich-Fargues.

\medskip

\noindent \textbf{0.13.} --- Le quatrième paragraphe est consacré aux énoncés globaux. On commence, après
avoir fait quelques rappels sur les faisceaux pervers \S \ref{intro}, par découper, proposition
(\ref{prop-so}), nos faisceaux pervers selon leur composantes isotypiques pour le sous-groupe d'inertie $I_o$
et on note ainsi $gr_{k,\pi_o}$ (resp. $gr_{k,\r_\oo,\pi_o}$) la composante $L_g(\pi_o)_{|I_o}$-isotypiques
du gradué $gr_k$ (resp. $gr_{k,\r_\oo}$) pour la filtration de monodromie du faisceau pervers
$R\Psi_{\eta_o}(\bar \Qm_l)[d-1]$ (resp. $R\Psi_{\eta_o}(\LC_{\r_\oo})[d-1]$ pour $\r_\oo$ une représentation
irréductible de $\bar D_\oo^\times$), où $\pi_o$ est une représentation irréductible cuspidale de $GL_g(F_o)$
avec $1 \leq g \leq d$. La représentation $L_g(\pi_o)_{|I_o}$ n'étant pas irréductible, $V_{\pi_o}$ est le
facteur direct de $V$ sur lequel l'action de $I_o$ se décompose en une somme de sous-représentations
irréductibles qui sont aussi des sous-représentations de $L_g(\pi_o)_{|I_o}$.

On introduit ensuite \S \ref{defi-fph}, certaines catégories de faisceaux pervers de Hecke qui fourniront le
cadre catégoriel des différents complexe de faisceaux que nous considérerons dans la suite et on donne,
définition (\ref{defi-type}), un certain nombre de notations attachées aux systèmes locaux d'Harris-Taylor
$\FC(g,l,\pi_o)_1$ et aux faisceaux induits $HT(g,l,\pi_o,\Pi_l)$ qui leurs sont associés sur la strate
$M_{I,s_o}^{=lg}$ où $\Pi_l$ est une représentation de $GL_{lg}(F_o)$ qui sera le plus souvent elliptique de
type $\pi_o$. On introduit alors le groupe de Grothendieck $\GF$ des faisceaux pervers de Hecke sur la tour
des $M_{I,s_o}$ munis d'une action compatible de $(D_\Am^\oo)^\times \times W_o$.

On énonce alors les théorèmes globaux qui précisent, théorème (\ref{theo-global0}) (resp.
(\ref{theo-global1})), les faisceaux pervers simples de $R\Psi_{\eta_o}(\LC_{\rho_\oo})[d-1]$ (resp. des
$gr_{k,\r_\oo,\pi_o}$) en termes des faisceaux pervers $j^{\geq lg}_{!*} HT(g,l,\pi_o,\st_l(\pi_o))[d-lg]
\otimes L_g(\pi_o) (-\frac{lg-1+k}{2})$ avec $1 \leq l \leq d/g$ (resp. $|k| < l \leq s_g$ et $l \equiv k-1
\mod 2$).

D'après le théorème de comparaison de Berkovich-Fargues, le théorème local (\ref{theo-ripsi-local}) se déduit
alors du calcul, théorème (\ref{theo-global2}), des faisceaux de cohomologie des $j^{\geq lg}_{!*}
\FC(g,l,\pi_o)$ et de la détermination, théorème (\ref{theo-ss}), de la suite spectrale de monodromie
associée.

On donne ensuite le schéma de la preuve qui procède par récurrence en supposant connu (\ref{theo-local-fil})
pour tout $d'<d$. On renvoie le lecteur à \S \ref{schema} pour un apercu des implications logiques entre les
divers énoncés locaux et globaux.

\medskip

\noindent \textbf{0.14.} --- Dans le cinquième paragraphe, on démontre le théorème (\ref{theo-global0}), i.e.
on donne l'image de $R\Psi_{\eta_o}(\bar \Qm_l)[d-1]$ dans le groupe de Grothendieck $\GF$ des faisceaux
pervers de Hecke sur la tour des $M_{I,s_o}$. La preuve procède en plusieurs étapes. Tout d'abord de la
description (\ref{iso1}) des restrictions aux strates $M_{I,s_o}^{=h}$ des $R^i \Psi_{\eta_o}(\bar \Qm_l)$ et
du calcul (\ref{somme-alternee}) de $\sum_i (-1)^i [\widetilde{\UC_{F_o}^{d,i}}]$, on en déduit, proposition
(\ref{prop-libre}), l'égalité suivante où on a posé $s_g=\lfloor \frac{d}{g} \rfloor$, la partie entière de
$d/g$:

\begin{multline} \label{egaliteR}
[R\Psi_{\eta_o}(\bar \Qm_l)[d-1]]=\sum_{g=1}^d \sum_{\pi_o \in \cusp_o(g)} \frac{1}{e_{\pi_o}}
\sum_{i=1}^{s_g} \sum_{l=i}^{s_g} (-1)^{l-i} \\
[j^{\geq lg}_! HT(g,l,\pi_o, [\overleftarrow{i-1},\overrightarrow{l-i}]_{\pi_o})[d-lg] \otimes L_g(\pi_o)
(-\frac{lg-2+2i-l}{2})]
\end{multline}
L'étape suivante consiste alors à exprimer l'image des faisceaux pervers qui interviennent dans le membre de
droite de l'égalité ci-dessus soit:
\begin{multline} \label{egalitej}
[j^{\geq lg}_! HT(g,l,\pi_o,[\overleftarrow{l-1}]_{\pi_o})[d-lg]]= \\ \sum_{i=0}^{s_g-l} j^{\geq l'g}_{!*}
HT(g,l+i,\pi_o,[\overleftarrow{l-1}]_{\pi_o} \overrightarrow{\times} [\overleftarrow{i-1}]_{\pi_o})[d-(l+i)g]
\otimes \Xi^{i(g-1)/2}
\end{multline}
Pour prouver cette dernière égalité, on raisonne par récurrence descendante sur la dimension des supports des
faisceaux pervers simples dans $j^{\geq lg}_! HT(g,l,\pi_o,[\overleftarrow{l-1}]_{\pi_o})[d-lg]$ en
réinjectant (\ref{egalitej}) pour les différents $l$, avec $g$ et $\pi_o$ fixés, dans (\ref{egaliteR}). On
argumente tout d'abord sur le fait que le membre de droite de (\ref{egaliteR}) doit s'écrire comme une somme
à coefficients positifs de faisceaux pervers simples, autoduale pour la dualité de Verdier.

On invoque ensuite le théorème de comparaison de Berkovich-Fargues afin d'utiliser l'hypothèse de récurrence
sur les modèles locaux en hauteur $h<d$, afin d'obtenir des renseignements sur les germes aux points non
supersinguliers de ces faisceaux pervers simples. On démontre alors, proposition (\ref{prop-p}), le résultat
hors des points supersinguliers au sens où les égalités du théorème précédent et de (\ref{egalitej}) sont
vraies si on rajoute une somme alternée de faisceaux pervers à support sur les points supersinguliers. Le
fait est qu'on utilise vraiment (\ref{iso1}) et pas seulement un calcul de somme alternée des restrictions
des faisceaux des cycles évanescents, ce qui explique l'indétermination au niveau des points supersinguliers.

Le paragraphe (\ref{fp-ponctuel}) est consacré à la détermination de ces faisceaux pervers ponctuels qui nous
manquent. Une idée naive est que pour connaître un faisceau ponctuel, on peut commencer par calculer son
groupe de cohomologie $H^0$. Étant donné un tel faisceau pervers ponctuel $\PC$ à déterminer, nous verrons en
fait que la connaissance de son $H^0$ suffit à le déterminer complètement: en effet grâce à
(\ref{egalite-groth}), nous montrerons que $\PC$ contient $\FC(g,s,\pi_o) \otimes (\Pi_l
\overrightarrow{\times} [\overleftarrow{s-l-1}]_{\pi_o}) \otimes \Xi^{(s-l)(g-1)/2}$ avec
$$H^0(\FC(g,s,\pi_o)\otimes\Pi_l \overrightarrow{\times} [\overleftarrow{s-l-1}]_{\pi_o}) \otimes \Xi^{(s-l)(g-1)/2}
=H^0(\PC)$$ de sorte que $\PC =\FC(g,s,\pi_o) \otimes \Pi_l \overrightarrow{\times}
[\overleftarrow{s-l-1}]_{\pi_o} \otimes \Xi^{(s-l)(g-1)/2}$.

On commence alors par calculer les groupes de cohomologie des faisceaux pervers d'Harris-Taylor, ou tout du
moins leur $\Pi^{\oo,o}$-partie, pour $\Pi$ automorphe vérifiant $\hyp(\oo)$ avec $\Pi_o=\st_s(\pi_o)$ de
sorte que, remarque (\ref{rema-hyp1}), la condition $\hyp(\oo)$ implique que $\Pi_\oo$ est isomorphe à
$\st_{s'}(\pi_\oo)$ pour un diviseur $s'$ de $d=s'g'$ et $\pi_\oo$ une représentation irréductible cuspidale
de $GL_{g'}(F_\oo)$. On montre à la proposition (\ref{prop-coho1}) que ceux-ci sont alors tous nuls de sorte
que, d'après la proposition (\ref{prop-p}), l'égalité (\ref{somme-alternee}) fournit, corollaire
(\ref{coro-pnul}), le $H^0$ des faisceaux ponctuels manquant ainsi que leur détermination. Le théorème
(\ref{theo-global0}) découle alors directement de ces résultats, cf. le corollaire (\ref{global0-ok}).

\medskip

\noindent \textbf{0.15.} --- Le sixième paragraphe est consacré, proposition (\ref{prop-hic}), au calcul des
faisceaux de cohomologies des faisceaux pervers $j^{\geq lg}_{!*} HT(g,l,\pi_o,\Pi_l)[d-lg]$ en fonction des
extensions par zéro $j^{\geq l'g}_! HT(g,l',\pi_o,\Pi_l \overrightarrow{\times} [\overrightarrow{l'-l-1}]_{\pi_o}) 
\otimes
\Xi^{(l'-l)(g-1)/2}$. D'après la proposition (\ref{prop-p}), on peut procéder par récurrence en utilisant la
suite spectrale (\ref{suite-spectrale}) dont l'aboutissement est connu sauf au niveau des points
supersinguliers.

Ainsi ces faisceaux de cohomologie ne sont pas complètement déterminés aux points supersinguliers mais il ne
reste qu'un nombre réduit de possibilités, cf. le lemme (\ref{lem-hij3}), que l'on peut obtenir, d'après
(\ref{prop-p}), par récurrence en utilisant par exemple la suite spectrale (\ref{ss-dualite}).

On donne ensuite, proposition (\ref{prop-grk-final}), une preuve de monodromie-poids. En égale
caractéristique ce résultat est connu et nous l'utilisons à maintes reprises, cependant, nous en démontrons
un cas particulier qui dans la situation de caractéristique mixte fournira une preuve complète du théorème de
monodromie-poids pour les variétés de Shimura de \cite{h-t}. On remarque tout d'abord que pour $s>2$, le
raisonnement par récurrence implique le résultat de sorte qu'il reste à initialiser la récurrence et donc à
prouver que pour $s=2$, la monodromie n'est pas triviale, en particulier sur la cohomologie. En égale
caractéristique le résultat est connu d'après \cite{lrs} tandis qu'en caractéristique il est connu pour $s=2$
et $g=1$ d'après \cite{ca}. Le principe, qui nous a été suggéré par M. Harris, est de s'y ramener via un
changement de base adéquat, cf. appendice A. \footnote{Notons par ailleurs que le cas général, i.e. $g=1$ et
$s$ quelconque a récemment été prouvé directement par Yoshida et Taylor.}

\medskip

\noindent \textbf{0.16.} --- Au septième paragraphe, on prouve les théorèmes globaux sous la proposition
(\ref{cas-m}). D'après le théorème de comparaison de Berkovich-Fargues la connaissance du théorème local
(\ref{theo-local-fil}) revient à celle des germes aux points supersinguliers des faisceaux de cohomologie des
faisceaux pervers d'Harris-Taylor ainsi que les flèches correspondantes dans la suite spectrale
(\ref{suite-spectrale}). Si nous disposions du théorème (\ref{theo-local-fil}) pour $d$, le raisonnement de
la preuve de la proposition (\ref{prop-hij1}) nous permettrait de déterminer complètement les faisceaux de
cohomologie des $gr_k$. Par souci d'efficacité nous montrons, proposition (\ref{prop-hij-ss-g}), qu'il
suffit, en utilisant l'opérateur de monodromie $N$, en fait de connaître les parties de poids $s(g-1)$ de
l'aboutissement de la suite spectrale de monodromie-locale de (\ref{theo-local-fil}). On est alors ramené à
prouver la proposition (\ref{cas-m}), i.e. à déterminer les parties de poids $s(g-1)$ des
$\widetilde{\UC_{F_o}^{d,i}}$.

\medskip

\noindent \textbf{0.17.} --- Le huitième paragraphe est consacré à la preuve de la proposition (\ref{cas-m}).
Celle-ci repose sur l'étude de la suite spectrale des cycles évanescents. Pour $\Pi$ une représentation
irréductible automorphe de $D_\Am^\times$ telle que $\Pi_o \simeq \st_s(\pi_o)$ avec $\pi_o$ cuspidale, la
$\Pi^{\oo,o}$-partie de son aboutissement est connue d'après \cite{lrs}, ou peut-être rapidement recalculée à
partir de la proposition (\ref{prop-coho1}). La détermination de cet aboutissement nous restreint, corollaire
(\ref{coro-1}), alors le nombre de possibilités pour les
$\widetilde{\UC_{F_o}^{d,i}}(\JL^{-1}(\st_s(\pi_o)))$.

En utilisant une propriété, théorème (\ref{involution}), d'invariance des $\UC_{F_o}^{d,\bullet}$ sous
l'involution de Zelevinski, il est alors possible de prouver la proposition (\ref{cas-m}). Cependant la
preuve de celle-ci repose sur tout ou partie du théorème de comparaison de Faltings à partir d'un énoncé
similaire du coté espace de Drinfeld, ce qui dépasse le cadre de ce texte. Par ailleurs la preuve, dans le
cas de la caractéristique mixte, de la conjecture de monodromie-poids version cohomologique demande une étude
des $\Pi^{\oo}$-parties des divers groupes de cohomologie, pour $\Pi$ automorphe vérifiant $\hyp(\oo)$ et
telle que $\Pi_o \simeq \speh_s(\pi_o)$. La $\Pi^{\oo,o}$-partie de la cohomologie de la fibre générique
n'est pas à priori connue. On se propose dans un premier temps de la calculer, ou tout du moins les bouts de
poids $s(g-1)$, le cas général étant traité, de manière indépendante, à la proposition (\ref{prop-lrs2}).

On calcule tout d'abord, proposition (\ref{prop-not}), les $\Pi^{\oo,o}$-parties des groupes de cohomologie
des faisceaux pervers d'Harris-Taylor. On en déduit alors, corollaire (\ref{coro-ssce-min}), la connaissance
des $\Pi^{\oo,o}$-parties de poids $s(g-1)$ des groupes de cohomologie de la fibre générique. Pour ce faire
on utilise le théorème de Lefschetz difficile. On étudie alors, proposition (\ref{prop-ssce-poids}), les
$\Pi^{\oo,o}$-parties des termes $E_2^{p,q}$ de la suite spectrale des cycles évanescents.

On montre, qu'à travers les suites spectrales associées à la stratification, qu'en ce qui concerne les bouts
de poids $s(g-1)$, seule la strate supersingulière contribue de sorte que l'on se retrouve dans une situation
similaire à celle de \cite{boy} dans le cas cuspidal, où les bouts de poids $s(g-1)$ des
$E_2^{p,q}[\Pi^{\oo,o}]$ sont nuls pour $p \neq 0$. Cette constatation découle du contrôle, lemme
(\ref{lem-rj-combi}), des contributions des strates non supersingulières. La proposition (\ref{cas-m}),
découle alors directement de la connaissance des $\Pi^{\oo,o}$-parties de l'aboutissement de la suite
spectrale des cycles évanescents et du fait qu'en ce qui concerne celles de poids $s(g-1)$, celles-ci ne
proviennent que des points supersinguliers.

\bigskip

\noindent \textit{Note:} pour les lecteurs qui n'ont pas lu les chapitres précédents, on leur conseille de
lire toutefois le \S \ref{rappel-ze-0} qui rappelle la combinatoire sur les représentations elliptiques de
$GL_d(F_o)$ due à Zelevinski (cf. \cite{ze}).


%% file: enonces.tex
\section{Rappels des données géométriques}

\begin{defisb}
- Soit $\cl:W_o \longto F_o^\times$ le morphisme de la théorie du corps de classe qui
envoie les frobenius
géométriques sur les uniformisantes, i.e. $\val(\cl (\fr_o))=-1$.

- Pour $c_o \in W_o$, on notera $\deg(c_o):=\val(\cl(c_o))$.

- Étant donnés une représentation complexe $\s_o$ (resp. $\pi_o$) de $W_o$ (resp. de
$GL_d(F_o)$) et un
entier $r$, on notera $\s_o(r)$ (resp. $\pi_o(r)$) la représentation $\s_o \otimes
|\cl|^{r}$ (resp. $\pi_o
\otimes |\det |^{r}$).
\end{defisb}

\subsection{Le modèle local de Deligne-Carayol}
\label{rapel-DC}

Pour tout $n \geq 0$, $\Def_n^d$ représente \footnote{cf. \cite{dr1} ou \cite{boy}} le
foncteur des
déformations de niveau $n$, par isogénies, du $\OC_o$-module formel de hauteur $d$ sur
$\bar \Fm_q$.

\marque Soit alors $\Psi_{F_o,n}^{d,i}$ le $\bar \Qm_l$-espace vectoriel de dimension
finie obtenu via la
théorie de Berkovich comme le $i$-ème foncteur des cycles évanescents associé au
morphisme structural
$$\spf \Def_n^d \longto \spf \hat \OC_o^{nr}.$$
Cet espace vectoriel est muni entre autre d'une action de $GL_d(\OC_o)$ qui se factorise
par le morphisme
surjectif naturel $GL_d(\OC_o) \longto GL_d(\OC_o/\MC_o^n)$ et on pose
$$\Psi_{F_o}^{d,i}= {\DS \lim_{\genfrac{}{}{0pt}{}{\longto}{n}}} \Psi_{F_o,n}^{d,i}$$
de sorte que pour $K_{o,n}:=\ker (\OC_o^\times \longto (\OC_o/\MC_o^n)^\times)$,
$\Psi_{F_o,n}^{d,i}=(\Psi_{F_o}^{d,i})^{K_{o,n}}$. On introduit le groupe $\NC_o$ (resp.
$\NC_o'$) défini
comme le noyau de
$$(g_o,\d_o,c_o) \in GL_d(F_o) \times D_{o,d}^\times \times W_o \mapsto
\val(\det(g_o^{-1})\rn(\d_o)\cl(c_o))
\in \Zm$$ (resp. composé avec la projection canonique $\Zm \longto \Zm/d\Zm$).

\rem Parfois on précisera la dimension en notant $\NC_o(d)$.

\marque Pour $\xi_o$ un caractère d'ordre fini de $F_o^\times$, on note
$\Psi_{F_o,\xi_o}^{d,i}$ la
$\xi_o'$-composante isotypique où $\xi_o'$ est la restriction de $\xi_o$ à
$\OC_o^\times$. Ainsi
$\Psi_{F_o,\xi_o}^{d,i}$ (resp. $\Psi_{F_o}^{d,i}$) est muni d'une action de $\NC_o'$
(resp. de $\NC_o$).

\marque Dans la définition de $\Def_n^d$, il est agréable de considérer plutôt les
déformations par
quasi-isogénies ce qui donne un schéma formel $\spf \Def_n^{d,\Zm} \simeq \coprod_\Zm
\spf \Def_n^d$ de sorte
que la construction précédente fourni des $\bar \Qm_l$-espaces vectoriels
$\UC_{F_o,n}^{d,i} \simeq
(\UC_{F_o}^{d,i})^{K_{o,n}}$ où
$$\UC_{F_o,\xi_o}^{d,i}:=\ind_{\NC_o}^{GL_d(F_o) \times D_{o,d}^\times \times
W_{F_o}}\Psi_{F_o,\xi_o}^{d,i}$$
est une représentation de $GL_d(F_o) \times D_{o,d}^\times \times W_{o}$.

\marque Pour toute représentation admissible irréductible $\t_o$ de $D_{o,d}^\times$ de
caractère central
$\xi_o$, la réciprocité de Frobenius donne que la composante isotypique
$\UC_{F_o,\xi_o}^{d,i}(\t_o)$ est
isomorphe à $\hom_{\DC_{o,d}^\times}(\res^{D_{o,d}^\times}_{\DC_{o,d}^\times}
\t_o,\Psi_{F_o,\xi_o^{-1}}^{d,i})$ où l'action de $(g_o^c,\s_o)$ est donnée par celle de
$(g_o,\d_o,\s_o) \in
\NC_o'$ pour $\d_o \in D_{o,d}^\times$ quelconque.

\subsection{Les variétés globales et les systèmes locaux d'Harris-Taylor}
\label{rappels-globaux}

Par souci de simplicité, les notations et le contexte sont ceux de \cite{lrs} et
\cite{boy}, i.e. le cas
d'égale caractéristique. Les résultats que nous rappelons sont une traduction de ceux de
\cite{h-t} et
peuvent être trouvés dans le chapitre I. Il se trouve que les preuves qui suivent
reposent formellement sur
les propriétés cohomologiques ci-après de sorte que les résultats que nous obtenons sont
également valides en
caractéristique mixte dans le cadre de \cite{h-t}; nous donnerons dans l'appendice A le
dictionnaire entre
nos notations et celle de \cite{h-t} et nous préciserons où trouver les propriétés
cohomologiques que nous
utilisons; pour de plus amples détails on renvoie à \cite{boy-duke}. Les données sont
alors:

\begin{itemize}
\item[(1)] une tour de schémas $(M_{I,o})_I$ définis sur le trait $\spec \OC_o$ et
indexés par les idéaux $I$ de
$A$; cette tour est munie d'une action, par correspondances, du groupe
$(D_\Am^\oo)^\times$;

\item[(2)] pour $\r_\oo$ une représentation irréductible du groupe des inversibles $\bar
D_\oo^\times$ de
l'algèbre à division centrale sur $F_\oo$ d'invariant $-1/d$, un système local
$\LC_{\r_\oo}$ sur les
$M_{I,o}$;

\item[(3)] une stratification de la fibre spéciale $M_{I,s_o}$, par des sous-schémas
fermés $M_{I,s_o}^{\geq h}$,
pour $1 \leq h \leq d$, de pure dimension $d-h$, telle que\footnote{C'est l'équivalent
du théorème de
Serre-Tate.} le complété de l'anneau local de $M_{I,o}$ en tout point géométrique de
$M_{I,s_o}^{=h}$ est
isomorphe à $\Def_n^h[[x_1^n,\cdots,x_{d-h}^n]]$. Par ailleurs les inclusions $j^{\geq
h}: M_{I,s_o}^{=h}
\hookrightarrow M_{I,s_o}^{\geq h}$ sont affines;

\item[(4)] les strates non supersingulières, i.e. $h \neq d$, sont induites, au sens où
il existe un sous-schéma
fermé $M_{I,s_o,1}^{=h}$ de $M_{I,s_o}^{=h}$ muni d'une action par correspondances de
$P_{h,d}(F_o)$ telle
que
$$M_{I^o\MC_o^n,s_o}^{=h}= M_{I^o\MC_o^n,s_o,1}^{=h}
\times_{P_{h,d}^{op}(\OC_o/\MC_o^n)} GL_d(\OC_o/\MC_o^n),$$
où l'action de $\left( \begin{array}{cc} g_o^c & 0 \\ * & g_o^{et} \end{array} \right) $
se factorise par
$(-\val(\det (g_o^c)),g_o^{et}) \in \Zm \times GL_{d-h}(F_o)$ et où $\frob_o$ agit via
$(-1,I_{d-h})$. Pour
tout $g \in GL_d(F_o)$ d'image $a \in GL_d(F_o)/P_{h,d}(F_o)$, on notera
$M_{I,s_o,a}^{=h}$ l'image par $g$
de $M_{I,s_o,1}^{=h}$ et $P_{h,d,a}$ le parabolique correspondant;

\item[(5)] pour $\t_o$ une représentation irréductible de $\bar D_{o,h}^\times$, un
système local $\FC_{\t_o,I,a}$ sur
$M_{I,s_o,a}^{=h}$; on note $\FC_{\t_o,I}:=\FC_{\t_o,I,a}
\times_{P_{h,d,a}(\OC_o/\MC_o^n)}
GL_d(\OC_o/\MC_o^n)$ où $n=\mult_o(I)$;

\item[(6)] Soit $\bar D$ l'algèbre à division centrale sur $F$ telle que $\bar
D^{\oo,o}_\Am \simeq
D^{\oo,o}_\Am$ et $\bar D_o$ (resp. $\bar D_\oo$) est l'algèbre à division centrale sur
$F_o$ (resp. sur
$F_\oo$) d'invariant $1/d$ (resp. $-1/d$). Pour tout idéal $I$ de $A$, l'ensemble des
points supersinguliers
$M_{I,s_o}^d(\bar \Fm_q)$ est en bijection avec le $(D_\Am^\oo)^\times \times
W_o$-ensemble
$$\bar D^\times \backslash [ (\bar D_\Am^{\oo,o})^\times / K_{\Am,I}^{\oo,o} \times
\Zm]$$
où l'action de $c_o \in W_o$ est donnée par la translation de $-\deg(c_o)$ sur la
composante $\Zm$ et où
l'action de $g^{\oo} \in (D_\Am^{\oo})^\times$ est donnée par la correspondance

$$\diagram & \bar D^\times \backslash [ \frac{(\bar
D_\Am^{\oo,o})^\times}{K_{\Am,J}^{\oo,o}} \times \Zm]
\dlto^{c_1} \drto^{c_2} \\
\bar D^\times \backslash [ \frac{(\bar D_\Am^{\oo,o})^\times}{K_{\Am,I}^{\oo,o}} \times
\Zm] & & \bar
D^\times \backslash [ \frac{(\bar D_\Am^{\oo,o})^\times}{K_{\Am,I}^{\oo,o}} \times \Zm]
\enddiagram$$
où $J$ est un idéal de $A$ tel que $K_{\Am,J}^{\oo,o} \subset K_{\Am,I}^{\oo,o} \cap
(g^{\oo,o})^{-1}
K_{\Am,I}^{\oo,o} g^{\oo,o}$, avec $c_1$ (resp. $c_2$) induit par l'inclusion
$K_{\Am,J}^{\oo,o} \subset
K_{\Am,I}^{\oo,o}$ (resp. par la multiplication à droite de $(g^{\oo,o})^{-1}$ sur
$(D_\Am^{\oo,o})^\times$)
et la translation de $-\val(\det(g_o))$ sur la composante $\Zm$.

\begin{defi} \label{defi-fg-s}
Pour tout diviseur $g$ de $d=sg$ et toute représentation irréductible cuspidale $\pi_o$
de $GL_g(F_o)$, on
notera $\FC(g,s,\pi_o,I)$ le faisceau concentré aux points supersinguliers
$$\diagram \bar D^\times \backslash [(\bar D_\Am^{\oo,o})^\times / K_{\Am,I}^{\oo,o}
\times (\Zm \times
(\JL^{-1}([\overleftarrow{s-1}]_{\pi_o}))^{K_o^n})] \dto \\
M_{I,s_o}^d(\bar \Fm_q)=\bar D^\times \backslash [(\bar D_\Am^{\oo,o})^\times /
K_{\Am,I}^{\oo,o} \times \Zm]
\enddiagram$$
avec $n=\mult_o(I)$ et où l'action diagonale de $\bar D^\times$ est donnée par
translation à droite sur
$(\bar D^{\oo,o}_\Am)^\times$, par translation de valeur $-\val \rn(.)$ sur $\Zm$ et par
l'action naturelle
sur $\JL^{-1}([\overleftarrow{s-1}]_{\pi_o})$.
\end{defi}
\end{itemize}

\section{Rappels des propriétés cohomologiques}
\label{rappels-coho}

\begin{defib} \label{defi-inert}
Soit $\t_o$ une représentation irréductible de $D_{o,h}^\times$, sa restriction à
$\DC_{o,h}^\times$ est une
somme de représentations irréductibles
$$\r_{o,1} \oplus \cdots \oplus \r_{o,e_{\t_o}}$$
et on notera $e_{\t_o}$ le nombre de celles ci. Étant donnée une représentation
irréductible $\r_o$ de
$\DC_{o,h}^\times$, soient alors $\t_o$ et $\t_o'$ des sous-représentations
irréductibles de l'induite de
$\DC_{o,h}^\times$ à $D_{o,h}^\times$ de $\r_o$: d'après la réciprocité de Frobenius, ce
sont exactement
celles telles que leur restriction à $\DC_{o,h}^\times$ contienne $\r_o$. On en déduit
alors que $\t_o$ et
$\t_o'$ sont inertiellement équivalentes, i.e. $\t_o' \simeq \t_o \otimes \xi_o$ avec
$\xi_o:\d \mapsto
x^{v(\det \d)}$ pour $x \in \bar \Qm_l^\times$. On note  $\CF_h$ l'ensemble des classes
d'équivalences
inertielles des représentations admissibles et irréductibles du groupe $D_{o,h}^\times$.
De la même façon,
deux représentations $\pi_o$ et $\pi_o'$ de $GL_g(F_o)$ seront dites inertiellement
équivalentes, s'il existe
un caractère $\xi:\Zm \longto \Qm_l^\times$ tel que $\pi_o \simeq \pi_o' \circ \xi \circ
\val \circ \det$ et
on notera $\cusp_g$ l'ensemble des classes d'équivalence inertielle des représentations
cuspidales de
$GL_g(F_o)$ ainsi que $e_{\pi_o}$ le cardinal de la classe d'équivalence inertielle de
$\pi_o$.
\end{defib}

\begin{defib} \label{wo-isotypique}
Soit $\s_o$ une représentation irréductible de $W_o$. Une représenta\-tion de $I_o$ sera
dite
$\s_o$-isotypique, si ses sous-représentations irréductibles sont aussi des
sous-représentations
irréductibles de la restriction à $I_o$ de $\s_o$.
\end{defib}

\rem Si une représentation de $I_o$ est $\s_o$-isotypique et $\s_o'$-isotypique alors
$\s_o$ et $\s_o'$ sont
inertiellement équivalentes.

\begin{lemb} \label{epio}
Pour tout entier $g$ et toute représentation irréductible cuspidale $\pi_o$ de
$GL_g(F_o)$,
$e_{\JL^{-1}([\overleftarrow{s-1}]_{\pi_o})}=e_{\pi_o}$.
\end{lemb}

\begin{proof} On rappelle que l'entier en question correspond au nombre de caractères
$\xi_o$ de $\Zm$ tels que
$\t_o \simeq \t_o \otimes \xi_o \circ \val \circ \det$ en notant
$\t_o=\JL^{-1}([\overleftarrow{s-1}]_{\pi_o})$; par Jacquet-Langlands c'est aussi le
nombre de caractères
$\xi_o$ de $\Zm$ tels que
$$[\overleftarrow{s-1}]_{\pi_o} \simeq [\overleftarrow{s-1}]_{\pi_o} \otimes (\xi_o
\circ \val \circ \det) \simeq
[\overleftarrow{s-1}]_{\pi_o \otimes (\xi_o \circ \val \circ \det)}$$ et donc au nombre
de caractères $\xi_o$
tels que $\pi_o \simeq \pi_o \otimes \xi_o \circ \val \circ \det$, d'où le résultat.

\end{proof}

\subsection{Sur le modèle local}

Soit $\t_o$ une représentation admissible irréductible de $D_{o,h}^\times$. L'espace que
l'on souhaite
étudier est le $(GL_h(F_o) \times W_{F_o})$-module
$$\hom_{D_{o,h}^\times}(\t_o,\UC_{F_o}^{h,i})=\UC_{F_o}^{h,i}(\t_o) \simeq
\Psi_{F_o}^{h,i}(\t_o):=\hom_{\DC_{o,h}^\times}(\t_o,\Psi_{F_o}^{h,i})$$

\marque Pour tout $\t_o$, on a un morphisme naturel de $\NC_o$-modules:
$$\UC_{F_o}^{h,i}(\t_o) \otimes \t_o \longto \Psi_{F_o}^{h,i}$$
qui envoie $f \otimes v$ sur $f(v)$. On note $\Psi_{F_o}^{h,i}[\t_o]$ l'image de ce
morphisme et soit
$\Psi_{F_o,m}^{h,i}[\t_o]$ la préimage de $\Psi_{F_o}^{h,i}[\t_o]$ dans
$\Psi_{F_o,m}^{h,i}$. Le sous-module
$\Psi_{F_o}^{h,i}[\t_o]$ ne dépend que de la classe d'équivalence inertielle de $\t_o$.
Le groupe
$\DC_{o,h}^\times$ étant compact, on a
$$\Psi_{F_o}^{h,i} = \bigoplus_{\t_o \in \CF_h} \Psi_{F_o}^{h,i}[\t_o]$$
Soit $\Delta_{\t_o}$ un ensemble d'\eles de $D_{o,h}^\times$ tel que les congruences des
$\val (\det \d)$
pour $\d \in \Delta_{\t_o}$ forment un système de représentants de $\Zm/e_{\t_o} \Zm$.
L'application
$$\begin{array}{rl} \label{action-tordue}
\UC_{F_o}^{h,i}(\t_o) \otimes \t_o & \longto \bigoplus_{\d \in \Delta_{\t_o}}
\Psi_{F_o}^{h,i}[\t_o]^\d \\
f \otimes v & \mapsto (f(\d^{-1}v))_\d
\end{array}$$
où $\Psi_{F_o}^{h,i}[\t_o]^\d$ est l'espace $\Psi_{F_o}^{h,i}[\t_o]$ muni de la
structure de $\NC_o$-module
où $(g_o,\d_o,c_o)$ agit via $(g_o,\d^{-1}\d_o \d ,c_o)$, est un isomorphisme de
$\NC_o$-modules.

\begin{defi}\label{action-modifiee}
On notera $\widetilde{\UC_{F_o}^{h,i}(\tau_o)}$ l'espace $\UC_{F_o}^{h,i}(\tau_o)$ où
l'action de $GL_h(F_o)$
est tordue par l'application $g_o \mapsto \lexp t g_o^{-1}$.
\end{defi}

\begin{theo} (cf. \cite{boy} théorème (3.2.4))
Pour toute représentation irréductible cuspidale $\pi_o$ de $GL_d(F_o)$,
$\widetilde{\UC_{F_o}^{d,i}}(\JL^{-1}(\pi_o))$ est nul pour $i \neq  d-1$ et
$$\widetilde{\UC_{F_o}^{d,d-1}}(\JL^{-1}(\pi_o)) \simeq \pi_o \otimes L_g(\pi_o)
(-\frac{d-1}{2})$$
\end{theo}

\begin{theo} (cf. \cite{h-t} théorème VII.1.5, ou le théorème (\ref{theo-calcul-psi}))
Pour tout diviseur $g$ de $d=sg$ et toute
représentation irréductible cuspidale $\pi_o$ de $GL_g(F_o)$, on a
\begin{multline} \label{inutile}
\sum_{i=0}^{d-1} (-1)^i
[\widetilde{\UC_{F_o}^{d,d-1-i}}(\JL^{-1}([\overleftarrow{s-1}]_{\pi_o}))]= \\
\sum_{i=1}^s (-1)^i [\overleftarrow{s-1-i},\overrightarrow{i}]_{\pi_o} \otimes
L_g(\pi_o)
(-\frac{d+s-2-2i}{2})
\end{multline}
ou de manière équivalente
\begin{multline} \label{inutile2}
\sum_{i=0}^{d-1} (-1)^i [\UC_{F_o}^{d,d-1-i}(\JL^{-1}([\overleftarrow{s-1}]_{\pi_o}))]=
\\ \sum_{i=1}^s
(-1)^i [\overrightarrow{i},\overleftarrow{s-1-i}]_{\pi_o} \otimes
L_g(\pi_o)(-\frac{d+s-2-2i}{2})
\end{multline}
\end{theo}

\subsection{Sur les systèmes locaux d'Harris-Taylor}
\label{ht-prop}

Dans la suite nous ne considérerons plus les schémas sur $\Fm_p$, $M_{I,s_o}$,
$M_{I,s_o}^{=h}$... mais
plutôt les $\bar \Fm_p$-schémas $M_{I,s_o} \times_{\Fm_p} \bar \Fm_p$, $M_{I,s_o}^{=h}
\times_{\Fm_p} \bar
\Fm_p$... Afin de ne pas alourdir encore plus les énoncés, nous garderons les mêmes
notations, par exemple
$M_{I,s_o}$ désignera ce que l'on devrait noter $M_{I,s_o} \times_{\Fm_p} \bar \Fm_p$.

Les systèmes locaux $\FC_{\t_o}$ d'Harris-Taylor sont tels que la restriction à
$M_{I,s_o,1}^{=h}$ du $i$-ème
faisceau des cycles évanescents $R^i \Psi_{\eta_o,I}(\bar \Qm_l)$ vérifie
\begin{equation} \label{iso1}
(R^i \Psi_{\eta_o,I}(\bar \Qm_l))_{M_{I,s_o,1}^{=h}}^h \simeq \bigoplus_{\t_o \in \CF_h}
(\FC_{\t_o,I}
\otimes \widetilde{\UC_{F_o,n}^{h,i}(\t_o)})^{h/e_{\t_o}},
\end{equation}
où $n=\mult_o(I)$ et où l'action l'action d'un élément $(g^{\oo,o},
g_o^{et},g_o^c,r,c_o) \in
(D_\Am^{\oo,o})^\times \times GL_{d-lg}(F_o) \times GL_{lg}(F_o) \times \Zm \times W_o$
est donnée par
l'action naturelle de
$$(g^{\oo,o},g_o^{et},r-\val(\det g_o^c)-\deg(c_o))) \in (D_\Am^{\oo,o})^\times \times
GL_{d-lg}(F_o) \times \Zm$$
sur la tour des $(\FC_{\t_o,I})_I$ au dessus de $M_{I,s_o,1}^{=lg}$, et par celle de
$(g_o^c,c_o)$ sur la
tour des $(\widetilde{\UC_{F_o,n}^{h,i}(\tau_o)})_n$. L'isomorphisme (\ref{iso1})
implique alors la
proposition suivante.

\begin{prop} \label{prop-hic}
Pour tout $i,j$, on a un isomorphisme canonique
$$H^j_c(M_{I,s_o,1}^{=h}, R^i\Psi_{\eta_o,I}(\LC_{\r_\oo}))^h \simeq \bigoplus_{\t_o \in
\CF_h}
(H^j_c(M_{I,s_o,1}^{=h},\FC_{\t_o,I} \otimes \LC_{\r_\oo}) \otimes
\widetilde{\UC_{F_o,m}^{h,i}(\t_o)})^{h/e_{\t_o}}$$ avec $m=\mult_o(I)$ et tel que
l'action de
$(g^{\oo,o},g_o^c,g_o^{et},\s_o) \in (D_\Am^{\oo,o})^\times \times GL_h(F_o) \times
GL_{d-h}(F_o) \times
W_{F_o}$ sur la limite inductive \footnote{indexée par les idéaux $I=I^o\MC_o^m$} du
membre de gauche, induit
l'action de
$$(g^{\oo,o},g_o^{et},-\val(\det g_o^c)-\deg(\s_o)) \otimes (g_o^c,\s_o)$$
sur la limite inductive du membre de droite.
\end{prop}

\rem Les systèmes locaux $\FC_{\t_o,I}$ ne sont pas irréductibles mais plutôt une somme
directe de $e_{\t_o}$
systèmes locaux irréductibles. La complexité de l'écriture de (\ref{iso1}), est la
contre-partie de la
simplicité de la description de l'action de $GL_d(F_o) \times W_o$ qui tient au fait que
l'on a fait
apparaître $\widetilde{\UC_{F_o}^{h,i}}(\t_o)$ plutôt que
$\hom_{\DC_{o,h}^\times}(\r_o,\widetilde{\Psi_{F_o}^{h,i}})$ avec $\r_o$ une
représentation irréductible de
$I_o$ de sorte que ce dernier est seulement muni d'une action de $\widetilde{\NC_o} \cap
(GL_h(F_o) \times
D_{o,h}^{\times})$.

\begin{defi} \label{defi-hyp}
Dans la suite, pour $\r_\oo$ une représentation irréductible de $\bar D_\oo^\times$, qui
est "l"'algèbre à
division centrale sur $F_\oo$ d'invariant $-1/d$, on considère pour une représentation
automorphe $\Pi$ de
$D_\Am^\times$, l'hypothèse $\hyp(\r_\oo)$ suivante: si $\r_\oo=JL^{-1}(\st_s(\pi_\oo)$
pour $\pi_\oo$ une
représentation irréductible cuspidale de $GL_g(F_\oo)$ pour $d=sg$, alors $\Pi_\oo$ est
soit isomorphe à
$\overleftarrow{(s)}_{\pi_\oo}$ soit à $\overrightarrow{(s)}_{\pi_\oo}$. Plus
généralement on notera
$\hyp(\oo)$ l'hypothèse sur $\Pi$ qu'il existe $\r_\oo$ tel que $\Pi$ vérifie
$\hyp(\r_\oo)$.
\end{defi}

\rem Dans le cas de caractéristique mixte la condition $\hyp(\oo)$ est rappelée dans
l'appendice A.

\medskip

On note ${\DS [H^*_{h,\r_\oo,\t_o}]:= \sum_i (-1)^i
\lim_{\genfrac..{0pt}{1}{\longto}{I}}
[H^i_c(M_{I,s_o,1}^{=h},\FC_{\t_o,I} \otimes \LC_{\r_\oo})]}$ dans le groupe de
Grothendieck des
représentations admissibles de $(D_\Am^{\oo,o})^\times \times GL_{d-h}(F_o) \times
\Zm$.

\begin{defi} \label{defi-carac}
Pour tout entier $h$, on considère l'identification canonique $D_{o,h}^\times /
\DC_{o,h}^\times \longto \Zm$
défini par la valuation de la norme réduite. On notera $\Xi:\Zm \longto \bar
\Qm_l^\times$ le caractère
défini par $\Xi(1)=\frac{1}{p}$
\end{defi}

\begin{prop} \label{prop-somme-alternee}
Pour $\Pi$ une représentation de $D_\Am^\times$, on a alors
\begin{equation} \label{somme-alternee}
[H^*_{h,\r_\oo,\t_o}(\Pi^{\oo,o})]= \left \{
  \begin{array}{ll}
  \e(\Pi)  m(\Pi) \Red_{\t_o}^h (\Pi_o) & \hbox{si } \Pi_\oo \hbox{ vérifie } \hyp(\oo)
\\
0 & \hbox{sinon}
  \end{array}
\right . \end{equation} où $\e(\Pi)$ est un signe qui dépend de $\Pi$, $m(\Pi)$ est la
multiplicité de $\Pi$
dans l'espace des formes automorphes et $\Red_{\t_o}^h: \groth(GL_d(F_o)) \longto
\groth(D_{o,h}^\times/\DC_{o,h}^\times \times GL_{d-h}(F_o))$ est défini comme la
composition des deux
homomorphismes suivant:
\begin{itemize}
\item en premier lieu, on a un homomorphisme
$$\begin{array}{l}
\groth(GL_d(F_o)) \longto \groth(GL_h(F_o) \times GL_{d-h}(F_o)) \\
~ [ \pi_o ] \mapsto [ J_{P_{h,d}}(\pi_o) \otimes \d_{P_{h,d}}^{1/2} ]
\end{array}$$

\item ensuite on a un homomorphisme
$$\begin{array}{l}
\groth(GL_h(F_o) \times GL_{d-h}(F_o)) \longto \groth(D_{o,h}^\times/\DC_{o,h}^\times
\times GL_{d-h}(F_o)) \\
~ [\a \otimes \b] \mapsto \sum_{\psi} \vol(D_{o,h}^\times/F_o^\times,d \bar h_o)^{-1}
\tr \a (\phi_{\JL(\t_o
\otimes \psi^{-1})})[\psi \otimes \b],
\end{array}$$
où $\psi$ décrit les caractères de $\Zm \simeq D_{o,h}^\times/\DC_{o,h}^\times$ tels que
$\a$ et $\t_o
\otimes \psi^{-1}$ ont le même caractère central et où l'on considère des mesures de
Haar associées sur
$GL_h(F_o)$ et $D_{o,h}^\times$.
\end{itemize}
\end{prop}

\rem En égale caractéristique on a $\e(\Pi)=1$ (resp. $\e(\Pi)=(-1)^{s-1}$) si $\Pi
\simeq \st_s(\pi_\oo)$
(resp. $\Pi_\oo \simeq \speh_s(\pi_\oo)$) pour $\pi_\oo$ une représentation irréductible
cuspidale de
$GL_g(F_\oo)$ avec $d=sg$. En caractéristique mixte, $\e(\Pi)=1$ (resp.
$\e(\Pi)=(-1)^{s-1}$) s'il existe une
place $x$ tel que $\Pi_x$ est de carré intégrable (resp. $\Pi_x \simeq \speh_s(\pi_x)$
pour $\pi_x$
irréductible cuspidale de $GL_g(F_x)$ avec $d=sg$). \footnote{L'auteur n'est pas certain
qu'en
caractéristique mixte, le cas $\Pi_x \simeq \speh_s(\pi_x)$ soit connu, de sorte qu'on
en donnera une
preuve.}

\medskip

Pour toute représentation irréductible $\r_\oo$ de $\bar D_\oo^\times$, on note
$\CC^\oo_{\bar D,\r_\oo}$ la
composante $\r_\oo$-isotypique de $C^\oo(\bar D^\times \backslash (\bar
D_\Am)^\times)$.

\begin{prop} \label{prop-h0-ss}
Pour tout $i$, on a un isomorphisme $(D_\Am^\oo)^\times \times W_o$-équivariant
\begin{equation} \label{h0-ss}
\lim_{\genfrac..{0pt}{1}{\to}{I}} H^0(M_{I,s_o}^d,R^i \Psi_{\eta_o,I}(\LC_{\r_\oo}))
\simeq \hom_{\bar
D_o^\times} ((\CC^\oo_{\bar D,\r_\oo})^\vee, \widetilde{\UC_{F_o}^{d,i}})
\end{equation}
\end{prop}

\section{Énonces des théorèmes locaux}

\subsection{Groupes de cohomologie des modèles de Deligne-Carayol}

Considérons le complexe à flèches nulles
\begin{multline*}
ML^\bullet (s):=([\overrightarrow{s-1}]_{1_o} \otimes |\cl|^{(s-1)/2},
[\overleftarrow{1},\overrightarrow{s-2}]_{1_o} \otimes |\cl|^{(s-3)/2} ,\cdots,
[\overleftarrow{s-1}]_{1_o} \otimes |\cl|^{(1-s)/2}) \\
\otimes \JL^{-1}(\st_s(1_o)) \otimes L_g(1_o)((1-s)/2)
\end{multline*}
où $[\overleftarrow{s-1}]_{1_o}$ est placé en degré $0$. Par définition, pour $\pi_o$
une représentation
cuspidale de $GL_g(F_o)$, on pose
\begin{multline*}
\pi_o \diamond ML^\bullet(s):=([\overrightarrow{s-1}]_{\pi_o} \otimes |\cl|^{(s-1)/2},
[\overleftarrow{1},
\overrightarrow{s-2}]_{\pi_o} \otimes |\cl|^{(s-3)/2},\cdots,
[\overleftarrow{s-1}]_{\pi_o} \otimes
|\cl|^{(1-s)/2}) \\
\otimes \JL^{-1}(\st_s(\pi_o)) \otimes L_g(\pi_o)((1-sg)/2)
\end{multline*}

\begin{theo} \label{theo-ripsi-local}
Pour tout $d$, on a
$$\widetilde{\UC_{F_o}^{d,d-1+\bullet}}=\bigoplus_{\atop{g|d}{d=sg}} \bigoplus_{\pi_o
\in \cusp_g}
\pi_o \diamond ML^\bullet(s)$$
où $\cusp_g$ désigne l'ensemble des classes d'équivalence des représentations
irréductibles cuspidales de
$GL_g(F_o)$.
\end{theo}

\marque Autrement dit, pour $\pi_o$ une représentation irréductible cuspidale de
$GL_g(F_o)$ on a:
$$\widetilde{\UC_{F_o}^{d,d-s+i}}(\JL^{-1}(\st_{s}(\pi_o))) \simeq
\left \{ \begin{array}{cl} L_g(\pi_o) (-\frac{d-s+2i}{2}) \otimes
[\overleftarrow{i},\overrightarrow{s-i-1}]_{\pi_o} & 0 \leq i < s \\ 0 & i < 0
\end{array} \right. $$ où
$[\overleftarrow{i},\overrightarrow{s-i-1}]_{\pi_o}$ est l'unique quotient irréductible
de l'induite
$$\ind_{P_{ig,sg}(F_o)^{op}}^{GL_{sg}(F_o)} \st_i(\pi_o(-\frac{(s-i)(g-1)}{2})) \otimes
\speh_{s-i}(\pi_o(\frac{i(g-1)}{2})).$$

\rem De manière équivalente on a
$$\UC_{F_o}^{d,d-s+i}(\JL^{-1}(\st_{s-1}(\pi_o))) \simeq L_g(\pi_o)(-\frac{d-s+2i}{2})
\otimes
[\overrightarrow{s-1-i},\overleftarrow{i}]_{\pi_o^\vee}$$

\subsection{Filtration de monodromie-locale}

On rappelle le théorème principal de l'appendice de \cite{boy-duke}, fournie par Laurent
Fargues.

\begin{theo-defi} \label{theo-defi-monodromie} (cf. le théorème principal de l'appendice
de \cite{boy-duke})
Soit $X \to \spec \OC_o$ un morphisme propre d'un schéma $X$ de type fini de dimension
$d$ sur le trait
$\spec \OC_o$. Le complexe des cycles évanescents $R\Psi_{\eta_o}(\Qm_l)[d-1]$ est un
faisceau pervers muni
d'une action de $W_o$ qui fournit une filtration de monodromie dont on notera $gr_{k}$
les gradués et on
considère la suite spectrale $E_1^{i,j}=h^{i+j} gr_{-i} \Rightarrow R^{i+j+d-1}
\Psi_{\eta_o}(\Qm_l)$. Pour
tout point géométrique $x$ de la fibre spéciale $X_{s_o}$, en considérant les germes en
$x$ des $h^i gr_{k}$,
on obtient une suite spectrale
$$E_{1,x}^{i,j}=(h^{i+j} gr_{-i})_x \Rightarrow R^{i+j+d-1} \Psi_{\eta_o}(\Qm_l)_x$$
dont la nature est purement locale de sorte que si $Y \to \spec \OC_o$ est un autre
schéma avec un point $y$
tel que le complété formel de l'anneau local de $Y$ en $y$ est isomorphe, en tant que
$\OC_o$-schéma formel,
au complété formel de l'anneau local de $X$ en $x$, alors pour tout $r \geq 1$, on a
$E_{r,x}^{i,j}=E_{r,y}^{i,j}$. La filtration ainsi obtenue sera dite de
monodromie-locale.
\end{theo-defi}

\marque Pour tout $s \geq 1$, on introduit une suite de bicomplexe $MLE_r^{i,j}(s)$
définit comme suit:

- pour $r=1$, $MLE_1^{i,j}(s)$ est nul pour $|j|\geq s$ ou $j \equiv s \mod 2$ ou $i+j >
0$ ou $i<(1-s-j)/2$.
Sinon pour $j=1-s+2r$ avec $0 \leq r < s$ et $i=s-1-2r-k$ avec $0 \leq k \leq(s-1-j)/2$,
on a
$$MLE_1^{s-1-2r-k,1-s+2r}(s)=[\overleftarrow{s-1-k}]_{1_o} \overrightarrow{\times}
[\overrightarrow{k-1}]_{1_o} \otimes
|\cl|^{(s-1-2r)/2}$$

- les flèches $d_1^{i,j}$ se déduisent des suites exactes courtes
(\ref{suites-exactes});

- pour tout $r \geq 2$, $MLE_r^{i,j}(s)=MLE_2^{i,j}(s)$ est nul pour $|j| \geq s$ ou $j
\equiv s \mod 2$ ou
$2i \neq 1-s-j$, et
$$MLE_2^{1-s+r,1-s+2r}=[\overleftarrow{s-1-r},\overrightarrow{r}]_{1_o} \otimes
|\cl|^{-(s-1-2r)/2}$$

Par exemple pour $s=4$ on a représenté à la figure (\ref{figure7}) les $MLE_1^{i,j}(4)$
et
$MLE_\oo^{i,j}(4)$.


Nous allons en fait prouver l'énoncé suivant dont découle directement le théorème
(\ref{theo-ripsi-local}).

\begin{theo} \label{theo-local-fil}
La filtration de monodromie-locale du complexe $\widetilde{\UC_{F_o,n}^{d,\bullet}}$ est
équivariante pour
l'action de l'algèbre de Hecke $\HC_n=\CC^\oo(K_{o,n} \backslash GL_d(F_o)/K_{o,n})$ et
celle de $I_o$; en
outre pour $n$ variant, les filtrations de monodromie locale forment des systèmes
inductifs compatibles à
l'action des correspondances de Hecke. Les termes $E_{r,gr,loc}^{i,j}$ de la suite
spectrale associée
$$E_{1,gr,loc}^{i,j}=h^{i+j} gr_{d,n,loc,-i} \Rightarrow
\widetilde{\UC_{F_o,n}^{d,d-1+i+j}}$$
vérifient
$$E_{\bullet,gr,loc}^{\bullet,d-1+\bullet}(\frac{d-1}{2})=\bigoplus_{g|d=sg}
\bigoplus_{\pi_o \in \cusp_g}
JL^{-1}(\st_s(\pi_o)) \otimes \bigl ( \pi_o \diamond MLE_\bullet^{\bullet,\bullet}(s)
\bigr )$$
\end{theo}

\marque Autrement dit, les gradués $gr_{d,n,loc,k}$ vérifient les propriétés suivantes:
\begin{itemize}
\item[(1)] pour tout $k$, $gr_{d,n,loc,k}$ est une somme directe
$$\bigoplus_{\genfrac{}{}{0pt}{}{g ~|~ d=sg}{\pi_o \in \cusp_g}}
gr_{d,n,loc,k,\pi_o},$$
où en tant que $I_o$-module, $gr_{d,n,loc,k,\pi_o}$ est $L_g(\pi_o)$-isotypique au sens
de la définition
(\ref{wo-isotypique}) \footnote{Sans utiliser le fait que les $L_g(\pi_o)$ décrivent
l'ensemble des
représentations irréductibles de $W_o$, l'énoncé dit que pour $\s_o$ qui ne serait pas
de la forme
$L_g(\pi_o)$, $gr_{d,loc,k,\s_o}$ est nul.}

\item[(2)] les $gr_{d,n,loc,k,\pi_o}$ sont nuls pour $|k| \geq s$ et pour $|k| < s$, ses
groupes
de cohomologie $h^i gr_{d,n,loc,k,\pi_o}$ sont nuls pour $i<d-s+|k|$ et pour $i \not
\equiv d-s+k \mod 2$;
pour $d-s+|k| \leq i=d-s+|k|+2r \leq d-1$, $h^{d-s+|k|+2r} gr_{d,loc,k,\pi_o}$ est
l'espace des vecteurs
invariants sous $K_{o,n} \subset GL_d(F_o)$ de l'espace suivant:
$$\JL^{-1}([\overleftarrow{s-1}]_{\pi_o}) \otimes [\overleftarrow{|k|+2r}]_{\pi_o}
\overrightarrow{\times}
[\overrightarrow{s-2-2r-|k|}]_{\pi_o} \otimes L_g(\pi_o)(-\frac{i+k}{2});$$

\item[(3)] la suite spectrale $E_{1,gr,loc,\pi_o}^{i,j}=h^{i+j} gr_{d,n,loc,-i,\pi_o}
\Rightarrow
\widetilde{\UC_{F_o,n}^{d,i+j}}$ dégénère en $E_2$ et ses $d_1^{i,j}$ se déduisent des
suites exactes courtes
\begin{multline} \label{suites-exactes}
0 \to [\overleftarrow{l-1},\overrightarrow{s-l}]_{\pi_o} \to
[\overleftarrow{l-1}]_{\pi_o}
\overrightarrow{\times} [\overrightarrow{s-l-1}]_{\pi_o} \to
[\overleftarrow{l},\overrightarrow{s-l-1}]_{\pi_o} \to 0
\end{multline}
\end{itemize}

\section{Énonces des théorèmes globaux}

\subsection{Rappels sur les faisceaux pervers}
\label{intro}

Pour $X$ un schéma séparé de type fini sur un corps, on considère la $t$-structure
perverse autoduale sur
$D_c^b(X,\bar \Qm_l)$. On reprend les notations et les résultats de \cite{ast}: en
particulier quand il n'y a
pas de confusion relativement à $X$, on note pour $a \in \Zm$, $\lexp p \t^{\leq a}$,
$\lexp p \t^{\geq a}$
les foncteurs de troncations de $D^b_c:=D_c^b(X,\bar \Qm_l)$ dans $D^{\leq
a}:=D_c^b(X,\bar \Qm_l)^{\leq a}$
et $D^{\geq a}:=D_c^b(X,\bar \Qm_l)^{\geq a}$, $\perv:=D^{\geq 0} \cap D^{\leq 0}$,
$\lexp p H^0:= \lexp p
\t^{\leq 0} \lexp p \t^{\geq 0}:D^b_c \longto \perv$.

\begin{defi} Soit $\s_o$ une représentation irréductible de $W_o$. Un faisceau pervers
$P$ irréductible, de la forme $j_{!*} \LC$ pour $j:U \hookrightarrow X$ localement fermé
et $\LC$ un système
local sur $U$, muni d'une action de $I_o$, sera dit $\s_o$-isotypique si $\LC$ l'est au
sens de la définition
(\ref{wo-isotypique}).
\end{defi}

\begin{lemm} Soient $\s_{o,1}$ et $\s_{o,2}$ des représentations irréductibles non
inertiellement équivalentes de $W_o$
et soient $P_1$, $P_2$ des faisceaux pervers muni d'une action de $I_o$, respectivement
$\s_{o,1}$ et
$\s_{o,2}$ isotypique. Si $\vphi:P_1 \longto P_2$ est un morphisme de faisceau pervers
$I_o$-équivariant,
pour tout $i$, on a $h^i \vphi=0$.
\end{lemm}

\begin{proof} C'est immédiat.

\end{proof}

Le lemme (5.3.6) de \cite{ast} donne alors la décomposition en somme directe suivante.

\begin{prop} \label{prop-so}
Tout faisceau pervers $P$ sur $X$ muni d'une action de $I_o$ s'écrit comme une somme
directe
$$ \bigoplus_{\s_o} P_{\s_o}$$
où la somme porte sur les classes d'équivalence inertielle des représentations
irréductibles de $W_o$ et où
$P_{\s_o}$ est $\s_o$-isotypique au sens de la définition (\ref{wo-isotypique}). En
particulier pour tout
point géométrique $x:\spec K \longto X$, $(P_{\s_o})_x$ est une représentation
$\s_o$-isotypique.
\end{prop}

\begin{coro} Soit $P$ un faisceau pervers sur $X$ muni d'une action compatible de $W_o$.
On considère sa filtration de
monodromie et on note $gr_k P$ son gradué de poids $k$. La suite spectrale
$$E_{1,gr,glob}^{i,j}=h^{i+j} gr_{-i} P \Longrightarrow h^{i+j} P$$
est la somme directe sur les classes d'équivalence inertielle des représentations
irréductibles $\s_o$ de
$W_o$ des suites spectrales
\begin{equation} \label{ss-gr}
E_{1,gr,glob,\sigma_o}^{i,j} =h^{i+j} gr_{-i} P_{\s_o} \Longrightarrow h^{i+j} P_{\s_o}
\end{equation}
\end{coro}

\marque Le complexe des cycles évanescents $R \Psi_{\eta_o,I}(\bar \Qm_l)[d-1]$ sur
$M_{I,s_o}$ est un
faisceau pervers muni d'une action de $W_{F_o}$; on considère alors sa filtration de
monodromie et on note
$gr_{k,I}$ son gradué\footnote{en égale caractéristique, les gradués $gr_{k,I}$ de la
filtration de
monodromie sont purs de poids $k+d-1$, i.e. les gradués de la filtration de monodromie
sont purs.} de niveau
$k$. On note $N$ le logarithme de la partie unipotente de la monodromie de sorte qu'en
particulier on a
$$N^k:gr_{k,I} \simeq gr_{-k,I}.$$
On pose $gr_k=(gr_{k,I})_I$, et on dispose alors de la suite spectrale
\begin{equation*}
E_{1,I,gr,glob}^{i,j}=h^{i+j} (gr_{-i,I}) \Rightarrow h^{i+j} (R \Psi_{\eta_o,I}(\bar
\Qm_l))[d-1]=R^{i+j+d-1} \Psi_{\eta_o,I}(\bar \Qm_l)
\end{equation*}
équivariante sous l'action de $(D_\Am^{\oo})^\times \times W_o$, qui, d'après le
corollaire précédent, se
scinde en une somme directe portant sur les classes d'équivalence inertielle des
représentations
irréductibles $\s_o$ de $W_o$:
\begin{equation} \label{suite-spectrale}
E_{1,I,gr,glob,\sigma_o}^{i,j}=h^{i+j} gr_{-i,I,\s_o} \Rightarrow h^{i+j} (R
\Psi_{\eta_o,I,\sigma_o}(\bar
\Qm_l))[d-1]=R^{i+j+d-1} \Psi_{\eta_o,I,\sigma_o}(\bar \Qm_l)
\end{equation}

\begin{defi} Pour $\r_\oo$ une représentation irréductible de $\bar D_\oo^\times$, on
notera $gr_{k,I,\r_\oo}$ le
gradué de $R\Psi_{\eta_o,I}(\LC_{\r_\oo})$ soit $gr_{k,I,\rho_\oo} \simeq gr_{k,I}
\otimes \LC_{\r_\oo}$.
Pour $\sigma_o$ une représentation irréductible de $W_o$, on notera
$gr_{k,I,\rho_\oo,\sigma_o}:=gr_{k,I,\sigma_o} \otimes \LC_{\rho_\oo}$. Pour $\pi_o$ une
représentation
irréductible cuspidale de $GL_g(F_o)$, on notera $gr_{k,I,\pi_o}$ (resp.
$gr_{k,I,\rho_\oo,\pi_o}$) pour
$gr_{k,I,L_g(\pi_o)}$ (resp. $gr_{k,I,\rho_\oo,L_g(\pi_o)}$) et de manière générale on
mettra un indice
$\pi_o$ en lieu et place de $L_g(\pi_o)$.
\end{defi}

\subsection{Définition de la catégorie des faisceaux pervers de Hecke sur $X_\IC$}
\label{defi-fph}

Il s'agit ici de donner le cadre catégoriel pour les énoncés et les preuves des
résultats des paragraphes
suivants; ce qui suit m'a été suggéré par Jean-François Dat.

\begin{defi}
On considère un groupe $G$ et on considère un ensemble $\IC$ tel qu'à tout $I \in \IC$
soit associé un
sous-groupe compact $\KC_I$ de $\tilde G$. On met sur $\IC$ une relation d'ordre partiel
$J \subset I$ si et
seulement si $\KC_J$ est un sous-groupe distingué de $\KC_I$. Un schéma de Hecke pour
$(G,\IC)$, est un
système projectif de schémas
$$X_\IC=(X_I)_{I \in \IC}$$
relativement à des morphismes finis $r_{J,I}:X_J \longto X_I$ de restriction du niveau,
tel que pour tout $g
\in G$ et tout $J \subset I$ tels que $g^{-1} \KC_J g \subset \KC_I$, on dispose d'un
morphisme fini de
schémas
$$[g]_{J,I}:X_J \longto X_I$$
vérifiant les propriétés suivantes:

- pour tout $g \in \KC_I$ et tout $J \subset I$, $[g]_{J,I}: X_J \longto X_I$ est égal
au morphisme de
restriction du niveau $r_{J,I}$;

- pour tout $g,g' \in G$, et tout $K \subset J \subset I$ tels que $g^{-1} \KC_J g
\subset \KC_I$ et
$(g')^{-1} \KC_K g' \subset \KC_J$, on a
$$[g]_{J,I} \circ [g']_{K,J}: X_K \longto X_J \longto X_I$$
est égal à $[gg']_{K,I}$.

- pour tout $J \subset I$, $X_J/\KC_I \simeq X_I$, où $g \in \KC_I$ agit sur $X_J$ via
$[g]_{J,J}$; autrement
dit $r_{J,I}:X_J \longto X_I$ est un bon quotient de $X_J$ par $\KC_I/\KC_J$.
\end{defi}

\noindent \textit{Exemples}: $\IC$ est l'ensemble des idéaux $I$ de $A$, avec
$X_I=M_{I,o},~M_{I,s_o}$ avec
$G=(D_\Am^\oo)^\times$ ou $X_I=M_{I,s_o,1}^{=h},~ M_{I,s_o,1}^{\geq h}$ avec
$G=(D_\Am^{\oo,o})^\times \times
P_{h,d}(F_o)$.

\begin{defi} Soit $X_\IC=(X_I)_I$ un schéma de Hecke pour $(G,\IC)$, on définit
alors la catégorie $\FPH_{G}(X_\IC)$ (resp. $\FH_{G}(X_\IC)$) des faisceaux pervers
(resp. des faisceaux) de
Hecke sur $X_\IC$ comme la catégorie dont les objets sont les systèmes $(\FC_I)_{I \in
\IC}$ où $\FC_I$ est
un faisceau pervers (resp. faisceau) sur $X_I$ tels que
\begin{itemize}
\item pour tout $g \in G$ et $J \subset I$ tel que $g^{-1} \KC_J g \subset \KC_I$, on
dispose d'un morphisme
de faisceau sur $X_I$:
$$u_{J,I}(g):\FC_I \longto [g]_{J,I,*} \FC_J$$
soumise à la condition de cocycle habituelle
$$u_{K,I}(g'g)=[g]_{J,I,*} (u_{K,J}(g')) \circ u_{J,I}(g)$$

\item pour tout $g \in \KC_I$, $u_{J,I}(g)$ se factorise par les $\KC_I$ invariants,
i.e. induit une flèche
$\FC_I \longto r_{J,I,*} \FC_J^{\KC_I}$ où $\FC_J^{\KC_I}$ désigne le sous-faisceau de
$\FC_J$ invariant par
tous les $u_{J,J}(g)$ où $g$ décrit $\KC_I$.
\end{itemize}

Les flèches sont alors les systèmes $(f_I)_{I \in \IC}$ avec $f_I: \FC_I \longto \FC'_I$
tel que le diagramme
suivant soit commutatif:
$$\diagram
\FC_I \rrto^{u_{J,I}(g)} \dto^{f_I} & & [g]_{J,I,*} (\FC_J) \dto^{[g]_{J,I,*} (f_J)} \\
\FC_I' \rrto^{u_{J,I}(g)} & & [g]_{J,I,*} (\FC_J')
\enddiagram$$
\end{defi}

\begin{prop} Pour tout schéma de Hecke $X_\IC$, la catégorie $\FH_G(X_\IC)$ (resp.
$\FPH_G(X_\IC)$) est abélienne
(resp. abélienne et artinienne).
\end{prop}

\begin{proof} Elle ne présente aucune difficulté. Montrons par exemple l'existence d'un
noyau pour
$f=(f_I:\FC_I \to \FC_I')_I$. Les catégories à niveau fini étant abélienne, soit pour
tout $I$, $\GC_I$ un
noyau de $f_I$. On a alors le diagramme commutatif suivant:
$$\xymatrix{
\GC_I \ar[r] \ar@{-->}[d] & \FC_I \ar[r]^{f_I} \ar[d]^{u_{J,I}(g)} & \FC_I'
\ar[d]^{u_{J,I}(g)} \\
[g]_{J,I,*} (\GC_J) \ar[r] \ar@{-->}[d] & [g]_{J,I,*} (\FC_J) \ar[r]^{[g]_{J,I,*} (f_J)}
\ar[d]^{[g]_{J,I,*}
(u_{K,J}(g'))} & [g]_{J,I,*} (\FC_J') \ar[d]^{[g]_{J,I,*}( u_{K,J}(g'))} \\
[g'g]_{K,J,*} (\GC_K) \ar[r] & [g'g]_{K,J,*} (\FC_K) \ar[r]^{[g'g]_{K,J,*}( f_K)} &
[g'g]_{K,J,*} (\FC_K')
}$$ De la propriété universelle du noyau $\GC_J$ (resp. $\GC_K$), on en déduit une
flèche
\begin{multline*}
u_{J,I}(g):\GC_I \longto [g]_{J,I,*} (\GC_J) \\ (\hbox{resp. }[g]_{K,J,*}(u_{K,J}(g')):
[g]_{J,I,*} (\GC_J)
\longto [g'g]_{K,J,*} (\GC_K))
\end{multline*}
la propriété de cocycle pour $\GC$ découle alors de celles pour $\FC$ et $\FC'$ via la
commutativité du
diagramme précédent. De la même façon, en utilisant la propriété universelle du noyau et
la commutativité du
diagramme suivant
$$\xymatrix{
\GC_I \ar[r] \ar@{-->}[d] & \FC_I \ar[r]^{f_I} \ar[d]^{u_{J,I}(g)} & \FC_I'
\ar[d]^{u_{J,I}(g)} \\
r_{J,I,*} (\GC_J)^{\KC_I} \ar[r] & r_{J,I,*} (\FC_J)^{\KC_I} \ar[r]^{r_{J,I,*} (f_J)} &
r_{J,I,*}
(\FC_J')^{\KC_I} }$$ on obtient que pour tout $g \in \KC_I$, $u_{J,I}(g)$ induit une
flèche $\GC_I \longto
r_{J,I,*} \GC_J^{\KC_I}$.

\end{proof}

\rem Dans la suite du texte, les faisceaux de Hecke que l'on considérera seront tels que
pour $g \in
\KC_{I}$, $u_{J,I}(g)$ induit un isomorphisme $\FC_I \longto r_{J,I,*} \FC_J^{\KC_I}$,
de sorte que par
exemple, les invariants sous $\KC_I$ de la cohomologie en niveau $J$, est isomorphe à la
cohomologie en
niveau $I$. Cette propriété est clairement conservée par passage au noyau mais ne l'est
pas pour les
conoyaux. Cela étant, on remarquera que les constructions des paragraphes suivants se
font toujours par des
noyaux à partir d'une flèche entre deux faisceaux vérifiant cette propriété
additionnelle de sorte que tous
les faisceaux que l'on construit la vérifient aussi.

\begin{nota} \label{nota-coho}
Pour $\FC=(\FC_I)_{I \in \IC}$ un objet de $\FPH_G(X_\IC)$ ou de $\FH_G(X_\IC)$, on
notera $H^i(\FC)$ pour la
limite inductive ${\DS \lim_{\genfrac..{0pt}{1}{\to}{I}}} H^i(X_I,\FC_I).$
\end{nota}

\begin{prop} \label{prop-def-fph}
Soient $\bar X_\IC$, $X_\IC$ et $Y_\IC$ des schémas de Hecke pour $(G,\IC)$ tels que
$j_\IC:X_\IC
\hookrightarrow \bar X_\IC$ soit un système projectif d'inclusions affines
$G$-équivariantes et $Y_\IC = \bar
X_\IC \backslash X_\IC$. On dispose alors des foncteurs
$$j_!,Rj_*,j_{!*}:\FPH_G(X_\IC) \longto \FPH_G(\bar X_\IC) \qquad \lexp p j^*=\lexp p
j^!:\FPH_G(\bar X_\IC) \longto
\FPH_G(X_\IC)$$ et pour $i:Y_\IC:=\bar X_\IC - X_\IC \hookrightarrow \bar X_\IC$ des
foncteurs
$$i_*=i_!:\FPH_G(Y_\IC) \longto \FPH_G(\bar X_\IC) \qquad \lexp p i^*, ~\lexp p
Ri^!:\FPH_G(\bar X_\IC) \longto \FPH_G(Y_\IC)$$
\end{prop}

\begin{proof} On commence par rappeler la proposition suivante tirée de \cite{ast}.

\begin{prop} \label{prop-ast0} (cf. proposition 1.3.17 de \cite{ast})
Pour tout foncteur exact $f$, on note $\lexp p f$ pour $\lexp p H^0 f$ comme foncteur
sur la catégorie des
faisceaux pervers.
\begin{itemize}
\item[(i)] Si $f$ est $t$-exact à gauche (resp. à droite), $\lexp p f$ est alors exact à
gauche (resp. à droite). Par ailleurs pour $\FC$
dans $D^{\geq 0}$ (resp. $D^{\leq 0}$), on a
$$\lexp p f (\lexp p H^0 \FC) \longmapright{\sim} \lexp p H^0 f(\FC) \quad (\hbox{resp.
} \lexp p H^0 f(\FC)
\longmapright{\sim} \lexp p f (\lexp p H^0 \FC)).$$

\item[(ii)] Soit $(f^+,f_+)$ une paire de foncteurs exacts adjoints; pour que $f^+$ soit
$t$-exact à droite, il faut et il suffit que $f_+$
soit $t$-exact à gauche et dans ce cas $(\lexp p f^+,\lexp p f_+)$ forment une paire de
foncteurs adjoints.

\item[(iii)] Si $f_1:D^b_{c,1} \to D^b_{c,_2}$ et $f_2:D^b_{c,2} \to D^b_{c,3}$ sont
$t$-exacts à gauche
(resp. à droite), alors $f_2 \circ f_1$ l'est aussi et $\lexp p (f_2 \circ f_1)= \lexp p
f_2 \circ f_1$.

\end{itemize}
\end{prop}

On rappelle que $j_!$, $Rj_*$, $j^*$, $i_*$ sont $t$-exacts et donc égaux à leur version
perverse alors que
$i^*$ est $t$-exact à droite. De même les morphismes $r_{J,I}$ et $[g]_{J,I}$ étant
finis, $r_{J,I}^*$ (resp.
$[g]_{J,I,*}$) est $t$-exact à droite (resp. $t$-exact) de sorte que $(\lexp p
r_{J,I}^*,\lexp p r_{J,I,*})$
et $(\lexp p [g]_{J,I,*}=[g]_{J,I,*},\lexp p [g]_{J,I}^!)$ forment des paires de
foncteurs adjoints avec
$\lexp p r_{J,I}^*$ et $\lexp p [g]_{J,I}^!$ (resp. $\lexp p r_{J,I,*}$, resp. $\lexp p
[g]_{J,I,*}$) exacts
à droite (resp. à gauche, resp. exact). Pour $J \subset I$, considérons le diagramme
suivant:
$$\xymatrix{ X_J \ar@{^{(}->}[rr]^{j_J} \ar[d]^{[g]_{J,I}} & & \bar X_J
\ar[d]_{\overline{[g]}_{J,I}}
& & Y_J \ar@{_{(}->}[ll]_{i_J} \ar[d]^{\widetilde{[g]}_{J,I}} \\
X_I \ar@{^{(}->}[rr]^{j_I} & & \bar X_I & & Y_I \ar@{_{(}->}[ll]_{i_I} }$$ Soient alors
$(\FC_I)_{I \in \IC}$
un objet de $\FPH(X_\IC)$ et $g,J,I$ avec $u_{J,I}(g): \FC_I \longto [g]_{J,I,*}
(\FC_I)$. Par application de
$Rj_{\IC,*}$ (resp. $j_{\IC,!}$), on obtient
$$Rj_{I,*} \FC_I \longmapright{Rj_{I,*}(u_{J,I}(g))} Rj_{I,*} [g]_{J,I,*}
(\FC_J)=\overline{[g]}_{J,I,*}
(Rj_{J,*} \FC_J)$$
$$(\hbox{resp. } j_{I,!} \FC_I \longmapright{j_{I,!}(u_{J,I}(g))} j_{I,!}
[g]_{J,I,!}(\FC_J)=\overline{[g]}_{J,I,*} (j_{J,!} \FC_J))$$ En ce qui concerne
$j_{\IC,!*}$, on considère le
diagramme commutatif suivant:
$$\xymatrix{
0 \to j_{I,!*} (\FC_I) \ar[r] \ar[ddddddr] \ar@{-->}[dddddd] & Rj_{I,*} (\FC_I)
\ar[d]_{Rj_{I,*}(u_{J,I}(g))}
\ar[r] & i_{I,*} \circ \lexp p i_I^* \circ Rj_{I,*} (\FC_I) \ar[d]_{i_{I,*}
\circ \lexp p i_I^* \circ Rj_{I,*}(u_{J,I}(g))} \to 0 \\
 & Rj_{I,*} \circ [g]_{J,I,*} (\FC_I) \ar[ddddd]^{=} \ar[r] &
 i_{I,*} \circ \lexp p i_I^* \circ Rj_{I,*} \circ [g]_{J,I,*} (\FC_I) \ar[d]^{=} \\
 & & i_{I,*} \circ \lexp p i_I^* \circ \overline{[g]}_{J,I,*} \circ Rj_{J,*} (\FC_J)
\ar[d]^{=} \\
 & & i_{I,*} \circ \lexp p (i_I^* \circ \overline{[g]}_{J,I,*}) \circ Rj_{J,*} (\FC_J)
\ar[d]^{(1)} \\
 & & i_{I,*} \circ \lexp p (\widetilde{[g]}_{J,I,*} \circ i_J^*) \circ Rj_{J,*} (\FC_J)
\ar[d]^{=} \\
 &  &
i_{I,*} \circ \widetilde{[g]}_{J,I,*} \circ \lexp p i_J^* \circ Rj_{J,*} (\FC_J)
\ar[d]^{=} \\
0 \to \overline{[g]}_{J,I,*} \circ j_{J,!*} (\FC_J) \ar[r] & \overline{[g]}_{J,I,*}
\circ Rj_{J,*} (\FC_J)
\ar[r] & \overline{[g]}_{J,I,*} \circ i_{J,*} \circ \lexp p i_J^* \circ Rj_{J,*} (\FC_J)
\to 0 }$$ où la
flèche notée (1) s'obtient à partir du diagramme commutatif suivant
$$\diagram Y_J \rto^{\widetilde{[g]}_{J,I}} \dto_{i_J} & Y_I \dto^{i_I} \\ X_J
\rto^{[g]_{J,I}} & X_I \enddiagram$$
comme le morphisme de changement de base $\bar i_I^* \circ [g]_{J,I,*} \longto
[g]_{J,I,*} \circ i_J^* $.

Les cas de $i_*$, $\lexp p i^*$ sont traités ci dessus et ceux de $\lexp p j^*$ et
$\lexp p Ri^!$ se traitent
de manière strictement similaire.

\end{proof}

\noindent \textit{Exemples}: étant donné un système local $\LC_\emptyset$ sur
$X_\emptyset$, soit pour $I \in
\IC$, $\LC_I=r_{I,\emptyset}^* \LC_\emptyset$; $(\LC_I[d])_{I \in \IC}$ est alors un
objet de
$\FPH_G(X_\IC)$, où $d$ est la dimension de $X_\IC$. Si en outre les $X_I$ sont la fibre
générique d'un
$S$-schéma $\XC_I$ de dimension relative $d-1$, pour $S$ un trait, alors
$(R\Psi_{\eta}(\LC_I)[d-1])_{I \in
\IC}$ est un objet de $\FPH_G(\bar X_\IC)$ où $\bar X_I$ désigne la fibre spéciale de
$X_I$.

\subsection{Notations}

\begin{defis} \label{defi-type}
\begin{itemize}
\item Pour $1 \leq g \leq d$, on notera dans la suite $s_g:=\lfloor \frac{d}{g}
\rfloor$, la partie entière
de $d/g$.

\item On introduit les injections
$$i^{h}:M_{I,s_o}^{\geq h} \hookrightarrow M_{I,s_o}, \qquad j^{ \geq h}: M_{I,s_o}^{=h}
\hookrightarrow
M_{I,s_o}^{\geq h} \qquad k^{=h}_1:M_{I,s_o,1}^{=h} \hookrightarrow M_{I,s_o}^{=h}$$ On
omet volontairement
l'indice $I$ dans les notations afin de ne pas encore alourdir le texte surtout qu'en
général le contexte est
clair.

\item On notera $\FPH(M_{s_o})$ (resp. $\FPH(M_o)$, resp. $\FPH(M_{s_o}^{=h})$, resp.
$\FPH(M_{s_o}^{\geq h})$)
la catégorie des faisceaux pervers de Hecke sur la tour des $M_{I,s_o}$ (resp.
$M_{I,o}$, resp.
$M_{I,s_o}^{=h}$, resp. $M_{I,s_o}^{\geq h}$) avec $G=(D_\Am^\oo)^\times$; on notera
aussi
$\FPH(M_{s_o,1}^{=h}$ (resp. $\FPH(M_{s_o,1}^{\geq h})$) la catégorie des faisceaux
pervers de Hecke sur la
tour des $M_{I,s_o,1}^{=h}$ (resp. $M_{I,s_o,1}^{\geq h}$) avec
$G=(D_\Am^{\oo,o})^\times \times
P_{h,d}(F_o)$. On utilisera des notations similaires pour les faisceaux de Hecke.

\item Un système projectif de faisceaux pervers ou pas, non nul $\FC=(\FC_I)_I \in
\FPH(M_{s_o}^{=h})$ (ou dans
$\FH(M_{s_o})$) sera dit \textrm{induit} s'il est de la forme
$$(k^{=h,*}_1 \FC_I)_{I \in \IC} \times_{P_{h,d}^{op}(F_o)} GL_d(F_o)$$
où $\times_{P_{h,d}^{op}(F_o)} GL_d(F_o):\FPH(M_{s_o,1}^{=h}) \longto
\FPH(M_{s_o}^{=h})$ est le foncteur qui
à un faisceau $(\FC_{0,I})_I$ associe $(\FC_{0,I} \times_{P_{h,d}^{op}(\OC_o)}
GL_d(\OC_o))_I$ de sorte que
pour tout $g \in GL_d(F_o)$, la correspondance associée induit un isomorphisme de $g^*
\FC_{\bar g}
\longmapright{\sim} \FC_0$ où $\FC_{\bar g}$ est la restriction de $\FC$ à la composante
$(M_{I,s_o,\bar
g}^{=h})_I$ où $\bar g$ est l'image de $g$ dans le quotient $GL_d(F_o)/P_{h,d}(F_o)$.

\item Pour $lg<d$, on notera $\FC(g,l,\pi_o)_1=(\FC(g,l,\pi_o,I)_1)_I$ le système local
$\FC_{\JL^{-1}([\overleftarrow{l-1}]_{\pi_o})}=(\FC_{\JL^{-1}([\overleftarrow{l-1}]_{\pi_o}),I})_I$
sur la
tour des $(M_{I,s_o,1}^{=lg})_I$ et soit $\FC(g,l,\pi_o)=(\FC(g,l,\pi_o,I))_I$ le
faisceau induit associé,
défini donc sur la tour des $(M_{I,s_o}^{=lg})_I$:
$$\FC(g,l,\pi_o)[d-lg] \in \FPH(M_{s_o}^{=lg})$$

\item Pour $1 \leq g \leq d$, un $GL_d(F_o) \times W_o$-faisceau sur la tour
$(M_{I,s_o}^{=lg})_I$ sera
dit de type $HT(g,l)$ s'il est induit, et si sa restriction à $(M_{I,s_o,1}^{=lg})_I$
est de la forme
$$(\FC(g,l,\pi_o,I)_1 \otimes \Pi_o^{K_{o,n}} \otimes \xi)_I$$
pour un certain triplet $(\Pi_o,\pi_o,\xi)$ où:

- $n=\mult_o(I)$,

- $\pi_o$ est une représentation irréductible cuspidale de $GL_g(F_o)$,

- $\Pi_o$ est une représentation de $GL_{lg}(F_o)$

- et $\xi$ un caractère de $\Zm$,

\noindent telle que l'action d'un élément $(g^{\oo,o}, g_o^{et},g_o^c,c_o) \in
(D_\Am^{\oo,o})^\times \times
GL_{d-lg}(F_o) \times GL_{lg}(F_o) \times W_o$ soit donnée par l'action naturelle de
$$(g^{\oo,o},g_o^{et},-\val(\det g_o^c) -\deg(c_o)) \in (D_\Am^{\oo,o})^\times \times
GL_{d-lg}(F_o) \times \Zm$$
sur $\FC_{\JL^{-1}([\overleftarrow{l-1}]_{\pi_o})}$ au dessus de
$(M_{I,s_o,1}^{=lg})_I$, et par celle de
$(g_o^c,c_o)$ sur $\Pi_o \otimes \xi$. On le notera alors
\begin{multline*}
HT(g,l,\pi_o,\Pi_o,\xi)=(HT(g,l,\pi_o,\Pi_o,\xi,I))_I, \\ HT(g,l,\pi_o,\Pi_o,\xi)[d-lg]
\in
\FPH(M_{s_o}^{=lg})
\end{multline*}
et on omettra $\xi$ si celui-ci est trivial. De même on notera
\begin{multline*}
HT_{\rho_\oo}(g,l,\pi_o,\Pi_o)=(HT_{\rho_\oo}(g,l,\pi_o,\Pi_o,I))_I, \\
HT_{\rho_\oo}(g,l,\pi_o,\Pi_o)[d-lg]
\in \FPH(M_{s_o}^{=lg})
\end{multline*}
le faisceau $HT(g,l,\pi_o,\Pi_o) \otimes \LC_{\rho_\oo}$.

\item Pour $\pi_o$ une représentation irréductible cuspidale de $GL_g(F_o)$ et $1 \leq l
\leq s_g$, on note
$$\PC(g,l,\pi_o)=(\PC(g,l,\pi_o,I))_I \in \FPH(M_{s_o})$$
le faisceau pervers défini comme l'extension intermédiaire $j^{\geq lg}_{!*}$ du
faisceau pervers sur
$M_{I,s_o}^{=lg}$, $HT(g,l,\pi_o,[\overleftarrow{l-1}]_{\pi_o})[d-lg] \otimes
L_g(\pi_o)$, qui est donc pur
de poids $d-lg$ pour $\pi_o$ unitaire \footnote{sinon on rajoute le poids de
$L_g(\pi_o)$} et qui ne dépend
que de la classe d'équivalence inertielle de $\pi_o$.
\end{itemize}
\end{defis}

\begin{rema} \label{rema-I} $\FC(g,l,\pi_o)$ et $HT(g,l,\pi_o,\Pi_o)$ (resp.
$\PC(g,l,\pi_o)$) se présentent
sous la forme d'une somme directe de $e_{\pi_o}$ systèmes locaux (resp. faisceaux
pervers) irréductibles:
$$\FC(g,l,\pi_o)=\bigoplus_{i=1}^{e_{\pi_o}} \FC(g,l,\rho_{o,i}) \qquad (\hbox{resp. }
\PC(g,l,\pi_o)=\bigoplus_{i=1}^{e_{\pi_o}} \PC(g,l,\rho_{o,i}))$$ donnés par la
restriction de
$\JL^{-1}([\overleftarrow{l-1}]_{\pi_o})$ à $\DC_{o,lg}^\times$ qui s'écrit comme une
somme directe
$\bigoplus_{i=1}^{e_{\pi_o}} \rho_{o,i}$ de représentations irréductibles, de sorte que
la différence entre
les faisceaux $\FC(g,l,\rho_{o,i})$ (resp. $\PC(g,l,\rho_{o,i})$) provient de l'action
de $\DC_{o,lg}^\times
\subset \NC_o$ donnée sur chacun comme dans la formule (\ref{action-tordue}).
\end{rema}

\begin{defi} \label{defi-gf}
Soit $\GF=(\GF_I)_I$ le système projectif de groupes de Grothendieck associé à
$\FPH(M_{I,s_o})$.
\end{defi}

\subsection{Filtration de monodromie}

\begin{theo} \label{theo-global0}
Soit $1 \leq g \leq d$ et $s_g$ la partie entière de $d/g$. Pour $\pi_o$ une
représentation irréductible
cuspidale de $GL_g(F_o)$, on a l'égalité suivante dans $\GF$:
$$e_{\pi_o} [R\Psi_{\eta_o,\pi_o}(\bar \Qm_l[d-1])]= \sum_{k=1-s_g}^{s_g-1}
\sum_{\genfrac{}{}{0pt}{}{|k| < l \leq s_g}
{l \equiv k-1 \mod 2}} \PC(g,l,\pi_o)(-\frac{lg+k-1}{2}),$$ où la torsion concerne
l'action de $W_o$.
\end{theo}

\marque Pour tout $s$, soit $MGr_k(s)$ le faisceau pervers nul pour $|k| \geq s$ et égal
à
$$\bigoplus_{\genfrac{}{}{0pt}{}{|k| <l \leq s}{l \equiv k-1 \mod 2}} j^{\geq l}_{!*}
\FC(1,l,1_o)[d-l] \otimes
[\overleftarrow{l-1}]_{1_o} \otimes L_g(1_o) (-(l-1+k)/2)$$ Concrètement on peut lire la
définition des
$MGr_\bullet(s)$ dans un repère $(l,k)$ en marquant sur la ligne $k=k_0$ les $l$ tels
que $j^{\geq l}_{!*}
\FC(1,l,1_o)[d-l] \otimes [\overleftarrow{l-1}]_{1_o} \otimes L_g(1_o) (-(l-1+k)/2)$
soit un constituant de
$MGr_{k_0}(s)$; par exemple pour $s=4$ on obtient la figure (\ref{figgg0}).

\begin{figure}[!h]
\input{figgure0.tex}

\caption{\label{figgg0} $MGr_\bullet(4)$}
\end{figure}

\marque Pour tout représentation cuspidale $\pi_o$ de $GL_g(F_o)$, on pose
\begin{multline*}
\pi_o \diamond \bigl ( j^{\geq l}_{!*} \FC(1,l,1_o)[d-l] \otimes
[\overleftarrow{l-1}]_{1_o} \otimes L_g(1_o)
(-(l-1+k)/2) \bigr ) := \\
j^{\geq lg}_{!*} \FC(g,l,\pi_o)[d-lg] \otimes [\overleftarrow{l-1}]_{\pi_o} \otimes
L_g(\pi_o) (-(lg-1+k)/2)
\end{multline*}
Le théorème (\ref{theo-global1}) se reformule alors sous la forme

\begin{theo} \label{theo-global1}
Soit $1 \leq g \leq d$ et $s_g$ la partie entière de $d/g$. Pour $\pi_o$ une
représentation irréductible
cuspidale de $GL_g(F_o)$, on a
$$e_{\pi_o} gr_{\bullet,\pi_o} = \pi_o \diamond MGr_\bullet(s_g)$$
\end{theo}

\marque Autrement dit, pour $1 \leq g \leq d$ et $\pi_o$ une représentation irréductible
cuspidale de
$GL_g(F_o)$, pour tout $|k| < s_g$, $e_{\pi_o} gr_{k,\pi_o}$ est égale à
$$\bigoplus_{\genfrac{}{}{0pt}{}{|k| < l \leq s_g}{l \equiv k-1 \mod 2}}
\PC(g,l,\pi_o)(-\frac{lg+k-1}{2}),$$
dans $\FPH(M_{s_o})$, où la torsion concerne l'action de $W_o$.

\begin{defi} \label{defi-pil}
Pour $1 \leq g \leq d$, et on fixe une représentation $\pi_o$ irréductible cuspidale de
$GL_g(F_o)$, que l'on
notera $\xi_o$ pour $g=1$. Pour tout $1 \leq l \leq s_g$, $\Pi_l$ désigne une
représentation quelconque de
$GL_{lg}(F_o)$, qui dans les applications sera elliptique de type $\pi_o$.
\end{defi}

\begin{theo} \label{theo-global2}
Pour $g$, $\pi_o$ et $\Pi_l$ comme ci-dessus, on considère la restriction à la tour des
$(M_{I,s_o}^{=h})_I$
de $h^i j^{\geq lg}_{!*} HT(g,l,\pi_o,\Pi_l)[d-lg]$;
\begin{itemize}
\item elle est nulle pour $h$ ne s'écrivant pas sous la forme $(l+a)g$ avec $0 \leq a
\leq s_g-l$;

\item pour $h=(l+a)g$ avec $0 \leq a \leq s_g-l$, elle est nulle pour $i \neq
lg-d+a(g-1)$ et sinon elle est
isomorphe dans $\FH(M_{s_o})$ à \footnote{cf. la définition (\ref{defi-carac})}
$$HT(g,l+a,\pi_o,\Pi_l \overrightarrow{\times} [\overrightarrow{a-1}]_{\pi_o}) \otimes
\Xi^{\frac{a(g-1)}{2}}$$
\end{itemize}
\end{theo}

On en déduit alors les corollaires suivant.

\begin{coro} Pour $g>1$, $h^{i} j^{\geq lg}_{!*} HT(g,l,\pi_o,\Pi_l)[d-lg]$
est nul pour $i$ ne s'écrivant pas sous la forme $lg-d+a(g-1)$ avec $0 \leq a \leq
s_g-l$ et pour un tel
$i=lg-d+a(g-1)$ il est isomorphe dans $\FH(M_{s_o})$ à
$$j^{\geq (l+a)g}_!HT(g,l+a,\pi_o,\Pi_l \overrightarrow{\times}
[\overrightarrow{a-1}]_{\pi_o}) \otimes
\Xi^{\frac{a(g-1)}{2}}.$$
\end{coro}

\begin{coro} \label{coro-higr}
Pour $\pi_o$ une représentation irréductible cuspidale de $GL_g(F_o)$ pour $1 \leq g
\leq d$, on a
\begin{itemize}
\item[(i)] pour $g>1$, $(h^i gr_{k,\pi_o})^{e_{\pi_o}}$ est une somme directe sur tous
les couples $(l,a)$ tels que
$-d+lg+a(g-1)=i$ et $l \equiv k-1 \mod 2$ des faisceaux induits dans $\FH(M_{s_o})$:
$$j^{\geq (l+a)g}_!HT(g,l+a,\pi_o,[\overleftarrow{l-1}]_{\pi_o} \overrightarrow{\times}
[\overrightarrow{a-1}]_{\pi_o})
\otimes L_g(\pi_o)(-\frac{lg-1+k+a(g-1)}{2})$$

\item[(ii)] pour $g=1$, la restriction de $h^i gr_{k,\xi_o}$ pour $i \equiv k-1 \mod 2$,
à la tour des $M_{I,s_o}^{=l}$
pour $l \geq d+i$ est égale au faisceau induit dans $\FH(M_{s_o})$
$$HT(1,l,\xi_o,[\overleftarrow{i-1}]_{\xi_o} \overrightarrow{\times}
[\overrightarrow{l-i-1}]_{\xi_o}) \otimes
\xi_o (-\frac{d+i-1+k}{2})$$ et nulle dans tous les autres cas.
\end{itemize}
\end{coro}

Le point (2) du théorème (\ref{theo-local-fil}) et donc l'énoncé de
(\ref{theo-ripsi-local}), découle alors
du résultat suivant. On introduit les bicomplexes
$$MGE_\bullet^{\bullet,\bullet}(s)=\bigoplus_{l=1}^s j^{\geq l}_! \FC(1,l,1_o)[d-l]
\otimes
MLE_\bullet^{\bullet,\bullet}(l)$$

\begin{theo} \label{theo-ss}
Soit $1 \leq g \leq d$ et $s_g$ la partie entière de $d/g$. Pour $\pi_o$ une
représentation irréductible
cuspidale de $GL_g(F_o)$, on a
$$e_{\pi_o} E_{\bullet,gr,glob,\pi_o}^{\bullet,\bullet}(\frac{d-1}{2})=
\pi_o \diamond MGE_\bullet^{\bullet,\bullet}(s_g)$$
\end{theo}

\marque Autrement dit la suite spectrale (\ref{suite-spectrale}) dégénère en $E_2$ et
les applications
$d_1^{i,j}$ sont induites par les suites exactes (\ref{suites-exactes}).


\subsection{Schéma de la preuve: par récurrence}
\label{schema}

Remarquons tout d'abord que le théorème (\ref{theo-global1}) découle directement de
(\ref{theo-global0}) en
utilisant monodromie-poids (MP), i.e. la pureté des gradués de la filtration de
monodromie. D'après le
théorème de comparaison de Berkovich-Fargues (BF), le théorème (\ref{theo-local-fil})
découle directement de
la détermination des germes en un point supersingulier des gradués de la filtration de
monodromie du complexe
des cycles évanescents de la variété de Drinfeld-Stuhler et de la suite spectrale
associée, soit donc des
théorèmes (\ref{theo-global1}), (\ref{theo-global2}) et (\ref{theo-ss}). Par ailleurs,
comme on l'a déjà
noté, le théorème (\ref{theo-ripsi-local}) découle directement de (\ref{theo-local-fil})
puisqu'il s'agit de
l'aboutissement de la suite spectrale associée à la filtration de monodromie-locale. On
résume ces
implications dans le tableau suivant.
\begin{figure}[!h]
\begin{center}
\begin{tabular}{|ll|}
\hline (\ref{theo-global0}) & \\ $\Downarrow$ MP & \\ (\ref{theo-global1}) + &
(\ref{theo-global2}) +
(\ref{theo-ss}) \\ \hline & $\Downarrow$ BF \\ \hfill (\ref{theo-local-fil})& en hauteur
$\leq d$ \\ \hline &
$\Downarrow$ \\ \hfill (\ref{theo-ripsi-local}) & en hauteur $\leq d$ \\ \hline
\end{tabular}
\end{center}
\end{figure}

Le principe est alors de raisonner par récurrence en supposant connu
(\ref{theo-local-fil}) pour tout $d'<d$.

\begin{rema} \label{rema-lim}
A priori, il nous faut travailler à niveau fini, pour tous les idéaux $I$ de $A$ et
montrer que les résultats
obtenus sont compatibles aux correspondances de Hecke; cependant devant l'immédiateté de
ces vérifications et
afin de ne pas alourdir encore les preuves, on raisonne souvent sur le système inductif
total.
\end{rema}

\marque Dans un premier temps le théorème de comparaison de Berkovich-Fargues (BF) et la
description (HT) des
faisceaux des cycles évanescents en fonction des systèmes locaux d'Harris-Taylor (cf.
(\ref{iso1})), nous
permet, \S \ref{glob0}, de montrer le théorème (\ref{theo-global0}).

\noindent --- On montre ensuite, proposition (\ref{prop-hij}), le théorème
(\ref{theo-global2}) hors des
points supersinguliers et on obtient un contrôle (C) sur ce qui peut se passer au niveau
de ces derniers.

\noindent --- Selon le même principe, si on disposait de (\ref{theo-local-fil}) pour la
hauteur $d$, on en
déduirait (\ref{theo-global2}) et (\ref{theo-ss}). On démontre, \S \ref{globss}, comment
la connaissance des
parties de poids minimal de (\ref{theo-ripsi-local}) permet d'obtenir
(\ref{theo-global2}), proposition
(\ref{prop-hij-ss-g}) \footnote{de sorte que
(\ref{theo-global2})=(\ref{prop-hij})+(\ref{prop-hij-ss-g})},
(\ref{theo-ss}), \S \ref{globss2} et même (\ref{theo-global1}) dans le cas Iwahori sans
utiliser
monodromie-poids, \S \ref{globalss3} ce qui, en caractéristique mixte, nous permettra de
prouver
monodromie-poids en toute généralité en se ramenant, par un changement de base
automorphe adéquat, au cas
Iwahori.

\noindent --- Reste alors, \S \ref{ssce}, à calculer les parties de poids minimal de
(\ref{theo-ripsi-local})
ce qui se fait via l'étude de la suite spectrale des cycles évanescents (SSCE) à travers
les suites
spectrales associées à la stratification (SSS) en utilisant le contrôle (C) sur ce qui
peut se passer au
niveau des points supersinguliers. Evidemment tout ceci passe par le calcul, ou au moins
du contrôle, des
groupes de cohomologie (CHT) des systèmes locaux d'Harris-Taylor, qui sont distillés au
fur et à mesure et
dont nous détaillons l'enchainement logique après le tableau récapitulatif suivant.
\begin{figure}[!h]
\begin{center}
\begin{tabular}{|ccl|}
\hline
\begin{tabular}{l} \S \ref{glob0}: (\ref{theo-local-fil}) \\ en hauteur $d'<d$
\end{tabular}
& \begin{tabular}{c} $\Longrightarrow$ \\ BF+HT \end{tabular} & (\ref{theo-global0})  \\
\hline
\begin{tabular}{l} \S \ref{glob2}: (\ref{theo-local-fil}) \\ en hauteur $d'<d$
\end{tabular}
& \begin{tabular}{c} BF+HT \\ $\Longrightarrow$ \\ (\ref{prop-hij}) \end{tabular} &
\begin{tabular}{l} (\ref{theo-global2}) hors des points supersinguliers \\ + C=contrôle
aux points supersinguliers \end{tabular} \\ \hline
\begin{tabular}{l} (\ref{theo-local-fil}) \\ en hauteur $d' \leq d$ \end{tabular}
& \begin{tabular}{c} $\Longrightarrow$ \\ (\ref{prop-hij-ss-g}) \end{tabular} &
(\ref{theo-global2})+ (\ref{theo-ss}) \\ \hline
\begin{tabular}{l} \S \ref{globss}: (\ref{theo-ripsi-local}) \\ en hauteur $d' < d$ \\
et en poids minimal \\
pour la hauteur $d$ \end{tabular} &
\begin{tabular}{c} $\Longrightarrow$ \\ C+N \end{tabular} &
$\left \{ \begin{array}{l} \hbox{(\ref{theo-global2}):\S \ref{globss1} proposition
(\ref{prop-hij-ss-g})} \\ \hbox{(\ref{theo-ss}): \S \ref{globss2}
} \\ \hbox{(\ref{theo-global1}): \S \ref{globalss3} sans MP}
\end{array} \right.$ \\ \hline
\S \ref{ssce0} & \begin{tabular}{c} $\Longrightarrow$ \\ SSCE+SSS \\ +C+CHT
\end{tabular} & \begin{tabular}{l} poids minimal de (\ref{theo-ripsi-local})\\ en
hauteur $d$
\end{tabular} \\ \hline
\end{tabular}
\end{center}
\end{figure}

\marque En ce qui concerne le rôle des divers calculs de groupes de cohomologie, voici
leur enchaînement
logique. Pour commencer on calcule exclusivement des composantes
$\Pi^{\oo,o}$-isotypiques pour $\Pi$ une
représentation irréductible automorphe de $D_\Am^\times$ vérifiant $\hyp(\oo)$, car dans
tous les autres cas
tous les calculs donnent des résultats nuls. On considère plus particulièrement deux
situations selon que
$\Pi_o$ est de la forme $[\overleftarrow{s-1}]_{\pi_o}$ ou
$[\overrightarrow{s-1}]_{\pi_o}$ pour $\pi_o$ une
représentation cuspidale de $GL_g(F_o)$ avec $d=sg$. Pour alléger les notations, on
notera $H^i(j^{\geq
lg}_!)$ pour la composante $\Pi^{\oo,o}$-isotypique du $i$-ème groupe de cohomologie du
complexe $j^{\geq
lg}_! \FC(g,l,\pi_o)[d-lg]$ et on utilisera une notation similaire pour $H^i(j^{\geq
lg}_{!*})$.

\marque La situation est largement plus simple dans le cas où
$\Pi_o=[\overleftarrow{s-1}]_{\pi_o}$. On
commence, proposition (\ref{prop-coho1}), par calculer les $H^i(j^{\geq lg}_{!*})$ et on
trouve que ceux-ci
sont tous nuls si $lg \neq d$. En utilisant (\ref{somme-alternee}), on en déduit alors,
corollaire
(\ref{coro-hij-nul}), le calcul des $H^i(j^{\geq lg}_!)$ qui sont nuls pour $i \neq 0$.
Ces calculs
permettent de compléter, corollaire (\ref{coro-pnul}), aux points supersinguliers, la
détermination des
constituants simples du faisceau pervers $j^{\geq lg}_! \FC(g,l,\pi_o)[d-lg]$ obtenu,
corollaire
(\ref{coro-sec-fp}), sur les autres strates grâce à la description des faisceaux des
cycles évanescents en
fonction des systèmes locaux d'Harris-Taylor ainsi qu'au théorème de comparaison de
Berkovich-Fargues. On en
déduit alors, corollaire (\ref{coro-hic}), le calcul des $\Pi^{\oo,o}$-parties des
groupes de cohomologie
$H^j_c(M_{I,s_o}^{=h},R^i\Psi_{\eta_o,I}(\LC_{\r_\oo}))$. A ce moment ci, on dispose du
théorème
(\ref{theo-global1}) de sorte que l'on peut calculer via la suite spectrale des poids
SSP, corollaire
(\ref{coho-global1}), les $\Pi^{\oo,o}$-parties des groupes de cohomologie
$H^i_{\eta_o,\r_\oo}$ de la fibre
générique. On en déduit alors, proposition (\ref{prop-hipsi}), les $\Pi^{\oo,o}$-parties
des termes
$E_2^{p,q}$ de la suite spectrale des cycles évanescents qui nous permet d'obtenir,
corollaire
(\ref{coro-1}), une première liste de possibilités pour les
$\widetilde{\UC_{F_o}^{d,i}}(\JL^{-1}(\st_s(\pi_o)))$.

\marque On considère alors $\Pi_o=[\overrightarrow{s-1}]_{\pi_o}$ et on calcule,
proposition
(\ref{prop-not}), les $H^i(j^{\geq lg}_{!*})$ qui sont non nuls pour $i \equiv s-l \mod
2$ et $|i| \leq s-l$.
De ce calcul et de la connaissance des constituants simples des $j^{\geq lg}_!
\FC(g,l,\pi_o)[d-lg]$, on en
déduit, lemme (\ref{lem-rj-combi}), un contrôle sur les $H^i(j^{\geq lg}_!)$. On étudie
ensuite la suite
spectrale des cycles évanescents. On commence par en déterminer l'aboutissement,
corollaire
(\ref{coro-ssce-min}) et proposition (\ref{prop-lrs2}), en utilisant le calcul des
$H^i(j^{\geq lg}_{!*})$
ainsi que le théorème (\ref{theo-global1}). En utilisant le théorème de Lefschetz
difficile ainsi que le
corollaire (\ref{coro-1}) qui contrôle ce qui se passe au niveau des points
supersinguliers, on en déduit,
proposition (\ref{prop-hic-poids}) et corollaire (\ref{strate-alt-p}), le calcul des
$H^i(j^{\geq lg}_!)$. Le
résultat est qu'en ce qui concerne les parties de poids minimal, $s(g-1)$, tout ce passe
au niveau des points
supersinguliers. L'étude de la suite spectrale des cycles évanescents fournit alors la
partie de poids
$s(g-1)$ du théorème (\ref{theo-ripsi-local}) qui d'après ce que l'on a vu, implique les
théorèmes globaux
(\ref{theo-global2}) et (\ref{theo-ss}) qui d'après le théorème de comparaison de
Berkovich-Fargues,
impliquent le théorème local (\ref{theo-local-fil}), ce qui complète la démonstration.
On résume la
discussion précédente dans le tableau suivant.
\begin{figure}[!ht]
\begin{center}
\begin{tabular}{|c|r|c|}
\hline $\Pi_o =\st_s(\pi_o)$ & & $\Pi_o=\speh_s(\pi_o)$ \\ \hline

$H^i(j^{\geq lg}_{!*})$ (\ref{prop-coho1}) & & \\
\begin{tabular}{c|c|} SSP $\Downarrow$ & $\Downarrow$ (\ref{somme-alternee}) \\
$H^i_{\eta_o}$ (\ref{coho-global1}) & $H^i(j^{\geq lg}_!)$ (\ref{coro-hij-nul}) \\
\hline  \end{tabular}
$\Rightarrow$ &
\begin{tabular}{c} $j^{\geq lg}_!$ \\ $=$ \\ $\sum j^{\geq (l+r)g}_{!*}$ \\
(\ref{coro-sec-fp}) \\  + \\  (\ref{coro-pnul})
\end{tabular} $\Rightarrow$ & $H^i(j^{\geq lg}_{!*})$ (\ref{prop-not}) \\

\begin{tabular}{lr} & $\Downarrow$ \end{tabular} & & \begin{tabular}{lr} $\Downarrow$
(\ref{theo-global1}) &
$\Downarrow \genfrac{}{}{0pt}{}{SSP}{LD}$ \end{tabular} \\

$H^j(M^{=h},R^i\Psi)~h \neq d$ (\ref{coro-hic}) & & \begin{tabular}{l|l} contrôle &
$H^i_{\eta_o}$ \\
$H^i(j^{\geq lg}_!)$ & (\ref{coro-ssce-min}) \\ (\ref{lem-rj-combi}) & (\ref{prop-lrs2})
\\ \hline \end{tabular} \\

$(\ref{coho-global1}) \Downarrow SSS$ & & \begin{tabular}{c} LD+\textbf{(\ref{coro-1})}
\\ $\Downarrow$ \end{tabular} \\

$H^j(R^i\Psi)~j \neq 0$ (\ref{prop-hipsi}) & & \begin{tabular}{c} \hline $H^i(j^{\geq
lg}_!)$ \\
(\ref{prop-hic-poids}) (\ref{strate-alt-p}) \\ \hline \end{tabular} \\

$\Downarrow$ SSCE & & \begin{tabular}{c} SSCE \\ $\Downarrow$ \end{tabular}  \\

contrôle $\Psi_{F_o}^{d,i}$ \textbf{(\ref{coro-1})} & \begin{tabular}{c|}
(\ref{theo-global2}) \\ (\ref{theo-ss}) \\
\hline \end{tabular} & \begin{tabular}{ll} $\Longleftarrow$ & $\Psi^{d,i}$
(\ref{theo-ripsi-local}) \end{tabular} \\

& \begin{tabular}{c} BF \\ $\Downarrow$ \end{tabular} & \\

& (\ref{theo-local-fil}) & \\ \hline
\end{tabular}
\end{center}
\end{figure}

\begin{rema} \label{rema-zw}
Afin de ne pas multiplier les situations, nous n'avons considéré dans nos énoncés que
les faisceaux induits
notés $HT(g,l,\pi_o,\Pi_l)$ à partir des systèmes locaux $\FC(g,l,\pi_o)_1 \otimes
\Pi_l$. Ceux-ci sont munis
d'une action par correspondances de $(D_\Am^\oo)^\times \times \Zm$; cependant quand on
les considère comme
des constituants de $R\Psi_{\eta_o}(\bar \Qm_l)$, on les voit comme munis d'une action
par correspondances de
$(D_\Am^\oo)^\times \times W_o$. Le lien entre les deux situations est donnée par $c_o
\in W_o \mapsto
-\deg(c_o) \in \Zm$. Afin de distinguer les deux situations, l'action de $W_o$ sera
toujours accompagnée d'un
$L_g(\pi_o)$ tandis que pour $\Zm$ il sera toujours question du caractère $\Xi$ (cf.
plus loin).

\noindent --- Par ailleurs nous ne considérons pas par la suite l'action totale de $W_o$
mais plutôt ses
frobenius semi-simplifiés.
\end{rema}

Avant de rentrer dans la preuve proprement dite, on rappelle la proposition
(\ref{prop-poids}) du chapitre I.

\begin{prop} \label{prop-poids2}
Soit $1 \leq lg <d$ et soit $\pi_o$ une représentation irréductible cuspidale de
$GL_g(F_o)$. Pour $1 \leq i
\leq l$ et pour tout $j$,
$$\lim_{\genfrac..{0pt}{1}{\to}{I}}
H^j_c(M_{I,s_o}^{=lg},HT_{\rho_\oo}(g,l,\pi_o,\Pi_l,I) \otimes L_g(\pi_o))$$
en tant que représentation de $GL_d(F_o) \times W_o$ est de la forme
$$\bigoplus_\xi(\ind_{P_{lg,d}^{\op}(F_o)}^{GL_d(F_o)} \Pi_l \circ
\xi(\val(\det))\otimes \pi_\xi)\otimes
L_g(\pi_o \circ \xi(\val(\det)))$$ où $\xi$ décrit les caractères $\Zm \longto \bar
\Qm_l^\times$ et
$\pi_\xi$ est une représentation de $GL_{d-lg}(F_o)$.
\end{prop}


%% file: figgure0.tex
\setlength{\unitlength}{.6cm} \centering
\begin{picture}(6,10)(0,-4)
\linethickness{.1pt}

\multiput(0,-3)(0,1){7}{\line(1,0){5}} \multiput(0,-3)(1,0){5}{\line(0,1){7}}

\put(5,0){\vector(1,0){1}} \put(6,.2){$l$}

\put(0,3){\vector(0,1){1}} \put(0.2,4){$k$}

\multiput(1,0)(1,1){4}{\circle{.2}}

\multiput(2,-1)(1,1){3}{\circle{.2}}

\multiput(3,-2)(1,1){2}{\circle{.2}}

\put(4,-3){\circle{.2}}

\put(0,0){\circle*{.1}}

\end{picture}

%% file: preuve.tex
\section{Étude dans le groupe de Grothendieck $\GF$}
\label{glob0}

\subsection{Groupe de Grothendieck des faisceaux pervers: généralités}

Soit $X$ un $\Fm_q$-schéma $X$. On rappelle que la catégorie $\perv(X)$ des faisceaux
pervers sur $X$ est
noethérienne et artinienne, ses objets simples étant de la forme $j_{!*} \LC$ où $\LC$
est un système local
irréductible sur $j: U \hookrightarrow X$, où $U$ est un ouvert d'un fermé de $X$.

Pour $P$ un objet de $\perv(X)$, sa dimension $n$ est par définition la plus grande
dimension des supports
$U$ des faisceaux pervers simples constituant $P$ de sorte que $-n= \min \{ i~/~ h^i P
\neq 0 \}$ et $h^{-n}
P$ a un support de dimension $n$.

\begin{defi} Étant donnée une catégorie triangulée $A$ localement petite, on considère
son groupe de Grothendieck $K(A)$ défini
comme le groupe libre engendré par les classes d'isomorphismes d'objets de $A$ quotienté
par les relations:

- $A[1]=-A$;

- $A=B+C$ pour tout triangle distingué $B \longto A \longto C \longmapright{+1}$.
\end{defi}

\begin{prop} Soit $(\DC^{\leq 0},\DC^{\geq 0})$ une catégorie dérivée d'une catégorie
abélienne $\AC$, munie d'une
$t$-structure non dégénérée: on note $\CC$ son coeur qui
est alors une catégorie abélienne de groupe de Grothendieck $\groth(\CC)$. L'application
qui à un objet $\FC$ de $\DC$ associe
$$\sum_i (-1)^i [\lexp p h^i \FC ] \in \groth(\CC)$$
induit un isomorphisme du groupe de Grothendieck $K(\DC)$ de la catégorie triangulée
$\DC$, sur $\groth(\CC)$.
\end{prop}

\begin{proof} Pour tout objet $\GC$ de $\DC$ et pour tout $n \in \Zm$, on a un triangle
distingué
$$\tau_{\leq n} \GC \longto \GC \longto \tau_{\geq n+1} \GC \longmapright{+1}$$
de sorte que si $\GC$ est un objet de $\DC^{\geq n}$, on a $[\GC]=(-1)^n [\lexp p h^n
\GC]+[\tau_{n+1} \GC]$
car $\tau_{\leq n} \GC=\tau_{\leq n} \tau_{\geq n} \GC=(\lexp p h^n \GC)[n]$. Soit alors
$a$ et $b$ tel que
$\FC$ soit un objet de $\DC^{[a,b]}$. En appliquant ce qui précède à $\tau_{\geq n} \FC$
pour $n$ variant de
$a$ à $b$, on obtient l'égalité
$$[\FC]=\sum_i (-1)^i [\lexp p h^i \FC],$$
de sorte que l'application $\sum_i \alpha_i [P_i] \in \groth(\CC) \mapsto \sum_i
\alpha_i [P_i] \in K(\DC)$ est inverse de celle de l'énoncé.
\end{proof}

\begin{coro}
Soit $D_b^c(X,\bar \Qm_l)$ la catégorie dérivée des $\bar \Qm_l$-complexes
constructibles sur un $\Fm_q$-schéma $X$ que l'on muni de la $t$-structure
de perversité autoduale (resp. de la $t$-structure triviale) de coeur la catégorie
$\perv(X)$ (resp. $\const(X)$) des faisceaux pervers (resp. des faisceaux
constructibles) sur $X$. Pour tout objet $\FC$ de $D_b^c(X,\bar \Qm_l)$, son image dans
$\groth(\perv(X))$ est déterminée par son image
dans $\groth(\const(X))$.
\end{coro}

\begin{lemm} Soit $j:U \hookrightarrow X$ l'inclusion d'un ouvert $U$ d'un
$\Fm_q$-schéma $X$; on note $i:Z \hookrightarrow X$ le fermé complémentaire.
L'image d'un objet $\FC$ de $D_b^c(X,\bar \Qm_l)$ dans $\groth(\const(X))$ est
déterminée par l'image de $j^* \FC$ dans $\groth(\const(U))$ et celle
de $i^* \FC$ dans $\groth(\const(Z))$.
\end{lemm}

\begin{proof} Le résultat découle directement de l'existence du triangle distingué
$j_! j^* \FC \longto \FC \longto i_*i^* \FC \longmapright{+1}$.

\end{proof}

On en déduit alors la proposition suivante.

\begin{prop} \label{prop-perv1}
Soit $X$ un schéma muni d'une stratification
$$X^{\geq d} \subset X^{\geq d-1} \subset \cdots \subset X^{\geq 1}=X$$
et soit $\PC$ un faisceau pervers sur $X$. Alors l'image de $\PC$ dans le groupe de
Grothendieck des faisceaux pervers sur $X$ est
déterminée par l'image de $\sum_i (-1)^i [h^i(\PC)_{|X^{=h}}$ dans le groupe de
Grothendieck des faisceaux localement constant sur
$X^{=h}:=X^{\geq h}-X^{\geq h-1}$ pour tout $1 \leq h \leq d$.
\end{prop}

\subsection{Image dans $\GF$ de $R\Psi_{\eta_o}(\bar \Qm_l)[d-1]$}

\begin{prop}
Les complexes $Rj^{\geq lg}_* \FC(g,l,\pi_o)_0[d-lg]$ et $j^{\geq lg}_!
\FC(g,l,\pi_o)_0[d-lg]$ sont des
faisceaux pervers.
\end{prop}

\begin{proof} Le résultat découle du fait que $j^{\geq lg}$ est une inclusion affine.

\end{proof}

\begin{prop} \label{prop-libre}
Pour $\pi_o$ une représentation irréductible cuspidale de $GL_g(F_o)$, on a l'égalité
suivante dans $\GF$:
\begin{multline} \label{egalite-groth}
e_{\pi_o}[R\Psi_{\eta_o,\pi_o}(\bar \Qm_l)[d-1]]= \sum_{i=1}^{s_g} \sum_{l=i}^{s_g}
(-1)^{l-i} \\
[j^{\geq lg}_! HT(g,l,\pi_o, [\overleftarrow{i-1},\overrightarrow{l-i}]_{\pi_o})
\otimes
L_g(\pi_o)(-\frac{lg-2+2i-l}{2})]
\end{multline}
\end{prop}

\begin{proof} L'isomorphisme (\ref{iso1}) et l'égalité (\ref{inutile}) déterminent pour
tout $0 \leq h < d$, la somme suivante
$$\sum_i (-1)^i [R^i\Psi_{v,|M_{I,s_o}^{=h}}]$$
dans le groupe de Grothendieck de la catégorie abélienne des faisceaux constructibles
sur $M_{I,s_o}^{=h}$. Le
résultat découle alors de la proposition (\ref{prop-perv1}).

\end{proof}

\subsection{Décomposition dans $\GF$ des $j_!^{\geq lg} \FC(g,l,\pi_o)[d-lg]$}
\label{bf1}

Rappelons le théorème de comparaison de Berkovich-Fargues

\begin{theo} Soit $z$ un point géométrique de $M_{I,s_o}^{=h}$. Le germe en $z$ de
$R\Psi_{\eta_o,I}(\bar
\Qm_l)$ est égal à $\widetilde{\Psi_{F_o,n}^{h}}$ en tant qu'objet de la catégorie
dérivée filtrée $\Dm^b
F(\bar \Qm_l)$, avec $n=\mult_o(I)$.
\end{theo}

Ainsi l'hypothèse de récurrence sur la filtration de monodromie-locale des modèles de
Deligne-Carayol de
hauteur strictement inférieure à $d$, donne le corollaire suivant.

\begin{coro} \label{coro-comparaison-grk}
Pour tout point géométrique $z_I$ de $M_{I,s_o}^{=h}$, le germe en $z_I$ de $h^i
gr_{k,\pi_o,I}$ vérifie les
propriétés suivantes:

\begin{itemize}
\item[(1)] il est nul si $h$ n'est pas divisible par $g$;

\item[(2)] pour $h=lg$, ils sont nuls pour $|k| \geq l$ quelque soit $i$;

\item[(3)] pour $h=lg$ et $|k| < l$:

(i) ils sont nuls pour $i < lg-d-l+1+|k|$ ou pour $i \not \equiv lg-d-l+1+k \mod 2$ ou
pour $i>lg-d$;

(ii) pour $lg-d-l+1+|k| \leq i=lg-d-l+1+|k|+2r \leq lg-d$, la fibre en $z_I$ de $h^i
gr_{k,\pi_o,I}$ est
naturellement munie d'une action de $\NC_o(h) \cap \HC_{o,n}(h)$ où $\NC_o(h) \subset
GL_h(F_o) \times
D_{o,h}^\times \times W_o$ est défini comme au \S \ref{rapel-DC}, où $GL_h(F_o)$ est vu
naturellement comme
un sous-groupe du Levi $GL_h(F_o) \times GL_{d-h}(F_o)$ de $GL_d(F_o)$ et où
$\HC_{o,n}(h)$ est l'algèbre de
Hecke associée à $K_{o,n}(h) \subset GL_h(\OC_o)$ avec $n=\mult_o(I)$. On obtient alors
que la fibre en $z_I$
de $(h^i gr_{k,\pi_o,I})^{e_{\pi_o}}$ est isomorphe aux invariants sous $K_{o,n}$ de
$$\JL^{-1}([\overleftarrow{l-1}]_{\pi_o}) \otimes [\overleftarrow{|k|+2r}]_{\pi_o}
\overrightarrow{\times}
[\overrightarrow{l-2r-|k|-2}]_{\pi_o} \otimes L_g(\pi_o)(-\frac{l(g-1)+|k|+k+2r}{2})$$
en tant que $\NC_o(h)
\cap \HC_{o,n}(h)$ module.

\end{itemize}
\end{coro}

\rem En ce qui concerne le dernier point du corollaire précédent, et donc le fait que la
filtration locale
soit équivariante, on peut la déduire du théorème de comparaison de Berkovich-Fargues
avec les strates des
variétés de Drinfeld-Stuhler en rang $d'<d$ en supposant connus, par récurrence, les
théorèmes globaux pour
celles-ci.

\begin{prop} \label{prop-p}
Pour tout $1 \leq l \leq s_g$, on a l'égalité dans $\GF$
$$j^{\geq lg}_! HT(g,l,\pi_o,\Pi_l)[d-lg] = \sum_{i=l}^{s_g} j^{\geq ig}_{!*}
HT(g,i,\pi_o,\Pi_l \overrightarrow{\times}
[\overleftarrow{i-l-1}]_{\pi_o})[d-ig] \otimes \Xi^{\frac{(l-i)(g-1)}{2}} \oplus
P_{!,l}$$ où $P_{!,l}$ est
une somme de faisceaux pervers concentrés aux points supersinguliers.
\end{prop}

\begin{proof} Dans $\GF$, on écrit $j^{\geq lg}_! HT(g,l,\pi_o,\Pi_l)[d-lg]$ sous la
forme
$\sum_i A_{\pi_o,l,i}(\Pi_l)$ où les $A_{\pi_o,l,i}(\Pi_l)$ sont une somme à
coefficients positifs de
faisceaux pervers simples de poids $d-lg-i$, que l'on appellera les constituants de
$j^{\geq lg}_!
HT(g,l,\pi_o,\Pi_l)[d-lg]$. On raisonne par récurrence descendante sur la dimension des
faisceaux pervers qui
interviennent, en traitant tous les $l$ par récurrence de $s_g$ à $1$ et en travaillant
sur les systèmes
inductifs comme indiqué dans la remarque (\ref{rema-lim}).

\begin{lemm} \label{lem-fait0}
Pour tout $1 \leq l \leq s_g$, $j^{\geq lg}_{!*} HT(g,l,\pi_o,\Pi_l)[d-lg]$ est le seul
constituant de
dimension $d-lg$ de $j^{\geq lg}_! HT(g,l,\pi_o,\Pi_l)[d-lg]$. Par ailleurs tous les
autres constituants
de $j^{\geq lg}_! HT(g,l,\pi_o,\Pi_l)[d-lg]$ sont de poids strictement inférieur à
$d-lg$.
\end{lemm}

\begin{proof} On a la suite exacte courte de faisceaux pervers
$$0 \to P_{\pi_o,l,0}(\Pi_l) \longto j^{\geq lg}_! HT(g,l,\pi_o,\Pi_l)[d-lg] \longto
j^{\geq lg}_{!*} HT(g,l,\pi_o,\Pi_l)[d-lg] \to 0$$ où $P_{\pi_o,l,0}$ est égale à
$i^{\geq lg,*}
j^{\geq lg}_{!*} HT(g,l,\pi_o,\Pi_l)[d-lg-1]$ qui est donc de poids inférieur ou égal à
$d-lg-1$ et
de dimension strictement inférieure à $d-lg$, d'où le résultat.

\end{proof}

D'après le lemme précédent, tous les faisceaux pervers qui interviennent dans l'écriture
de $j^{\geq lg}_!
HT(g,l,\pi_o,\Pi_l)[d-lg]$ sont de dimension inférieure ou égale à $d-lg$. Ainsi pour
$h>d-g$, il n'y a aucun
constituant de dimension $h$ dans les $j^{\geq lg}_! HT(g,l,\pi_o,\Pi_l)[d-lg]$ pour $1
\leq l \leq s_g$ et
pour $h=d-g$, $j^{\geq g}_{!*} HT(g,1,\pi_o,\Pi_1)[d-g]$ est le seul constituant de
dimension $d-g$ pour
$l=1$. Supposons donc que pour tout $h > h_0$, les constituants des $j^{\geq lg}_!
HT(g,l,\pi_o,\Pi_l)[d-lg]$
pour $1 \leq l \leq s_g$, sont ceux prédits par l'énoncé et traitons la dimension
$h_0$.

\begin{lemm} \label{lem-hok}
Supposons la proposition (\ref{prop-p}) vérifiée pour les faisceaux pervers de dimension
strictement
supérieure à $0<h_0 < d-g$, i.e. pour tout $1 \leq l \leq s_g$, les constituants de
dimension strictement
plus grande que $h_0$ des $j^{\geq lg}_! HT(g,l,\pi_o,\Pi_l)[d-lg]$ sont les
$$j^{\geq (l+r)g}_{!*} HT(g,l+r,\pi_o,\Pi_l \overrightarrow{\times}
[\overleftarrow{r-1}]_{\pi_o})
[d-lg] \otimes \Xi^{\frac{r(g-1)}{2}}$$ avec $d-(l+r)g>h_0$. Pour tout $k$, l'image dans
$\GF$ de $e_{\pi_o}
gr_{k,\pi_o}$ est alors égale à
$$\sum_{\genfrac{}{}{0pt}{}{|k| < l < (d-h_0)/g}{l \equiv k-1 \mod 2}} \PC(g,l,\pi_o)
(-\frac{lg+k-1}{2}) +P_{k,h_0}$$
où $P_{k,h_0}=(P_{k,h_0,I})_I$ est une somme de faisceaux pervers simples de dimension
inférieure ou égale à
$h_0$ à support dans la tour des $(M_{I,s_o}^{\geq d-h_0})_I$, i.e. la fibre de
$h^iP_{k,h_0,I}$ en tout
point géométrique de $M_{I,s_o}^{=h}$ est nulle pour $h<d-h_0$ et pour tout $i$.
\end{lemm}

\begin{proof} D'après (\ref{prop-libre}), chacun\footnote{On rappelle, cf. la remarque
(\ref{rema-I}), que les
$\PC(g,l,\pi_o)$ ne sont pas irréductibles; cependant les arguments fonctionnent de
manière strictement
identique pour tous ses constituants simples car la seule différence entre ceux-ci
provient de l'action de
$\DC_{o,lg}^\times$.} des $\PC(g,l,\pi_o) (-\frac{lg+k-1}{2})$, pour $d-lg>h_0$ est un
constituant d'un
$e_{\pi_o}gr_{r,\pi_o}$; il s'agit alors de montrer que $r=k$. Le résultat découle
directement via le
théorème de comparaison de Berkovich-Fargues, de l'hypothèse de récurrence sur le
théorème local
(\ref{theo-local-fil}). En effet, pour $lg<d$ les $\pi_o$-parties des $h^0
gr_{lg,loc,k}$ sont pures de poids
$lg-1+k$. L'équivalent du théorème de Serre-Tate et le théorème de Berkovich-Fargues
impliquent alors que la
fibre en tout point géométrique $z$ de la strate $lg$ de $h^{lg-d} gr_{k,\pi_o}$ est
pure de poids $lg-1+k$.
Or la fibre en un point géométrique de la strate $lg$ de $h^{lg-d}
\PC(g,l,\pi_o)(-\frac{lg+k-1}{2})$ est de
poids $lg-1+k$. Ainsi pour $l=1$, $\PC(g,1,\pi_o)((1-g)/2)$ est un constituant de
$gr_{0,\pi_o}$ car tous les
autres faisceaux pervers qui interviennent sont de dimension strictement plus petite.

On raisonne alors par récurrence sur $l$; on suppose que pour tout $l <l_0<(d-h_0)/g$,
et tout $|k|<l-1$,
$\PC(g,l,\pi_o)(-\frac{lg+k-1}{2})$ est un constituant de $gr_{k,\pi_o}$ et traitons le
cas de $l_0$. D'après
ce qui précède on en déduit que pour tout $|k|<l_0$,
$\PC(g,l_0,\pi_o)(-\frac{l_0g-1+k}{2})$ est un
constituant de $gr_{r,\pi_o}$ avec $r \geq k$; en effet dans le cas contraire la
filtration par le poids de
$gr_r$ donnerait une suite exacte courte
$$0 \to P \longto gr_{r,\pi_o} \longto \PC(g,l_0,\pi_o)(-\frac{l_0g-1+k}{2}) \to 0$$
modulo des faisceaux pervers de dimension strictement plus petite, ce qui donnerait une
suite exacte longue
$$ \cdots h^{l_0g-d} gr_{r,\pi_o} \longto h^{l_0g-d}
\PC(g,l_0,\pi_o)(-\frac{l_0g-1+k}{2}) \longto h^{l_0g-d+1} P
\cdots$$ or $h^{lg-d+1} P$ a un support de dimension strictement plus petite que
$d-l_0g$ de sorte qu'il
existerait un point géométrique de $M_{I,s_o}^{=l_0g}$ tel que la fibre de $h^{l_0g-d}
gr_{r,\pi_o}$ serait
de poids $l_0g-1+k$ ce qui n'est pas. Ainsi en particulier pour $k=l_0-1$, il existe $r
\geq l_0-1$ tel que
$gr_{r,\pi_o}$ soit de dimension $d-l_0g$. Pour la même raison que ci-dessus, on en
déduit qu'il s'agit de
$r=l_0-1$. On utilise alors l'opérateur $N$ qui finit de placer les
$\PC(g,l_0,\pi_o)(-\frac{l_0g+k-1}{2})$
sur les $gr_{k,\pi_o}$, d'où la récurrence.

Par ailleurs on remarque qu'en tout point géométrique $z$ de $M_{I,s_o}^{\leq d-h_0}$,
la fibre en $z$ des
faisceaux de cohomologie de $gr_{k,\pi_o,I}-P_{k,h_0,I}$ est égale à celle du modèle
local i.e., d'après le théorème de comparaison de
Berkovich-Fargues à celle de $gr_{k,\pi_o,I}$ de sorte que celle de $P_{k,h_0,I}$ est
nulle.

\end{proof}

\noindent \textit{Retour à la preuve de la proposition (\ref{prop-p})}: \\ le principe
est d'étudier
l'égalité (\ref{egalite-groth}) en utilisant le lemme précédent ainsi que
(\ref{lem-fait0}) pour $j^{\geq
l_0g}_{!} HT(g,l_0,\pi_o,\Pi_{l_0})[d-l_0g]$ dont le seul constituant de dimension
supérieur ou égal à
$d-l_0g$ est $j^{\geq l_0g}_{!*} HT(g,l_0,\pi_o,\Pi_{l_0})[d-l_0g].$

\noindent \textbf{Note:} \textit{on suppose, afin de simplifier la rédaction, que
$\Pi_l$ est elliptique de
type $\pi_o$.}

\begin{lemm} \label{lem-fp-rpsi}
Pour tout $1 \leq l \leq s_g$, soit $A_{\pi_o,l,r}(\Pi_l)$ un constituant de
$$j^{\geq lg}_! HT(g,l,\pi_o,\Pi_l)[d-lg]$$
de poids $d-lg-r$ et de dimension $h_0$. Il existe alors des entiers $k$ et $t \geq 0$
ainsi que $\Pi_l'$ une
représentation elliptique de type $\pi_o$, et donc de même support cuspidal que $\Pi_l$,
telle que
$A_{\pi_o,l,r}(\Pi_l') \otimes L_g(\pi_o)(-\frac{l(g+1)-2(t+1)}{2})$ soit un constituant
de $e_{\pi_o}
gr_{k,\pi_o}$ \footnote{En utilisant monodromie-poids, on obtient $k=l-1-2t$.}.
\end{lemm}

\begin{proof} Soit $A_{\pi_o,l,r}(\Pi_l)$ un constituant de $j^{\geq lg}_!
HT(g,l,\pi_o,\Pi_l)[d-lg]$ de sorte que
$A_{\pi_o,l,r}([\overleftarrow{l-1}]_{\pi_o}) \otimes L_g(\pi_o)(-\frac{l(g+1)-2}{2})$
apparaît dans le
membre de droite de l'égalité de la proposition (\ref{prop-libre}). On considère alors
$\d$ maximal tel qu'il
existe une représentation elliptique $\Pi_l'$ de type $\pi_o$ de $GL_{lg}(F_o)$ tel que
$A_{\pi_o,l,r}(\Pi_l')
\otimes L_g(\pi_o)(-\frac{\d}{2})$ apparaisse dans le membre de droite de l'égalité de
la proposition
(\ref{prop-libre}). On remarque alors qu'il y apparaît avec un coefficient positif et
aucun négatif de sorte
qu'il reste dans la somme ce qui implique le résultat. En effet s'il apparaissait avec
un coefficient
négatif, ce serait comme constituant d'un $j^{\geq l'g}_!
HT(g,l',\pi_o,[\overleftarrow{l'-r-1},\overrightarrow{r}]_{\pi_o}) [d-l'g] \otimes
L_g(\pi_o)
(-\frac{l'(g+1)-2(r+1)}{2})$ pour $r$ impair positif et on note que $j^{\geq l'g}_!
HT(g,l',\pi_o,[\overleftarrow{l'-1}]_{\pi_o})[d-l'g] \otimes
L_g(\pi_o)(-\frac{l'(g+1)-2}{2})$ contiendrait
un $A_{\pi_o,l,r}(\Pi_l'') \otimes L_g(\pi_o)(-\frac{\d+2r}{2})$ pour $\Pi_l''$ une
représentation de même
support cuspidal que $\Pi_l'$, contredisant la maximalité de $\d$.

\end{proof}

- Supposons dans un premier temps que $h_0$ n'est pas de la forme $d-lg$ et qu'il existe
$l,r$ tel que
$A_{\pi_o,l,r}(\Pi_l)$ contienne un faisceau pervers simple de dimension $h_0$. Le lemme
précédent implique
alors qu'il existe $k \geq 0$ tel que, avec les notations du lemme (\ref{lem-hok}),
$P_{k,h_0}$ ait un
constituant simple, disons $B_{\pi_o,k}$, de dimension $h_0$ et de poids $d-1+t$ pour un
certain entier $t$.
On raisonne alors comme dans la preuve du lemme (\ref{lem-hok}). Étant de dimension
inférieure ou égale à
$h_0$, $P_{k,h_0}$ est alors de dimension $h_0$ de sorte que $h^{-h_0} P_{k,h_0}$ est de
dimension $h_0$. En
outre comme $h^{-h_0} P_{k,h_0,I}$, d'après le lemme (\ref{lem-hok}), est à support dans
$M_{I,s_o}^{\geq
d-h_0}$, on en déduit que pour tout point générique de $M_{I,s_o}^{=d-h_0}$, la fibre en
ce point de
$h^{-h_0} P_{k,h_0,I}$ est non nulle de sorte qu'il existe un point géométrique $z_I$ de
$M_{I,s_o}^{=d-h_0}$
tel que la fibre en $z_I$ de $h^{-h_0} P_{k,h_0,I}$ soit non nulle. Par ailleurs, les
$\pi_o$-parties des
$h^{-h_0}gr_{k,loc}$ du modèle local de hauteur $d-h_0$ sont toutes nulles, on en
déduit, d'après le théorème
de comparaison de Berkovich-Fargues, que la fibre en $z$ de $h^{-h_0} gr_{k,\pi_o}$ est
nulle ce qui implique
que $t \neq k$: en effet sinon $e_{\pi_o} gr_{k,\pi_o}$ est une extension de
$$(\bigoplus_{\genfrac{}{}{0pt}{}{|k| < l < (d-h_0)/g}{l \equiv k-1 \mod 2}}
\PC(g,l,\pi_o) (-\frac{lg+k-1}{2}))
 \oplus B_{k,\pi_o}$$
par des faisceaux pervers de dimension inférieure ou égale à $h_0$ ce qui implique que
$(h^{-h_0}
gr_{k,\pi_o})^{e_{\pi_o}}$ admettrait $h^{-h_0} B_{k,h_0}$ comme facteur direct. Selon
le même principe, on
doit même avoir $t <k$. Or d'après la pureté des $gr_k$, on doit avoir $t=k$ d'où le
résultat. Sans utiliser
monodromie-poids, on raisonne alors comme suit.

\begin{lemm} \label{lem-dualite}
Soit $\PC$ un constituant de $R\Psi_{\eta_o,\pi_o}(\bar \Qm_l[d-1])$. On en
déduit alors que $D \PC(1-d)$ est un constituant de $R\Psi_{\eta_o,\pi_o^\vee}(\bar
\Qm_l[d-1])$.
\end{lemm}

\begin{proof} Le résultat découle simplement de la compatibilité de $R\Psi_{\eta_o}$
avec la dualité de Verdier,
soit
$$D R\Psi_{\eta_o}(\bar \Qm_l[d-1]) \simeq R\Psi_{\eta_o}(D \bar \Qm_l[d-1])$$
et du fait que la fibre générique de $M_{I,o}$ est lisse de sorte que $D \bar \Qm_l[d-1]
\simeq \bar
\Qm_l[d-1](1-d)$.

\end{proof}

On considère alors un constituant simple de dimension $h_0$ de poids maximal, disons
$d-1+r$ de
$R\Psi_{\eta_o,\pi_o}(\bar \Qm_l)[d-1]$. C'est alors un constituant de $gr_{k,\pi_o}$
pour $k > r$ d'après ce
qui précède. Par ailleurs en utilisant le lemme (\ref{lem-dualite}), on obtient un
constituant simple de
dimension $h_0$ de $gr_{-k,\pi_o}$ et de poids $d-1-r$, de sorte que l'opérateur de
monodromie $N^k$ fournit
un constituant simple de dimension $h_0$ de $gr_{k,\pi_o}$ et de poids $d-1-r+2k >
d-1+r$, d'où la
contradiction par maximalité de $r$.

- Supposons désormais que $h_0=d-l_0g$ et supposons avoir montré par récurrence que pour
tout $l_1 \leq l
\leq s_g$, le seul constituant de dimension $d-l_0g$ de $j^{\geq lg}_!
HT(g,l,\pi_o,\Pi_l)[d-lg]$, est
$j^{\geq l_0g}_{!*} HT(g,l_0g,\pi_o,\Pi_l \overrightarrow{\times}
[\overleftarrow{l_0-l-1}]_{\pi_o})[d-l_0g]
\otimes \Xi^{\frac{(l_0-l)(g-1)}{2}}$. Le résultat est vérifié pour $l_1=l_0$, supposons
le donc vérifié
jusqu'au rang $l_1$ et traitons le cas de $l_1-1$. On étudie alors les faisceaux pervers
de dimension
$d-l_0g$ dans le membre de droite de (\ref{egalite-groth}), en particulier ceux de poids
$d-l_0$, ce qui
donne $j^{\geq l_0g}_{!*} HT(g,l_0,\pi_o,\Pi)[d-l_0g] \otimes L_g(\pi_o)
(-\frac{l_0(g-1)}{2})$ avec
$$(-1)^{l_0-1}\Pi=[\overrightarrow{l_0-1}]_{\pi_o} - \sum_{i=1}^{l_0-l_1} (-1)^i
[\overrightarrow{l_0-i-1}]_{\pi_o}
\overrightarrow{\times} [\overleftarrow{i-1}]_{\pi_o}$$ soit $\Pi=(-1)^{l_1-1}
[\overrightarrow{l_1-1},\overleftarrow{l_0-l_1}]_{\pi_o}$. Or ce dernier ne peut pas
apparaître dans le
membre de gauche de (\ref{egalite-groth}). En effet si $l_1$ est pair c'est évident car
à gauche il ne peut y
avoir que des coefficients positifs. Sinon de manière générale, on raisonne comme suit.
Le lemme
(\ref{lem-dualite}) donnerait que $j^{\geq l_0g}_{!*} HT(g,l_0,\pi_o^\vee,
[\overleftarrow{l_0-l_1},\overrightarrow{l_1-1}]_{\pi_o^\vee})[d-l_0g] \otimes
L_g(\pi_o^\vee)(-\frac{l_0(g+1)-2}{2})$ serait un constituant de
$e_{\pi_o}[R\Psi_{\eta_o,\pi_o}(\bar
\Qm_l)[d-1]]$ qui, vu le poids, ne pourrait provenir que de
$$j^{\geq l_0g}_! HT(g,l_0,\pi_o^\vee,[\overleftarrow{l_0-1}]_{\pi_o^\vee})[d-l_0g]
\otimes L_g(\pi_o^\vee)
(-\frac{l_0(g+1)-2}{2}),$$ ce qui n'est pas d'après le lemme (\ref{lem-fait0}).
\footnote{En utilisant la
pureté des gradués de la filtration de monodromie, on peut argumenter comme suit: le
faisceau pervers en
question serait alors un constituant de $e_{\pi_o}gr_{1-l_0,\pi_o}$ et donc $j^{\geq
l_0g}_{!*}
HT(g,l_0,\pi_o, [\overrightarrow{l_1-1},\overleftarrow{l_0-l_1}]_{\pi_o}) \otimes
L_g(\pi_o)
(-\frac{l_0(g+1)-2}{2})$ serait un constituant de $e_{\pi_o}[gr_{l_0-1,\pi_o}]$ qui ne
pourrait être obtenu
que via $j^{\geq l_0g}_! HT(g,l,\pi_o, [\overleftarrow{l_0-1}]_{\pi_o})[d-l_0g] \otimes
L_g(\pi_o)
(-\frac{l_0(g+1)-2}{2})$ ce qui ne se peut pas d'après le lemme (\ref{lem-fait0})}.

Ainsi $j^{\geq l_0g}_{!*}
HT(g,l_0,\pi_o,[\overrightarrow{l_1-1},\overleftarrow{l_0-l_1}]_{\pi_o})[d-l_0g]
\otimes L_g(\pi_o) (-\frac{l_0(g-1)}{2})$ doit être un constituant d'un $j^{\geq lg}_!
HT(g,l,\pi_o,[\overleftarrow{i},\overrightarrow{l-i-1}]_{\pi_o})[d-lg] \otimes
L_g(\pi_o)(-\frac{l(g-1)+2i}{2})$ pour $0 \leq i < l$ et $l>l_0$ et $l \equiv l_1-1 \mod
2$. Le résultat,
i.e. $l=l_1-1$ et $r=0$, découle alors des trois lemmes suivants.

\begin{lemm} \label{lem-induite}
Pour un point géométrique $z$ de $M_{I,s_o,1}^{=l_0g}$, les fibres en $z$ des faisceaux
de cohomologies des
$j^{\geq lg}_{!*} HT(g,l,\pi_o,\Pi_l)[d-lg] \otimes L_g(\pi_o)$, pour $l \geq l_0$,
sont, en tant que
représentation de $GL_{l_0g}(F_o) \times GL_{d-l_0g}(F_o) \times W_o$ de la forme
$$\bigoplus_\xi (\Pi_l(\xi) \times \pi_\xi) \otimes \pi_\xi' \otimes L_g(\pi_o)(\xi)$$
où $\xi$ décrit les caractères de $\Zm$ et $\pi_\xi$ (resp. $\pi_\xi'$) est une
représentation de
$GL_{(l_0-l)g}(F_o)$ (resp. $GL_{d-l_0g}(F_o)$).
\end{lemm}

\begin{proof} C'est évident en utilisant que les strates sont induites, i.e.
$$j^{\geq lg}_{!*} HT(g,l,\pi_o,\Pi_l)= j^{\geq lg}_{1,!*} \FC(g,l,\pi_o)_1 \otimes
\Pi_l
\times_{P_{lg,l_0g,d}(F_o)} P_{l_0g,d}(F_o)$$ en tant que $P_{l_0g,d}(F_o) \times
\Zm$-module. Les torsions
découlent alors de l'action telle qu'elle est décrite au \S \ref{ht-prop}.

\end{proof}

\begin{lemm} Supposons que pour tout $1 \leq l \leq l_0$, les constituants simples de
dimension strictement supérieure
à $d-l_0g$ des $j^{\geq lg}_{!*} HT(g,l,\pi_o,\Pi_l)[d-lg]$ sont ceux prévus par la
proposition
(\ref{prop-p}) et supposons qu'il existe $l<l_0$ tel que
$$j^{\geq lg}_! HT(g,l,\pi_o,[\overleftarrow{r},\overrightarrow{l-r-1}]_{\pi_o})[d-lg]
\otimes
L_g(\pi_o)(-\frac{l(g-1)+2r}{2})$$ admette $j^{\geq l_0g}_{!*}
HT(g,l_0,\pi_o,[\overrightarrow{l_1-1},\overleftarrow{1},\overleftrightarrow{l_0-l_1-1}]_{\pi_o})[d-l_0g]
\otimes L_g(\pi_o) (-\frac{l_0(g-1)}{2})$ comme constituant. On a alors $l \leq l_1$ et
$r=0$.
\end{lemm}

\begin{proof} On considère la filtration par le poids de
$$j^{\geq lg}_!
HT(g,l,\pi_o,[\overleftarrow{r},\overrightarrow{l-r-1}]_{\pi_o})[d-lg] \otimes
L_g(\pi_o)(-\frac{l(g-1)+2r}{2})$$ et la suite spectrale qui s'en déduit

\begin{multline*}
E_1^{i,j}=h^{i+j} gr_{-i}(l,r) \Rightarrow  h^{i+j} j^{\geq lg}_!
HT(g,l,\pi_o,[\overleftarrow{r},
\overrightarrow{l-r-1}]_{\pi_o})[d-lg] \otimes L_g(\pi_o)(-\frac{l(g-1)+2r}{2})
\end{multline*}
où $gr_k(l,r)$ est le gradué de poids $k$ de $j^{\geq lg}_!
HT(g,l,\pi_o,[\overleftarrow{r},\overrightarrow{l-r-1}]_{\pi_o})[d-lg] \otimes
L_g(\pi_o)(-\frac{l(g-1)+2r}{2})$. D'après le lemme précédent toutes les fibres aux
points géométriques de
$M_{I,s_o,1}^{=l_0g}$ des $E_1^{i,j}$ avec $i+j<l_0g-d$ sont, en tant que
$GL_{l_0g}(F_o) \times W_o$-module,
de la forme
$$\bigoplus_\xi [\overleftarrow{r}, \overrightarrow{l-r-1}]_{\pi_o(\xi)} \times \pi_\xi
\otimes
L_g(\pi_o(\xi))(-\frac{l(g-1)+2r}{2})$$ alors que celles de $j^{\geq l_0g}_{!*}
HT(g,l_0,\pi_o,[\overrightarrow{l_1-1},\overleftarrow{1},\overleftrightarrow{l_0-l_1-1}]_{\pi_o})[d-l_0g]
\otimes L_g(\pi_o) (-\frac{l_0(g-1)}{2})$ sont de la forme
$$\bigoplus_\xi
[\overrightarrow{l_1-1},\overleftarrow{1},\overleftrightarrow{l_0-l_1-1}]_{\pi_o(\xi)}
\times
\pi'_\xi \otimes L_g(\pi_o(\xi))(-\frac{l(g-1)}{2})$$ En remarquant que les
$E_1^{i,l_0g-d+1-i}$ ont un
support de dimension strictement inférieur à $d-l_0g$ et que les $E_\oo^n$ sont nuls
pour $n \neq lg-d$, on
en déduit alors que $r=0$ et $l \leq l_1$.

\end{proof}

\begin{lemm} \label{lem-combi}
Pour tout $1 \leq l < l_0$ et toute représentation $\Pi$ elliptique de type $\pi_o$ de
$GL_{(l_0-l)g}(F_o)$
distincte de $[\overleftarrow{l_0-l-1}]_{\pi_o}$, $j^{\geq lg}_!
HT(g,l,\pi_o,\Pi_l)[d-lg]$ ne contient pas
$$j^{\geq l_0g}_{!*} HT(g,l_0,\pi_o,\Pi_l \overrightarrow{\times} \Pi)[d-lg] \otimes
\Xi^{(l_0-l)(g-1)/2}.$$
\end{lemm}

\begin{proof} Dans le cas contraire considérons $l$ minimal pour cette propriété. En
l'appliquant à
$\Pi_l=[\overrightarrow{l-1}]_{\pi_o}$, on en déduit que
$$j^{\geq l_0g}_{!*}
HT(g,l_0,\pi_o,[\overrightarrow{l-1},\overleftarrow{1},\overleftrightarrow{l_0-l-1})[d-lg]
\otimes
L_g(\pi_o)(-\frac{l_0(g-1)}{2})$$ reste dans le membre de droite de
(\ref{egalite-groth}), où
$[\overleftrightarrow{l_0-l-1}]_{\pi_o}$ désigne $\Pi$. En effet il y apparaît via
$j^{\geq lg}_{!*}
HT(g,l,\pi_o,[\overrightarrow{l-1}]_{\pi_o})[d-lg] \otimes
L_g(\pi_o)(-\frac{l(g-1)}{2})$ et n'est pas
compensé car d'après le lemme précédent ce ne pourrait qu'être pour un $l'<l$ ce qui
contredirait la
minimalité de $l$. Si le signe est négatif on obtient directement la contradiction,
sinon on argument comme
suit: en utilisant la pureté (resp. sans utiliser la pureté), par application de
$N^{l_0-1}$ (resp. de la
dualité de Verdier), on en déduit que $e_{\pi_o}[gr_{l_0-1,\pi_o}]$ (resp. $e_{\pi_o}
[gr_{k,\pi_o^\vee}]$
pour un certain $k$) devrait contenir
$$j^{\geq l_0g}_{!*} HT(g,l_0,\pi_o,\lceil
[\overrightarrow{l-1}]_{\pi_o} \overrightarrow{\times} \Pi_l \rceil)[d-l_0g] \otimes L_g(\pi_o)
(-\frac{l_0(g+1)-2}{2})$$ (resp. $j^{\geq l_0g}_{!*} HT(g,l_0,\pi_o^\vee,\lceil
[\overrightarrow{l-1}]_{\pi_o} \overrightarrow{\times} \Pi_l \rceil^\vee)[d-l_0g] \otimes L_g(\pi_o^\vee)
(-\frac{l_0(g+1)-2}{2})$) qui est de poids $d+l_0-2$. Or tous les constituants de dimension strictement
supérieur à $d-l_0g$ de $R\Psi_{\eta_o,\pi_o}(\bar \Qm_l)[d-1]$ (resp. de $R\Psi_{\eta_o,\pi_o^\vee}(\bar
\Qm_l)[d-1]$ ) sont de poids strictement inférieur à $d+l_0-2$ de sorte qu'il existerait un point géométrique
de $M_{I,s_o}^{=l_0g}$ tel que la fibre de $h^{l_0g-d}gr_{l_0-1,\pi_o}$ (resp. $h^{l_0g-d}gr_{k,\pi_o^\vee}$)
en ce point admettrait un facteur direct de poids $d+l_0-2$ de la forme $\lceil
[\overrightarrow{l-1}]_{\pi_o} \overrightarrow{\times} \Pi_l \rceil$ (resp. $\lceil
[\overrightarrow{l-1}]_{\pi_o} \overrightarrow{\times} \Pi_l \rceil^\vee$) ce qui n'est pas d'après le
corollaire (\ref{coro-comparaison-grk}).

\end{proof}

\noindent \textit{Fin de la preuve de la proposition (\ref{prop-p})}: pour tout $k$, on
pose dans $\GF$:
$$Q_{k,\pi_o,l_0}:=e_{\pi_o} gr_{k,\pi_o}-\sum_{\genfrac{}{}{0pt}{}{|k| < l \leq l_0}{l
\equiv k-1 \mod 2}} \PC(g,l,\pi_o)
(-\frac{lg+k-1}{2})$$ Le corollaire (\ref{coro-comparaison-grk}) et ce qui précède,
prouvent alors que pour
tout point générique $z$ de dimension supérieure ou égale à $d-l_0g$, $\sum_i (-1)^i
(h^i Q_{k,\pi_o,l_0})_z$
est nulle. On en déduit donc d'après la proposition (\ref{prop-perv1}), que
$Q_{k,\pi_o,l_0}$ est de
dimension strictement inférieure à $d-l_0g$, ce qui conclut la récurrence.

\end{proof}

\begin{defi}
Pour $l \leq l'$ et $\Pi_l$ une représentation de $GL_{lg}(F_o)$, on introduit les
faisceaux pervers purs de
$\FPH(M_{s_o})$ \footnote{Comme dans la remarque (\ref{rema-I}), ceux-ci ne sont pas
simples mais plutôt une
somme directe de $e_{\pi_o}$ faisceaux pervers simples où la différence entre eux
provient de l'action de
$\DC_{o,lg}^\times$. Par ailleurs pour ce qui est des poids cf. la remarque
(\ref{rema-zw}).}
$$\PC_-(g,l',\pi_o,l,\Pi_l):= j^{\geq l'g}_{!*} HT(g,l',\pi_o,\Pi_l
\overrightarrow{\times}
[\overleftarrow{l'-l-1}]_{\pi_o}) [d-l'g] \otimes \Xi^{\frac{(l'-l)(g - 1)}{2}}$$ purs
de poids $d-l'g
-2(l'-l)$ et
$$\PC_+(g,l',\pi_o,l,\Pi_l):= j^{\geq l'g}_{!*} HT(g,l',\pi_o,\Pi_l
\overleftarrow{\times}
[\overleftarrow{l'-l-1}]_{\pi_o}) [d-l'g] \otimes \Xi^{\frac{(l'-l)(g + 1)}{2}}$$ purs
de poids $d-l'g +
2(l'-l)$.
\end{defi}

La filtration par le poids donne alors le corollaire suivant.

\begin{coro} \label{coro-sec-fp}
\begin{itemize}

\item[(i)] Il existe, pour $0 \leq i \leq s_g-l$, des faisceaux pervers
$P_{\pi_o,l,i}(\Pi_l)\in \FPH(M_{s_o})$
tels que l'on ait, dans $\FPH(M_{s_o})$, les suites exactes suivantes:

$$0 \to P_{\pi_o,l,0}(\Pi_l) \longto j^{\geq lg}_! HT(g,l,\pi_o,\Pi_l)[d-lg] \longto
j^{\geq lg}_{!*}
HT(g,l,\pi_o,\Pi_l)[d-lg] \to 0$$
$$0 \to P_{\pi_o,l,1}(\Pi_l) \longto P_{\pi_o,l,0}(\Pi_l) \longto
\PC_-(g,l+1,\pi_o,l,\Pi_l) \oplus A_{\pi_o,l,0}(\Pi_l)
\to 0$$
$$\cdots$$
$$0 \to P_{\pi_o,l,r}(\Pi_l) \longto P_{\pi_o,l,r-1}(\Pi_l) \longto
\PC_-(g,l+r,\pi_o,l,\Pi_l) \oplus A_{\pi_o,l,r-1}(\Pi_l)
 \to 0$$
$$\cdots$$
$$0 \to P_{\pi_o,l,s_g-l}(\Pi_l) \longto P_{\pi_o,l,s_g-l-1}(\Pi_l) \longto
\PC_-(g,s_g,\pi_o,l,\Pi_l) \oplus
A_{\pi_o,l,s_g-l-1}(\Pi_l) \longto 0$$ avec $A_{\pi_o,l,i}(\Pi_l)$, pour $1 \leq i \leq
s_g-l$, des faisceaux
pervers concentrés aux points supersinguliers purs de poids $d-lg-i-1$ et
$P_{\pi_o,l,s_g-l}(\Pi_l)$ un
faisceau pervers concentré aux points supersinguliers de poids inférieur ou égal à
$d-lg-s_g+l$.

\item[(ii)] Dualement pour la dualité de Verdier, le complexe $Rj^{\geq lg}_*
HT(g,l,\pi_o,\Pi_l)[d-lg]$
est un objet de $\FPH(M_{s_o})$ qui s'insère dans les suites exactes courtes
\begin{multline*}
0 \to j^{\geq lg}_{!*} HT(g,l,\pi_o,\Pi_l)[d-lg] \longto Rj^{\geq lg}_*
HT(g,l,\pi_o,\Pi_l)[d-lg] \\
\longto DP_{\pi_o,,l,0}(\Pi_l)(lg-d) \to 0
\end{multline*}
\begin{multline*} 0 \to DA_{\pi_o,l,0}(\Pi_l)(lg-d) \oplus \PC_+(g,l+1,\pi_o,l,\Pi_l)
\longto \\
DP_{\pi_o,l,0}(\Pi_l)(lg-d) \longto DP_{\pi_o,l,1}(\Pi_l)(lg-d) \to 0
\end{multline*}
$$\cdots$$
\begin{multline*} 0 \to DA_{\pi_o,l,r-1}(\Pi_l)(lg-d) \oplus \PC_+(g,l+r,\pi_o,l,\Pi_l)
\longto \\
DP_{\pi_o,l,r-1}(\Pi_l)(lg-d) \longto DP_{\pi_o,l,r}(\Pi_l)(lg-d) \to 0
\end{multline*}
$$\cdots$$
\begin{multline*} 0 \to DA_{\pi_o,l,s_g-l-1}(\Pi_l)(lg-d) \oplus
\PC_+(g,s_g,\pi_o,l,\Pi_l) \longto \\
DP_{\pi_o,l,s_g-l-1}(\Pi_l)(lg-d) \longto DP_{\pi_o,l,s_g-l}(\Pi_l)(lg-d) \to 0
\end{multline*}

\end{itemize}
\end{coro}

\rem Hors des points supersinguliers, nous avons montré que les $gr_k$ étaient purs de
poids $d-1+k$, ce qui
n'a rien d'impressionnant puisque finalement on l'a déduite de la pureté locale qui nous
est donnée d'après
l'hypothèse de récurrence. Nous verrons à la proposition (\ref{prop-mono}) comment la
démontrer pour les
faisceaux pervers supportés par les points supersinguliers.

\begin{coro} \label{coro-grk}
Pour tout $|k| < s_g$, l'image de $e_{\pi_o} gr_{k,\pi_o}$ dans $\GF$ est égale à
$$(\sum_{\genfrac{}{}{0pt}{}{|k| < l \leq s_g}{l \equiv k-1 \mod 2}}
\PC(g,l,\pi_o)(-\frac{lg+k-1}{2})) + P_k$$
où comme ci-dessus, $P_k$ est une somme de faisceaux pervers concentrés aux points
supersinguliers. Dans tous
les autres cas $gr_{k,\pi_o}$ est de dimension nulle, concentré aux points
supersinguliers.
\end{coro}

\subsection{Étude aux points supersinguliers des $j_!^{\geq lg} \FC(g,l,\pi_o)[d-lg]$}
\label{fp-ponctuel}

Le but de ce paragraphe est de déterminer les faisceaux pervers ponctuels non précisés
dans la proposition
(\ref{prop-p}). On rappelle que le raisonnement du paragraphe précédent ne s'appliquait
pas au niveau des
points supersinguliers car nous n'y connaissons pas l'aboutissement de
(\ref{suite-spectrale}). Une idée
naive est que pour connaître un faisceau ponctuel, on peut commencer par calculer son
groupe de cohomologie
$H^0$.

\begin{prop} \label{prop-coho1}
Soit $1 \leq g < d$ ne divisant pas $d$ (resp. $g$ divisant $d=sg$), et soit $\Pi$ une
représentation globale
de $D_\Am^\times$ telle que $\Pi$ vérifie $\hyp(\r_\oo)$ avec $\Pi_o \simeq
\st_s(\pi_o)$ pour $\pi_o$ une
représentation irréductible cuspidale de $GL_g(F_o)$. Pour tout $i$ et $1 \leq l \leq
s_g$ (resp. $1 \leq l <
s$), la composante $\Pi^{\oo,o}$-isotypique des groupes de cohomologie des faisceaux
pervers $j^{\geq
lg}_{!*} HT_{\rho_\oo}(g,l,\pi_o,\Pi_l)[d-lg]$ est nulle, soit avec les notations de
(\ref{nota-coho}):
$$H^i(j^{\geq lg}_{!*} HT_{\rho_\oo}(g,l,\pi_o,\Pi_l)[d-lg])[\Pi^{\oo,o}]=0.$$
\end{prop}

\begin{proof} On raisonne par récurrence pour $l$ variant de $s-1$ à $1$;
l'initialisation de la récurrence se fait
d'elle même dans la preuve de l'induction qui suit. On reprend les suites exactes
courtes de faisceaux
pervers du corollaire (\ref{coro-sec-fp}). On montre tout d'abord, par récurrence sur
$i$ de $s-l-1$ à $0$,
que pour tout $j \neq 0$, les groupes de cohomologie $H^j(P_{\pi_o,l,i}(\Pi_l)\otimes
\LC_{\r_\oo})[\Pi^{\oo,o}]$ sont nuls et que
\begin{multline*}
H^0(P_{\pi_o,l,i}(\Pi_l)\otimes \LC_{\r_\oo})[\Pi^{\oo,o}] = \\
\sum_{k=i}^{s_g-l-1} H^0(A_{\pi_o,l,k}(\Pi_l) \otimes \LC_{\rho_\oo})[\Pi^{\oo,o}] +
H^0(P_{\pi_o,l,s_g-l}(\Pi_l) \otimes \LC_{\rho_\oo})[\Pi^{\oo,o}]
\end{multline*}
Le résultat est clairement vrai pour $i=s-l-1$ car $P_{\pi_o,l,s-l-1}(\Pi_l)$ est un
faisceau pervers
ponctuel. Supposons donc le résultat acquis jusqu'au rang $i+1$ et montrons le au rang
$i$. La suite exacte
longue de cohomologie associée à
$$0 \to P_{\pi_o,l,i}(\Pi_l) \longto P_{\pi_o,l,i-1}(\Pi_l) \longto
\PC_-(g,l+i,\pi_o,l,\Pi_l) \oplus
A_{\pi_o,l,i-1}(\Pi_l) \to 0$$ fournit, pour $j \neq 0$, les isomorphismes
$$H^j(P_{\pi_o,l,i}(\Pi_l)\otimes \LC_{\r_\oo})[\Pi^{\oo,o}]
\simeq H^j(P_{\pi_o,l,i-1}(\Pi_l)\otimes \LC_{\r_\oo})[\Pi^{\oo,o}]$$ car d'après
l'hypothèse de récurrence
portant sur les $l$,
$$H^j(\PC(g,l+i,\pi_o,l,\Pi_l)\otimes \LC_{\r_\oo})[\Pi^{\oo,o}]$$
est nul ce qui implique la nullité de
$$H^j(\PC_-(g,l+i,\pi_o,l,\Pi_l)\otimes \LC_{\r_\oo})[\Pi^{\oo,o}]$$
Par ailleurs pour $j=0$, on a la suite exacte courte
\begin{multline*}
0 \to H^0(P_{\pi_o,l,i}(\Pi_l) \otimes \LC_{\r_\oo})[\Pi^{\oo,o}] \longto
H^0(P_{\pi_o,l,i-1}(\Pi_l) \otimes \LC_{\rho_\oo}) \\
\longto H^0(A_{\pi_o,l,i-1}(\Pi_l) \otimes \LC_{\rho_\oo}) \to 0
\end{multline*}
d'où le résultat.

On considère alors la suite exacte longue de cohomologie associée à
$$0 \to P_{\pi_o,l,0}(\Pi_l) \longto j^{\geq lg}_! HT(g,l,\pi_o,\Pi_l)[d-lg] \longto
j^{\geq lg}_{!*} HT(g,l,\pi_o,\Pi_l)[d-lg] \to 0$$ qui s'écrit
$$0 \to H^{-1}(j^{\geq lg}_!) \longto H^{-1}(j^{\geq lg}_{!*}) \longto
H^0(P_{\pi_o,l,s-l-1})
\longto H^0(j^{\geq lg}_!) \longto H^0(j^{\geq lg}_{!*}) \to 0$$ et pour tout $i \neq
-1,0$, $H^i(j^{\geq
lg}_!) \simeq H^i(j^{\geq lg}_{!*})$ où pour alléger les notations, on a omis d'écrire
$HT_{\rho_\oo}(g,l,\pi_o,\Pi_l)[d-lg]$ ainsi que $[\Pi^{\oo,o}]$. On en déduit alors que
pour $i\neq 0$,
$H^i(j^{\geq lg}_!)$ est pur de poids $d-1+i$ alors que $H^0(j^{\geq lg}_!)$ est mixte
de poids inférieur ou
égal à $d-1$. Le calcul de la somme alternée $\sum_i (-1)^i H^i(j^{\geq lg}_!)$,
laquelle est nulle pour $g
\nmid d$ et pour $g | d$ est constituée d'un seul terme de poids $d-1-(s-l)$, implique
alors la nullité des
$H^i(j^{\geq lg}_!)$ pour $i >0$ et celle des $H^i(j^{\geq lg}_{!*}
HT_{\rho_\oo}(g,l,\pi_o,\Pi_l)[d-lg])[\Pi^{\oo,o}]$ pour $i \geq 0$ et pour tout
$\pi_o$. La dualité de
Poincaré donne alors la nullité des $H^i(j^{\geq lg}_{!*}
HT_{\rho_\oo^\vee}(g,l,\pi_o^\vee,\Pi_l^\vee)[d-lg])[\Pi^{\oo,o}]$ pour tout $\pi_o$ et
tout $i \leq 0$ et
donc finalement la nullité des $H^i(j^{\geq
lg}_{!*}HT_{\rho_\oo}(g,l,\pi_o,\Pi_l)[d-lg])[\Pi^{\oo,o}]$ pour
tout $i$ et tout $\pi_o$.

\end{proof}

\rem Il est possible de faire des calculs strictement similaires pour $\Pi_o$
quelconque. Par exemple, au
lemme (\ref{prop-mono}) on traite le cas de $\Pi_o \simeq \st_{n_1}(\xi_1) \boxplus
\cdots \boxplus
\st_{n_r}(\xi_r)$, le cas général étant traité dans la preuve du théorème
(\ref{theo-mono-glob}).

\medskip

\marque D'après le corollaire (\ref{coro-sec-fp}) et en remarquant qu'un faisceau
pervers ponctuel n'a de la
cohomologie qu'en degré zéro, les corollaires suivants découlent directement de
(\ref{somme-alternee}).

\begin{coro} \label{coro-hij-nul}
Pour $g$ divisant $d=sg$ (resp. $1 \leq g < d$ ne divisant pas $d$), soit $\Pi$ une
représentation
irréductible de $D_{\Am}^\times$ vérifiant $\hyp(\r_\oo)$ telle que $\Pi_o \simeq
\st_s(\pi_o)$. On a alors:

(i) pour tout $i \neq d-lg$ (resp. pour tout $i$),
$H^i_c(M_{I,s_o,1}^{=lg},\FC(g,l,\pi_o)_1 \otimes
\LC_{\r_\oo})[\Pi^{\oo,o}]$ et donc $H^i_c(M_{I,s_o}^{=lg}, \FC(g,l,\pi_o,I)\otimes
\LC_{\r_\oo})[\Pi^{\oo,o}]$ sont nuls et pour $i=d-lg$
\begin{multline*}
\lim_{\genfrac..{0pt}{1}{\to}{I}} H^{d-lg}_c(M_{I,s_o,1}^{=lg},\FC(g,l,\pi_o,I)_1
\otimes \LC_{\rho_\oo})[\Pi^{\oo,o}]= \\
 m(\Pi) [\overleftarrow{s-l-1}]_{\pi_o(l(g-1)/2)} \otimes
(\Xi^{\frac{(s-l)(g-1)}{2}} \bigoplus_{\xi \in \AF(\pi_o)} \xi^{-1})
\end{multline*}
\begin{multline*} \lim_{\genfrac..{0pt}{1}{\longto}{I}} H^{d-lg}_c(M_{I,\bar
s_o}^{=lg},HT_{\r_\oo}(g,l,\pi_o,\Pi_l,I))
[\Pi^{\oo,o}] \\ = m(\Pi) (\Pi_l \overrightarrow{\times}
[\overleftarrow{s-l-1}]_{\pi_o}) \otimes
(\Xi^{(s-l)(g-1)/2} \bigoplus_{\xi \in \AF(\pi_o)} \xi^{-1})
\end{multline*}
en tant que représentation de $GL_{(s-l)g}(F_o) \times \Zm$ et de $GL_d(F_o) \times
\Zm$, où
$\AF(\pi_o)$ est l'ensemble des caractères $\xi:\Zm \longto \Qm_l^\times$, tels que
$\pi_o \otimes \xi^{-1}
\circ \val(\det) \simeq \pi_o$ et $m(\Pi)$ est la multiplicité de $\Pi$ dans l'espace
des formes automorphes.

(ii) pour tout $1 \leq l \leq s$, on a
\begin{multline*}
H^0(j^{\geq lg}_! HT(g,l,\pi_o,\Pi_l)[d-lg] \otimes L_g(\pi_o)) \\ =m(\Pi)
e_{\pi_o}(\Pi_l
\overrightarrow{\times} [\overleftarrow{s-l-1}]_{\pi_o}) \otimes L_g(\pi_o)
|-|^{-(s-l)(g-1)/2}
\end{multline*}

\end{coro}

\begin{rema} \label{rema-hyp1}
Si $\Pi$ vérifie $\hyp(\oo)$, la condition $\Pi_o \simeq [\overleftarrow{s-1}]_{\pi_o}$
pour $\pi_o$ une
représentation irréductible cuspidale de $GL_g(F_\oo)$, implique qu'il existe une
représentation irréductible
cuspidale $\pi_\oo$ de $GL_{g'}(F_\oo)$ avec $d=s'g'$ telle que $\Pi_\oo \simeq
[\overleftarrow{s'-1}]_{\pi_\oo}$ et donc $\Pi$ vérifie $\hyp(\r_\oo)$ avec
$\r_\oo=\JL^{-1}([\overleftarrow{s'-1}]_{\pi_\oo})$.
\end{rema}

\marque On rappelle que $\bar D$ est une algèbre à division centrale sur $F$ dont les
invariants locaux en
les places $x$ de $F$ distinctes de $\oo$ et $o$ sont égaux à ceux de $D$ et $\bar D_o
\simeq D_{o,d}$ (resp.
$\bar D_\oo$) est l'algèbre à division centrale sur $F_o$ (resp. $F_\oo$) d'invariant
$1/d$ (resp. $-1/d$).
Pour la preuve de l'énoncé suivant on considère une représentation automorphe
irréductible $\Pi$ de
$D_\Am^\times$ telle que:
\begin{itemize}
\item $\Pi_\oo \simeq [\overleftarrow{s-1}]_{\pi_\oo}$ pour $\pi_\oo$ une représentation
cuspidale de
$GL_d(F_\oo)$ et $\Pi_o \simeq [\overleftarrow{s-1}]_{\pi_o}$;

\item $m(\Pi)=1$;

\item l'ensemble $\AF_{\bar D^\times}(\Pi)$ des représentations irréductibles
automorphes $\bar \Pi$ de $\bar
D_\Am^\times$ telle que $\bar \Pi^{\oo,o} \simeq \Pi^{\oo,o}$ est réduit à un élément
avec $m(\bar \Pi)=1$,
$\bar \Pi_o \simeq \JL^{-1}(\Pi_o)$ et $\bar \Pi_\oo \simeq
\JL^{-1}([\overleftarrow{s-1}]_{\pi_\oo})$.
\end{itemize}

\noindent L'existence d'une telle représentation $\Pi$ est assurée par Henniart (cf.
\cite{he2} annexe A-4).

\begin{coro} \label{coro-pnul}
Pour tout $1 \leq l \leq s_g$, les faisceaux pervers $A_{\pi_o,l,i}(\Pi_l)$ du
corollaire (\ref{coro-sec-fp})
sont tous nuls et pour $g$ ne divisant pas $d$ (resp. $g|d=sg$)
$P_{\pi_o,l,s_g-l}(\Pi_l)$ est nul (resp.
$P_{\pi_o,l,s-l-1}(\Pi_l)$ est le faisceau pervers ponctuel de support l'ensemble des
points supersinguliers
tel que $P_{\pi_o,l,s-l-1}(\Pi_l)$ est isomorphe à
$\FC(g,s,\pi_o)(-\frac{(s-l)(g-1)}{2}) \otimes (\Pi_l
\overrightarrow{\times} [\overleftarrow{s-l-1}]_{\pi_o})$).
\end{coro}

\begin{proof} D'après la proposition précédente on a
\begin{multline*}
\sum_{i=0}^{s_g-l-1} \lim_{\genfrac..{0pt}{1}{\to}{I}}
H^0(M_{I,s_o},A_{\pi_o,l,i}(\Pi_l) \otimes
\LC_{\r_\oo})[\Pi^{\oo,o}]
+ \\
 \lim_{\genfrac..{0pt}{1}{\to}{I}} H^0(M_{I,s_o},P_{\pi_o,l,s_g-l}(\Pi_l)\otimes
\LC_{\r_\oo})[\Pi^{\oo,o}] = \\
\lim_{\genfrac..{0pt}{1}{\to}{I}} H^{d-lg}(M_{I,s_o},j^{\geq lg}_!
HT_{\rho_\oo}(g,l,\pi_o,\Pi_l,I))[\Pi^{\oo,o}].
\end{multline*}
Le membre de droite est d'après (\ref{somme-alternee}) soit nul, pour $g$ ne divisant
pas $d$, soit de poids
$d-lg-s_g+l$ pour $g$ divisant $d$.

- Les faisceaux pervers $A_{\pi_o,l,i}(\Pi_l)$, pour $0 \leq i < s_g-l$ étant de
dimension zéro et purs de
poids $d-lg-i$, la nullité de $H^0(M_{I,s_o,1}^{\geq lg},A_{\pi_o,l,i}(\Pi_l)\otimes
\LC_{\r_\oo})[\Pi^{\oo,o}]$ implique celle des $A_{\pi_o,l,i}(\Pi_l)$.

- De même pour $g$ ne divisant pas $d$, $P_{\pi_o,l,s-l}$, faisceau pervers concentré
aux points
supersinguliers, est nul car son groupe de cohomologie en degré $0$ l'est.

- Pour $g$ divisant $d$, on rappelle, cf. le corollaire (\ref{coro-hij-nul}), que
$$H^i_c(M_{I,s_o,1}^{=lg},\FC(g,l,\pi_o,\Pi_l,I)_1 \otimes \LC_{\r_\oo})[\Pi^{\oo,o}]$$
est nul pour $i \neq d-lg$, et
\begin{multline*}
\lim_{\genfrac..{0pt}{1}{\to}{I}} H^{d-lg}_c(M_{I,s_o,1}^{=lg},\FC(g,l,\pi_o)_1 \otimes
\LC_{\rho_\oo}))[\Pi^{\oo,o}]= \\
 m(\Pi) \bigoplus_{\xi \in \AF(\pi_o)} [\overleftarrow{s-l-1}]_{\pi_o(l(g-1)/2)}
\otimes
\Xi^{\frac{(s-l)(g-1)}{2}} \xi
\end{multline*}
où $\AF(\pi_o)$ est l'ensemble des caractères $\xi:\Zm \longto \Qm_l^\times$, tels que
$\pi_o \otimes
\xi^{-1} \circ \val(\det) \simeq \pi_o$. Ainsi pour toute représentation elliptique
$\Pi_l$ de type $\pi_o$
de $GL_{lg}(F_o)$,
\begin{multline*}
\lim_{\genfrac..{0pt}{1}{\to}{I}} H^0(M_{I,s_o}^{\geq lg},j^{\geq lg}_!
HT_{\rho_\oo}(g,l,\pi_o,\Pi_l,I)[d-lg]) [\Pi^{\oo,o}]=
\\ m(\Pi) e_{\pi_o} \Pi_l \overrightarrow{\times} [\overleftarrow{s-l-1}]_{\pi_o}
\otimes \Xi^{\frac{(s-l)(g-1)}{2}}
\end{multline*}
en tant que représentation de $GL_d(F_o) \times \Zm$. On étudie alors comme précédemment
l'égalité
(\ref{egalite-groth}) et tout particulièrement les faisceaux pervers de dimension nulle.
En ce qui concerne
les faisceaux pervers simples de poids $s(g+1)-4$, le membre de droite de
(\ref{egalite-groth}) est égal à
\begin{multline*}
(P_{\pi_o,s-1,1}([\overleftarrow{s-2}]_{\pi_o})(-\frac{(s-1)(g+1)-2}{2}) - \\
\FC(g,s,\pi_o) \otimes [\overleftarrow{s-2},\overrightarrow{1}]_{\pi_o}) \otimes
L_g(\pi_o)
(-\frac{s(g+1)-4}{2})
\end{multline*}
de sorte que $P_{\pi_o,s-1,1}([\overleftarrow{s-2}]_{\pi_o})$ contient
$\FC(g,s,\pi_o)(-\frac{g-1}{2})
\otimes [\overleftarrow{s-2},\overrightarrow{1}]_{\pi_o}$ et donc, vu que les strates
sont induites,
$\FC(g,s,\pi_o) (-\frac{g-1}{2}) \otimes [\overleftarrow{s-2}]_{\pi_o}
\overrightarrow{\times}
[\overleftarrow{0}]_{\pi_o}$. Par ailleurs d'après (\ref{h0-ss}), on a le lemme suivant

\begin{lemm} \label{lem-pts-ss}
Pour toute représentation $\Pi_s$ de $GL_d(F_o)$, on a
$$\lim_{\genfrac..{0pt}{1}{\to}{I}} H^0(M_{I,s_o}^d,\FC(g,s,\pi_o,I) \otimes
\LC_{\r_\oo} \otimes \Pi_{s})=
e_{\pi_o} \CC^{\oo}_{\bar D,\r_\oo} [\JL^{-1}(\overleftarrow{(s)}_{\pi_o})] \otimes
\Pi_{s}$$ de sorte que sa
partie $\bar \Pi^{\oo,o}$-isotypique est $e_{\pi_o} m(\bar \Pi) \Pi_{s}$.
\end{lemm}

En appliquant ce lemme à $\Pi_s=\Pi_{s-1} \overrightarrow{\times}
[\overleftarrow{0}]_{\pi_o}$, et en
utilisant les propriétés imposées à $\Pi$, on obtient que la partie $\bar
\Pi^{\oo,o}$-isotypique de ${\DS
\lim_{\genfrac..{0pt}{1}{\to}{I}}} H^0(M_{I,s_o}^d,\FC(g,s,\pi_o) \otimes \Pi_{s-1}
\overrightarrow{\times}
[\overleftarrow{0}]_{\pi_o})$ est égale à la partie $\Pi^{\oo,o}$-isotypique de ${\DS
\lim_{\genfrac..{0pt}{1}{\to}{I}}}
H^0_c(M_{I,s_o}^{=(s-1)g},HT(g,s-1,\pi_o,\Pi_{s-1})[g]$ qui d'après ce qui
précède est égale à celle de ${\DS \lim_{\genfrac..{0pt}{1}{\to}{I}}}
H^0(M_{I,s_o},P_{\pi_o,s-1,1}(\Pi_{s-1}))$. On en déduit donc que
$P_{\pi_o,s-1,1}(\Pi_{s-1})$ est égal à
$\FC(g,s,\pi_o)(-\frac{g-1}{2}) \otimes \Pi_{s-1} \overrightarrow{\times}
[\overleftarrow{0}]_{\pi_o}$.

On raisonne alors par récurrence sur $l$ de $s-1$ à $1$, en supposant que pour tout $l >
l_0$,
$P_{\pi_o,l,s-l-1}(\Pi_l)$ est isomorphe à $\FC(g,s,\pi_o)(-\frac{(s-l)(g-1)}{2})
\otimes \Pi_l
\overrightarrow{\times} [\overleftarrow{s-l-1}]_{\pi_o}$. On regarde alors les faisceaux
pervers de dimension
nulle et de poids $s(g+1)-2(s-l_0+1)$ dans l'égalité (\ref{egalite-groth}), soit
\begin{multline*}
P_{\pi_o,l_0,s-l_0-1} ( [\overleftarrow{l_0-1}]_{\pi_o}) \otimes
(-\frac{l_0(g+1)-2}{2})+ \FC(g,s,\pi_o) \otimes \\
(\sum_{i=1}^{s-l_0} (-1)^i [\overleftarrow{l_0-1},\overrightarrow{i}]_{\pi_o}
\overrightarrow{\times}
[\overleftarrow{s-l_0-i-2}]_{\pi_o}) \otimes L_g(\pi_o) (-\frac{s(g+1)-2(s-l_0+1)}{2})
\end{multline*}
ce qui donne
$$P_{\pi_o,l_0,s-l_0-1}([\overleftarrow{l_0-1}]_{\pi_o}) - \FC(g,s,\pi_o) \otimes
[\overleftarrow{l_0-1},\overrightarrow{1},\overleftarrow{s-l_0-1}]_{\pi_o} \otimes
L_g(\pi_o)
(-\frac{(s-l_0)(g-1)}{2})$$ positif. On en déduit alors que
$P_{\pi_o,l_0,s-l_0-1}([\overleftarrow{l_0-1}]_{\pi_o})$ contient
$$\FC(g,s,\pi_o) (-\frac{(s-l_0)(g-1)}{2})
\otimes [\overleftarrow{l_0-1},\overrightarrow{1},\overleftarrow{s-l_0-1}]_{\pi_o},$$ et
donc vu le caractère
induit des strates, $P_{\pi_o,l_0,s-l_0-1}(\Pi_{l_0})$ contient
$$\FC(g,s,\pi_o)(-\frac{(s-l_0)(g-1)}{2})
\otimes \Pi_{l_0} \overrightarrow{\times} [\overleftarrow{s-l_0-1}]_{\pi_o}.$$ L'égalité
provient alors de
l'égalité des $H^0$ comme dans le cas $l_0=s-1$, traité ci-avant, d'où le résultat.

\end{proof}

\begin{coro} \label{global0-ok}
Le théorème (\ref{theo-global0}) est vrai.
\end{coro}

\begin{proof} En effet le résultat découle directement de la proposition
(\ref{prop-libre}) en réinjectant les
égalités de la proposition (\ref{prop-p}) où les $P_{!,l}$ sont nuls d'après le
corollaire (\ref{coro-pnul}).
On obtient ainsi pour tout $1 \leq i \leq s_g$:
\begin{multline*}
\sum_{l=i}^{s_g} (-1)^{l-i} j^{\geq lg}_!
HT(g,l,\pi_o,[\overleftarrow{i-1},\overrightarrow{l-i}]_{\pi_o})[d-lg] \otimes
L_g(\pi_o)
(-\frac{lg-2+2i-l}{2}) = \\ \sum_{l=i}^{s_g} (-1)^{l-i} \sum_{r=0}^{s_g-l} j^{\geq
(l+r)g}_{!*}
HT(g,l+r,\pi_o,[\overleftarrow{i-1},\overrightarrow{l-i}]_{\pi_o}
\overrightarrow{\times}
[\overleftarrow{r-1}]_{\pi_o})[d-(l+r)g] \\ \otimes L_g(\pi_o)
(-\frac{(l+r)(g-1)+2(i-1)}{2})
\end{multline*}
laquelle somme est alors égale à
$$\sum_{l=i}^{s_g} j^{\geq lg}_{!*} HT(g,l,\pi_o,\Pi_{l,i})[d-lg] \otimes L_g(\pi_o)
(-\frac{l(g-1)+2(i-1)}{2})$$ où
$$\Pi_{l,i}:=\sum_{r=i}^l (-1)^{r-i} [\overleftarrow{i-1},\overrightarrow{r-i}]_{\pi_o}
\overrightarrow{\times} [\overleftarrow{l-r-1}]_{\pi_o} \\ =
[\overleftarrow{l-1}]_{\pi_o} $$ avec
$$\begin{array}{rl} e_{\pi_o}[R\Psi_{\eta_o,\pi_o}(\bar \Qm_l)[d-1]] = &
\sum_{i=1}^{s_g} \sum_{l=i}^{s_g}
j^{\geq lg}_{!*}
HT(g,l,\pi_o,\Pi_{l,i})[d-lg] \otimes \\ & \hfill L_g(\pi_o) (-\frac{l(g-1)+2(i-1)}{2})
\\
= & \sum_{l=1}^{s_g} j^{\geq lg}_{!*} HT(g,l,\pi_o,[\overleftarrow{l-1}]_{\pi_o})
\otimes \\ & \hfill
L_g(\pi_o)
(-\frac{lg-1}{2}) (\sum_{\genfrac{}{}{0pt}{}{|k| < l}{k \equiv l-1 \mod 2}}
(-\frac{k}{2}) \\
 = & \sum_{k=1-s}^{s-1} \sum_{\genfrac{}{}{0pt}{}{|k| < l \leq s_g}{l \equiv k-1 \mod
2}} j^{\geq lg}_{!*}
HT(g,l,\pi_o,[\overleftarrow{l-1}]_{\pi_o})[d-lg] \\
& \hfill \otimes L_g(\pi_o) (-\frac{lg+k-1}{2})
\end{array}$$
d'où le résultat.

\end{proof}

\begin{coro} \label{coro-hic}
Pour tout $g|d=sg$, et $\pi_o$ une représentation cuspidale de $GL_g(F_o)$, on considère
les groupes de
cohomologie $H^j_c(M_{I,\bar s_o}^{=h},R^i \Psi_{\eta_o,I,\pi_o}(\LC_{\r_\oo}))$ ainsi
que leur limite
inductive sur tous les idéaux $I$ de $A$, limite que l'on notera
$H^j_{=h,i,\pi_o,\r_\oo}$. Pour $\Pi$ une
représentation automorphe de $D_\Am^\times$ vérifiant $\hyp(\r_\oo)$, on a les résultats
suivant:
\begin{itemize}
    \item[(1)] pour $g$ ne divisant pas $d$, les $H^j_{=h,i,\pi_o,\r_\oo}[\Pi^{\oo,o}]$
sont nuls pour tous $j,h,i$;

    \item[(2)] pour $g$ divisant $d=sg$ et $\Pi_o \simeq [\overleftarrow{s-1}]_{\pi_o}$,
les
    $H^j_{=h,i,\pi_o,\r_\oo}[\Pi^{\oo,o}]$ sont nuls pour $h$ qui n'est pas de la forme
    $lg$ avec $1 \leq l \leq s$;

    \item[(3)] pour $g$ divisant $d=sg$ et $\Pi_o \simeq [\overleftarrow{s-1}]_{\pi_o}$,
les
    $H^j_{=lg,i,\pi_o,\r_\oo}[\Pi^{\oo,o}]$ sont nuls pour $j \neq d-lg$;

    \item[(4)] pour $g$ divisant $d=sg$ et $\Pi_o \simeq [\overleftarrow{s-1}]_{\pi_o}$,
les
    $H^{d-lg}_{=lg,i,\pi_o,\r_\oo}[\Pi^{\oo,o}]$ sont nuls pour $i$ ne vérifiant pas
    $l(g-1) \leq i \leq lg-1$. Si $1 \leq l < s$ et $i=lg-r$ avec $1 \leq r \leq l$,
    $H^{d-lg}_{=lg,lg-r,\pi_o,\r_\oo}[\Pi^{\oo,o}]$ est isomorphe à
    $$m(\Pi)
    ([\overleftarrow{l-r},\overrightarrow{r-1}]_{\pi_o} \overrightarrow{\times}
[\overleftarrow{s-l-1}]_{\pi_o})
    \otimes L_g(\pi_o) (-\frac{s(g-1)-2(r-l)}{2})$$ en tant que représentation de
$GL_d(F_o) \times W_o$, où $m(\Pi)$
    est la multiplicité de $\Pi$ dans l'espace des formes automorphes.

\end{itemize}
\end{coro}

\begin{proof} On rappelle (cf. la proposition (\ref{prop-hic})) que pour tout $i,j$, on
a un isomorphisme canonique
$$H^j_c(M_{I,s_o,1}^{=h}, R^i\Psi_{\eta_o}(\LC_{\r_\oo}))^h \simeq \bigoplus_{\t_o \in
\CF_h}
(H^j_c(M_{I,s_o,1}^{=h},\FC_{\t_o} \otimes \LC_{\r_\oo}) \otimes
\widetilde{\UC_{F_o,m}^{h,i}(\t_o)})^{h/e_{\t_o}}$$ de sorte que les résultats découle
directement du
corollaire (\ref{coro-hij-nul}).

\end{proof}

\section{Faisceaux de cohomologie des $j^{\geq lg}_{!*} \FC(g,l,\pi_o)$}
\label{glob2}

Rappelons, cf \S \ref{schema}, que d'après le théorème de comparaison de
Berkovich-Fargues, la filtration de
monodromie-locale du complexe $\Psi_{F_o}^d$ sera donnée par les germes en un point
supersinguliers des
faisceaux de cohomologie des $gr_k$. On propose alors de calculer tous les faisceaux de
cohomologie des
$gr_k$. D'après la proposition (\ref{prop-p}) précisée par le corollaire
(\ref{coro-pnul}), il nous suffit de
déterminer les faisceaux de cohomologie des $j^{\geq lg}_{!*} \FC(g,l,\pi_o)[d-lg]$.
Nous verrons que cela ne
pose aucun problème en dehors des points supersinguliers car on peut utiliser
l'hypothèse de récurrence dans
le théorème (\ref{theo-local-fil}). Au niveau des points supersinguliers on ne dispose
d'aucun renseignement
local, cependant la proposition (\ref{prop-p}) nous permet de restreindre efficacement
les possibilités pour
les germes en un point supersingulier de ces faisceaux de cohomologie.

\subsection{Une écriture dans $\GF$ de $j^{\geq lg}_{!*} \FC(g,l,\pi_o)[d-lg]$}

\begin{lemm} \label{lem-hij0}
Pour tout $1 \leq l \leq s_g$, on a l'égalité dans $\GF$:
\begin{multline*}
\PC(g,l,\pi_o)= \sum_{r=0}^{s_g-l} (-1)^r j^{\geq (l+r)g}_!
HT(g,l+r,\pi_o,[\overleftarrow{l-1}]_{\pi_o}
\overrightarrow{\times} [\overrightarrow{r-1}]_{\pi_o})[d-(l+r)g] \\
\otimes L_g(\pi_o) (-\frac{r(g-1)}{2})
\end{multline*}
\end{lemm}

\begin{proof} On démontre le résultat par récurrence sur $l$ de $s_g$ à $1$, en
utilisant les corollaires
(\ref{coro-sec-fp}) et (\ref{coro-pnul}). Le cas $l=s_g$ est directement donné par loc.
cit. Supposons le
résultat acquis jusqu'au rang $l+1$ et traitons le cas de $l$. D'après loc. cit. on a
\begin{multline} \label{egalite-hij}
j^{\geq lg}_! HT(g,l,\pi_o,[\overleftarrow{l-1}]_{\pi_o})[d-lg] \otimes L_g(\pi_o) =
\PC(g,l,\pi_o)
\\ + \sum_{i=1}^{s-l} j^{\geq (l+i)g}_{!*} HT(g,l+i,\pi_o,[\overleftarrow{l-1}]_{\pi_o}
\overrightarrow{\times} [\overleftarrow{i-1}]_{\pi_o})[d-(l+i)g] \otimes L_g(\pi_o)
(-\frac{ig-i}{2})
\end{multline}

D'après l'hypothèse de récurrence on a
\begin{multline*}
j^{\geq (l+i)g}_{!*} HT(g,l+i,\pi_o,[\overleftarrow{l-1}]_{\pi_o}
\overrightarrow{\times}
[\overleftarrow{i-1}]_{\pi_o})[d-(l+i)g]
\otimes L_g(\pi_o) (-\frac{ig-i}{2}) = \\
\sum_{r=0}^{s-l-i} (-1)^r j^{\geq (l+i+r)g}_!
HT(g,l+i+r,\pi_o,[\overleftarrow{l-1}]_{\pi_o}
\overrightarrow{\times} [\overleftarrow{i-1}]_{\pi_o} \overrightarrow{\times}
[\overrightarrow{r-1}]_{\pi_o})[d-(l+i+r)g] \\ \otimes L_g(\pi_o)
(-\frac{(i+r)g-i-r}{2})
\end{multline*}
de sorte que (\ref{egalite-hij}) s'écrit
\begin{multline*}
\PC(g,l,\pi_o)-j^{\geq lg}_! HT(g,l,\pi_o,[\overleftarrow{l-1}]_{\pi_o})[d-lg] \otimes
L_g(\pi_o)= \\
\sum_{r=1}^{s-l} (-1)^{r} j^{\geq (l+r)g}_!
HT(g,l+r,\pi_o,[\overleftarrow{l-1}]_{\pi_o}
\overrightarrow{\times} \Pi_r) [d-(l+r)g] \otimes L_g(\pi_o) (-\frac{(rg-r}{2})
\end{multline*}
où $\Pi_r=(-1)^{r-1} [\overleftarrow{r-1}]_{\pi_o} + \sum_{i=1}^{r-1} (-1)^{i-1}
[\overleftarrow{i-1}]_{\pi_o} \overrightarrow{\times}
[\overrightarrow{r-i-1}]_{\pi_o}=[\overrightarrow{r-1}]_{\pi_o}$ d'où le résultat.

\end{proof}

\subsection{Filtration de monodromie-locale en hauteur non maximale}
\label{bf2}

On rappelle que d'après l'hypothèse de récurrence, le théorème (\ref{theo-local-fil})
est connu pour les
modèles locaux de Deligne-Carayol de hauteur strictement inférieure à $d$. Ainsi
l'aboutissement de la suite
spectrale (\ref{suite-spectrale}) est connue en dehors des points supersinguliers tandis
qu'on ne connaît que
les germes en des points non supersinguliers des termes $E_1$ (cf. le
théorème-définition
(\ref{theo-defi-monodromie})).

\begin{prop} \label{prop-hij}
Pour $g \neq 1$ et $1 \leq l \leq s_g$, les faisceaux de cohomologie $h^i
\PC(g,l,\pi_o)$ sont nuls pour tout
$i < l-s_g$ et $i$ qui n'est pas de la forme $lg-d+a(g-1)$. Pour $i=lg-d+a(g-1)$ avec $0
\leq a < s_g-l$, ils
sont égaux dans $\FH(M_{s_o})$ à
$$j^{\geq (l+a)g}_! HT(g,l+a,\pi_o, [\overleftarrow{l-1}]_{\pi_o}
\overrightarrow{\times}
[\overrightarrow{a-1}]_{\pi_o}) \otimes L_g(\pi_o) (-\frac{a(g-1)}{2}).$$
\end{prop}

\begin{proof} Pour $i < lg-d$, les $h^i \PC(g,l,\pi_o)$ sont tous nuls car
$\PC(g,l,\pi_o)$ est de dimension
$d-lg$. On a la suite exacte courte de faisceaux pervers
$$0 \to P_{\pi_o,l,0}(\pi_o) \otimes L_g(\pi_o) \longto j^{\geq g}_!
\FC(g,l,\pi_o,\pi_o)[d-lg] \otimes L_g(\pi_o)
\longto \PC(g,l,\pi_o) \to 0$$
où $P_{\pi_o,1,0}(\pi_o)$ est un faisceau pervers de dimension $d-(l+1)g$ de sorte que
$$h^{lg-d} j^{\geq g}_{!*} \FC(g,1,\pi_o)[d-g]=j^{\geq g}_! \FC(g,1,\pi_o).$$
Le principe est alors d'utiliser le théorème de
comparaison de Berkovich-Fargues couplé avec le lemme (\ref{lem-hij0}).

Remarquons tout d'abord que d'après le corollaire (\ref{coro-grk}), pour tout $k$,
$gr_{k,\pi_o}$ est pur
hors des points supersinguliers, de sorte qu'en ce qui concerne les strates non
supersingulières, on a

$$e_{\pi_o} h^i gr_{k,\pi_o}=\bigoplus_{\genfrac{}{}{0pt}{}{|k| < l \leq s_g}{l \equiv
k-1 \mod 2}} h^i \PC(g,l,\pi_o)(-\frac{lg-1+k}{2})$$
Ainsi d'après le corollaire (\ref{coro-comparaison-grk}), on en déduit le lemme
suivant.

\begin{lemm} \label{lem-hij1}
Les $h^i \PC(g,l,\pi_o)$ vérifient les propriétés suivantes:
\begin{itemize}
\item hormis les points supersinguliers, ils sont à support dans les strates
$M_{I,s_o}^{=l'g}$ pour
$l \leq l' \leq s_g$ avec $l'g-d-l'+l-1' \leq i \leq l'g-d$ et $i \equiv l'g-d-l'+l-1
\mod 2$;

\item pour $i=l'g-d-l'+l-1+2r$, la fibre en un point géométrique de $M_{I,s_o}^{=l'g}$
de $h^i
\PC(g,l,\pi_o)$ est un facteur direct de
$$[\overleftarrow{l+2r-1}]_{\pi_o} \overrightarrow{\times}
[\overrightarrow{l'-l-2r-1}]_{\pi_o} \otimes L_g(\pi_o)
(-\frac{l'(g-1)+2(l-1)+2r}{2})$$
\end{itemize}
\end{lemm}

Soit $z$ un point géométrique de $M_{I,s_o}^{=l'g}$. Les strates étant induites, on en
déduit que la fibre en
$z$ de $h^i \PC(g,l,\pi_o)(-\frac{l(g-1)}{2})$ est de la forme:
$$\bigoplus_{\xi} ([\overleftarrow{l-1}]_{\pi_o \circ \xi (-(l'-l)(g-1)/2)}
\overrightarrow{\times} \pi_\xi)
\otimes L_g(\pi_o \circ \xi) (-\frac{l(g-1)}{2})$$ où $\pi_\xi$ est une représentation
de
$GL_{(l'-l)g}(F_o)$. Pour $\xi=\Xi^{((l'-l)(g-1)+2r)/2}$ avec $r>0$, on remarque que si
$[\overleftarrow{l-1}]_{\pi_o(r)} \overrightarrow{\times} \pi_\xi$ contient un des deux
constituants de
$$[\overleftarrow{l+2r-1}]_{\pi_o} \overrightarrow{\times}
[\overrightarrow{l'-l-2r-1}]_{\pi_o}$$
alors il contient aussi tous les constituants de
$$([\overleftarrow{l-1}]_{\pi_o} \overrightarrow{\times}
[\overrightarrow{l'-l-2r-1}]_{\pi_o})
\overleftarrow{\times} [\overleftarrow{2r-1}]_{\pi_o}$$ alors que, par exemple
$[\overleftarrow{2r-1},\overrightarrow{1},\overleftarrow{l-1},\overrightarrow{l'-l-2r}]_{\pi_o}
\otimes
L_g(\pi_o) (-\frac{l'(g-1)+2r}{2})$ n'est pas un constituant de $e_{\pi_o}
h^{l'g-d-l'+l-1+2r}
gr_{1-l,\pi_o}$ d'où la contradiction.

On en déduit ainsi que $h^{l'g-d-l'+l-1} \PC(g,l,\pi_o)$ est le seul faisceau de
cohomologie ayant un support
d'intersection non vide avec la strate $l'g$. Le raisonnement étant valide pour tout
$l'$, on en déduit aussi
que le support de $h^{l'g-d-l'+l-1} \PC(g,l,\pi_o)$ est contenu dans la strate $l'g$. Le
lemme
(\ref{lem-hij0}) donne alors que la restriction à la strate $l'g$ de $h^{l'g-d-l'+l-1}
\PC(g,l,\pi_o)$ est
isomorphe à $HT(g,l',\pi_o, [\overleftarrow{l-1}]_{\pi_o} \overrightarrow{\times}
[\overrightarrow{l'-l-1}]_{\pi_o}) \otimes L_g(\pi_o) (-\frac{(l'-l)(g-1)}{2})$, de
sorte que
\begin{multline*}
h^{l'g-d-l'+l-1} \PC(g,l,\pi_o) \simeq j^{\geq l'g}_! HT(g,l',\pi_o,
[\overleftarrow{l-1}]_{\pi_o}
\overrightarrow{\times} [\overrightarrow{l'-l-1}]_{\pi_o}) \\
\otimes L_g(\pi_o) (-\frac{(l'-l)(g-1)}{2}).
\end{multline*}

\end{proof}

\begin{prop} \label{prop-hij1}
Pour $g=1$, les $h^i \PC(1,l,\xi_o)$, sont nuls pour $i \neq l-d$ et pour $i=l-d$ leur
restriction à
$M_{I,s_o}^{=l+r}$, pour $0 \leq r < d-l$ est égale dans $\FH(M_{s_o}^{=l+r})$ à
$HT(1,l+r,\xi_o,[\overleftarrow{l-1}]_{\xi_o} \overrightarrow{\times}
[\overrightarrow{r-1}]_{\xi_o})$.
\end{prop}

\begin{proof} Elle est strictement identique à celle pour $g \neq 1$, en considérant à
chaque étape les restrictions
aux strates $M_{I,s_o}^{=l+r}$.

\end{proof}

\begin{prop} \label{prop-hij-ss}
Pour $g$ ne divisant pas $d$, $h^{s_g g-d-(s_g-l)} \PC(g,l,\pi_o)$ est, dans
$\FH(M_{s_o})$, isomorphe à
$$j^{\geq s_g g}_! HT(g,s_g,\pi_o,([\overleftarrow{l-1}]_{\pi_o}
\overrightarrow{\times}
[\overrightarrow{s_g-l-1}]_{\pi_o})) \otimes L_g(\pi_o)(-\frac{(s_g-l)(g-1)}{2}).$$ Par
ailleurs les $h^i
\PC(g,l,\pi_o)$ sont nuls pour $s_g g -d -(s_g-l) < i \leq 0$.
\end{prop}

\begin{proof} Tant qu'on est en dehors des points supersinguliers, les arguments
précédents s'adaptent, en utilisant
le théorème de comparaison de Berkovich-Fargues avec la connaissance de la filtration de
monodromie-locale
des modèles de Deligne-Carayol en hauteur strictement inférieure à $d$. En ce qui
concerne les points
supersinguliers, on raisonne par récurrence sur $l$ de $s_g$ à $1$. Pour $l=s_g$, le
corollaire
(\ref{coro-sec-fp}), joint au corollaire (\ref{coro-pnul}), donne
$$\PC(g,s_g,\pi_o)=j^{\geq s_g g}_!
HT(g,s_g,\pi_o,[\overleftarrow{s_g-1}]_{\pi_o})[d-s_gg] \otimes L_g(\pi_o)$$
d'où le résultat. Supposons donc le résultat acquis jusqu'au rang $l+1$ et traitons le
cas de $l$. On
considère les suites exactes courtes du corollaire (\ref{coro-sec-fp}), où l'on rappelle
que d'après
(\ref{coro-pnul}), les $A_{\pi_o,l,g}$ sont nuls:
\begin{multline*}
0 \to P_{\pi_o,l,0}([\overleftarrow{l-1}]_{\pi_o}) \longto j^{\geq lg}_! HT(g,l,\pi_o,
[\overleftarrow{l-1}]_{\pi_o})[d-lg] \longto \\ j^{\geq lg}_{!*}
HT(g,l,\pi_o,[\overleftarrow{l-1}]_{\pi_o})[d-lg] \to 0
\end{multline*}
$$\cdots$$
$$0 \to P_{\pi_o,l,s_g-l}([\overleftarrow{l-1}]_{\pi_o}) \longto
P_{\pi_o,l,s_g-l-1}([\overleftarrow{l-1}]_{\pi_o}) \longto
\PC_-(g,s_g,\pi_o,l,[\overleftarrow{l-1}]_{\pi_o})
\to 0$$ où $P_{\pi_o,l,s_g-l}([\overleftarrow{l-1}]_{\pi_o})$ est le faisceau pervers
nul. On démontre alors
par récurrence sur $r$, de $s_g-l$ à $0$, que le germe en un point supersingulier de
$h^i
P_{\pi_o,l,r}([\overleftarrow{l-1}]_{\pi_o})$ est nul pour tout $i$, d'où le résultat.

\end{proof}

En direction du théorème (\ref{theo-global2}), on a le résultat suivant.

\begin{lemm} \label{lem-hij3}
Pour $g \neq 1$ divisant $d=sg$ (resp. $g=1$), la fibre en un point supersingulier de
$h^i \PC(g,l,\pi_o)$
(resp. de $h^i \PC(1,l,\xi_o)$) est nulle pour $i<l-s$. Par ailleurs les $h^{l-s+i}
\PC(g,l,\pi_o)$ sont
concentrés aux points supersinguliers pour $0 \leq i < s-l$, de fibre (resp. la fibre en
un point
supersingulier de $h^{l-s+i} \PC(1,l,\xi_o)$ est), en tant que $GL_d(F_o) \times
W_o$-module, un sous-espace,
éventuellement nul, de
$$[\overleftarrow{l-1}]_{\pi_o} \overrightarrow{\times}
[\overrightarrow{s-l-2-i}]_{\pi_o}
\overrightarrow{\times} [\overleftarrow{i}]_{\pi_o} \otimes L_g(\pi_o)
(-\frac{(s-l)(g-1)}{2}).$$
\end{lemm}

\begin{proof} Nous ne traitons que le cas $g \neq 1$, le raisonnement pour $g=1$ étant
strictement identique. On
raisonne par récurrence sur $l$ de $s$ à $1$, le cas $l=s$ étant trivial. Supposons donc
le résultat acquis
jusqu'au rang $l+1$ et traitons le cas de $l$. On considère alors les suites exactes
courtes de faisceaux
pervers de (\ref{coro-sec-fp}). On remarque tout d'abord que toutes les fibres en un
point supersingulier des
faisceaux de cohomologie des $\PC_-(g,l+r,\pi_o,l,\Pi_l)$ sont de poids $(s-l)(g-1)$ de
sorte qu'il en est de
même pour les $P_{\pi_o,l,r}(\Pi_l)$ et donc pour $\PC(g,l,\pi_o)$. Le résultat découle
alors de l'étude des
suites exactes longues associées. Une façon plus visuelle et plus directe d'obtenir le
résultat est de
considérer l'isomorphisme $Dj^{\geq lg}_{!*} \FC(g,l,\pi_o)[d-lg] \simeq j^{\geq
lg}_{!*}
\FC(g,l,\pi_o^\vee)[d-lg](d-lg)$ et de regarder la suite spectrale
\begin{equation} \label{ss-dualite}
E_2^{p,q}=R^p Hom(h^{-q} \PC(g,l,\pi_o),K_{s_o}) \Rightarrow h^{p+q}
\PC(g,l,\pi_o^\vee)(d-1)
\end{equation}
où $K_{s_o}$ désigne le complexe dualisant sur $M_{I,s_o}$ \footnote{On pourra voir les termes $E_2^{p,q}$ à
la figure (\ref{figure1}) (resp. (\ref{figure2})), pour $g=7$, $s=5$ et $l=4$ (resp. $l=3$).}

On rappelle que pour un faisceau $\LC$ sur $M_{I,s_o}^{=(l+r)g}$, par adjonction on a
$$RHom(i^{(l+r)g}_* j^{\geq (l+r)g}_! \LC,f^! \bar \Qm_l) \simeq  i^{(l+r)g}_* Rj^{\geq
(l+r)g}_*
RHom(\LC,j^{\geq (l+r)g,!} i^{(l+r)g,!} f^! \bar \Qm_l)$$ et comme $M_{I,s_o}^{=(l+r)g}$
est lisse, on a
$$j^{\geq (l+r)g,!} i^{(l+r)g,!} f^! \bar \Qm_l \simeq \bar
\Qm_l[2(d-(l+r)g)](d-(l+r)g)$$
soit pour $p \geq -2(d-(l+r)g)$,
\begin{multline} \label{etoile}
R^p Hom (j^{\geq (l+r)g}_! HT(g,l+a,\pi_o, \\ [\overleftarrow{l-1}]_{\pi_o}
\overrightarrow{\times}
[\overleftarrow{a-1}]_{\pi_o}) \otimes |\cl|^{-a(g-1)/2}, K_{s_o}) \\ \simeq
(i^{(l+a)g}_* R^{p+2(d-(l+a)g)}
j^{\geq (l+a)g}_* HT(g,l+a,\pi_o^\vee, \\ [\overleftarrow{l-1}]_{\pi_o^\vee}
\overleftarrow{\times}
[\overleftarrow{a-1}]_{\pi_o^\vee}) \otimes |\cl|^{-a(g+1)/2})(lg-d)
\end{multline}

Ainsi pour $z$ un point supersingulier, $(E_2^{p,q})_z$, pour $q$ de la forme
$d-lg-a(g-1)$ (resp. $q=a+1$)
avec $0 \leq a < s-l$, et $p=-(s-l-a)g-i$ avec $0 \leq i \leq s-l-a$ (resp. $p=0$), s'il
est non nul, est
mixte de poids $(s-l-a-r)(g-1)/2+(a+r)(g+1)/2$ avec $0 \leq r <s-l-a$ si $p<-(s-l-a)g$
et pur de poids
$(s-l)(g+1)/2$ si $p=(s-l-a)g$. Dans ce dernier cas on obtient alors
$$[\overleftarrow{l-1}]_{\pi_o} \overrightarrow{\times}
[\overrightarrow{s-l-a-1}]_{\pi_o} \overleftarrow{\times}
[\overleftarrow{a-1}]_{\pi_o} \otimes |\cl|^{-(s-l)(g+1)/2}$$ (resp. $(h^{-q} j^{\geq
lg}_{!*}
\FC(g,l,\pi_o)_0)_z^\vee (-lg) \otimes |\cl|^{-(s-l)(g+1)/2}$ où le dual est pris en
tant que représentation
de $GL_{(s-l)g}(F_o)$). Dans la figure (\ref{fig-dualite2}) on représente ces
$(E_2^{p,q})_z$ de poids
$(s-l)(g+1)$. Le résultat découle alors du fait que les $E_\oo^i$ sont tous nuls pour $i
\geq 0$.

\input{figg1.tex}

\end{proof}

\section{Preuve des théorèmes globaux sous IV.7.1.1}
\label{globss}

On rappelle, cf. \S \ref{schema}, que le théorème (\ref{theo-local-fil}) découle,
d'après le théorème de comparaison de Berkovich-Fargues,
des théorèmes globaux (\ref{theo-global1}), (\ref{theo-global2}) et (\ref{theo-ss}). En
outre le théorème
(\ref{theo-local-fil}) implique, grâce à (\ref{theo-ss}) via Berkovich-Fargues, le
théorème (\ref{theo-ripsi-local}).
Par ailleurs le principe de la preuve de la proposition (\ref{prop-hij}) montre que
(\ref{theo-local-fil}) implique les théorèmes
globaux. Le but de ce paragraphe est de montrer que le théorème (\ref{theo-ripsi-local})
implique les théorèmes globaux
(\ref{theo-global1}), (\ref{theo-global2}) et (\ref{theo-ss}) et donc le théorème local
(\ref{theo-local-fil}).
Par souci d'efficacité, on montrera, en utilisant l'opérateur $N$, qu'il
nous suffit en fait de connaître les parties de poids $s(g-1)$ de
(\ref{theo-ripsi-local}).

\subsection{Preuve de (\ref{theo-global2}) sous (\ref{cas-m})}
\label{globss1}

\begin{prop} \label{cas-m}
Pour tout diviseur $g$ de $d=sg$ et $\pi_o$ une représentation irréductible cuspidale de
$GL_g(F_o)$, la
partie de poids $d-s$ de
$\widetilde{\UC_{F_o,\xi_o}^{d,d-s+i}}(\JL^{-1}([\overleftarrow{s-1}]_{\pi_o}))$ est
nulle pour $0 < i <s$ et égale à $[\overrightarrow{s-1}]_{\pi_o} \otimes L_g(\pi_o)
(-\frac{d-s}{2})$ pour
$i=0$.
\end{prop}

La preuve de cette proposition sera donnée au paragraphe suivant; montrons comment le
théorème
(\ref{theo-local-fil}) en découle. On commence par prouver le théorème
(\ref{theo-global2}) qui d'après la
proposition (\ref{prop-hij}) découle alors de la proposition suivante.

\begin{prop} \label{prop-hij-ss-g}
Sous le résultat de la proposition (\ref{cas-m}), pour $g \neq 1$ divisant $d$,
$h^{l-s+i} \PC(g,l,\pi_o)$
est nul pour $i \neq 0$ et sinon est concentré aux points supersinguliers de fibre
isomorphe à
$[\overleftarrow{l-1}]_{\pi_o} \overrightarrow{\times} [\overrightarrow{s-l-1}]_{\pi_o}
\otimes L_g(\pi_o)
(-\frac{(s-l)(g-1)}{2})$ pour $i=0$.
\end{prop}

\begin{proof} D'après le lemme (\ref{lem-hij3}) tous les germes en un point
supersingulier des $\PC(g,l,\pi_o)$,
sont de poids $(s-l)(g-1)$. Le raisonnement de la preuve de la proposition
(\ref{prop-hij}) s'applique alors
tel quel en étudiant les $\PC(g,l,\pi_o)(-\frac{l(g-1)}{2})$ pour $1 \leq l \leq s$ et
en utilisant la
connaissance des parties de poids $s(g-1)$ de (\ref{theo-local-fil}).

\end{proof}

\subsection{Preuve de (\ref{theo-ss}) sous (\ref{cas-m})}
\label{globss2}

On considère la suite spectrale
(\ref{suite-spectrale}) associée à la filtration de monodromie. D'après la proposition
précédente et le point
(ii) de la proposition (\ref{prop-p}), le germe en un point supersingulier $z$ de
$E_1^{i,j}$ vérifie les
propriétés suivantes (cf. la figure (\ref{figure7}) dans le cas $s=4$ et $g=3$):
\begin{itemize}
\item il est nul pour $i,j$ ne vérifiant pas la condition suivante: il existe $0 \leq k
\leq s-1$ tel
que $j=1-s+2k$ et $-k \leq i \leq s-1-2k$;

\item pour $j=1-s+2k$ avec $0 \leq k \leq s-1$, le germe en $z$ de
$E_1^{s-1-2k-r,1-s+2k}$ pour $0 \leq r \leq s-1-k$
est isomorphe à
$$[\overleftarrow{s-r-1}]_{\pi_o} \overrightarrow{\times} [\overrightarrow{r-1}]_{\pi_o}
\otimes L_g(\pi_o)
(-\frac{s(g-1)+2k}{2}).$$
\end{itemize}

En outre d'après (\ref{cas-m}), pour $j=1-s$ et pour tout $0 \leq i < s-1$,
l'application
$d_1^{i,1-s}:E_1^{i,1-s} \longto E_1^{i+1,1-s}$ induit en $z$ la flèche non triviale
$$[\overleftarrow{i}]_{\pi_o} \overrightarrow{\times} [\overrightarrow{s-i-2}]_{\pi_o}
\longto
[\overleftarrow{i+1}]_{\pi_o} \overrightarrow{\times} [\overrightarrow{s-i-3}]_{\pi_o}$$
dont le noyau est
$[\overleftarrow{i},\overrightarrow{s-i-1}]_{\pi_o}$ et le conoyau est
$[\overleftarrow{i+2},\overrightarrow{s-i-3}]_{\pi_o}$. Le logarithme $N^k$ de la partie
unipotente de la
monodromie induit un diagramme commutatif
$$\xymatrix{ E_1^{i,1-s} \rrto^{d_1^{i,1-s}} & & E_1^{i+1,1-s} \\
E_1^{i-2k,1-s+2k} \uto^{N^k} \rrto^{d_1^{i-2k,1-s+2k}} & & E_{i-2k+1,1-s+2k}
\uto^{N^k}}$$ En remarquant que
pour $0 \leq k \leq -i$, $N^k$ induit des isomorphismes $E_1^{i-2k,1-s+2k} \simeq
E_1^{i,1-s}$ et
$E_1^{i+1-2k,1-s+2k} \simeq E_1^{i+1,16s}$, on en déduit que la flèche
$d_1^{i-2k,1-s+2k}$ est aussi induite
par la flèche non triviale
$$[\overleftarrow{i}]_{\pi_o} \overrightarrow{\times} [\overrightarrow{s-i-2}]_{\pi_o}
\longto
[\overleftarrow{i+1}]_{\pi_o} \overrightarrow{\times}
[\overrightarrow{s-i-1}]_{\pi_o}$$

Finalement pour tout $1 \leq l' < l \leq s$ et tout $k \equiv l'-1 \mod 2$ avec $|k|
\leq l'-1$, les flèches
$$h^{-d+lg-(l-l')} gr_{k,\pi_o} \longto h^{-d+lg-(l-l'-1)} gr_{k-1,\pi_o}$$
de la suite spectrale (\ref{suite-spectrale}), sont non nulles, et se déduisent des
suites exactes suivantes
\begin{multline*} \label{suites-exactes}
0 \to [\overleftarrow{l'-1},\overrightarrow{l-l'}]_{\pi_o} \to
[\overleftarrow{l'-1}]_{\pi_o}
\overrightarrow{\times} [\overrightarrow{l-l'-1}]_{\pi_o} \to  \\
\to [\overleftarrow{l'}]_{\pi_o} \overrightarrow{\times}
[\overrightarrow{l-l'-2}]_{\pi_o} \to
[\overleftarrow{l'+1},\overrightarrow{l-l'-2}]_{\pi_o} \to 0
\end{multline*}
pour $l' < l-1$ et pour $l'=l-1$ de la suite exacte courte:
\begin{equation*}
0 \to [\overleftarrow{l-2},\overrightarrow{1}]_{\pi_o} \to
[\overleftarrow{l-2}]_{\pi_o}
\overrightarrow{\times} [\overrightarrow{0}]_{\pi_o} \to [\overleftarrow{l-1}]_{\pi_o}
\to 0.
\end{equation*}
d'où le résultat.

\end{proof}

\rem On verra plus loin (cf. les remarques (\ref{rema-N}) et (\ref{rema-N2})) qu'on
pourrait en fait montrer
(\ref{theo-local-fil}) directement, sans utiliser l'opérateur de monodromie $N$.

\begin{coro}
Pour $0 \leq i \leq d-lg$ et $g >1$, $R^ij^{\geq lg}_* HT(g,l,\pi_o,\Pi_l)[d-lg]$ est,
dans $\FH(M_{s_o})$,
une somme directe sur tous les couples $(n,r)$ tels que $ng+r(g-1)=i$ des faisceaux de
type $HT(g,l+n+r)$

$$j^{\geq (l+n+r)g}_!HT(g,l+n+r,\pi_o,(\Pi_l \overrightarrow{\times}
[\overleftarrow{n-1}]_{\pi_o})
\overleftarrow{\times} [\overrightarrow{r-1}]_{\pi_o}) \otimes
\Xi^{\frac{n(g+1)+r(g-1)}{2}}$$ Pour $g=1$,
$R^ij^{\geq h}_* HT(1,h,\xi_o,\Pi_h)[d-h]$ est à support dans la tour des
$M_{I,s_o}^{\geq h+i}$ et sa
restriction à $M_{I,s_o}^{=h+i+r}$ est le faisceau de type $HT(1,h+i+r)$:
$$HT(1,h+i+r,\xi_o, (\Pi_h \overrightarrow{\times} [\overleftarrow{i-1}]_{\xi_o})
\overleftarrow{\times}
[\overleftarrow{r-1}]_{\xi_o}) \otimes \Xi^{i}$$
\end{coro}

\subsection{Pureté de la filtration de monodromie de $R\Psi_{\eta_o}(\bar \Qm_l)$}
\label{globalss3}

La proposition suivante correspond au théorème (\ref{theo-global1}) qui rappelons le
découle, en utilisant la
pureté de la filtration de monodromie, du théorème (\ref{theo-global0}) qui est prouvé
au corollaire
(\ref{global0-ok}). Ce paragraphe n'a d'intérêt que par son application au cas de la
caractéristique mixte.

\begin{prop} \label{prop-grk-final}
Pour tout $|k| < s_g$, $e_{\pi_o} gr_{k,\pi_o}$ est, dans $\FPH(M_{s_o})$, égal à
$$\sum_{\genfrac{}{}{0pt}{}{|k| < l \leq s_g}{l \equiv k-1 \mod 2}} \PC(g,l,\pi_o)
(-\frac{lg+k-1}{2})$$
Dans tous les autres cas $gr_{k,\pi_o}$ est nul.
\end{prop}

\begin{proof} Par rapport au corollaire (\ref{coro-grk}), il s'agit de traiter les
points supersinguliers et donc,
d'après le corollaire (\ref{coro-pnul}) de montrer que pour tout $r \equiv s-1 \mod 2$
et $|r| \leq s-1$,
$\PC(g,s,\pi_o)(-\frac{sg+r-1}{2})$ est un constituant de $e_{\pi_o} gr_{r,\pi_o}$.
\footnote{On remarquera
que le raisonnement suivant est valable pour toutes les strates de sorte qu'en
raisonnant par récurrence
l'argument du corollaire (\ref{coro-grk}) en utilisant le théorème de comparaison de
Berkovich-Fargues
n'était pas strictement nécessaire.}

On commence par remarquer, en utilisant la partie unipotente $N$ de la monodromie, qu'il
suffit en fait de
montrer que $\PC(g,s,\pi_o) (-\frac{s(g-1)}{2})$ est un constituant de $e_{\pi_o}
gr_{1-s,\pi_o}$.
Considérons alors $\PC(g,s,\pi_o) (-\frac{s(g-1)}{2})$ et soit $r$ tel qu'il soit un
constituant de
$e_{\pi_o} gr_{r,\pi_o}$. Soit $z_I$ un point géométrique de $M_{I,s_o}^{=d}$. D'après
la proposition
(\ref{prop-hij-ss-g}), les germes en $z_I$ des $h^{-1} gr_{k,I,\pi_o}$ sont tous de
poids strictement plus
grand que $s(g-1)$ sauf pour $k=2-s$. Par ailleurs le germe en $z_I$ de $R^{d-1}
\Psi_{\eta_o,I,\pi_o}(\bar{{\mathbb Q}_l})$ n'a aucun constituant de poids $s(g-1)$ de
sorte que $r \leq
2-s$.

\noindent - Si on avait $r \leq -s$, la partie unipotente $N$ de la monodromie donnerait
que
$$\PC(g,s,\pi_o)(-\frac{s(g-1)-2r}{2})$$
serait un constituant de $e_{\pi_o} R\Psi_{\eta_o,\pi_o}(\bar{{\mathbb Q}_l})[d-1]$ ce
qui n'est pas d'après
(\ref{prop-libre}) (\ref{prop-p}) et (\ref{coro-pnul}), car $1-s-2r > s-1$.

\noindent - Supposons $r=2-s$: pour $s>2$, par application de la dualité de Verdier et
de l'opérateur de monodromie $N$, on en
déduit que $\PC(g,s,\pi_o) \otimes L_g(\pi_o)(-\frac{s(g+1)-2}{2})$ et $\PC(g,s,\pi_o)
(-\frac{s(g+1)-4}{2})$ sont des constituants de $e_{\pi_o}gr_{s-2,\pi_o}$. Ainsi par une
nouvelle application de $N$,
$\PC(g,s,\pi_o) (-\frac{s(g+1)-4}{2})$ est aussi un constituant de $e_{\pi_o}
gr_{s-5,\pi_o}$ et donc devrait
apparaître deux fois ce qui n'est pas. Dans le cas $s=2$, la situation défavorable
correspondrait à
$\PC(g,2,\pi_o)(-g)$ et $\PC(g,2,\pi_o)(1-g)$ constituants de $e_{\pi_o} gr_{0,\pi_o}$
de sorte que dans la
cohomologie globale, la monodromie serait nulle ce qui est en contradiction avec la
proposition suivante pour
$s=2$.

\end{proof}

\begin{prop} \label{prop-mono}
Pour tout diviseur $g$ de $d=sg$ et toute représentation irréductible cuspidale de
$GL_g(F_o)$, il existe une
représentation irréductible automorphe $\Pi$ de $D_\Am^\times$ vérifiant $\hyp(\oo)$
telle que $\Pi_o \simeq
[\overleftarrow{s-1}]_{\pi_o}$ telle que l'opérateur de monodromie $N$ sur
$H^{d-1}_{\eta_o}[\Pi^\oo]$ soit
non nul.
\end{prop}

\begin{proof} En égale caractéristique, le résultat est connu d'après \cite{lrs}. En
caractéristique mixte, le
principe, qui nous a été suggéré par M. Harris, est de se ramener au cas Iwahori par
changement de base
résoluble. Ledit changement de base est expliqué dans l'appendice A, traitons alors le
cas Iwahori. Notons
que récemment Yoshida et Taylor, cf. \cite{y-t}, ont rédigé ce résultat qui par ailleurs
m'avait été expliqué
par A. Genestier. Cependant dans notre cas, vu que l'on ne s'intéresse qu'au cas $s=2$
et donc $g=1$, le
résultat a déjà été prouvé par Carayol, cf. \cite{ca}.
\end{proof}

En particulier on obtient la proposition suivante.

\begin{prop} \label{prop-mono-1}
Soit $\Pi$ une représentation irréductible automorphe de $D_\Am^\times$ telle que
$\Pi_o$ est l'induite
irréductible $[\overleftarrow{s_1-1}]_{\xi_1} \boxplus \cdots \boxplus
[\overleftarrow{s_r-1}]_{\xi_r}$ où
les $\xi_i$ sont des caractères de $F_o^\times$. En tant que représentation de $W_o$,
$H^{d-1}_{\eta_o}[\Pi^\oo]$ est la somme directe des
$$\bigoplus_{i=1}^r \sp_{s_i} \otimes \xi_i$$
où $\sp_s$ est la représentation de dimension $s$, $|-|^{(1-s)/2} \oplus \cdots \oplus
|-|^{(s-1)/2}$ où
l'indice de nilpotence de l'opérateur de monodromie $N$ est égal à $s$.
\end{prop}

\begin{proof} Le cas $s=2$ et $g=1$ de la proposition (\ref{prop-mono}) est prouvé par
Carayol dans \cite{ca} de
sorte que la proposition précédente est vraie pour $g=1$ et $s$ quelconque. On étudie
ensuite la suite
spectrale
\begin{equation} \label{ss-poids-hi1}
E_1^{i,j}[\Pi^\oo]=H^{i+j}(gr_{-i})[\Pi^\oo] \Rightarrow H^{i+j}_{\eta_o}[\Pi^\oo]
\end{equation}
Cette étude est faite dans le cas général au théorème (\ref{theo-mono-glob}); le point
est qu'en ce qui
concerne les $\Pi^\oo$-parties, nous ne devons considérer que les représentations
$\pi_o$ de $GL_1(F_o)$,
i.e. le cas $g=1$, qui rappelons le est connu grâce à Carayol, d'où le résultat.

\end{proof}

\rem Pour des résultats généraux sur les composantes locales $\Pi_o$ ainsi que sur la
partie
$\Pi^\oo$-isotypique des groupes de cohomologie de la fibre générique, on renvoie
respectivement aux
théorèmes (\ref{prop-compo-locale}) et (\ref{theo-mono-glob}).

\section{Étude de la suite spectrale des cycles évanescents}
\label{ssce0}

Le but de ce paragraphe est de prouver la proposition (\ref{cas-m}). Le principe de la
preuve est
d'étudier la suite spectrale des cycles évanescents en y
intégrant, via les suites spectrales associées à la stratification, la connaissance des
$H^i_c(M_{I,s_o}^{=lg}, \FC(g,l,\pi_o)\otimes \LC_{\r_\oo})$.

\subsection{Cas où $\Pi_o \simeq \st_s(\pi_o)$}

Soit $1 \leq g \leq d$ et $\pi_o$ une représentation irréductible cuspidale de
$GL_g(F_o)$. On fixe dans la
suite une représentation automorphe irréductible $\Pi$ de $D_\Am^\times$ vérifiant
$\hyp(\rho_\oo)$ et telle que
$\Pi_o \simeq [\overleftarrow{s-1}]_{\pi_o}$ pour $\pi_o$ une représentation
irréductible cuspidale de
$GL_g(F_o)$ et on note $m(\Pi)$ la multiplicité de $\Pi$ dans l'espace des formes
automorphes. On commence
par un résultat déjà présent dans \cite{lrs}.

\begin{coro} \label{coho-global1}
Les groupes de cohomologie de la fibre générique
$H^i_{\eta_o,\rho_\oo}[\Pi^{\oo,o}]$ sont nuls pour $i \neq d-1$ et
$$H^{d-1}_{\eta_o,\rho_\oo}[\Pi^{\oo,o}]=[\overleftarrow{s-1}]_{\pi_o} \otimes
L_d([\overleftarrow{s-1}]_{\pi_o}) (-\frac{sg-1}{2})$$
\end{coro}

\begin{proof} On considère la suite spectrale
$$E_{1,\rho_\oo}^{i,j}[\Pi^{\oo,o}]:=H^{i+j}(gr_{-i,\rho_\oo})[\Pi^{\oo,o}] \Rightarrow
H^{i+j+d-1}_{\eta_o,\rho_\oo}[\Pi^{\oo,o}]$$
On rappelle qu'en utilisant la pureté des $gr_k$, le corollaire (\ref{global0-ok})
implique le théorème
(\ref{theo-global1}). Sans utiliser monodromie-poids, en caractéristique mixte, le
théorème
(\ref{theo-global1}) découle de la proposition (\ref{prop-mono}). La proposition
(\ref{prop-coho1}) donne
alors la nullité de $E_{1,\rho_\oo}^{i,j}[\Pi^{\oo,o}]$ pour: $i+j \neq 0$, ou $|i| \geq
s$ ou $i \equiv s \mod 2$. Pour
$|i|<s$ et $i=s-1-2r$, on a
$$E_{\oo,\rho_\oo}^{s-1-2r,2r+1-s}[\Pi^{\oo,o}]=E_{1,\rho_\oo}^{s-1-2r,2r+1-s}[\Pi^{\oo,o}]=[\overleftarrow{s-1}]_{\pi_o}
\otimes
L_g(\pi_o) (-\frac{s(g-1)+2r}{2})$$ d'où le résultat.

\end{proof}

\begin{prop} \label{prop-hipsi}
Pour tout $0 \leq i \leq d-1$, les $GL_d(F_o) \otimes W_o$-modules
$$H^j(R^i\Psi_{\eta_o,\pi_o,\rho_\oo}(\bar \Qm_l))[\Pi^{\oo,o}]$$
vérifient les propriétés suivantes:
\begin{itemize}
\item[(1)] ils sont nuls si $g$ n'est pas un diviseur de $d$;

\item[(2)] pour $g$ un diviseur de $d=sg$, ils sont nuls si $j$ n'est pas de la forme
$d-lg$ pour $1 \leq l
\leq s$;

\item[(3)] pour $g$ un diviseur de $d=sg$ et $j=d-lg$ avec $1 \leq l \leq s$, ils sont
nuls si $i$ n'est pas
de la forme $lg-r$ avec $1 \leq r \leq l$;

\item[(4)] pour $g$ un diviseur de $d=sg$ et $1 \leq l < s$,
$H^{d-lg}(R^{lg-r} \Psi_{\eta_o,I,\pi_o,\rho_\oo}(\bar \Qm_l))[\Pi^{\oo,o}]$ est
isomorphe à
$$ m(\Pi) [\overleftarrow{l-r},
\overrightarrow{r-1}]_{\pi_o} \overrightarrow{\times} [\overleftarrow{s-l-1}]_{\pi_o}
\otimes L_g(\pi_o)
(-\frac{s(g-1)+2(l-r)}{2})$$

\end{itemize}
\end{prop}

\begin{proof} On utilise la suite spectrale associée à la stratification
\begin{multline} \label{sss}
E_{1,I,\pi_o,\rho_\oo}^{p,q;i}=H_c^{p+q}(M_{I,\bar s_o}^{=p-1},R^i
\Psi_{\eta_o,I,\pi_o,\rho_\oo}(\bar
\Qm_l)) \\ \Rightarrow H_c^{p+q}(M_{I,\bar s_o},R^i \Psi_{\eta_o,I,\pi_o,\rho_\oo}(\bar
\Qm_l))
\end{multline}
On rappelle que d'après le corollaire (\ref{coro-hic}),
$E_{1,I,\pi_o,\rho_\oo}^{p,q;i}[\Pi^{\oo,o}]$ est non nul si
et seulement si:
\begin{itemize}
\item $g$ est un diviseur de $d=sg$,

\item $p-1=lg$ pour $1 \leq l \leq s$,

\item $p+q=d-lg$,

\item $i=lg-r$ avec $1 \leq r \leq l$.
\end{itemize}
Les points (1), (2) et (3) en découlent alors directement. Par ailleurs on a
\begin{multline*}
\lim_{\genfrac..{0pt}{1}{\longto}{I}}
E_{1,I,\pi_o,\rho_\oo}^{lg+1,d-2lg-1;lg-r}[\Pi^{\oo,o}]\simeq m(\Pi)
[\overleftarrow{l-r},
\overrightarrow{r-1}]_{\pi_o} \overrightarrow{\times} [\overleftarrow{s-l-1}]_{\pi_o}
\\
\otimes L_g(\pi_o) (-\frac{s(g-1)-2(r-l)}{2})
\end{multline*}
de sorte que pour tout $k\geq 1$, les flèches
$d_k^{p,q;i}:E_{k,\pi_o,\rho_\oo}^{p,q;i}[\Pi^{\oo,o}] \longto
E_{k,\pi_o,\rho_\oo}^{p+k,q+k-1}[\Pi^{\oo,o}]$ de (\ref{sss}) sont toutes nulles. En
effet pour que
$E_{k,\pi_o,\rho_\oo}^{p,q;i}[\Pi^{\oo,o}]$ (resp.
$E_{k,\pi_o,\rho_\oo}^{p+k,q+k-1;i}[\Pi^{\oo,o}]$) soit non nul, il faut
qu'il existe $1 \leq l_1 \leq s$ et $1 \leq r_1 \leq l_1$ (resp. $1 \leq l_2 \leq s$ et
$1 \leq r_2 \leq
l_2$) avec
$$(p,q,i)=(l_1g+1,d-2l_1g-1,l_1g-r_1)$$
$$(\hbox{resp. } (p+k,q+k-1,i)=(l_2g+1,d-2l_2g-1,l_2g-r_2)).$$
Ce qui donne $1=3g(l_2-l_1)$; on voit alors que pour tout $k\geq 1$,
$E_{k,\pi_o,\rho_\oo}^{p,q;i}[\Pi^{\oo,o}]$ et
$E_{k,\pi_o,\rho_\oo}^{p+k,q+k-1}[\Pi^{\oo,o}]$ ne peuvent pas être tous deux non nuls
de sorte que
$E_{\oo,I,\pi_o,\rho_\oo}^{p,q;i}[\Pi^{\oo,o}]=E_{1,I,\pi_o,\rho_\oo}^{p,q;i}[\Pi^{\oo,o}]$
d'où le résultat d'après le
corollaire (\ref{coro-hic}).

\end{proof}

\begin{coro} \label{coro-1}
La partie de poids $d-s$ de
$\widetilde{\UC_{F_o}^{d,d-i}}(\JL^{-1}([\overleftarrow{s-1}]_{\pi_o}))$ est un
constituant de $\Pi_i \otimes L_g(\pi_o) (-\frac{d-s}{2})$ où
$$\Pi_i= \left \{ \begin{array}{cl} [\overleftarrow{s-1}]_{\pi_o} & \hbox{pour } i=1 \\
{[}\overrightarrow{0}]_{\pi_o} \overrightarrow{\times} [\overleftarrow{s-2}]_{\pi_o} &
\hbox{pour } i=2 \\
\cdots \\
{[}\overrightarrow{i-2}]_{\pi_o} \overrightarrow{\times} [\overleftarrow{s-i}]_{\pi_o} &
\hbox{pour } i \\
\cdots \\
{[}\overrightarrow{s-2}]_{\pi_o} \overrightarrow{\times} [\overleftarrow{0}]_{\pi_o} &
\hbox{pour } i=s
\end{array} \right.$$
En outre $\widetilde{\UC_{F_o}^{d,d-s}}(\JL^{-1}([\overleftarrow{s-1}]_{\pi_o}))$
contient
$[\overrightarrow{s-1}]_{\pi_o} \otimes L_g(\pi_o) (-\frac{s(g-1)}{2})$ et la partie de
poids $s(g-1)$ de
$\sum_{i=1}^s (-1)^i
\widetilde{\UC_{F_o}^{d,d-i}}(\JL^{-1}([\overleftarrow{s-1}]_{\pi_o}))$ est égale à
$(-1)^s [\overrightarrow{s-1}]_{\pi_o} \otimes L_g(\pi_o) (-\frac{s(g-1)}{2})$.

\end{coro}

\begin{proof} On étudie la $\Pi^{\oo,o}$-partie de la suite spectrale des cycles
évanescents pour $\Pi$ vérifiant
les propriétés du début de ce paragraphe:
\begin{multline} \label{ssce}
E_{2,\pi_o,\rho_\oo}^{p,q}[\Pi^{\oo,o}]=\lim_{\genfrac..{0pt}{1}{\longto}{I}}
H^p(M_{I,\bar
s_o},R^q\Psi_{\eta_o,I,\pi_o,\rho_\oo}(\bar \Qm_l))[\Pi^{\oo,o}] \\ \Rightarrow
H^{p+q}_{\eta_o,\pi_o,\rho_\oo}[\Pi^{\oo,o}]
\end{multline}

\begin{lemm} \label{lemm-prec}
Pour tout $1 < r < s$, la partie de poids $s(g-1)$ de
$E_{2,\pi_o,\rho_\oo}^{0,d-r}[\Pi^{\oo,o}]$
est un constituant de $E_{2,\pi_o,\rho_\oo}^{d-(r-1)g,(r-1)(g-1)}[\Pi^{\oo,o}]$, i.e.
de
$$[\overrightarrow{r-2}]_{\pi_o} \overrightarrow{\times} [\overleftarrow{s-r}]_{\pi_o}
\otimes L_g(\pi_o)
(-\frac{s(g-1)}{2});$$ celle de $E_{2,\pi_o}^{0,d-1}[\Pi^{\oo,o}]$ est nulle.
\end{lemm}

\begin{proof} D'après la proposition précédente, les
$E_{2,\pi_o,\rho_\oo}^{p,q}[\Pi^{\oo,o}]$ de poids $s(g-1)$ pour $p \neq
0$, non nuls, sont
$$E_{2,\pi_o,\rho_\oo}^{d-rg,r(g-1)}[\Pi^{\oo,o}] \simeq [\overrightarrow{r-1}]_{\pi_o}
\overrightarrow{\times}
[\overleftarrow{s-r-1}]_{\pi_o} \otimes L_g(\pi_o) (-\frac{s(g-1)}{2})$$ Or, d'après le
corollaire
(\ref{coho-global1}), pour $1 \leq r < s$, la partie de poids $s(g-1)$ de
$E_{\oo,\pi_o,\rho_\oo}^{d-r}[\Pi^{\oo,o}]$
est nulle d'où le résultat.

\end{proof}

\begin{lemm}
Pour $g >1$, le terme $E_{2,\pi_o}^{0,d-r}[\Pi^{\oo,o}]$ de la suite spectrale des
cycles évanescents est
égal à $\hom_{\bar D_o^\times}((\CC^\oo_{\bar
D})^\vee,\widetilde{\UC^{d,d-r}_{F_o}})[\Pi^{\oo,o}]$.
\end{lemm}

\begin{proof} On étudie la suite spectrale associée à la stratification pour
$R^{d-r} \Psi_{\eta_o,\pi_o,\rho_\oo}(\bar \Qm_l)$. D'après la proposition précédente, pour tout $1 \leq l <
s$, $H^0_c(M_{I,s_o}^{=lg},R^{d-r} \Psi_{\eta_o,I,\pi_o,\rho_\oo}(\bar \Qm_l))[\Pi^{\oo,o}]$ est nul; pour $g
\neq 1$, il en est de même de $H^1_c(M_{I,s_o}^{=lg},R^{d-r} \Psi_{\eta_o,I,\pi_o,\rho_\oo}(\bar
\Qm_l))[\Pi^{\oo,o}]$, d'où le résultat.

\end{proof}

\begin{lemm} Pour $g=1$, $\hom_{\bar D_o^\times}((\CC^\oo_{\bar
D})^\vee,\widetilde{\UC^{d,d-r}_{F_o}})[\xi_o)]$
est un constituant de
$$[\overrightarrow{r-2}]_{\xi_o} \overrightarrow{\times} [\overleftarrow{d-r}]_{\xi_o}
\otimes \xi_o.$$
\end{lemm}

\begin{proof} Pour $g=1$, on remarque que le seul $H^1_c(M_{I,s_o}^{=h},R^{d-r}
\Psi_{\eta_o,I,\xi_o,\rho_\oo}(\bar \Qm_l))$
ayant une partie de poids $0$ non nulle, est pour $h=d-1$ et $r=d$, de
$\Pi^{\oo,o}$-composante isotypique
égale à $[\overrightarrow{d-2}]_{\xi_o} \overrightarrow{\times}
[\overleftarrow{0}]_{\xi_o} \otimes \xi_o$ ce
qui correspond à la contribution de $E_{2,\xi_o,\rho_\oo}^{1,0}[\Pi^{\oo,o}]$ dans le
lemme (\ref{lemm-prec}), d'où le
résultat.

\end{proof}

Ainsi pour tout $1 \leq i \leq s$, la partie de poids $d-s$ de
$\widetilde{\UC_{F_o}^{d,d-i}}(\JL^{-1}([\overleftarrow{s-1}]_{\pi_o}))$ est un
constituant de $\Pi_i \otimes
L_g(\pi_o) (-\frac{d-s}{2})$ où $\Pi_i$ est comme dans l'énoncé. Par ailleurs on
remarque que
$E_{2,\pi_o,\rho_\oo}^{d-s+1,(s-1)(g-1)}[\Pi^{\oo,o}]$ contient
$[\overrightarrow{s-1}]_{\pi_o} \otimes L_g(\pi_o)
(-\frac{s(g-1)}{2})$ alors que $E_{2,\pi_o}^{d-s+2,(s-2)(g-1)}[\Pi^{\oo,o}]$ non; on en
déduit alors que
$\widetilde{\UC_{F_o}^{d,d-i}}(\JL^{-1}(\pi_o))$ contient
$[\overrightarrow{s-1}]_{\pi_o} \otimes L_g(\pi_o)
(-\frac{s(g-1)}{2})$. Le calcul sur la somme alternée correspond à (\ref{inutile}).

\end{proof}

\subsection{Involution de Zelevinski et première preuve de (\ref{cas-m})}

Il est possible de prouver la proposition (\ref{cas-m}) sans plus d'étude cohomologique
en utilisant le résultat suivant.

\begin{theo} \label{involution}
Pour toute représentation irréductible cuspidale $\pi_o$ de $GL_g(F_o)$, on a, pour tout
$s \geq 1$ et pour tout $i$, un isomorphisme canonique
$$\left ( \UC_{F_o,l,sg}^{sg-1+\bullet}(\JL^{-1}(\st_s(\pi_o))) \right
)^{\vee,\iota}(d-1) \simeq \UC{F_o,l,sg}^{sg+s-2-2
\bullet}(\JL^{-1}(\st_s(\pi_o^\vee)))$$
où dans le membre de droite l'exposant $\vee,\iota$ désigne le dual composé avec
l'involution
de Zelevinski, $\iota$, sur $GL_{sg}(F_o)$.
\end{theo}

\rem Laurent Fargues a une preuve de ce résultat qui utilise tout ou partie
l'isomorphisme de Faltings, à partir d'un résultat similaire du coté de
l'espace de Drinfeld.

Ainsi avec les notations du corollaire (\ref{coro-1}), s'il existait $1 \leq i \leq s$
tel que $\Pi_i$ admette un constituant autre que
$[\overrightarrow{s-1}]_{\pi_o}$, on en déduirait que
$\UC_{F_o,l,sg}^{sg-s-1+i}(\JL^{-1}([\overleftarrow{s-1}]_{\pi_o}))$ aurait un
constituant de poids
$sg+s-2$ autre que $[\overleftarrow{s-1}]_{\pi_o^\vee}$. Or, d'après les corollaires
(\ref{coro-sec-fp}) et (\ref{coro-pnul}), pour tout $k$, les germes
aux points supersinguliers des $h^i gr_{k,\pi_o^\vee}$ sont tous de poids
strictement plus petit que $sg+s-2$ ou alors égales à
$[\overleftarrow{s-1}]_{\pi_o^\vee} \otimes \rec_{F_o}(\pi_o) (-\frac{sg+s-2}{2})$, d'où
le résultat.

\subsection{Cas  $\Pi_o \simeq \speh_s(\pi_o)$}

Nous allons prouver la proposition (\ref{cas-m}) sans utiliser le théorème
(\ref{involution}). Pour cela nous revenons à l'étude de la suite spectrale
des cycles évanescents.
On peut remarquer que la connaissance des $E_{2,\pi_o,\rho_\oo}^{p,q}[\Pi^{\oo,o}]$ de
la suite spectrale (\ref{ssce}) ne nous
fournit pas le théorème local. La figure (\ref{figu1}) illustre ce fait dans le cas
$s=2$, où pour $\pi_o$
une représentation irréductible cuspidale de $GL_{d/2}(F_o)$, l'on n'arrive pas à
exclure le cas:
\begin{itemize}
\item
$\widetilde{\UC_{F_o}^{d,d-1}}(\JL^{-1}([\overleftarrow{1}]_{\pi_o}))=[\overleftarrow{1}]_{\pi_o}
\otimes
(L_g(\pi_o)(-\frac{d-1}{2}) \otimes \sp_2)$;

\item $\widetilde{\UC_{F_o}^{d,d-2}}(\JL^{-1}([\overleftarrow{1}]_{\pi_o}))=
[\overleftarrow{0}]_{\pi_o}
\overrightarrow{\times} [\overleftarrow{0}]_{\pi_o} \otimes L_g(\pi_o)
(-\frac{d-2}{2})$.
\end{itemize}

\input{figu1.tex}

\noindent Si la flèche indiquée est un isomorphisme, on obtient bien le bon
aboutissement.

\medskip

\marque Dans le cas où l'on considère une représentation automorphe $\Pi$ de
$D_\Am^\times$ vérifiant
$\hyp(\rho_\oo)$ et telle que $\Pi_o \simeq [\overrightarrow{s-1}]_{\pi_o}$, on
calculera les termes
$E_{2,\rho_\oo}^{p,q}[\Pi^{\oo,o}]$ de la suite spectrale des cycles évanescents et on
montrera alors que le
théorème local en découle. Dans la figure (\ref{figu2}), on illustre comment on exclut
le cas défavorable
ci-avant, en remarquant que $[\overleftarrow{0}]_{\pi_o} \overrightarrow{\times}
[\overleftarrow{0}]_{\pi_o}$
n'est pas isomorphe à $[\overleftarrow{0}]_{\pi_o} \overleftarrow{\times}
[\overleftarrow{0}]_{\pi_o}$,
contredisant le théorème de Lefschetz difficile. On pourra aussi se référer aux figures
(\ref{figure5}) et
(\ref{figure6}) qui détaillent pour $s=4$ et $g=2$, les $E_{2,\pi_o}^{p,q}[\Pi^{\oo,o}]$
respectivement dans
les cas $\Pi_\oo=[\overleftarrow{s-1}]_{\pi_\oo}$ et
$\Pi_\oo=[\overrightarrow{s-1}]_{\pi_\oo}$.

\input{figu2.tex}

\begin{rema} \label{rema-hyp2}
Si $\Pi$ vérifie $\hyp(\rho_\oo)$, la condition $\Pi_o \simeq
[\overrightarrow{s-1}]_{\pi_o}$ pour $\pi_o$ une
représentation irréductible cuspidale de $GL_g(F_o)$, implique qu'il existe une
représentation irréductible
cuspidale $\pi_\oo$ de $GL_{g'}(F_\oo)$ avec $d=s'g'$ telle que $\Pi_\oo \simeq
[\overleftarrow{s'-1}]_{\pi_\oo}$ et donc $\Pi$ vérifie $\hyp(\r_\oo)$ avec
$\r_\oo=\JL^{-1}([\overleftarrow{s'-1}]_{\pi_\oo})$ et $s=s'$.
\end{rema}

\marque On considère \footnote{cf. aussi les conditions qui précèdent le corollaire
(\ref{coro-pnul})} dans
la suite $\Pi$ une représentation irréductible automorphe de $D_\Am^\times$ vérifiant
$\hyp(\rho_\oo)$ et
telle que $\Pi_o \simeq [\overrightarrow{s-1}]_{\pi_o}$, pour $\pi_o$ une représentation
irréductible
cuspidale de $GL_g(F_o)$. On suppose par ailleurs que l'ensemble des représentations
irréductibles
automorphes $\bar \Pi$ de $\bar D_\Am^\times$ telles que $\bar \Pi^{\oo,o} \simeq
\Pi^{\oo,o}$ est réduite à
une unique représentation, en particulier $m(\Pi)=m(\bar \Pi)=1$, et $\bar \Pi_o \simeq
\JL^{-1}([\overleftarrow{s-1}]_{\pi_o})$ et $\bar \Pi_\oo \simeq
\JL^{-1}([\overleftarrow{s-1}]_{\pi_\oo})$.
L'objectif est de calculer les termes $E_{2,\rho_\oo}^{p,q}[\Pi^{\oo,o}]$ de la suite
spectrale des cycles
évanescents.

\begin{prop} \label{prop-not}
Pour tout $1 \leq l \leq s$,
$$H^i(j^{\geq lg}_{!*}
HT_{\rho_\oo}(g,l,\pi_o,[\overleftarrow{l-1}]_{\pi_o})[d-lg])[\Pi^{\oo,o}]$$ est nul
pour $|i|
> s-l$ ou $i \not \equiv s-l \mod 2$ et sinon, il est isomorphe à
$$([\overleftarrow{l-1}]_{\pi_o} \overrightarrow{\times}
[\overrightarrow{\frac{s-l-i}{2}-1}]_{\pi_o})
\overleftarrow{\times} [\overrightarrow{\frac{s-l+i}{2}-1}]_{\pi_o} \otimes
(\Xi^{\frac{(s-l)g+i}{2}} \otimes
\bigoplus_{\xi \in \AF(\pi_o)} \xi^{-1})$$ en tant que représentation de $GL_d(F_o)
\times \Zm$, où
$\AF(\pi_o)$ est l'ensemble des caractères $\xi:\Zm \longto \Qm_l^\times$, tels que
$\pi_o \otimes \xi^{-1}
\circ \val(\det) \simeq \pi_o$.
\end{prop}

\begin{proof} On raisonne par récurrence pour $l$ variant de $s$ à $1$, le cas $l=s$ est
donné par le lemme
(\ref{lem-pts-ss}) et les conditions imposées ci-dessus à $\Pi$. Les suites exactes
courtes
\begin{multline*}
0 \to P_{\pi_o,l,0}([\overleftarrow{l-1}]_{\pi_o}) \longto j^{\geq lg}_! HT(g,l,\pi_o,
[\overleftarrow{l-1}]_{\pi_o})[d-lg] \\ \longto j^{\geq lg}_{!*}
HT(g,l,\pi_o,[\overleftarrow{l-1}]_{\pi_o})[d-lg] \to 0
\end{multline*}
$$0 \to P_{\pi_o,l,i}([\overleftarrow{l-1}]_{\pi_o}) \longto
P_{\pi_o,l,i-1}([\overleftarrow{l-1}]_{\pi_o})
\longto \PC_-(g,l+i,\pi_o,l,[\overleftarrow{l-1}]_{\pi_o}) \to 0$$ pour $1 \leq i \leq
s-l-1$, fournissent
\begin{multline} \label{egalite-purete}
\sum_i (-1)^i H^i(j_{!*}^{\geq lg}
HT_{\r_\oo}(g,l,\pi_o,[\overleftarrow{l-1}]_{\pi_o})[d-lg])[\Pi^{\oo,o}]= \\
\sum_{r=1}^{s-l} (-1)^r \sum_i (-1)^i H^i([\overleftarrow{l-1}]_{\pi_o}
\overrightarrow{\times}
[\overleftarrow{r-1}]_{\pi_o})[\Pi^{\oo,o}] \otimes \Xi^{r/2}+ \\
\sum_i (-1)^i H^i(j_!^{\geq lg}
HT_{\r_\oo}(g,l,\pi_o,[\overleftarrow{l-1}]_{\pi_o})[d-lg])[\Pi^{\oo,o}]
\end{multline}
où par simplification, on écrit $H^i([\overleftarrow{l-1}]_{\pi_o}
\overrightarrow{\times}
[\overleftarrow{r-1}]_{\pi_o})$ pour
$$\lim_{\genfrac..{0pt}{1}{\to}{I}} H^i(M_{I,s_o}^{\geq (l+r)g}, j^{\geq (l+r)g}_{!*}
HT_{\rho_\oo}(g,l+r,\pi_o,
[\overleftarrow{l-1}]_{\pi_o} \overrightarrow{\times}
[\overleftarrow{r-1}]_{\pi_o})[d-(l+r)g])$$

D'après l'hypothèse de récurrence et (\ref{somme-alternee}), pour tout $i \equiv 1 \mod
2$, la partie de
poids $(s-l)(g-1)+i$ du membre de droite de (\ref{egalite-purete}) est nulle; il en est
de même pour
$i>2(s-l)$ ou $i<0$, tandis que pour $0 \leq i=2k < 2(s-l)$ (resp. $i=2(s-l)$) celle-ci
est égale à

$$\sum_{r=1}^{s-l-k} H^{l-s+2k+r}([\overleftarrow{l-1}]_{\pi_o} \overrightarrow{\times}
[\overleftarrow{r-1}]_{\pi_o})[\Pi^{\oo,o}] \otimes \Xi^{r/2}$$ (resp. à la partie de
poids $(s-l)(g+1)$ de
(\ref{somme-alternee})), ce qui donne $\pi \otimes (\Xi^{\frac{s(g-1)+2k}{2}} \otimes
\bigoplus_{\xi \in
\AF(\pi_o)} \xi^{-1})$ où
$$\begin{array}{rl} \pi = & \bigl ( [\overleftarrow{l-1}]_{\pi_o}
\overrightarrow{\times} ((-1)^{s-l-k}
[\overleftarrow{s-l-k-1}]_{\pi_o} \\ & +\sum_{r=1}^{s-l-k-1} (-1)^{r-1}
[\overleftarrow{r-1}]_{\pi_o}
\overrightarrow{\times} [\overrightarrow{s-l-r-k-1}]_{\pi_o}) \bigr )
\overleftarrow{\times} [\overrightarrow{k-1}]_{\pi_o} \\
= & ([\overleftarrow{l-1}]_{\pi_o} \overrightarrow{\times}
[\overrightarrow{s-l-k-1}]_{\pi_o})
\overleftarrow{\times} [\overrightarrow{k-1}]_{\pi_o} \\
& (\hbox{resp. } [\overleftarrow{l-1}]_{\pi_o} \overleftarrow{\times}
[\overrightarrow{s-l-1}]_{\pi_o}).
\end{array}$$
Or $H^i(j_{!*}^{\geq lg} HT_{\r_\oo}(g,l,\pi_o,[\overleftarrow{l-1}]_{\pi_o}) [d-lg])$
est pur de poids
$(s-l)g+i$ de sorte que son semi-simplifié est égal à celui de l'énoncé. On conclut
alors à l'égalité des
représentations, et pas seulement de leur semi-simplifiée, en remarquant que les strates
étant induites,
l'espace précédent, en tant que représentation de $GL_d(F_o)$ est de la forme
$$\ind_{P_{lg,d}^{op}(F_o)}^{GL_d(F_o)} [\overleftarrow{l-1}]_{\pi_o((s-l)g+i)/2}
\otimes \pi'$$
pour une certaine représentation $\pi'$ de $GL_{d-lg}(F_o)$.

\end{proof}

\begin{coro} \label{coro-ssce-min}
Pour tout $i \neq d-s$, la partie de poids $s(g-1)$ de
$H^i_{\eta_o,\r_\oo}[\Pi^{\oo,o}]$ est nulle alors que
pour $i=d-s$ elle est égale à $[\overrightarrow{s-1}]_{\pi_o} \otimes L_g(\pi_o)
(-\frac{s(g-1)}{2})$.
\end{coro}

\begin{rema} \label{rema-N} En fait on peut à ce stade déterminer complètement les
$H^i_{\eta_o,\r_\oo}$, cependant comme on l'a
déjà remarqué seule la connaissance des parties de poids $s(g-1)$ nous est nécessaire.
Pour un énoncé
complet, on pourra voir la proposition (\ref{prop-lrs2}).
\end{rema}

\begin{proof} On écrit $i$ sous la forme $d-1-\d$ et on étudie la suite spectrale
\begin{equation} \label{ss-poids}
E_{1,\r_\oo}^{i,j}:=H^{i+j}(gr_{-i,\r_\oo}) \Rightarrow H^{d-1+i+j}_{\eta_o,\r_\oo}
\end{equation}
qui, d'après la pureté, dégénère en $E_2$. On pourra se référer à la figure
(\ref{figure8}) où l'on a
représenté les $H^i(gr_{k,\r_\oo})[\Pi^{\oo,o}]$ pour $s=4$. D'après les propositions
(\ref{prop-p}) (ii) et
(\ref{prop-grk-final}), on a
$$e_{\pi_o} H^i(gr_{k,\r_\oo})[\Pi^{\oo,o}]=\bigoplus_{\genfrac{}{}{0pt}{}{|k|<l \leq
s}{l \equiv k+1 \mod 2}} H^i(\PC(g,l,\pi_o)
(-\frac{lg-1+k}{2}))[\Pi^{\oo,o}]$$ Ainsi d'après la proposition (\ref{prop-not}), pour
$\d>0$, la partie de
poids $s(g-1)$ de $H^{d-1+\d}_{\eta_o,\r_\oo}[\Pi^{\oo,o}]$ est nulle; précisément pour
$\d >0$, les poids de
$H^{d-1+\d}_{\eta_o,\r_\oo}[\Pi^{\oo,o}]$ sont parmi les $k=s(g-1)+2 \d +2r$ avec $0
\leq r < s-\d$ et sa
partie de poids $s(g-1)+2\d$, que l'on notera avec un indice, est un quotient de
\begin{multline*}
E_{1,\r_\oo}^{1-s+\d,s-1}[\Pi^{\oo,o}]_{s(g-1)+2\d}=\frac{1}{e_{\pi_o}}
H^{\d}(\PC(g,s-\d,\pi_o)
(-\frac{(s-\d)(g-1)}{2})) \\
=[\overleftarrow{s-\d-1}]_{\pi_o} \overleftarrow{\times} [\overrightarrow{\d-1}]_{\pi_o}
\otimes L_g(\pi_o)
(-\frac{s(g-1)+2\d}{2})
\end{multline*}

Par ailleurs, pour $\d>0$, d'après (\ref{prop-not}), la partie de poids $s(g-1)$ de
$H^{d-s-\d}_{\eta_o,\r_\oo}[\Pi^{\oo,o}]$ est un sous-quotient de
\begin{multline*}
E_{1,\r_\oo}^{1-s+\d,s-1-2\d}[\Pi^{\oo,o}]_{s(g-1)}=\frac{1}{e_{\pi_o}}
H^{-\d}(\PC(g,s-\d,\pi_o)
(-\frac{(s-\d)(g-1)}{2})) \\
=[\overleftarrow{s-\d-1}]_{\pi_o} \overrightarrow{\times}
[\overrightarrow{\d-1}]_{\pi_o} \otimes L_g(\pi_o)
(-\frac{s(g-1)}{2})
\end{multline*}

On utilise alors le théorème de Lefschetz difficile dont on rappelle l'énoncé ci-après.

\begin{theo} (cf. \cite{lrs} 14.19) Il existe une classe $h \in H^2_{\eta_o,\r_\oo}
(1)$,
invariante sous les actions de $(D_\Am^\oo)^\times$ et $W_o$, telle que les applications
itérées du cup
produit
$$h^i: H^{d-1-i}_{\eta_o,\r_\oo} \longto H^{d-1+i}_{\eta_o,\r_\oo} (i)$$
sont des isomorphismes.
\end{theo}

- Ainsi pour $s>2$, on observe que si la partie de poids $s(g-1)$ de
$H^{d-1-\d}_{\eta_o}[\Pi^{\oo,o}]$, pour
$1 < \d <s-1$, est non nulle, ses constituants sont de la forme
$[\overleftrightarrow{s-2},\overrightarrow{1}]_{\pi_o} \otimes L_g(\pi_o)
(-\frac{s(g-1)}{2})$ alors qu'un
constituant non nul de poids $s(g-1)+2\d$ de $H^{d-1+2\d}_{\eta_o,\r_\oo}[\Pi^{\oo,o}]$,
s'il existe, est de
la forme $[\overleftrightarrow{s-2},\overleftarrow{1}]_{\pi_o} \otimes L_g(\pi_o)
(-\frac{s(g-1)+2\d}{2})$.
La contradiction découle alors du théorème de Lefschetz difficile et de l'observation de
l'orientation de la
dernière flèche.

- Pour $\d=1$ et $s>2$, la partie de poids $s(g-1)+2$ de
$H^{d}_{\eta_o,\r_\oo}[\Pi^{\oo,o}]$ est un quotient
de $[\overleftarrow{s-2}]_{\pi_o} \overleftarrow{\times} [\overrightarrow{0}]_{\pi_o}
\otimes L_g(\pi_o)
(-\frac{s(g-1)+2}{2})$ alors que la partie de poids $s(g-1)$ de
$H^{d-2}_{\eta_o,\r_\oo}[\Pi^{\oo,o}]$ est un
constituant de $[\overleftarrow{s-2}]_{\pi_o} \overrightarrow{\times}
[\overrightarrow{0}]_{\pi_o} \otimes
L_g(\pi_o) (-\frac{s(g-1)}{2})$. La contradiction découle alors du théorème de Lefschetz
difficile et du fait
que $[\overleftarrow{s-1}]_{\pi_o}$ n'est pas un quotient de
$[\overleftarrow{s-2}]_{\pi_o}
\overleftarrow{\times} [\overrightarrow{0}]_{\pi_o}$.

- Pour $\d=s-1$ et $s \geq 2$, on observe que
$E_{2,\r_\oo}^{d-s}[\Pi^{\oo,o}]=H^{d-s}_{\eta_o,\r_\oo}[\Pi^{\oo,o}]$ est un
sous-espace de

$$E_{1,\r_\oo}^{0,1-s}[\Pi^{\oo,o}]=\frac{1}{e_{\pi_o}} H^{1-s}(\PC(g,1,\pi_o))
=[\overleftarrow{0}]_{\pi_o} \overrightarrow{\times} [\overrightarrow{s-2}]_{\pi_o}
\otimes L_g(\pi_o)
(-\frac{s(g-1)}{2})$$ alors que
$E_{2,\r_\oo}^{d-s}[\Pi^{\oo,o}]=H^{d+s-2}_{\eta_o,\r_\oo}[\Pi^{\oo,o}]$ est
un quotient de
$$E_{1,\r_\oo}^{0,s-1}[\Pi^{\oo,o}]=\frac{1}{e_{\pi_o}} H^{s-1}(\PC(g,1,\pi_o))
=[\overleftarrow{0}]_{\pi_o} \overleftarrow{\times} [\overrightarrow{s-2}]_{\pi_o}
\otimes L_g(\pi_o)
(-\frac{s(g+1)-2}{2})$$ Ainsi d'après Lefschetz difficile, s'ils sont non nuls, ils
doivent être égaux à
$[\overrightarrow{s-1}]_{\pi_o}$. Par ailleurs on remarque que ce dernier n'est pas un
constituant de
$H^{d-s+1}_{\eta_o,\r_\oo}[\Pi^{\oo,o}]$, ni de
$H^{d+s-3}_{\eta_o,\r_\oo}[\Pi^{\oo,o}]$, de sorte que
$[\overrightarrow{s-1}]_{\pi_o} \otimes L_g(\pi_o) (-\frac{s(g-1)}{2})$ est
effectivement un constituant de
$H^{d-s}_{\eta_o,\r_\oo}$, d'où le résultat.

\end{proof}

\begin{prop} \label{prop-ssce-poids} Pour tout $p \neq 0$ et tout $q$, les parties de
poids $s(g-1)$ des
$E_{2,\r_\oo}^{p,q}[\Pi^{\oo,o}]$ de la suite spectrale des cycles évanescents
\begin{equation} \label{ssce2}
E_{2,\r_\oo}^{p,q}=\lim_{\genfrac..{0pt}{1}{\longto}{I}} H^p(M_{I,\bar
s_o},R^q\Psi_{\eta_o,I}(\LC_{\r_\oo}))
\Rightarrow H^{p+q}_{\eta_o,\r_\oo}
\end{equation}
sont nulles.
\end{prop}

\begin{proof} Le principe est d'étudier les $E_{2,\r_\oo}^{p,q}[\Pi^{\oo,o}]$ via les
suites spectrales associées à la
stratification
\begin{equation} \label{sss2}
E_{1,I,\r_\oo}^{p,q;i}=H_c^{p+q}(M_{I,\bar s_o}^{=p-1},R^i
\Psi_{\eta_o,I}(\LC_{\r_\oo})) \Rightarrow
H_c^{p+q}(M_{I,\bar s_o},R^i \Psi_{\eta_o,I}(\LC_{\r_\oo}))
\end{equation}
Ainsi le résultat découle simplement de la proposition suivante.

\begin{prop} \label{prop-hic-poids}
Pour tout $1 \leq l < s$ et pour tout $i$, la partie de poids $(s-l)(g-1)$ de
$H^i_c(M_{I,s_o}^{=lg},\FC(g,l,\pi_o,I) \otimes \LC_{\r_\oo})[\Pi^{\oo,o}]$ est nulle.
\end{prop}

\begin{proof} Il s'agit dans un premier temps d'étudier les parties de poids
$(s-l)(g-1)$ des
$H^i_c(M_{I,s_o}^{=lg},HT_{\r_\oo}(g,l,\pi_o,\Pi_l,I))[\Pi^{\oo,o}]$ pour une
représentation $\Pi_l$
quelconque de $GL_{lg}(F_o)$.

\begin{lemm} \label{lem-rj-combi}
Pour tout $1 \leq l < s$, les
$$\lim_{\genfrac..{0pt}{1}{\to}{I}} H^i(M_{I,s_o}^{\geq lg},(j^{\geq lg}_!
HT_{\r_\oo}(g,l,\pi_o,\Pi_l,I)[(s-l)g]
))[\Pi^{\oo,o}]$$ vérifient les propriétés suivantes:

\begin{itemize}
\item[(i)] ils sont de la forme $\bigoplus_\xi(\ind_{P_{lg,d}^{\op}(F_o)}^{GL_d(F_o)}
(\Pi_l \otimes \xi) \otimes \pi_\xi)
\otimes \xi^{-1}$, en tant que représentation de $GL_d(F_o) \times \Zm$, où $\xi$ décrit
les caractères $\Zm
\longto \bar \Qm_l^\times$ et $\pi_\xi$ est une représentation de $GL_{(s-l)g}(F_o)$;

\item[(ii)] ils sont nuls pour $|i| \geq s-l+1$;

\item[(iii)] ils sont en général mixtes de poids $(s-l)(g+1)-2(k-1)$ pour $1 \leq k \leq
s-l+1$ vérifiant
$s-l-2(k-1) \leq i \leq s-l-(k-1)$;

\item[(iv)] soit $l-s \leq i_0 < 0$, le plus petit indice $i$ tel que la partie de
poids
$(s-l)(g-1)$ de
$$\lim_{\genfrac..{0pt}{1}{\to}{I}} H^i(M_{I,s_o}^{\geq lg},j^{\geq lg}_!
HT_{\r_\oo}(g,l,\pi_o,\Pi_l,I)[(s-l)g])
[(s-l)g])[\Pi^{\oo,o}]$$ soit non nulle\footnote{Si un tel $i_0$, n'existe pas l'énoncé
est vide.}. Cette
dernière est alors égale à
$$\Pi_l \overrightarrow{\times}
[\overleftarrow{i_0-l+s-1},\overrightarrow{-i_0}]_{\pi_o} \otimes
(\Xi^{\frac{(s-l)(g-1)}{2}} \otimes
\bigoplus_{\xi \in \AF(\pi_o)} \xi^{-1}).$$
\end{itemize}
\end{lemm}

\begin{proof} Le point (i) correspond à la proposition (\ref{prop-poids}) qui découle
directement de l'action de
$GL_{lg}(F_o)$ sur la strate via $\val(\det)$ et du fait que les strates non
supersingulières sont induites.

Pour les points (ii)-(iii), on considère les suites exactes
\begin{equation} \label{secfpi}
0 \to P_{\pi_o,l,i+1}(\Pi_l) \longto P_{\pi_o,l,i}(\Pi_l) \longto
\PC_-(g,l+i+1,\pi_o,l,\Pi_l) \to 0
\end{equation}
et on reprend les notations simplifiées de la preuve de la proposition (\ref{prop-not})
en notant
$P_{\pi_o,l,-1,\r_\oo}(\Pi_l):=j^{\geq lg}_! HT_{\r_\oo}(g,l,\pi_o,\Pi_l)[d-lg]$. Par
facilité on notera
$P_{\pi_o,l,i}(\Pi_l)$ pour $P_{\pi_o,l,i}(\Pi_l) \otimes \LC_{\r_\oo}$. Le résultat
découle alors
directement du cas $i=-1$ dans le lemme suivant.

\begin{lemm} \label{lem-utile}
Pour $-1 \leq i \leq s-l-1$, les ${\DS \lim_{\genfrac..{0pt}{1}{\to}{I}}}
H^j(M_{I,s_o}^{\geq
(l+i+1)g},P_{\pi_o,l,i}(\Pi_l) \otimes \LC_{\rho_\oo})[\Pi^{\oo,o}]$ vérifient les
points suivants:
\begin{itemize}

\item[(1)] ils sont nuls pour $|j| \geq s-l-i$;

\item[(2)] pour $|j| < s-l-i$, ils sont mixtes de poids $(s-l)(g+1)-2(k+i)$ pour $1 \leq
k \leq s-l-i$ vérifiant
$s-l-i-1-2(k-1) \leq j \leq s-l-i-1-(k-1)$;

\end{itemize}
\end{lemm}

\begin{proof} (1)-(2) On raisonne par récurrence descendante, le cas $i=s-l-1$ étant
évident car
$P_{\pi_o,l,s-l-1}(\Pi_l)$ est le faisceau concentré aux points supersinguliers
$$\FC(g,s,\pi_o) \otimes \Pi_l(-(s-l)/2) \times [\overleftarrow{s-l-1}]_{\pi_o(l/2)}
\otimes |\cl|^{-\frac{(s-l)(g-1)}{2}}$$
Supposons donc le résultat acquis jusqu'au rang $i+1$ et traitons le cas de
$P_{\pi_o,l,i}(\Pi_l)$. On
considère la suite exacte longue de cohomologie associée à la suite exacte courte
(\ref{secfpi}). D'après
l'hypothèse de récurrence les $H^j(P_{\pi_o,l,i+1}(\Pi_l))$ sont nuls pour $|j| \geq
s-l-i-1$ et d'après la
proposition (\ref{prop-not}), les $H^j(\PC_-(g,l+i+1,\pi_o,l,\Pi_l))$ sont nuls pour
$|j| \geq s-l-i$, d'où
le point (1). En ce qui concerne le point (2), la suite exacte longue en question
s'écrit alors
\begin{multline} \label{sel1}
0 \to H^{-(s-l-i-1)}(P_{\pi_o,l,i}(\Pi_l)) \longto H^{-(s-l-i-1)}(j^{\geq
(l+i+1)g}_{!*}) \longto \\
H^{-(s-l-i-2)}(P_{\pi_o,l,i+1}(\Pi_l))
\longto H^{-(s-l-i-2)} P_{\pi_o,l,i}(\Pi_l) \longto 0 \cdots \\
\cdots 0 \to H^{(s-l-i-1)-2r}(P_{\pi_o,l,i+1}(\Pi_l)) \longto
H^{(s-l-i-1)-2r}(P_{\pi_o,l,i}(\Pi_l)) \longto
\\ H^{(s-l-i-1)-2r}(j^{\geq (l+i+1)g}_{!*}) \longto
 H^{(s-l-i-1)-2r+1}(P_{\pi_o,l,i+1}(\Pi_l)) \\ \longto
H^{(s-l-i-1)-2r+1}(P_{\pi_o,l,i}(\Pi_l)) \to 0 \\
\cdots \\
0 \to H^{(s-l-i-1)-2}(P_{\pi_o,l,i+1}(\Pi_l)) \longto
H^{(s-l-i-1)-2}(P_{\pi_o,l,i}(\Pi_l)) \longto \\
H^{(s-l-i-1)-2}(j^{\geq (l+i+1)g}_{!*})  \longto
H^{s-l-i-2}(P_{\pi_o,l,i+1}(\Pi_l)) \longto H^{s-l-i-2}(P_{\pi_o,l,i}(\Pi_l)) \to 0 \\
0 \to H^{s-l-i-1}(P_{\pi_o,l,i}(\Pi_l)) \longto H^{s-l-i-1}(j^{\geq (l+i+1)g}_{!*}) \to
0
\end{multline}
On rappelle que $H^j(j^{\geq (l+i+1)g}_{!*})$ est pur de poids $(s-l)g-(i+1)+j$ de sorte
qu'en utilisant
l'hypothèse de récurrence, les $H^j(P_{\pi_o,l,i}(\Pi_l))$ sont de poids
$(s-l)(g+1)-2(i+1+k-1)$ avec $1 \leq
k \leq s-l-i$.

Soit alors $j$ de la forme $s-l-1-i-(2r+1)$ avec $0 \leq 2r+1 \leq 2(s-l-1-i)$; la suite
exacte longue
ci-dessus montre alors que les poids de $H^j(P_{\pi_o,l,i}(\Pi_l))$ sont ceux de
$H^j(P_{\pi_o,l,i+1}(\Pi_l))$, i.e. $(s-l)(g+1)-2(i+1+k)$ pour $1 \leq k \leq s-l-i-1$
vérifiant $-2(k-1)
\leq j-(s-l-2-i) \leq -(k-1)$. Le changement de variable $k'=k+1$ donne alors le
résultat, i.e.
$H^j(P_{\pi_o,l,i}(\Pi_l))$ est de poids $(s-l)(g+1)-2(i+k')$ avec $1 \leq k' \leq
s-l-i$ vérifiant
$-2(k'-1)+1 \leq j-(s-l-1-i) \leq -(k'-1)$ soit ce qui est prévu car $j-(s-l-1-i)$ est
impair\footnote{Le cas
$j-(s-l-1-i)=-2(k'-1)$ n'est pas à considérer.}.

Pour $j$ de la forme $s-l-1-i-2r$ avec $0 \leq 2r \leq 2(s-l-1-i)$, la suite exacte
longue précédente montre
que les poids de $H^j(P_{\pi_o,l,i}(\Pi_l))$ sont, à priori, ceux de
$H^j(P_{\pi_o,l,i+1}(\Pi_l))$ ainsi que
celui de $H^j(j^{\geq (l+i+1)g}_{!*})$ soit $(s-l)(g+1)-2(i+1+r)$. D'après l'hypothèse
de récurrence,
$H^j(P_{\pi_o,l,i}(\Pi_l))$ est de poids $(s-l)(g+1)-2(i+1+k)$ avec $1 \leq k \leq
s-l-i-1$ vérifiant
$-2(k-1) \leq j-(s-l-2-i) \leq -(k-1)$, de sorte que $H^j(P_{\pi_o,l,i}(\Pi_l))$ est de
poids
$(s-l)(g+1)-2(i+k')$ avec $1 \leq k' \leq s-l-i$ vérifiant $-2(k'-1)+1 \leq j-(s-l-1-i)
\leq -(k'-1)$ soit ce
qui est prévu car le cas $j-(s-l-1-i)=-2(k'-1)$ est justement donné par le poids de
$H^j(j^{\geq
(l+i+1)g}_{!*})$ soit $(s-l)(g+1)-2(i+1+r)$.

\end{proof}

\noindent \textit{Suite de la preuve du lemme (\ref{lem-rj-combi})}: (iv) Dans la suite
on ne considère que
les parties de poids $(s-l)(g-1)$. D'après le lemme (\ref{lem-utile}),
$H^{l-s}(P_{\pi_o,l,0}(\Pi_l))$ est
nul de sorte que la suite exacte longue (\ref{sel1}) s'écrit
\begin{multline*}
0 \to H^{l-s}(j^{\geq lg}_!) \longto \Pi_l \overrightarrow{\times}
[\overrightarrow{s-l-1}]_{\pi_o} \otimes
(\Xi^{(s-l)(g-1)/2} \otimes \bigoplus_{\xi \in \AF(\pi_o)} \xi^{-1})\\
\longto H^{1+l-s}(P_{\pi_o,l,0}(\Pi_l)) \longto H^{1+l-s}(j^{\geq lg}_!) \to 0 \to
\cdots
\end{multline*}
Le cas $i_0=l-s$ découle alors du fait que les strates sont induites, i.e. si
$H^{l-s}(j^{\geq lg}_!)$ est un
sous-espace de $\Pi_l \overrightarrow{\times} [\overrightarrow{s-l-1}]_{\pi_o} \otimes
(\Xi^{\frac{(s-l)(g-1)}{2}} \otimes \bigoplus_{\xi \in \AF(\pi_o)} \xi^{-1})$ alors il
est égal à tout
l'espace.

Pour $i_0=1+l-s$, (\ref{sel1}) s'écrit
\begin{multline*}
0 \to \Pi_l \overrightarrow{\times} [\overrightarrow{s-l-1}]_{\pi_o} \otimes
(\Xi^{\frac{(s-l)(g-1)}{2}}
\otimes \bigoplus_{\xi \in \AF(\pi_o)} \xi^{-1}) \longto \\
H^{1+l-s}(P_{\pi_o,l,0}(\Pi_l)) \longto H^{1+l-s}(j^{\geq lg}_!) \to 0
\end{multline*}
tandis que celle associée à la suite exacte courte
$$0 \to P_{\pi_o,l,1}(\Pi_l) \longto P_{\pi_o,l,0}(\Pi_l) \longto
\PC_-(g,l+1,\pi_o,l,\Pi_l) \to 0$$
s'écrit
$$0 \to H^{1+l-s}(P_{\pi_o,l,0}(\Pi_l)) \longto \Pi_l \overrightarrow{\times}
[\overrightarrow{0}]_{\pi_o}
\overrightarrow{\times} [\overrightarrow{s-l-2}]_{\pi_o} \otimes
(\Xi^{\frac{(s-l)(g-1)}{2}} \otimes
\bigoplus_{\xi \in \AF(\pi_o)} \xi^{-1}) \to \cdots$$ Ainsi si $H^{1+l-s}(j^{\geq
lg}_!)$ est non nul, alors
$H^{1+l-s}(P_{\pi_o,l,0}(\Pi_l))$ contient strictement
$$\Pi_l \overrightarrow{\times}
\overrightarrow{(s-l)}_{\pi_o} \otimes (\Xi^{\frac{(s-l)(g-1)}{2}} \otimes
\bigoplus_{\xi \in \AF(\pi_o)}
\xi^{-1})$$ et étant de la forme $\Pi_l \overrightarrow{\times} \pi \otimes
(\Xi^{\frac{(s-l)(g-1)}{2}}
\otimes \bigoplus_{\xi \in \AF(\pi_o)} \xi^{-1})$ ainsi qu'un sous-espace de
$$\Pi_l \overrightarrow{\times} [\overrightarrow{0}]_{\pi_o} \overrightarrow{\times}
[\overrightarrow{s-l-2}]_{\pi_o} \otimes (\Xi^{\frac{(s-l)(g-1)}{2}} \otimes
\bigoplus_{\xi \in \AF(\pi_o)}
\xi^{-1}),$$ on en déduit qu'il est égal à ce dernier de sorte que $H^{1+l-s}(j^{\geq
lg}_!)$ est isomorphe à
$$\Pi_l \overrightarrow{\times}
[\overleftarrow{1},\overrightarrow{s-l-1}]_{\pi_o} \otimes (\Xi^{\frac{(s-l)(g-1)}{2}}
\otimes \bigoplus_{\xi
\in \AF(\pi_o)} \xi^{-1}),$$ d'où le résultat.

Par ailleurs on remarque de la même façon que si $i_0 >1+l-s$, la partie de poids
$(s-l)(g-1)$ de
$H^{1+l-s}(P_{\pi_o,l,0}(\Pi_l))$ est alors égale à
$$\Pi_l \overrightarrow{\times} [\overrightarrow{s-l-1}]_{\pi_o} \otimes
(\Xi^{\frac{(s-l)(g-1)}{2}} \otimes
\bigoplus_{\xi \in \AF(\pi_o)} \xi^{-1}).$$

Supposons alors $i_0 \geq 2+l-s$. La suite exacte longue (\ref{sel1}) donne l'égalité
pour tout $i \geq 2$,
des parties de poids $(s-l)(g-1)$ de $H^{l-s+i}(P_{\pi_o,l,0}(\Pi_l))$ et de
$H^{l-s+i}(j^{\geq lg}_!)$. On
va montrer que, pour tout $2 \leq i \leq i_0-(l-s)$, la partie de poids $(s-l)(g-1)$ de
$H^{l-s+r}(P_{\pi_o,l,i-2}(\Pi_l))$ est nulle pour $i \leq r < i_0-l+s$ et que celle de
$H^{i+l-s-1}(P_{\pi_o,l,i-2}(\Pi_l))$ est égale à $\Pi_l \overrightarrow{\times}
[\overleftarrow{i-1},\overrightarrow{s-l-i}]_{\pi_o} \otimes (\Xi^{\frac{(s-l)(g-1)}{2}}
\otimes
\bigoplus_{\xi \in \AF(\pi_o)} \xi^{-1})$. D'après ce que l'on vient de voir, c'est vrai
pour $i=2$.
Supposons donc le résultat acquis jusqu'au rang $i$ et traitons le cas de $i+1$. La
suite exacte longue de
cohomologie associée à la suite exacte courte de faisceaux pervers
$$0 \to P_{\pi_o,l,i-1}(\Pi_l) \longto P_{\pi_o,l,i-2}(\Pi_l) \longto
\PC_-(g,l+i-1,\pi_o,l,\Pi_l) \to 0$$
s'écrit
\begin{multline} \label{sel2}
0 \to H^{i-1+l-s}(P_{\pi_o,l,i-2}(\Pi_l)) \longto \Pi_l \overrightarrow{\times}
[\overleftarrow{i-2}]_{\pi_o}
\overrightarrow{\times} [\overrightarrow{s-l-i}]_{\pi_o} \otimes
(\Xi^{\frac{(s-l)(g-1)}{2}} \otimes
\bigoplus_{\xi \in \AF(\pi_o)} \xi^{-1})
\\ \longto H^{i+l-s}(P_{\pi_o,l,i-1}(\Pi_l)) \longto H^{i+l-s}(P_{\pi_o,l,i-2}(\Pi_l))
\to 0 \cdots
\end{multline}
ainsi que l'égalité des parties de poids $(s-l)(g-1)$ des espaces
$H^{i+l-s+r}(P_{\pi_o,l,i-1}(\Pi_l))$ et
$H^{i+l-s+r}(P_{\pi_o,l,i-2}(\Pi_l))$ pour tout $r>0$. La nullité de la partie de poids
$(s-l)(g-1)$ de
$H^{i+l-s}(P_{\pi_o,l,i-2}(\Pi_l))$ et l'isomorphisme
$$H^{i+l-s-1}(P_{\pi_o,l,i-2}(\Pi_l)) \simeq \Pi_l \overrightarrow{\times}
[\overleftarrow{i-2},
\overrightarrow{s-l-i+1}]_{\pi_o} \otimes (\Xi^{\frac{(s-l)(g-1)}{2}} \otimes
\bigoplus_{\xi \in \AF(\pi_o)}
\xi^{-1})$$ donne le cas $i+1$.

Traitons désormais le cas $i=i_0-l+s$. On pose $i_1=i_0-l+s-2$ et considérons alors la
suite exacte longue de
cohomologie associée à
$$0 \to P_{\pi_o,l,i_1+1}(\Pi_l) \longto P_{\pi_o,l,i_1}(\Pi_l) \longto
\PC_-(g,l+i_1+1,\pi_o,l,\Pi_l) \to 0$$
qui s'écrit
\begin{multline} \label{sel3}
0 \to H^{i_0-1}(P_{\pi_o,l,i_1}(\Pi_l)) \longto \Pi_l \overrightarrow{\times}
[\overleftarrow{i_1}]_{\pi_o}
\overrightarrow{\times} [\overrightarrow{s-l-i_1-2}]_{\pi_o} \otimes
(\Xi^{\frac{(s-l)(g-1)}{2}} \otimes
\bigoplus_{\xi \in \AF(\pi_o)} \xi^{-1}) \\ \longto H^{i_0}(P_{\pi_o,l,i_1+1}(\Pi_l))
\longto
H^{i_0}(P_{\pi_o,l,i_1}(\Pi_l)) \to 0 \cdots
\end{multline}
avec $H^{i_0}(P_{\pi_o,l,i_1}(\Pi_l))=H^{i_0}(j^{\geq lg}_!)$ non nul par hypothèse. La
suite exacte longue
associée à
$$0 \to P_{\pi_o,l,i_1+2}(\Pi_l) \longto P_{\pi_o,l,i_1+1}(\Pi_l) \longto
\PC_-(g,i_1+2,\pi_o,l,\Pi_l) \to 0$$
s'écrit
\begin{multline*}
0 \to H^{i_0}(P_{\pi_o,l,i_1+1}(\Pi_l)) \longto \Pi_l \overrightarrow{\times}
[\overleftarrow{i_1+1}]_{\pi_o}
\overrightarrow{\times} [\overrightarrow{s-l-i_1-3}]_{\pi_o} \\ \otimes
(\Xi^{\frac{(s-l)(g-1)}{2}} \otimes
\bigoplus_{\xi \in \AF(\pi_o)} \xi^{-1}) \longto \cdots
\end{multline*}
Ainsi si on veut que $H^{i_0}(P_{\pi_o,l,i_1}(\Pi_l))$ soit non nul, il faut que
$H^{i_0}(P_{\pi_o,l,i_1+1}(\Pi_l))$ soit égal à $\Pi_l \overrightarrow{\times}
[\overleftarrow{i_1+1}]_{\pi_o} \overrightarrow{\times}
[\overrightarrow{s-l-i_1-3}]_{\pi_o} \otimes
(\Xi^{\frac{(s-l)(g-1)}{2}} \otimes \bigoplus_{\xi \in \AF(\pi_o)} \xi^{-1})$ et donc
$$H^{i_0}(P_{\pi_o,l,i_1+1}(\Pi_l)) \simeq \Pi_l \overrightarrow{\times}
[\overleftarrow{i_1+2},\overrightarrow{s-l-i_1-3}]_{\pi_o} \otimes
(\Xi^{\frac{(s-l)(g-1)}{2}} \otimes
\bigoplus_{\xi \in \AF(\pi_o)} \xi^{-1})$$ d'où le résultat.

\end{proof}

\noindent \textit{Retour sur la preuve de la proposition (\ref{prop-hic-poids}):}

\rem Si on savait que les strates $M_{I,s_o}^{=h}$ étaient affines, la proposition
(\ref{prop-hic-poids})
découlerait directement du point (iii) du lemme (\ref{lem-rj-combi}) pour $k=s-l+1$ car
seul
$H^{d-lg}_c(M_{I,s_o}^{=h},\FC(g,l,\pi_o,I) \otimes \LC_{\r_\oo})[\Pi^{\oo,o}]$ pourrait
être de poids $l-s$
ce qui serait contradictoire avec (\ref{somme-alternee}). Ne disposant pas de l'affinité
des strates
$M_{I,s_o}^{=h}$, on entre plus précisément dans la combinatoire.

On raisonne par récurrence sur $l$ de $1$ à $s_g$ puis pour $l$ fixé, par récurrence sur
$i$ de $l-s$ à $0$.
L'initialisation de la récurrence pour $l=1$ se traite comme le passage de $l$ à $l+1$;
on suppose alors le
résultat acquis pour tout $1 \leq l' < l$ \footnote{Pour $l=1$, l'hypothèse de
récurrence est vide.}.

- Pour $i=l-s$, si l'espace en question était non nul, on aurait d'après le point (iv)
du lemme
(\ref{lem-rj-combi}),
\begin{multline*}
\lim_{\genfrac..{0pt}{1}{\to}{I}}
H^{(s-l)(g-1)}_c(M_{I,s_o},HT_{\rho_\oo}(g,l,\pi_o,\Pi_l))[\Pi^{\oo,o}] =
\\ \Pi_l \overrightarrow{\times} [\overrightarrow{s-l-1}]_{\pi_o} \otimes
(\Xi^{\frac{(s-l)(g-1)}{2}} \otimes
\bigoplus_{\xi \in \AF(\pi_o)} \xi^{-1}).
\end{multline*}
On considère alors la suite spectrale (\ref{sss}) pour $i=l(g-1)$. On remarque que
toutes les parties de
poids $s(g-1)$ de $E_1^{p,q;l(g-1)}[\Pi^{\oo,o}]$ sont nulles pour $p-1 \neq lg$; en
effet pour $p-1 < lg$
cela découle de l'hypothèse de récurrence et pour $p-1 > lg$ du fait que
$l'(g-1)>l(g-1)$ pour $p-1=l'g$. On
obtient alors
\begin{multline*}
\lim_{\genfrac..{0pt}{1}{\to}{I}} H^{(s-l)(g-1)}(M_{I,\bar
s_o},R^{l(g-1)}\Psi_{\eta_o,I,\pi_o}(\LC_{\r_\oo}))[\Pi^{\oo,o}]
\simeq \\
[\overleftarrow{l-1}]_{\pi_o} \overrightarrow{\times} [\overrightarrow{s-l-1}]_{\pi_o}
\otimes L_g(\pi_o)
(-\frac{s(g-1)}{2})
\end{multline*}
On considère alors la suite spectrale (\ref{ssce}), de sorte que
$E_{2,\r_\oo}^{(s-l)(g-1),l(g-1)}[\Pi^{\oo,o}]$ est isomorphe à l'espace ci-dessus. Or
comme d'après
l'hypothèse de récurrence toutes les parties de poids $s(g-1)$ des
$E_{2,\r_\oo}^{p,q}[\Pi^{\oo,o}]$ pour
$q<l(g-1)$ sont nulles et que d'après le lemme (\ref{lem-rj-combi}) il en est de même
pour $p+q < d-s$, on en
déduit que $E_{\oo,\r_\oo}^{d-s}[\Pi^{\oo,o}]$ admettrait
$[\overleftarrow{l-1}]_{\pi_o}
\overrightarrow{\times} [\overrightarrow{s-l-1}]_{\pi_o} \otimes L_g(\pi_o)
(-\frac{s(g-1)}{2})$ comme
sous-espace ce qui n'est pas d'après la proposition (\ref{prop-ssce-poids}).

- Supposons donc le résultat vérifié pour tout $l-s \leq i=l-s+\d \leq i_0 < -1$ et
traitons le cas de $i_0$
D'après le point (iv) du lemme (\ref{lem-rj-combi}), on obtiendrait
$$[\overleftarrow{l-1}]_{\pi_o} \overrightarrow{\times}
[\overleftarrow{\d},\overrightarrow{s-l-1-\d}]_{\pi_o} \otimes
(\Xi^{\frac{(s-l)(g-1)}{2}} \otimes
\bigoplus_{\xi \in \AF(\pi_o)} \xi^{-1})$$ Comme précédemment, la suite spectrale
(\ref{sss}) pour
$i=l(g-1)$, donne que la partie de poids $s(g-1)$ de
$E_{2,\r_\oo}^{(s-l)(g-1)+\d,l(g-1)}[\Pi^{\oo,o}]$ dans
la suite spectrale (\ref{ssce}) est isomorphe à
$$[\overleftarrow{l-1}]_{\pi_o} \overrightarrow{\times}
[\overleftarrow{\d},\overrightarrow{s-l-1-\d}]_{\pi_o} \otimes L_g(\pi_o)
(-\frac{s(g-1)}{2}).$$ On remarque
alors que $E_{\oo,\r_\oo}^{d-s+\d}[\Pi^{\oo,o}]$ admettrait
$[\overleftarrow{l+\d},\overrightarrow{s-l-1-\d}]_{\pi_o} \otimes L_g(\pi_o)
(-\frac{s(g-1)}{2})$ comme
sous-espace car ce dernier n'apparaît pas dans les $E_{2,\r_\oo}^{p,q}[\Pi^{\oo,o}]$
pour $q<l(g-1)$ d'après
l'hypothèse de récurrence, ni pour $q=l'(g-1)>l(g-1)$ d'après le lemme
(\ref{lem-rj-combi}) (i) ainsi que le
corollaire (\ref{coro-1}).

\end{proof}

\subsection{Preuve de \ref{cas-m}}

Il s'agit donc de prouver la proposition (\ref{cas-m}). En étudiant les parties de poids
$s(g-1)$ des
composantes $\Pi^{\oo,o}$ isotypiques pour $\Pi$ irréductible automorphe tel que $\Pi_o
\simeq
\speh_s(\pi_o)$, on s'est ramené d'après le paragraphe précédent, à une situation
similaire à celle de
\cite{boy}, i.e. dans la suite spectrale des cycles évanescents, seuls les points
supersinguliers
contribuent. La fin de la preuve procède alors exactement comme dans loc. cit.

Précisément, la partie de poids $s(g-1)$ de
$\widetilde{\UC^{d,d-s}_{F_o}}(\JL^{-1}([\overleftarrow{s-1}]_{\pi_o}))$, d'après le
corollaire
(\ref{coro-1}), est un constituant de $[\overrightarrow{s-2}]_{\pi_o}
\overrightarrow{\times}
[\overrightarrow{0}]_{\pi_o} \otimes L_g(\pi_o) (-\frac{s(g-1)}{2})$. On considère alors
la suite spectrale
(\ref{sss2}) pour $i=d-s$, où, d'après le corollaire (\ref{coro-hic}), tous les
$E_{1,\r_\oo}^{p,q;d-s}$ sont
nuls pour $p+q \neq 0$ ou $p-1 \neq d$ de sorte que le terme
$E_{2,\r_\oo}^{0,d-s}[\Pi^{\oo,o}]$ de
(\ref{ssce2}), et donc $E_{\oo,\r_\oo}^{0,d-s}[\Pi^{\oo,o}]$ d'après
(\ref{prop-hic-poids}), est égal à cet
espace qui est donc, d'après (\ref{prop-ssce-poids}), égal à
$[\overrightarrow{s-1}]_{\pi_o} \otimes
L_g(\pi_o) (-\frac{s(g-1)}{2})$.

On suppose avoir montré par récurrence que pour tout $0 \leq r < r_0< s-1$, les parties de poids $s(g-1)$ de
$\widetilde{\UC_{F_o}^{d,d-s+r}}(\JL^{-1}([\overleftarrow{s-1}]_{\pi_o}))$ sont comme prévues, i.e. nulles
pour $r \neq 0$, et égales à $[\overrightarrow{s-1}]_{\pi_o} \otimes L_g(\pi_o) (-\frac{s(g-1)}{2})$ pour
$r=0$. D'après le corollaire (\ref{coro-1}), la partie de poids $s(g-1)$ de
$\widetilde{\UC_{F_o}^{d,d-s+r_0}}(\JL^{-1}([\overleftarrow{s-1}]_{\pi_o}))$ est un constituant de
$[\overrightarrow{s-2-r_0}]_{\pi_o} \overrightarrow{\times} [\overleftarrow{r_0}]_{\pi_o} \otimes L_g(\pi_o)
(-\frac{s(g-1)}{2})$. Or comme la partie de poids $s(g-1)$ de $\sum_{i=0}^{s-1} (-1)^i
\widetilde{\UC_{F_o}^{d,d-s+i}}(\JL^{-1}([\overleftarrow{s-1}]_{\pi_o}))$ est égale à
$[\overrightarrow{s-1}]_{\pi_o} \otimes L_g(\pi_o) (-\frac{s(g-1)}{2})$, en remarquant, d'après le corollaire
(\ref{coro-1}), que
$$[\overrightarrow{s-1-r_0},\overleftarrow{r_0+1}]_{\pi_o} \otimes L_g(\pi_o) (-\frac{s(g-1)}{2})$$
ne peut pas être un constituant de
$\widetilde{\UC_{F_o}^{d,d-s+r}}(\JL^{-1}([\overleftarrow{s-1}]_{\pi_o}))$ pour $r>r_0$, on en déduit que si
la partie de poids $s(g-1)$ de $\widetilde{\UC_{F_o}^{d,d-s+r_0}}(\JL^{-1}([\overleftarrow{s-1}]_{\pi_o}))$
est non nulle, elle est alors égale à $[\overrightarrow{s-2-r_0},\overleftarrow{r_0+2}]_{\pi_o} \otimes
L_g(\pi_o) (-\frac{s(g-1)}{2})$. Comme précédemment ce dernier espace est aussi égal à la partie de poids
$s(g-1)$ de $E_{\oo,\r_\oo}^{d-s+r_0}[L_g(\pi_o)][\Pi^{\oo,o}]$ ce qui n'est pas d'après
(\ref{prop-ssce-poids}).

Finalement le cas $r_0=s-1$ est donné en utilisant que la partie de poids $s(g-1)$ de
$\sum_{i=0}^{s-1}
(-1)^i \widetilde{\UC_{F_o}^{d,d-s+i}}(\JL^{-1}([\overleftarrow{s-1}]_{\pi_o}))$ est
égale à
$[\overrightarrow{s-1}]_{\pi_o} \otimes L_g(\pi_o) (-\frac{s(g-1)}{2})$.

\end{proof}

\begin{rema} \label{rema-N2}
Comme dans la remarque (\ref{rema-N}), on aurait pu traiter tous les poids des
$\widetilde{\UC_{F_o}^{d,i}}$.
Il suffit pour cela de montrer un analogue du point (iv) de la proposition
(\ref{prop-hic-poids}) pour tous
les poids, ce qui ne pose aucune difficulté sauf qu'on ne se retrouve pas dans une
situation aussi favorable
où seuls les points supersinguliers contribuent.
\end{rema}


%% file: figg1.tex
\begin{figure}[!h]
\setlength{\unitlength}{.17cm} \centering
\begin{picture}(50,53)(-5,-50)
\linethickness{.1pt}

\put(0,0){\circle{.4}} \put(10,0){\makebox(0,0){$\overleftarrow{l-1} \overleftarrow{\times}
\overleftarrow{s-l-1}$}}

\put(0,-2){\vector(0,1){2}} \put(0,-2.5){\makebox(0,0){$p=-q=-(s-l)g$}} \put(0,0){\line(1,-1){14}}

\put(7,-7){\circle{.1}} \put(7,-7){\line(0,1){1}} \put(7,-6){\circle{.4}} \put(7,-6){\line(1,-1){7}}
\put(7,-9){\vector(0,1){2}} \put(5,-9){\makebox(0,0){$p=-q=-(s-l-1)g$}}

\put(17,-6){\makebox(0,0){$\overleftarrow{l-1} \overleftarrow{\times} \overrightarrow{0}
\overleftarrow{\times} \overleftarrow{s-l-2}$}} \put(14,-14){\circle{.1}} \put(14,-14){\line(0,1){2}}
\put(12,-14){\vector(1,0){2}} \put(2,-14){\makebox(0,0){$p=-q=-(s-l-2)g$}} \dashline{1}(14,-14)(28,-28)
\put(14,-13){\circle{.1}} \dashline{1}(14,-13)(28,-27) \put(14,-12){\circle{.4}}

\put(24,-12){\makebox(0,0){$\overleftarrow{l-1} \overleftarrow{\times} \overrightarrow{1}
\overleftarrow{\times} \overleftarrow{s-l-3}$}} \dashline{1}(14,-12)(28,-26) \put(28,-28){\line(0,1){2}}
\put(26,-28){\vector(1,0){2}} \put(20,-28){\makebox(0,0){$p=-q=-3g$}} \dashline{1}(28,-26)(28,-24)
\multiput(28,-28)(0,1){3}{\circle{.1}} \multiput(28,-28)(0,1){3}{\line(1,-1){21}} \put(28,-24){\circle{.4}}
\put(28,-24){\line(1,-1){21}}

\put(38,-24){\makebox(0,0){$\overleftarrow{l-1} \overleftarrow{\times} \overrightarrow{s-l-4}
\overleftarrow{\times} \overleftarrow{2}$}} \put(35,-30){\circle{.4}} \put(35,-30){\line(1,-1){14}}
\put(35,-35){\line(0,1){2}} \put(33,-35){\vector(1,0){2}} \put(28,-35){\makebox(0,0){$p=q=-2g$}}
\multiput(35,-35)(0,1){3}{\circle{.1}} \put(35,-31){\circle{.1}} \put(35,-31){\line(0,1){1}}

\put(45,-30){\makebox(0,0){$\overleftarrow{l-1} \overleftarrow{\times} \overrightarrow{s-l-3}
\overleftarrow{\times} \overleftarrow{1}$}} \put(42,-36){\circle{.4}} \put(42,-36){\line(1,-1){7}}
\put(42,-42){\line(0,1){2}} \multiput(42,-42)(0,1){3}{\circle{.1}} \put(40,-42){\vector(1,0){2}}
\put(35,-42){\makebox(0,0){$p=-q=-g$}} \put(42,-38){\line(0,1){2}} \multiput(42,-38)(0,1){2}{\circle{.1}}

\put(52,-36){\makebox(0,0){$\overleftarrow{l-1} \overleftarrow{\times} \overrightarrow{s-l-2}
\overleftarrow{\times} \overleftarrow{0}$}} \put(49,-49){\line(0,1){2}} \put(47,-49){\vector(1,0){2}}
\put(42,-49){\makebox(0,0){$p=q=0$}} \dashline{1}(49,-47)(49,-45) \put(49,-45){\line(0,1){3}}
\multiput(49,-48)(0,1){2}{\circle{.4}} \multiput(49,-45)(0,1){4}{\circle{.4}}

\end{picture}
\caption{\label{fig-dualite2} parties de poids $(s-l)(g+1)$ de la fibre aux points supersinguliers des
$E_2^{p,q}$ de la suite spectrale (\ref{ss-dualite})}
\end{figure}

%% file: figu1.tex
\begin{figure}[!t]
\setlength{\unitlength}{.7cm}
\begin{picture}(7,12)(-5,-2)
\linethickness{.1pt}

\put(0,2){\vector(0,1){5}}

\dashline{3}(0,0)(0,2)

\put(0,0){\vector(1,0){6}}

\put(0,2){\line(1,0){5}}

\put(7,2){\makebox(0,0){$j=g-1$}}

\put(0,5){\circle{.3}}

\put(-2,5){\makebox(0,0){$\genfrac{}{}{0pt}{}{[\overleftarrow{0}]_{\pi_o} \overrightarrow{\times}
[\overrightarrow{0}]_{\pi_o}}{\otimes L_g(\pi_o)|-|^{-2(g-1)/2}}$}}

\put(0,6){\circle{.3}}

\put(2,6){\vector(-1,0){2}}

\put(2,6){$i+j=d-1$}

\put(2,8){\vector(-1,-1){2}}

\put(2,9){\makebox(0,0){$\genfrac{}{}{0pt}{}{[\overleftarrow{1}]_{\pi_o} \otimes L_g(\pi_o)}{(
|-|^{-2(g+1)/2} \oplus |-|^{-2(g-1)/2})}$}}

\put(0,6){\line(1,-1){4}}

\put(4,2){\circle{.3}}

\put(5,3){\vector(-1,-1){1}}

\put(5,4){\makebox(0,0){$\genfrac{}{}{0pt}{}{[\overleftarrow{0}]_{\pi_o} \overrightarrow{\times}
[\overrightarrow{0}]_{\pi_o}}{\otimes L_g(\pi_o)|-|^{-2(g-1)/2}}$}}

\put(0,5){\vector(4,-3){4}} \dashline{3}(4,0)(4,2)

\put(4,-1){\vector(0,1){1}}

\put(4,-1.5){\makebox(0,0){$i=d-g$}}

\put(0,-.5){\makebox(0,0){$(0,0)$}}

\multiput(0,2)(1,0){5}{\circle*{.1}}

\multiput(0,2)(0,1){5}{\circle*{.1}}

\end{picture}
\caption{\label{figu1} $E_2^{p,q}[\Pi^{\oo,o}]$ de (\ref{ssce}) pour $\Pi_o \simeq \st_s(\pi_o)$ compatible
avec l'aboutissement}
\end{figure}

%% file: figu2.tex
\begin{figure}[!t]
\setlength{\unitlength}{.7cm}
\begin{picture}(7,12)(-4,-2)
\linethickness{.1pt}

\put(0,2){\vector(0,1){5}}
\dashline{3}(0,0)(0,2)

\put(0,0){\vector(1,0){6}}

\put(0,2){\line(1,0){6}}

\put(8,2){\makebox(0,0){$j=g-1$}}

\put(0,5){\circle{.3}}

\put(-2,5){\makebox(0,0){$\genfrac{}{}{0pt}{}{[\overleftarrow{0}]_{\pi_o} \overrightarrow{\times}
[\overrightarrow{0}]_{\pi_o}}{\otimes L_g(\pi_o) |-|^{-2(g-1)/2}}$}}

\put(0,6){\circle{.3}}

\put(2,6){\vector(-1,0){2}}

\put(2,6){$i+j=d-1$}

\put(2,8){\vector(-1,-1){2}}

\put(2,9){\makebox(0,0){$\genfrac{}{}{0pt}{}{[\overleftarrow{1}]_{\pi_o} \otimes L_g(\pi_o)}{(|-|^{-2(g-1)/2}
\oplus |-|^{-2(g+1)/2})}$}}

\put(0,6){\line(1,-1){4}}

\put(5,2){\circle{.3}}

\put(6,3){\vector(-1,-1){1}}

\put(6,4){\makebox(0,0){$\genfrac{}{}{0pt}{}{[\overleftarrow{0}]_{\pi_o} \overleftarrow{\times}
[\overrightarrow{0}]_{\pi_o}}{\otimes L_g(\pi_o) |-|^{-2(g+1)/2}}$}}

\put(0,6){\line(5,-4){5}}

\put(5,2){\vector(1,-1){0}} \dashline{3}(5,0)(5,2)

\put(5,-1){\vector(0,1){1}}

\put(5,-1.5){\makebox(0,0){$i=d-g+1$}}

\put(0,-.5){\makebox(0,0){$(0,0)$}}

\multiput(0,2)(1,0){5}{\circle*{.1}}

\multiput(0,2)(0,1){5}{\circle*{.1}}

\end{picture}
\caption{\label{figu2} $E_{2,\pi_o}^{p,q}[\Pi^{\oo,o}]$ de (\ref{ssce}) pour $\Pi_\oo \simeq
\speh_s(\pi_\oo)$ non compatible avec l'aboutissement}
\end{figure}

%% file: introduction-chap5.tex
On donne quelques applications des calculs du chapitre précédent.

\medskip

\noindent \textbf{0.1.} --- On commence, proposition (\ref{prop-lrs2}), par calculer les
$\Pi^{\oo,o}$-parties des groupes de cohomologie de la fibre générique pour $\Pi$ automorphe tel que $\Pi_o
\simeq \speh_s(\pi_o)$, alors que nous n'avions traité, corollaire (\ref{coro-ssce-min}), que celles de poids
$s(g-1)$. \footnote{Ce résultat correspond à la conjecture 14.21 de \cite{lrs}.}

\medskip

\noindent \textbf{0.2.} --- On en déduit ensuite, proposition (\ref{theo-jl}), une correspondance de
Jacquet-Langlands entre les représentations automorphes de $D_\Am^\times$ et celles de $\bar D_\Am^\times$:
d'après les experts ce résultat s'obtient aisément à partir de la formule des traces simples.

\medskip

\noindent \textbf{0.3.} --- On décrit, proposition (\ref{prop-compo-locale}), les composantes locales des
représentations automorphes de $D_\Am^\times$: au final en égale caractéristique, corollaire
(\ref{compo-nulle}), on obtient ce que l'on déduirait des résultats de Moeglin-Waldspurger sur les
représentations de carré intégrable de $GL_d(\Am)$ et de la conjecture de Ramanujan-Peterson, prouvée par
Lafforgue dans ce cadre, si on disposait d'une correspondance de Jacquet-Langlands globale entre $GL_d(\Am)$
et $D_\Am^\times$. \footnote{En particulier, on en déduit que l'hypothèse 14.23 de \cite{lrs} est vérifiée.}

\medskip

\noindent \textbf{0.4.} --- On conclut enfin, théorème (\ref{theo-mono-glob}), par la preuve, qui n'a
d'intérêt qu'en caractéristique mixte, de la conjecture de monodromie-poids globale, i.e. les gradués de la
filtration de monodromie des groupes de cohomologie de la fibre générique sont purs.

%% file: applications.tex
\section{Preuve de la conjecture (14.21) de [LRS]}

\begin{propb} \label{prop-lrs2}
(cf. \cite{lrs} (14.21)) Soient $g$ un diviseur de $d=sg$ et $\pi_o$ (resp. $\pi_\oo$) une représentation
irréductible cuspidale unitaire de $GL_g(F_o)$ (resp. de $GL_g(F_\oo)$). On considère une représentation
automorphe $\Pi$ de $D_\Am^\times$ telle que $\Pi_o \simeq [\overrightarrow{s-1}]_{\pi_o}$ (resp. $\Pi_\oo
\simeq [\overrightarrow{s-1}]_{\pi_\oo}$). Le $GL_d(F_o) \times W_o$-module
$H^{d+s-2-2i}_{\eta_o,\r_\oo}[\Pi^{\oo,o}]$ est alors isomorphe à
$$[\overrightarrow{s-1}]_{\pi_o} \otimes L_g(\pi_o) (-\frac{d+s-2-2i}{2})$$
pour $0 \leq i <s$ et $\r_\oo:=\JL^{-1}([\overleftarrow{s-1}]_{\pi_\oo})$.
\end{propb}

\begin{proof} La proposition (\ref{prop-ssce-poids}) joint au théorème de Lefschetz difficile implique que pour tout
$0 \leq r \leq s-1$, $H^{d-s+2r}_{\eta_o,\r_\oo}[\Pi^{\oo,o}]$ admet $[\overrightarrow{s-1}]_{\pi_o} \otimes
L_g(\pi_o) (-\frac{s(g-1)+2r}{2})$ comme facteur direct. On reprend la suite spectrale (\ref{ss-poids}) en
utilisant la proposition (\ref{prop-not}) ainsi que le théorème de Lefschetz difficile. On pourra se reporter
à la figure (\ref{figure8}) où l'on a représenté les $H^i(gr_{k,\r_\oo})[\Pi^{\oo,o}]$ dans le cas $s=4$. On
rappelle tout d'abord que pour $1< \d <s$, on a
\begin{multline} \label{higrk}
H^{-\d}(gr_{k,\r_\oo})[\Pi^{\oo,o}]=\bigoplus_{\genfrac{}{}{0pt}{}{|k|<l \leq s}{\genfrac{}{}{0pt}{}{l \equiv
k+1 \mod 2}{l \equiv s+ \d \mod 2}}}  ([\overleftarrow{l-1}]_{\pi_o} \overrightarrow{\times}
[\overrightarrow{\frac{s-l+\d-2}{2}}]_{\pi_o}) \overleftarrow{\times}
[\overrightarrow{\frac{s-l-\d-2}{2}}]_{\pi_o} \\ \otimes L_g(\pi_o) (-\frac{d-1-\d+k}{2})
\end{multline}
\begin{multline}
H^{\d}(gr_{k,\r_\oo})[\Pi^{\oo,o}]=\bigoplus_{\genfrac{}{}{0pt}{}{|k|<l \leq s}{\genfrac{}{}{0pt}{}{l \equiv
k+1 \mod 2}{l \equiv s+ \d \mod 2}}}  ([\overleftarrow{l-1}]_{\pi_o} \overrightarrow{\times}
[\overrightarrow{\frac{s-l-\d-2}{2}}]_{\pi_o}) \overleftarrow{\times}
[\overrightarrow{\frac{s-l+\d-2}{2}}]_{\pi_o} \\ \otimes L_g(\pi_o) (-\frac{d-1+\d+k}{2})
\end{multline}
On remarque ainsi que tous les $[\overrightarrow{s-1}]_{\pi_o}$ des $E_{1,\r_\oo}^{p,q}[\Pi^{\oo,o}]$ de
(\ref{ss-poids}) restent dans l'aboutissement; il nous faut alors montrer que tous les autres disparaissent.
Dans la suite on ne considérera plus les $[\overrightarrow{s-1}]_{\pi_o}$.

Prenons dans un premier temps $\d>1$, de sorte qu'outre les $[\overrightarrow{s-1}]_{\pi_o}$, les éventuels
constituants non nuls de $E_{\oo,\r_\oo}^{d-s-\d}[\Pi^{\oo,o}]$ (resp. de
$E_{\oo,\r_\oo}^{d-s+\d}[\Pi^{\oo,o}]$) sont de la forme
$$[\overrightarrow{r},\overleftarrow{l-1},\overrightarrow{s-l-r}]_{\pi_o} \otimes L_g(\pi_o) (-\frac{d-1-\d-k}{2})$$
avec $s-l-r=r+\d,r+\d+1,r+\d-1$ (resp. $s-l-r=r-\d,r-\d-1,r-\d+1$) et certains entiers $k$ qu'il n'est pas
nécessaire de préciser. Pour $\d \geq 2$, on remarque que les ensembles $\{ \d,\d+1,\d-1 \}$ et $\{
-\d,-\d+1,-\d-1 \}$ sont disjoints de sorte que d'après le théorème de Lefschetz difficile, l'éventuel
constituant est forcément nul. \footnote{On peut aussi argumenter en remarquant que ce ne sont pas des
constituants locaux licites d'une représentation automorphe alors que dans l'aboutissement seules celles ci
doivent apparaître.}

Pour $\d=1$, le raisonnement est plus fin et demande de distinguer les sous-espaces des quotients dans nos
induites. Si on reprend le raisonnement précédent dans le cas $\d=1$, le théorème de Lefschetz difficile
impose que les éventuels constituants de $H^{d-1 \pm 1}_{\eta_o,\r_\oo}[\Pi^{\oo,o}]$ sont de la forme
\begin{equation} \label{forme}
[\overrightarrow{\frac{s-l}{2}},\overleftarrow{l-1},\overrightarrow{\frac{s-l}{2}}]_{\pi_o}  \otimes
L_g(\pi_o) (-\frac{d-1-k}{2})
\end{equation}
pour $1 \leq l \leq s$, $0 \leq r \leq (s-l)/2$ et certains entiers $k$ qu'il n'est pas nécessaire de
préciser. On revient sur (\ref{higrk}) et sur la suite spectrale (\ref{ss-poids}). Ainsi on a
\begin{multline} \label{e11}
E_{1,\r_\oo}^{-k,-1+k}[\Pi^{\oo,o}]=  \bigoplus_{\genfrac{}{}{0pt}{}{|k|<l \leq s}{\genfrac{}{}{0pt}{}{l
\equiv k+1 \mod 2}{l \equiv s-1 \mod 2}}}([\overleftarrow{l-1}]_{\pi_o} \overrightarrow{\times}
[\overrightarrow{\frac{s-l-1}{2}}]_{\pi_o}) \overleftarrow{\times} [\overrightarrow{\frac{s-l-3}{2}}]_{\pi_o}
\\ \otimes L_g(\pi_o) (-\frac{d-2+k}{2})
\end{multline}
\begin{multline} \label{e12}
E_{1,\r_\oo}^{-k-1,-1+k}[\Pi^{\oo,o}]=  \bigoplus_{\genfrac{}{}{0pt}{}{|k|<l \leq s}{\genfrac{}{}{0pt}{}{l
\equiv k+1 \mod 2}{l \equiv s- 2 \mod 2}}}([\overleftarrow{l-1}]_{\pi_o} \overrightarrow{\times}
[\overrightarrow{\frac{s-l}{2}}]_{\pi_o}) \overleftarrow{\times} [\overrightarrow{\frac{s-l-4}{2}}]_{\pi_o}
\\ \otimes L_g(\pi_o) (-\frac{d-2+k}{2})
\end{multline}
\begin{multline} \label{e10}
E_{1,\r_\oo}^{-k+1,-1+k}[\Pi^{\oo,o}]=  \bigoplus_{\genfrac{}{}{0pt}{}{|k|<l \leq s}{\genfrac{}{}{0pt}{}{l
\equiv k+1 \mod 2}{l \equiv s \mod 2}}}([\overleftarrow{l-1}]_{\pi_o} \overrightarrow{\times}
[\overrightarrow{\frac{s-l-2}{2}}]_{\pi_o}) \overleftarrow{\times} [\overrightarrow{\frac{s-l-2}{2}}]_{\pi_o}
\\ \otimes L_g(\pi_o) (-\frac{d-2+k}{2})
\end{multline}

Pour $k=s-1$, on remarque que la partie de poids $s(g+1)-2$ de $H^{d-2}_{\eta_o,\r_\oo}[\Pi^{\oo,o}]$ est un
sous-espace de $[\overleftarrow{s-2}]_{\pi_o} \overrightarrow{\times} [\overrightarrow{0}]_{\pi_o} \otimes
L_g(\pi_o) (-\frac{s(g+1)-2}{2})$. Or $[\overleftarrow{s-1}]_{\pi_o}$ n'étant pas un sous espace de
$[\overleftarrow{s-2}]_{\pi_o} \overrightarrow{\times} [\overrightarrow{0}]_{\pi_o}$, on en déduit, en
utilisant (\ref{forme}), que la partie de poids $s(g+1)-2$ de $H^{d-2}_{\eta_o,\r_\oo}[\Pi^{\oo,o}]$ est
nulle.

Pour $0 \leq k \leq s-3$, on a:
\begin{itemize}
\item $E_{1,\r_\oo}^{-k,-1+k}[\Pi^{\oo,o}]=E_{1,\r_\oo}^{-k-2,k+1}[\Pi^{\oo,o}] \otimes |-|^{2} \oplus
V_{k,-1} \otimes L_g(\pi_o) (-\frac{d-2+k}{2}) $, avec
$$V_{k,-1}=([\overleftarrow{k}]_{\pi_o}
\overrightarrow{\times} [\overrightarrow{\frac{s-k-2}{2}}]_{\pi_o}) \overleftarrow{\times}
[\overrightarrow{\frac{s-k-4}{2}}]_{\pi_o};$$

\item $E_{1,\r_\oo}^{-k-1,-1+k}[\Pi^{\oo,o}]=E_{1,\r_\oo}^{-k-3,k+1}[\Pi^{\oo,o}] \otimes |-|^{2} \oplus
V_{k,-2} \otimes L_g(\pi_o) (-\frac{d-2+k}{2})$, avec
$$V_{k,-2}=([\overleftarrow{k+1}]_{\pi_o}
\overrightarrow{\times} [\overrightarrow{\frac{s-k-2}{2}}]_{\pi_o}) \overleftarrow{\times}
[\overrightarrow{\frac{s-k-6}{2}}]_{\pi_o};$$

\item $E_{1,\r_\oo}^{-k+1,-1+k}[\Pi^{\oo,o}]=E_{1,\r_\oo}^{-k-1,k+1}[\Pi^{\oo,o}] \otimes |-|^{2} \oplus
V_{k,0} \otimes L_g(\pi_o) (-\frac{d-2+k}{2})$ avec $V_{0,0}=0$ et
$$V_{k,0}=([\overleftarrow{k-1}]_{\pi_o}
\overrightarrow{\times} [\overrightarrow{\frac{s-k-4}{2}}]_{\pi_o}) \overleftarrow{\times}
[\overrightarrow{\frac{s-k-4}{2}}]_{\pi_o} \hbox{ pour }0 < k \leq s-3.$$
\end{itemize}

On suppose alors par récurrence que
$$E_\oo^{-k-2,k+1}[\Pi^{\oo,o}]=(\ker d_1^{-k-2,k+1} / \im d_1^{-k-3,k+1})[\Pi^{\oo,o}]$$
est nulle; l'opérateur de monodromie $N$ implique alors que
$E_\oo^{-k,-1+k}[\Pi^{\oo,o}]$ est égal à $\ker d_+ /\im d_- \otimes L_g(\pi_o) (-\frac{d-2+k}{2})$ avec
$$V_{k,-2} \longmapright{d_-} V_{k,-1} \longmapright{d_+} V $$
Ainsi d'après (\ref{forme}), si $E_\oo^{-k,-1+k}[\Pi^{\oo,o}]$ était non nul, il serait égal à
$$[\overrightarrow{\frac{s-k}{2}},\overleftarrow{k-1},\overrightarrow{\frac{s-k}{2}}]_{\pi_o} \otimes
L_g(\pi_o) (-\frac{d-2+k}{2}).$$ En remarquant qu'aucun des constituants de
$[\overrightarrow{\frac{s-k}{2}},\overleftarrow{k-1}]_{\pi_o} \overrightarrow{\times}
[\overrightarrow{\frac{s-k-2}{2}}]_{\pi_o}$ n'est un constituant de $V_{k,-2}$, on en déduit que quelque soit
$d_-$, on a une surjection
$$V_{k,-1} / \im d_- \twoheadrightarrow [\overrightarrow{\frac{s-k}{2}},\overleftarrow{k-1}]_{\pi_o}
\overrightarrow{\times} [\overrightarrow{\frac{s-k-2}{2}}]_{\pi_o} \otimes L_g(\pi_o) (-\frac{d-2+k}{2})$$
alors que $[\overrightarrow{\frac{s-k}{2}},\overleftarrow{k-1},\overrightarrow{\frac{s-k}{2}}]_{\pi_o}$ n'est
pas un sous-espace de $[\overrightarrow{\frac{s-k}{2}},\overleftarrow{k-1}]_{\pi_o} \overrightarrow{\times}
[\overrightarrow{\frac{s-k-2}{2}}]_{\pi_o}$ de sorte que $E_\oo^{-k,-1+k}[\Pi^{\oo,o}]$ est nul.

Pour $1-s \leq k \leq 0$, on peut, dans le cas d'égale caractéristique, conclure en utilisant les cas $0 \leq
k \leq s-1$ et la pureté des $H^i_{\eta_o,\r_\oo}$, dont les gradués pour la filtration de monodromie sont
purs. Sinon on raisonne de manière strictement identique en étudiant les suites exactes
$$E_{1,\r_\oo}^{-k-1,k+1} \longmapright{d_1^{-k-1,k+1}} E_{1,\r_\oo}^{-k,k+1} \longmapright{d_1^{-k,k+1}}
E_{1,\r_\oo}^{-k+1,k+1}$$ et en remplaçant l'étude des sous-espaces par celle des quotients.

\end{proof}

On donne alors l'amélioration suivante de (\ref{somme-alternee}) qui constitue le pendant du corollaire
(\ref{coro-hij-nul}) dans le cas $\Pi_\oo=[\overrightarrow{s-1}]_{\pi_\oo}$.

\begin{corob} \label{strate-alt-p}
Soit $\Pi$ une représentation irréductible automorphe de $D_\Am^\times$ telle que $\Pi_\oo \simeq
[\overrightarrow{s-1}]_{\pi_\oo}$ où $\pi_\oo$ est une représentation cuspidale de $GL_g(F_\oo)$) avec
$d=sg$. On suppose que $\Pi_o \simeq [\overrightarrow{s-1}]_{\pi_o}$ où $\pi_o$ est une représentation
irréductible cuspidale unitaire de $GL_g(F_o)$. Si
$$[H^i_{h,\JL^{-1}([\overleftarrow{s-1}]_{\pi_\oo}),\t_o}(\Pi^{\oo,o})]$$
muni de son action naturelle de $GL_d(F_o)$ comme dans \cite{lrs}, est non nul dans le groupe de Grothendieck
$\groth(GL_{d-h}(F_o) \times (D_{o,h}^\times/\DC_{o,h}^\times))$, on a alors
\begin{itemize}

\item $\t_o=\JL^{-1}([\overleftarrow{s-1}]_{\pi_o})$;

\item $h=lg$ avec $1 \leq l \leq s$ et $i=(s-l)(g+1)$.
\end{itemize}
Dans ce cas, ils sont donnés par
\begin{multline} \label{formule-modif}
\lim_{\genfrac..{0pt}{1}{\to}{I}} H^{(s-l)(g+1)}_c(M_{I,s_o,1}^{=lg},\FC_{\t_o} \otimes \LC_{\r_\oo})[\Pi^{\oo,o}] 
\simeq \\
m(\Pi) [\overrightarrow{s-l-1}]_{\pi_o(l(g+1)/2)} \otimes (\Xi^{\frac{(s-l)(g+1)}{2}} \bigoplus_{\xi \in
\DF(\pi_o)} \xi^{-1})
\end{multline}
en tant que représentation de $GL_{(s-l)g}(F_o) \times (D_{o,h}^\times/\DC_{o,h}^\times)$, où $\DF(\pi_o)$
est l'ensemble des caractères $\xi$ de $\Zm \simeq D_{o,h}^\times /\DC_{o,h}^\times$ tels que
$$\JL^{-1}([\overleftarrow{s-1}]_{\pi_o}) \otimes \xi^{-1} \simeq \JL^{-1}([\overleftarrow{s-1}]_{\pi_o}).$$
\end{corob}

\begin{proof} On rappelle que d'après le lemme (\ref{lem-rj-combi}), les composantes $\Pi^{\oo,o}$-isotypiques
des $H^i_c(M_{I,s_o}^{=lg},\FC(g,l,\pi_o) \otimes \LC_{\r_\oo}) \otimes \Pi_l$ sont mixtes de poids
$(s-l)(g-1)+2k$ avec $0 \leq k < s-l$ et de la forme
$$(\Pi_l \overrightarrow{\times} \pi_+) \overleftarrow{\times} \pi_- \otimes L_g(\pi_o)
(-\frac{(s-l)(g-1)+2k}{2})$$ où $\pi_+$ (resp. $\pi_-$) est une représentation elliptique de
$GL_{(s-l-k)g}(F_o)$ (resp. de $GL_{kg}(F_o)$). On raisonne alors par récurrence sur $k$, de $0$ à $s-l-2$,
afin de montrer que pour tout $i$, les parties de poids $(s-l)(g-1)+2k$ des
$H^i_c(M_{I,s_o}^{=lg},\FC(g,l,\pi_o,I) \otimes \LC_{\r_\oo})[\Pi^{\oo,o}]$ sont nulles. Le cas $k=0$ a été
traité dans la proposition (\ref{prop-hic-poids}). Supposons le résultat acquis jusqu'au rang $k<s-l-2$ et
montrons le au rang $k+1$. On raisonne alors par récurrence sur $l$ de $1$ à $s-1$. L'initialisation de la
récurrence se prouve comme le cas général; on suppose donc le résultat acquis jusqu'au rang $l<s-2$ et
prouvons le au rang $l+1$.

Supposons qu'il existe $j$ tel que la partie de poids $(s-l)(g-1)+2(k+1)$ de
$H^j_c(M_{I,s_o}^{=lg},\FC(g,l,\pi_o,I) \otimes \LC_{\r_\oo})[\Pi^{\oo,o}]$ soit non nulle. On étudie ensuite
la suite spectrale (\ref{sss2}) associée à la stratification pour $i=l(g-1)$. Ainsi la partie de poids
$m_k:=s(g-1)+2(k+1)$ de $E_{1,\r_\oo}^{lg+1,j-lg-1;l(g-1)}[\Pi^{\oo,o}]$ est non nulle et de la forme
$([\overrightarrow{l-1}]_{\pi_o} \overrightarrow{\times} \pi_+) \overleftarrow{\times} \pi_- \otimes
L_g(\pi_o) (-\frac{m_k}{2})$, où $\pi_+$ (resp. $\pi_-$) est une représentation elliptique de
$GL_{(s-l-k)g}(F_o)$ (resp. de $GL_{kg}(F_o)$). Or d'après l'hypothèse de récurrence pour tout $1 \leq l'
<l$, les parties de poids $m_k=(s-l')(g-1)+l'(g-1)+2(k+1)$ des
$E_{1,\r_\oo}^{l'g+1,j'-l'g-1;l(g-1)}[\Pi^{\oo,o}]$ sont nulles; en effet celles-ci sont données par les
parties de poids $(s-l')(g-1)+2k'$ avec $s(g-1)+2(k+1)=(s-l')(g-1)+l'(g-1)+2\d+2k'$ avec $l(g-1)=l'(g-1)+\d$,
$0 < \d < l'<l$, soit $k'=k+1-\d$. Pour $l'>l$, les parties de poids $m_k$ des
$E_{1,\r_\oo}^{l'g+1,j'-l'g-1;l(g-1)}[\Pi^{\oo,o}]$ sont nulles car $l(g-1)$ ne s'écrit pas sous la forme
$l'(g-1)+\d$ avec $0 \leq \d < l'$.

On étudie alors la suite spectrale (\ref{ssce2}) des cycles évanescents. D'après ce qui précède, la partie de
poids $m_k$ de $E_{2,\r_\oo}^{j,l(g-1)}[\Pi^{\oo,o}]$ est de la forme $([\overrightarrow{l-1}]_{\pi_o}
\overrightarrow{\times} \pi_+) \overleftarrow{\times} \pi_- \otimes L_g(\pi_o) (-\frac{m_k}{2})$. Par
ailleurs les parties de poids $m_k$ des $E_{2,\r_\oo}^{j+r+1,l(g-1)-r}[\Pi^{\oo,o}]$, pour $r \geq 1$ sont
nulles; en effet celles-ci proviendraient à travers la suite spectrale (\ref{sss2}), des parties de poids
$(s-l')(g-1)+2k'$ de la composante $\Pi^{\oo,o}$-isotypique de
$$H^{j+r-1}_c(M_{I,s_o}^{=l'g},\FC(g,l',\pi_o) \otimes \LC_{\r_\oo}) \otimes [\overleftarrow{k'},
\overrightarrow{l-1-k'}]_{\pi_o},$$ pour $1 \leq l' <l$ avec $s(g-1)+2(k+1)=(s-l')(g-1)+l'(g-1)+2\d+2k'$,
soit $k'=k+1-\d$ avec $0 < \d < l'$ qui sont nulles d'après l'hypothèse de récurrence. En ce qui concerne les
parties de poids $m_k$ des $E_{2,\r_\oo}^{j-r-1,l(g-1)+r}[\Pi^{\oo,o}]$ pour $r >0$, elles proviennent à
nouveau des parties de poids $(s-l')(g-1)+2k'$ des composantes $\Pi^{\oo,o}$-isotypiques des
$H^{j-r-1}_c(M_{I,s_o}^{=l'g},\FC(g,l',\pi_o) \otimes \LC_{\r_\oo}) \otimes
[\overleftarrow{k'},\overrightarrow{l-1-k'}]_{\pi_o}$ pour $1 \leq l' \leq s$ avec $k'=k+1-\d$ et $0 \leq
\d=(l-l')(g-1)+r < l'$, ce qui impose d'après l'hypothèse de récurrence $l'>l$, avec $\d=0$ et $k'=k+1$. On
remarque alors que
$$[\overleftrightarrow{k},\overrightarrow{l-1},\overleftarrow{1},\overleftrightarrow{s-l-k-1}]_{\pi_o}
\otimes L_g(\pi_o) (-\frac{s(g-1)+2(k+1)}{2})$$ serait un constituant de $E_{\oo,\r_\oo}^j[\Pi^{\oo,o}]$, ce
qui n'est pas d'après la proposition (\ref{prop-lrs2}).

\end{proof}

\rem On pourra voir la figure (\ref{figure6}), où l'on a représenté pour $s=4$ et $g=2$, les
$E_2^{p,q}[\Pi^{\oo,o}]$ de la suite spectrale des cycles évanescents.

\section{Correspondances de Jacquet-Langlands globales}

\textit{Badulescu m'a expliqué que les résultats de ce paragraphe peuvent s'obtenir aisément via la formule des traces 
simples. Nous
avons tout de même souhaité montrer comment ils découlaient de nos calculs.}

\medskip

On rappelle que $D$ et $\bar D$ sont des algèbres à division sur un corps global $F$ d'égale caractéristique
$p$, telles que:
\begin{itemize}
\item pour tout place $x \neq \oo,o$, $D_x \simeq \bar D_x$;

\item $D_\oo^\times \simeq GL_d(F_\oo)$ et $\inv_\oo \bar D^\times=-1/d$;

\item $D_o^\times \simeq GL_d(F_o)$ et $\inv_o \bar D^\times=1/d$.
\end{itemize}

\begin{propb} \label{theo-jl}
On fixe un diviseur $g$ de $d=sg$ ainsi qu'une représentation irréductible cuspidale unitaire $\pi_\oo$ de
$GL_g(F_\oo)$. Il existe alors une bijection dite de Jacquet-Langlands entre:

\noindent - les représentations irréductibles automorphes $\bar \Pi$ de $\bar D_\Am^\times$ telles que $\bar
\Pi_\oo \simeq \JL^{-1}([\overleftarrow{s-1}]_{\pi_\oo})$,

\noindent - les représentations irréductibles automorphes $\Pi$ de $D_\Am^\times$ vérifiant l'une des deux
conditions suivantes:
\begin{itemize}
\item[(a)] $\Pi_\oo$ est isomorphe à $[\overleftarrow{s-1}]_{\pi_\oo}$ et $\Pi_o$ est une représentation
de carré intégrable;

\item[(b)] $\Pi_\oo \simeq [\overrightarrow{s-1}]_{\pi_\oo}$ et $\Pi_o \simeq [\overrightarrow{s-1}]_{\pi_o}$
pour $\pi_o$ une représentation irréductible cuspidale de $GL_g(F_o)$.
\end{itemize}
compatible aux correspondances de Jacquet-Langlands locales, soit $\bar \Pi^{\oo,o} \simeq \Pi^{\oo,o}$ et
$\Pi_o \simeq \JL(\bar \Pi_o)$ dans le cas (a) et $\Pi_o \simeq \lexp t \JL(\bar \Pi_o)$ dans le cas (b) où
$\lexp t \pi$ désigne la représentation duale associée à $\pi$ pour la dualité de Zelevinski. En outre on a
$m(\Pi)=m(\bar \Pi)$.

Par ailleurs soit $\Pi^o$ une représentation de $(D_\Am^o)^\times$ telle que pour toute représentation
$\Pi_o$ de $GL_d(F_o)$ avec $\Pi:=\Pi^o \Pi_o$ vérifiant $\hyp(\oo)$, $\Pi_o$ n'est pas de la forme
$[\overleftarrow{s-1}]_{\pi_o}$ ou $[\overrightarrow{s-1}]_{\pi_o}$ pour $\pi_o$ une représentation
irréductible cuspidale de $GL_g(F_o)$ avec $d=sg$. Alors il n'existe pas de représentation irréductible
automorphe $\bar \Pi$ de $\bar D_\Am^\times$ telle que $\bar \Pi^{\oo,o} \simeq \Pi^{\oo,o}$.
\end{propb}

\begin{proof} Posons $\r_\oo=\JL^{-1}([\overleftarrow{s-1}]_{\pi_\oo})$.

- Soit $\Pi$ une représentation automorphe de $D_\Am^\times$ telle que $\Pi_\oo \simeq
[\overleftarrow{s-1}]_{\pi_\oo}$ et $\Pi_o \simeq [\overleftarrow{s'-1}]_{\pi_o}$ pour une certaine
représentation cuspidale unitaire $\pi_o$ de $GL_{g'}(F_o)$. D'après la proposition (\ref{prop-coho1}), pour
tout $1 \leq l < s'$, on a
\begin{multline*}
\lim_{\genfrac..{0pt}{1}{\to}{I}} H^0(M_{I,s_o}^{\geq lg'},j^{\geq lg'}_!
HT_{\rho_\oo}(g',l,\pi_o,\Pi_l,I)[d-lg']) [\Pi^{\oo,o}]=
\\ m(\Pi) \Pi_l \overrightarrow{\times} [\overleftarrow{s'-l-1}]_{\pi_o} \otimes (\Xi^{\frac{(s'-l)(g'-1)}{2}}
\otimes \bigoplus_{\xi \in \AF(\pi_o)} \xi^{-1})
\end{multline*}
Par ailleurs ce dernier est aussi égal à
\begin{multline*}
\lim_{\genfrac..{0pt}{1}{\to}{I}} H^0(M_{I,s_o}^d,\FC(g',s',\pi_o,I) \otimes \LC_{\r_\oo} \otimes \Pi_l
\overrightarrow{\times} [\overleftarrow{s'-l-1}]_{\pi_o})[\Pi^{\oo,o}] \\
\otimes (\Xi^{\frac{(s'-l)(g'-1)}{2}} \otimes \bigoplus_{\xi \in \AF(\pi_o)} \xi^{-1})
\end{multline*}
qui d'après le lemme (\ref{lem-pts-ss}) est égal à
$$\sum_{\bar \Pi} \Pi_l \overrightarrow{\times} [\overleftarrow{s'-l-1}]_{\pi_o} \otimes 
(\Xi^{\frac{(s'-l)(g'-1)}{2}}
\otimes \bigoplus_{\xi \in \AF(\pi_o)} \xi^{-1})$$ où $\bar \Pi$ décrit les représentations irréductibles
automorphes de $\bar D_\Am^\times$ telles que $\bar \Pi^{\oo,o} \simeq \Pi^{\oo,o}$, $\bar \Pi_\oo \simeq
\JL^{-1}([\overleftarrow{s-1}]_{\pi_\oo})$ et $\bar \Pi_o \simeq \JL^{-1}([\overleftarrow{s'-1}]_{\pi_o})$.

- De la même façon, soit $\Pi$ telle que $\Pi_\oo \simeq [\overrightarrow{s-1}]_{\pi_\oo}$ et $\Pi_o \simeq
[\overrightarrow{s'-1}]_{\pi_o}$. On reprend la preuve de la proposition (\ref{prop-not}). Pour $l=s'$, le
lemme (\ref{lem-pts-ss}) donne
\begin{multline*}
\lim_{\genfrac..{0pt}{1}{\to}{I}} H^0(M_{I,s_o}^{s'g'},\FC(g',s',\pi_o,I) \otimes
[\overleftarrow{s'-1}]_{\pi_o})[\Pi^{\oo,o}]= \\
\Bigl ( \sum_{\bar \Pi \in \UF_{\bar D}(\Pi^{\oo,o})} m(\bar \Pi) \Bigr ) [\overleftarrow{s'-1}]_{\pi_o}
\otimes \bigoplus_{\xi \in \AF(\pi_o)} \xi^{-1}
\end{multline*}
où $\UF_{\bar D}(\Pi^{\oo,o})$ désigne l'ensemble des $\bar \Pi$ automorphes telles que $\bar \Pi^{\oo,o}
\simeq \Pi^{\oo,o}$, $\bar \Pi_\oo \simeq \JL^{-1}([\overleftarrow{s-1}]_{\pi_\oo})$ et $\bar \Pi_o \simeq
[\overleftarrow{s'-1}]_{\pi_o}$. Par ailleurs pour $l=s'-1$, on trouve que
\begin{multline*}
\sum_i (-1)^i H^i(j^{\geq (s'-1)g'}_{!*})=\biggl ( -\Bigl ( \sum_{\bar \Pi \in \UF_{\bar D}(\Pi^{\oo,o})}
m(\bar \Pi) \Bigr ) [\overleftarrow{s'-2}]_{\pi_o} \overrightarrow{\times} [\overleftarrow{0}]_{\pi_o}
\otimes \Xi^{-1/2} \\
- \Bigl ( \sum_{\Pi \in \UF_D(\Pi^{\oo,o})} m(\Pi) \bigr ) [\overleftarrow{s'-2}]_{\pi_o}
\overleftarrow{\times} [\overleftarrow{0}]_{\pi_o} \otimes \Xi^{1/2} \biggr ) ( \Xi^{g'/2} \otimes
\bigoplus_{\xi \in \AF(\pi_o)} \xi^{-1})
\end{multline*}
D'après la pureté, le terme de gauche de l'égalité ci-dessus doit être égal à
$$-H^{-1}(j^{\geq (s'-1)g'}_{!*})-H^1(j^{\geq (s'-1)g'}_{!*})$$
de sorte que la dualité de Verdier donne le résultat.

- Supposons que $\Pi_o$ ne soit ni de la forme $[\overleftarrow{s'-1}]_{\pi_o}$ ni
$[\overrightarrow{s'-1}]_{\pi_o}$ pour tout diviseur $s'$ de $d=s'g'$ et toute représentation cuspidale
$\pi_o$ de $GL_{g'}(F_o)$. Supposons qu'il existe une représentation $\bar \Pi$ telle que $\bar \Pi^{\oo,o}
\simeq \Pi^{\oo,o}$ et soit $s'$, $\pi_o$ tel que $\bar \Pi_o \simeq
\JL^{-1}([\overleftarrow{s'-1}]_{\pi_o})$. Le lemme (\ref{lem-pts-ss}) donne comme précédemment
$$H^0(\FC(g,s',\pi_o) \otimes
[\overleftarrow{s'-1}]_{\pi_o}) \simeq \bigl ( \sum_{\bar \Pi \in \UF_{\bar D}(\Pi^{\oo,o})} m(\bar \Pi)
\bigr ) [\overleftarrow{s'-1}]_{\pi_o} \otimes \bigoplus_{\xi \in \AF(\pi_o)} \xi^{-1}$$ de sorte que
\begin{multline*}
\sum_i (-1)^i H^i(j^{\geq (s'-1)g'}_{!*})[\Pi^{\oo,o}]=\biggl ( - \Bigl ( \sum_{\Pi \in \UF_D(\Pi^{\oo,o})}
m(\Pi)
\Bigr ) (\ind_{P_{d-g',d}(F_o)}^{GL_d(F_o)} [\overleftarrow{s'-2}]_{\pi_o} \otimes \\
\red_{\JL^{-1} [\overleftarrow{s'-2}]_{\pi_o}}^{d-g'}(\Pi_o))  +  \Bigl ( \sum_{\bar \Pi \in \UF_{\bar
D}(\Pi^{\oo,o})} m(\bar \Pi) \Bigr ) [\overleftarrow{s'-2}]_{\pi_o} \overrightarrow{\times}
[\overleftarrow{0}]_{\pi_o} \otimes \Xi^{(g'-1)/2} \biggr ) \otimes \bigoplus_{\xi \in \AF(\pi_o)} \xi^{-1})
\end{multline*}
La condition de pureté et la compatibilité à la dualité de Verdier, impose alors que $\Pi_o$ est soit
isomorphe à $[\overleftarrow{s'-1}]_{\pi_o}$ soit à $[\overrightarrow{s'-1}]_{\pi_o}$ d'où la contradiction.

Ce dernier raisonnement montre en outre qu'étant donné une représentation $\bar \Pi$ comme dans l'énoncé, il
lui correspond une représentation $\Pi$ vérifiant les conditions de l'énoncé. Dans le cas où $\Pi_o \simeq
[\overrightarrow{s'-1}]_{\pi_o}$, l'égalité $s=s'$ découle du corollaire (\ref{coro-compo-locale}).

 \end{proof}

\section[Composantes locales]{Composantes locales des représentations automorphes vérifiant $\hyp(\oo)$}

\textit{Les résultats ci-dessous sont bien connus dans le cas où il existe une correspondance de
Jacquet-Langlands vers $GL_d(F)$, à partir des travaux de Moeglin-Waldspurger. C'est le cas par exemple s'il
existe une place $x$ telle que $D_x$ est une algèbre à division sur $F_x$.}

\medskip

\begin{propb} \label{prop-compo-locale}
Soit $\Pi$ une représentation irréductible automorphe de $D_\Am^\times$ vérifiant $\hyp(\rho_\oo)$. Pour tout
$0 \leq i < d$, il existe alors un réel $\d$ ainsi que des entiers $n_i$ \footnote{$n_i$ représente la
dimension de la partie primitive de la représentation galoisienne de $H^i_{\eta_o}$} tels que:

\begin{itemize}
\item $\sum_{i=0}^{d-1} n_i(d-i)=d$,

\item pour $i$ tel que $n_i \neq 0$, $n_j=0$ pour $j \equiv i+1 \mod 2$, tels que pour toute place $o$ telle que
$D_o^\times \simeq GL_d(F_o)$, il existe:
\begin{itemize}
\item des entiers $t_j>1,~g_j$, $1 \leq j \leq u$, et
$t'_k,~g'_k$, $1 \leq k \leq u'$, avec $\sum_{j=1}^u t_jg_j+\sum_{k=1}^{u'} t'_k g'_k=d$,

\item des représentations irréductibles cuspidales $\pi_{o,j}$ et $\pi_{o,k}'$ de respectivement $GL_{g_j}(F_o)$ et 
$GL_{g'_k}(F_o)$
telles que $\pi_{o,j} (\delta)$ et $\pi'_{o,k} (\delta)$ soient unitaires
\end{itemize}
vérifiant les conditions suivantes:

\begin{itemize}
\item pour tout $0 \leq i < d-1$, $\sum_{k~/~ t'_k=d-1-i} g'_k=n_i$;

\item $\Pi_o$ est l'induite irréductible $[\overleftarrow{t_1-1}]_{\pi_{o,1}} \times \cdots \times
[\overleftarrow{t_u-1}]_{\pi_{o,u}} \times [\overrightarrow{t_1'-1}]_{\pi'_{o,1}} \times \cdots \times
[\overrightarrow{t'_{u'}-1}]_{\pi'_{o,u'}}$, avec:

\begin{itemize}
\item si $u >0$ alors pour tout $j$, $s_j \equiv 1 \mod 2$;

\item si $u=0$ alors tous les $s_i$ ont la même parité donnée par $(-1)^{s_i-1}=\e(\Pi)$.
\end{itemize}
\end{itemize}
\end{itemize}
\end{propb}

On en déduit en particulier les corollaires suivant.

\begin{corob} Pour toute représentation automorphe de $D_\Am^\times$ vérifiant $\hyp(\rho_\oo)$, la conjecture de
Ramanujan-Peterson est vérifiée.
\end{corob}

\begin{corob} \label{coro-compo-locale}
Soit $\Pi$ une représentation irréductible automorphe de $D_\Am^\times$ vérifiant $\hyp(\rho_\oo)$ et supposons
qu'il existe une place $o_0$ telle que $D_{o_0}^\times \simeq GL_d(F_{o_0})$ et $\Pi_{o_0}$ tempérée, i.e.
avec les notations du théorème précédent $t'_i=1$ pour tout $1 \leq i \leq u'$ (resp. $\Pi_{o_0} \simeq
[\overrightarrow{s-1}]_{\pi_{o_0}}$ pour $\pi_{o_0}$ une représentation cuspidale de $GL_g(F_{o_0})$ avec
$d=sg$). On en déduit alors, en utilisant les notations du théorème précédent, que pour toute place $o$ non
ramifiée, pour tout $1 \leq k \leq u'$ (resp. pour tout $1 \leq j \leq u$) $t'_k=0$ (resp. $t_j=0$), i.e.
$\Pi_o$ est tempérée et $n_i=0$ pour tout $0 \leq i <d-1$ (resp. $n_i=0$ pour $i \neq s$ et $\Pi_o \simeq
[\overrightarrow{s-1}]_{\pi'_{o,1} \times \cdots \times \pi'_{o,u'}}$). Dans le cas d'égale caractéristique,
on en déduit en outre que $\Pi_\oo \simeq [\overleftarrow{s'-1}]_{\pi_\oo}$ (resp. $\Pi_\oo \simeq
[\overrightarrow{s-1}]_{\pi_\oo}$) avec $\rho_\oo$ cuspidale de $GL_{g'}(F_\oo)$ (resp. de $GL_g(F_\oo)$ i.e.
le même $s$ et $g$ que celui en la place $o_0$) avec $d=s'g'$ (resp. $d=sg$).
\end{corob}

\begin{corob} \label{compo-nulle}
En égale caractéristique, soit $\Pi$ une représentation irréductible automorphe de $D_\Am^\times$ telle que
$\Pi_\oo$ est de la forme $[\overleftarrow{s-1}]_{\pi_\oo}$ (resp. $[\overrightarrow{s-1}]_{\pi_\oo}$) pour
$\pi_\oo$ une représentation irréductible cuspidale de $GL_g(F_\oo)$ avec $d=sg$. Pour toute place $o$ où
$D_o^\times \simeq GL_d(F_o)$, la composante locale $\Pi_o$ est tempérée, i.e. de la forme
$[\overleftarrow{t_1-1}]_{\pi_{o,1}} \times \cdots \times [\overleftarrow{t_r-1}]_{\pi_{o,u}}$ (resp. de la
forme $[\overrightarrow{s-1}]_{\pi_{o,1}} \times \cdots \times [\overrightarrow{s-1}]_{\pi_{o,u'}}$).
\end{corob}

\rem On en déduit en particulier que l'hypothèse 14.23 de \cite{lrs} est vérifiée.

\begin{proof} On reprend les notations de la proposition (\ref{prop-compo-locale}). Dans le cas $\Pi_\oo \simeq
[\overleftarrow{s-1}]_{\pi_\oo}$, si $u >0$ (resp. $u=0$) la représentation galoisienne
$$\s(\Pi):=|\cl |^{(d-1)/2} \bigoplus_{i \equiv 1 \mod 2} H^{d-1+i}_{\eta_o,\rho_\oo}[\Pi^\oo]$$
$$(\hbox{resp. } \s(\Pi):=|\cl |^{(d-1)/2} \bigoplus_{i \equiv \e(\Pi)+1 \mod 2} 
H^{d-1+i}_{\eta_o,\rho_\oo}[\Pi^\oo])$$
est telle qu'en toute place $x \neq \oo$ telle que $D_x^\times \simeq GL_d(F_x)$, $\s(\Pi)_x$ est la
représentation de Langlands $L_d(\Pi_x)$ associée à la composante locale $\Pi_x$ de $\Pi$. Il est alors connu
qu'il en est de même à la place $\oo$ de sorte que $\s(\Pi)$ est irréductible ce qui implique qu'il existe un
unique $i$ tel que $H^{d-1+i}_{\eta_o,\rho_\oo}[\Pi^\oo]$ soit non nul: la dualité implique qu'il s'agit de $i=0$,
d'où le résultat.

Dans le cas où $\Pi_\oo \simeq [\overrightarrow{s-1}]_{\pi_\oo}$, on considère $\s(\Pi):= \bigoplus_{i \equiv
s-1 \mod 2} H^{d-1+i}_{\eta_o,\rho_\oo}[\Pi^\oo]$. Le même raisonnement que précédemment s'applique où cette
fois $\s(\Pi)_\oo$ est la somme directe de $s$ représentation irréductible de dimension $g$. Or d'après le
théorème de Cebotarev, $\s(\Pi)=\bigoplus_{i \equiv s-1 \mod 2} V_{d-1+i}$ où
$V_{d-1+i}=\bigoplus_{\genfrac{}{}{0pt}{}{0 \leq r \leq d-1}{r \equiv s-1 \mod 2}} W_r$ où $W_r$ est de
dimension $n_r$. D'après la pureté on en déduit alors que chacun des $V_{d-1+i}$ pour $i \equiv s-1 \mod 2$
et $|i|< s$ est de dimension $g$ de sorte que tous les $n_j$ sont nuls sauf $n_{d-s}=g$, d'où le résultat.

\end{proof}

\rem Dans le cas où l'on dispose de la correspondance de Jacquet-Langlands globale entre $D_\Am^\times$ et
$GL_d(\Am)$, ces résultats sont compatibles avec ceux de Moeglin-Waldspurger. En effet soit $\JL(\Pi)$ la
représentation de $GL_d$ correspondante qui est donc de carré intégrable modulo le centre. Si $\JL(\Pi)_\oo
\simeq [\overleftarrow{s-1}]_{\pi_\oo}$ pour $\pi_\oo$ irréductible cuspidale, alors $\JL(\Pi)$ est cuspidale
de sorte que ses composantes locale $\JL(\Pi)_o$ sont génériques et même tempérée d'après la conjecture de
Ramanujan-Peterson prouvée dans ce cadre par Lafforgue. Si $\JL(\Pi)_\oo \simeq
[\overrightarrow{s-1}]_{\pi_\oo}$ alors $\JL(\Pi)$ ne peut être cuspidale car $[\overrightarrow{s-1}]_{\oo}$
n'est pas générique de sorte que d'après Moeglin-Waldspurger, $\JL(\Pi)$ est de la forme
$[\overrightarrow{s-1}]_{\tilde \Pi}$ pour $\tilde \Pi$ une représentation cuspidale de $GL_g(\Am)$.

\begin{proof} \textit{de la proposition (\ref{prop-compo-locale}):} A torsion par un caractère près,
$\Pi_o$ est unitaire et donc d'après la classification de Tadic de la forme
$$[\overrightarrow{s_1-1}]_{[\overleftarrow{t_1-1}]_{\pi_{o,1}(\lambda_1)}} \times \cdots \times
[\overrightarrow{s_u-1}]_{[\overleftarrow{t_u-1}]_{\pi_{o,u}(\lambda_u)}}$$ où $\pi_{o,i}$ sont des
représentations irréductibles cuspidales unitaires de $GL_{g_i}(F_o)$ avec $\lambda_i \in ]
\frac{-1}{2},\frac{1}{2}[$ et $\sum_{i=1}^u g_is_it_i=d$ et \textit{où on induit par rapport aux paraboliques
opposés aux paraboliques standards}. L'induite étant irréductible, l'ordre des facteurs n'importe pas. Soit
alors $r \geq 1$ tel que $\pi_o:=\pi_{o,1}(\lambda_1) \simeq \cdots \simeq \pi_{o,r}(\lambda_r)$ et
$\pi_{o,i}(\lambda_i)$ n'est pas isomorphe à $\pi_o$ pour $i
>r$. On note $l_0:=max_{1 \leq i \leq r} \{ s_i,t_i \}$. On note $g$ l'entier tel que $\pi_o$ est une
représentation cuspidale de $GL_g(F_o)$. Par ailleurs soit $0\leq k_1$ (resp. $k_1 \leq k_2$, resp. $k_2 \leq
k_3 \leq r$) tel que pour tout $1 \leq i \leq k_1$ (resp. $k_1 < i \leq k_2$, resp. $k_2 < i \leq k_3$) on
ait $s_i=t_i=l_0$ (resp. $t_i < s_i=l_0$, resp. $s_i < t_i=l_0$) et $\max_{k_3 < i \leq r} \{ t_i,s_i \} <
l_0$.

\begin{lemb} \label{lem-red1}
Pour tout $\max \{ s ,t \}< l \leq st$, $\red_{\JL^{-1}([\overleftarrow{l-1}]_{\pi_o})}
([\overrightarrow{s-1}]_{[\overleftarrow{t-1}]_{\pi_o}})=(0)$.
\end{lemb}

\begin{proof} Le résultat découle essentiellement de \cite{ze} et en particulier du lemme suivant.

\begin{lemb} \label{lem-ze}
Pour tout $s,t$, les constituants de $J_{N_{g,2g,\cdots,stg}}
([\overrightarrow{s-1}]_{[\overleftarrow{t-1}]_{\pi_o}})$ sont de la forme $\pi_o(\frac{t+s-2}{2}-\s(1))
\otimes \cdots \otimes \pi_o(\frac{t+s-2}{2}-\s(st))$ où le multi-ensemble $R:=\{ \s(1),\cdots, \s(st) \}$
est tel que $\{ \pi_o(\frac{t+s-2}{2}-r)~/~r \in R \}$ soit égal au support cuspidal de
$[\overrightarrow{s-1}]_{[\overleftarrow{t-1}]_{\pi_o}}$. Pour tout $0 \leq i \leq st-1$, on pose
$\s^{-1}(i)=\{ n_i(1) < n_i(2) < \cdots \}$. On a alors les conditions suivantes:

\begin{itemize}
\item[(i)] pour tout $0 \leq k <s-1$ et $0 \leq i < t$, on a $n_{t-1+k-i}(\min(k+1,i+1)) > n_{t+k-i}(\min(k+2,i+1))$;

\item[(ii)] pour tout $0 \leq k < s$ et tout $0 \leq i < j <t$, on a $n_{t-1+k-i}(\min(k+1,i+1)) >
n_{t-1+k-j}(\min(k+1,j+1))$.
\end{itemize}
\end{lemb}

\rem Quand on regarde le support cuspidal de
$$\overbrace{[\overleftarrow{t-1}]_{\pi_o(\frac{s-1}{2})} \boxplus \cdots
[\overleftarrow{t-1}]_{\pi_o(\frac{1-s}{2})}}^s$$ dans l'ordre de gauche à droite, pour tout $0 \leq k < s-1$
et $0 \leq i < t$, le $\pi_o(\frac{t+s-2}{2}-(t-1+k-i))$ du $k$-ème facteur est le $\min(k+1,i+1)$-ème. La
condition (ii) du lemme affirme que dans un même facteur, on doit prendre les $\pi_o(i)$ de gauche à droite,
tandis que la condition (i) précise que le $\pi_o(i)$ du $k$-ème facteur arrive après le $\pi_o(i+1)$ du
$k+1$-ème facteur.

\begin{proof} En remarquant que $([\overrightarrow{s-1}]_{[\overleftarrow{t-1}]_{\pi_o}})$ s'injecte dans les
induites
$$[\overleftarrow{t-1}]_{\pi_o(\frac{s-1}{2})} \times
[\overrightarrow{s-2}]_{[\overleftarrow{t-1}]_{\pi_o(-1/2)}} \qquad
[\overleftarrow{t-1}]_{\pi_o(\frac{s-1}{2})} \boxplus [\overleftarrow{t-1}]_{\pi_o(\frac{s-3}{2})} \times
\overrightarrow{(s-2)}_{[\overleftarrow{t-1}]_{\pi_o(-1)}}$$ on se ramène aisément, par récurrence sur $s$,
au cas $s=2$ soit à étudier $[\overleftarrow{t-1}]_{\pi_o(1/2)} \boxplus
[\overleftarrow{t-1}]_{\pi_o(-1/2)}$. Or dans le groupe de Grothendieck, cette dernière est égale à
$$[\overleftarrow{t-1}]_{\pi_o(1/2)} \times [\overleftarrow{t-1}]_{\pi_o(-1/2)} - [\overleftarrow{t}]_{\pi_o}
\times [\overleftarrow{t-2}]_{\pi_o}$$ La condition (ii) constitue l'unique restriction sur les $\s$ pour
l'induite totale $[\overleftarrow{t-1}]_{\pi_o(1/2)} \times [\overleftarrow{t-1}]_{\pi_o(-1/2)}$. Il suffit
alors de remarquer que les $\s$ vérifiant (i) et pas (ii) sont obtenues à partir de
$[\overleftarrow{t}]_{\pi_o} \times [\overleftarrow{t-2}]_{\pi_o}$, d'où le résultat.

\end{proof}

Ainsi si on obtenait un constituant de la forme $[\overleftarrow{l-1}]_{\pi_o} \otimes \pi$, on en déduirait,
par transitivité du foncteur de Jacquet, avec les notations du lemme précédent, que la bijection $\s_0$ telle
que $\s_0(i)=\s(1)-i+1$ pour $1 \leq i \leq l$ vérifierait les conditions (i) et (ii) ce qui impose $l \leq
t$.

De même si on obtenait $[\overrightarrow{l-1}]_{\pi_o} \otimes \pi$, la bijection $\s_0$ telle que
$\s_0(i)=\s_0(1)+i-1$ pour $1 \leq i \leq l$ vérifierait les conditions (i) et (ii) ce qui impose $l \leq s$,
d'où le résultat.

\end{proof}

On en déduit alors que pour tout $l> l_0$, $\red_{\JL^{-1}([\overleftarrow{l-1}]_{\pi_o})} (\Pi_o)$ est nul
de sorte que d'après (\ref{somme-alternee}), que pour tout $l_0 < l \leq s$, $\sum_i (-1)^i H^i(j^{\geq lg}_!
HT_{\rho_\oo}(g,l,\pi_o,\Pi_l))[\Pi^\oo]=0$. D'après le lemme (\ref{lem-hij0}), on en déduit alors que
$H^i(j^{\geq lg}_{!*} HT_{\rho_\oo}(g,l,\pi_o,\Pi_l))[\Pi^\oo]$ est nul, pour tout $i$.

\begin{lemb} \label{lem-red2}
Soit $l:=\max \{ s,t \}$. Si $l>t$ (resp. $l>s$), alors
$$\red_{\JL^{-1}([\overleftarrow{l-1}]_{\pi_o})} ([\overrightarrow{s-1}]_{[\overleftarrow{t-1}]_{\pi_o}})$$
est isomorphe à
$$(-1)^{l-1} \Xi^{\frac{(t-1)(sg-1)}{2}} \otimes 
[\overrightarrow{s-1}]_{[\overleftarrow{t-2}]_{\pi_o(\frac{sg-1}{2})}}$$
$$(\hbox{resp. } \Xi^{\frac{(s-1)(tg+1)}{2}} \otimes 
[\overrightarrow{s-2}]_{[\overleftarrow{t-1}]_{\pi_o(\frac{tg+1}{2})}})$$
Dans le cas $l=s=t$, on obtient la somme des deux termes ci-dessus.
\end{lemb}

\begin{proof} L'argument, d'après le lemme (\ref{lem-ze}), est le même que celui de la preuve du lemme
(\ref{lem-red1}).

\end{proof}

Ainsi d'après (\ref{somme-alternee}), on en déduit que
$$\e(\Pi) \sum_i (-1)^i H^i (j^{\geq l_0g}_!
HT_{\rho_\oo}(g,l_0,\pi_o,\Pi_l)[d-l_0g])[\Pi^\oo]$$
qui d'après ce qui précède est égal à $\e(\Pi)\sum_i (-1)^i
H^i(j^{\geq l_0g}_{!*} HT_{\rho_\oo}(g,l,\pi_o,\Pi_l)[d-l_0g])$, est donné par
\begin{multline*}
\Xi^{\frac{d-l_0g}{2}} ((\Xi^{\frac{l_0-1}{2}} (-1)^{l_0-1} \Pi_l (\frac{1-l_0}{2}) \times
([\overrightarrow{l_0-2}]_{[\overleftarrow{l_0-1}]_{\pi_o(1/2)}})^{k_1}+ \\
\Xi^{-\frac{l_0-1}{2}} (\Pi_l (\frac{l_0-1}{2}) \times
[\overrightarrow{l_0-1}]_{[\overleftarrow{l_0-2}]_{\pi_o(-1/2)}})^{k_1}) \times
[\overrightarrow{s_2-1}]_{[\overleftarrow{t_2-1}]_{\pi_{o,2}(\lambda_2)}} \times \cdots \times
[\overrightarrow{s_u-1}]_{[\overleftarrow{t_u-1}]_{\pi_{o,u}(\lambda_u)}} + \\
\sum_{i=k_1+1}^{k_2} \Xi^{-\frac{t_i-1}{2}} (-1)^{l_0-1}
([\overrightarrow{s_1-1}]_{[\overleftarrow{t_1-1}]_{\pi_{o,1} (\lambda_1)}} \times \cdots \times
[\overrightarrow{s_{i-1}-1}]_{[\overleftarrow{t_{i-1}-1}]_{\pi_{o,i-1}(\lambda_{i-1})}} \times \\
\Pi_l (\frac{t_i-1}{2}) \times [\overrightarrow{l_0-1}]_{[\overleftarrow{t_i-2}]_{\pi_o(1/2)}} \times \\
[\overrightarrow{s_{i+1}-1}]_{[\overleftarrow{t_{i+1}-1}]_{\pi_{o,i+1}(\lambda_{i+1})}} \times \cdots \times
[\overrightarrow{s_u-1}]_{[\overleftarrow{t_u-1}]_{\pi_{o,u}(\lambda_u)}}) + \\
\sum_{i=k_2+1}^{k_3} \Xi^{\frac{s_i-1}{2}}
(([\overrightarrow{s_1-1}]_{[\overleftarrow{t_1-1}]_{\pi_{o,1}(\lambda_1)}} \times \cdots \times
[\overrightarrow{s_{i-1}-1}]_{[\overleftarrow{t_{i-1}-1}]_{\pi_{o,i-1}(\lambda_{i-1})}}
\times \\
\Pi_l (\frac{1-s_i}{2}) \times [\overrightarrow{s_i-2}]_{[\overleftarrow{l_0-1}]_{\pi_o(-1/2)}}) \times \\
[\overrightarrow{s_{i+1}-1}]_{[\overleftarrow{t_{i+1}-1}]_{\pi_{o,i+1}(\lambda_{i+1})}} \times \cdots \times
[\overrightarrow{s_u-1}]_{[\overleftarrow{t_u-1}]_{\pi_{o,u}(\lambda_u)}})
\end{multline*}

En utilisant que pour tout $i$, $H^i(j^{\geq l_0g}_{!*} HT_{\rho_\oo}(g,l_0,\pi_o,\Pi_l)[d-l_0g])$ est pur de poids
$d-l_0g+i$, on en déduit alors que ceux-ci sont alors isomorphes à
\begin{multline*} \Pi_l (\frac{1-l_0}{2}) \times
([\overrightarrow{l_0-1}]_{[\overleftarrow{l_0-2}]_{\pi_o(1/2)}})^{k_1} \times
[\overrightarrow{s_2-1}]_{[\overleftarrow{t_2-1}]_{\pi_{o,2}(\lambda_2)}} \times \cdots \times
[\overrightarrow{s_u-1}]_{[\overleftarrow{t_u-1}]_{\pi_{o,u}(\lambda_u)}} \\ \hbox{ pour } i=1-l_0
\end{multline*}
\begin{multline*}
\Pi_l (\frac{l_0-1}{2}) \times ([\overrightarrow{l_0-2}]_{[\overleftarrow{l_0-1}]_{\pi_o(-1/2)}})^{k_1}
\times [\overrightarrow{s_2-1}]_{[\overleftarrow{t_2-1}]_{\pi_{o,2}(\lambda_2)}} \times \cdots \times
[\overrightarrow{s_u-1}]_{[\overleftarrow{t_u-1}]_{\pi_{o,u}(\lambda_u)}} \\
\hbox{ pour } i=l_0-1
\end{multline*}
\begin{multline*}
\bigoplus_{k_1<k \leq k_2~/~1-t_k=i} (-1)^{l_0-1}
([\overrightarrow{s_1-1}]_{[\overleftarrow{t_1-1}]_{\pi_{o,1}(\lambda_1)}} \times \cdots \times
[\overrightarrow{s_{i-1}-1}]_{[\overleftarrow{t_{i-1}-1}]_{\pi_{o,i-1}(\lambda_{i-1})}} \times \\ \Pi_l
(\frac{t_i-1}{2}) \times [\overrightarrow{l_0-1}]_{[\overleftarrow{t_i-2}]_{\pi_o(1/2)}} \times
[\overrightarrow{s_{i+1}-1}]_{[\overleftarrow{t_{i+1}-1}]_{\pi_{o,i+1}(\lambda_{i+1})}} \times \cdots \times
[\overrightarrow{s_u-1}]_{[\overleftarrow{t_u-1}]_{\pi_{o,u}(\lambda_u)}}) \\ \hbox{ pour } 0 \leq -i<l_0-1
\end{multline*}
\begin{multline*}
\bigoplus_{k_2<k \leq k_3 ~/~s_k-1=i}
(([\overrightarrow{s_1-1}]_{[\overleftarrow{t_1-1}]_{\pi_{o,1}(\lambda_1)}} \times
\cdots \times [\overrightarrow{s_{i-1}-1}]_{[\overleftarrow{t_{i-1}-1}]_{\pi_{o,i-1}(\lambda_{i-1})}} \times \\
\Pi_l (\frac{1-s_i}{2}) \times [\overrightarrow{s_i-2}]_{[\overleftarrow{l_0-1}]_{\pi_o(-1/2)}}) \times
[\overrightarrow{s_{i+1}-1}]_{[\overleftarrow{t_{i+1}-1}]_{\pi_{o,i+1}(\lambda_{i+1})}} \times \cdots \times
[\overrightarrow{s_u-1}]_{[\overleftarrow{t_u-1}]_{\pi_{o,u}(\lambda_u)}}) \\ \hbox{ pour } 0 \leq i <l_0-1
\end{multline*}

Si on veut que ces résultats soient compatibles à la dualité de Verdier, il faut alors que pour tout $i$,
$\min \{ s_i,t_i \}=1$. En outre s'il existe $i$ tel que $t_i>1$ alors $\e(\Pi)=1$ et pour tout $j$, $s_j
\equiv 1 \mod 2$. Si pour tout $i$, $t_i=1$ alors tous les $s_i$ ont la même parité donnée par
$(-1)^{s_i-1}=\e(\Pi)$.

Par ailleurs les $H^i(j^{\geq lg}_{!*} HT_{\rho_\oo}(g,l,\pi_o,[\overleftarrow{l-1}]_{\pi_o})[d-lg]) \otimes 
L_g(\pi_o)
(1-lg)$ servent à calculer la cohomologie de la fibre générique via la suite spectrale des poids. On en
déduit alors que tous les poids $L_g(\pi_{o,i}(\lambda_i))$ s'obtiennent à partir de ceux obtenus par la
cohomologie de la fibre générique par une translation entière, d'où le résultat.

\end{proof}

\section{Pureté de la filtration de monodromie des $H^i_{\eta_o}$}

\begin{theob} \label{theo-mono-glob}
Pour toute représentation irréductible $\Pi$ et pour tout $i$, les gradués, $gr_{k,\rho_\oo}[\Pi^\oo]$ de la
filtration de monodromie de $H^i_{\eta_o,\rho_\oo}[\Pi^\oo]$ sont purs de poids $d-1+k$ \footnote{plus le
poids de $\rho_\oo$ que l'on supposera nul}. En fait pour $i \neq d-1$, l'opérateur de monodromie $N$ agit
trivialement sur les $H^i_{\eta_o,\rho_\oo}$.
\end{theob}

\rem Il s'agit donc de prouver la conjecture de monodromie-poids, dite globale, résultat déjà connu en égale
caractéristique d'après Deligne. On ne cherche pas ici à prouver ce résultat dans cette situation, puisqu'à
de nombreuses reprises nous avons déjà utilisé la pureté des gradués de la filtration de monodromie. La
preuve qui suit n'a donc de sens qu'en caractéristique mixte.

\begin{proof} Soit $o$ une place non ramifiée. Remarquons tout d'abord que s'il n'existe pas de $\Pi_o'$
\footnote{On peut en outre supposer que ce dernier soit de même support cuspidal que $\Pi_o$} tel que $\Pi^o
\Pi_o'$ ne vérifie pas $\hyp(\oo)$, alors tous les $H^i_{\eta_o}[\Pi^\oo]$ sont nuls. Considérons alors $\Pi$
irréductible vérifiant $\hyp(\oo)$ de sorte que, d'après ce qui précède, $\Pi_o$ est de la forme
$$[\overrightarrow{s_1-1}]_{\pi_{o,1}} \times \cdots \times [\overrightarrow{s_u-1}]_{\pi_{o,u}} \times
[\overleftarrow{t_1-1}]_{\pi_{o,1}'} \times \cdots \times [\overleftarrow{t_v-1}]_{\pi_{o,v}'}$$ avec
$\pi_{o,i}$ et $\pi_{o,j}'$ irréductibles cuspidales unitaires de respectivement $GL_{g_i}(F_o)$ et
$GL{g_j'}(F_o)$. \footnote{En outre tous les $s_i$ ont la même congruence modulo $2$.} On étudie alors la
suite spectrale
\begin{equation} \label{ss-gr-3}
E_1^{i,j}=H^{i+j}(gr_{-i,\r_\oo}) \Rightarrow H^{i+j+d-1}_{\eta_o,\r_\oo}
\end{equation}
D'après le théorème (\ref{theo-global1}), les $H^r(gr_{k,\r_\oo})[\Pi^{\oo,o}]$ se calculent à partir des
$H^r(\PC(g_i,l,\pi_{o,i}))[\Pi^{\oo,o}]$ pour $1 \leq i \leq u$ et $1 \leq lg_i \leq d$ ainsi que des
$H^r(\PC(g_i',l,\pi_{o,i}'))[\Pi^{\oo,o}]$ pour $1 \leq i \leq u'$ et $1 \leq lg_i' \leq d$. Le calcul de ces
derniers s'obtient par récurrence descendante sur $l$ via le calcul des
$\red_{\JL^{-1}([\overleftarrow{l-1}]_{\pi_{o,i}})}^{lg_i}(\Pi_o)$.

A priori, pour une preuve rigoureuse on devrait introduire des notations comme dans la preuve de la
proposition (\ref{prop-compo-locale}) en regroupant les $\pi_{o,i}$ qui sont isomorphes, puis en les triant
selon leur $s_i$. Cependant comme le foncteur de Jacquet est additif sur les induites et vu que les
$H^i(j^{\geq lg}_{!*} HT(g,l,\pi_o,\Pi_l))[\Pi^{\oo,o}]$ se calculent à partir de ceux-ci, par récurrence, on
remarque alors que les divers $\pi_{o,i}$ n'interagissent pas entre eux de sorte que chacun des $H^i(j^{\geq
lg}_{!*} HT(g,l,\pi_o,\Pi_l)[d-lg]) [\Pi^{\oo,o}]$ sera de la forme
$$\bigoplus_{\genfrac{}{}{0pt}{}{k/\pi_{o,k} \simeq \pi_o}{s_k \geq l}} H(i,\pi_o,l,s_k,-) \oplus
\bigoplus_{\genfrac{}{}{0pt}{}{k/\pi_{o,k}' \simeq \pi_o}{t_k \geq l}} H(i,\pi_o,l,t_k,+)$$ où
$H(i,\pi_o,l,s_k,-)$ (resp. $H(i,\pi_o,l,t_k,+)$) correspond au calcul de
$$H^i(j^{\geq lg}_{!*} HT(g,l,\pi_o,\Pi_l)[d-lg])[\Pi^{\oo,o}]$$
effectué comme si $k$ est tel que pour tout $1 \leq i \neq k \leq u$ et pour tout $1 \leq j \leq u'$ (resp.
pour tout $1 \leq i \leq u$ et $1 \leq j \neq k \leq u'$), $\pi_{o,i}$ et $\pi_{o,j}'$ ne sont pas isomorphes
à $\pi_o$. Plaçons nous alors dans une telle situation, i.e. $\Pi_o \simeq [\overleftarrow{s-1}]_{\pi_o}
\times \pi_{o,1}$ (resp. $\Pi_o \simeq [\overrightarrow{s-1}]_{\pi_o} \times \pi_{o,1}$) où $\pi_{o,1}$ est
une représentation irréductible de $GL_{d-sg}(F_o)$ de support cuspidal disjoint de la droite associée à
$\pi_o$ irréductible cuspidale de $GL_g(F_o)$. Les calculs sont alors exactement similaires à ceux déjà
effectués, rappelons rapidement comment cela s'articule.

- pour $l > s$, les $H^i(j^{\geq lg}_{!*} HT(g,l,\pi_o,\Pi_l))[\Pi^{\oo,o}]$ sont nuls: on procède par
récurrence descendante sur $l$ de $s_g$ à $s+1$. Le cas $l=s_g$ quand $s_g g =d$ découle du lemme
(\ref{lem-pts-ss}) en utilisant le théorème (\ref{theo-jl}). Dans les autres situations le résultat découle
du fait que $\red_{\JL^{-1}([\overleftarrow{l-1}]_{\pi_o})}^{lg} (\Pi_o)$ est nul;

- pour $l=s$, on a $H^0(j^{\geq s'g}_{!*} HT(g,s',\pi_o,[\overleftarrow{s-1}]_{\pi_o})[d-sg])[\Pi^{\oo,o}]
\simeq \Pi_o$ \footnote{Dans le cas respé on utilise que $(-1)^{s-1}=\e(\Pi)$.}: le résultat découle
directement du point précédent et du fait que
$\red_{\JL^{-1}([\overleftarrow{s-1}]_{\pi_o})}^{sg}(\Pi_o)=\pi_{o,1}(sg/2) \otimes \Xi^{\frac{d-sg}{2}}$;

- pour $1 \leq l < s$ et $\Pi_o \simeq [\overleftarrow{s-1}]_{\pi_o} \times \pi_{o,1}$, tous les $H^i(j^{\geq
lg}_{!*} HT(g,l,\pi_o,[\overleftarrow{l-1}]_{\pi_o}))[\Pi^{\oo,o}]$ sont nuls: la preuve est strictement
identique à celle de la proposition (\ref{prop-coho1});

- pour $1 \leq l < s$ et $\Pi_o \simeq [\overrightarrow{s-1}]_{\pi_o} \times \pi_{o,1}$, les $H^i(j^{\geq
lg}_{!*} HT(g,l,\pi_o,[\overleftarrow{l-1}]_{\pi_o})[d-lg])[\Pi^{\oo,o}]$ sont nuls pour $|i| > s-l$ ou $i
\equiv s-l-1 \mod 2$. Dans les autres cas, il est isomorphe à
$$([\overleftarrow{l-1}]_{\pi_o} \overrightarrow{\times} [\overrightarrow{\frac{s-l+i - 4}{2}}]_{\pi_o})
\overleftarrow{\times} [\overrightarrow{\frac{s-l-i-4}{2}}]_{\pi_o} \otimes (\Xi^{\frac{(s-l)g+i}{2}} \otimes
\bigoplus_{\xi \in \AF(\pi_o)} \xi^{-1})$$ La preuve est identique à celle de la proposition
(\ref{prop-not}).

On revient alors à l'étude de la suite spectrale (\ref{ss-gr-3}) dont on vient de déterminer les termes
$E_1^{i,j}[\Pi^{\oo,o}]$. La détermination des flèches $d_1^{i,j}$ se fait alors comme dans la preuve de la
proposition (\ref{prop-lrs2}) car à nouveau les divers $\pi_{o,i}$ n'interagissent pas entre eux dans les
arguments de dualité. Considérons en effet le cas de deux contributions: si celles-ci correspondent à des
Steinberg généralisées alors tout est concentré sur les $H^0$ de sorte que toutes les flèches sont nulles.
Supposons donc avoir à traiter le cas de deux $\speh$ \footnote{Le cas d'une $\speh$ avec une $\st$ se traite
de la même façon.}. Si elles sont de même longueur, on se retrouve alors dans la situation de la preuve de
(\ref{prop-lrs2}) mais avec des multiplicité $2$, ce qui ne change en rien les arguments. Supposons les alors
de longueur distinctes, $l'<l$ et reprenons les arguments de la preuve de (\ref{prop-lrs2}). Tant qu'on
étudie des $H^i(gr_k)$ avec $\min \{ |i|, |k| \}  >l$, il n'y a rien à rajouter car n'intervient que la
composante de plus grande longueur. Dans les autres cas, l'argument consiste à comparer par dualité les
termes $H^i$ et $H^{-i}$. Pour $i>0$, $H^i$ (resp. $H^{-i}$) est une somme directe constituée de
représentations de la forme
\begin{equation} \label{type1}
[\overrightarrow{l-1}]_{\pi_o} \times ([\overrightarrow{r_1-1}]_{\pi_o} \overrightarrow{\times}
[\overleftarrow{l'-r_1-r_2-1}]_{\pi_o} \overrightarrow{\otimes} [\overrightarrow{r_2-1}]_{\pi_o})
\end{equation}
avec $r_2 \leq l'/2$ (resp. $r_1 \leq l'/2$) \footnote{On peut être plus précis, cf. la proposition
(\ref{prop-not}).}
\begin{equation} \label{type2}
([\overrightarrow{r_1-1}]_{\pi_o} \overrightarrow{\times} [\overleftarrow{l-r_1-r_2-1}]_{\pi_o}
\overrightarrow{\times} [\overrightarrow{r_2-1}]_{\pi_o} ) \times [\overrightarrow{l'-1}]_{\pi_o}
\end{equation}
avec $r_2 \leq l/2$ (resp. $r_1 \leq l/2$). On remarque alors que les composants de (\ref{type1}) sont
disjoints de ceux de (\ref{type2}), sauf éventuellement $[\overrightarrow{l-1}]_{\pi_o} \times
[\overrightarrow{l'-1}]_{\pi_o}$. Ce sont donc les seuls qui peuvent rester dans l'aboutissement \footnote{On
le savait déjà car seules les représentations automorphes subsistent dans l'aboutissement.}. Pour voir que
ces derniers restent effectivement dans l'aboutissement il suffit de remarquer qu'ils n'apparaissent que dans
les termes $H^i(gr_0)[\Pi^{\oo,o}]$ de sorte qu'il ne peut y avoir de simplifications, d'où le résultat.

\end{proof}


%% file: appendiceB.tex
\chapter{Cas de la caractéristique mixte}
\label{appendiceB}


Il s'agit dans un premier temps de donner la correspondance entre nos notations et celles de \cite{h-t}.

\begin{figure}[!h]
\begin{center}
\begin{tabular}{|c | c |}
\hline égale caractéristique:  & caractéristique mixte: \cite{h-t} \\ \hline $F$ & $F$ \\ $o$ & $\omega$
\\ $\Psi_{F_o,m}^{d,i}$ & $\Psi_{F_\omega,l,d,m}^i$ \\ $L_g$ & $\rec_{F_\omega}^\vee$
\\ $\cl$ & $\art_{F_\omega}^{-1}$ \\
$D_\Am^\times$ & $G_\tau(\Am)$ \\ $(D_\Am^\oo)^\times$ & $G(\Am^\oo)=G_\tau(\Am^\oo)$ \\ $M_{I,o}$ & $X_U$ \\ 
$M_{I,s_o}$ & $\bar X_U$ \\
$M_{I,s_o}^{=h}$ & $\bar X_U^{(d-h)}$ \\ $P_{h,d}$ & $P_{d-h}$ \\ $M_{I^o \MC_o^n,s_o,1}^{=h}$ & $\bar 
X_{U^p,m,M}^{(d-h)}$ \\
$\FC_{\t_o}$ & $\FC_{\t_o^\vee}$ \\ $(D_{\Am}^{\oo,o})^\times$ & $G(\Am^{\oo,p}) \times \Qm_p^\times \times
\prod_{i=2}^r (B_{\omega_i}^{op})^\times$ \\ $\r_\oo$ & $\xi$ \\ $P_{h,d}$ & $P_h^{\op}$ \\ \hline
\end{tabular}
\end{center}
\end{figure}
où on a
$$G(\Qm_p) \simeq \Qm_p^\times \times GL_d(F_\omega) \times \prod_{i=2}^r (B_{\omega_i}^{op})^\times$$
$$G(\Am^\oo) =G(\Am^{\oo,p}) \times \Qm_p^\times \times GL_d(F_\omega) \times \prod_{i=2}^r
(B_{\omega_i}^{op})^\times$$ L'algèbre $\bar D$ associée aux points supersinguliers est $H_{z}$ (cf.
\cite{h-t} p.153, qui dans \cite{ha} p.33 et p.62 est noté $I_x$) où $z$ est un point de la strate basique
(groupe des quasi-isogénies du triplet $(A,\l,i)$ associé à $z$, cf. loc. cit.).

\medskip

\marque \textit{Notes sur les induites}: dans le but de garder la même combinatoire avec les mêmes notations,
il convient de faire les modifications suivantes aux définitions du \S \ref{defi-induite}.

\begin{itemize}

\item comme précisé plus haut, on considère les paraboliques standards (en caractéristique
mixte) en lieu et place des paraboliques opposés (en égale caractéristique) en particulier dans la
définition de $\red_{\tau_o}^h$;

\item pour $\pi_1$ (resp. $\pi_2$) une représentation de $GL_{n_1}(F_o)$ (resp. $\pi_1 \times \pi_2$ désigne
l'induite
$$\ind_{P_{n_1,n_1+n_2}(F_o)}^{GL_{n_1+n_2}(F_o)} \pi_1(n_2) \otimes \pi_2(-n_1).$$
Dans le cas où $n_1=l_1g$ et $n_2=l_2g$, $\pi_1 \overrightarrow{\times} \pi_2$ (resp. $\pi_1 \overleftarrow{\times} 
\pi_2$)
désigne l'induite normalisée $\pi_1(-l_2) \times \pi_2(l_1)$ (resp. $\pi_1(l_2) \times \pi_2(-l_1)$);

\item soit $g$ un diviseur de $d=sg$ et $\pi_o$ une représentation irréductible cuspidale de
$GL_g(F_o)$. Pour $1 \leq h \leq d$, le foncteur de Jacquet vérifie les propriétés suivantes:
\begin{itemize}
\item si $g$ ne divise pas $h$, alors
$$J_{N_{h,d}^{\op}}([\overleftarrow{s-1}]_{\pi_o}) = J_{N_{h,d}^{\op}}([\overrightarrow{s-1}]_{\pi_o})=(0)$$

\item si $h=lg$ alors
$$J_{N_{lg,sg}^{\op}}([\overleftarrow{s-1}]_{\pi_o})=[\overleftarrow{l-1}]_{\pi_o((l-s)/2)} \otimes
[\overleftarrow{s-l-1}]_{\pi_o(l/2)}$$
$$J_{N_{lg,sg}^{\op}}([\overrightarrow{s-1}]_{\pi_o})=[\overrightarrow{l-1}]_{\pi_o((s-l)/2)} \otimes
[\overrightarrow{s-l-1}]_{\pi_o(-l/2)}$$
\end{itemize}

\item $V(\pi_o,s)$ est ici l'induite $\pi_o(\frac{1-s}{2}) \times \cdots \times \pi_o(\frac{s-1}{2})$. La 
représentation de Steinberg
généralisée $\st_s(\pi_o)$ (resp. $\speh_s(\pi_o)$) est encore notée $[\overleftarrow{s-1}]_{\pi_o}$ (resp. 
$[\overrightarrow{s-1}]_{\pi_o}$)
et de manière générale on reprend les notations de \S \ref{defi-induite} de sorte que par rapport à la remarque 
(\ref{rema-numero}),
les flèches correspondent aux orientations du graphe
$$\circ^{\frac{1-s}{2}} \to \circ \cdots \to \circ^{\frac{s-1}{2}}$$
numéroté donc dans l'ordre inverse par rapport à loc. cit. de sorte
que les orientations sont cette fois-ci compatibles aux foncteurs de Jacquet pour les paraboliques standards. On 
notera par ailleurs que ce l'on note,
par exemple, $[\overleftarrow{l-1},\overrightarrow{s-l}]_{\pi_o}$ dans ce paragraphe, est noté 
$[\overrightarrow{s-l},\overleftarrow{l-1}]_{\pi_o}$
au \S \ref{defi-induite} de sorte que, par exemple, les énoncés locaux en égale caractéristique et en caractéristique 
mixte, bien que noté de la
manière similaire avec nos notations, sont en fait duaux.
\end{itemize}

\bigskip

\marque On a ainsi: $\bar X_{U^p,m}^{(d-h)}=\bar X_{U^p,m,M}^{(d-h)} \times_{P_{h,d}(\OC_o/\MC_o^n)}
GL_d(\OC_o/\MC_o^n)$ et donc
$$\lim_{\genfrac{}{}{0pt}{}{\to}{U}} H^i_c(\bar X_{U^p,m}^{(d-h)},\LC_\xi) = \ind_{P_{h,d}(F_o)}^{GL_d(F_o)}
\lim_{\genfrac{}{}{0pt}{}{\to}{U}} H^i_c(\bar X_{U^p,m,M}^{(d-h)},\LC_\xi)$$ où les $\bar
X_{U^p,m,M}^{(d-h)}$ \footnote{On rappelle que pour tout $U$, il y a un morphisme radiciel $I_{U^p,m}^{(h)}
\longto \bar X_{U^p,m,M}^{(d-h)}$.} sont munis par correspondances d'une action de $\Zm \times GL_{d-h}(F_o)$
de sorte que l'action de $(\left (
\begin{array}{cc} g_o^c & * \\ 0 & g_o^{et} \end{array} \right ) , c_o) \in P_{h,d}(F_o) \times W_o$
est donnée par l'action de $(\val(\det g_o^c)- \deg c_o,g_o^{et}) \in \Zm \times GL_{d-h}(F_o)$.

\bigskip

\marque En ce qui concerne les propriétés cohomologiques associées aux systèmes locaux d'Harris-Taylor du
paragraphe (\ref{ht-prop}), elles sont présentes dans \cite{h-t}, soit:

(a) l'isomorphisme (\ref{iso1}) avec l'action de $G(\Am^\oo)$, est donné au corollaire IV.2.3 et IV.2.4 en
utilisant la formule du bas de la page 138. L'action de $\Zm$ dans l'énoncé de (\ref{iso1}) découle de
l'action naturelle de $\Zm$ sur les variétés d'Igusa de première espèce, cf. loc. cit. p.122;

(b) la proposition (\ref{prop-hic}) correspond au corollaire IV.2.8 de \cite{h-t}, où donc l'action d'un
élément
\begin{multline*}
(g^p,g_{p,0},g^0_\omega,g^{et}_\omega,g_{\omega_i},\s) \in \\ G(\Am^{\oo,p}) \times \Qm_p^\times \times
GL_{n-h}(F_\omega) \times GL_h(F_\omega) \times \prod_{i=2}^r (B_{\omega_i}^{op})^\times \times W_{F_\omega}
\end{multline*}
est donnée par l'action de
$$(g^p,g_{p,0} p^{f_1(w(w(\s)-\det g^0_\omega))},\d,g^{et}_\omega,g_{\omega_i}) \times (g^0_\omega,\s) \in
G^{(h)}(\Am^\oo)/\OC_{D_{F_\omega,n-h}}^\times  \times GL_{d-h}(F_\omega) \times W_{F_\omega}$$ avec
$$G^{(h)}(\Am^\oo)=G(\Am^{\oo,o}) \times \Qm_p^\times \times D_{F_\omega,n-h}^\times \times GL_h(F_\omega)
\times \prod_{i=2}^r (B_{\omega_i}^{op})^\times$$ où $\d \in D_{F_\omega,n-h}^\times$ est tel que $w(\det
\d)=w(\s)-w(\det g^0_\omega)$ et $\OC_{D_{F_\omega,n-h}}$ est l'ordre maximal de $D_{F_\omega,n-h}$;

\rem On notera que les signes opposés dans $w(\s)-w(\det g^0_\omega)$ sont compatibles à l'utilisation des
paraboliques standards $P_{h,d}$ dans les formules d'induction tout au long de la preuve, cf. par exemple la
proposition (\ref{prop-poids}). Les preuves sont par ailleurs strictement similaires en modifiant les calculs
de $\red_{\t_o}^h$. Dans certaines formules où n'apparaît pas les notations $\overleftarrow{\times}$ ou
$\overrightarrow{\times}$ qui intègrent les modifications à effectuer entre égale et inégale caractéristique,
il faut rajouter un signe devant la torsion pour l'action du groupe linéaire; en particulier la formule
(\ref{formule-modif}) doit se modifier comme suit:

\begin{multline*}
\lim_{\genfrac..{0pt}{1}{\to}{I}} H^{d-lg}_c(M_{I,s_o,1}^{=lg},\FC(g,l,\pi_o)_1 \otimes \LC_{\rho_\oo})[\Pi^{\oo,o}]= 
\\
 m(\Pi) \bigoplus_{\xi \in \AF(\pi_o)} [\overleftarrow{s-l-1}]_{\pi_o(-l(g-1)/2)} \otimes \Xi^{(s-l)(g-1)/2}
\end{multline*}
de sorte qu'en tant que $GL_d(F_o) \times W_o$-module, on a
\begin{multline*}
\lim_{\genfrac..{0pt}{1}{\to}{I}}
H^{d-lg}_c(M_{I,s_o}^{=lg},HT_{\rho_\oo}(g,l,\pi_o,[\overleftarrow{l-1}]_{\pi_o}))[\Pi^{\oo,o}]
\otimes \Xi^{\frac{lg-1}{2}}= \\
 m(\Pi)e_{\pi_o} \ind_{P_{lg,d}(F_o)}^{GL_d(F_o)} [\overleftarrow{l-1}]_{\pi_o((s-l)(g-1)/2)} \otimes
[\overleftarrow{s-l-1}]_{\pi_o(-l(g-1)/2)} \otimes \Xi^{\frac{d-s+l-1}{2}} \\
\simeq m(\Pi) e_{\pi_o} [\overleftarrow{l-1}]_{\pi_o} \overrightarrow{\times} [\overleftarrow{s-l-1}]_{\pi_o}
\otimes \Xi^{\frac{d-s+l-1}{2}}
\end{multline*}
\textbf{où l'on notera bien que} $\FC(g,l,\pi_o)$ est le faisceau induit à partir du système local
d'Harris-Taylor, $\FC_{\JL^{-1}(\st_s(\pi_o)^\vee)}$, de même que $\PC(g,l,\pi_o)$ est l'extension
intermédiaire de $\FC(g,l,\pi_o)[d-lg] \otimes [\overleftarrow{l-1}]_{\pi_o} \otimes L_g(\pi_o)^\vee$. De
manière générale il faut aussi remplacer $\JL^{-1}(\bullet)$ par $\JL^{-1}(\bullet)^\vee$.

(c) la condition $\hyp(\oo)$ est la suivante:
$\Pi_\oo$ est cohomologique pour une certaine représentation algébrique $\xi$ sur $\Cm$ de la
restriction des scalaires de $F$ à $\Qm$ de $GL_g$, i.e. il existe $i$ tel que

$$H^i((\lie G_\tau(\Rm)) \otimes_\Rm \Cm,U_\tau,\Pi_\oo \otimes \xi') \neq (0)$$
où $U_\tau$ est un sous-groupe compact modulo le centre de $G_\tau(\Rm)$, maximal, cf. \cite{h-t} p.92, et où $\xi'$
est le caractère sur $\Cm$ associé à $\xi$ via un isomorphisme $\bar \Qm_l \simeq \Cm$ fixé.
On dira que $\Pi$ vérifie $\hyp(\xi)$ si $\Pi$ vérifie $\hyp(\oo)$ pour $\xi$.

(d) l'équivalent de la proposition (\ref{prop-somme-alternee}) correspond à la deuxième identité
fondamentale, théorème V.5.4, où $\Red$ est définie p.182.

\bigskip

\begin{proof}[Retour sur la preuve de la proposition (\ref{prop-mono})]: (d'après une idée de M. Harris) on
rappelle que dans \cite{h-t}, la somme alternée de la cohomologie de la variété globale à valeur dans
$\LC_\xi$, dans le groupe de Grothendieck correspondant, est écrite sous la forme
$$\sum_{\pi^\oo} \pi^\oo \otimes [R_\xi(\pi^\oo)] \qquad [R_\xi(\pi^\oo)]=\sum_i (-1)^i R^i_\xi(\pi^\oo)$$
où $\pi^\oo$ décrit l'ensemble des représentations irréductibles de $G(\Am^\oo)$. Par ailleurs, corollaire
VI.2.7 de loc. cit., s'il existe une représentation irréductible $\pi_\oo$ de $G_\tau({\mathbb R})$ telle que
$\pi:= \pi^\oo \otimes \pi_\oo$ vérifie $\hyp(\oo)$ avec $BC(\pi)=(\psi,\Pi)$ où $\JL(\Pi)$ est cuspidale,
avec les notations de loc. cit., alors $R^i_\xi(\pi^\oo)$ est nul pour $i \neq d-1$.

\noindent Soit alors $\s_\omega=L_g(\pi_\omega)$ et soit $F'_\omega$ une extension de $F_\omega$ telle que la
restriction de $\s_\omega$ à $F'_\omega$ soit non ramifiée. On rappelle que l'extension
$F'_\omega/F_\omega$ est résoluble; cf.\cite{se} IV-2. On globalise alors la situation comme dans \cite{h-t}: soit 
$F=F^+M$ où $M$ est
un corps quadratique imaginaire dans lequel $p$ se décompose et soit $(F')^+/F^+$ une extension résoluble de
corps totalement réels telle que:

-la place $\omega$ de $F^+$ soit inerte dans $(F')^+$;

- l'extension $(F')^+_\omega / F^+_\omega$ est isomorphe à l'extension
$F'_\omega / F_\omega$.

\noindent On pose alors $F'=M(F')^+$. D'après le corollaire VI.2.6 de loc. cit., on peut choisir
une représentation automorphe cuspidale $\Pi$ de $GL_d(\Am_F)$ telle que $\Pi$ vérifie:

- $\Pi^c \simeq \Pi^\vee$;

- $\Pi_\oo$ a le même caractère central qu'une représentation algébrique de $Res^G_\Qm (GL_d)$;

- $\Pi_\omega \simeq \st_s(\pi_\omega \otimes \psi)$ pour un certain caractère $\psi$ de $F'_\omega$.

\noindent Soit alors le changement de base $\Pi':=BC_{F'/F}(\Pi)$; d'après loc. cit. (théorèmes VI.1.1 et
VI.2.9) on associe à $\Pi$ et $\Pi'$ des représentations $\pi$ et $\pi'$ de respectivement $G_\tau(\Am_F)$ et
$G_\tau(\Am_{F'})$ qui vérifient $\hyp(\oo)$ avec $\pi'_\omega \simeq [\overleftarrow{s-1}]_{\zeta} \boxplus
\cdots \boxplus [\overleftarrow{s-1}]_{\zeta}$, pour $\zeta$ un caractère de $F'_o$. Par ailleurs d'après le
théorème VII.1.9 de loc. cit. appliqué aux bonnes places de $F$, et en utilisant le théorème de densité de
Cebotarev, on en déduit que la représentation galoisienne $R_{\xi'}(\pi^{',\oo})$ est isomorphe à
$R_\xi(\pi^\oo)_{|\gal(\bar F'/F')}$ avec
$$R_{\xi}(\pi^{',\oo})_\omega \simeq L_d(\pi'_\omega)^\vee=(\sp_s \otimes \zeta^{-1})^g$$
de sorte que $R_\xi(\pi^\oo) \simeq \sp_s \otimes L_g(\pi_o)^\vee$, d'où le résultat.

\end{proof}


%% file: recap.tex
\section{sur les faisceaux}


\begin{itemize}

\item[(1)] - $R\Psi_{\eta_o}(\LC_{\r_\oo}) \simeq R \Psi_{\eta_o}(\bar \Qm_l) \otimes \LC_{\r_\oo}$;

- $R \Psi_{\eta_o}(\bar \Qm_l)[d-1]$ est un objet de $\FPH(M_{s_o})$;

- $R \Psi_{\eta_o}(\bar \Qm_l)[d-1]=\bigoplus_{\pi_o \in \cusp(d)} R \Psi_{\eta_o,\pi_o}(\bar \Qm_l)[d-1]$
avec les notations de la proposition (\ref{prop-so}), où $\cusp(d)$ désigne l'ensemble des classes
d'équivalence inertielle des représentations irréductibles cuspidales de $GL_g(F_o)$ pour tout $1 \leq g \leq
d$;

- $R \Psi_{\eta_o,\pi_o}(\bar \Qm_l)[d-1]$ est mixte de poids $d-1$; sa filtration de monodromie est égale à
celle par le poids au décalage de $d-1$ près, on note $gr_{k,\pi_o}$ le gradué de poids $k$ de sa filtration
de monodromie de sorte qu'en particulier $N^k: gr_{k,\pi_o} \longmapright{\sim} gr_{-k,\pi_o}$;

- pour $1 \leq g \leq d$ et $\pi_o$ une représentation irréductible cuspidale de $GL_g(F_o)$; pour tout $|k|
< s_g:=\lfloor d/g \rfloor$, dans $\FPH(M_{s_o})$ on a (cf. \ref{theo-global1}):
$$e_{\pi_o} gr_{k,\pi_o} = \bigoplus_{\genfrac{}{}{0pt}{}{|k| < l \leq s_g}{l \equiv k-1 \mod 2}} 
\PC(g,l,\pi_o)(-\frac{lg+k-1}{2})$$
où $\PC(g,l,\pi_o):=j^{\geq lg}_{!*} HT(g,l,\pi_o,[\overleftarrow{l-1}]_{\pi_o})[d-lg] \otimes L_g(\pi_o)$;

- la suite spectrale des poids (\ref{suite-spectrale}) dégénère en $E_2$ et les applications $d_1^{i,j}$ sont
induites par les suites exactes (\ref{suites-exactes}).

\item[(2)] $j^{\geq lg}_{!*} HT(g,l,\pi_o,\Pi_l)[d-lg] \in \FPH(M_{s_o})$ est pur de poids $d-lg$ pour $\pi_o$ 
unitaire;
ses faisceaux de cohomologies $h^i j^{\geq lg}_{!*} HT(g,l,\pi_o,\Pi_l)[d-lg]$ sont, pour $g>1$
\footnote{pour $g=1$ cf. (\ref{theo-global2})}, nuls pour $i$ ne s'écrivant pas sous la forme $lg-d+a(g-1)$
avec $0 \leq a \leq s_g-l$ et pour un tel $i=lg-d+a(g-1)$, isomorphe dans $\FH(M_{s_o})$ à
$$j^{\geq (l+a)g}_!HT(g,l+a,\pi_o,\Pi_l \overrightarrow{\times} [\overrightarrow{a-1}]_{\pi_o}) \otimes
\Xi^{\frac{a(g-1)}{2}}$$

\item[(3)] Pour tout $1 \leq l \leq s_g$, on a, cf. (\ref{prop-p}), l'égalité dans $\GF$
$$j^{\geq lg}_! HT(g,l,\pi_o,\Pi_l)[d-lg] = \sum_{i=l}^{s_g} j^{\geq ig}_{!*} HT(g,i,\pi_o,\Pi_l 
\overrightarrow{\times}
[\overleftarrow{i-l-1}]_{\pi_o})[d-ig] \otimes \Xi^{\frac{(l-i)(g-1)}{2}}$$ Dualement on a

$$Rj^{\geq lg}_* HT(g,l,\pi_o,\Pi_l)[d-lg] = \sum_{i=l}^{s_g} j^{\geq ig}_{!*} HT(g,i,\pi_o,\Pi_l 
\overleftarrow{\times}
[\overleftarrow{i-l-1}]_{\pi_o})[d-ig] \otimes \Xi^{\frac{(l-i)(g+1)}{2}} $$ On en déduit alors que, au moins
pour $g>1$,

\begin{multline*}
R^i j^{\geq lg}_* HT(g,l,\pi_o,\Pi_l)[d-lg]=\bigoplus_{(r,a)~/~i=(l+r+a)g-d-a} \\
 j^{\geq (l+r+a)g}_! HT(g,l+r+a,\pi_o,(\Pi_l \overleftarrow{\times} [\overleftarrow{r-1}]_{\pi_o})
\overrightarrow{\times} [\overrightarrow{a-1}]_{\pi_o}) \otimes \Xi^{((r+a)(g-1)+2r)/2}
\end{multline*}

\end{itemize}

\section{de groupe de cohomologie}


\noindent \textit{- Cas où $\Pi_o$ est de la forme $\st_s(\pi_o)$}:

\begin{itemize}
\item[(1)] d'après la proposition (\ref{prop-coho1}),
$H^i(j^{\geq lg}_{!*} HT_{\rho_\oo}(g,l,\pi_o,\Pi_l))[\Pi^\oo]$ est nul pour tout $i$.

\item[(2)] d'après le corollaire (\ref{coro-hij-nul}),
$H^i(j^{\geq lg}_{!} HT_{\rho_\oo}(g,l,\pi_o,\Pi_l))[\Pi^\oo]$ est donné par
$$\left \{ \begin{array}{cl} 0 & \hbox{si } i \neq 0 \\
m(\Pi) \Pi_l \overrightarrow{\times} [\overleftarrow{s-l-1}]_{\pi_o} \otimes ( \Xi^{\frac{(s-l)(g-1)}{2}}
\otimes \bigoplus_{\xi \in \AF(\pi_o)} \xi^{-1}) & i=0 \end{array} \right.$$ où $\AF(\pi_o)$ est l'ensemble
des caractères $\xi:\Zm \longto \Qm_l^\times$, tels que $\pi_o \otimes \xi^{-1} \circ \val(\det) \simeq
\pi_o$.

\item[(3)] par application de la dualité de Verdier, $H^i(Rj^{\geq lg}_{*} HT_{\rho_\oo}(g,l,\pi_o,\Pi_l))[\Pi^\oo]$
est donné par
$$\left \{ \begin{array}{cl} 0 & \hbox{si } i \neq 0 \\
m(\Pi)\Pi_l \overleftarrow{\times} [\overleftarrow{s-l-1}]_{\pi_o} \otimes (\Xi^{(s-l)(g+1)/2} \otimes
\bigoplus_{\xi \in \AF(\pi_o)} \xi^{-1}) & i=0
\end{array} \right.$$
\end{itemize}

\bigskip

\noindent \textit{- Cas où $\Pi_o$ est de la forme $\speh_s(\pi_o)$}:

\begin{itemize}
\item[(1)] d'après la proposition (\ref{prop-not}), $H^i(j^{\geq lg}_{!*}
HT_{\rho_\oo}(g,l,\pi_o,\Pi_l))[\Pi^\oo]$ est donné par

$$\left \{
\begin{array}{lr} m(\Pi) \Pi_l \overrightarrow{\times} [\overrightarrow{s-l-1}]_{\pi_o} \otimes
(\Xi^{\frac{(s-l)(g-1)}{2}} \otimes \bigoplus_{\xi \in \AF(\pi_o)} \xi^{-1}) & i=l-s \\
m(\Pi) (\Pi_l \overrightarrow{\times} [\overrightarrow{\frac{s-l-i}{2}-1}]_{\pi_o}) \overleftarrow{\times}
[\overrightarrow{\frac{s-l+i}{2}-1}]_{\pi_o} \otimes (\Xi^{\frac{(s-l)g+i}{2}} \otimes \bigoplus_{\xi \in
\AF(\pi_o)} \xi^{-1}) & l-s < i < s-l \\ & i \equiv l-s \mod 2 \\ m(\Pi) \Pi_l \overleftarrow{\times}
[\overrightarrow{s-l-1}]_{\pi_o} \otimes (\Xi^{\frac{(s-l)(g+1)}{2}} \otimes \bigoplus_{\xi \in \AF(\pi_o)}
\xi^{-1}) & i=s-l
\end{array} \right. $$

\item[(2)] d'après le corollaire (\ref{strate-alt-p}),
$H^i(j^{\geq lg}_{!} HT_{\rho_\oo}(g,l,\pi_o,\Pi_l))[\Pi^\oo]$ est donné par
$$\left \{ \begin{array}{cl} 0 & \hbox{si } i \neq s-l \\
m(\Pi) \Pi_l \overleftarrow{\times} [\overrightarrow{s-l-1}]_{\pi_o} \otimes (\Xi^{\frac{(s-l)(g+1)}{2}}
\bigoplus_{\xi \in \AF(\pi_o)} \xi^{-1}) & i=s-l \end{array} \right.$$

\item[(3)] par application de la dualité de Verdier, $H^i(Rj^{\geq lg}_{*} HT_{\rho_\oo}(g,l,\pi_o,\Pi_l))[\Pi^\oo]$ 
est donné par
$$\left \{ \begin{array}{cl} 0 & \hbox{si } i \neq 0 \\
m(\Pi) \Pi_l \overleftarrow{\times} [\overleftarrow{s-l-1}]_{\pi_o} \otimes (\Xi^{(s-l)(g+1)/2} \otimes
\bigoplus_{\xi \in \AF(\pi_o)} \xi^{-1}) & i=0
\end{array} \right.$$
\end{itemize}

\section{sur la suite spectrale des cycles évanescents}


A partir du calcul des groupes de cohomologie à support compact des $HT(g,l,\pi_o,\Pi_l)$, l'isomorphisme
(cf. la proposition (\ref{prop-hic}))
$$H^j_c(M_{I,s_o,1}^{=h}, R^i\Psi_{\eta_o}(\LC_{\r_\oo}))^h \simeq \bigoplus_{\t_o \in \CF_h}
(H^j_c(M_{I,s_o,1}^{=h},\FC_{\t_o} \otimes \LC_{\r_\oo}) \otimes
\widetilde{\UC_{F_o,m}^{h,i}(\t_o)})^{h/e_{\t_o}}$$ permet de calculer les
$$H^j_{=h,i,\pi_o,\r_\oo}[\Pi^{\oo,o}]:=\lim_{\genfrac..{0pt}{1}{\longto}{I}} H^j_c(M_{I,\bar s_o}^{=h},R^i 
\Psi_{\eta_o,I,\pi_o}
(\LC_{\r_\oo}))$$ Les suites spectrales associées à la stratification
\begin{equation}
E_{1,I,\pi_o}^{p,q;i}=H_c^{p+q}(M_{I,\bar s_o}^{=p-1},R^i \Psi_{\eta_o,I,\pi_o}(\bar \Qm_l)) \Rightarrow
H_c^{p+q}(M_{I,\bar s_o},R^i \Psi_{\eta_o,I,\pi_o}(\bar \Qm_l))
\end{equation}
permettent alors de calculer les termes $E_2^{i,j}[\Pi^{\oo,o}]$ de la suite spectrale des cycles
évanescents. On résume les résultats obtenus dans

\bigskip

- \noindent \textit{le cas où $\Pi_o \simeq \st_s(\pi_o)$}: le corollaire (\ref{coro-hic}) donne tout
d'abord:

\begin{coro} Les $H^j_{=h,i,\pi_o,\r_\oo}[\Pi^{\oo,o}]$ vérifient les propriétés suivantes:
\begin{itemize}
    \item[(1)] pour $g$ divisant $d=sg$ et $\Pi_o \simeq [\overleftarrow{s-1}]_{\pi_o}$, les
    $H^j_{=h,i,\pi_o,\r_\oo}[\Pi^{\oo,o}]$ sont nuls pour $h$ qui n'est pas de la forme
    $lg$ avec $1 \leq l \leq s$;

    \item[(2)] pour $g$ divisant $d=sg$ et $\Pi_o \simeq [\overleftarrow{s-1}]_{\pi_o}$, les
    $H^j_{=lg,i,\pi_o,\r_\oo}[\Pi^{\oo,o}]$ sont nuls pour $j \neq d-lg$;

    \item[(3)] pour $g$ divisant $d=sg$ et $\Pi_o \simeq [\overleftarrow{s-1}]_{\pi_o}$, les
    $H^{d-lg}_{=lg,i,\pi_o,\r_\oo}[\Pi^{\oo,o}]$ sont nuls pour $i$ ne vérifiant pas
    $l(g-1) \leq i \leq lg-1$. Si $1 \leq l < s$ et $i=lg-r$ avec $1 \leq r \leq l$, on a
    \begin{multline*} H^{d-lg}_{=lg,lg-r,\pi_o,\r_\oo}[\Pi^{\oo,o}] \simeq  \\ m(\Pi)
    ([\overleftarrow{l-r},\overrightarrow{r-1}]_{\pi_o} \overrightarrow{\times} [\overleftarrow{s-l-1}]_{\pi_o})
    \otimes L_g(\pi_o) (-\frac{s(g-1)-2(r-l)}{2})
    \end{multline*}
    en tant que représentation de $GL_d(F_o) \times W_o$, où $m(\Pi)$
    est la multiplicité de $\Pi$ dans l'espace des formes automorphes.

\end{itemize}
\end{coro}

Ainsi, cf. la preuve de la proposition (\ref{prop-hipsi}), les suites spectrales associées à la
stratification sont triviales, i.e. dégénèrent en $E_1$ de sorte que les aboutissements sont donnés par la
proposition (\ref{prop-hipsi}) dont on rappelle l'énoncé ci-dessous:

\begin{prop}
Pour tout $0 \leq i \leq d-1$, les $H^j(R^i\Psi_{\eta_o,\pi_o}(\bar \Qm_l))[\Pi^{\oo,o}]$, en tant que
$GL_d(F_o) \otimes W_o$-modules, vérifient les propriétés suivantes:
\begin{itemize}
\item[(1)] ils sont nuls si $g$ n'est pas un diviseur de $d$;

\item[(2)] pour $g$ un diviseur de $d=sg$, ils sont nuls si $j$ n'est pas de la forme $d-lg$ pour $1 \leq l
\leq s$;

\item[(3)] pour $g$ un diviseur de $d=sg$ et $j=d-lg$ avec $1 \leq l \leq s$, ils sont nuls si $i$ n'est pas
de la forme $lg-r$ avec $1 \leq r \leq l$;

\item[(4)] pour $g$ un diviseur de $d=sg$ et $1 \leq l < s$,
$H^{d-lg}(R^{lg-r} \Psi_{\eta_o,I,\pi_o}(\bar \Qm_l))[\Pi^{\oo,o}]$ est isomorphe à
$$ m(\Pi) [\overleftarrow{l-r},
\overrightarrow{r-1}]_{\pi_o} \overrightarrow{\times} [\overleftarrow{s-l-1}]_{\pi_o} \otimes L_g(\pi_o)
(-\frac{s(g-1)-2(r-l)}{2})$$

\end{itemize}
\end{prop}

Enfin la suite spectrale des cycles évanescents dégénère en $E_3$ et l'aboutissement
$H^i_{\eta_o}[\Pi^{\oo,o}]$ est donné par le corollaire (\ref{coho-global1}):

\begin{coro}
Les groupes de cohomologie de la fibre générique $H^i_{\eta_o}[\Pi^{\oo,o}]$ sont nuls pour $i \neq d-1$ et
$$H^{d-1}_{\eta_o}[\Pi^{\oo,o}]=[\overleftarrow{s-1}]_{\pi_o} \otimes L_g([\overleftarrow{s-1}]_{\pi_o}) 
(-\frac{sg-1}{2})$$
\end{coro}

\bigskip

- \noindent \textit{Cas où $\Pi_o \simeq \speh_s(\pi_o)$}:

\begin{coro} Les $H^j_{=h,i,\pi_o,\r_\oo}[\Pi^{\oo,o}]$ vérifient les propriétés suivantes:
\begin{itemize}
    \item[(1)] pour $g$ divisant $d=sg$ et $\Pi_o \simeq [\overrightarrow{s-1}]_{\pi_o}$, les
    $H^j_{=h,i,\pi_o,\r_\oo}[\Pi^{\oo,o}]$ sont nuls pour $h$ qui n'est pas de la forme
    $lg$ avec $1 \leq l \leq s$;

    \item[(2)] pour $g$ divisant $d=sg$ et $\Pi_o \simeq [\overrightarrow{s-1}]_{\pi_o}$, les
    $H^j_{=lg,i,\pi_o,\r_\oo}[\Pi^{\oo,o}]$ sont nuls pour $j \neq (s-l)(g+1)$;

    \item[(3)] pour $g$ divisant $d=sg$ et $\Pi_o \simeq [\overrightarrow{s-1}]_{\pi_o}$, les
    $H^{(s-l)(g+1)}_{=lg,i,\pi_o,\r_\oo}[\Pi^{\oo,o}]$ sont nuls pour $i$ ne vérifiant pas
    $l(g-1) \leq i \leq lg-1$. Si $1 \leq l < s$ et $i=lg-r$ avec $1 \leq r \leq l$, on a
    \begin{multline*} H^{d-lg}_{=lg,lg-r,\pi_o,\r_\oo}[\Pi^{\oo,o}] \simeq \\ m(\Pi)
    ([\overleftarrow{l-r},\overrightarrow{r-1}]_{\pi_o} \overleftarrow{\times} [\overrightarrow{s-l-1}]_{\pi_o})
    \otimes L_g(\pi_o) (-\frac{s(g+1)-2r}{2})
    \end{multline*}
    en tant que représentation de $GL_d(F_o) \times W_o$, où $m(\Pi)$
    est la multiplicité de $\Pi$ dans l'espace des formes automorphes.

\end{itemize}
\end{coro}

\begin{proof} Le raisonnement est identique à celui de la preuve de la proposition (\ref{prop-hipsi}). On remarque
alors que pour tout $k\geq 1$, les flèches
$$d_k^{p,q;i}:E_{k,\pi_o}^{p,q;i}[\Pi^{\oo,o}] \longto E_{k,\pi_o}^{p+k,q+k-1}[\Pi^{\oo,o}]$$
des suites spectrales associées à la stratification sont toutes nulles. En effet, d'après le corollaire
(\ref{strate-alt-p}), pour que $E_{k,\pi_o}^{p,q;i}[\Pi^{\oo,o}]$ (resp.
$E_{k,\pi_o}^{p+k,q+k-1;i}[\Pi^{\oo,o}]$) soit non nul, il faut qu'il existe $1 \leq l_1 \leq s$ et $1 \leq
r_1 \leq l_1$ (resp. $1 \leq l_2 \leq s$ et $1 \leq r_2 \leq l_2$) avec
$$(p,q,i)=(l_1g+1,d-2l_1g-1+s-l_1,l_1g-r_1)$$
$$(\hbox{resp. } (p+k,q+k-1,i)=(l_2g+1,d-2l_2g-1+s-l_2,l_2g-r_2)).$$
Ce qui donne $1=(3g+1)(l_2-l_1)$; on voit alors que pour tout $k\geq 1$, $E_{k,\pi_o}^{p,q;i}[\Pi^{\oo,o}]$
et $E_{k,\pi_o}^{p+k,q+k-1}[\Pi^{\oo,o}]$ ne peuvent pas être tous deux non nuls de sorte que
$E_{\oo,I,\pi_o}^{p,q;i}[\Pi^{\oo,o}]=E_{1,I,\pi_o}^{p,q;i}[\Pi^{\oo,o}]$ d'où le résultat d'après le
corollaire précédent.

\end{proof}

Enfin la suite spectrale des cycles évanescents dégénère en $E_3$ et l'aboutissement
$H^i_{\eta_o}[\Pi^{\oo,o}]$ est donnée par la proposition (\ref{prop-lrs2}) dont on rappelle l'énoncé.

\begin{prop} Les groupes de cohomologie de la fibre générique $H^i_{\eta_o,\rho_\oo}[\Pi^{\oo,o}]$
sont nuls pour $|d-1-i| \geq s$ ou pour $i \not \equiv d-s \mod 2$ et sinon
$$H^{d+s-2-2i}_{\eta_o,\r_\oo}[\Pi^{\oo,o}] \simeq [\overrightarrow{s-1}]_{\pi_o} \otimes L_g(\pi_o) 
(-\frac{d+s-2-2i}{2})$$
\end{prop}


%% file: figure7.tex
\begin{figure}[!h]
\setlength{\unitlength}{1.3cm}
\begin{picture}(9,16)(0,-7)
\linethickness{.1pt}

\put(4,1){\line(0,1){7}} \put(0,5){\line(1,0){8}}

\multiput(4,2)(0,1){7}{\circle*{.1}} \multiput(1,5)(1,0){7}{\circle*{.1}}

\put(1,8){\line(1,-1){6}} \multiput(1,8)(2,-2){4}{\circle{.3}}

\multiput(1,8)(1,-1){7}{\circle*{.1}}

\put(2,6){\line(1,-1){4}} \multiput(2,6)(2,-2){3}{\circle{.3}}

\multiput(2,6)(1,-1){5}{\circle*{.1}}

\put(3,4){\line(1,-1){2}} \multiput(3,4)(2,-2){2}{\circle{.3}}

\multiput(3,4)(1,-1){3}{\circle*{.1}}

\put(4,2){\circle{.3}}

\put(4,2){\vector(1,0){1}} \put(5,2){\vector(1,0){1}} \put(6,2){\vector(1,0){1}}

\put(3,4){\vector(1,0){1}} \put(4,4){\vector(1,0){1}}

\put(2,6){\vector(1,0){1}}

\put(4,7.5){\vector(0,1){1} $j$} \put(7,5){\vector(1,0){1} $i$}

\put(0.5,8){\makebox(0,0){$\overleftarrow{3}$}}

\put(1,6){\makebox(0,0){$\overleftarrow{2} \overrightarrow{\times} \overrightarrow{0}$}}

\put(2.5,4.5){\makebox(0,0){$\overleftarrow{1} \overrightarrow{\times} \overrightarrow{1}$}}

\put(3,2){\makebox(0,0){$\overleftarrow{0} \overrightarrow{\times} \overrightarrow{2}$}}

\multiput(3,6.5)(2,-2){3}{\makebox(0,0){$\overleftarrow{3}$}}

\put(5,2.5){\makebox(0,0){$\overleftarrow{1} \overrightarrow{\times} \overrightarrow{1}$}}

\put(6,2.5){\makebox(0,0){$\overleftarrow{2} \overrightarrow{\times} \overrightarrow{0}$}}

\put(4,4.5){\makebox(0,0){$\overleftarrow{2} \overrightarrow{\times} \overrightarrow{0}$}}

\put(7,1){\vector(-1,1){1}} \put(7,1){\vector(-2,1){2}} \put(7,1){\vector(-3,1){3}}
\put(7,1){\vector(0,1){1}}

\put(7,.7){\makebox(0,0){poids=$-3$}}

\put(3,3){\vector(0,1){1}} \put(3,3){\vector(1,1){1}} \put(3,3){\vector(2,1){2}}

\put(3,2.7){\makebox(0,0){poids=$-1$}}

\put(2,5){\vector(0,1){1}} \put(2,5){\vector(1,1){1}} \put(2,4.9){\makebox(0,0){poids=$1$}}

\put(1,7){\vector(0,1){1}} \put(1,6.9){\makebox(0,0){poids=$3$}}

\put(4,-7){\line(0,1){7}} \put(0,-3){\line(1,0){8}}

\multiput(4,-6)(0,1){7}{\circle*{.1}} \multiput(1,-3)(1,0){7}{\circle*{.1}}

\put(1,0){\line(1,-1){6}} \put(1,0){\circle{.3}}

\multiput(1,0)(1,-1){7}{\circle*{.1}}

\put(2,-2){\line(1,-1){4}} \put(2,-2){\circle{.3}}

\multiput(2,-2)(1,-1){5}{\circle*{.1}}

\put(3,-4){\line(1,-1){2}} \put(3,-4){\circle{.3}}

\multiput(3,-4)(1,-1){3}{\circle*{.1}}

\put(4,-6){\circle{.3}}

\put(4,-.5){\vector(0,1){1} $j$} \put(7,-3){\vector(1,0){1} $i$}

\put(0.5,0){\makebox(0,0){$\overleftarrow{3}$}}

\put(1,-2){\makebox(0,0){$[\overleftarrow{2},\overrightarrow{1}]$}}

\put(2.5,-3.5){\makebox(0,0){$[\overleftarrow{1},\overrightarrow{2}]$}}

\put(3,-6){\makebox(0,0){$\overrightarrow{3}$}}

\end{picture}

\caption{\label{figure7} Termes $MLE_\bullet^{i,j}$ pour $s=4$}

\end{figure}

\clearpage



%% file: figure8.tex
\begin{figure}[!h]
\setlength{\unitlength}{1cm}
\begin{picture}(15,20)(0,-10)
\linethickness{.1pt}

\multiput(1,3)(0,1){5}{\line(1,0){8}} \multiput(2,3)(1,0){7}{\line(0,1){4}}

\multiput(5,2)(0,2){4}{\circle{.3}}

\multiput(4,3)(0,2){3}{\circle{.2}} \multiput(6,3)(0,2){3}{\circle{.2}}

\multiput(3,4)(2,0){3}{\circle{.4}} \multiput(3,6)(2,0){3}{\circle{.4}}

\multiput(2,5)(2,0){4}{\circle*{.1}}

\put(8,5){\vector(1,0){1} $i$} \put(5,8){\vector(0,-1){7}} \put(4.7,1){$k$}

\multiput(4.5,1.5)(-1,1){4}{\line(1,1){4}}

\multiput(6.5,8.5)(1,-1){4}{\vector(-1,0){1}}

\put(6.7,8.5){poids=$s(g-1)$} \put(7.7,7.5){poids=$s(g-1)+2$} \put(8.7,6.5){poids=$s(g-1)+4$}
\put(9.4,6.1){$=s(g+1)-4$} \put(9.7,5.5){poids=$s(g+1)-2$}

\qbezier(3,4)(3.5,5)(3,6) \put(3.1,5.8){\vector(-1,2){0}} \qbezier(7,4)(7.5,5)(7,6)
\put(7.1,5.8){\vector(-1,2){0}}

\multiput(3.3,5.6)(1,1){3}{$N$} \multiput(4.3,4.6)(1,1){3}{$N$} \multiput(5.3,3.6)(1,1){3}{$N$}

\qbezier(4,3)(4.5,4)(4,5) \put(4.1,4.8){\vector(-1,2){0}} \qbezier(4,5)(4.5,6)(4,7)
\put(4.1,6.8){\vector(-1,2){0}} \qbezier(6,3)(6.5,4)(6,5) \put(6.1,4.8){\vector(-1,2){0}}
\qbezier(6,5)(6.5,6)(6,7) \put(6.1,6.8){\vector(-1,2){0}}

\qbezier(5,2)(5.5,3)(5,4) \put(5.1,3.8){\vector(-1,2){0}}\qbezier(5,4)(5.5,5)(5,6)
\put(5.1,5.8){\vector(-1,2){0}}\qbezier(5,6)(5.5,7)(5,8) \put(5.1,7.8){\vector(-1,2){0}}

\put(4,8.5){\vector(2,-1){1}} \put(4,8.8){\makebox(0,0){$\overleftarrow{3}$}}

\put(3,8.5){\vector(2,-3){1}} \put(3,8.8){\makebox(0,0){$\overleftarrow{2} \overrightarrow{\times}
\overrightarrow{0}$}}

\put(2,7.5){\vector(2,-3){1}} \put(2,7.7){\makebox(0,0){$\overleftarrow{1} \overrightarrow{\times}
\overrightarrow{1}$}}

\put(1,5.5){\vector(2,-1){1}} \put(1,5.7){\makebox(0,0){$\overleftarrow{0} \overrightarrow{\times}
\overrightarrow{2}$}}

\put(3,1){\vector(1,4){1}} \put(3,.5){\makebox(0,0){$\overleftarrow{2} \overrightarrow{\times}
\overrightarrow{0} \oplus$}} \put(3,0){\makebox(0,0){$(\overleftarrow{0} \overrightarrow{\times}
\overrightarrow{1}) \overleftarrow{\times} \overrightarrow{0}$}}

\put(7,1){\vector(-1,4){1}} \put(7,.5){\makebox(0,0){$\overleftarrow{2} \overleftarrow{\times}
\overrightarrow{0} \oplus$}} \put(7,0){\makebox(0,0){$(\overleftarrow{0} \overrightarrow{\times}
\overrightarrow{0}) \overleftarrow{\times} \overrightarrow{1}$}}

\put(6,2){\vector(-1,2){1}} \put(6,1.7){\makebox(0,0){$\overleftarrow{3} \oplus$}}
\put(6,1.2){\makebox(0,0){$(\overleftarrow{1} \overrightarrow{\times} \overrightarrow{0})
\overleftarrow{\times} \overrightarrow{0}$}}

\put(8,2){\vector(-1,2){1}} \put(8,1.7){\makebox(0,0){$\overleftarrow{1} \overleftarrow{\times}
\overrightarrow{1}$}}

\put(9,4.5){\vector(-2,1){1}} \put(9.7,4.5){\makebox(0,0){$\overleftarrow{0} \overleftarrow{\times}
\overrightarrow{2}$}}

\multiput(1,-7)(0,1){5}{\line(1,0){8}} \multiput(2,-7)(1,0){7}{\line(0,1){4}}

\put(3,-2){\vector(-1,-3){1}} \put(3,-1.5){\makebox(0,0){$\overrightarrow{3}$ poids=s(g-1)}}

\put(3,-8){\vector(1,3){1}} \put(3,-8.5){\makebox(0,0){$\overrightarrow{3}$ poids=s(g-1)+2}}

\put(7,-2){\vector(-1,-3){1}} \put(7,-1.5){\makebox(0,0){$\overrightarrow{3}$ poids=s(g-1)+4}}

\put(7,-8){\vector(1,3){1}} \put(7,-8.5){\makebox(0,0){$\overrightarrow{3}$ poids=s(g+1)-2}}

\multiput(2,-5)(2,0){4}{\circle*{.1}}

\put(8,-5){\vector(1,0){1} $i$} \put(5,-2){\vector(0,-1){7}} \put(4.7,-9){$k$}

\end{picture}

\caption{\label{figure8} $H^i(gr_{k,\r_\oo})[\Pi^{\oo,o}]$ pour $\Pi_{\oo} \simeq \speh_4(\pi_\oo)$ et $\Pi_o
\simeq \speh_4(\pi_o)$ avec $\pi_\oo$ (resp. $\pi_o$) une représentation cuspidale de $GL_g(F_\oo)$ (resp.
$GL_g(F_o)$). Le dessin du bas représente le terme $E_2$ et donc l'aboutissement de la suite spectrale.}

\end{figure}


%% file: figure5.tex
\setlength{\unitlength}{1cm}
\begin{picture}(9,19)(0,-10)
\linethickness{.1pt}

\put(0,0){\vector(1,0){9}}
\put(0,0){\vector(0,1){9}}
\put(0,0){\circle*{.1}}

\put(0,1){\line(1,-1){1}}
\multiput(0,1)(1,-1){2}{\circle*{.1}}

\put(0,2){\line(1,-1){2}}
\multiput(0,2)(1,-1){3}{\circle*{.1}}

\put(0,3){\line(1,-1){3}}
\multiput(0,3)(1,-1){4}{\circle*{.1}}

\put(0,4){\line(1,-1){4}}
\multiput(0,4)(1,-1){5}{\circle*{.1}}

\put(0,5){\line(1,-1){5}}
\multiput(0,5)(1,-1){6}{\circle*{.1}}

\put(0,6){\line(1,-1){6}}
\multiput(0,6)(1,-1){7}{\circle*{.1}}

\put(0,7){\line(1,-1){7}}
\multiput(0,7)(1,-1){8}{\circle*{.1}}

\put(1,8){\vector(-1,-1){1}}
\put(1,8){\vector(-1,-2){1}}
\put(1,8){\vector(-1,-3){1}}
\put(1,8){\vector(-1,-4){1}}
\put(2,8.5){\makebox(0,0){$h=4 \times 2$}}

\put(0,7){\circle{.3}} \put(-.5,7){\makebox(0,0){$\genfrac{}{}{0pt}{}{\overleftarrow{3}}{10}$}}

\put(0,6){\circle{.2}} \put(-.5,6){\makebox(0,0){$\genfrac{}{}{0pt}{}{[\overleftarrow{2},\overrightarrow{1}]}
{8}$}} \put(0,6){\vector(2,-1){2}}

\put(0,5){\circle{.2}}
\put(-.5,5){\makebox(0,0){$\genfrac{}{}{0pt}{}{[\overleftarrow{1},\overrightarrow{2}]}{6}$}}
\put(0,5){\vector(2,-1){2}}

\put(0,4){\circle{.2}} \put(-.5,4){\makebox(0,0){$\genfrac{}{}{0pt}{}{\overrightarrow{3}}{4}$}}
\put(0,4){\vector(2,-1){2}}

\put(3,6){\vector(-1,-1){1}}
\put(3,6){\vector(-1,-2){1}}
\put(3,6){\vector(-1,-3){1}}
\put(3,6.5){\makebox(0,0){$h=3 \times 2$}}

\put(2,5){\line(0,-1){2}} \put(2,5){\circle{.3}}
\put(1.5,3){\makebox(0,0){$\genfrac{}{}{0pt}{}{\overrightarrow{2} \overrightarrow{\times}
\overrightarrow{0}]}{4}$}}

\put(2,4){\circle{.2}} \put(1.5,4){\makebox(0,0){$\genfrac{}{}{0pt}{}{[\overleftarrow{1},\overrightarrow{1}]
\overrightarrow{\times} \overrightarrow{0}}{6}$}} \put(2,4){\vector(2,-1){2}}

\put(2,3){\circle{.2}} \put(1.5,5){\makebox(0,0){$\genfrac{}{}{0pt}{}{\overleftarrow{2}
\overrightarrow{\times} \overrightarrow{0}}{8}$}} \put(2,3){\vector(2,-1){2}}

\put(5,4){\vector(-1,-1){1}}
\put(5,4){\vector(-1,-2){1}}
\put(5,4.5){\makebox(0,0){$h=2 \times 2$}}

\put(4,3){\line(0,-1){1}} \put(4,3){\circle{.3}}
\put(4.7,3){\makebox(0,0){$\genfrac{}{}{0pt}{}{\overleftarrow{1} \overrightarrow{\times}
\overleftarrow{1}}{6}$}}

\put(4,2){\circle{.2}} \put(4.7,2){\makebox(0,0){$\genfrac{}{}{0pt}{}{\overrightarrow{1}
\overrightarrow{\times} \overleftarrow{1}}{4}$}} \put(4,2){\vector(2,-1){2}}

\put(7,2){\vector(-1,-1){1}}
\put(7,2.5){\makebox(0,0){$h=1 \times 2$}}

\put(6,1){\circle{.3}} \put(6.5,1){\makebox(0,0){$\genfrac{}{}{0pt}{}{\overleftarrow{0}
\overrightarrow{\times} \overleftarrow{2}}{4}$}}

\put(0,-10){\vector(1,0){9}} \put(0,-10){\vector(0,1){9}} \put(0,-10){\circle*{.1}}

\put(0,-9){\line(1,-1){1}} \multiput(0,-9)(1,-1){2}{\circle*{.1}}

\put(0,-8){\line(1,-1){2}} \multiput(0,-8)(1,-1){3}{\circle*{.1}}

\put(0,-7){\line(1,-1){3}} \multiput(0,-7)(1,-1){4}{\circle*{.1}}

\put(0,-6){\line(1,-1){4}} \multiput(0,-6)(1,-1){5}{\circle*{.1}}

\put(0,-5){\line(1,-1){5}} \multiput(0,-5)(1,-1){6}{\circle*{.1}}

\put(0,-4){\line(1,-1){6}} \multiput(0,-4)(1,-1){7}{\circle*{.1}}

\put(0,-3){\line(1,-1){7}} \multiput(0,-3)(1,-1){8}{\circle*{.1}}

\multiput(0,-3)(2,-2){4}{\circle{.3}}

\put(1,-2){\vector(-1,-1){1}} \put(2,-1.5){\makebox(0,0){$\overleftarrow{3}$ poids=s(g+1)-2}}

\put(3,-4){\vector(-1,-1){1}} \put(3,-3.5){\makebox(0,0){$\overleftarrow{3}$ poids=s(g+1)-4}}

\put(5,-6){\vector(-1,-1){1}} \put(5,-5.5){\makebox(0,0){$\overleftarrow{3}$ poids=s(g-1)+2}}

\put(7,-8){\vector(-1,-1){1}} \put(7,-7.5){\makebox(0,0){$\overleftarrow{3}$ poids=s(g-1)}}

\end{picture}

%% file: figure6.tex
\setlength{\unitlength}{1cm}
\begin{picture}(9,19)(0,-10)
\linethickness{.1pt}

\put(0,0){\vector(1,0){11}}
\put(0,0){\vector(0,1){9}}
\put(0,0){\circle*{.1}}

\put(0,1){\line(1,-1){1}}
\multiput(0,1)(1,-1){2}{\circle*{.1}}

\put(0,2){\line(1,-1){2}}
\multiput(0,2)(1,-1){3}{\circle*{.1}}

\put(0,3){\line(1,-1){3}}
\multiput(0,3)(1,-1){4}{\circle*{.1}}

\put(0,4){\line(1,-1){4}}
\multiput(0,4)(1,-1){5}{\circle*{.1}}

\put(0,5){\line(1,-1){5}}
\multiput(0,5)(1,-1){6}{\circle*{.1}}

\put(0,6){\line(1,-1){6}}
\multiput(0,6)(1,-1){7}{\circle*{.1}}

\put(0,7){\line(1,-1){7}}
\multiput(0,7)(1,-1){8}{\circle*{.1}}

\put(1,8){\vector(-1,-1){1}}
\put(1,8){\vector(-1,-2){1}}
\put(1,8){\vector(-1,-3){1}}
\put(1,8){\vector(-1,-4){1}}
\put(2,8.5){\makebox(0,0){$h=4 \times 2$}}

\put(0,7){\circle{.2}} \put(-.5,7){\makebox(0,0){$\genfrac{}{}{0pt}{}{\overleftarrow{3}}{10}$}}
\put(0,7){\vector(3,-2){3}}

\put(0,6){\circle{.2}}
\put(-.5,6){\makebox(0,0){$\genfrac{}{}{0pt}{}{[\overleftarrow{2},\overrightarrow{1}]}{8}$}}
\put(0,6){\vector(3,-2){3}}

\put(0,5){\circle{.2}}
\put(-.5,5){\makebox(0,0){$\genfrac{}{}{0pt}{}{[\overleftarrow{1},\overrightarrow{2}]}{6}$}}
\put(0,5){\vector(3,-2){3}}

\put(0,4){\circle{.3}} \put(-.5,4){\makebox(0,0){$\genfrac{}{}{0pt}{}{\overrightarrow{3}}{4}$}}

\put(4,6){\vector(-1,-1){1}} \put(4,6){\vector(-1,-2){1}} \put(4,6){\vector(-1,-3){1}} \put(4,6.5){\makebox(0,0){$h=3
\times 2$}}

\put(3,5){\line(0,-1){2}} \put(3,5){\circle{.2}}
\put(2.5,5){\makebox(0,0){$\genfrac{}{}{0pt}{}{\overleftarrow{2} \overleftarrow{\times}
\overrightarrow{0}}{10}$}} \put(3,5){\vector(3,-2){3}} \put(3,5){\line(1,-1){5}}
\multiput(3,5)(1,-1){6}{\circle*{.1}}

\put(3,4){\circle{.2}} \put(2.5,4){\makebox(0,0){$\genfrac{}{}{0pt}{}{[\overleftarrow{1},\overrightarrow{1}]
\overleftarrow{\times} \overleftarrow{0}}{8}$}} \put(3,4){\vector(3,-2){3}}

\put(3,3){\circle{.3}} \put(2.5,3){\makebox(0,0){$\genfrac{}{}{0pt}{}{\overrightarrow{2}
\overleftarrow{\times} \overleftarrow{0}}{6}$}}

\put(7,4){\vector(-1,-1){1}}
\put(7,4){\vector(-1,-2){1}}
\put(7,4.5){\makebox(0,0){$h=2 \times 2$}}

\put(6,3){\line(0,-1){1}} \put(6,3){\circle{.2}}
\put(6.5,3){\makebox(0,0){$\genfrac{}{}{0pt}{}{\overleftarrow{1} \overleftarrow{\times} \overrightarrow{1}}
{10}$}} \put(6,3){\vector(3,-2){3}} \put(6,3){\line(1,-1){3}} \multiput(6,3)(1,-1){4}{\circle*{.1}}

\put(6,2){\circle{.3}} \put(6.5,2){\makebox(0,0){$\genfrac{}{}{0pt}{}{\overrightarrow{1}
\overleftarrow{\times} \overrightarrow{1}}{8}$}}

\put(9,2){\vector(0,-1){1}}
\put(9,2.5){\makebox(0,0){$h=1 \times 2$}}

\put(9,1){\circle{.3}} \put(9.5,1){\makebox(0,0){$\genfrac{}{}{0pt}{}{\overleftarrow{0}
\overleftarrow{\times} \overrightarrow{2}}{10}$}} \put(9,1){\line(1,-1){1}}
\multiput(9,1)(1,-1){2}{\circle*{.1}}

\put(0,-10){\vector(1,0){11}} \put(0,-10){\vector(0,1){9}} \put(0,-10){\circle*{.1}}

\put(0,-9){\line(1,-1){1}} \multiput(0,-9)(1,-1){2}{\circle*{.1}}

\put(0,-8){\line(1,-1){2}} \multiput(0,-8)(1,-1){3}{\circle*{.1}}

\put(0,-7){\line(1,-1){3}} \multiput(0,-7)(1,-1){4}{\circle*{.1}}

\put(0,-6){\line(1,-1){4}} \multiput(0,-6)(1,-1){5}{\circle*{.1}}

\put(0,-5){\line(1,-1){5}} \multiput(0,-5)(1,-1){6}{\circle*{.1}}

\put(0,-4){\line(1,-1){6}} \multiput(0,-4)(1,-1){7}{\circle*{.1}}

\put(0,-3){\line(1,-1){7}} \multiput(0,-3)(1,-1){8}{\circle*{.1}}

\put(1,-2){\vector(-1,-4){1}} \put(2,-1.5){\makebox(0,0){$\overrightarrow{3}$ poids=s(g-1)}}

\put(4,-4){\vector(-1,-3){1}} \put(4,-3.5){\makebox(0,0){$\overrightarrow{3}$ poids=s(g-1)+2}}

\multiput(0,-6)(3,-1){4}{\circle{.3}}

\put(3,-5){\line(1,-1){5}} \multiput(3,-5)(1,-1){6}{\circle*{.1}}

\put(7,-6){\vector(-1,-2){1}} \put(7,-5.5){\makebox(0,0){$\overrightarrow{3}$ poids=s(g+1)-4}}

\put(6,-7){\line(1,-1){3}} \multiput(6,-7)(1,-1){4}{\circle*{.1}}

\put(9,-8){\vector(0,-1){1}} \put(9,-7.5){\makebox(0,0){$\overrightarrow{3}$ poids=s(g+1)-2}}

\put(9,-9){\line(1,-1){1}} \multiput(9,-9)(1,-1){2}{\circle*{.1}}

\end{picture}

%% file: figure1.tex
\setlength{\unitlength}{.5cm}
\begin{picture}(17,9)(-1,-1)
\linethickness{.1pt} \put(0,7){\circle{.2}}
\put(-.5,7.5){\makebox(0,0){$\genfrac{}{}{0pt}{}{\overleftarrow{3}}{0}$}}

\put(6,7){\circle*{.2}} \put(5.5,6){\makebox(0,0){$\genfrac{}{}{0pt}{}{\overleftarrow{3}
\overrightarrow{\times} \overrightarrow{0}}{0}$}}

\put(7,7){\circle{.2}} \put(7,7){\line(1,-1){7}}
\put(7.5,7.5){\makebox(0,0){$\genfrac{}{}{0pt}{}{\overleftarrow{3} \overleftarrow{\times}
\overrightarrow{0}}{1}$}}

\put(14,1){\circle*{.2}} \put(14.5,1.5){\makebox(0,0){$\genfrac{}{}{0pt}{}{\overleftarrow{3}
\overleftarrow{\times} \overrightarrow{0}}{0}$}}

\put(14,0){\circle*{.1}} \put(14,0){\line(0,1){1}} \put(14,-1){\makebox(0,0){$(0,0)$}}
\put(0,6){\makebox(0,0){$(-14,7)$}} \multiput(0,7)(1,0){8}{\circle*{.1}}
\multiput(7,7)(0,-1){8}{\circle*{.1}} \multiput(7,0)(1,0){8}{\circle*{.1}}

\qbezier(0,7)(3,8)(6,7)

\end{picture}

%% file: figure2.tex
\setlength{\unitlength}{.5cm}
\begin{picture}(30,16)(-1,1)
\linethickness{.1pt}

\put(0,14){\circle{.2}} \put(-.5,14.5){\makebox(0,0){$\genfrac{}{}{0pt}{}{\overleftarrow{2}}{0}$}}
\qbezier(0,14)(3,15)(6,14)

\put(6,14){\circle*{.2}} \put(5.5,13){\makebox(0,0){$\genfrac{}{}{0pt}{}{0}{\overleftarrow{2}
\overrightarrow{\times} \overrightarrow{0}}$}}

\put(7,14){\circle{.2}} \put(7,14){\line(1,-1){7}}
\put(7,15){\makebox(0,0){$\genfrac{}{}{0pt}{}{\overleftarrow{2} \overleftarrow{\times}
\overrightarrow{0}}{1}$}}

\qbezier(6,14)(9,13)(12,14)

\put(12,14){\circle*{.2}} \put(11.5,13){\makebox(0,0){$\genfrac{}{}{0pt}{}{0}{\overleftarrow{2}
\overrightarrow{\times} \overrightarrow{1}}$}}

\qbezier(7,14)(10,15)(13,14) \put(13,14){\circle*{.2}}
\put(13,15){\makebox(0,0){$\genfrac{}{}{0pt}{}{(\overleftarrow{2} \overleftarrow{\times} \overrightarrow{0})
\overrightarrow{\times} \overrightarrow{0}}{1}$}}

\put(14,14){\circle{.2}} \put(16,14){\makebox(0,0){$\genfrac{}{}{0pt}{}{\overleftarrow{2}
\overleftarrow{\times} \overleftarrow{1}}{2}$}} \put(14,14){\line(1,-1){14}}

\multiput(0,14)(1,0){15}{\circle*{.1}}
\multiput(7,14)(0,-1){8}{\circle*{.1}}
\multiput(7,7)(1,0){15}{\circle*{.1}}

\put(14,8){\circle*{.2}} \put(14,9){\makebox(0,0){$\genfrac{}{}{0pt}{}{\overleftarrow{2}
\overleftarrow{\times} \overrightarrow{0}}{0}$}}

\qbezier(14,8)(17,9)(20,8)

\put(20,8){\circle*{.2}} \put(19.5,7){\makebox(0,0){$\genfrac{}{}{0pt}{}{0}{(\overleftarrow{2}
\overleftarrow{\times} \overrightarrow{0}) \overrightarrow{\times} \overrightarrow{0}}$}}

\put(21,8){\circle{.2}} \put(21.5,9){\makebox(0,0){$\genfrac{}{}{0pt}{}{(\overleftarrow{2}
\overleftarrow{\times} \overleftarrow{0}) \overleftarrow{\times} \overleftarrow{0}}{1}$}}
\put(21,8){\line(1,-1){7}}

\multiput(14,8)(1,0){8}{\circle*{.1}}
\multiput(21,7)(0,-1){8}{\circle*{.1}}
\multiput(21,0)(1,0){8}{\circle*{.1}}

\put(14,7){\line(0,1){1}}
\put(28,0){\line(0,1){2}}

\put(28,2){\circle{.2}} \put(28,3){\makebox(0,0){$\genfrac{}{}{0pt}{}{\overleftarrow{2}
\overleftarrow{\times} \overrightarrow{1}}{0}$}}

\put(28,1){\circle{.1}}

\put(14,6){\makebox(0,0){$(-14,7)$}}
\put(0,13){\makebox(0,0){$(-28,14)$}}
\put(28,-1){\makebox(0,0){$(0,0)$}}

\end{picture}